\documentclass[11pt]{scrartcl}

\usepackage[utf8]{inputenc}

\usepackage{config}
\usepackage{titling}

\usepackage{geometry}
\usepackage[inline]{enumitem}

\newcommand{\descriptionfont}[1]{{\bfseries\sffamily #1}}

\setlist[itemize]{topsep=0pt,partopsep=0pt,itemsep=0pt,parsep=0pt}
\setlist[itemize,1]{label={\small\textbullet}}
\setlist[itemize,2]{label={\tiny\textbullet}}
\setlist[itemize,3]{label=$\cdot$}
\setlist[enumerate]{topsep=0pt,partopsep=0pt,itemsep=0pt,parsep=0pt}
\setlist[enumerate,1]{label=\roman*)}
\setlist[enumerate,2]{label=\alph*)}
\setlist[enumerate,3]{label=\arabic*)}

\hypersetup{
	colorlinks=true,
	linkcolor=AO!65!black,
	citecolor=AO!65!black,
	urlcolor=AppleGreen!65!black,
	bookmarksopen=true,
	bookmarksnumbered,
	bookmarksopenlevel=2,
	bookmarksdepth=3
      }

\newcommand{\SC}{\omega}

\title{Packing even directed circuits quarter-integrally}
\predate{}
\date{}
\postdate{}

\preauthor{}
\DeclareRobustCommand{\authorthing}{
	\begin{center}
		Maximilian Gorsky\footnote{\href{mailto:m.gorsky@pm.me}{m.gorsky@pm.me}} \\
		{\small Technische Universit\"at Berlin, Germany} \\
		  \medskip
		Ken-ichi Kawarabayashi\thanks{Supported by the JST ERATO Kawarabayashi Large Graph Project JPMJER1201 and JSPS KAKENHI JP18H05291.}~~\!\!\footnote{\href{mailto:k\_keniti@nii.ac.jp}{k\_keniti@nii.ac.jp}} \\
		{\small National Institute of Informatics, Japan} \\
		\medskip
		Stephan Kreutzer\footnote{\href{mailto:sk@stephankreutzer.net}{sk@stephankreutzer.net}} \\
		{\small Technische Universit\"at Berlin, Germany} \\
		\medskip
		Sebastian Wiederrecht\thanks{Supported by the Institute for Basic Science (IBS-R029-C1).}~~\!\footnote{\href{mailto:sebastian.wiederrecht@gmail.com}{sebastian.wiederrecht@gmail.com}} \\
		{\small Institute for Basic Science, South Korea} \\
		
\end{center}}
\author{\authorthing}
\postauthor{}

\begin{document}
\maketitle

\begin{abstract}
We prove the existence of a computable function $f\colon\mathbb{N}\to\mathbb{N}$ such that for every integer $k$ and every digraph $D$ either contains a collection $\mathcal{C}$ of $k$ directed cycles of even length such that no vertex of $D$ belongs to more than four cycles in $\mathcal{C}$, or there exists a set $S\subseteq V(D)$ of size at most $f(k)$ such that $D-S$ has no directed cycle of even length.
Moreover, we provide an algorithm that finds one of the two outcomes of this statement in time $g(k)n^{\mathcal{O}(1)}$ for some computable function $g\colon \mathbb{N}\to\mathbb{N}$.

Our result unites two deep fields of research from the algorithmic theory for digraphs: The study of the \textit{Erd\H{o}s-P{\'o}sa property} of digraphs and the study of the \textit{Even Dicycle Problem}.
The latter is the decision problem which asks if a given digraph contains an even dicycle and can be traced back to a question of P\'olya from 1913.
It remained open until a polynomial time algorithm was finally found by Robertson, Seymour, and Thomas (Ann.\ of Math.\ (2) 1999) and, independently, McCuaig (Electron.\ J.\ Combin.\ 2004; announced jointly at STOC 1997).
The Even Dicycle Problem is equivalent to the recognition problem of \textit{Pfaffian} bipartite graphs and has applications even beyond discrete mathematics and theoretical computer science.
On the other hand, Younger's Conjecture (1973), states that dicycles have the Erd\H{o}s-P\'osa property.
The conjecture was proven more than two decades later by Reed, Robertson, Seymour, and Thomas (Combinatorica 1996) and opened the path for structural digraph theory as well as the algorithmic study of the directed feedback vertex set problem.
Our approach builds upon the techniques used to resolve both problems and combines them into a powerful structural theorem that yields further algorithmic applications for other prominent problems.
\end{abstract}
\let\sc\itshape
\thispagestyle{empty}

\newpage

\setcounter{page}{1}

\section{Introduction} \label{sec:intro}

Cycles are fundamental objects in combinatorics.
So much so, that their study, in one form or another, is integral to almost all aspects of graph theory.
Two particularly interesting properties of cycles are the duality between \textit{hitting} and \textit{packing} cycles as proven by Erd\H{o}s and P\'osa \cite{ErdosPosa1965OnIndependetCircuits} and the structural implications of forbidding cycles with a fixed parity.
Studying the combination of both, namely a half-integral version of the \textit{Erd\H{o}s-P\'osa property}\footnote{A collection $\mathcal{F}$ of subgraphs of a graph $G$ is said to be a \emph{packing} if no two members of $\mathcal{F}$ share a vertex and $\mathcal{F}$ is called a \emph{half-integral packing} if every vertex of $G$ belongs to at most two graphs in $\mathcal{F}$.
We say that a class of graphs $\mathcal{C}$ has the \emph{(half-integral) Erd\H{o}s-P\'osa property} if there exists a function $f\colon\mathbb{N}\to\mathbb{N}$ such that for every $k\in\N$ and every graph $G$ it holds that, if $G-S$ contains some member of $\mathcal{C}$ as a subgraph for every $S\subseteq V(G)$ with $|S|\leq f(k)$, then $G$ contains a (half-integral) packing of members from $\mathcal{C}$.} for odd cycles has lead to a variety of structural and algorithmic results \cite{Reed2004FindingOddCycleTransversals,KawarabayashiReed2010OddCyclePacking,KawarabayashiReed2010AlmostLinearTimeOddCycleTransversal,Kakimura2012ErdosPosaParityConstratins,Joos2016ParityLinkage,Huynh2019UnifiedErdosPosaConstrained,Conforti2020ExtendedFormulationsStableSet,Conforti2020StableSetGenusOddCycles,Jansen2021VertexDeletingAndEvenLess,gollin2021unified,gollin2022unified,ThomasYoo2023PackingCyclesUndirected}, spearheaded by Reed's seminal proof for the half-integral Erd\H{o}s-P\'osa property of odd cycles \cite{Reed1999Fruitsalad}.
Notice that the focus lies on packing and hitting \textit{odd} cycles in particular.
One reason for this can be found in the comparably simple structure of graphs without even cycles; every $2$-connected graph without an even cycle must be an odd cycle\footnote{This property of even cycles is sometimes referred to as the \emph{theta-property} \cite{Zaslavsky1989BiasedGraphsI}.}.
This implies that large \textit{walls} contain many pairwise vertex-disjoint even cycles and thus, by the \textit{Grid Theorem} of Robertson and Seymour \cite{RobertsoSeymour1986GraphMinorsV}, any graph of sufficiently large treewidth contains many vertex-disjoint even cycles.
As a consequence, even cycles have the Erd\H{o}s-P\'osa property and the related computational problems are relatively straight forward.

\vspace{-11pt}
\paragraph{Digraphs}
Originally conjectured by Younger \cite{younger1973graphs}, a directed analogue of the original result by Erd\H{o}s and P\'osa was proven by Reed, Robertson, Seymour, and Thomas \cite{Reed1996PackingDirectedCirctuis}.
The way Reed et al.\@ prove their theorem strongly resembles similar strategies which have emerged from the theory of graph minors for undirected graphs.
Indeed, the insights and techniques gained from their work on Younger's conjecture lead to the inception of \textit{directed treewidth} by Johnson, Robertson, Seymour, and Thomas \cite{Johnson2001DirectedTreewidth} and provided key ideas for the proof of the directed grid theorem by Kawarabayashi and Kreutzer \cite{KawarabayashiKreutzer2015DirectedGrid}.
Thus in the directed setting studying the Erd\H{o}s-P\'osa property laid the groundwork for structure theory in directed graphs.
Several variants of the Erd\H{o}s-P\'osa property for dicycles have been proven around the time Younger's conjecture was resolved and continue to be studied to this day, further refining the techniques involving directed treewidth \cite{Seymour1995PackingDirectedCircuitsFractionally,Seymour1996PackingDirectectedCircuitsEulerian,GueninThomas2011PackingDirectedCircuits,KakimuraKawarabayashi2012PackingDirectedCircuits,Kawarabayashi2012DirectedCyclesFixedVertexSet,Masarik2019PackingDirectedCircuits}.
The investigation of the Erd\H{o}s-P\'osa property for dicycles with parity constraints has only gained traction very recently with the $\mathsf{W}[1]$-hardness of \textsc{Directed Odd Cycle Transversal} by Lokshtanov, Ramanujan, Saurabh, and Zehavi \cite{Lokshtanov2020ComplexityOddCycle} and the half-integral Erd\H{o}s-P\'osa property for odd dicycles by Kawarabayashi, Kreutzer, Kwon, and Qiqin \cite{Kawarabyashi2023HalfInegralOddDirectedCycles}.

\vspace{-11pt}
\paragraph{Even dicycles}
In contrast to their undirected relatives, the subject of \textit{even dicycles} is surprisingly deep.
For a long time even the computational complexity of deciding whether a digraph contains an even dicycle was unknown \cite{BermondT1981Cycles}.
Vazirani and Yannakakis \cite{VaziraniYannakakis1989PFaffian} showed that the decision problem \textsc{Even Dicycle} is polynomial time equivalent to the recognition of \textit{bipartite Pfaffian graphs}, a problem that can the traced back to an innocent question by P\'olya \cite{polya1913aufgabe} from 1913 which has become known as \textit{P\'olya's Permanent Problem} and is deeply tied to the computational complexity of the permanent \cite{valiant1979complexity,Thomassen1992EvenDicycles}.
A partial solution in terms of structural matching theory was given by Little who described the unavoidable \textit{matching minors} of non-Pfaffian bipartite graphs \cite{Little1975Convertible}.
Indeed, a digraph $D$ has no even directed cycle if and only if there exists an orientation $\vec{B}$ of a bipartite graph $B$ with bipartition $A$ and $B$ a perfect matching $M$ such that (i) exactly the edges in $E(B)\setminus M$ are oriented from $A$ to $B$, (ii) $\vec{B}$ is a \textit{Pfaffian orientation} of $B$, and (iii) $D$ is exactly the digraph obtained from $\vec{B}$ by contracting all edges in $M$.
See the survey in \cite{McCuaig2004Polya} for an explanation of the above observation.
This illustrates that digraphs without even dicycles are exactly those that encode the simplest possible Pfaffian orientations of bipartite graphs.
In 1987, Seymour and Thomassen \cite{SeymourThomassen1987EvenDirectedGraphs} found an analogue of Little's theorem for digraphs and gave a precise description of those \textit{butterfly minors} in a digraph which force an even dicycle to appear.
Coincidentally, this is the first paper using, but not yet naming, butterfly minors as a concept.
A solution for \textsc{Even Dicycle} was independently found by McCuaig \cite{McCuaig2004Polya} and Robertson, Seymour, and Thomas \cite{Robertson1999PermanentsPfaffianOrientations} and jointly announced in 1997 \cite{mccuaig1997permanents}.
Most recently, a first polynomial-time algorithm for \textsc{Shortest Even Dicycle} was found by Bj\"{o}rklund, Husfeldt, and Kaski \cite{Bjorklund2022ShortestEvenDicycle} using algebraic methods and randomisation.
Notably, in stark contrast to the undirected case, large directed treewidth does not force a digraph to contain an even dicycle.

\vspace{-11pt}
\paragraph{Our Result}
In this paper we prove the \textit{first} result on the (fractional) Erd\H{o}s-P\'osa property of even dicycles.
Our proofs are constructive and yield an algorithm with running time $f(k)n^{\mathcal{O}(1)}$ for some computable function $f$ that, given a digraph $D$, either finds a set $S\subseteq V(D)$ of size at most $f(k)$ such that $D-S$ has no even dicycle, or concludes that no such set of size at most $k$ can exist.

Let $D$ be a digraph and $\mathcal{C}=\{ C_1,\dots,C_k\}$ be a family of even dicycles in $D$.
We say that $\mathcal{C}$ is a \emph{quarter-integral packing} of even dicycles if every vertex of $D$ belongs to at most four members of $\mathcal{C}$.
With this, our main theorem reads as follows.

\begin{theorem}\label{thm:mainthm1}
There exists a computable function $f\colon\mathbb{N}\to\mathbb{N}$ such that for every integer $k$ and every digraph $D$ either
\begin{enumerate}
    \item $D$ contains a quarter-integral packing of $k$ even dicycles, or
    \item there exists a set $S\subseteq V(D)$ with $|S|\leq f(k)$ such that $D-S$ does not contain an even dicycle.
\end{enumerate}
Moreover, there exists a computable function $g\colon\mathbb{N}\to\mathbb{N}$ and an algorithm that finds one of the outcomes above in time $g(k)|V(D)|^{\mathcal{O}(1)}$.
\end{theorem}

\begin{figure}[ht]
    \centering
    \scalebox{0.7}{
    \begin{tikzpicture}[scale=1]

        \pgfdeclarelayer{background}
		\pgfdeclarelayer{foreground}
			
		\pgfsetlayers{background,main,foreground}
			
        \begin{pgfonlayer}{main}
        \node (C) [v:ghost] {};

        \foreach \x in {0,...,5}
        {
            \foreach \y in {1,...,12}
            {
                \node (v_\x_\y) [v:mainempty,position=27.69*\y+90:32+5*\x mm from C] {};
            }
        }

        \foreach \x in {0,...,5}
        {
            \node (v_\x_0) [v:main,minimum width=6pt,color=CornflowerBlue,position=27.69*0+90:32+5*\x mm from C] {};
        }

        \foreach \x in {0,...,5}
        {
            \node (w_\x) [v:main,minimum width=6pt,color=CarrotOrange,position=14.19+117.69+55.38*\x:42mm from C] {};
        }

        \end{pgfonlayer}{main}
        
        \begin{pgfonlayer}{foreground}
        \end{pgfonlayer}{foreground}

        \begin{pgfonlayer}{background}

        \draw[e:main,->,bend right=10] (v_0_0.center) to (v_0_1.center);
        \draw[e:main,->,bend right=10] (v_0_1.center) to (v_0_2.center);
        \draw[e:main,->,bend right=10] (v_0_2.center) to (v_0_3.center);
        \draw[e:main,->,bend right=10] (v_0_3.center) to (v_0_4.center);
        \draw[e:main,->,bend right=10] (v_0_4.center) to (v_0_5.center);
        \draw[e:main,->,bend right=10] (v_0_5.center) to (v_0_6.center);
        \draw[e:main,->,bend right=10] (v_0_6.center) to (v_0_7.center);
        \draw[e:main,->,bend right=10] (v_0_7.center) to (v_0_8.center);
        \draw[e:main,->,bend right=10] (v_0_8.center) to (v_0_9.center);
        \draw[e:main,->,bend right=10] (v_0_9.center) to (v_0_10.center);
        \draw[e:main,->,bend right=10] (v_0_10.center) to (v_0_11.center);
        \draw[e:main,->,bend right=10] (v_0_11.center) to (v_0_12.center);
        \draw[e:main,->,bend right=10] (v_0_12.center) to (v_0_0.center);

        \draw[e:main,->,bend right=10] (v_1_0.center) to (v_1_1.center);
        \draw[e:main,->,bend right=10] (v_1_1.center) to (v_1_2.center);
        \draw[e:main,->,bend right=10] (v_1_2.center) to (v_1_3.center);
        \draw[e:main,->,bend right=10] (v_1_3.center) to (v_1_4.center);
        \draw[e:main,->,bend right=10] (v_1_4.center) to (v_1_5.center);
        \draw[e:main,->,bend right=10] (v_1_5.center) to (v_1_6.center);
        \draw[e:main,->,bend right=10] (v_1_6.center) to (v_1_7.center);
        \draw[e:main,->,bend right=10] (v_1_7.center) to (v_1_8.center);
        \draw[e:main,->,bend right=10] (v_1_8.center) to (v_1_9.center);
        \draw[e:main,->,bend right=10] (v_1_9.center) to (v_1_10.center);
        \draw[e:main,->,bend right=10] (v_1_10.center) to (v_1_11.center);
        \draw[e:main,->,bend right=10] (v_1_11.center) to (v_1_12.center);
        \draw[e:main,->,bend right=10] (v_1_12.center) to (v_1_0.center);

        \draw[e:main,->,bend right=10] (v_2_0.center) to (v_2_1.center);
        \draw[e:main,->,bend right=10] (v_2_1.center) to (v_2_2.center);
        \draw[e:main,->,bend right=10] (v_2_2.center) to (v_2_3.center);
        \draw[e:main,->,bend right=10] (v_2_3.center) to (v_2_4.center);
        \draw[e:main,->,bend right=10] (v_2_4.center) to (v_2_5.center);
        \draw[e:main,->,bend right=10] (v_2_5.center) to (v_2_6.center);
        \draw[e:main,->,bend right=10] (v_2_6.center) to (v_2_7.center);
        \draw[e:main,->,bend right=10] (v_2_7.center) to (v_2_8.center);
        \draw[e:main,->,bend right=10] (v_2_8.center) to (v_2_9.center);
        \draw[e:main,->,bend right=10] (v_2_9.center) to (v_2_10.center);
        \draw[e:main,->,bend right=10] (v_2_10.center) to (v_2_11.center);
        \draw[e:main,->,bend right=10] (v_2_11.center) to (v_2_12.center);
        \draw[e:main,->,bend right=10] (v_2_12.center) to (v_2_0.center);

        \draw[e:main,->,bend right=10] (v_3_0.center) to (v_3_1.center);
        \draw[e:main,->,bend right=10] (v_3_1.center) to (v_3_2.center);
        \draw[e:main,->,bend right=10] (v_3_2.center) to (v_3_3.center);
        \draw[e:main,->,bend right=10] (v_3_3.center) to (v_3_4.center);
        \draw[e:main,->,bend right=10] (v_3_4.center) to (v_3_5.center);
        \draw[e:main,->,bend right=10] (v_3_5.center) to (v_3_6.center);
        \draw[e:main,->,bend right=10] (v_3_6.center) to (v_3_7.center);
        \draw[e:main,->,bend right=10] (v_3_7.center) to (v_3_8.center);
        \draw[e:main,->,bend right=10] (v_3_8.center) to (v_3_9.center);
        \draw[e:main,->,bend right=10] (v_3_9.center) to (v_3_10.center);
        \draw[e:main,->,bend right=10] (v_3_10.center) to (v_3_11.center);
        \draw[e:main,->,bend right=10] (v_3_11.center) to (v_3_12.center);
        \draw[e:main,->,bend right=10] (v_3_12.center) to (v_3_0.center);

        \draw[e:main,->,bend right=10] (v_4_0.center) to (v_4_1.center);
        \draw[e:main,->,bend right=10] (v_4_1.center) to (v_4_2.center);
        \draw[e:main,->,bend right=10] (v_4_2.center) to (v_4_3.center);
        \draw[e:main,->,bend right=10] (v_4_3.center) to (v_4_4.center);
        \draw[e:main,->,bend right=10] (v_4_4.center) to (v_4_5.center);
        \draw[e:main,->,bend right=10] (v_4_5.center) to (v_4_6.center);
        \draw[e:main,->,bend right=10] (v_4_6.center) to (v_4_7.center);
        \draw[e:main,->,bend right=10] (v_4_7.center) to (v_4_8.center);
        \draw[e:main,->,bend right=10] (v_4_8.center) to (v_4_9.center);
        \draw[e:main,->,bend right=10] (v_4_9.center) to (v_4_10.center);
        \draw[e:main,->,bend right=10] (v_4_10.center) to (v_4_11.center);
        \draw[e:main,->,bend right=10] (v_4_11.center) to (v_4_12.center);
        \draw[e:main,->,bend right=10] (v_4_12.center) to (v_4_0.center);

        \draw[e:main,->,bend right=10] (v_5_0.center) to (v_5_1.center);
        \draw[e:main,->,bend right=10] (v_5_1.center) to (v_5_2.center);
        \draw[e:main,->,bend right=10] (v_5_2.center) to (v_5_3.center);
        \draw[e:main,->,bend right=10] (v_5_3.center) to (v_5_4.center);
        \draw[e:main,->,bend right=10] (v_5_4.center) to (v_5_5.center);
        \draw[e:main,->,bend right=10] (v_5_5.center) to (v_5_6.center);
        \draw[e:main,->,bend right=10] (v_5_6.center) to (v_5_7.center);
        \draw[e:main,->,bend right=10] (v_5_7.center) to (v_5_8.center);
        \draw[e:main,->,bend right=10] (v_5_8.center) to (v_5_9.center);
        \draw[e:main,->,bend right=10] (v_5_9.center) to (v_5_10.center);
        \draw[e:main,->,bend right=10] (v_5_10.center) to (v_5_11.center);
        \draw[e:main,->,bend right=10] (v_5_11.center) to (v_5_12.center);
        \draw[e:main,->,bend right=10] (v_5_12.center) to (v_5_0.center);


        \draw[e:main,->] (v_0_1.center) to (v_1_1.center);
        \draw[e:main,->] (v_1_1.center) to (v_2_1.center);
        \draw[e:main,->] (v_2_1.center) to (v_3_1.center);
        \draw[e:main,->] (v_3_1.center) to (v_4_1.center);
        \draw[e:main,->] (v_4_1.center) to (v_5_1.center);

        \draw[e:main,->] (v_0_3.center) to (v_1_3.center);
        \draw[e:main,->] (v_1_3.center) to (v_2_3.center);
        \draw[e:main,->] (v_2_3.center) to (v_3_3.center);
        \draw[e:main,->] (v_3_3.center) to (v_4_3.center);
        \draw[e:main,->] (v_4_3.center) to (v_5_3.center);

        \draw[e:main,->] (v_0_5.center) to (v_1_5.center);
        \draw[e:main,->] (v_1_5.center) to (v_2_5.center);
        \draw[e:main,->] (v_2_5.center) to (v_3_5.center);
        \draw[e:main,->] (v_3_5.center) to (v_4_5.center);
        \draw[e:main,->] (v_4_5.center) to (v_5_5.center);

        \draw[e:main,->] (v_0_7.center) to (v_1_7.center);
        \draw[e:main,->] (v_1_7.center) to (v_2_7.center);
        \draw[e:main,->] (v_2_7.center) to (v_3_7.center);
        \draw[e:main,->] (v_3_7.center) to (v_4_7.center);
        \draw[e:main,->] (v_4_7.center) to (v_5_7.center);

        \draw[e:main,->] (v_0_9.center) to (v_1_9.center);
        \draw[e:main,->] (v_1_9.center) to (v_2_9.center);
        \draw[e:main,->] (v_2_9.center) to (v_3_9.center);
        \draw[e:main,->] (v_3_9.center) to (v_4_9.center);
        \draw[e:main,->] (v_4_9.center) to (v_5_9.center);

        \draw[e:main,->] (v_0_11.center) to (v_1_11.center);
        \draw[e:main,->] (v_1_11.center) to (v_2_11.center);
        \draw[e:main,->] (v_2_11.center) to (v_3_11.center);
        \draw[e:main,->] (v_3_11.center) to (v_4_11.center);
        \draw[e:main,->] (v_4_11.center) to (v_5_11.center);


        \draw[e:main,->] (v_5_2.center) to (v_4_2.center);
        \draw[e:main,->] (v_4_2.center) to (v_3_2.center);
        \draw[e:main,->] (v_3_2.center) to (v_2_2.center);
        \draw[e:main,->] (v_2_2.center) to (v_1_2.center);
        \draw[e:main,->] (v_1_2.center) to (v_0_2.center);

        \draw[e:main,->] (v_5_4.center) to (v_4_4.center);
        \draw[e:main,->] (v_4_4.center) to (v_3_4.center);
        \draw[e:main,->] (v_3_4.center) to (v_2_4.center);
        \draw[e:main,->] (v_2_4.center) to (v_1_4.center);
        \draw[e:main,->] (v_1_4.center) to (v_0_4.center);

        \draw[e:main,->] (v_5_6.center) to (v_4_6.center);
        \draw[e:main,->] (v_4_6.center) to (v_3_6.center);
        \draw[e:main,->] (v_3_6.center) to (v_2_6.center);
        \draw[e:main,->] (v_2_6.center) to (v_1_6.center);
        \draw[e:main,->] (v_1_6.center) to (v_0_6.center);

        \draw[e:main,->] (v_5_8.center) to (v_4_8.center);
        \draw[e:main,->] (v_4_8.center) to (v_3_8.center);
        \draw[e:main,->] (v_3_8.center) to (v_2_8.center);
        \draw[e:main,->] (v_2_8.center) to (v_1_8.center);
        \draw[e:main,->] (v_1_8.center) to (v_0_8.center);

        \draw[e:main,->] (v_5_10.center) to (v_4_10.center);
        \draw[e:main,->] (v_4_10.center) to (v_3_10.center);
        \draw[e:main,->] (v_3_10.center) to (v_2_10.center);
        \draw[e:main,->] (v_2_10.center) to (v_1_10.center);
        \draw[e:main,->] (v_1_10.center) to (v_0_10.center);

        \draw[e:main,->] (v_5_12.center) to (v_4_12.center);
        \draw[e:main,->] (v_4_12.center) to (v_3_12.center);
        \draw[e:main,->] (v_3_12.center) to (v_2_12.center);
        \draw[e:main,->] (v_2_12.center) to (v_1_12.center);
        \draw[e:main,->] (v_1_12.center) to (v_0_12.center);


        \draw[e:main,->,bend right=25] (v_5_1.center) to (w_0.center);
        \draw[e:main,->,bend right=25] (w_0.center) to (v_5_2.center);

        \draw[e:main,->,bend right=25] (v_5_3.center) to (w_1.center);
        \draw[e:main,->,bend right=25] (w_1.center) to (v_5_4.center);

        \draw[e:main,->,bend right=25] (v_5_5.center) to (w_2.center);
        \draw[e:main,->,bend right=25] (w_2.center) to (v_5_6.center);

        \draw[e:main,->,bend right=25] (v_5_7.center) to (w_3.center);
        \draw[e:main,->,bend right=25] (w_3.center) to (v_5_8.center);

        \draw[e:main,->,bend right=25] (v_5_9.center) to (w_4.center);
        \draw[e:main,->,bend right=25] (w_4.center) to (v_5_10.center);

        \draw[e:main,->,bend right=25] (v_5_11.center) to (w_5.center);
        \draw[e:main,->,bend right=25] (w_5.center) to (v_5_12.center);
        
        \end{pgfonlayer}{background}
        
    \end{tikzpicture}
    }
    \caption{A strongly planar digraph without two vertex-disjoint even dicycles but a large hitting set for all even dicycles.
    Each dicycle contains exactly one of the \textcolor{CornflowerBlue}{blue} vertices and a dicycle is even if and only if it contains an odd number of the \textcolor{CarrotOrange}{orange} vertices.
    It is not difficult to see that this digraph does not contain two disjoint even dicycles.
    At the same time, one needs to delete at least $k=6$ vertices to remove all even dicycles in the graph.
    Moreover, it is straight forward to scale this example for any integer $k$.}
    \label{fig:planarcounterexample}
\end{figure}
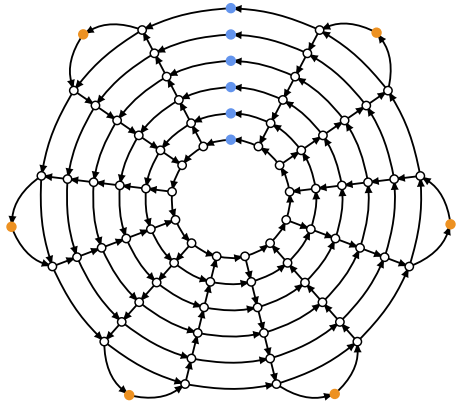

\Cref{thm:mainthm1} says that even dicycles have the quarter-integral Erd\H{o}s-P\'osa property.
Notice this is almost best possible since even dicycles do not have the \textit{integral} Erd\H{o}s-P\'osa property as illustrated by the example in \cref{fig:planarcounterexample}.

Our results imply an approximation algorithm for a generalised notion of directed treewidth in the following sense.
An \emph{odd directed tree-decomposition} for a digraph $D$ is a four tuple $\mathcal{T}=(T,\alpha,\beta,\gamma)$ where $\alpha,\beta\colon V(T)\to 2^{V(G)}$, we call $\beta$ the \emph{bags}, and $\gamma\colon E(T)\to 2^{V(D)}$ where $\{ \beta(t) \mid t\in V(T)\}\setminus\{ \emptyset\}$ is a partition of $V(D)$, for every $t\in V(T)$, $\alpha(t)\subseteq \Gamma(t)\coloneqq \beta(t)\cup\bigcup_{e\sim t}\gamma(t)$, for every $t\in V(T)$, no even dicycle of $D-\alpha(t)$ contains a vertex of $\Gamma(t)\setminus\alpha(t)$, and for every $e\in V(T)$, $\gamma(t)$ contains a vertex of every dicycle (even or not) which contains vertices from bags of both components of $T-e$.
We call $\mathcal{T}$ \emph{strong} if for any $t\in V(T)$ no strong component of $D-\alpha(t)$ contains a vertex of $\Gamma(t)$ and $D-\Gamma(t)$ at the same time.
The \emph{width} of $\mathcal{T}$ is the maximum size of $\alpha(t)$ over all $t\in V(T)$ and $\gamma(e)$ over all $e\in E(T)$, and the \emph{odd directed treewidth}, denoted by $\mathsf{odtw}(D)$, of $D$ is the minimum width over all odd directed tree-decompositions for $D$.

By only considering odd directed tree-decomposition where $\alpha(t)=\Gamma(t)$ for all $t\in V(T)$, we recover the notion of \emph{directed treewidth}\footnote{Up to a $-1$ in the original definition in \cite{Johnson2001DirectedTreewidth}.}, denoted by $\mathsf{dtw}(D)$.
Hence, $\mathsf{odtw}(D)\leq\mathsf{dtw}(D)$ for all digraphs $D$.
As part of our proof for \cref{thm:mainthm1} we obtain the following.

\begin{theorem}\label{thm:mainthm2}
There exist functions $f,g\colon\mathbb{N}\to\mathbb{N}$ and an algorithm that provides, for every digraph $D$ and every integer $k$, one of the following
\begin{enumerate}
    \item A strong odd directed tree-decomposition for $D$ of width at most $f(k)$, or
    \item the correct conclusion that $\mathsf{odtw}(D)>k$.
\end{enumerate}
This algorithm runs in time $g(k)|V(D)|^{\mathcal{O}(1)}$.
\end{theorem}

Moreover, as a consequence of \cref{cor:ddpp} below and \cref{thm:mainthm2}, we can extend the algorithm for the $t$-\textsc{Directed Disjoint Paths Problem} ($t$-DDPP) of Johnson et al.\@ \cite{Johnson2001DirectedTreewidth} for digraphs of bounded $\mathsf{dtw}$ to digraphs of bounded $\mathsf{odtw}$.
The $t$-DDPP asks, given a set of $t$ vertex pairs $(s_i,t_i)$ in a digraph $D$, for $t$ pairwise internally disjoint directed paths $P_1,\dots,P_t$ such that $P_i$ starts in $s_i$ and ends in $t_i$ for all $i\in\{1,\dots,t \}$.
This problem is known to be $\mathsf{NP}$-complete even for $t=2$ \cite{Fortune1980DirectedSubgraphHomeomorphism} and has been the focus of a long line of research by various authors \cite{Even1975Timetable,Fortune1980DirectedSubgraphHomeomorphism,schrijver1994finding,Johnson2001DirectedTreewidth,slivkins2010parameterized,cygan2013planar,steiner2020parameterized,giannopoulou2020canonical}.

\begin{theorem}\label{thm:mainthm3}
For any integers $t,k$ and any digraph $D$ with $\mathsf{odtw}(D)\leq k$, the $t$-DDPP can be solved in time $n^{\mathcal{O}(t+f(k))}$, where $f$ is the function from \Cref{thm:mainthm2}.
\end{theorem}

\Cref{thm:mainthm3} is a \textit{strict} generalisation of the algorithm of Johnson et al.\@ for the $t$-DDPP to classes more general than those of bounded directed treewidth.
This is the first direct parameterisation strictly generalising this result.
Since $2$-DDPP is $\mathsf{NP}$-complete on general digraphs \cite{Fortune1980DirectedSubgraphHomeomorphism}, a full generalisation of the Graph Minors Algorithm of Robertson and Seymour \cite{GraphMinorsXIII} is unlikely to exist.

\subsection{Digraphs without even dicycles and algorithmic implications of our result}

The algorithm of \cref{thm:mainthm1} allows to find in polynomial time a set $S$ of at most $f(k)$ vertices in a digraph $D$ such that $D-S$ does not contain an even dicycle or to conclude that one needs to delete at least $k+1$ vertices to achieve this.

One explanation for the high level of interest in Reed's proof of the half-integral Erd\H{o}s-P\'osa property of odd undirected cycles \cite{Reed1996PackingDirectedCirctuis} and the strategies for finding small hitting sets for odd undirected cycles that followed, is that deleting such a small \textit{odd cycle transversal} leaves a bipartite graph.
Bipartite graphs themselves have strong algorithmic properties, many of which can be lifted to graphs with small odd cycle transversal.
This includes fundamental computational problems such as \textsc{Minimum Vertex Cover}, \textsc{Maximum Independent Set}, and \textsc{Maximum Cut} (see for example \cite{jaffke2023dynamic}).
Somewhat surprisingly, digraphs without even dicycles resemble bipartite undirected graphs in many ways as we will explain in the next paragraph.

A digraph $D$ is said to be \emph{non-even} if it can be turned into a digraph $D'$ that does not contain an even dicycle by subdividing some of its edges.
Non-even digraphs first came up because the decision problem \textsc{Even Dicycle} can be seen to be polynomial-time equivalent to the recognition problem for non-even digraphs \cite{SeymourThomassen1987EvenDirectedGraphs}.
Indeed, as an immediate consequence of the definition one can observe that, if $D-S$ does not have an even dicycle, it surely is non-even.

A characterisation for non-even digraphs had already been sought for a long time, in particular in connection to the \textsc{Even Dicycle} problem, and they played a major role it its resolution by Robertson, Seymour, and Thomas, and McCuaig \cite{mccuaig1997permanents}.
Whilst McCuaig \cite{McCuaig2004Polya} worked in the matching theoretic setting, Robertson, Seymour, and Thomas \cite{Robertson1999PermanentsPfaffianOrientations} directly\footnote{They later opted for a matching theoretic presentation of their results as it simplifies some of the arguments.} showed that, with the exception of one special graph, the strongly 2-connected non-even digraphs, also called \emph{dibraces}, could be constructed from \emph{strongly planar}\footnote{Strongly planar digraphs form an important subclass of non-even and planar digraphs.} dibraces via a simple sum operation, leaving the summands to interact on at most four vertices in a very restricted fashion.
Thus, when narrowing our view to dibraces, non-even digraphs are structurally close to being strongly planar. 

Analysing the particularities of this class, several nice algorithmic properties have been found.
\begin{enumerate}
\item Where for bipartite graphs the vertex cover number and the matching number are equal, for non-even digraphs, with the exception of a unique example, the \textit{feedback vertex number}\footnote{The \emph{feedback vertex number} of a digraphs is the minimum size of a vertex set whose deletion leave an acyclic digraph.} equals the \textit{cycle packing number}\footnote{The \emph{cycle packing number} of a digraph $D$ is the maximum size of a packing of dicycles in $D$.} \cite{GueninThomas2011PackingDirectedCircuits}.
\item The notoriously hard $t$-DDPP \cite{Even1975Timetable,Fortune1980DirectedSubgraphHomeomorphism,slivkins2010parameterized} can be solved in polynomial time (if $t$ is fixed) on non-even digraphs \cite{Giannopoulou2023ExcludingSingleCrossing}.
\item Moreover, a well known construction shows that every non-even digraph corresponds directly to a bipartite graph $B$ together with a perfect matching $M$, where $B$ is \emph{Pfaffian} \cite{McCuaig2004Polya,Robertson1999PermanentsPfaffianOrientations,Giannopoulou2024ExcludingPlanar}.
Bipartite Pfaffian graphs are of significance because they admit so called Pfaffian orientations that allow for the efficient computation of the number of perfect matchings (see the survey by Thomas for more on Pfaffian orientations \cite{Thomas2006PfaffianSurvey}) which naturally corresponds to the fundamental $\#\mathsf{P}$-complete problem of computing the permanent \cite{valiant1979complexity}.
\end{enumerate}

We can define a simple parameter which, if bounded, yields algorithmic results relating to each one of the above points as direct corollaries to \Cref{thm:mainthm1}, or at least heavily suggests an avenue through which these problems should become more tractable.
Given a digraph $D$, let us call a set $S\subseteq V(D)$ an \emph{even dicycle transversal} if $D-S$ does not contain an even dicycle.
The \emph{even dicycle transversal number} of $D$ is the minimum size of an even dicycle transversal in $D$.


\begin{corollary}\label{cor:ddpp}
    For any integers $t,k$ and any digraph $D$ with even dicycle transversal number at most $k$, the $t$-DDPP can be solved in time $n^{\mathcal{O}(t+f(k))}$, where $f$ is the function from \Cref{thm:mainthm1}.
\end{corollary}
\begin{proof}[Proof of \Cref{cor:ddpp} (Rough Sketch)]
    The algorithm for the $t$-DDPP from \cite{Giannopoulou2023ExcludingSingleCrossing} has a running time of the form $n^{\mathcal{O}(t)}$.
    This means that one may simply enumerate all $\mathcal{O}(n^{2k})$ possible involvements of the even dicycle transversal in potential solutions to a given instance of the $t$-DDPP.
    This yields $\mathcal{O}(n^{2k})$ instances of the $(t+2k)$-DDPP on non-even digraphs.
    Hence, whenever the even dicycle transversal number of a digraph is bounded by some number $k$, one can use our algorithm from \cref{thm:mainthm1} to obtain an even dicycle transversal of size at most $f(3k)$ and can therefore solve the $t$-DDPP in time $n^{\mathcal{O}(t+f(k))}$.
    Notice that a better running-time is unlikely to be achievable since $t$-DDPP is $\mathsf{W}[1]$-hard already on acyclic digraphs \cite{Even1975Timetable,Fortune1980DirectedSubgraphHomeomorphism}.
\end{proof}

\begin{corollary}\label{cor:countingperfectmatchings}
    For any integer $k$ and any bipartite graph $B$ such that there exists a set $F$ of at most $k$ edges in $B$ for which $B - V(F)$ is Pfaffian, one can determine the number of perfect matchings of $B$ in time $n^{\mathcal{O}(|F|)}$.
\end{corollary}
\begin{proof}[Proof of \Cref{cor:countingperfectmatchings} (Rough Sketch)]
    Given some bipartite graph $B$ with a perfect matching $M$ and the corresponding digraph $D$, one can observe that any even dicycle transversal $S$ of $D$ corresponds to a set $F$ of $|S|$ edges of $M$ such that deleting the endpoints of the edges in $F$ from $B$ yields a Pfaffian graph.
    By enumerating all $|V(B)|^{2|S|}$ possible intersections of perfect matchings of $B$ with the endpoints of the edges in $F$, this allows to count the perfect matchings in $B$ in time $|V(B)|^{\mathcal{O}(|S|)}$.
    Again the key is the algorithm from \cref{thm:mainthm1} which makes it possible to find such a set $S$ or determine that no small even dicycle transversal can exist in $D$.
    Since counting perfect matchings is $\mathsf{W}[1]$-hard on graphs which become planar after deleting at most $k$ vertices \cite{CurticapeanXia2015ParameterizingThePermanent}, it is unlikely that the dependency on the size of $S$ in the degree of the polynomial of the running time can be removed.
\end{proof}

The case of \textsc{Minimum Feedback Vertex Set} and \textsc{Maximum Cycle Packing} in non-even digraphs is less immediate.
The proofs of Guenin and Thomas \cite{GueninThomas2011PackingDirectedCircuits} imply polynomial time algorithms for both problems on strongly planar digraphs via the Lucchesi-Younger Theorem \cite{LucchesiYounger1978Minmax}.
However, it is not clear if this can be extended to solve any of the two problems on the class of non-even digraphs in polynomial time.
If this is the case, it seems conceivable that such an algorithm could be extended to graphs of bounded even dicycle transversal number.
We leave this as an open problem.

\begin{question}\label{que:questionregardinghittingandpackingdicycles}
What is the computational complexity of the problems \textsc{Maximum Cycle Packing} and \textsc{Minimum Feedback Vertex Set} on non-even digraphs?
\end{question}

The original paper on directed treewidth \cite{Johnson2001DirectedTreewidth} and a recent result of Giannopoulou, Kreutzer, and Wiederrecht \cite{Giannopoulou2024ExcludingPlanar} show that the strategies above for the $t$-DDPP and counting perfect matchings in bipartite graphs also apply for digraphs of bounded directed treewidth.
However, digraphs without even dicycles may have arbitrarily large directed treewidth (consider for an example the digraphs from \cref{fig:planarcounterexample} when omitting the \textcolor{CarrotOrange}{orange} vertices).
It can be observed that the underlying undirected graphs of digraphs without even dicycles may have arbitrarily large clique minors \cite{Giannopoulou2023ExcludingSingleCrossing} and can therefore also be far away from being planar.
We believe that this class of digraphs has considerable untapped algorithmic potential and hope that this paper will help with its realisation.

\subsection{The established strategy: A case study of odd cycles}\label{subsec:undirectedstrategy}

To get a grasp on the type of difficulties one encounters specifically when dealing with even dicycles we use this subsection to sketch the approach Reed took to prove the half-integral Erd\H{o}s-P\'osa property of odd cycles \cite{Reed1999Fruitsalad}.
On an abstract level, this approach is somewhat universal and was adapted for proving extensions of this result to the group labelled setting \cite{Huynh2019UnifiedErdosPosaConstrained,ThomasYoo2023PackingCyclesUndirected}.
Another line of arguments replaces some of the deeper structural tools of Reed's approach with the more abstract notion of \textit{highly linked sets} and \textit{tangles}.
This type of argument was used for more recent generalisations of Reed's original result to group labelled graphs \cite{gollin2021unified,gollin2022unified} as well as to prove the half-integral Erd\H{o}s-P\'osa property for odd dicycles \cite{Kawarabyashi2023HalfInegralOddDirectedCycles}.
Both approaches share a common idea at their core and use this to enhance suitable parts of the graph minors theory of Robertson and Seymour.
We focus here on the more involved structural approach of Reed\footnote{Reed's original paper did not use the \textit{odd $A$-path theorem} directly, but his proof already bears a remarkable resemblance to the modern strategies.} as this is the one closer in spirit to our own strategy.

\vspace{-11pt}
\paragraph{The $A$-path theorem(s)}
The main tool which is used to augment the graph minors theory of Robertson and Seymour is an appropriate variation of Mader's theorem, sometimes referred to as the $A$-path theorem \cite{Mader1978Kreuzungsfrei}.
Roughly speaking, given some graph $G$ and a set $A\subseteq V(G)$, Mader's theorem asserts that there are either many pairwise disjoint paths, each with exactly its endpoints in $A$, so called \emph{$A$-paths}, or there is a small set of vertices whose deletion removes all such paths from $G$.
A variant of Mader's theorem for $A$-paths of odd length was proven in the context of the odd-minor variant of Hadwiger's conjecture \cite{Geelen2009OddHadwiger} and has since been generalised in various ways \cite{Chudnovsky2006PackingAPaths,Chudnovsky2008APathAlgorithm,Wollan2010PackingAPaths,Yamaguchi2019PackingAPaths}.

\vspace{-11pt}
\paragraph{Bipartite tree-decompositions}
Once an $A$-path theorem for paths of odd length has been established, one may use this to augment the structure theory of Robertson and Seymour as follows.
Let us fix some integer $k$ and suppose we want to either find a half-integral packing of $k$ odd cycles in a graph $G$, or delete a small set of vertices $S$ such that $G-S$ is bipartite.
\begin{enumerate}
    \item In the case where $G$ has small treewidth a standard argument may be applied to either find an integral packing of odd cycles or a small odd cycle transversal.
    \item If $G$ does not have small treewidth it contains some witness\footnote{Typically such a witness would be a \textit{tangle} or a \textit{highly linked set}.} of large treewidth.
    This witness can take one of two forms.
    \begin{enumerate}
        \item It might be that there exists a large, compared to $k$, clique minor.
        If this clique minor does not already contain many pairwise disjoint odd cycles, it contains a bipartite subgraph which still has a large clique minor.
        Using the odd $A$-path theorem either yields a large packing of odd cycles, or a small set $S\subseteq V(G)$ such that the \textit{block}\footnote{A \emph{block} in a graph is a maximal subgraph without a cut vertex.} of $G-S$ which contains the remaining clique minor is bipartite.
        See \cite{Geelen2009OddHadwiger} for more details.
        \item If there is no clique minor, then one may use the \textit{Flat Wall Theorem} of Robertson and Seymour \cite{GraphMinorsXIII}.
        This theorem, in the absence of a clique minor, yields a small set $A\subseteq V(G)$ and a large ``piece'' of $G-A$ which behaves like a planar graph and is supported by a large grid-like infrastructure.
        This \textit{flat wall} now either contains a large bipartite part, still with a large wall, or we have found our packing of odd cycles.
        Finally, applying the odd $A$-path theorem to the boundary of the wall yields either a half-integral packing, or a small set $S$ such that the block of $G-(A \cup S)$ containing the wall is bipartite.
    \end{enumerate}
    \item Finally, these observations may be used to find a tree-decomposition of $G$ into pieces that are either small or become bipartite after deleting a small set of vertices.
    See \cite{KawarabayashiReed2010OddCyclePacking,KawarabayashiReed2010AlmostLinearTimeOddCycleTransversal,jaffke2023dynamic} for more in-depth descriptions of this structure.
    Such a tree-decomposition yields the half-integral Erd\H{o}s-P\'osa property for odd cycles via a simple dynamic programming argument.
\end{enumerate}
The reason why this only provides a half-integral version of the Erd\H{o}s-P\'osa property can be found in one of the possible outcomes of applying the odd $A$-path theorem to the flat wall.
Namely the \emph{Escher-wall}, which is a bipartite wall $W$ together with a set of parity-breaking $A$-paths where $A$ is the vertex set of the boundary of the wall such that that the entire graph can be embedded in the projective plane in a way such that every face is bounded by an even cycle.

\subsection{Challenges in the directed setting}
While some of the important building blocks for the strategy presented in \cref{subsec:undirectedstrategy} have been extended to digraphs, one key aspect seems almost insurmountable.
That is, while \textsc{Even Dicycle} belongs to $\mathsf{P}$, the related problem \textsc{Even Directed $s$-$t$-Path} of deciding if there exists a directed $s$-$t$-path of even length in a digraph is known to be $\mathsf{NP}$-complete \cite{LaPaughPapadimitriou1984EvenPathProblem}.
This means that a directed analogue of the $A$-path theorem for even directed paths is unlikely to exist, or at least be tractible.
A similar problem occurs for odd directed $A$-paths.
To prove their half-integral Erd\H{o}s-P\'osa theorem for odd dicycles, Kawarabayashi, Kreutzer, Kwon, and Qiqin \cite{Kawarabyashi2023HalfInegralOddDirectedCycles} proved a half-integral variant for odd directed $A$-\textit{walks}.
For odd dicycles this suffices because any odd directed $s$-$t$-walk contains either an odd dicycle, or an odd directed $s$-$t$-path.
So in the cases where the odd directed $A$-walk theorem does not generate many odd directed $A$-walks, it necessarily produces a large half-integral packing of odd dicycles.
Sadly, this strategy is not feasible for even dicycles as an even directed $s$-$t$-walk can possibly be broken into an odd directed $s$-$t$-path and an odd dicycle.

The lack of an even directed $A$-path theorem creates a ripple effect that forces us to rethink the entire approach.
For instance: A priori we do not have a way to handle the situation in which we find a large complete digraph as a butterfly minor.
This is, because in the established strategy one would either find many cycles of the desired parity or a large subclique where all cycles are of the other parity.
We would then use the $A$-path theorem for paths of the correct parity to either remove all ``bad'' cycles which attach to the clique, or find the desired packing.
Without an even directed $A$-path theorem, cliques need to be dealt with in an entirely different way.
Luckily for us, non-even digraphs are indeed butterfly minor closed which means that any large enough clique minor contains many copies of obstructions for non-even digraphs \cite{SeymourThomassen1987EvenDirectedGraphs}.
By definition, each of these obstructions needs to contain an even dicycle.
We elaborate on this in \cref{subsec:approach}.

To make things worse, the \textit{Directed Flat Wall Theorem} \cite{Giannopoulou2020DirectedFlatWall} does not fully provide the nice topological properties which are guaranteed by its undirected counterpart.
This means that the ``area'' defined by the directed flat wall cannot be treated as a planar graph in the way it is being done in the undirected setting.
We will discuss how to overcome this problem in \cref{subsec:approach}.
For now let us assume we can obtain a stronger version of the Directed Flat Wall Theorem that lets us deal with this situation.
Due to the lack of the even directed $A$-path theorem, we cannot easily remove all even dicycles which might interact with our wall.
Instead, a promising strategy could be to try to \textit{expand} the domain which is governed by our flat structure.

An established and elegant way of increasing the territory where we can fully capture the structure of our graph through ``flatness'' is described in the new proof for the Graph Minor Structure Theorem, originally conceived by Robertson and Seymour \cite{robertson2003graph}, by Kawarabayashi, Thomas, and Wollan \cite{kawarabayashi2020quickly}.
This strategy starts from a large flat wall by iteratively attaching huge linkages to the wall, thereby creating new walls which can be treated with the Flat Wall Theorem.
However, to maintain control over this process, Kawarabayashi et al.\ make use of the ability of undirected path families crossing through a large family of cycles to be \textit{normalised}.
This normalisation is a rerouting process which becomes impossible in certain situations in directed graphs as the directions of paths and cycles need to be respected.
Another particularly nice detail of the proof of Kawarabayashi et al.\ is that they can use large linkages that run through a wall to split their wall into two disjoint walls, each describing a new area which then needs to be further flattened over the course of their procedure.
This however does not work in the directed setting, which can easily be seen by considering a directed cycle with a chord.
In this configuration there exists a unique dicycle using the chord, a property that is not true in the undirected setting. 
This forces us to diversify our structure away from simply considering walls.

\subsection{Our approach}\label{subsec:approach}
Recall the discussion on Reed's strategy for the half-integral Erd\H{o}s-P\'osa property for odd cycles.
The crucial step of this proof was to either fully describe how many odd cycles attach to a given wall, or to isolate the wall, and therefore the associated area of ``high connectivity'' from all odd cycles that occur in the graph by deleting a small set of vertices.
While, as we have laid out in detail above, we cannot expect to be successful as easily as in the proofs for odd (di)cycles, we can still follow this very abstract guideline.
Thus, our strategy will be to decompose any given digraph into pieces which are, after the deletion of a small vertex set, completely disjoint from any even dicycle in the digraph.
As before, the crucial part of this is to describe how such an outcome can be achieved in areas of the digraph which witness large directed treewidth: a \textit{local structure theorem}.
The rest of this introduction is dedicated to providing insight into how we resolve this problem.

The first part of this subsection is dedicated to providing a light introduction to the several steps and tools necessary to establish the structural result at the heart of our strategy.
While the second part strives to provide more detailed insights.
The one advantage we have, despite the many challenges imposed by the directed setting, is the fact that the exclusion of \textit{odd bicycles}\footnote{A \emph{bicycle} is a digraph obtained from an undirected cycle by replacing every undirected edge with two anti-parallel directed edges. It is odd if it was obtained from an odd cycle.} provides a structure that is topological in nature, as we explained when addressing the results of \cite{McCuaig2004Polya} and \cite{Robertson1999PermanentsPfaffianOrientations}.
Moreover, and somewhat ironically, odd bicycles are exactly the obstruction to being non-even and as such, they must contain even dicycles.
Hence, in many ways, our proof is a structural study of odd bicycles in disguise.
So the main theme of the entire strategy is to ``corner'' the even dicycles in our graph inside disks where the graph behaves in a non-planar way, while everywhere else it is embedded in the plane in a quite structured manner.

\vspace{-11pt}
\paragraph{Trapping even dicycles: Our main tools}
To perform the necessary surgery on the pieces of our digraph supported by a large grid-like infrastructure, we import one major tool from structural matching theory which was introduced by Giannopoulou and Wiederrecht \cite{giannopoulou2021two} and then further refined by them in joint work with Thilikos \cite{Giannopoulou2023ExcludingSingleCrossing}.
This tool is a matching theoretic variant of the \textbf{Two Paths Theorem} (see for example \cite{kawarabayashi2020quickly} for a discussion on the role the original Two Paths Theorems plays in the Graph Minors Series of Robertson and Seymour) and provides a link between topology, i.\@e.\@ allows is to find our desired (almost) embedding in the first place, and the existence of an odd bicycle, i.\@e.\@ the presence of an even dicycle with additional topological information.

The second major tool we use to finally confine all even dicycles in the graph into controlled areas of our (almost) embedding is a \textbf{shifting argument} for even dicycles in (almost) planar graphs.
Informally speaking, our argument says that, given a digraph $D$ which is the union of an odd and an even dicycle, with the odd dicycle embedded in the plane without crossings,
we find an even dicycle in $D$ that is entirely contained in one of the two disks bounded by the odd cycle.

\subsubsection{A brief overview on how to obtain our local structure theorem}

Before we (quite literally) dive into a more detailed discussion of our technique, let us describe, on a more abstract level, the different stages of our approach.
The structural tools developed for our proof are fairly intricate and only bare superficial semblance to those of the Graph Minors Series by Robertson and Seymour.
This part is meant as a gentle introduction and to improve the reader's orientation throughout the following sections while avoiding heavy notation and is therefore kept intentionally vague.

The overall strategy is to iteratively construct an (almost\footnote{We say \textit{almost} because not the entire graph is planar. However, those pieces that are not still behave, essentially, in a planar way.}) embedding of our digraph in the plane by adding a new \textit{strip}, i.\@e.\@ a newly embedded piece of the digraph, at a time.

\vspace{-11pt}
\paragraph{Diving from a wall.}
The first step is to observe that any large enough wall-like infrastructure can be tamed in the sense that it either hosts a large subwall which (almost) embeds into the plane, or a large integral packing of even dicycles is attached to it.
This procedure is repeated throughout our proof in many ways as the final goal of our procedure is to \textit{embed everything that is highly linked to the initial wall in a similar way.}

\begin{figure}[ht!]
    \centering
    \scalebox{0.7}{
    \begin{tikzpicture}[scale=1]

        \pgfdeclarelayer{background}
		\pgfdeclarelayer{foreground}
			
		\pgfsetlayers{background,main,foreground}
			
        \begin{pgfonlayer}{main}
        \node (C) [v:ghost] {};

        \node (L) [v:ghost,position=0:0mm from C] {};

        \node (Lpicture) [v:ghost,position=0:0mm from L] {
            \begin{tikzpicture}[scale=1]

                \pgfdeclarelayer{background}
		          \pgfdeclarelayer{foreground}
			
		          \pgfsetlayers{background,main,foreground}
			
                \begin{pgfonlayer}{main}
                    \node (C) [v:ghost] {};

                    \node (1) [v:ghost,position=140:40mm from C] {};
                    \node (2) [v:ghost,position=270:30mm from 1] {};
                    \node (3) [v:ghost,position=180:5mm from 2] {};
                    \node (4) [v:ghost,position=270:7mm from 2] {};
                    \node (5) [v:ghost,position=315:20mm from 3] {};
                    \node (25) [v:ghost,position=0:3mm from 2] {};
                    \node (6) [v:ghost,position=45:20mm from 5] {};
                    \node (7) [v:ghost,position=180:5mm from 6] {};
                    \node (8) [v:ghost,position=180:2mm from 7] {};
                    \node (9) [v:ghost,position=270:7mm from 7] {};
                    \node (10) [v:ghost,position=270:33mm from 4] {};
                    \node (11) [v:ghost,position=270:15mm from 9] {};
                    \node (12) [v:ghost,position=0:130mm from 10] {};
                    \node (13) [v:ghost,position=0:93.5mm from 11] {};
                    \node (14) [v:ghost,position=90:30mm from 12] {};
                    \node (15) [v:ghost,position=0:5mm from 14] {};
                    \node (16) [v:ghost,position=180:3mm from 14] {};
                    \node (17) [v:ghost,position=90:7mm from 14] {};
                    \node (26) [v:ghost,position=135:20mm from 15] {};
                    \node (18) [v:ghost,position=225:20mm from 26] {};
                    \node (19) [v:ghost,position=90:12mm from 13] {};
                    \node (20) [v:ghost,position=0:2mm from 19] {};
                    \node (21) [v:ghost,position=90:7mm from 19] {};
                    \node (22) [v:ghost,position=90:33mm from 17] {};
                    \node (23) [v:ghost,position=90:15mm from 21] {};
                    \node (24) [v:ghost,position=180:93.5mm from 23] {};
                    \node (30) [v:ghost,position=180:37.6mm from 23] {};
                    \node (29) [v:ghost,position=180:18.3mm from 30] {};
                    \node (27) [v:ghost,position=90:18mm from 30] {};
                    \node (28) [v:ghost,position=180:18.3mm from 27] {};
                    \node (35) [v:ghost,position=270:12mm from 29] {};
                    \node (32) [v:ghost,position=270:12mm from 30] {};
                    \node (31) [v:ghost,position=0:5mm from 32] {};
                    \node (33) [v:ghost,position=180:2mm from 32] {};
                    \node (34) [v:ghost,position=0:2mm from 35] {};
                    \node (36) [v:ghost,position=180:5mm from 35] {};
                    \node (39) [v:ghost,position=315:20mm from 36] {};
                    \node (38) [v:ghost,position=270:7mm from 35] {};
                    \node (37) [v:ghost,position=270:7mm from 32] {};
                    \node (40) [v:ghost,position=180:37.6mm from 13] {};
                    \node (41) [v:ghost,position=180:18.3mm from 40] {};
                    \node (42) [v:ghost,position=270:18mm from 41] {};
                    \node (43) [v:ghost,position=270:18mm from 40] {};

                    \node (l1) [v:ghost,position=135:6.8mm from 24] {};
                    \node (l4) [v:ghost,position=225:6.8mm from 11] {};
                    \node (l9) [v:ghost,position=315:6.8mm from 41] {};
                    \node (l12) [v:ghost,position=45:6.8mm from 29] {};

                    \node (l2) [v:ghost,position=90:6pt from l1] {};
                    \node (l3) [v:ghost,position=180:6pt from l1] {};
                    \node (l5) [v:ghost,position=270:6pt from l4] {};
                    \node (l6) [v:ghost,position=180:6pt from l4] {};
                    \node (l7) [v:ghost,position=0:21mm from l4] {};
                    \node (l8) [v:ghost,position=0:5mm from l7] {};
                    \node (l10) [v:ghost,position=0:6pt from l9] {};
                    \node (l11) [v:ghost,position=270:6pt from l9] {};
                    \node (l13) [v:ghost,position=0:6pt from l12] {};
                    \node (l14) [v:ghost,position=90:6pt from l12] {};
                    \node (l15) [v:ghost,position=180:21mm from l12] {};
                    \node (l16) [v:ghost,position=180:5mm from l15] {};

                    \node (r1) [v:ghost,position=135:6.8mm from 30] {};
                    \node (r4) [v:ghost,position=225:6.8mm from 40] {};
                    \node (r9) [v:ghost,position=315:6.8mm from 13] {};
                    \node (r12) [v:ghost,position=45:6.8mm from 23] {};

                    \node (r2) [v:ghost,position=90:6pt from r1] {};
                    \node (r3) [v:ghost,position=180:6pt from r1] {};
                    \node (r5) [v:ghost,position=270:6pt from r4] {};
                    \node (r6) [v:ghost,position=180:6pt from r4] {};
                    \node (r7) [v:ghost,position=0:21mm from r4] {};
                    \node (r8) [v:ghost,position=0:5mm from r7] {};
                    \node (r10) [v:ghost,position=0:6pt from r9] {};
                    \node (r11) [v:ghost,position=270:6pt from r9] {};
                    \node (r13) [v:ghost,position=0:6pt from r12] {};
                    \node (r14) [v:ghost,position=90:6pt from r12] {};
                    \node (r15) [v:ghost,position=180:21mm from r12] {};
                    \node (r16) [v:ghost,position=180:5mm from r15] {};

                \end{pgfonlayer}{main}

                \begin{pgfonlayer}{background}

                \draw [color=CornflowerBlue,opacity=0.4,line width=12pt] (l2.center) to (l5.center);

                \draw [color=HotMagenta,opacity=0.4,line width=12pt] (r2.center) to (r5.center);

                \draw [color=CornflowerBlue,opacity=0.4, line width=12pt] (l7.center) to (l10.center);

                \node (lbox1) [color=white,regular polygon, regular polygon sides=3,fill,minimum width=36pt,position=180:2.25mm from l8,shape border rotate=270] {};

                \node (lbox2) [color=CornflowerBlue,opacity=0.4,regular polygon, regular polygon sides=3,fill,minimum width=30pt,position=180:2.38mm from l8,shape border rotate=270] {};

                \draw [color=CornflowerBlue,opacity=0.4,line width=12pt] (l6.center) to (l7.center);

                \draw [color=HotMagenta,opacity=0.4, line width=12pt] (r7.center) to (r10.center);

                \node (lbox1) [color=white,regular polygon, regular polygon sides=3,fill,minimum width=36pt,position=180:2.25mm from r8,shape border rotate=270] {};

                \node (lbox2) [color=HotMagenta,opacity=0.4,regular polygon, regular polygon sides=3,fill,minimum width=30pt,position=180:2.38mm from r8,shape border rotate=270] {};

                \draw [color=HotMagenta,opacity=0.4,line width=12pt] (r6.center) to (r7.center);

                \draw [color=CornflowerBlue,opacity=1,line width=12pt] (l11.center) to (l14.center);

                \draw [color=CornflowerBlue,opacity=0.4, line width=12pt] (l15.center) to (l3.center);

                \node (lbox3) [color=white,regular polygon, regular polygon sides=3,fill,minimum width=36pt,position=0:2.25mm from l16,shape border rotate=90] {};

                \node (lbox4) [color=CornflowerBlue,opacity=0.4,regular polygon, regular polygon sides=3,fill,minimum width=30pt,position=0:2.38mm from l16,shape border rotate=90] {};

                \draw [color=CornflowerBlue,opacity=0.4,line width=12pt] (l13.center) to (l15.center);

                \draw [color=HotMagenta,opacity=0.4,line width=12pt] (r11.center) to (r14.center);

                \draw [color=HotMagenta,opacity=0.4, line width=12pt] (r15.center) to (r3.center);

                \node (lbox3) [color=white,regular polygon, regular polygon sides=3,fill,minimum width=36pt,position=0:2.25mm from r16,shape border rotate=90] {};

                \node (lbox4) [color=HotMagenta,opacity=0.4,regular polygon, regular polygon sides=3,fill,minimum width=30pt,position=0:2.38mm from r16,shape border rotate=90] {};

                \draw [color=HotMagenta,opacity=0.4,line width=12pt] (r13.center) to (r15.center);

                \node (m) [v:ghost,position=90:10mm from 39] {};

                \node (l17) [v:ghost,position=90:0.39mm from m] {};
                \node (l18) [v:ghost,position=180:1.8mm from l17] {};
                \node (l19) [v:ghost,position=180:1.8mm from m] {};
                \node (r17) [v:ghost,position=90:0.39mm from m] {};
                \node (r18) [v:ghost,position=0:1.8mm from r17] {};
                \node (r19) [v:ghost,position=0:1.8mm from m] {};

                \node (lbox5) [color=white,isosceles triangle,
	isosceles triangle apex angle=90,fill,rotate=45,minimum width=47pt,position=180:1.5mm from l19] {};

                 \node (rbox5) [color=white,isosceles triangle,
	isosceles triangle apex angle=90,fill,rotate=135,minimum width=47pt,position=0:1.5mm from r19] {};

                \node (lbox6) [opacity=1,color=CornflowerBlue,isosceles triangle,
	isosceles triangle apex angle=90,fill,rotate=45,minimum width=42pt,position=180:1.5mm from l18] {};

                \node (rbox6) [opacity=0.4,HotMagenta,isosceles triangle,
	isosceles triangle apex angle=90,fill,rotate=135,minimum width=42pt,position=0:1.5mm from r18] {};

                \draw [e:main] (2.center) to (1.center) to (28.center) to (27.center) to (22.center) to (17.center);

                \draw [e:main] (7.center) to (24.center) to (29.center) to (30.center) to (23.center) to (21.center);

                \draw [e:main] (25.center) to (2.center) to (3.center) to (5.center) to (6.center) to (7.center) to (8.center);

                \draw [e:main] (4.center) to (10.center) to (42.center) to (43.center) to (12.center) to (14.center);

                \draw [e:main] (9.center) to (11.center) to (41.center) to (40.center) to (13.center) to (19.center);

                \draw [e:main] (20.center) to (19.center) to (18.center) to (26.center) to (15.center) to (14.center) to (16.center);

                \draw [e:main] (29.center) to (35.center);

                \draw [e:main] (30.center) to (32.center);

                \draw [e:main] (38.center) to (41.center);

                \draw [e:main] (37.center) to (40.center);

                \draw [e:main] (34.center) to (35.center) to (36.center) to (39.center) to (31.center) to (32.center) to (33.center);

                \draw [e:main,dashed] (28.center) to (29.center);
                \draw [e:main,dashed] (27.center) to (30.center);
                \draw [e:main,dashed] (41.center) to (42.center);
                \draw [e:main,dashed] (40.center) to (43.center);

                
                \end{pgfonlayer}{background}
        
                \begin{pgfonlayer}{foreground}
                    
                \end{pgfonlayer}{foreground}
        
            \end{tikzpicture}
        };


			
			
                    

                


                

        \end{pgfonlayer}{main}
        
        \begin{pgfonlayer}{foreground}
        \end{pgfonlayer}{foreground}

        \begin{pgfonlayer}{background}
        \end{pgfonlayer}{background}
        
    \end{tikzpicture}
    }
    \caption{A large collection of concentric dicycles in the plane together with a transaction on the boundary of the disk defined by the dicycles.
    The transaction splits the disk into two new disks, $\Delta_1$ on the left, and $\Delta_2$ on the right.
    Notice that, when forgetting the orientations of the edges, we can surround both disks with their own family of concentric cycles.
    Moreover, the family surrounding $\Delta_2$, depicted in \textcolor{HotMagenta}{magenta}, is again a family of concentric dicycles.
    Hence, $\Delta_2$ is a \textit{circle}.
    The family surrounding $\Delta_1$ on the left, depicted in \textcolor{CornflowerBlue}{blue}, however, now consists of cycles that decompose into pairs of parallel directed paths sharing their start and endpoint.
    Thus, $\Delta_1$ is a \textit{diamond}.}
    \label{fig:diamondsandcircles}
\end{figure}
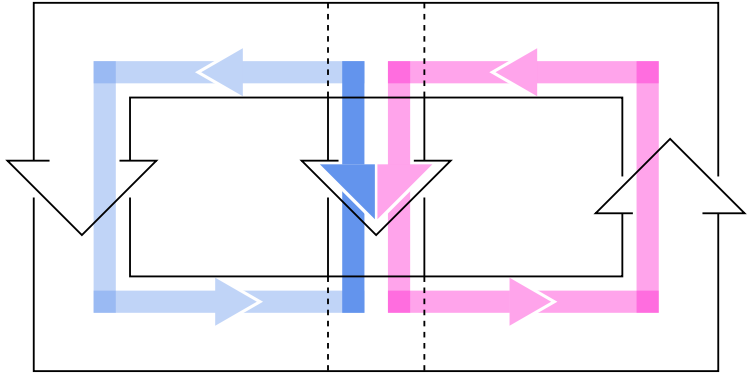

\vspace{-11pt}
\paragraph{Different shapes of currents: Diamonds and circles.}
As we mentioned above, a chord in a dicycle splits the cycle into exactly one smaller dicycle and an object that, if we forget about the orientations, is still a cycle, but now it consists of two internally disjoint directed paths with a common start and a common endpoint.
A similar situation is encountered when instead of a chord and a single dicycle, we deal with a large family of concentric cycles in the plane and a large linkage crossing from one part of this cycle family to another.
When seen as an annulus that gets split by a ribbon, this yields another annulus (or a new family of concentric cycles): a \textit{circle} on one side, and a \textit{diamond} shape on the other.
See \cref{fig:diamondsandcircles} for a schematic illustration of such a situation.
As it turns out, these two types of objects exhibit different behaviour and need to be handled separately.

\vspace{-11pt}
\paragraph{The outlines of maelstroms.}
What both, circles and diamonds, have in common is that they enclose disks in whose interior the graph is not necessarily plane.
These disks are where the even dicycles will be trapped.
In the Graph Minors Series such a disk is called a \textit{vortex}.
However, because of the unfriendly behavior of directed linkages, a single circle or diamond may host a plethora of such disks which are surrounded by smaller and smaller ``whirls'' of infrastructure.
This forces us to maintain a large backlog of infrastructure to ensure that any even dicycle inside still stays trapped: the \textit{maelstrom outlines}.

\vspace{-11pt}
\paragraph{A glimpse of the abyss.}
As discussed above, maelstroms come in two different shapes: circles and diamonds.
Circles are somewhat well behaved in the sense that their behaviour is very close to the behaviour of ``nests'' and vortices from the undirected graph minors theory.
Hence, here a lot of intuition applies.
Diamonds, however, are a completely different beast.
Their structure is inherently \textit{oriented} by the shape of the diamond itself which makes it necessary to distinguish two different types of \textit{transactions}\footnote{\emph{Transactions} are linkages that start and end on different sections of the boundary of a disk.}: Those that agree with the orientation of the diamond and those that do not.
The latter pose a major problem as they form, together with a piece of the diamond infrastructure, new circles which are now trapped inside the diamond.
This induces a somewhat recursive structure as these circles might trap additional (or possibly all) even dicycles within their boundaries.

\vspace{-11pt}
\paragraph{Charting the currents.}
With all necessary tools to handle the situations above at hand, we can now begin to fully ``map out'' the landscape of our digraph to identify the regions which contain even dicycles.
Broadly speaking, all of these tools have two possible outcomes: Either they provide the pieces necessary to continue the construction of our (almost) embedding, or they yield a large quarter-integral packing of $t$ even dicycles, where $t \in \N$ is given to us beforehand.
In the latter case we may immediately terminate the entire process as we have found the object we desired.
Assuming we never encounter the large packing, our process is as follows:
\begin{enumerate}[label=\arabic*)]
\item Decide that the directed treewidth of $D$ is small or find a large \textit{odd wall} and define the first two circle maelstroms.
\item Given any maelstrom outline decide if there exists a large \textit{transaction} on its boundary. If yes continue, if no jump to step 4).
\item Within the large transaction find an \textit{odd strip} and augment the current (almost) embedding.
Some cases may appear.
\begin{enumerate}
    \item The strip splits a \textit{circle maelstrom} into \textbf{two areas}, each with an even dicycle.
    \item The strip splits a \textit{diamond maelstrom} into \textbf{two areas}, each with an even dicycle, and it \textit{agrees} with the orientation of the diamond.
    \item The strip splits a \textit{diamond maelstrom} into \textbf{two areas}, each with an even dicycle, but it \textit{disagrees} with the orientation of the diamond.
    \item The strip splits a \textit{circle maelstrom} into exactly \textbf{one area} with an even dicycle and an odd area.
    \item The strip splits a \textit{diamond maelstrom} into exactly \textbf{one area} with an even dicycle and an odd area, it agrees with the orientation of the diamond, and the area with the even dicycle is itself a diamond.
    \item The strip splits a \textit{diamond maelstrom} into exactly \textbf{one area} with an even dicycle and an odd area, and it either agrees with the orientation of the diamond but the area with the even dicycle is a circle, or it disagrees with the orientation of the diamond.
    \item After finding the strip, \textbf{no even dicycle remains}.
\end{enumerate}
The ways we recurse in any of these cases differ in varying degrees but, holistically speaking, in each of the cases a) to f) we define one (or two) new maelstroms, now enclosing a smaller area, and revisit this maelstrom starting with step 2).
In case g) however, we have fully exhausted the maelstrom and may forget it completely.
Notice that, at any given point in time, the number of ``active'' maelstroms equals the size of a quarter-integral packing of even dicycles we can guarantee.
This means that, if we ever reach $t$ or more active maelstroms our procedure terminates as we have found our packing.
\item Finally, the maelstrom's infrastructure is exhausted and no new transactions of large enough order can be found.
This means we have reached a situation which, in the Graph Minors Series of Robertson and Seymour, is called \emph{bounded depth}.
\end{enumerate}

\vspace{-11pt}
\paragraph{The death of a vortex.}
The previous step has now either found a large quarter-integral packing of even dicycles and has thus terminated, or it has produced an (almost) embedding of the digraph into the plane with a small number of ``active'' maelstroms, each of which has bounded depth.
What remains is to deal with this situation.
Here we employ a strategy that was first used by Thilikos and Wiederrecht in their dichotomy result for counting perfect matchings \cite{thilikos2022killing}.
As each maelstrom has bounded depth, any \textit{segment} on its boundary can be separated (in both directions) from the rest of the boundary via deleting a small number of vertices.
This means that we can now either isolate many segments, each hosting their own private even dicycle and thus yielding a large quarter-integral packing, or we can delete a bounded number of vertices to separate the entire maelstrom boundary from all even dicycles inside the maelstrom.
This is the final step in obtaining the local structure theorem and completes the difficult parts of the proof.

\subsubsection{A more detailed discussion.}
The skeleton of our proof is inspired by the proof of the Graph Minor Structure Theorem constructed by Kawarabayashi, Thomas, and Wollan \cite{kawarabayashi2020quickly}.
Though, aside from the proof of the global structure theorem in the very last section, we have to replace almost every major tool from their proof entirely and in the few places where a section of our paper clearly maps to one of theirs, their proofs make heavy use of the niceties of the undirected setting, forcing us to come up with entirely new proof strategies and develop and whole new set of tools.

\paragraph{From a local to a global structure theorem.}
The proof of our global structure theorem in \Cref{sec:globalstructure} uses the many already developed abstract tools for digraph structure theory and thanks to the successes there, the proof will seem very familiar to readers who know the literature surrounding the Graph Minors Structure Theorem (in fact, the same basic idea has already appeared in \cite{robertson1991graph}).
The main idea here is to use a balanced separator type argument to construct a directed tree-decomposition in a way very similar to Reed's seminal algorithm for approximating (undirected) treewidth \cite{reed1992finding}.
That is, we iteratively build a tree-decomposition by splitting a separator set in a balanced way.
Whenever we encounter a situation where such a balanced split is impossible we have found a witness of large directed treewidth.
Ideally, this allows us to evoke our local structure theorem to which the remaining graph only attaches via small interfaces.
In our case, all of these interfaces are governed by the apex set provided by our local structure theorem.
This is the reason why the resulting odd directed tree-decomposition will be \textit{strong} in the end.
Although some arguments require translation from the undirected to the directed setting, this procedure is fairly standard in structural graph theory and does not pose a big obstacle.
Hence, similar to the situation in \cite{kawarabayashi2020quickly}, the difficulty of our project lies in proving a theorem which provides structure in relation to a fixed wall, commonly known as a local structure theorem.

Once all pieces are in place, the global structure theorem reads as follows.

\begin{theorem}\label{thm:mainthm4}
There exist a functions $f,g\colon\mathbb{N}\to\mathbb{N}$ such that for every integer $k$ and every digraph $D$ either
\begin{enumerate}
    \item $D$ has a quarter integral packing of $k$ even dicycles, or
    \item $D$ has a strong odd directed tree-decomposition of width at most $f(k)$.
\end{enumerate}
Moreover, there exists an algorithm that finds one of the two outcomes above in time $g(k)|V(D)|^{\mathcal{O}(1)}$.
\end{theorem}

\paragraph{Almost planar embeddings: \textit{Odd decompositions}.}
We can rely on the Directed Grid Theorem by Kawarabayashi and Kreutzer \cite{KawarabayashiKreutzer2015DirectedGrid} to find a large wall\footnote{Proving that even dicycles have the half-integral Erd\H{o}s-P\'osa property in digraphs with low directed treewidth is quite easy.}.
Starting there, we would like to find a ``flat'' wall, but as mentioned prior, the existing definitions of flatness in the directed setting do not actually function like planar graphs.
To circumvent this, we define a framework called an \emph{odd decomposition} in \Cref{sec:notation}.
The basis of the odd decomposition is formed by a part of the graph that is more or less free of even dicycles and already has an embedding that is strongly planar.
We then allow for three types of non-planarity in certain areas, two of which replace the small separations used to define flat embeddings in \cite{kawarabayashi2020quickly} and the third represents the unruly areas that we need to decompose further, which we call \emph{maelstroms}, with the analogous areas being called vortices in \cite{kawarabayashi2020quickly}.

\paragraph{A relative from matching theory: Obstructions to odd decompositions yield even dicycles.}
As in the undirected setting, we want an odd decomposition to behave like a planar embedding of a graph with respect to the part of the graph that does not lie in the non-planar areas.
This is guaranteed by the Two Paths Theorem \cite{jung1970verallgemeinerung,seymour1980disjoint,shiloach1980polynomial,thomassen19802,robertson1990graph} in the undirected setting, but there is no analogous result in the directed setting.
The fact that an odd decomposition nonetheless exhibits the behaviour we need thus already has to rely on some very non-trivial facts.
In particular we use a result analogous to the Two Paths Theorem which relates to non-even digraphs and is implied by the results of Giannopoulou and Wiederrecht \cite{giannopoulou2021two}.
We note that their result was proven in the matching theoretic setting and the connection between (structural) matching theory and structural digraph theory is in many ways the engine behind this entire project.
Whilst we formulate most of our proofs for digraphs, much of the conceptual work behind them stems from insights we gained through the use of matching theory.

\paragraph{An odd wall and the first two maelstroms.}
We then reprove major parts of the Directed Flat Wall Theorem \cite{Giannopoulou2020DirectedFlatWall} in \Cref{sec:oddwall} to either find a large half-integral packing of even dicycles or a set of vertices that we can delete to find a large wall with an associated odd decomposition in which this wall lies essentially outside of all maelstroms.
This proof does not only resemble the directed version of the Flat Wall Theorem, it also evokes ideas and techniques from the matching theoretic analogue of Giannopoulou and Wiederrecht \cite{giannopoulou2021flat}.
Once we have obtained the odd decomposition of the wall, we then consider the two sides of the graph that the wall intuitively separates, since it is cylindrical.
Within the wall we choose two large disjoint walls which each surround one of these two sides and then try to decompose the graph starting from there.

The \textit{odd wall theorem} reads as follows.
We provide a formal definition of an \textit{odd wall} later on.
Informally, an odd wall is a wall for which the strong component of $D$, after deleting the two boundary cycles of the wall, which contains the interior of the wall has an (almost) embedding in the plane and is free of even dicycles.

\begin{theorem}\label{thm:mainthm5}
There exist functions\footnote{We write $\mathbb{N}^{(2)}$ to indicate that some of the functions need two arguments instead of one.} $a,f,g\colon\mathbb{N}^{(2)}\to\mathbb{N}$ such that for all positive integers $r,k$ and all digraphs $D$ that contain a cylindrical $f(r,k)$-wall $W$, there either exists an integral packing of $k$ even dicycles, or a set $A\subseteq V(D)$ of size at most $a(k)$ and a cylindrical $r$-wall $W'\subseteq W-A$ such that $W'$ is an odd wall in $D-A$.

Moreover, there exists an algorithm that, given $D$ and $W$, finds either the packing or $A$ and $W'$ in time $g(k)|V(D)|^{\mathcal{O}(1)}$.
\end{theorem}

\paragraph{Refining a maelstrom.}
As in \cite{kawarabayashi2020quickly}, we take a large linkage running through the wall and the undecomposed part, decompose this linkage in a way that is compatible with the existing odd decomposition in \Cref{sec:transaction}, and split the maelstrom into several new maelstroms in \Cref{sec:buildoutline}.
However, unlike in \cite{kawarabayashi2020quickly}, the structure for the new maelstroms is not uniformly the same, forcing us to introduce several types of valid structures around a maelstrom.
Additionally, we need to show that the number of maelstroms we produce whilst we keep refining the decomposition cannot grow arbitrarily.

Formulated in an abstract way, this leads to us proving in \Cref{sec:shifting} that even dicycles are local with respect to the structure we find.
For planar graphs this is not too difficult, but the proof that this is still true for odd decompositions is surprisingly involved.
Aside from this, we also need to show that this is true even if the even dicycle behaves in a non-planar way with respect to our existing decomposition, which involves having to invoke an even more specific version of the two paths theorem for non-even digraph implied by a result of Giannopoulou, Thilikos, and Wiederrecht \cite{Giannopoulou2023ExcludingSingleCrossing}.
Large parts of \Cref{sec:shifting} take place in the matching theoretic setting and finding direct proofs in the directed setting would likely involve even more unpleasant case analysis.

With these tools in hand, we integrate the linkage we find into the odd decomposition and, using the decomposed part of the linkage, we then split our maelstrom into several new maelstroms, each with their own infrastructure.
Both of these steps incur sets of vertices that need to be removed to proceed and owing to the complications of the directed setting, the infrastructure we find for the different maelstroms is not guaranteed to be disjoint, but interacts at worst in a third-integral fashion, making it amenable to our project of proving the quarter-integral Erd\H{o}s-P\'osa property for even dicycles.

\paragraph{Closing the rifts at the bottom of the abyss.}
It would at this point be possible for us to prove a local structure theorem that leaves us with an odd decomposition of whatever is attached to the wall we start with, once we remove a set of vertices of bounded size, whilst leaving a small number of maelstroms on which no large linkage can be found\footnote{We are able to always find a decomposition on the sphere and can avoid having to increase the genus of the surface we find our decomposition in, contrasting again with \cite{kawarabayashi2020quickly}.}.
However, we go a step further and also decompose these maelstroms using ideas by Thilikos and Wiederrecht \cite{thilikos2022killing} in \Cref{sec:killmaelstroms} to ultimately show in \Cref{sec:localstructure} that one can remove a small number of vertices and in return make the entire graph attached to what remains of the initial wall odd.
While the ideas in \cite{thilikos2022killing} are quite elegant, the directed setting and the structure we ultimately find in our proof is so complicated that a translation of their methods only resolves the base case in our setting and we are forced to provide substantially more analysis to actually finish decomposing our maelstroms.
Let us formulate a rough version of our local structure theorem without many of the more involved technical definitions necessary to facilitate to full proof.

\begin{theorem}\label{thm:mainthm6}
There exist functions $a,f,g\colon\mathbb{N}^{(2)}\to\mathbb{N}$ such that for all pairs of positive integers $k$ and $r$ and all digraphs $G$ with a cylindrical $f(r,k)$-wall $W$ there either exists
\begin{enumerate}
    \item a quarter-integral packing of $k$ even dicycles, or
    \item a set $A\subseteq V(D)$ of size at most $a(k)$ and a cylindrical $r$-wall $W'\subseteq W-A$ such that the strong component of $D-A$ which contains $W'$ does not contain an even dicycle.
\end{enumerate}
Moreover, there exists an algorithm that, given a digraph $D$ and a wall $W$ as above, finds one of the two outcomes in time $g(k)|V(D)|^{\mathcal{O}(1)}$.
\end{theorem}

\subsection{From the global structure to quarter-integral Erd\H{o}s-P\'osa}

For the sake of completeness let us also provide a short proof that explains how one can obtain \cref{thm:mainthm1} from \cref{thm:mainthm4}.

\begin{proof}[Proof of \Cref{thm:mainthm1}]
Fix some integer $k$ and let $D$ be any digraph.
\cref{thm:mainthm4} now either provides a quarter-integral packing of $k$ even dicycles, in which case we are done, or it finds a strong odd directed tree-decomposition $\mathcal{T}=(T,\alpha,\beta,\gamma)$ of bounded width.
Let $f'\colon\mathbb{N}\to\mathbb{N}$ denote the function from $\cref{thm:mainthm4}$.
As $\mathcal{T}$ is strong we may assume that $\gamma(dt)\subseteq \alpha(t)$ for all $dt\in E(T)$ since $\alpha(t)$, by definition, already meets all dicycles which contain a vertex of $\Gamma(t)$ and a vertex that does not belong to $\Gamma(t)$.

For every $e\in E(T)$ the digraph $D$ is split into digraphs $D_1$ and $D_2$ corresponding to the two digraphs induced by the vertices in the bags of the two components of $T-e$ after deleting the vertices of $\gamma(e)$.
Suppose there exists some $e\in E(T)$ such that both, $D_1$ and $D_2$, contain an even dicycle.
In this case we may proceed by induction for $k-1$ on both $D_i$.
If both $D_i$ return an even dicycle transversal $X_i$ of order $f(k-1)$ we have found a set $S\coloneqq X_1\cup X_2\cup \gamma(e)$ of size at most $2f(k-1)+f'(k)$ such that $D-S$ has no even dicycle.
Otherwise, at least one of $D_1$ or $D_2$ contains a quarter-integral packing of $k-1$ even dicycles.
Together with the one even dicycle guaranteed to exist in the other digraph this provides a quarter-integral packing of $k$ even dicycles.

Similarly, we are done immediately if there exists some $e\in E(T)$ such that $D-\gamma(e)$ does not contain an even dicycle.

Hence, for each $e\in E(T)$ exactly one of $D_1$ and $D_2$ contains an even dicycle.
Let $\vec{T}$ be an orientation of $T$ obtained by orienting every edge of the tree towards the side that contains the even dicycle.
Observe that we may assume that there cannot exist a $t\in V(T)$ together with edges $(t,d)$ and $(t,d')\in E(\vec{T})$ as otherwise deleting $\gamma(dt)\cup \gamma(d't)$ would remove all even dicycles from $D$.
Thus, there exists a unique $t\in V(T)$ which is a \textit{sink} of $\vec{T}$, i.\@e.\@ $t$ is the unique vertex of $V(T)$ with only incoming edges.
This implies that $D-\alpha(t)$ does not contain an even dicycle and completes our proof.
\end{proof}

\subsection{Odd directed treewidth}

In this subsection we provide proof sketches for the other two main theorems, that is the approximation algorithm for odd directed treewidth, namely \cref{thm:mainthm2} and the algorithm for $t$-DDPP on digraphs of bounded odd directed treewidth, that is \cref{thm:mainthm3}.

\subsubsection{A sketch for \Cref{thm:mainthm2}}

Notice that the second outcome of \cref{thm:mainthm4} is exactly a strong odd directed tree-decomposition of bounded width.
This means that all we need to discuss is that the first outcome of \cref{thm:mainthm4} indeed provides a proof that $\mathsf{odtw}(D)>k$.

For this we need to dive back into the proof of the \textit{local} structure theorem.
Hence, we will only provide a rough sketch of the proof for the sake of brevity.
The proof of the local structure theorem, as discussed above, iteratively constructs an (almost) embedding of our digraph into the plane.
This is done by starting from a large cylindrical wall, which is a witness of large directed treewidth, and then iteratively attaching large linkages to the wall.
Overall there are four potential sources of quarter-integral packings of even dicycles throughout this process:
\begin{enumerate}
    \item An application of the \textit{odd wall theorem} (or one of its relatives), that is \cref{thm:mainthm5}, finds, as a butterfly minor, a large complete digraph.
    In the statement of \cref{thm:mainthm5} this outcome is already refined as a large complete digraph contains a large integral packing of odd bicycles which, in turn, contains a large integral packing of even dicycles.
    However, the clique itself is also a witness of large odd directed treewidth and thus this outcome would yield the lower bound for \cref{thm:mainthm2}.
    \item The second source for a quarter-integral packing originates from a substep of finding new \textit{odd transactions} and further ``splitting'' the current maelstrom.
    These even dicycles appear if we have found our candidate for the odd transaction but now find a large number of disjoint paths that ``jump'' over the transaction from one side of the maelstrom to the other.
    Such a jumping transaction yields a large quarter-integral packing of odd bicycles.
    Moreover, since this situation arises from iteratively attaching large linkages to a huge wall $W$, for each of these odd bicycles $C$ there exists a large number of disjoint paths from $C$ to $W$ and from $W$ to $C$.
    Hence, for any small set of vertices $S$, at least one of the concentric cycles of the wall and at least one of the odd bicycles must be left intact and, moreover, both of them must belong to the same strong component of $D-S$.
    This implies that the even dicycles within the odd bicycles cannot be separated from the wall by a small set of vertices and thus, $\mathsf{odtw}(D)$ must be large.
    \item The third source for our quarter-integral packing would be if, at any point in time, we have at least $k$ ``active'' maelstroms.
    Similar to case ii) above, each maelstrom, and therefore each even dicycle within, is highly connected to the wall and thus, by reiterating the arguments above, one can see that a small set of vertices can never separate the wall from all these even dicycles which yields the desired lower bound.
    \item Finally, the fourth source would be the step of ``killing the vortex''.
    In this step a single maelstrom may yield a large quarter-integral packing.
    However, the even dicycles found in this step all have the property, by construction, that they cannot be separated from the maelstrom-infrastructure by a small set of vertices.
    As discussed above, the maelstrom itself is highly linked to the wall $W$ we started with and thus, so must be all the even dicycles we have found.
    Hence, as before, we obtain the desired lower bound on the odd directed treewidth.
\end{enumerate}
The discussion above shows that, in fact, the proof of our local structure theorem finds, whenever it produces a quarter-integral packing of even dicycles, an object usually called a \textit{bramble} with the additional property that each element of the bramble also contains an even dicycle.
Brambles are among the most well studied objects which yield lower bounds to (directed) treewidth (see for example \cite{Reed1997Brambles,masarik2022constant}).
This particular type of bramble can now be seen to yield a lower bound for odd directed treewidth.

\subsubsection{A sketch for \Cref{thm:mainthm3}}

We conclude by briefly describing our algorithm for the $t$-DDPP for graphs of bounded odd directed treewidth.
In light of \cref{thm:mainthm2} this algorithm becomes a straight forward generalisation of the original algorithm for the $t$-DDPP by Johnson et al.\@ \cite{Johnson2001DirectedTreewidth} for bounded directed treewidth.
So let us begin by briefly summarising this algorithm instead.

Given a directed tree-decomposition\footnote{As explained earlier, if we require $\alpha(t)=\Gamma(t)$ for all $t\in V(T)$ in an odd directed tree-decomposition, we obtain exactly the notion of directed treewidth. In this case we can simply omit the $\alpha$-sets as they do not provide any additional information. The result is indeed the original directed tree-decomposition.} $(T,\beta,\gamma)$ of bounded width we compute partial solutions to the $t$-DDPP ``bottom up''.
To be able to do this we need a few observations.
First, let us assume $(T,\beta,\gamma)$ to be rooted at some node $r\in V(T)$.
Second, we may also assume that the digraph induced by the bags in each subtree is strongly connected.
This can be achieved by further refining the decomposition and ``splitting off'' new branches for each strong component within the subtree while simply copying the guards, i.\@e.\@ the $\gamma$-sets.
\begin{enumerate}
    \item Since $|\gamma(dt)|\leq k$ for all edges $dt\in E(T)$, the total number of times a potential solution to $t$-DDPP can enter and leave a subtree is bounded by $k$.
    Hence, the problem of enumerating possible terminal sets for partial solutions to the $t$-DDPP within the subtree is also bounded.
    However, this bound is not simply $k+t$ since $\gamma(dt)$ does not fully separate the vertices within a subtree from those outside.
    Instead, it blocks any path that \textit{leaves and returns}.
    Hence, the total number of terminals for partial solutions is roughly $|V(D)|^{t+k}$.
    \item As $|\Gamma(t)|\leq k$ we can solve each instance of a partial solution exhaustively in the leaves.
    \item Finally, since $|\Gamma(t)|\leq k$ and the number of partial solutions in each subtree is at most $|V(D)|^{t+k}$, we are able to merge partial solutions of siblings and then extend them to partial solutions for the entire subtree in time $|V(D)|^{\mathcal{O}(k+t)}$ again by exhaustively trying all possible extensions through the bag of the new root in the subtree.
\end{enumerate}

Since \cref{thm:mainthm2} provides us with a \textit{strong} odd directed tree-decomposition, we may assume that $\gamma(dt)\subseteq \alpha(t)$ for all $dt\in E(T)$.
This assumption is very strong and allows to to maintain almost all properties that enabled the algorithm for $t$-DDPP for ordinary directed tree-decompositions.
The only difference is not that instead of $\Gamma(t)\setminus\alpha(t)$ being empty as before, we now have that $\Gamma(t)\setminus\alpha(t)$ consists of strong components of $D-\alpha(t)$ which do not contain even dicycles.
Hence, we cannot solve, nor extend, the partial solutions within the subtree exhaustively.
But fear not, as the odd directed tree-decomposition $(T,\alpha,\beta,\gamma)$ provided by \cref{thm:mainthm2} is strong, let $k$ denote its width, we may still enumerate all possible choices for terminal pairs of partial solutions within the subtree in time $|V(D)|^{\mathcal{O}(t+k)}$.
We may then use the algorithm from \cite{Giannopoulou2023ExcludingSingleCrossing} to solve each instance of such a partial problem in time $|V(D)|^{\mathcal{O}(t+k)}$.
By replacing each exhaustive search in the algorithm above with this subroutine, our algorithm is complete.

\subsection{The structure of the paper}

We conclude this introduction with a comprehensive guide to the rest of the paper.

The key definitions together with an extended introduction to the \textsc{Even Dicycle Problem} can be found in \cref{sec:notation}.
Here we also provide formal definitions of our almost embeddings, odd decompositions, maelstroms, and their outlines.

With the definitions of \cref{sec:notation} we take the time to revisit the main points of this introduction in a more formal way and provide a second, slightly different proof for \cref{thm:mainthm1} in \cref{sec:statementandoverview} within a more formal framework.

Within \cref{sec:oddwall} we provide a full proof of the \textit{Odd Wall Theorem}, namely \cref{thm:mainthm5}.
This proof requires several highly technical definitions and tools which are also to be found within this section.

The proof for the Odd Wall Theorem has no direct need for the more involved tools from structural matching theory.
For the rest of our proof, however, we require the two main tools mentioned earlier:
The matching theoretic \textbf{Two Paths Theorem} and the \textbf{Shifting Argument} which allows us to trap even dicycles within certain areas of our almost embedding.
In \Cref{sec:shifting} we provide a short introduction to structural matching theory and provide the necessary material and proofs to establish these two essential tools.
This section also contains the proofs necessary to ``trap'' the even dicycles within our maelstroms.

The outline of a maelstrom is somewhat unwieldy to handle, but sadly a necessary consequence of some later arguments.
However in most situation, we can provide a simpler, more refined structure to which several of the tools from \Cref{sec:shifting} can be applied.
This refined structure is called the \emph{rim} of a maelstrom.
In \Cref{sec:localising} we show that such a rim always exists in an outline, though this requires sacrificing some infrastructure.

The next crucial step is to show that, whenever we find a large transaction on the boundary (or \textit{society}) of a maelstrom outline, we either find a large quarter-integral packing of even dicycles or, possibly after deleting a small number of vertices, a new strip that can be added to our almost embedding.
The purpose of \Cref{sec:transaction} is to prove the corresponding statements.

Once we have found a new strip to be incorporated into our almost embedding, we have to discuss how a maelstrom can be split by this strip and how new outlines for the resulting maelstroms can be constructed.
This is by far the most technical part of the proof and contained in \Cref{sec:buildoutline}.

By now we know how to find new strips within large transactions on the society of a maelstrom outline and we have developed the tools to embed these strips into our almost embedding, thereby dividing the maelstrom into smaller ones.
The final piece of the puzzle is now to deal with maelstrom outlines of bounded depth.
This is the ``vortex killing'' part and takes place in \Cref{sec:killmaelstroms}.

With this, all pieces for the proof of the local structure theorem are in place.
The proof proceeds by induction where our progress measure is, in part, the number of distinct ``active'' maelstroms we have constructed so far.
Two last obstacles need to be overcome at this point:
\begin{enumerate}
    \item In the case where we split a circle maelstrom and only one side contains an even dicycle, or where we split a diamond maelstrom and only the diamond contains an even dicycle, we still ``pay'' a small amount of vertices that go into the apex set.
    Since, in these situations we do not make progress in the sense that we can guarantee an additional even dicycle for our packing, we need to ensure that this situation does not keep on adding apex vertices indefinitely.
    This requires additional technical arguments.
    \item Finally, if we split a diamond maelstrom, but the only even dicycle remaining attaches to a transaction that traverses the diamond ``in reverse'' (we call such an object an \emph{eddy}), we have to readjust the entire odd decomposition of the maelstrom.
    In this case we must open up the recursive structure of our almost embedding and get drawn closer to the \textit{abyss}.
\end{enumerate}
The necessary arguments to deal with the two situations above, and the complete proof of \cref{thm:mainthm6} can be found in \Cref{sec:localstructure}.

All that is left is to leverage the findings of \Cref{sec:localstructure} into the global structure theorem.
The proof of \cref{thm:mainthm4} can be found in \cref{sec:globalstructure}.

\newpage
\tableofcontents

\thispagestyle{empty}
\newpage

\section{Notation and an introduction to structural digraph theory}\label{sec:notation}

\subsection{Preliminaries}\label{subsec:prelim}

Given two integers $i,j\in\Z$, we denote by $[i,j]$ the set $\CondSet{z\in\Z}{i\leq z\leq j}$.
Please note that $[i,j]=\emptyset$ if $i>j$.
If we need an interval $[x,y]\subseteq\R$, we will explicitly write $[x,y]_{\R}$ to avoid ambiguity.
In case $x=1$, we abbreviate the above notation to $[y]$.

We mostly deal with digraphs, but in some contexts we also make use of undirected graphs.
To avoid confusion we stick to the convention of naming our digraphs $D$ and our undirected graphs $G$.
Whenever we break this convention due to a lack of symbols we make it clear explicitly.
All graphs and digraphs considered in this paper are simple, that is they do not contain loops nor parallel edges.
The edges $(u,v)$ and $(v,u)$ in a digraph do not count as parallel and we call the digraph $(\Set{u,v},\Set{(u,v),(v,u)})$ a \emph{digon}.

If $G$ is a (di)graph we denote by $2^G$ the set of all subgraphs of $G$.

A \emph{path} $P$ of \emph{length} $\ell\geq 0$ in a digraph $D$ is a sequence of distinct vertices $(v_0,v_1,\dots,v_{\ell})$ such that for each $i\in[0,\ell-1]$ the edge $(v_i,v_{i+1})$ exists in $\E{D}$.
Note that in case $\ell=0$ the path $P$ consists of exactly one vertex.
We call the vertex $v_0$ the \emph{tail} of $P$ and $v_{\ell}$ its \emph{head}.
In a slight abuse of notation, we identify the path $P$ and the subgraph $(\Set{v_0,v_1,\dots,v_{\ell}},\CondSet{(v_i,v_{i+1})}{i\in[0,\ell-1]})$ of $D$.
If $P$ is a path and $u,v\in\V{P}$, we denote by $uPv$ the unique subpath of $P$ starting in $u$ and ending in $v$.
In case $u$ comes after $v$ along $P$ the path $uPv$ is empty.
We will let the transitive closure of the edge set of a directed path $P$ define the \emph{order of a directed path} and say that a vertex \emph{$u$ appears before $v$ on $P$} to say that $u$ is lower in the order of $P$ than $v$.

A \emph{directed cycle}, or \emph{dicycle} $C$ of \emph{length} $\ell\geq 2$ in $D$, is a sequence $(v_0,v_1,\dots,v_{\ell-1})$ of distinct vertices such that for each $i\in[0,\ell-1]$ the edge $(v_i,v_{i+1})$ exists, where $v_{\ell}=v_0$.
As with paths, we identify the vertex sequence $C$ with its corresponding subgraph in $D$.
In case $\ell$ is even, we say that $C$ is an \emph{even dicycle}, otherwise $C$ is an \emph{odd dicycle}.

\begin{definition}[Fractional packing]
    Let $D$ be a digraph and let $t, n \geq 1$ be integers.
    A \emph{$\nicefrac{1}{n}$-integral packing} of $t$ even dicycles is a family of $t$ even dicycles such that any vertex of $D$ is contained in at most $n$ even dicycles in the family.
    We use the same conventions for packings of arbitrary digraphs.
    For small $n$ it is common to write out the fraction, such as with \emph{half-, third-, and quarter-integral packings}.
\end{definition}

Let $D$ be a digraph and $X,Y\subseteq\V{D}$ a path $P$ is an \emph{$X$-$Y$-path} if its tail lies in $X$, its head lies in $Y$, and it is otherwise disjoint from $X\cup Y$.
Similarly, if $H_1$ and $H_2$ are, not necessarily distinct, subgraphs of $D$ then $P$ is an \emph{$H_1$-$H_2$-path} if it is a $\V{H_1}$-$\V{H_2}$-path and no edge of $P$ is contained in $\E{H_1}\cup \E{H_2}$.
If $H_1 = H_2 = H$, we sometimes call an $H_1$-$H_2$-path an $H$-path.

\begin{definition}[Linkage]\label{def:linkage}
Let $D$ be a digraph and $k \in \N$ be positive.
A \emph{$k$-linkage} is a family of $k$ pairwise disjoint directed paths in $D$.
Given two vertex sets $X$ and $Y$, a family of $k$ pairwise disjoint $X$-$Y$-paths is called a \emph{$X$-$Y$-linkage} and a family of $k$ $X$-$Y$-paths such that no three paths share a vertex is a called a \emph{half-integral $X$-$Y$-linkage}.
\end{definition}

In yet another slight abuse of notation, we identify the linkage $\mathcal{L}$ and the subgraph of $D$ obtained from the union of all paths in $\mathcal{L}$.
In general, given a set of graphs $\mathcal{H}$, we define $V(\mathcal{H}) = \cup_{H \in \mathcal{H}} V(H)$ and $E(\mathcal{H}) = \cup_{H \in \mathcal{H}} E(H)$.

\subsection{The Even Dicycle Problem}\label{sec:evendicycleproblem}

The \emph{Even Dicycle Problem} is the decision problem of determining if a digraph contains an even dicycle.
It is equivalent to many algorithmic problems involving a plethora of algebraic and combinatorial objects (see \cite{McCuaig2004Polya} for an overview).
The oldest known equivalent problem is known as \emph{P\'olya's Permanent Problem} which can be traced back to a question about the permanent of a matrix \cite{polya1913aufgabe}.
The complexity of the Even Dicycle Problem was open for almost a century until McCuaig \cite{McCuaig2004Polya} and Robertson, Seymour, and Thomas \cite{Robertson1999PermanentsPfaffianOrientations} independently announced a solution in 1999 \cite{mccuaig1997permanents}.
This solution is built upon a powerful structural theorem which was obtained using the theory of matching minors, but it can also be stated exclusively in the language of digraphs.

\begin{definition}[Non-even digraphs]\label{def:nonevendigraphs}
Let $D$ be a digraph.
A $0,1$-edge weighting $\omega\colon\E{D}\rightarrow\{0,1\}$ is called \emph{even} if there exists a dicycle $C$ in $D$ such that $\sum_{e\in\E{C}}\omega(e)$ is even.
A digraph $D$ is called \emph{even} if every $0,1$-edge weighting of $D$ is even, otherwise it is called \emph{non-even}.
\end{definition}

\begin{definition}[Butterfly minor]\label{def:butterflyminor}
Let $D$ be a digraph.
An edge $(u,v)\in\E{D}$ is \emph{contractible} if it is the only outgoing edge of $u$ or the only incoming edge at $v$.
The operation of contracting a contractible edge, that is identifying both of its endpoints and removing all resulting parallel edges and loops, is called a \emph{butterfly contraction}.
A digraph $H$ is a \emph{butterfly minor} of $D$ if it can be obtained from $D$ by a sequence of vertex deletions, edge deletions, and butterfly contractions.
\end{definition}

A \emph{butterfly minor model} of a digraph $H$ in a digraph $D$ is a mapping $\mu\colon \V{H}\cup\E{H}\rightarrow 2^D$, where $2^D$ denotes the collection of all subgraphs of $D$, such that
\begin{enumerate}
    \item for all $v\in \V{H}$, $\mu(v)$ consists of an in-arborescence
      $T^-_v$ with sink $r_v$ and an out-arborescence $T^+_v$ with
      source $r_v$ such that $T^-_v \cap T^+_v = r_v$,
    \item if $u, v\in\V{H}$ and $u \neq v$ then $\mu(u)$ and $\mu(v)$ are disjoint,
    \item for each $e=(u,v)\in\E{H}$, $\mu(e)$ is a directed path with tail in $T^+_u$ and head in $T^-_v$ which is internally disjoint from $\bigcup_{v\in V(H)}\V{\mu(v)}$, and
    \item if $e, e'\in\E{H}$ and $e \ne e'$ then $\mu(e)$ and $\mu(e')$ are internally disjoint.
\end{enumerate}
Note that a model $\mu$ of $H$ in a digraph $D$ describes a subgraph $H'$ of $D$ that can be transformed into $H$ be repeated applications of butterfly contractions.
In the following we will sometimes use the term `butterfly minor model' to refer to this subgraph instead of the mapping $\mu$.

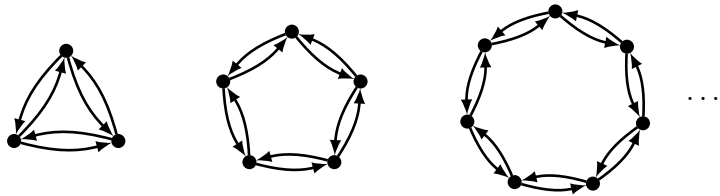
\begin{figure}[h!]
	\begin{center}
				\begin{tikzpicture}
		
		\node (D1) [v:ghost] {};
		\node (D2) [v:ghost,position=0:30mm from D1] {};
		\node (D3) [v:ghost,position=0:35mm from D2] {};
		\node (D4) [v:ghost,position=0:20mm from D3] {$\dots$};

		\node (C1) [v:ghost,position=0:0mm from D1] {
			
			\begin{tikzpicture}[scale=0.8]
			
			\pgfdeclarelayer{background}
			\pgfdeclarelayer{foreground}
			
			\pgfsetlayers{background,main,foreground}
			
			\begin{pgfonlayer}{main}
			
			\node (C) [] {};

			

			
			
			\node (v1) [v:main,position=90:10mm from C] {};
			\node (v2) [v:main,position=210:10mm from C] {};
			\node (v3) [v:main,position=330:10mm from C] {};
			
			
			
			

			

			
			
			\draw (v1) [e:main,->,bend right=15] to (v2);
			\draw (v2) [e:main,->,bend right=15] to (v3);
			\draw (v3) [e:main,->,bend right=15] to (v1);
			
			\draw (v1) [e:main,->,bend right=15] to (v3);
			\draw (v2) [e:main,->,bend right=15] to (v1);
			\draw (v3) [e:main,->,bend right=15] to (v2);

			

			
			
			\end{pgfonlayer}
			

			\begin{pgfonlayer}{background}
			
			\end{pgfonlayer}	
			
			\begin{pgfonlayer}{foreground}

			\end{pgfonlayer}
			\end{tikzpicture}
			
		};
		
		\node (C2) [v:ghost,position=0:0mm from D2] {
			
			\begin{tikzpicture}[scale=0.8]
			
			\pgfdeclarelayer{background}
			\pgfdeclarelayer{foreground}
			
			\pgfsetlayers{background,main,foreground}
			
			\begin{pgfonlayer}{main}
			
			\node (C) [] {};

			

			
			
			\node (v1) [v:main,position=90:12mm from C] {};
			\node (v2) [v:main,position=162:12mm from C] {};
			\node (v3) [v:main,position=234:12mm from C] {};
			\node (v4) [v:main,position=306:12mm from C] {};
			\node (v5) [v:main,position=18:12mm from C] {};
			
			
			
			

			

			
			
			\draw (v1) [e:main,->,bend right=15] to (v2);
			\draw (v2) [e:main,->,bend right=15] to (v3);
			\draw (v3) [e:main,->,bend right=15] to (v4);
			\draw (v4) [e:main,->,bend right=15] to (v5);
			\draw (v5) [e:main,->,bend right=15] to (v1);
			
			\draw (v1) [e:main,->,bend right=15] to (v5);
			\draw (v2) [e:main,->,bend right=15] to (v1);
			\draw (v3) [e:main,->,bend right=15] to (v2);
			\draw (v4) [e:main,->,bend right=15] to (v3);
			\draw (v5) [e:main,->,bend right=15] to (v4);
			
			

			
			
			\end{pgfonlayer}
			

			\begin{pgfonlayer}{background}
			
			\end{pgfonlayer}	
			
			\begin{pgfonlayer}{foreground}

			\end{pgfonlayer}
			\end{tikzpicture}
			
		};
		
		\node (C3) [v:ghost,position=0:0mm from D3] {
			
			\begin{tikzpicture}[scale=0.8]
			
			\pgfdeclarelayer{background}
			\pgfdeclarelayer{foreground}
			
			\pgfsetlayers{background,main,foreground}
			
			\begin{pgfonlayer}{main}
			
			\node (C) [] {};

			

			
			
			\node (v1) [v:main,position=90:15mm from C] {};
			\node (v2) [v:main,position=141.4:15mm from C] {};
			\node (v3) [v:main,position=192.8:15mm from C] {};
			\node (v4) [v:main,position=243.2:15mm from C] {};
			\node (v5) [v:main,position=294.6:15mm from C] {};
			\node (v6) [v:main,position=346.1:15mm from C] {};
			\node (v7) [v:main,position=37.5:15mm from C] {};
			
			
			
			

			

			
			
			\draw (v1) [e:main,->,bend right=15] to (v2);
			\draw (v2) [e:main,->,bend right=15] to (v3);
			\draw (v3) [e:main,->,bend right=15] to (v4);
			\draw (v4) [e:main,->,bend right=15] to (v5);
			\draw (v5) [e:main,->,bend right=15] to (v6);
			\draw (v6) [e:main,->,bend right=15] to (v7);
			\draw (v7) [e:main,->,bend right=15] to (v1);
			
			\draw (v2) [e:main,->,bend right=15] to (v1);
			\draw (v3) [e:main,->,bend right=15] to (v2);
			\draw (v4) [e:main,->,bend right=15] to (v3);
			\draw (v5) [e:main,->,bend right=15] to (v4);
			\draw (v6) [e:main,->,bend right=15] to (v5);
			\draw (v7) [e:main,->,bend right=15] to (v6);
			\draw (v1) [e:main,->,bend right=15] to (v7);

			

			
			
			\end{pgfonlayer}
			

			\begin{pgfonlayer}{background}
			
			\end{pgfonlayer}	
			
			\begin{pgfonlayer}{foreground}

			\end{pgfonlayer}
			\end{tikzpicture}
			
		};
		
	\end{tikzpicture}
    \end{center}
\caption{The family of odd bicycles.}
\label{fig:oddbicycles}
\end{figure}

Let $G$ be an undirected graph.
We denote by $\Bidirected{G}$ the \emph{bidirection} of $G$, that is the digraph obtained from $G$ by replacing every edge $uv\in G$ with the two directed edges $(u,v)$ and $(v,u)$.
In case $C$ is an undirected cycle we call $\Bidirected{C}$ a \emph{bicycle}.
Moreover, in case $C$ is of odd length we say that $\Bidirected{C}$ is an \emph{odd bicycle} (see \cref{fig:oddbicycles}), otherwise $\Bidirected{C}$ is said to be an \emph{even bicycle}.

We are now ready to state a first key result towards the resolution of the Even Dicycle Problem which was also the first appearance of the concept of butterfly minors in the literature.

\begin{theorem}[Seymour and Thomassen \cite{SeymourThomassen1987EvenDirectedGraphs}]\label{thm:nonevendigraphs}
A digraph $D$ is non-even if and only if it does not contain an odd bicycle as a butterfly minor.
\end{theorem}

Following \cite{SeymourThomassen1987EvenDirectedGraphs}, we will also call the butterfly minor of an odd bicycle a \emph{weak odd bicycle}.
The existence of \cref{thm:nonevendigraphs} has several interesting implications.
The first is that the property of being non-even is closed under butterfly minors and the second is that every weak odd bicycle has a dicycle with even total weight under every $0,1$-edge weighting.
Hence we may establish the following fundamental observation.

\begin{observation}\label{obs:oddbicycleevendicycle}
Any weak odd bicycle contains an even dicycle.
\end{observation}

Let $G$ be an undirected graph.
A \emph{block} of $G$ is a maximal $2$-connected subgraph of $G$.
In undirected graphs it is straight forward to see that a $2$-connected minor must be contained in one of its blocks.
For digraphs it is not even clear how to define a block in the first place since there exists a discrepancy between the existence of a cut-vertex that just destroys strong connectivity and a cut-vertex which acts as an undirected separator.
The following is a possible extension of the idea of blocks and the block decomposition to digraphs which leans heavily into the matching theoretic roots of the Even Dicycle Problem.

\begin{definition}[(Directed) separation]\label{def:separation}
	Let $D$ be a digraph.
	A \emph{directed separation} is a tuple $(X,Y)$ such that $X\cup
        Y=\V{D}$ and there is no directed edge with tail in $Y \setminus X$
        and head in $X \setminus Y$.
	The set $X \cap Y$ is called the \emph{separator} and the
        \emph{order} of $(X,Y)$ is $\Abs{X\cap Y}$. 
	A tuple $(X,Y)$ is called a \emph{separation}, or \emph{undirected separation}, if we want to emphasise this fact, if $(X,Y)$ and $(Y,X)$ both are directed separations.
	If $G$ is an undirected graph, then $(X,Y)$ is a separation if
        it is a separation in $\Bidirected{G}$.

        A (directed) separation $(X, Y)$ is \emph{trivial} if $X \subseteq Y$
        or $Y \subseteq X$. Otherwise it is \emph{non-trivial}.
\end{definition}

\begin{definition}[Directed tight cut contraction]\label{def:directedtightcut}
  Let $D$ be a digraph, and $\Brace{X,Y}$ be a non-trivial
  directed separation of order $1$ in $D$. 
  Let $\Set{v}=X\cap Y$, we set
  \begin{align*}
    \ContractsTo{D}{X}{v_X}\coloneqq
    D-X+v_X+&\CondSet{\Brace{y,v_X}}{\Brace{y,v}\in\E{D}\text{ and }y\in
              Y}\\
    +&\CondSet{\Brace{v_X,y}}{\Brace{x,y}\in\E{D},~x\in X\text{, and }y\in Y}\text{, and}\\
    \ContractsTo{D}{Y}{v_Y}\coloneqq
    D-Y+v_Y+&\CondSet{\Brace{v_Y,x}}{\Brace{v,x}\in\E{D}\text{ and }x\in
              X}\\
    +&\CondSet{\Brace{x,v_Y}}{\Brace{x,y}\in\E{D},~x\in X\text{, and }y\in Y}.
  \end{align*}
  The two digraphs $\ContractsTo{D}{X}{v_X}$ and
  $\ContractsTo{D}{Y}{v_Y}$ are called the
  \emph{$\Brace{X,Y}$-contractions} of $D$. 
\end{definition}
It is
straightforward to observe that any directed tight cut contraction of a strongly connected digraph $D$ is in fact a butterfly minor of $D$.
Moreover, the butterfly contraction of an edge can be seen as a special case of a directed tight cut contraction\footnote{Tight cuts in undirected graphs will be discussed in \Cref{sec:shifting}.}.

We may select any non-trivial directed separation of order one in a given digraph $D$ and form its two contractions.
Then, in case any of the two contractions still has a non-trivial directed separation of order one we can continue with this process.
Eventually this yields a list of strongly $2$-connected digraphs uniquely determined\footnote{This follows from the famous result of Lovasz \cite{lovasz1987matching} on the uniqueness of the tight cut decomposition.} by $D$.
We call a maximal set of pairwise laminar non-trivial directed separations of order one together with the corresponding collection of dibraces a \emph{directed tight cut decomposition} of $D$.
These strongly $2$-connected digraphs are called the \emph{dibraces}\footnote{Dibraces are specialisations of a more complex matching theoretic concept known as `braces' which we will describe in \cref{sec:shifting}.} of $D$.

Besides the directed tight cut decomposition, the structure theorem for non-even digraphs also has a topological component.

\begin{definition}(Drawing in a surface)\label{def:drawing}
A \emph{drawing} (with crossings) \emph{in a surface $\Sigma$} is a triple $(U,V,E)$ such that
\begin{itemize}
	\item $V$ and $E$ are finite,
	\item $V\subseteq U\subseteq\Sigma$,
	\item $V\cup\bigcup_{e\in E}e=U$ and $V\cap(\bigcup_{e\in E}e)=\emptyset$,
	\item for every $e\in E$, either $e=h((0,1))$, where $h\colon[0,1]_{\R}\rightarrow U$ is a homeomorphism onto its image with $h(0),h(1)\in V$, or $e=h(\mathbb{S}^2-(1,0))$, where $h\colon\mathbb{S}^2\rightarrow U$ is a homeomorphism onto its image with $h(0,1)\in V$, and
	\item if $e,e'\in E$ are distinct, then $\Abs{e\cap e'}$ is finite.
\end{itemize}
If for any two distinct edges $e,e'\in E$, we have that $e$ and $e'$ intersect at most in their endpoints, we say that $(U,V,E)$ is \emph{cross-free}.
We call a drawing $(U,V,E)$ in $\Sigma$ a \emph{drawing of $D$ (or $G$)} if $V$ is the vertex set of $D$ (or $G$) and $E$ are the edges of $D$ (or $G$).
Given a drawing $\Gamma$ of a digraph, we write $\V{\Gamma}$ to reference its vertex set and $\E{\Gamma}$ to reference the edge set.
\end{definition}
By a \emph{surface} we mean a compact $2$-dimensional manifold with or without boundary.
By the classification theorem of surfaces every surface is homeomorphic to the sphere with $h$ handles and $c$ cross-caps added, and the interior of $d$ disjoint open disks removed, in which case the \emph{Euler genus} of the surface is defined to be $2h+c$.
A cross-free drawing on the sphere, possibily with some open holes, is called \emph{planar}.

Let $\Gamma$ be a drawing of the digraph $D$ in a surface $\Sigma$ and $v\in\V{G}$ be a vertex.
We say that a disk $\Delta\subseteq\Sigma$ is \emph{centred at $v$} if $v$ is drawn in its interior by $\Gamma$ and $\Delta$ contains no other vertices in its interior.
Let  $\Delta \subseteq \Sigma$ be a disk centred at $v$ such that every edge of $D$
incident with $v$ contains exactly one point from the boundary
$\Boundary{\Delta}$. 
Let $F\subseteq\E{D}$ be the edges incident with $v$ and $\Set{F_1,F_2}$ be a
bipartition of $F$. 
For each $f\in F$ let $p_f\in \Boundary{\Delta}$ be the point that $f$ has on
the boundary of $\Delta$. 
We say that $\Brace{F_1,F_2}$ is a \emph{butterfly in $\Delta$} if there
exists a curve $\gamma$ through $v$ in $\Sigma$ with both endpoints, $x$ and
$y$, on $\Delta$ such that we can number the two internally disjoint curves
$\delta_1\subseteq\Boundary{\Delta}$ and $\delta_2\subseteq\Boundary{\Delta}$ with endpoints $x$ and $y$
to obtain $\CondSet{p_f}{f\in F_i}\subseteq \delta_i$ for both $i\in[2]$. 

\begin{definition}[Strong drawing]
  Let $D$ be a digraph and $\Sigma$ be a surface.
  A cross-free drawing $\Gamma$ of $D$ in $\Sigma$ is \emph{strong} if for
  every vertex $v\in\V{D}$ there exists a non-trivial disk $\Delta\subseteq\Sigma$ centred
  at $v$ such that 
  \[
    \Brace{\CondSet{\Brace{u,v}}{\Brace{u,v}\in\E{D}},\CondSet{\Brace{v,u}}{\Brace{v,u}\in\E{D}}}
  \]
  is a butterfly in $\Delta$.

The smallest integer $g\in\N$ such that $D$ can be strongly drawn in a surface of Euler genus $g$ is called the \emph{strong genus} of $D$.
We denote the strong genus of $D$ by $\StrongGenus{D}$.
If $\StrongGenus{D}=0$, $D$ is said to be \emph{strongly planar}.
\end{definition}
Note that the strong genus of a digraph $D$ is closed under vertex and edge deletion, and under butterfly contractions.

The final ingredient we need is a directed analogue of small clique sums.
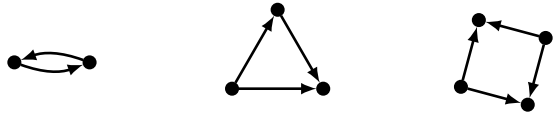
\begin{figure}[!ht]
	\centering
	    \begin{tikzpicture}
		\pgfdeclarelayer{background}
		\pgfdeclarelayer{foreground}
		\pgfsetlayers{background,main,foreground}
		
		\node (mid) [v:ghost] {};
		\node (left) [v:ghost,position=180:30mm from mid] {};
		\node (right) [v:ghost,position=0:30mm from mid] {};

		\node(u1) [v:main,position=180:5mm from left] {};
		\node(u2) [v:main,position=0:5mm from left] {};
		
		\draw [e:main,->,bend right=20] (u1) to (u2);
		\draw [e:main,->,bend right=20] (u2) to (u1);
		
		\node (v1) [v:main,position=90:7mm from mid] {};
		\node (v2) [v:main,position=210:7mm from mid] {};
		\node (v3) [v:main,position=330:7mm from mid] {};

		\draw [e:main,<-] (v1) to (v2);
		\draw [e:main,->] (v2) to (v3);
		\draw [e:main,<-] (v3) to (v1);
		
		\node (w1) [v:main,position=30:6.5mm from right] {};
		\node (w2) [v:main,position=120:6.5mm from right] {};
		\node (w3) [v:main,position=210:6.5mm from right] {};
		\node (w4) [v:main,position=300:6.5mm from right] {};
		
		\draw [e:main,->] (w1) to (w2);
		\draw [e:main,->] (w1) to (w4);
		\draw [e:main,->] (w3) to (w2);
		\draw [e:main,->] (w3) to (w4);
		
	\end{tikzpicture}
	\caption{The subgraphs necessary for the small-cycle-sum operation.}
	\label{fig:smallcyclesum}
\end{figure}

\begin{definition}[Small-cycle-sum]\label{def:smallcyclesum}
	Let $D_0$ be a digraph, let $u,v\in\V{D_0}$, and let $(u,v),(v,u)\in\E{D_0}$.
	Let $D_1$ and $D_2$ be such that $D_1\cup D_2=D_0$,
        $\V{D_1}\cap\V{D_2}=\Set{u,v}$, $\V{D_1}\setminus\V{D_2}\neq\emptyset$,
        $\V{D_2}\setminus\V{D_1}\neq\emptyset$, and $\E{D_1}\cap\E{D_2} = \{ (u,v), (v, u)\}$.
	Let $D$ be obtained from $D_0$ by deleting some (possibly neither) of the edges $(u,v)$, $(v,u)$.
	We say that $D$ is a \emph{$2$-sum} of $D_1$ and $D_2$.
	
	Let $D_0$ be a digraph, let $u,v,w\in\V{D_0}$ and $(w,u), (u,v),
        (w,v) \in
        \E{D_0}$, and assume that $D_0$ has a directed cycle containing the edge $(w,v)$, but not the vertex $u$.
	Let $D_1$ and $D_2'$ be such that $D_1\cup D_2' = D_0$, $\V{D_1}\cap\V{D_2'}=\Set{u,v,w}$, $\V{D_1}\setminus\V{D_2'}\neq\emptyset$, $\V{D_2'}\setminus\V{D_1}\neq\emptyset$, and $\E{D_1}\cap\E{D_2'}=\Set{(w,u),(u,v),(w,v)}$.
	Let $D_2'$ have no edge with tail $v$, and no edge with head $w$.  Let $D_2$ be the digraph obtained from $D_2'$ by contracting $(w,v)$. (Note that $(w,v)$ is not butterfly contractible in $D_2'$.)
	Let $D$ be obtained from $D_0$ by deleting some (possibly none) of the edges $(w,u), (u,v),(w,v)$.
	We say that $D$ is a \emph{$3$-sum} of $
	D_1$ and $D_2$.
	
	Let $D_0$ be a digraph, let $x,y,u,v\in\V{D_0}$ as well as $(x,y),(x,v),(u,y),(u,v) \in \E{D_0}$, and assume that $D_0$ has a directed cycle containing precisely two of the edges $(x,y),(x,v),(u,y),(u,v)$.
	Let $D_1$ and $D_2'$ be such that $D_1\cup D_2'=D_0$, $\V{D_1}\cap\V{D_2'}=\Set{x,y,u,v}$, $\V{D_1}\setminus\V{D_2'}\neq\emptyset$, $\V{D_2'}\setminus\V{D_1}\neq\emptyset$, and $\E{D_1}\cap\E{D_2'}=\Set{(x,y),(x,v),(u,y),(u,v)}$.
	Let $D_2'$ have no edge with tail $y$ or $v$, and no edge with head $x$ or $u$.
	Let $D_2$ be the digraph obtained from $D_2'$ by contracting the edges $(x,y)$ and $(u,v)$.
    (Again, note that the edges $(x,y)$ and $(u,v)$ are not butterfly contractible.)
	Finally, let $D$ be obtained from $D_0$ by deleting some (possibly none) of the edges $(x,y),(x,v),(u,y),(u,v)$.
	We say that $D$ is a \emph{$4$-sum} of $D_1$ and $D_2$.
	
	We say that a digraph $D$ is a \emph{small cycle sum} of two digraphs $D_1$ and $D_2$ if it is an $i$-sum of $D_1$ and $D_2$ for some $i\in[2,4]$.
\end{definition}

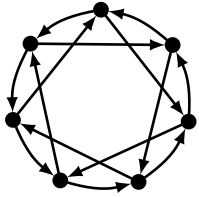
\begin{figure}[h!]
	\begin{center}
		        \begin{tikzpicture}[scale=0.8]
		
		\pgfdeclarelayer{background}
		\pgfdeclarelayer{foreground}
		
		\pgfsetlayers{background,main,foreground}
		
		\tikzstyle{v:main} = [draw, circle, scale=0.5, thick,fill=black]
		\tikzstyle{v:tree} = [draw, circle, scale=0.3, thick,fill=black]
		\tikzstyle{v:border} = [draw, circle, scale=0.75, thick,minimum size=10.5mm]
		\tikzstyle{v:mainfull} = [draw, circle, scale=1, thick,fill]
		\tikzstyle{v:ghost} = [inner sep=0pt,scale=1]
		\tikzstyle{v:marked} = [circle, scale=1.2, fill=CornflowerBlue,opacity=0.3]
		
		\tikzset{>=latex} 
		\tikzstyle{e:marker} = [line width=9pt,line cap=round,opacity=0.2,color=DarkGoldenrod]
		\tikzstyle{e:colored} = [line width=1.2pt,color=BostonUniversityRed,cap=round,opacity=0.8]
		\tikzstyle{e:coloredthin} = [line width=1.1pt,opacity=0.8]
		\tikzstyle{e:coloredborder} = [line width=2pt]
		\tikzstyle{e:main} = [line width=1pt]
		\tikzstyle{e:extra} = [line width=1.3pt,color=LavenderGray]
		
		\begin{pgfonlayer}{main}
		
		\node (C) [] {};
		
		\node (C1) [v:ghost, position=180:25mm from C] {};
		
		\node (C2) [v:ghost, position=0:0mm from C] {};
		
		\node (C3) [v:ghost, position=0:25mm from C] {};

		

		
		
		\node (v1) [v:main,position=90:15mm from C2] {};
		\node (v2) [v:main,position=141.4:15mm from C2] {};
		\node (v3) [v:main,position=192.8:15mm from C2] {};
		\node (v4) [v:main,position=243.2:15mm from C2] {};
		\node (v5) [v:main,position=294.6:15mm from C2] {};
		\node (v6) [v:main,position=346.1:15mm from C2] {};
		\node (v7) [v:main,position=37.5:15mm from C2] {};
		
		
		
		

		

		
		
		\draw (v1) [e:main,->,bend right=15] to (v2);
		\draw (v2) [e:main,->,bend right=15] to (v3);
		\draw (v3) [e:main,->,bend right=15] to (v4);
		\draw (v4) [e:main,->,bend right=15] to (v5);
		\draw (v5) [e:main,->,bend right=15] to (v6);
		\draw (v6) [e:main,->,bend right=15] to (v7);
		\draw (v7) [e:main,->,bend right=15] to (v1);
		
		\draw (v1) [e:main,->] to (v6);
		\draw (v2) [e:main,->] to (v7);
		\draw (v3) [e:main,->] to (v1);
		\draw (v4) [e:main,->] to (v2);
		\draw (v5) [e:main,->] to (v3);
		\draw (v6) [e:main,->] to (v4);
		\draw (v7) [e:main,->] to (v5);
		
		

		
		
		\end{pgfonlayer}
		

		\begin{pgfonlayer}{background}
		
		\end{pgfonlayer}	
		
		\begin{pgfonlayer}{foreground}

		\end{pgfonlayer}
		\end{tikzpicture}
	\end{center}
	\caption{The digraph $F_7$.}
	\label{fig:F7}
\end{figure}

Using this tool, Robertson, Seymour, and Thomas were able to provide a structural description of non-even digraphs.
Whilst they chose to pursue this goal via the theory of digraphs, independently McCuaig worked on the matching theoretic equivalent of their result (see the Acknowledgement section of \cite{mccuaig1997permanents} for reference).

\begin{theorem}[Robertson, Seymour, and Thomas \cite{Robertson1999PermanentsPfaffianOrientations}, and McCuaig \cite{McCuaig2004Polya}]\label{thm:nonevenstructure}
A digraph $D$ is non-even if and only if all of its dibraces are either $F_7$ (see \cref{fig:F7}), or can be obtained from strongly planar dibraces by means of small cycle sums.
\end{theorem}

This structural description of odd bicycle free digraphs immediately implies a polynomial time algorithm.

\begin{theorem}[Robertson, Seymour, and Thomas \cite{Robertson1999PermanentsPfaffianOrientations}, and McCuaig \cite{McCuaig2004Polya}]\label{thm:nonevenalgo}
There exists an algorithm that decides in polynomial time if a given digraph $D$ is non-even.
\end{theorem}

Finally, as a corollary, this implies a polynomial time algorithm to solve the Even Dicycle Problem.

\begin{corollary}[Seymour and Thomassen \cite{SeymourThomassen1987EvenDirectedGraphs}, Robertson, Seymour, and Thomas \cite{Robertson1999PermanentsPfaffianOrientations}, and McCuaig \cite{McCuaig2004Polya}]\label{cor:evendicyclealgo}
There exists an algorithm that decides in polynomial time if a given digraph $D$ contains an even dicycle.
\end{corollary}

\subsection{Societies, odd decompositions, and maelstroms}\label{subsec:drawing}

Large portions of our proof are inspired by the work of Kawarabayashi, Wollan, and Thomas on simplified methods to prove the global structure theorem of $H$-minor free graphs \cite{kawarabayashi2020quickly}.
The main concept upon which they build their proof is the idea of societies.

\begin{definition}[Society]\label{def:society}
Let $\Omega$ be a cyclic permutation of the elements of some set and let $\V{\Omega}$ denote this set.
A \emph{society} is a pair $(G,\Omega)$, where $G$ is a (di)graph, and $\Omega$ is a cyclic permutation with $\V{\Omega}\subseteq\V{G}$.

A \emph{cylindrical society} is a tuple $(G,\Omega_1,\Omega_2)$, where $G$ is a (di)graph and $\Omega_1$, $\Omega_2$ are cyclic permutations with $\V{\Omega_i}\subseteq\V{G}$ for both $i\in[2]$ and $\V{\Omega_1}\cap\V{\Omega_2}=\emptyset$.
\end{definition}

\begin{definition}[Plane decomposition]\label{def:planedecomposition}
A \emph{plane decomposition} of a digraph $D$ is a tuple $\delta=(\Gamma,\mathcal{V},\mathcal{D})$, where
\begin{itemize}
    \item $\Gamma$ is a drawing of $D$ in a sphere $\Sigma$ possibly with crossings and some open holes,
    
    \item $\mathcal{V}$ and $\mathcal{D}$ are collections of closed disks, each a subset of $\Sigma$, we call the disks in $\mathcal{V}$ \emph{big vertices}, and
    
    \item $\Gamma$, $\mathcal{V}$, and $\mathcal{D}$ satisfy axioms \descriptionfont{PD1} to \descriptionfont{PD7}.
\end{itemize}
Given a big vertex $\mathsf{v} \in \mathcal{V}$, the edges of $D$ drawn in $\Gamma$ such that they have at least one point on either side of $\Boundary{\mathsf{v}}$ are called the \emph{crossing edges of $\mathsf{v}$}.
We say that the boundary of a disk $\Delta$ \emph{intersects $\Gamma$ in a big vertex} $\mathsf{v}\in\mathcal{V}$ if $\Boundary{\Delta}$ meets $\Boundary{\mathsf{v}}$ in exactly two points, $\Delta\cap\mathsf{v}$ does not contain a vertex drawn in the interior of $\Delta$ or $\mathsf{v}$, and $\Delta$ intersects the edges of $\Gamma$ which have at least one point drawn in the interior of $\mathsf{v}$ in crossing edges only.
\begin{description}
    \item[PD1] The disks in $\mathcal{V}$ (respectively $\mathcal{D}$) have pairwise disjoint interiors.
    
    \item[PD2] The boundary of a disk $d$ in $\mathcal{D}$ intersects $\Gamma$ only in vertices, big vertices, and crossing edges of big vertices it intersects and $d$ must contain a point drawn by $\Gamma$.
    
    \item[PD3] For every $\mathsf{v}\in\mathcal{V}$ there exists a unique vertex $v_{\mathsf{v}}$, called the \emph{door of $\mathsf{v}$}, such that $v_{\mathsf{v}}$ is the only vertex drawn on the boundary of $\mathsf{v}$, $(X\cup\Set{v_{\mathsf{v}}},Y\cup\Set{v_{\mathsf{v}}})$ or $(Y\cup\Set{v_{\mathsf{v}}},X\cup\Set{v_{\mathsf{v}}})$ is a directed separation in $D$ of order one where $X$ is the collection of all vertices drawn in the interior of $\mathsf{v}$, and none of the edges crossing $\Boundary{\mathsf{v}}$ are incident with $v_{\mathsf{v}}$.

    \item[PD4] If $\Delta_1,\Delta_2\in\mathcal{D}$ are distinct, then $\Delta_1\cap \Delta_2\subseteq\V{\Gamma}$.
    
    \item[PD5] If $\mathsf{u},\mathsf{v}\in\mathcal{V}$ are distinct, let us denote by $E(\mathsf{u},\mathsf{v})$ the set of edges which are crossing for both $\mathsf{u}$ and $\mathsf{v}$.
    Then there exists a disk $\Delta\in\mathcal{D}$ intersecting $\mathsf{u}$ and $\mathsf{v}$ but neither intersecting nor containing another big vertex, nor another vertex of $\Gamma$, such that $\Delta$ is internally disjoint from $\mathsf{u}$ and $\mathsf{v}$.
    Moreover, $\Delta$ contains the drawings of all edges from $E(\mathsf{u},\mathsf{v})$ in $\Gamma$ minus the interiors of $\mathsf{u}$ and $\mathsf{v}$.
    
    \item[PD6] If $\mathsf{v}\in\mathcal{V}$ is a big vertex and $u$ a vertex of $\Gamma$ drawn on the boundary of some disk in $\mathcal{D}$ such that the set $E(\mathsf{v},u)$ of crossing edges of $\mathsf{v}$ which are incident with $u$ is non-empty, then there exists a disk $\Delta\in\mathcal{D}$ such that $\Delta$ intersects $\mathsf{v}$ and $u$ is drawn on the boundary of $\Delta$, but it neither intersects nor contains another big vertex, nor another vertex of $\Gamma$, such that $\Delta$ is internally disjoint from $\mathsf{v}$.
    Moreover, $\Delta$ contains the part of the drawings of all edges from $E(\mathsf{v},u)$ which are drawn by $\Gamma$ in the plane minus the interior of $\mathsf{v}$.
    
    \item[PD7] Every edge $e\in \E{\Gamma}$ either belongs to the interior of one of the disks in $\mathcal{D}\cup\mathcal{V}$ or there exist $\mathsf{v}\in\mathcal{V}$ with either $e\in E(\mathsf{u},\mathsf{v})$ or $e\in E(u,\mathsf{v})$ for some $\mathsf{u}\in\mathcal{V}$ or vertex $u$ of $\Gamma$ and $e \cap \Boundary{\mathsf{v}}$ is a single point.
\end{description}
We define containment and isomorphism of plane decompositions in the natural way.
Restrictions can similarly defined naturally, though a small complication may occur.
We highlight this problem in \Cref{rem:restriction} and explain why it does not impede our efforts.
Let $N$ be the set of all big vertices and all vertices of $\Gamma$ that do not belong to the interior of any of the disks in $\mathcal{D}\cup\mathcal{V}$ and are not the door of some big vertex.
We will refer to the elements of $N$ as the \emph{nodes} of $\delta$.
Note that any big vertex can be identified into a single vertex using a directed tight cut contraction, which justifies us treating the elements of $\mathcal{V}$ as vertices for the purpose of a plane decomposition.
If $\Delta\in\mathcal{D}$, we refer to the set $\Delta-N$ as a \emph{cell} of $\delta$.
We denote the set of nodes of $\delta$ by $N(\delta)$ and the set of cells by $C(\delta)$.
For a cell $c\in C(\delta)$, we denote by $\widetilde{c}$ the set of vertices of $\Gamma$ that lie on the boundary of the closure of $c$, together with all big vertices that are intersected by the closure of $c$.
Thus the cells $c$ of $\delta$ with $\widetilde{c}\neq\emptyset$ form the hyperedges of a hypergraph with vertex set $N(\delta)$, where $\widetilde{c}$ is the set of vertices incident with $c\in C(\delta)$.
For a cell $c\in C(\delta)$ (or any other closed disk $c$), we define $\sigma_{\delta}(c)$ (or $\sigma(c)$ when $\delta$ is understood from the context) to be the subgraph of $D$ consisting of all vertices and edges drawn in the union of the closure of $c$ and all big vertices it intersects.
We define $\pi_{\delta}\colon N(\delta)\rightarrow\V{D}\cup (\V{D}\times 2^{V(D)})$ to be the mapping that assigns to every vertex of $\Gamma$ in $N(\delta)$ its identity and to every big vertex $\mathsf{v}\in\mathcal{V}\cap N(\delta)$ a tuple $(v,X)$, where $v$ is the door of $\mathsf{v}$ and $X$ is the set of vertices of $D$ corresponding to the vertices contained in the closure of $\mathsf{v}$.

Let $x\in N(\delta)$.
If $x$ is a big vertex with the door $v$, we say that an edge is \emph{incident with $x$} if it is a crossing edge of $x$ or it is an edge incident with $v$ that is not drawn in the interior of $x$.
In this case the \emph{outgoing edges} of $x$ are those incident with $x$ whose tail is either $v$ or lies in the interior of $x$, while its \emph{incoming edges} are those incident with $x$ whose head is either $v$ or lies in the interior of $x$.
We say that a disk $\Delta$ is \emph{centred at $x$} if $x$ is a vertex of $\Gamma$ and $\Delta$ is centred at $x$, or $x$ is a big vertex and $\Delta$ contains $x$ in its interior, but is otherwise disjoint from the vertices of $\Gamma$.
A plane decomposition is \emph{strong} if for every $x\in N(\delta)$ there exists a disk $\Delta$ centred at $x$ such that we can partition the incoming edges at $x$ and the outgoing edges of $x$ into the sets $F_1$ and $F_2$ in such a way that $(F_1,F_2)$ forms a butterfly in $\Delta$.
\end{definition}

\begin{definition}[Cross]\label{def:cross}
    Let $D$ be a digraph, let $\delta=(\Gamma,\mathcal{V},\mathcal{D})$ be a strong plane decomposition of $D$, let $C$ be a (not necessarily directed) cycle $C$ in $D$, and let $P_1, P_2$ be two disjoint paths in $D$ starting and ending on $C$.
    For $i \in [2]$, let $u_i$ be the tail of $P_i$ and let $v_i$ be the head of $P_i$.
    
    We say that $P_1$ and $P_2$ form a \emph{cross over $C$} if $u_1,u_2,v_1,v_2$ appear in the given order on $C$ and within $\Gamma$ the interiors of both paths $P_1, P_2$ are drawn entirely within one of the two disks bounded by the trace of $C$.
    If $S = (D,\Omega)$ is a society, then we say that $P_1$ and $P_2$ form \emph{cross over $S$} if $u_1,u_2,v_1,v_2$ appear in the given order in $V(\Omega)$, with respect to the order provided by $\Omega$.
\end{definition}

\begin{definition}[Maelstrom]\label{def:maelstrom}
    Let $D$ be a digraph and let $\delta=(\Gamma,\mathcal{V},\mathcal{D})$ be a strong plane decomposition of $D$.
    Let $c\in C(\delta)$ be a cell of $\delta$.
    In case $\Abs{\widetilde{c}} = 2$ and there is no vertex of $\Gamma$ drawn inside of $c$, we say that $c$ is a \emph{common cell}.
    Suppose $\Abs{\widetilde{c}}\leq 4$.
    For every big vertex $\mathsf{v}\in\widetilde{c}$, let $\pi(\mathsf{v})=(v_{\mathsf{v}},X_{\mathsf{v}})$.
    Let $D'$ be the graph obtained from $D$ by contracting for every big vertex $\mathsf{v}\in\widetilde{c}$ the set $X_{\mathsf{v}}$ into the vertex $v_\mathsf{v}$ and removing all resulting parallel edges and loops.
    Let $Z\coloneqq \CondSet{\pi(v)}{v\in\widetilde{c}\setminus\mathcal{V}}\cup \CondSet{v_{\mathsf{v}}}{\mathsf{v}\in\widetilde{c}\cap\mathcal{V}}$.
    If there exist digraphs $D_1$ and $D_2$ such that $D_2-Z=\sigma(c)-Z-\bigcup_{\mathsf{v}\in\widetilde{c}\cap\mathcal{V}}X_{\mathsf{v}}$, $D_1-Z=D-\V{\sigma(c)}$, $D_2$ is non-even, and $D'$ is a small cycle sum of $D_1$ and $D_2$, then we call $c$ a \emph{conjunction cell}.
    
    A cell $c\in C(\delta)$ that is neither a common cell nor a conjunction cell is called a \emph{maelstrom}. 
\end{definition}

\begin{remark}\label{rem:restriction}
    We now note that the restriction of a plane decomposition may behave strangely when either deleting the door of a big vertex or a vertex on the boundary of a conjunction cell.
    In the first case, we cut off a part of the graph that may be non-planar and may no longer be decomposable according to our definition.
    However, in these cases, this part of the graph occupies a different strong component than the grounded part of the graph and we can simply remove this entire component from the decomposition.

    Similarly, removing a vertex on the boundary of a conjunction cell is also troublesome.
    Here we opt to simply remove all vertices on the boundary of a conjunction cell each time we have to remove one.
    This may blow up the size of the set of vertices hitting all even dicycles that we may find by at most a factor of 4.
    We can therefore ignore this problem for the sake of our arguments. 
\end{remark}

Notice, that a maelstrom as defined above can be seen as a directed form of a ``vortex'' as defined in \cite{kawarabayashi2020quickly}.
However, due to the richer structure of digraphs, we need to impose further structure onto our maelstroms.

\begin{definition}[Odd decomposition]\label{def:odddecomposition}
    A strong plane decomposition $\delta=(\Gamma,\mathcal{V},\mathcal{D})$ of a digraph $D$ is an \emph{odd decomposition} of $D$ if for every even dicycle $C$ in $D$ there exists either
    \begin{itemize}
        \item a big vertex $\mathsf{v}$ such that $C$ is drawn entirely in $\mathsf{v}$,
        \item a conjunction cell $c\in C(\delta)$ such that $C$ is completely drawn within $c$, or
        \item a maelstrom cell $c\in C(\delta)$ such that $C$ contains a vertex drawn within $c$.
    \end{itemize}
    An odd decomposition is \emph{pure} if no big vertex and no conjunction cell contains an even dicycle.
    
    An odd decomposition is \emph{maelstrom-free} if no cell in $C(\delta)$ is a maelstrom.
\end{definition}

It is straightforward to observe that the existence of a pure odd decomposition is necessary for the absence of even dicycles in a digraph.
For this recall \cref{obs:oddbicycleevendicycle} together with \cref{thm:nonevendigraphs} and \cref{thm:nonevenstructure}.

\begin{observation}\label{obs:puredecomposition}
A digraph $D$ has a maelstrom-free and pure odd decomposition if and only if it does not contain an even dicycle.
\end{observation}

In \cref{sec:shifting} we introduce several tools that will allow us to reduce the existence of even dicycles to areas bounded by a single closed curve, provided that sufficient infrastructure, in the form of concentric dicycles that are more or less aligned with the closed curve, exists.
This approach allows us to apply most of our arguments to a more local setting.
Towards a formalisation of this method, we consider odd decompositions with respect to societies and cylindrical societies.

\begin{definition}[Renditions]\label{def:rendition}
    Let $(D,\Omega)$ be a society and $\Delta$ be a closed disk, possibly with an open hole.
    A \emph{rendition in $\Delta$} of $D$ is a strong plane decomposition $\delta=(\Gamma,\mathcal{V},\mathcal{D})$ of $D$ such that 
    \begin{enumerate}
        \item $\V{\Gamma}\subseteq \Delta$,
        \item for all $\Delta'\in\mathcal{D}$ and all $\mathsf{v}\in\mathcal{V}$, we have $\Delta',\mathsf{v}\subseteq \Delta$,
        \item for each big vertex $\mathsf{v}\in\mathcal{V}$, $\Boundary{\Delta}$ meets $\mathsf{v}$ either in exactly a single segment, or is completely disjoint from $\mathsf{v}$,
        \item for each $\mathsf{v}\in\mathcal{V}$ which is intersected by $\Boundary{\Delta}$, exactly the image under $\pi_{\delta}$ of the door of $\mathsf{v}$ lies in $\V{\Omega}$ and no other vertex drawn in $\mathsf{v}$ lies in $V(\Omega)$,
        \item $\CondSet{\pi_{\delta}(v)}{v\in\Boundary{\Delta} \cap ( N(\delta) \setminus \mathcal{V} )} \subseteq \V{\Omega}$, and
        \item one of the cyclic orders of $\Boundary{\Delta}$ maps the images under $\pi_{\delta}$ of the vertices in $\Boundary{\Delta} \cap ( N(\delta) \setminus \mathcal{V} )$ and the doors of the big vertices intersected by $\Boundary{\Delta}$ to the order of $\Omega$.
    \end{enumerate}
    
    Now let $\Delta'$ be obtained from a closed disk $\Delta''$ by removing an open disk disjoint from the boundary of $\Delta''$ and let $(D,\Omega_1,\Omega_2)$ be a cylindrical society.
    Let $B_1$ and $B_2$ be the two closed curves in $\Delta'$ whose union equals $\Boundary{\Delta'}$.
    Observe that for each $i\in[2]$ the curve $B_i$ bounds a closed disk $\Delta_i'$ with an open hole that contains $\Delta'$.
    A \emph{rendition in $\Delta'$} is an odd decomposition $\delta=(\Gamma,\mathcal{V},\mathcal{D})$ such that for both $i\in[2]$, $\delta$ is a rendition of $(D,\Omega_i)$ in the disk $\Delta_i'$.
    We say that $(D,\Omega_1,\Omega_2)$ has a \emph{rendition in the disk} if there exists $\Delta'$ as above such that $(D,\Omega_1,\Omega_2)$ has a rendition in $\Delta'$.
    
    We call a rendition \emph{non-even} if it is maelstrom-free, and \emph{odd} if $\delta$ is a maelstrom-free, pure odd decomposition.
\end{definition}

\begin{definition}[Grounded subgraphs and traces]\label{def:trace}
Let $\delta=(\Gamma,\mathcal{V},\mathcal{D})$ be a maelstrom-free strong plane decomposition of $D$.

If $Q$ is a path in $D$ and $\mathsf{v}\in N(\delta)$ is a big vertex, we say $\Fkt{\pi_{\delta}}{\mathsf{v}}$ is an \emph{endpoint} of $Q$ if an endpoint of $Q$ is drawn within $\mathsf{v}$ by $\Gamma$.  
Now let $P$ be either a (not necessarily directed) cycle or a path in $D$ that uses no edge of $\sigma(c)$ of any maelstrom $c\in C(\delta)$.
We say that $P$ is \emph{grounded in $\delta$} if either $P$ is a non-zero length path with both endpoints in $\pi_{\delta}(N(\delta))$, or $P$ is a cycle and it uses edges of $\sigma(c_1)$ and $\sigma(c_2)$ for two distinct cells $c_1,c_2\in C(\delta)$.

For every big vertex $\mathsf{v}\in\mathcal{V}$ and every crossing edge $e$ of $\mathsf{v}$, we select one of the components $t_1$ of $\Boundary{\mathsf{v}}-v_{\mathsf{v}}-e$.
Additionally, for every common cell $c \in C(\delta)$ intersecting a big vertex $\mathsf{v}$, with $v \in \widetilde{c}$ being the other vertex on the boundary, and any crossing edge $e$ for $\mathsf{v}$ with $v$ as an endpoint, we select a component $t_2$ of $\Boundary{c}-v-e$.
Similarly, for every common cell $c' \in C(\delta)$ and every conjunction cell $c'' \in C(\delta)$ with $\Abs{\widetilde{c''}}=2$, we select one of the components $t_3$ of $\Boundary{c'}-\widetilde{c'}$ and a component $t_4$ of $\Boundary{c''}-\widetilde{c''}$ respectively.
We refer to $(t_1 \cup v_{\mathsf{v}} \cup e) \cap \Boundary{\mathsf{v}}$, $( t_2 \cup \widetilde{c} \cup e) \cap \Boundary{c}$, $t_3 \cup \widetilde{c'}$, and $t_4 \cup \widetilde{c''}$, all of which form connected components, as \emph{tiebreakers} in $\delta$, and we assume, without any further mention, that every strong plane decomposition comes equipped with tiebreakers.

If $P$ is grounded, we define the \emph{trace} of $P$ as follows.
Let $Q_1,\dots,Q_k$ be the distinct non-trivial, maximal subpaths of $P$ such that $Q_i$ is a subgraph of $\sigma(c)$ for some cell $c$, or subgraph of $\sigma(\mathsf{v})$ for some big vertex $\mathsf{v}\in\mathcal{V}$.
Fix an index $i \in [k]$.
The maximality of $Q_i$ and the fact that it is non-trivial implies that it has two endpoints, any of which is either a non-door vertex in $N(\delta)$, the door $v_{\mathsf{v}}$ of some big vertex $\mathsf{v}\in\mathcal{V}$, or a vertex drawn in the interior of some big vertex $\mathsf{u}\in N(\delta)$.

Should one endpoint of $Q_i$ be an endpoint of $P$ that is drawn in the interior of a big vertex $\mathsf{v}$, then there exists a unique crossing edge $e \in E(P)$ of $\mathsf{v}$ and we let $L_i$ be the tiebreaker for $e$ and $\mathsf{v}$.
If one endpoint of $Q_i$ is an endpoint of $P$ that is drawn in the interior of the conjunction cell, we define $L_i$ to be the tiebreaker of $c$, if $|\widetilde{c}| = 2$.
Should $|\widetilde{c}| = 3$, let $a \in \widetilde{c} \cap V(Q_i)$, then we define $L_i$ to be the shortest curve in $c$ that connects $a$ and the segment of $\Boundary{c} - a$ that does not contain $a$.
Finally, if $|\widetilde{c}| = 4$, let $a \in \widetilde{c} \cap V(Q_i)$ and let $b \in \widetilde{c}$ be chosen such that both components of $\Boundary{c} - a - b$ contain an element from $\widetilde{c}$, we define $L_i$ to be the shortest curve in $c$ that connects $a$ and $b$.\footnote{This definition does not accurately reflect the vertices which $P$ contains, but is sufficiently representative for our purposes.}

Now we can assume that $Q_i$ does not contain an endpoint of $P$.
If $c$ is a common cell, with $e$ being the crossing edge in $E(P)$ if $\Abs{\widetilde{c}}=1$, or a conjunction cell with $\Abs{\widetilde{c'}}=2$, we define $L_i$ as the respective tiebreaker provided by $\delta$.
For any big vertex $\mathsf{v}$, let $e \in E(P)$ be the crossing edge with an endpoint on $Q_i$, we again define $L_i$ as the tiebreaker provided by $\delta$.
If $Q_i$ is contained in a conjunction cell $c'\in C(\delta)$ with $\Abs{\widetilde{c'}}=3$, we define $L_i$ to be the union of $x,y$ and the component of $\Boundary{c'}-x-y$ which is disjoint from $\widetilde{c'}\setminus\Set{x,y}$, where $x,y\in N(\delta)$ are the two points of $Q_i$ drawn on $\Boundary{c'}$.

The case where $Q_i$ is contained in a conjunction cell $c'\in C(\delta)$ with $\Abs{\widetilde{c'}}=4$ is a bit more complicated.
Let $z_1,z_2\in N(\delta)$ be the two nodes of $\widetilde{c'}$ which are not part of $Q_i$.
It is possible that each of the two components of $\Boundary{c'}-z_1-z_2$ contains a part of $Q_i$.
In this case we define $L_i$ to be a shortest curve within $c'$ connecting the two nodes of $\widetilde{c'}$ which are in $Q_i$.
It follows from the results in \cite{giannopoulou2021two} that there can never be two disjoint directed paths forming a cross over $\Boundary{c'}$ in the situation where $\delta$ is a maelstrom-free, odd decomposition, and thus we deem this choice for $L_i$ to be justified.
If, however, there does exists a component of $\Boundary{c'}-z_1-z_2$ which contains both endpoints of $Q_i$, we define $L_i$ to be the shortest curve connecting the endpoints of $Q_i$ within this component.

The \emph{trace} of $P$ in $\delta$, denoted by $\Trace_\delta(P)$, is defined to be $\bigcup_{i\in[k]}L_i$.
Since a trace is unique up to homeomorphic shifts, we can speak of the trace of $P$ without serious ambiguity.
If $\delta$ is clear from the context, we speak simply of the trace of $P$, denoted as $\Trace(P)$.
Thus the trace of a cycle is the homeomorphic image of the unit circle, and the trace of a path is an arc with both endpoints on the boundaries of members of $N(\delta)$.

If $H$ is a subgraph of $D$, we say that $H$ is \emph{grounded} if every path and cycle in $D$ is grounded and the \emph{trace of $H$} is defined to be the union of the traces of all paths and cycles in $H$.
\end{definition}

\begin{definition}[Concentric and homogeneous cycles]\label{def:concentric}
    Let $D$ be a digraph and let $\delta=(\Gamma,\mathcal{V},\mathcal{D})$ be a strong plane decomposition of $D$.
    We say that a family of pairwise disjoint cycles $\mathcal{C} = \{ C_1, C_2, \ldots , C_k \}$ is \emph{concentric in} $\delta$ if for every $i \in[2,k-1]$ the trace of $C_i$ separates $C_{i-1}$ from $C_{i+1}$.
    A family of pairwise disjoint, concentric dicycles $\mathcal{C} = \{ C_1, C_2, \ldots , C_k \}$ is called \emph{homogeneous in $\delta$} if all of them are oriented either clockwise or counter-clockwise.
    If $\delta$ is clear from the context, we simply call a family of cycles concentric or homogeneous.
\end{definition}

\subsection{The inner workings of a maelstrom}

The concept of maelstroms is not exactly the same as the idea of vortices in the Graph Minors Series.
Indeed, it is possible to further subdivide a maelstrom into several vortices.
However, a maelstrom still represents an area in an odd decomposition where we cannot guarantee strong planarity, at least not in the sense which is required for renditions.
Instead, we make use of a slightly more relaxed concept.

Let $(D,\Omega)$ be a society and let $I\subseteq\V{\Omega}$.
We say that $I$ is a \emph{segment} of $\Omega$ if there do not exist distinct elements $x_1,x_2\in I$ and $y_1,y_2\in\V{\Omega}\setminus I$ such that $x_1,y_1,x_2,y_2$ occur in $\Omega$ in the order listed.
A vertex $x$ is an \emph{endpoint} of the segment $I$ if there is a vertex $y\in\V{\Omega}\setminus I$ which immediately succeeds or precedes $x$ in the order $\Omega$.
Note that every segment $I$ that is a non-empty proper subset of $\V{\Omega}$ has endpoints, one of which is \emph{first} while the other one is \emph{last} with respect to the order of $I$ given by $\Omega$.
Given two vertices $x,y\in\V{\Omega}$, we denote by $x\Omega y$ the uniquely determined segment with first vertex $x$ and last vertex $y$.
In case $y$ immediately precedes $x$, we define $x\Omega y$ to be the trivial segment $\V{\Omega}$.

\begin{definition}[Transaction]\label{def:transaction}
	Let $(D,\Omega)$ be a society.
	A \emph{transaction} on $(D,\Omega)$ is a linkage $\mathcal{L}$ in $D$ for which two disjoint segments $X,Y\subseteq\V{\Omega}$ exist such that $\mathcal{L}$ is an $X$-$Y$-linkage and the paths in $\mathcal{L}$ are internally disjoint from $\V{\Omega}$.
\end{definition}

An \emph{anti-directed cycle} is a cycle $C$ in a multi-digraph where each vertex is either the head of both of its incident edges, or their tail.
A \emph{subdivided anti-directed cycle} is obtained from an anti-directed cycle by replacing every edge by a directed path of length at least one.
Note that this definition allows for an anti-directed cycle of length two, made up of two parallel directed edges with the same head and the same tail.
In particular, a subdivided anti-directed cycle obtained from an anti-directed cycle of length two is called a \emph{diamond}.

\begin{definition}($d$-disks of grounded cycles)\label{def:cdisk}
    Let $D$ be a digraph and let $\delta=(\Gamma,\mathcal{V},\mathcal{D})$ be a strong plane decomposition of $D$.
    Let $C$ be a (not necessarily directed) cycle in $D$ which is grounded in $\delta$ and let $d$ be some disk whose interior is disjoint from the trace of $C$.
    The unique disk $\Delta$ which is bounded by the trace of $C$ and contains $d$ is called the \emph{$d$-disk} of $C$.
    On the other hand, the disk bounded by the trace of $C$ that does not contain $d$ is called the \emph{non-$d$-disk} of $C$.
    Given some subgraph $H\subseteq D$ such that $H$ contains a (not necessarily directed) cycle which is grounded in $\delta$, we say that a cycle $C\subseteq H$ is the \emph{$d$-tight cycle of $H$ in $\delta$} if $C$ is grounded in $\delta$ and for every cycle $C'\subseteq H$, either $C'$ is not grounded in $\delta$, or the $d$-disk of $C$ is properly contained in the $d$-disk of $C'$.
\end{definition}

Observe that the $m$-tight cycle of $H$ in the definition above is uniquely determined up to identity of traces.
Our next definition is the foundation of our structural results.

\begin{definition}[Outline of a maelstrom]\label{def:outline}
	Let $D$ be a digraph and let $\delta=(\Gamma,\mathcal{V},\mathcal{D})$ be an odd decomposition of $D$.
	Let $m\in C(\delta)$ be a maelstrom of $\delta$ and let $t,\theta$ be two integers, where $\theta$ is positive and even, and $t$ is non-negative.
	A tuple $\mathfrak{M}=(H,\mathcal{C},\mathfrak{E},\mathfrak{V})$, where $H$ is a subgraph of $D$, $\mathcal{C}=\Set{C_1,\dots,C_{\theta}}$ is a homogeneous family of dicycles, $\mathfrak{E}=\Set{\mathcal{E}_1,\dots,\mathcal{E}_{t}}$ is a collection of $t$ linkages, with $\mathcal{E}_i=\Set{E_1^i,\dots,E_{\theta}^i}$ for all $i\in[t]$, and $\mathfrak{V}$ is a non-empty collection of closed disks found in the interior of $m$, is called a \emph{$\theta$-outline} of $m$ of \emph{degree $t$} if it meets the following requirements.

	\begin{description}
		\item[M1] The graph $\bigcup_{i \in [\theta]}C_i \cup \bigcup_{i \in [t]}\bigcup_{j \in [\theta]} E^i_j$ is a subgraph of $H$.
		\item[M2] Let $H_0\coloneqq H$.
        The society $(H_0,\Omega_0)$ has a non-even rendition $\rho=(\Gamma_{H_0},\mathcal{V}_{H_0},\mathcal{D}_{H_0})$ in a disk $\Delta_0$, where $\V{\Omega_0}=\V{C_{\theta}}$ and all cycles in $\mathcal{C}$ and all paths in $\bigcup_{i \in [t]} \mathcal{E}_i$ are grounded in $\rho$.
		\item[M3] For all $i \in [0,t]$, let $O_i$ be the $m$-tight cycle of $C_{\nicefrac{\theta}{2}} \cup \bigcup_{j\in[i]}E^j_{\nicefrac{\theta}{2}}$ in $\rho$, we denote the $m$-disk of $O_i$ by $c_i$.
		Moreover, we require that $c_t = m$.
		\item[M4] For all $i \in [t]$, let $Q_i$ be the $m$-tight cycle of $C_{\theta} \cup \bigcup_{j\in[i]}E^j_{\theta}$ in $\rho$, we denote its $m$-disk by $\Delta_i$.
		We define the society $(H_i, \Omega_i)$, with $V(\Omega_i) = V(Q_i)$ and $H_i \subseteq H$ as the graph induced by all vertices and big vertices of $H$ drawn on $\Delta_i$ by $\Gamma_{H_0}$.
		\item[M5] For all $i \in [t]$, the linkage $\mathcal{E}_i$ is a transaction on $(H_{i-1},\Omega_{i-1})$.
		\item[M6] For every $i\in[t]$ and $j \in [2,\theta -1]$, the trace of $E_j^i$ separates $E_{j-1}^i$ from $E_{j+1}^i$ in $\Delta_i$.
		\item[M7] 
        The restriction of $\rho$ to all nodes drawn on and all conjunction cells intersecting $\bigcup_{i \in [0,t]} (\Delta_i \setminus (c_i \setminus \Boundary{c_i}))$ is an odd rendition that is a restriction of $\delta$. 
        \item[M8] We require that for all $i \in [t]$, the cycle $Q_i$ is either a diamond or a dicycle.

        \item[M9] For all $i \in [t]$, let $\Delta_i'$ be the disk in $\Delta_{i-1} - \Trace_{\delta}(E^i_\theta)$ that does not contain $m$.
        If $Q_i$ is a dicycle, the intersection of $c_{i-1}$ and $\Delta_i'$ must be a disk or empty.\footnote{Note that, since we do not demand this for diamond outlines, the paths in $\mathcal{E}_i$ may interact with the existing infrastructure at an arbitrary depth on the side of $Q_i$ that is turned away from the maelstrom.
		However, they cannot intersect $Q_{i-1}$ there, since this cycle forms their society.}

		\item[M10] The vertices of $\Fkt{\sigma}{m}-\V{H}$ are drawn by $\Gamma$ into the disks from $\mathfrak{V}$, which are called the \emph{vortex regions} or simply \emph{vortices} of $\mathfrak{M}$.
		Moreover, no vertex of $H$ is drawn in the interior of some disk from $\mathfrak{V}$ and every disk of $\mathfrak{V}$ is a subset of $\Delta_t$.
	\end{description}
    We call $\mathfrak{M}$ a \emph{diamond outline} if $Q_t$ is a diamond, and if $Q_t$ is a dicycle, $\mathfrak{M}$ is a \emph{circle outline}.
	The tuple $(H_{\mathfrak{M}},\Omega_{\mathfrak{M}})$ is the \emph{$\mathfrak{M}$-society} of $m$, where $\Omega_{\mathfrak{M}}\coloneqq \Omega_t$ and $H_{\mathfrak{M}}\coloneqq H_t\cup\sigma(m)$.
\end{definition}

The notion of maelstroms arises from the problem that we need to ensure sufficient infrastructure to support the continued decomposition of our digraph throughout our proofs.
In \cite{kawarabayashi2020quickly} this is done through the notions of nests.
These are large families of concentric cycles such that the outermost cycle is connected to the innermost by a large linkage.
Since \cite{kawarabayashi2020quickly} deals with undirected graphs, it suffices to ask for a large transaction on the society defined by the innermost cycle of a nest to ensure sufficiently large connectivity of the transaction with the society defined by the outermost cycle of the nest.
This has the advantage that one can ensure that no path which has entered the vortex will leave the vortex again and return a second time.
However, in the setting of digraphs these kinds of arguments are impossible, since we have to take the direction of our linkages into account.
What would be a single vortex can be split into several regions by a transaction only returning slightly below the surface of the surrounding ``nest'' of dicycles.
Since there might still be edges between these regions, we cannot fully integrate the infrastructure surrounding our vortex-analogues into the odd decompositions.
Hence the need for maelstroms to account for these problems.

\begin{figure}
    \centering
    \scalebox{0.85}{
    \begin{tikzpicture}[scale=1]

        \pgfdeclarelayer{background}
		\pgfdeclarelayer{foreground}
			
		\pgfsetlayers{background,main,foreground}
			
        \begin{pgfonlayer}{main}
        \node (C) [v:ghost] {};

            \pgftext{\includegraphics[width=10cm]{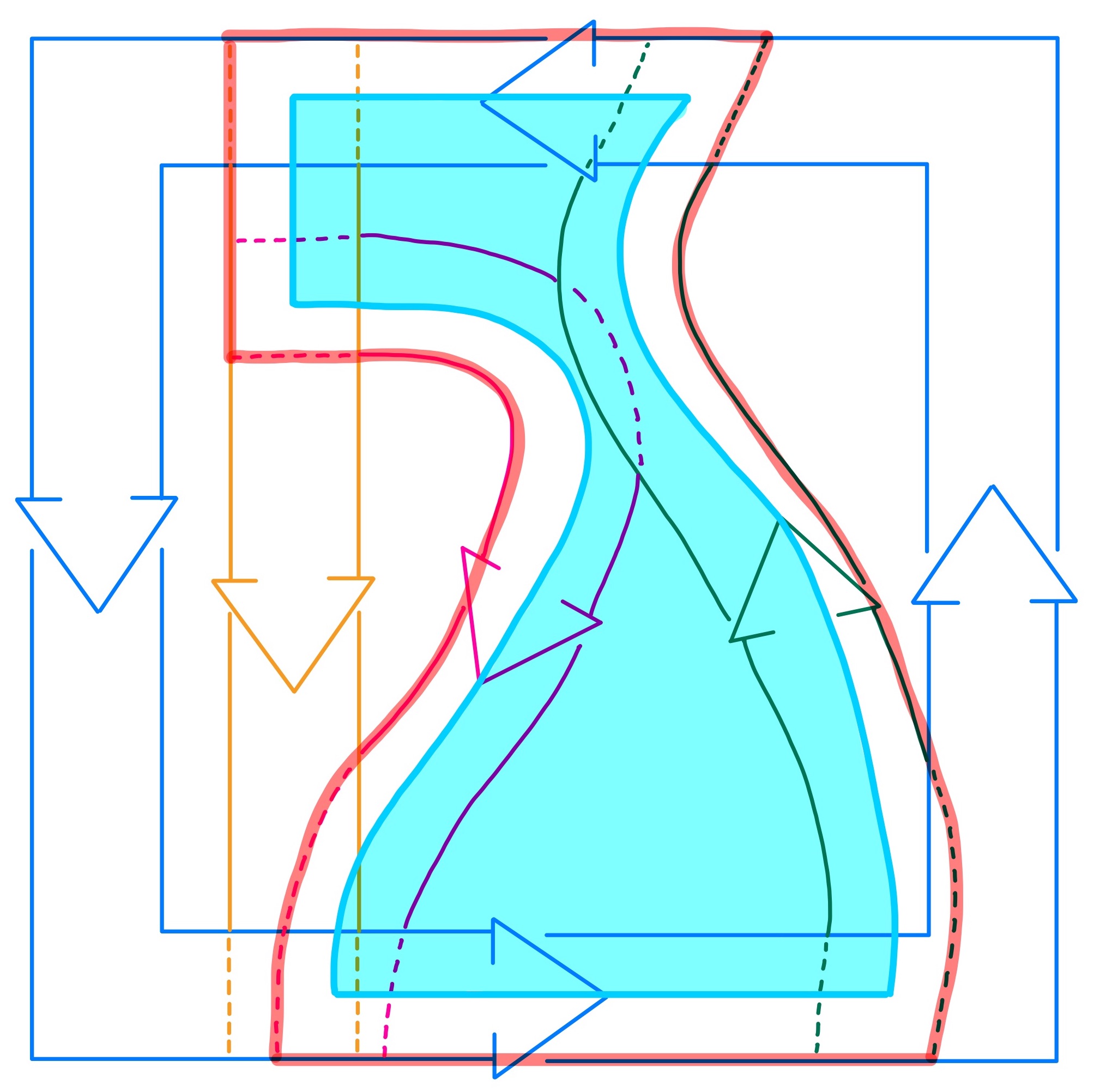}} at (C.center);

            \node (society) [v:ghost,position=229:53mm from C] {$(H_{\mathfrak{M}},\Omega_{\mathfrak{M}})$};

            \node (society2) [v:ghost,position=75:45.2mm from C] {$O_{\mathfrak{M}}$};

            \node (maelstrom) [v:ghost,position=340:10mm from C] {$c_0$};
            
        \end{pgfonlayer}{main}
        
        \begin{pgfonlayer}{foreground}
        \end{pgfonlayer}{foreground}

        \begin{pgfonlayer}{background}
        \end{pgfonlayer}{background}
        
    \end{tikzpicture}
    }
    \caption{A circle maelstrom together with the original wall in blue, and the transactions that contain its outline. This illustrates that transactions may intersect each other even if they lie on opposite sides of the outline, as seen here with the magenta transaction and the green transaction. However, this interaction is not allowed to pass the halfway point of the paths in either linkage. Thus $O_\mathfrak{M}$ defines a single disk that describes the maelstrom.}
    \label{fig:maelstrom}
\end{figure}
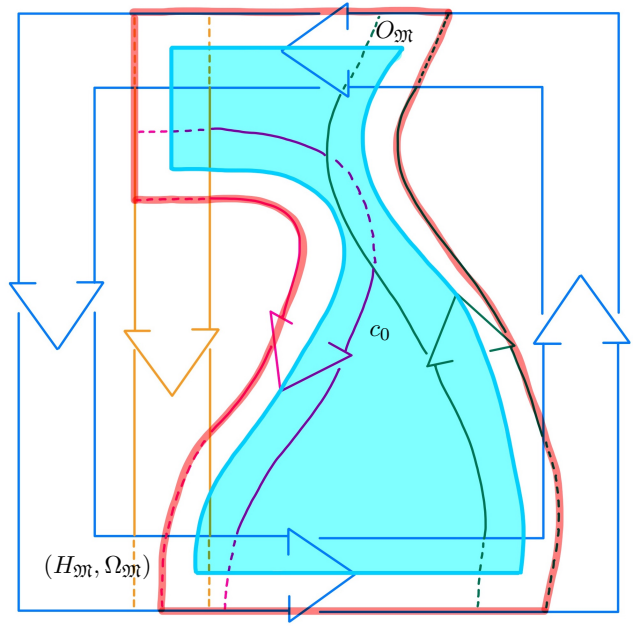

\subsection{Depth of a maelstrom outline}\label{subsec:nestsandtransactions}

Since our maelstrom outlines are a tool for induction, we need to describe the situation in which the inductive step may be applied to a fixed maelstrom outline.
This will be true whenever we are able to find a large enough transaction on the maelstrom society of the respective outline.
However, in the case of diamond maelstroms this is not enough as some of these transactions are not compatible with the inherent direction of the diamond itself.
Hence, we need to fix definitions of ``depth'', i.\@e.\@ the measure of a largest transaction on the respective maelstrom society which may be used to continue the construction, which differentiate between circle and diamond maelstroms.

\begin{definition}[Depth of a circle maelstrom]\label{def:depthofacircle}
    Let $\theta$ be a positive integer, let $D$ be a digraph, let $\delta=(\Gamma,\mathcal{V},\mathcal{D})$ be an odd decomposition of $D$, and let $c_0 \in C(\delta)$ be a maelstrom of $\delta$, with a \textbf{circle} $\theta$-outline $\mathfrak{M}=(H,\mathcal{C},\mathfrak{E},\mathfrak{V})$.
	
    The \emph{depth of $\mathfrak{M}$} is the largest order of a transaction on the maelstrom society $(H_{\mathfrak{M}},\Omega_{\mathfrak{M}})$ of $\mathfrak{M}$.
\end{definition}

\begin{definition}[Streams and the directional depth of a diamond maelstrom]\label{def:streamsanddirectionaldepth}
    Let $\theta$ be a positive integer, let $D$ be a digraph, let $\delta=(\Gamma,\mathcal{V},\mathcal{D})$ be an odd decomposition of $D$, and let $c_0 \in C(\delta)$ be a maelstrom of $\delta$, with a \textbf{diamond} $\theta$-outline $\mathfrak{M}=(H,\mathcal{C},\mathfrak{E},\mathfrak{V})$.
    Let $(H_{\mathfrak{M}},\Omega_{\mathfrak{M}})$ be the maelstrom society of $\mathfrak{M}$ and let $P_1,P_2$ be the two directed paths which together form the diamond corresponding to $\Omega_{\mathfrak{M}}$.
    
	Let $\mathcal{L}$ be a transaction on $(H_{\mathfrak{M}},\Omega_{\mathfrak{M}})$ starting on the segment $X\subseteq\V{\Omega_{\mathfrak{M}}}$ and ending on the segment $Y\subseteq\V{\Omega_{\mathfrak{M}}}$.
	If there exists an $i\in[2]$ such that $X\cup Y\subseteq\V{P_i}$ and the vertices of $Y$ appear on $P_i$ before the vertices of $X$, we say that $\mathcal{L}$ is \emph{choppy}.
	On the other hand, if $\mathcal{L}$ contains no choppy transaction\footnote{Notice that a single path moving ``backwards'' along one of the $P_i$ is considered to be choppy transaction.}, then $\mathcal{L}$ is called a \emph{stream} on $(H_{\mathfrak{M}},\Omega_{\mathfrak{M}})$.

	The \emph{directional depth} of $\mathfrak{M}$ is the largest order of a stream on $(H_{\mathfrak{M}},\Omega_{\mathfrak{M}})$.
\end{definition}

Once the directional depth of the tight outlines of a diamond maelstrom is bounded, we may start with the second phase of refining the plane decomposition $\delta$.
In the previous step it was possible that a new stream would split $c_0$ into several new maelstroms.
This, however, will no longer occur in this second step.
We will now exclusively attach choppy transactions to $(H_{\mathfrak{M}},\Omega_{\mathfrak{M}})$, until no large choppy transaction can be found anymore.

Afterwards there might still exist large transactions on the ``refined'' maelstrom society, but for any of the remaining transactions we can guarantee that they will be pushed into a part of $\delta$ which we can treat differently.
Hence, for our purposes this second step, although it does not yield a bound on the order of all transactions, is good enough.

\begin{definition}[Refined diamond outlines and eddies]\label{def:refinedoutline}
    Let $t, \theta$ be positive integers and let $k, g, h$ be non-negative integers, where $\theta$ is even and $h \geq k$.
    Let $D$ be a digraph, let $\delta=(\Gamma,\mathcal{V},\mathcal{D})$ be an odd decomposition of $D$, and let $m \in C(\delta)$ be a maelstrom of $\delta$, with a \textbf{diamond} $\theta$-outline $\mathfrak{M}=(H,\mathcal{C},\mathfrak{E},\mathfrak{V})$ of degree $t$,
    with the maelstrom society $(H_\mathfrak{M}, \Omega_{\mathfrak{M}})$.
    
    A \emph{refined diamond $\theta$-outline of} $\mathfrak{M}$ around a closed disk $n \subseteq m$ of \emph{turbulence $h$}, \emph{surplus degree $k$}, and \emph{roughness} $g$ is a tuple $\mathfrak{N} = (H^*,\mathcal{C},\mathfrak{E},\mathfrak{B},\mathfrak{V}^*,\mathfrak{D})$.
    Moreover, $\mathfrak{B}=\Set{\mathcal{B}_1,\dots,\mathcal{B}_k}$ is a possibly empty collection of $k$ linkages such that, for each $i\in[k]$, $\mathcal{B}_i= \{B_{1}^{i},\dots,B_{\theta(g+1)} \}$, and $\mathfrak{D}$ is a collection of closed disks that are disjoint from $n$, such that
    \begin{description}
        \item[E1] The graph $H\cup\bigcup_{i=1}^k\bigcup_{j=1}^{\theta(g+1)} B^i_j$ is a subgraph of $H^*$.
        
        \item[E2] Let $H^*_0=H^*$, and let $\Delta_0$ and $(H_0,\Omega_0)$ be as in the definition of $\mathfrak{M}$.
        The society $(H^*_0,\Omega_0)$ has a non-even rendition $\rho' = (\Gamma'_{H^*_0}, \mathcal{V}'_{H^*_0}, \mathcal{D}'_{H^*_0})$ in $\Delta_0$ which contains the rendition $\rho$ of $(H_0,\Omega_0)$, such that all paths of the linkages in $\mathfrak{B}$ are grounded in $\rho'$.
        
        \item[E3] For all $i \in [0,t]$, let $O_i$ and $c_i$ be defined as in $\mathfrak{M}$.
        For all $i \in [t+1,t+k]$, let $O_i$ be the $n$-tight cycle of $C_{\nicefrac{\theta}{2}} \cup \bigcup_{j\in[t]}E^j_{\nicefrac{\theta}{2}} \cup \bigcup_{j \in [i-t]}B^j_{\nicefrac{\theta(g+1)}{2}}$ in $\rho'$.
        We denote the $n$-disk of $O_i$ by $c_i$.
        Notice that $O_i$ is grounded in $\delta$ and $\rho'$, and we require that $c_{t+k} = n$.

        \item[E4] For all $i \in [t]$, let $Q_i$, $\Delta_i$, and $(H_i,\Omega_i)$ be defined as in $\mathfrak{M}$.
        For all $i \in [t+1,t+k]$, let $Q_i$ be the $n$-tight cycle of $C_{\theta} \cup \bigcup_{j\in[t]}E^j_{\theta} \cup \bigcup_{j \in [i-t]}B^j_{\theta(g+1)}$ in $\rho'$.
        We denote the $n$-disk of $Q_i$ by $\Delta_i$.
        For $i\in[t+1,t+k]$, we define the society $(H_i,\Omega_i)$, such that $\V{\Omega_i} = \V{Q_i}$ and $H_i \subseteq H^*$ is the graph induced by all vertices and big vertices of $H^*$ drawn on $\Delta_i$ by $\Gamma'_{H^*_0}$.
        
        \item[E5] $\mathcal{B}_1$ is a choppy transaction on $(H_t,\Omega_t)$.
        
        \item[E6] Let $P_1^1,P_2^1$ be the two directed paths which together make up the diamond corresponding to $\Omega_{\mathfrak{M}}=\Omega_t$ and let $u$ be the common tail, $v$ be the common head of the two paths.
        For every $i \in [2, k]$ and $j \in [2]$ we define the paths $P^i_j$ inductively as follows:
        \begin{itemize}
            \item $P^i_j$ has endpoints $u$ and $v$, and we treat $P^i_j$ as if it was oriented from $u$ to $v$ (although it is not necessarily a directed path).
            \item In case $B^{i-1}_{\nicefrac{\theta(g+1)}{2}}$ has at most one endpoint on $P^{i-1}_j$, we set $P^i_j\coloneqq P^{i-1}_j$.
            \item Otherwise, $B^{i-1}_{\nicefrac{\theta(g+1)}{2}}$ has both endpoints, say $x$ and $y$, on $P^{i-1}_j$, and they occur, with respect to the imposed orientation, in the order listed.
            In this case, we set $P^i_j\coloneqq uP^{i-1}_jxB^{i-1}_{\nicefrac{\theta(g+1)}{2}}yP^{i-1}_jv$.
        \end{itemize}
        
        For every $i\in[1,k]$, $\mathcal{B}_i$ is a transaction on $(H_{t+i-1},\Omega_{t+i-1})$.
        Let $X_i$ be the minimal segment of $\Omega_{t+i-1}$ containing all tails of the paths in $\mathcal{B}_i$ and $Y_i$ be the minimal segment of $\Omega_{t+i-1}$ containing all heads of the paths in $\mathcal{B}_i$.
        Then there exists $j\in[2]$ such that $X_i\cup Y_i\subseteq V(P^i_j)$ and the vertices of $Y_i$ appear on $P^i_j$ before the vertices of $X_i$ with respect to the orientation of $P^i_j$ imposed above\footnote{Notice that this means that each of the $\mathcal{B}_i$ may be extended to a choppy transaction on $(H_t,\Omega_t)$}.
        We say that $\mathcal{B}_i$ \emph{attaches} to $P^i_j$.
        
        \item[E7] For every $i \in [k]$, $j\in[2,\theta(g+1)-1]$ the trace of $B_j^i$ separates $B_{j-1}^i$ from $B_{j+1}^i$ in $\Delta_{t+i}$.

        

        \item[E8] 
        The restriction of $\rho'$ to all nodes drawn on common and conjunction cells intersecting $\Delta_0 - ( (n - \Boundary{n}) \cup \bigcup_{d \in \mathfrak{D}} (d - \Boundary{d} ) )$ is an odd rendition that is a restriction of $\delta$.
        
        
        \item[E9] The vertices of $\sigma(n)-\V{H^*}$ are drawn by $\Gamma$ into the disks from $\mathfrak{V}^*$, which are called the \emph{vortex regions} or simply \emph{vortices} of $\mathfrak{N}$.
        Moreover, the vertices of $\sigma(m-(V{(H^*)}\cup \sigma(n)))$ are drawn by $\Gamma$ into the disks from $\mathfrak{D}$.
		Finally, no vertex of $H^*$ is drawn in the interior of some disk from $\mathfrak{V}^* \cup \mathfrak{D}$.
        We call the disks in $\mathfrak{D}$ the \emph{eddy regions} or simply \emph{eddies} of $\mathfrak{N}$.

        \item[E10]For every pair $d,d'\in\mathfrak{D}\cup\{ n\}$, every (not necessarily directed) path $P$ in $D$ with one endpoint in $\sigma(d)$ and the other in $\sigma(d')$ contains a vertex of $H^* - (\sigma(d) \cup \sigma(d'))$.
    \end{description}
    For each $i \in [2]$, the path $P^{t+k}_i$ can be split into subpaths $S_i^1, \ldots , S_i^{h(i)}$ by removing all vertices from $P^{t+k}_i$ that are exclusive to $V(\Omega_t)$ and found in none of the linkages in $\mathfrak{B}$, such that $h = h(1) + h(2)$.
    We call $S_1^1, \ldots , S_1^{h(1)}, S_2^1, \ldots , S_2^{h(2)}$ the \emph{eddy segments} of $\Omega_\NN$.

    We call $(H_\NN,\Omega_\NN)$ the \emph{$\mathfrak{N}$-society} of $n$, where $H_\NN = H_{t+k} \cup \sigma(n)$ and $\Omega_\NN = \Omega_{t+k}$. 
    Finally, $\mathfrak{N}$ is said to be \emph{proper} if $H_{\mathfrak{M}}-\bigcup_{d\in\mathfrak{D}}\sigma(d)$ contains an even dicycle.
\end{definition}

The reason the additional parameters for the refined outline are all allowed to be zero is that this later allows us to combine proofs for diamond outlines and refined diamond outlines, since each diamond outline $\mathfrak{M}=(H,\mathcal{C},\mathfrak{E},\mathfrak{V})$ can be seen as a refined diamond outline $\mathfrak{N}=(H,\mathcal{C},\mathfrak{E}, { \mathfrak{B}= \emptyset } ,\mathfrak{V},\mathfrak{D} = \emptyset)$.

\begin{figure}
    \centering
    \scalebox{0.85}{
    \begin{tikzpicture}[scale=1]

        \pgfdeclarelayer{background}
		\pgfdeclarelayer{foreground}
			
		\pgfsetlayers{background,main,foreground}
			
        \begin{pgfonlayer}{main}
        \node (C) [v:ghost] {};

            \pgftext{\includegraphics[width=10cm]{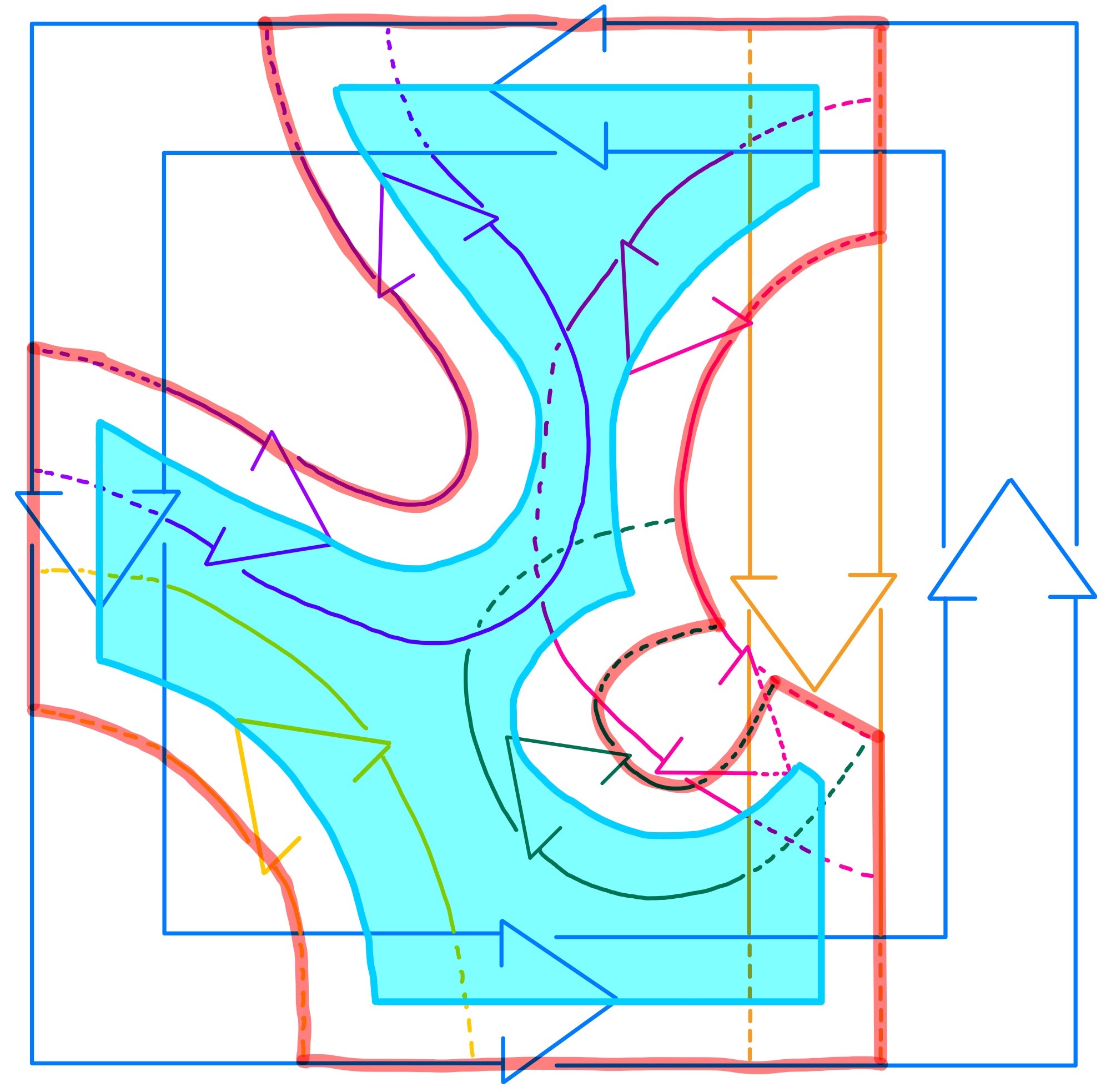}} at (C.center);

            \node (society) [v:ghost,position=232:52mm from C] {$(H_{\mathfrak{M}},\Omega_{\mathfrak{M}})$};

            \node (society2) [v:ghost,position=61:50.2mm from C] {$O_{\mathfrak{M}}$};

            \node (maelstrom) [v:ghost,position=225:18mm from C] {$c'_0$};

            \node (a) [v:ghost,position=225:43mm from C] {$(a)$};
            \node (b) [v:ghost,position=141:30mm from C] {$(b)$};
            \node (c) [v:ghost,position=309:19mm from C] {$(c)$};
            
        \end{pgfonlayer}{main}
        
        \begin{pgfonlayer}{foreground}
        \end{pgfonlayer}{foreground}

        \begin{pgfonlayer}{background}
        \end{pgfonlayer}{background}
        
    \end{tikzpicture}
    }
    \caption{A refined diamond maelstrom in blue and its outline.
    The areas $(a)$, $(b)$, and $(c)$ are the eddies of the maelstrom.
    As depicted, the transactions defining the eddies may overlap and create vortices as in \Cref{def:outline}.
    However, no transaction in the infrastructure may travel through more than half of the paths of any other transaction.
    This gives us a single disk described by the infrastructure that we choose as our maelstrom.}
    \label{fig:refineddiamond}
\end{figure}
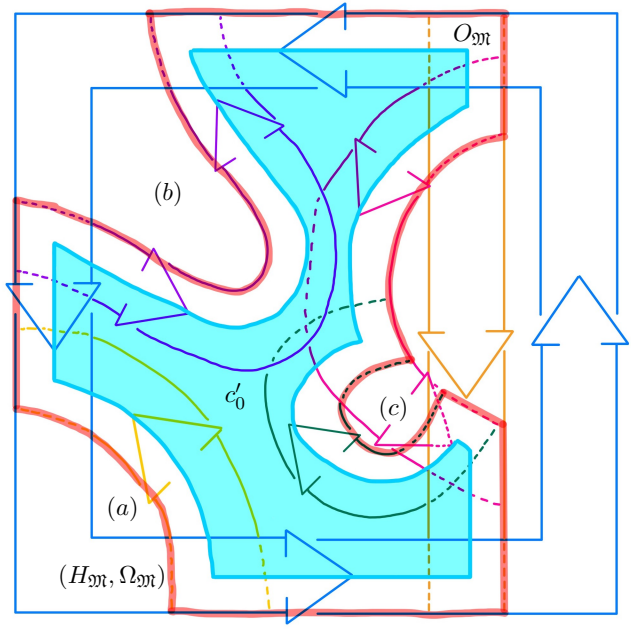

\begin{definition}[Calm transactions, whirls, and the depth of a refined diamond outline]\label{def:calmtransaction}
    Let $\theta$ be a positive integer and let $h,k$ be non-negative integers, with $\theta$ being even and $h \geq k$.
    Let $\delta=(\Gamma,\mathcal{V},\mathcal{D})$ be an odd decomposition of $D$, and let $m \in C(\delta)$ be a maelstrom of $\delta$, with a diamond $\theta$-outline $\mathfrak{M}=(H,\mathcal{C},\mathfrak{E},\mathfrak{V})$.
    Furthermore, let $\mathfrak{N} = (H^*,\mathcal{C},\mathfrak{E},\mathfrak{B},\mathfrak{V}^*,\mathfrak{D})$ be a \textbf{refined diamond} $\theta$-outline in $\mathfrak{M}$ around $n \subseteq m$ with turbulence $h$.
    Let $(H_\mathfrak{N}, \Omega_{\mathfrak{N}})$ be the $\mathfrak{N}$-society and $S_1, \ldots , S_h$ be the eddy segments of $\Omega_{\mathfrak{N}}$.

    A directed $\V{\Omega_{\mathfrak{N}}}$-path $P$ in $H_\mathfrak{N}$ is called a \emph{whirl} if there exists an $i \in [h]$, such that both endpoints of $P$ are on $S_i$.
    A transaction $\mathcal{P}$ on $(H_\mathfrak{N}, \Omega_{\mathfrak{N}})$ is called \emph{whirly} if some path in $\mathcal{P}$ is a whirl and if $\mathcal{P}$ is not whirly, then we call it \emph{calm}.
    
    The \emph{depth} of $\mathfrak{N}$ is the maximum order of a calm transaction on $(H_\mathfrak{N}, \Omega_{\mathfrak{N}})$.
\end{definition}

Observe that any large calm transaction on $(H_\mathfrak{N}, \Omega_{\mathfrak{N}})$ must either contain (or at least be extendable to) a large stream on the original maelstrom society $(H_{\mathfrak{M}},\Omega_{\mathfrak{M}})$, or a large order choppy transaction on $(H_{\mathfrak{M}},\Omega_{\mathfrak{M}})$.
Hence, the depth of $\mathfrak{N}$ is definitely an upper bound on the directional depth of $\mathfrak{M}$.

\subsection{Decomposing a digraph}\label{subsec:decomposing}

The directed tight cut decomposition is a way to decompose a given digraph in a tree-like manner into its dibraces.
Similarly, \cref{thm:nonevenstructure} yields a way to further decompose the dibraces of a non-even digraph, whilst maintaining a tree structure.
One possible way to generalise this idea can be found in the notion of directed treewidth \cite{Johnson2001DirectedTreewidth}.

Let $D$ be a digraph and $X,Y\subseteq\V{D}$.
A directed walk $W$ is a \emph{directed $X$-walk} if it starts and ends in $X$, and contains a vertex of $\V{D-X}$.
We say that $Y$ \emph{strongly guards} $X$ if every directed $X$-walk in $D$ contains a vertex of $Y$.
The set $Y$ \emph{weakly guards} $X$ if every directed $X$-$\V{D-X}$-path contains a vertex of $Y$, 

An \emph{arborescence} is a digraph $\vec{T}$ obtained from a tree $T$ by selecting a \emph{root} $r\in\V{T}$ and orienting all edges of $T$ away from $r$.
If $e$ is a directed edge and $v$ is an endpoint of $e$ we write $v\sim e$.

\begin{definition}[Directed Treewidth]\label{def:directedtreewidth}
	Let $D$ be a digraph.
	A \emph{directed tree decomposition} for $D$ is a tuple $\Brace{T,\beta,\gamma}$ where $T$ is an arborescence, $\beta\colon\Fkt{V}{T}\rightarrow 2^{V(D)}$ is a function that partitions $\Fkt{V}{D}$ into sets called the \emph{bags}\footnote{We do allow empty bags. Thus, this means $\CondSet{\Fkt{\beta}{t}}{t\in\V{T}}\setminus\{ \emptyset\}$ is a partition of $\V{D}$ into non-empty sets.}, and $\gamma\colon\Fkt{E}{T}\rightarrow 2^{V(D)}$ is a function, giving us sets called the \emph{guards}, satisfying the following requirement:
	\begin{enumerate}
		\item[] For every $\Brace{d,t}\in\Fkt{E}{T}$, $\Fkt{\gamma}{d,t}$ strongly guards $\Fkt{\beta}{T_t}\coloneqq\bigcup_{t'\in V(T_t)}\Fkt{\beta}{t'}$.
	\end{enumerate}
	Here $T_t$ denotes the subarboresence of $T$ with root $t$.
	For every $t\in\Fkt{V}{T}$ let $\Fkt{\Gamma}{t}\coloneqq\Fkt{\beta}{t}\cup\bigcup_{t\sim e}\Fkt{\gamma}{e}$.
	The \emph{width} of $\Brace{T,\beta,\gamma}$ is defined as
	\begin{align*}
		\Width{T,\beta,\gamma}\coloneqq\max_{t\in V(T)}\Abs{\Fkt{\Gamma}{t}}-1.
	\end{align*}
	The \emph{directed treewidth} of $D$, denoted by $\dtw{D}$, is the minimum width over all directed tree decompositions for $D$.
\end{definition}

\section{Statement of the structure theorems and overview of the proof}\label{sec:statementandoverview}

Overall our proof is structured much like similar proofs of structure theorems (see for example \cite{kawarabayashi2020quickly}), though the details are often quite distinct.
The bulk of our effort is dedicated to finding structure with respect to a given wall, which is often referred to as local structure.
This first part is the one that requires many new tools and the directed setting forces us to abandon known proof strategies in favour of alternative approaches.
Once we actually have local structure, the global structure theorem can largely be derived using existing tools for digraph structure theory and this part mercifully follows more along the lines of established approaches in graph structure theory.

\subsection{Towards local structure}

To find structure with respect to a wall, we first need to know what kind of wall we are searching for and how to find said wall.
Luckily both of these problems have already been figured out.

Let $k\in\N$ be a positive integer.
The \emph{cylindrical grid of order $k$} is the digraph obtained from the dicycles $C_1,\dots C_k$, with $C_i=\Brace{v_0^i,e^i_0,v_1^i,e_1^i,\dots,e_{2k-3}^i,v_{2k-2}^i,e_{2k-2}^i,v_{2k-1}^i,e_{2k-1}^i,v_0^i}$
for each $i\in[k]$, by adding the directed paths $P_i = v_{2i}^1v_{2i}^2\dots v_{2i}^{k-1}v_{2i}^k$ and $Q_i = v_{2i+1}^kv_{2i+1}^{k-1}\dots v_{2i+1}^2v_{2i+1}^1$ for every $i\in[0,k-1]$.

\begin{figure}[!h]
	\centering
	    \begin{tikzpicture}[scale=0.9]
		\pgfdeclarelayer{background}
		\pgfdeclarelayer{foreground}
		\pgfsetlayers{background,main,foreground}
		
		\node (mid) [v:ghost] {};
		
		\node (u11) [v:main,position=30:10mm from mid] {};
		\node (u21) [v:main,position=60:10mm from mid] {};
		\node (u31) [v:main,position=90:10mm from mid] {};
		\node (u41) [v:main,position=120:10mm from mid] {};
		\node (u51) [v:main,position=150:10mm from mid] {};
		\node (u61) [v:main,position=180:10mm from mid] {};
		\node (u71) [v:main,position=210:10mm from mid] {};
		\node (u81) [v:main,position=240:10mm from mid] {};
		\node (u91) [v:main,position=270:10mm from mid] {};
		\node (u101) [v:main,position=300:10mm from mid] {};
		\node (u111) [v:main,position=330:10mm from mid] {};
		\node (u121) [v:main,position=0:10mm from mid] {};
		
		\node (u12) [v:main,position=30:16mm from mid] {};
		\node (u22) [v:main,position=60:16mm from mid] {};
		\node (u32) [v:main,position=90:16mm from mid] {};
		\node (u42) [v:main,position=120:16mm from mid] {};
		\node (u52) [v:main,position=150:16mm from mid] {};
		\node (u62) [v:main,position=180:16mm from mid] {};
		\node (u72) [v:main,position=210:16mm from mid] {};
		\node (u82) [v:main,position=240:16mm from mid] {};
		\node (u92) [v:main,position=270:16mm from mid] {};
		\node (u102) [v:main,position=300:16mm from mid] {};
		\node (u112) [v:main,position=330:16mm from mid] {};
		\node (u122) [v:main,position=0:16mm from mid] {};
		
		\node (u13) [v:main,position=30:22mm from mid] {};
		\node (u23) [v:main,position=60:22mm from mid] {};
		\node (u33) [v:main,position=90:22mm from mid] {};
		\node (u43) [v:main,position=120:22mm from mid] {};
		\node (u53) [v:main,position=150:22mm from mid] {};
		\node (u63) [v:main,position=180:22mm from mid] {};
		\node (u73) [v:main,position=210:22mm from mid] {};
		\node (u83) [v:main,position=240:22mm from mid] {};
		\node (u93) [v:main,position=270:22mm from mid] {};
		\node (u103) [v:main,position=300:22mm from mid] {};
		\node (u113) [v:main,position=330:22mm from mid] {};
		\node (u123) [v:main,position=0:22mm from mid] {};
		
		\node (u14) [v:main,position=30:28mm from mid] {};
		\node (u24) [v:main,position=60:28mm from mid] {};
		\node (u34) [v:main,position=90:28mm from mid] {};
		\node (u44) [v:main,position=120:28mm from mid] {};
		\node (u54) [v:main,position=150:28mm from mid] {};
		\node (u64) [v:main,position=180:28mm from mid] {};
		\node (u74) [v:main,position=210:28mm from mid] {};
		\node (u84) [v:main,position=240:28mm from mid] {};
		\node (u94) [v:main,position=270:28mm from mid] {};
		\node (u104) [v:main,position=300:28mm from mid] {};
		\node (u114) [v:main,position=330:28mm from mid] {};
		\node (u124) [v:main,position=0:28mm from mid] {};
		
		\node (u15) [v:main,position=30:34mm from mid] {};
		\node (u25) [v:main,position=60:34mm from mid] {};
		\node (u35) [v:main,position=90:34mm from mid] {};
		\node (u45) [v:main,position=120:34mm from mid] {};
		\node (u55) [v:main,position=150:34mm from mid] {};
		\node (u65) [v:main,position=180:34mm from mid] {};
		\node (u75) [v:main,position=210:34mm from mid] {};
		\node (u85) [v:main,position=240:34mm from mid] {};
		\node (u95) [v:main,position=270:34mm from mid] {};
		\node (u105) [v:main,position=300:34mm from mid] {};
		\node (u115) [v:main,position=330:34mm from mid] {};
		\node (u125) [v:main,position=0:34mm from mid] {};
		
		\node (u16) [v:main,position=30:40mm from mid] {};
		\node (u26) [v:main,position=60:40mm from mid] {};
		\node (u36) [v:main,position=90:40mm from mid] {};
		\node (u46) [v:main,position=120:40mm from mid] {};
		\node (u56) [v:main,position=150:40mm from mid] {};
		\node (u66) [v:main,position=180:40mm from mid] {};
		\node (u76) [v:main,position=210:40mm from mid] {};
		\node (u86) [v:main,position=240:40mm from mid] {};
		\node (u96) [v:main,position=270:40mm from mid] {};
		\node (u106) [v:main,position=300:40mm from mid] {};
		\node (u116) [v:main,position=330:40mm from mid] {};
		\node (u126) [v:main,position=0:40mm from mid] {};
		
		\begin{pgfonlayer}{background}
			
			\draw[e:main,bend right=15,->] (u21) to (u31);
			\draw[e:main,bend right=15,->] (u41) to (u51);
			\draw[e:main,bend right=15,->] (u61) to (u71);
			\draw[e:main,bend right=15,->] (u81) to (u91);
			\draw[e:main,bend right=15,->] (u101) to (u111);
			\draw[e:main,bend right=15,->] (u121) to (u11);
			
			\draw[e:main,bend right=15,->] (u12) to (u22);
			\draw[e:main,bend right=15,->] (u32) to (u42);
			\draw[e:main,bend right=15,->] (u52) to (u62);
			\draw[e:main,bend right=15,->] (u72) to (u82);
			\draw[e:main,bend right=15,->] (u92) to (u102);
			\draw[e:main,bend right=15,->] (u112) to (u122);
			
			\draw[e:main,bend right=15,->] (u23) to (u33);
			\draw[e:main,bend right=15,->] (u43) to (u53);
			\draw[e:main,bend right=15,->] (u63) to (u73);
			\draw[e:main,bend right=15,->] (u83) to (u93);
			\draw[e:main,bend right=15,->] (u103) to (u113);
			\draw[e:main,bend right=15,->] (u123) to (u13);
			
			\draw[e:main,bend right=15,->] (u14) to (u24);
			\draw[e:main,bend right=15,->] (u34) to (u44);
			\draw[e:main,bend right=15,->] (u54) to (u64);
			\draw[e:main,bend right=15,->] (u74) to (u84);
			\draw[e:main,bend right=15,->] (u94) to (u104);
			\draw[e:main,bend right=15,->] (u114) to (u124);
			
			\draw[e:main,bend right=15,->] (u25) to (u35);
			\draw[e:main,bend right=15,->] (u45) to (u55);
			\draw[e:main,bend right=15,->] (u65) to (u75);
			\draw[e:main,bend right=15,->] (u85) to (u95);
			\draw[e:main,bend right=15,->] (u105) to (u115);
			\draw[e:main,bend right=15,->] (u125) to (u15);
			
			\draw[e:main,bend right=15,->] (u16) to (u26);
			\draw[e:main,bend right=15,->] (u36) to (u46);
			\draw[e:main,bend right=15,->] (u56) to (u66);
			\draw[e:main,bend right=15,->] (u76) to (u86);
			\draw[e:main,bend right=15,->] (u96) to (u106);
			\draw[e:main,bend right=15,->] (u116) to (u126);
			
			\draw[e:main,->] (u11) to (u12);
			\draw[e:main,->] (u12) to (u13);
			\draw[e:main,->] (u13) to (u14);
			\draw[e:main,->] (u14) to (u15);
			\draw[e:main,->] (u15) to (u16);
			
			\draw[e:main,->] (u26) to (u25);
			\draw[e:main,->] (u25) to (u24);
			\draw[e:main,->] (u24) to (u23);
			\draw[e:main,->] (u23) to (u22);
			\draw[e:main,->] (u22) to (u21);
			
			\draw[e:main,->] (u31) to (u32);
			\draw[e:main,->] (u32) to (u33);
			\draw[e:main,->] (u33) to (u34);
			\draw[e:main,->] (u34) to (u35);
			\draw[e:main,->] (u35) to (u36);
			
			\draw[e:main,->] (u46) to (u45);
			\draw[e:main,->] (u45) to (u44);
			\draw[e:main,->] (u44) to (u43);
			\draw[e:main,->] (u43) to (u42);
			\draw[e:main,->] (u42) to (u41);
			
			\draw[e:main,->] (u51) to (u52);
			\draw[e:main,->] (u52) to (u53);
			\draw[e:main,->] (u53) to (u54);
			\draw[e:main,->] (u54) to (u55);
			\draw[e:main,->] (u55) to (u56);
			
			\draw[e:main,->] (u66) to (u65);
			\draw[e:main,->] (u65) to (u64);
			\draw[e:main,->] (u64) to (u63);
			\draw[e:main,->] (u63) to (u62);
			\draw[e:main,->] (u62) to (u61);
			
			\draw[e:main,->] (u71) to (u72);
			\draw[e:main,->] (u72) to (u73);
			\draw[e:main,->] (u73) to (u74);
			\draw[e:main,->] (u74) to (u75);
			\draw[e:main,->] (u75) to (u76);
			
			\draw[e:main,->] (u86) to (u85);
			\draw[e:main,->] (u85) to (u84);
			\draw[e:main,->] (u84) to (u83);
			\draw[e:main,->] (u83) to (u82);
			\draw[e:main,->] (u82) to (u81);
			
			\draw[e:main,->] (u91) to (u92);
			\draw[e:main,->] (u92) to (u93);
			\draw[e:main,->] (u93) to (u94);
			\draw[e:main,->] (u94) to (u95);
			\draw[e:main,->] (u95) to (u96);
			
			\draw[e:main,->] (u106) to (u105);
			\draw[e:main,->] (u105) to (u104);
			\draw[e:main,->] (u104) to (u103);
			\draw[e:main,->] (u103) to (u102);
			\draw[e:main,->] (u102) to (u101);
			
			\draw[e:main,->] (u111) to (u112);
			\draw[e:main,->] (u112) to (u113);
			\draw[e:main,->] (u113) to (u114);
			\draw[e:main,->] (u114) to (u115);
			\draw[e:main,->] (u115) to (u116);
			
			\draw[e:main,->] (u126) to (u125);
			\draw[e:main,->] (u125) to (u124);
			\draw[e:main,->] (u124) to (u123);
			\draw[e:main,->] (u123) to (u122);
			\draw[e:main,->] (u122) to (u121);
			
		\end{pgfonlayer}
	\end{tikzpicture}
	\caption{The elementary directed $3$-wall.}
	\label{fig:directedwall1}
\end{figure}
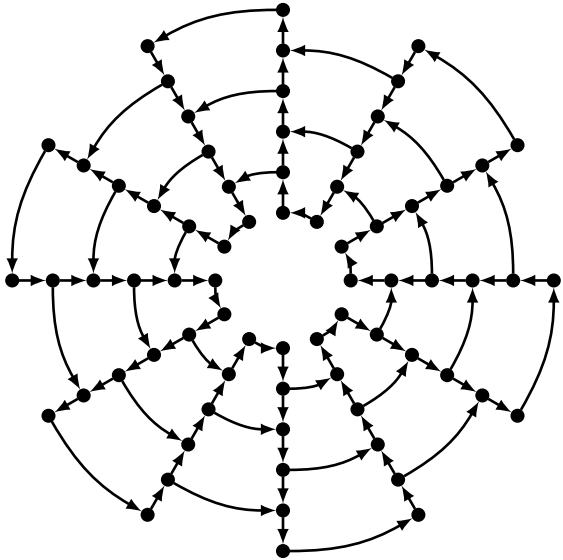

Let $k\in\N$ be a positive integer.
An \emph{elementary cylindrical $k$-wall} $W$ is the digraph obtained from the cylindrical grid $G$ of order $2k$ by deleting the edges $\Brace{v_{2i}^{2j},v_{2i+1}^{2j}}$ and $\Brace{v_{2i+1}^{2j+1},v_{2i+2}^{2j+1}}$ for every $i\in[0,2k-1]$ and every $j\in[0,k-1]$.
A \emph{cylindrical $k$-wall} is a subdivision of $W$.
See \cref{fig:directedwall1} for an illustration.

As in the undirected case, we can translate the Directed Grid Theorem \cite{KawarabayashiKreutzer2015DirectedGrid} into a Directed Wall Theorem by simply doubling the necessary quantities.
Since this is just a linear factor we may use the same function as the directed grid theorem and simply double the argument.

\begin{theorem}[Kawarabayashi and Kreutzer \cite{KawarabayashiKreutzer2015DirectedGrid}]\label{cor:directedwall}
	There exists a function $\DirectedGridNoArg \colon \N\rightarrow \N$ such that for every $k\in\N$ and every digraph $D$ we have $\dtw{D} \leq \DirectedGrid{2k}$, or $D$ contains a cylindrical $k$-wall.
\end{theorem}

The cylindrical wall plays a similar role to the wall in undirected graphs.
However, unlike walls in undirected graphs, which contain all planar graphs as minors, the cylindrical wall does not contain all planar graphs and not even all strongly planar graphs as butterfly minors.
In particular odd bicycles are not strongly planar and thus cannot be found in the cylindrical grid.

However, it turns out that if a digraph $D$ does not contain a large half-integral packing of even dicycles, then it either has small directed treewidth or contains a wall that can act as the starting point of an odd decomposition.
We will address the case of $D$ having small directed treewidth when we talk about the global structure theorem, but let us get more explicit about the second case.

\begin{definition}[Perimeter]\label{def:perimeter}
	Let $k\in\N$ be a positive integer, let $W$ be a $k$-wall, and let $Q_1, Q_k \subseteq W$ be the unique pair of distinct dicycles in $W$ whose deletion does not disconnect the underlying undirected graph of $W$. 
    We define the union of $Q_1$ and $Q_k$ to be the \emph{perimeter} of $W$ and write $\Perimeter{W}=Q_1\cup Q_k$.
\end{definition}

\begin{definition}[Odd wall]\label{def:oddwall}
Let $D$ be a digraph and $A \subseteq V(D)$ be a set of vertices.
Let $W \subseteq D-A$ be a cylindrical wall with the perimeter $Q_1 \cup Q_k$.
We say that $W$ is an \emph{odd wall under $A$} if 
\begin{enumerate}
	\item there is an undirected separation $(X,Y)$ of $D$ such that $X\cap Y=A\cup\Perimeter{W}$, $W \subseteq Y$, and every vertex in $Y$ reaches a vertex of $W-\Perimeter{W}$ or is reachable from it, and
	\item the cylindrical society $(\InducedSubgraph{D}{Y - A},\Omega_1,\Omega_2)$ with $\V{\Omega_1}=\V{Q_1}$ and $\V{\Omega_2}=\V{Q_k}$ has an odd rendition in the disk.
\end{enumerate} 
\end{definition}

This allows us to state the main theorem of \Cref{sec:oddwall}, which gives us the duality between either a large half-integral packing of even dicycles or a large odd wall that we mentioned above, provided the graph has high directed treewidth.

\begin{theorem}\label{thm:oddwall}
	There exist functions $\OddWallOrderNoArg \colon \N\times\N\rightarrow\N$ and $\OddWallApexNoArg \colon \N\rightarrow\N$ such that for all integers $r,t\geq 1$ and all digraphs $D$ there exists an integral packing of $t$ even dicycles in $D$ or for every $\OddWallOrder{r,t}$-wall $W$ in $D$ there exist a set $A\subseteq\V{D}$ with $\Abs{A}\leq \OddWallApex{t}$ and an $r$-wall $W'\subseteq W-A$ which is an odd wall under $A$.
\end{theorem}

Please note that the existence of an odd wall is a much stronger outcome than the notion of flatness used in \cite{Giannopoulou2020DirectedFlatWall}, as the notion used there does in fact allow for the existence of directed crosses within the flat wall\footnote{Although the occurrence of these crosses is highly controlled.}.
Nevertheless it is possible to derive \cref{thm:oddwall} directly using the techniques from \cite{Giannopoulou2020DirectedFlatWall}.
Since in some places the proofs need some adjustments propagating through many non-trivial steps of the arguments, we present an adjusted version of this proof in \Cref{sec:oddwall}.

We are now ready to state the local structure theorem.
Let $D$ be a digraph, let $W$ be a cylindrical wall in $D$, and let $A$ be a set of vertices in $D$.
If $W'$ is the union of all maximal walls within $W$ that do not contain vertices of $A$, then we call the strong component of $D - A$ that contains $W'$ the \emph{$W$-component of $D$ under $A$}.

\begin{theorem}[Local structure theorem]\label{thm:localstructureSec3}
    There exist two functions $\LocalStructureOutlineNoArg \colon \N \rightarrow \N$ and $\LocalStructureApexNoArg \colon \N \rightarrow \N$ such that for any integer $t$ and any digraph $D$ containing a cylindrical $\LocalStructureOutline{t}$-wall $W$ one of the following holds:
    \begin{itemize}
        \item $D$ contains a quarter-integral packing of $t$ even dicycles, or

        \item there exists a set of vertices $A \subseteq V(D)$, with $|A| \leq \LocalStructureApex{t}$, such that the $W$-component of $D$ under $A$ has a maelstrom-free, pure odd decomposition and is thus odd.
    \end{itemize}
\end{theorem}

Before moving on to the global structure theorem, we briefly discuss the other tools aside from the odd wall theorem that lead us to the proof of \Cref{thm:localstructure}.

Of particular importance are the contents of \Cref{sec:shifting}, where we start with an elementary proof of the relation of an odd and an even dicycle drawn in the plane.
It turns out that one can always shift the even dicycle to lie in the closure of one of the two disks bounded by the odd dicycle.
This is quite a powerful tool when trying to find a packing of even dicycles, as it allows us to move the even dicycles through the graph, provided we identified odd dicycles in strategic places.
In a very abstract sense, this will be the major tactic we use throughout our proof strategy.

However, an odd decomposition is not actually a planar drawing of a graph and this difference first of all forces us to prove a version of the planar shifting lemma we just described for largely grounded dicycles in an odd decomposition.
This is the one spot in our proof where we prominently use matching theory and despite the intuitive simplicity of the problem, the proof is quite arduous.
Once we are over that hill, we must contend with even dicycles that are not even planar in the sense of being grounded in an odd decomposition.
Intuitively, these even dicycles are the ones that occupy the maelstroms, though the contents of a maelstrom are also allowed to be planar as long as they contain an even dicycle.

For the purpose of shifting even dicycles to lie somewhere in the area surrounding a maelstrom, we use a powerful theorem by Giannopoulou and Wiederrecht \cite{Giannopoulou2023ExcludingSingleCrossing} that acts as a replacement for the two paths theorem, which does not hold in the directed setting in its straightforward translation.
At this point we will no longer be using a single odd dicycle to move our even dicycle, but instead a considerably large, but constant number of concentric, homogeneous odd dicycles.
Such a set of concentric odd dicycles can of course be found in an odd wall and forms the basis of the definition of a maelstrom outline.

The non-planar shifting lemmas also mark the first instance of us having to prove at least two different versions of a tool, because circle outlines, diamond outlines, and refined diamond outlines all behave slightly differently.
It also marks the first instance in which our proof turns from half-integral all the way to only guaranteeing a quarter-integral packing of even dicycles due to the structure of refined diamond outlines.

This section also introduces a more orderly substructure found within outlines, called \emph{rim}, whose existence we argue for in \Cref{sec:localising}.
We also want to highlight that \Cref{sec:shifting} contains the proof of \Cref{lem:detour}, a tool which allows us to reroute two linkages that intersect in an odd decomposition without loss.
This result will be used in almost every proof following it, up to the proof of \Cref{thm:localstructure}.
This highlights the difference to the undirected setting, where such a rerouting argument holds quite obviously, whilst here it takes a very specific set of circumstances and can only reroute our paths in a predetermined direction.

Once we reach \Cref{sec:transaction}, we will enter more familiar territory for structural graph theorists.
As in \cite{kawarabayashi2020quickly}, our intermediate goal before proving \Cref{thm:localstructure} will be to refine maelstroms and their associated outlines until they have bounded depth.
For this purpose, we need to show that given a huge transaction, we can either find a large quarter-integral packing of even dicycles, or we can integrate a large chunk of the given transaction into the already existing decomposition of our graph, whilst sacrificing a few vertices in the process.
The result itself is comparable to results in \cite{kawarabayashi2020quickly}, but our proof of this fact nonetheless takes some substantial detours compared to the methods in \cite{kawarabayashi2020quickly}.

With a large transaction that extends our decomposition in hand, the goal of \Cref{sec:buildoutline} will be to split the outline on whose society we found the transaction into more outlines, or find a large quarter-integral packing of even dicycles.
Due to the fact that we have three different types of outlines to work with and the complexities of the directed setting, this section is quite technical, but ultimately yields the expected result, again costing us a small set of vertices that we must delete to ensure that the new outlines and new maelstroms we find are actually separated.

Since the results of \Cref{sec:buildoutline} suggest that we are able to continuously split a maelstrom until we have found maelstrom outlines of bounded depth, we dedicate the rather short \Cref{sec:killmaelstroms} to the act of getting rid of the even dicycles in maelstroms whose outlines have bounded depth.
The proof in this section is inspired by the methods Thilikos and Wiederrecht present in \cite{thilikos2022killing}, hence the name of the section.
Though refined diamond outlines refuse to be dealt with entirely in this section, providing us with more work in the next section.

Our na{\"i}ve expectation that splitting maelstroms should be easy is then challenged in \Cref{sec:localstructure} once we try to put our tools to use to prove \Cref{thm:localstructure}, requiring us to first prove a technical statement on how splitting a maelstrom outline should actually function.
This forces us to show that, if we do not find more than one maelstrom hosting an even dicycle, we can iteratively refine the maelstrom further until one of several beneficial outcomes are reached.
The challenge of this part is that whilst we continue refining the outline, we must keep the apex set under a fixed bound, since we are not actually finding any even dicycles that we can separate.

After this is dealt with, we have to deal with an unfortunate edge case that occurs when killing a maelstrom with a refined diamond outline of bounded depth.
Since bounded depth for refined diamond outlines still leaves the possibility for certain large transactions to exist on eddy segments, we are forced to keep building our outline beyond the refined diamond outline and prove that this always yields a favourable outcome for us, whilst again keeping a fixed bounded on the apex set we incur in this process.
Once that is dealt with however, we finally show that \Cref{thm:localstructure} holds by combining the tools from the previous sections.

\subsection{Towards global structure}

To prove the global structure theorem we can mostly rely on existing results and tools for structural digraph theory from the literature.
The global structure theorem itself is then implied by a stronger version of it, which we can prove via induction in \Cref{sec:globalstructure}.
Much like the graph minor structure theorem for undirected graphs, the global structure of a digraph without a large quarter-integral packing of even dicycles takes the form of a (directed) tree decomposition.

\begin{theorem}\label{thm:globalstructureSec3}
    There exists a function $\GlobalStructureNoArg \colon \N \rightarrow \N$ such that for every digraph $D$ and every $k \in \N$ either
    \begin{enumerate}
        \item $D$ has a quarter-integral packing of $k$ even dicycles, or

        \item $D$ has a directed tree-decomposition $(T, \beta , \gamma )$ and a map $\alpha \colon V(T) \rightarrow 2^{V(D)}$ such that
        \begin{enumerate}
            \item $|\gamma(e)| \leq \GlobalStructure{k}$ for all $e \in E(T)$,

            \item $\alpha(t) \subseteq \Gamma(t)$ for all $t \in V(T)$,

            \item $|\alpha(t)| \leq \GlobalStructure{k}$ for all $t \in V(T)$, and

            \item no even dicycle of $D - \alpha(t)$ contains a vertex of $\Gamma(t) \setminus \alpha(t)$.
        \end{enumerate}
    \end{enumerate}
\end{theorem}

\subsection{Proof of the main theorem}

Using \Cref{thm:globalstructureSec3}, the proof of \Cref{thm:mainthm1} is a fairly standard instance of dynamic programming.
To illustrate this, and since we have not addressed digraphs with low directed treewidth yet, we give the analogous proof that graphs of low directed treewidth have the quarter-integral \EP-property for even dicycles before proving \Cref{thm:mainthm1}.

\begin{lemma}\label{lem:lowdtw}
    Let $k, t$ be integers and let $D$ be a digraph, then if $\dtw{D} < k$ there either exists an integral packing of $t$ even dicycles in $D$, or there exists a set of vertices $A \subseteq V(D)$ such that $D - A$ is odd and $|A| \leq k(t-1)$. 
\end{lemma}
\begin{proof}
    Since $D$ has bounded directed treewidth, there exists a directed tree decomposition $(T, \beta, \gamma)$ for $D$.
    (This decomposition can be found efficiently, as we will see in \Cref{thm:fptdirectedwall}, which we take from \cite{campos2022adapting}.)
    Choose $t \in V(T)$ such that $\bigcup_{s \in V(T_t)} \Gamma(s)$ contains an even dicycle and $T_t$ is minimal, in the sense that there exist no subarboresences within $T_t$ whose induced subgraph of $D$ also contains an even dicycle.
    Due to the minimality of $T_t$, $\bigcup_{s \in V(T_t)} \Gamma(s) - \Gamma(t)$ is odd and since $\dtw{D} \leq k$, we know that $|\Gamma(t)| \leq k$.
    
    We now iterate this procedure on $T - V(T_t)$ and the directed tree decomposition induced by this arboresence in $(T, \beta, \gamma)$.
    Clearly, if we can iterate this process $t$ times, we find an integral packing of $t$ even dicycles and otherwise, the union of all bags and guards we deleted contains at most $k (t-1)$ vertices.
\end{proof}

An alternative proof for \Cref{thm:mainthm1}, contrasting with the one we presented in the introduction, proceeds along the same lines as the proof of \Cref{lem:lowdtw}.

\begin{proof}[Proof of \Cref{thm:mainthm1} assuming \Cref{thm:globalstructureSec3}]
    Suppose that $D$ does not contain a quarter-integral packing of $k$ even dicycles and let $(T, \beta, \gamma)$ be the tree decomposition guaranteed by \Cref{thm:globalstructureSec3} together with the function $\alpha$.
    We again choose $V(T)$ such that $\bigcup_{s \in V(T_t)} \Gamma(s)$ contains an even dicycle and $T_t$ is minimal.
    Due to the minimality of $T_t$ and the fact that no even dicycle of $D - \alpha(t)$ contains a vertex of $\Gamma(t) \setminus \alpha(t)$ according to \Cref{thm:globalstructureSec3}, we know that $\bigcup_{s \in V(T_t)} \Gamma(s) - ( \alpha(t) \cup \gamma((d,t)) )$ is odd, where $(d,t) \in E(T) \setminus E(T_t)$.
    Thanks to \Cref{thm:globalstructureSec3}, we have $| \alpha(t) \cup \gamma((d,t)) | \leq 2\LocalStructureApex{t}$.

    We can then proceed as in the proof of \Cref{lem:lowdtw} either yielding a quarter-integral packing of $k$ even dicycles or an apex set of total size at most $(t-1) 2\LocalStructureApex{t}$, which proves the theorem.
\end{proof}

\section{Obtaining an odd wall}\label{sec:oddwall}

To prove \Cref{thm:oddwall} we will use a simple tool that shows that even small occurrences of non-planarity attached to the infrastructure offered by a cylindrical wall will create butterfly minor models of odd bicycles.
For a cylindrical grid $G$ of order $k$ obtained from the dicycles $C_1, \ldots , C_k$, we call a $C_1$-$C_k$- or $C_k$-$C_1$-path $P$ a \emph{perimeter jump} if it is internally disjoint from $G$.

\begin{lemma}\label{lem:gridplusjumpmakesoddbicycle}
Let $G$ be a cylindrical grid of order $k$, with $k \geq 3$, and let $P$ be a perimeter jump for $G$.
$G \cup P$ contains a weak odd bicycle and thus an even dicycle.
\end{lemma}
\begin{proof}
In fact, we can prove that given three integers $1 \leq i < j < \ell \leq k$, we can build a butterfly minor model of a weak odd bicycle using $C_1, C_i, C_j, C_{\ell}, C_k$ and only four paths from $P_1, \ldots , P_k, Q_1, \ldots , Q_k$.
W.l.o.g.\ we can assume that $P$ is a $C_1$-$C_k$-path.
Let $a \in \V{C_1}$ be the head of $P$ and let $x \in \V{C_k}$ be the tail of $P$.
Furthermore, let $s \in [k]$ and $b \in \V{Q_s}$ be chosen such that $bC_1a$ is of minimal length and let $t \in [k] \setminus \{ s \} $ and $y \in \V{Q_t}$  be chosen such that $xC_ky$ is of minimal length.
Let $b_s$ be the unique vertex in the intersection of $C_{\ell}$ and $Q_s$, and  let $y_t$ be the unique vertex in the intersection of $C_{\ell}$ and $Q_t$.
Let $r \in [k] \setminus \{ s,t \}$ and let $u_i, v_i, u_j, v_j \in \V{P_r \cup Q_r}$ be chosen such that $u_i$ and  $u_j$ are the unique vertices in $\V{P_r} \cap \V{C_i}$ and  $\V{P_r} \cap \V{C_i}$, respectively, and let $v_i$ and $v_j$ be defined analogously for $Q_r$.
We now define the following three dicycles:
\begin{align*}
    C_1' \coloneqq & \ v_iC_iu_i \cup u_iP_ru_j \cup u_jC_jv_j \cup v_jQ_rv_i , \qquad C_2' \coloneqq C_j , \text{ and} \\
    C_3' \coloneqq & \ y_tC_{\ell}b_s \cup b_sQ_sb \cup bC_1a \cup P \cup xC_ky \cup yQ_ty_t .
\end{align*}
It is easy to now verify that $C_1' \cup C_2' \cup C_3'$ form a butterfly minor model of the odd bicycle of order three.
\end{proof}

This observation will allow us to either find a large half integral packing of weak odd bicycles, and therefore the desired packing of even dicycles, or we will find, after removing a bounded size set of apex vertices, a large wall that has an odd cylindrical rendition.

\subsection{More definitions for the proof of \Cref{thm:oddwall}}

We will have to interact with the infrastructure of a wall quite closely and thus will need a bit more notation before we proceed.
Suppose $\mathcal{P}=\Set{P_1,\dots,P_k}$ is a family of pairwise disjoint directed paths, and $Q$ is a directed path that meets all paths in $\mathcal{P}$ such that $P_i\cap Q$ is a directed path.
\begin{itemize}
	\item We say that the paths $P_1,\dots,P_k$ \emph{appear in this order on $Q$} if for all $i\in[k-1]$, $P_i\cap Q$ occurs on $Q$ strictly before $P_{i+1}\cap Q$ with respect to the orientation of $Q$.
	
	\item In this case, for $i,j\in[k]$ with $i<j$, we denote by $\InducedSubgraph{Q}{P_i,\dots,P_j}$ the minimal directed subpath of $Q$ containing all vertices of $Q\cap P_{\ell}$ for $\ell\in[i,j]$.
\end{itemize}

\begin{figure}[!h]
	\centering
		\begin{tikzpicture}[scale=0.9]
		\pgfdeclarelayer{background}
		\pgfdeclarelayer{foreground}
		\pgfsetlayers{background,main,foreground}
		
		\node(v11) [v:main] {};
		\node(v12) [v:main,position=0:13mm from v11] {};
		\node(v13) [v:main,position=0:13mm from v12] {};
		\node(v14) [v:main,position=0:13mm from v13] {};
		\node(v15) [v:main,position=0:13mm from v14] {};
		\node(v16) [v:main,position=0:13mm from v15] {};
		\node(v17) [v:main,position=0:13mm from v16] {};
		\node(v18) [v:main,position=0:13mm from v17] {};
		
		\node(v01) [v:ghost,position=135:8mm from v11] {};
		\node(v03) [v:ghost,position=135:8mm from v13] {};
		\node(v05) [v:ghost,position=135:8mm from v15] {};
		\node(v07) [v:ghost,position=135:8mm from v17] {};
		
		\node(v21) [v:main,position=270:8mm from v11] {};
		\node(v22) [v:main,position=0:13mm from v21] {};
		\node(v23) [v:main,position=0:13mm from v22] {};
		\node(v24) [v:main,position=0:13mm from v23] {};
		\node(v25) [v:main,position=0:13mm from v24] {};
		\node(v26) [v:main,position=0:13mm from v25] {};
		\node(v27) [v:main,position=0:13mm from v26] {};
		\node(v28) [v:main,position=0:13mm from v27] {};
		
		\node(v31) [v:main,position=270:8mm from v21] {};
		\node(v32) [v:main,position=0:13mm from v31] {};
		\node(v33) [v:main,position=0:13mm from v32] {};
		\node(v34) [v:main,position=0:13mm from v33] {};
		\node(v35) [v:main,position=0:13mm from v34] {};
		\node(v36) [v:main,position=0:13mm from v35] {};
		\node(v37) [v:main,position=0:13mm from v36] {};
		\node(v38) [v:main,position=0:13mm from v37] {};
		
		\node(v41) [v:main,position=270:8mm from v31] {};
		\node(v42) [v:main,position=0:13mm from v41] {};
		\node(v43) [v:main,position=0:13mm from v42] {};
		\node(v44) [v:main,position=0:13mm from v43] {};
		\node(v45) [v:main,position=0:13mm from v44] {};
		\node(v46) [v:main,position=0:13mm from v45] {};
		\node(v47) [v:main,position=0:13mm from v46] {};
		\node(v48) [v:main,position=0:13mm from v47] {};
		
		\node(v51) [v:main,position=270:8mm from v41] {};
		\node(v52) [v:main,position=0:13mm from v51] {};
		\node(v53) [v:main,position=0:13mm from v52] {};
		\node(v54) [v:main,position=0:13mm from v53] {};
		\node(v55) [v:main,position=0:13mm from v54] {};
		\node(v56) [v:main,position=0:13mm from v55] {};
		\node(v57) [v:main,position=0:13mm from v56] {};
		\node(v58) [v:main,position=0:13mm from v57] {};
		
		\node(v61) [v:main,position=270:8mm from v51] {};
		\node(v62) [v:main,position=0:13mm from v61] {};
		\node(v63) [v:main,position=0:13mm from v62] {};
		\node(v64) [v:main,position=0:13mm from v63] {};
		\node(v65) [v:main,position=0:13mm from v64] {};
		\node(v66) [v:main,position=0:13mm from v65] {};
		\node(v67) [v:main,position=0:13mm from v66] {};
		\node(v68) [v:main,position=0:13mm from v67] {};
		
		\node(v71) [v:main,position=270:8mm from v61] {};
		\node(v72) [v:main,position=0:13mm from v71] {};
		\node(v73) [v:main,position=0:13mm from v72] {};
		\node(v74) [v:main,position=0:13mm from v73] {};
		\node(v75) [v:main,position=0:13mm from v74] {};
		\node(v76) [v:main,position=0:13mm from v75] {};
		\node(v77) [v:main,position=0:13mm from v76] {};
		\node(v78) [v:main,position=0:13mm from v77] {};
		
		\node(v81) [v:main,position=270:8mm from v71] {};
		\node(v82) [v:main,position=0:13mm from v81] {};
		\node(v83) [v:main,position=0:13mm from v82] {};
		\node(v84) [v:main,position=0:13mm from v83] {};
		\node(v85) [v:main,position=0:13mm from v84] {};
		\node(v86) [v:main,position=0:13mm from v85] {};
		\node(v87) [v:main,position=0:13mm from v86] {};
		\node(v88) [v:main,position=0:13mm from v87] {};
		
		\node(v91) [v:ghost,position=225:8mm from v81] {};
		\node(v93) [v:ghost,position=225:8mm from v83] {};
		\node(v95) [v:ghost,position=225:8mm from v85] {};
		\node(v97) [v:ghost,position=225:8mm from v87] {};
		
		\node(Q1) [v:ghost,position=90:4mm from v01] {$Q_1$};
		\node(Q2) [v:ghost,position=90:4mm from v03] {$Q_2$};
		\node(Q3) [v:ghost,position=90:4mm from v05] {$Q_3$};
		\node(Q4) [v:ghost,position=90:4mm from v07] {$Q_4$};
		
		\node(P1_1) [v:ghost,position=0:5mm from v18] {$P_1^1$};
		\node(P1_2) [v:ghost,position=0:5mm from v28] {$P_1^2$};
		
		\node(P2_1) [v:ghost,position=0:5mm from v38] {$P_2^1$};
		\node(P2_2) [v:ghost,position=0:5mm from v48] {$P_2^2$};
		
		\node(P3_1) [v:ghost,position=0:5mm from v58] {$P_3^1$};
		\node(P3_2) [v:ghost,position=0:5mm from v68] {$P_3^2$};
		
		\node(P4_1) [v:ghost,position=0:5mm from v78] {$P_4^1$};
		\node(P4_2) [v:ghost,position=0:5mm from v88] {$P_4^2$};
		
		\begin{pgfonlayer}{background}
			
			\draw[e:main,->,line width=1.9pt] (v11) to (v12);
			\draw[e:main,->] (v12) to (v13);
			\draw[e:main,->] (v13) to (v14);
			\draw[e:main,->] (v14) to (v15);
			\draw[e:main,->] (v15) to (v16);
			\draw[e:main,->] (v16) to (v17);
			\draw[e:main,->,line width=1.9pt] (v17) to (v18);
			
			\draw[e:main,<-,line width=1.9pt] (v21) to (v22);
			\draw[e:main,<-] (v22) to (v23);
			\draw[e:main,<-] (v23) to (v24);
			\draw[e:main,<-] (v24) to (v25);
			\draw[e:main,<-] (v25) to (v26);
			\draw[e:main,<-] (v26) to (v27);
			\draw[e:main,<-,line width=1.9pt] (v27) to (v28);
			
			\draw[e:main,->,line width=1.9pt] (v31) to (v32);
			\draw[e:main,->] (v32) to (v33);
			\draw[e:main,->] (v33) to (v34);
			\draw[e:main,->] (v34) to (v35);
			\draw[e:main,->] (v35) to (v36);
			\draw[e:main,->] (v36) to (v37);
			\draw[e:main,->,line width=1.9pt] (v37) to (v38);
			
			\draw[e:main,<-,line width=1.9pt] (v41) to (v42);
			\draw[e:main,<-] (v42) to (v43);
			\draw[e:main,<-] (v43) to (v44);
			\draw[e:main,<-] (v44) to (v45);
			\draw[e:main,<-] (v45) to (v46);
			\draw[e:main,<-] (v46) to (v47);
			\draw[e:main,<-,line width=1.9pt] (v47) to (v48);
			
			\draw[e:main,->,line width=1.9pt] (v51) to (v52);
			\draw[e:main,->] (v52) to (v53);
			\draw[e:main,->] (v53) to (v54);
			\draw[e:main,->] (v54) to (v55);
			\draw[e:main,->] (v55) to (v56);
			\draw[e:main,->] (v56) to (v57);
			\draw[e:main,->,line width=1.9pt] (v57) to (v58);
			
			\draw[e:main,<-,line width=1.9pt] (v61) to (v62);
			\draw[e:main,<-] (v62) to (v63);
			\draw[e:main,<-] (v63) to (v64);
			\draw[e:main,<-] (v64) to (v65);
			\draw[e:main,<-] (v65) to (v66);
			\draw[e:main,<-] (v66) to (v67);
			\draw[e:main,<-,line width=1.9pt] (v67) to (v68);
			
			\draw[e:main,->,line width=1.9pt] (v71) to (v72);
			\draw[e:main,->] (v72) to (v73);
			\draw[e:main,->] (v73) to (v74);
			\draw[e:main,->] (v74) to (v75);
			\draw[e:main,->] (v75) to (v76);
			\draw[e:main,->] (v76) to (v77);
			\draw[e:main,->,line width=1.9pt] (v77) to (v78);
			
			\draw[e:main,<-,line width=1.9pt] (v81) to (v82);
			\draw[e:main,<-] (v82) to (v83);
			\draw[e:main,<-] (v83) to (v84);
			\draw[e:main,<-] (v84) to (v85);
			\draw[e:main,<-] (v85) to (v86);
			\draw[e:main,<-] (v86) to (v87);
			\draw[e:main,<-,line width=1.9pt] (v87) to (v88);
			
			\draw[e:main,->,line width=1.9pt] (v12) to (v22);
			\draw[e:main,->] (v14) to (v24);
			\draw[e:main,->] (v16) to (v26);
			\draw[e:main,->,line width=1.9pt] (v18) to (v28);
			
			\draw[e:main,->,line width=1.9pt] (v21) to (v31);
			\draw[e:main,->] (v23) to (v33);
			\draw[e:main,->] (v25) to (v35);
			\draw[e:main,->,line width=1.9pt] (v27) to (v37);
			
			\draw[e:main,->,line width=1.9pt] (v32) to (v42);
			\draw[e:main,->] (v34) to (v44);
			\draw[e:main,->] (v36) to (v46);
			\draw[e:main,->,line width=1.9pt] (v38) to (v48);
			
			\draw[e:main,->,line width=1.9pt] (v41) to (v51);
			\draw[e:main,->] (v43) to (v53);
			\draw[e:main,->] (v45) to (v55);
			\draw[e:main,->,line width=1.9pt] (v47) to (v57);
			
			\draw[e:main,->,line width=1.9pt] (v52) to (v62);
			\draw[e:main,->] (v54) to (v64);
			\draw[e:main,->] (v56) to (v66);
			\draw[e:main,->,line width=1.9pt] (v58) to (v68);
			
			\draw[e:main,->,line width=1.9pt] (v61) to (v71);
			\draw[e:main,->] (v63) to (v73);
			\draw[e:main,->] (v65) to (v75);
			\draw[e:main,->,line width=1.9pt] (v67) to (v77);
			
			\draw[e:main,->,line width=1.9pt] (v72) to (v82);
			\draw[e:main,->] (v74) to (v84);
			\draw[e:main,->] (v76) to (v86);
			\draw[e:main,->,line width=1.9pt] (v78) to (v88);
			
			\draw[e:main,bend left=30,->,line width=1.9pt] (v01) to (v11);
			\draw[e:main,bend left=30,->] (v03) to (v13);
			\draw[e:main,bend left=30,->] (v05) to (v15);
			\draw[e:main,bend left=30,->,line width=1.9pt] (v07) to (v17);
			
			\draw[e:main,bend left=30,->,line width=1.9pt] (v81) to (v91);
			\draw[e:main,bend left=30,->] (v83) to (v93);
			\draw[e:main,bend left=30,->] (v85) to (v95);
			\draw[e:main,bend left=30,->,line width=1.9pt] (v87) to (v97);
			
		\end{pgfonlayer}
	\end{tikzpicture}
	\caption{The elementary cylindrical $4$-wall. The thick edges of the cycles $Q_1$ and $Q_4$ mark its perimeter.}
	\label{fig:cylindrical4wall}
\end{figure}
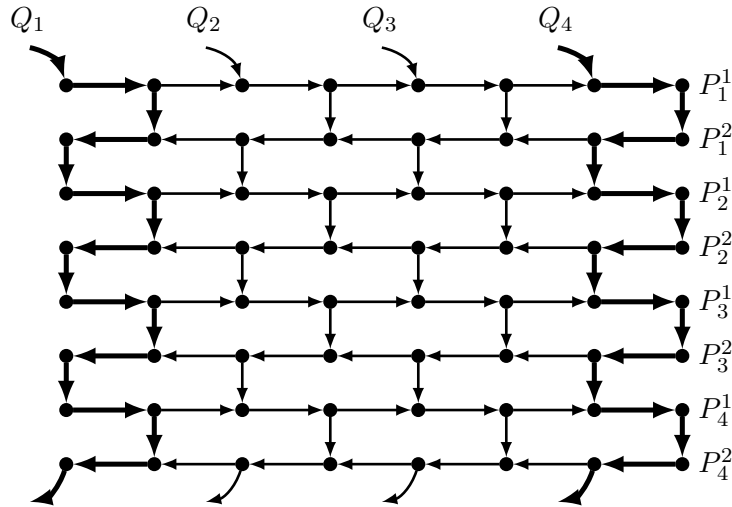

\begin{definition}[Vertical and horizontal paths]
	Let $k\in\N$ be a positive integer and $W$ be a cylindrical $k$-wall.
	
	We denote the vertical paths of $W$ by $Q_1,\dots,Q_k$, ordered from left to right.
	Let $\CondSet{P_j^i}{i\in[2],~j\in[k]}$ be the horizontal directed paths such that the paths $P^1_j$, $j\in[k]$, are oriented from left to right and the paths $P^2_j$, $j\in[k]$, are oriented from right to left such that $P^i_j$ is above $P^{i'}_{j'}$ whenever $j<j'$ and $P^1_j$ is above $P^2_j$ for all $j\in[k]$.
	The top line is $P_1^1$.
	
	By $\hat{P}_j$ we denote the disjoint union of $P_j^1$ and $P_j^2$ for all $j\in[k]$.
	
	Two horizontal paths $P_j^i$ and $P_{j'}^{i'}$ are \emph{consecutive} if $i\neq i'$, and $j'\in[j-1,j+1]$ or if $P_j^i+P_{j'}^{i'}=\hat{P}_j$.
	A family $\mathcal{P}\subseteq \CondSet{P_j^i}{i\in[2],~j\in[k]}$ is said to be \emph{consecutive} if there do not exist paths $P_1,P_2\in\mathcal{P}$, and $P_3\in\CondSet{P_j^i}{i\in[2],~j\in[k]}\setminus\mathcal{P}$ such that there is no directed path from $P_1$ to $P_2$ in $W-P_3$.
	We extend our notation for $\InducedSubgraph{P_j^i}{Q_p,\dots,Q_q}$ for $p<q$ in the natural way for $\InducedSubgraph{\mathcal{P}}{Q_p,\dots,Q_q}$ and, in a slight abuse of notation, identify $\hat{P}_i$ and $\Set{P^1_i,P^2_i}$.
	
	For more convenience we write ``Let $W=\Brace{Q_1,\dots,Q_k,\hat{P}_1,\dots,\hat{P}_k}$ be a cylindrical $k$-wall.'' to fix the embedding and naming of the vertical cycles and horizontal paths as explained above and depicted in \cref{fig:cylindrical4wall} for the case $k=4$.
\end{definition}

It is convenient for us to imagine cylindrical walls as illustrated in \cref{fig:cylindrical4wall}.
So we think of the concentric cycles as paths from the top to the bottom with an additional edge from the lowest row back to the top.
The directed paths that alternately go in and out of the wall are then seen as the horizontal paths.
Since the underlying undirected graph is planar and $3$-connected Whitney's Theorem \cite{whitney1992congruent} ensures a unique planar embedding.
In particular, we will refer to the faces of the embedding as the faces of the graph itself.
A face that is bounded by a cycle that is not the perimeter is called a \emph{brick} of the wall.

In what follows we introduce the remaining definitions necessary for our proof.

\begin{definition}[$W$-distance]\label{def:Wdistance}
	Let $k\in\N$ be a positive integer and $W=\Brace{Q_1,\dots,Q_k,\hat{P}_1,\dots,\hat{P}_k}$ be a cylindrical $k$-wall.
	Given two vertices $u,v\in\V{W}$, we say that they have \emph{$W$-distance} at least $i$ if there exist $i$ distinct vertical or $i$ distinct horizontal paths whose removal separates $u$ and $v$ in $W$.
\end{definition}

\begin{definition}[Slice]\label{def:slice}
	Let $k\in\N$ be a positive integer and $W=\Brace{Q_1,\dots,Q_k,\hat{P}_1,\dots,\hat{P}_k}$ be a cylindrical $k$-wall.
	A \emph{slice} $W'$ of $W$ is a cylindrical wall containing the vertical paths $Q_i,\dots,Q_{i+\ell}$ for some $i\in[k]$ and some $\ell\in[k-i]$, and the horizontal paths $\InducedSubgraph{P_1^1}{Q_i,\dots,Q_{i+\ell}},\dots,\InducedSubgraph{P_k^2}{Q_i,\dots,Q_{i+\ell}}$.
	We say that $W'$ is the \emph{slice of $W$ between $Q_i$ and $Q_{i+\ell}$} and that it is of \emph{width $\ell+1$}.
\end{definition}

\begin{definition}[Strip]\label{def:strip}
	Let $k\in\N$ be a positive integer and $W=\Brace{Q_1,\dots,Q_k,\hat{P}_1,\dots,\hat{P}_k}$ be a cylindrical $k$-wall.
	A \emph{strip of height $j-i+1$ between $i$ and $j$} of $W$ is the subgraph of $W$ induced by the horizontal paths $\hat{P}_i,\dots,\hat{P}_j$ for some $i<j\in[k]$ and the subpaths $\InducedSubgraph{Q_{\ell}}{\hat{P}_i,\dots,\hat{P}_j}$ for $\ell\in[k]$.
\end{definition}

\begin{figure}[!h]
	\centering
	\input{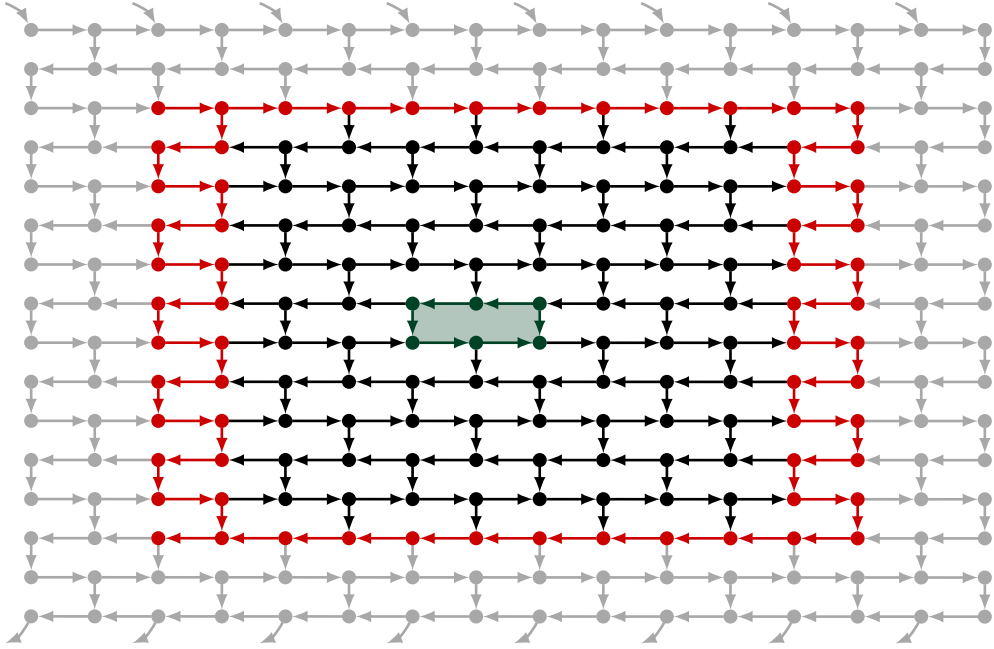}
	\caption{The tile $T$ at $\Brace{2,2}$ of width $2$ in a cylindrical $8$-wall.
		The green brick is the centre, the red paths are the perimeter of $T$, and the vertices that belong to $T$, but neither to the red paths, nor the green brick, are the internal vertices of $T$.}
	\label{fig:atile}
\end{figure}

\begin{definition}[Tiles]\label{def:tile}
	Let $k\in\N$ be a positive integer and $W=\Brace{Q_1,\dots,Q_k,\hat{P}_1,\dots,\hat{P}_k}$ be a cylindrical $k$-wall.
	Let $i,j\in[k]$ and $d\in\N$ be positive integers.
	The \emph{tile $T$ of $W$ at $\Brace{i,j}$ of width $d$} is defined as the subgraph of $W$ induced by
	\begin{align*}
		\bigcup_{\ell\in[i,i+2d+1]}\InducedSubgraph{Q_{\ell}}{\hat{P}_j,\dots,\hat{P}_{j+2d+1}}\cup\bigcup_{\ell\in[j,j+2d+1]}\InducedSubgraph{\hat{P}_{\ell}}{Q_i,\dots,Q_{i+2d+1}}.
	\end{align*}
	We call $i$ the \emph{column index} of the tile, $j$ the \emph{row index} of the tile, and say that the $j$-th row has a tile $T$ if the row index of $T$ is $j$.
	To make the notation a bit more compact we write $T_{i,j,d}$ for the tile of $W$ at $\Brace{i,j}$ of width $w$.
	
	Since $W$ is a cylindrical wall, there exist subgraphs of $W$ that technically also form tiles, but that do not necessarily fit into our \emph{parametrisation} of $W$.
	To overcome this, we agree for $\ell>k$ to set $P_{\ell}\coloneqq P_{\Brace{\Brace{\ell-1}\mod k}+1}$.
	This means that tiles that start near the bottom are allowed to continue at the top.
	Indeed, the notions of top and bottom are only present because of the way we parametrised the wall, and thus even those tiles are well defined.
	
	The \emph{perimeter} of the tile $T$ is $T\cap\Brace{Q_i\cup Q_{i+2d+1}\cup P_j^1\cup P_{j+2d+1}^2}$.
	We call $Q_i$ the \emph{left path} of the perimeter, $Q_{i+2d+1}$ its \emph{right path}, $P_j^1$ the \emph{upper path} of the perimeter, and finally $P_{j+2d+1}^2$ is its \emph{lower path}.
	
	The \emph{corners} of a tile are the vertices $a,b,c,d\in\V{T}$ where $a$, the \emph{upper left corner}, is the common starting point of $T\cap Q_i$ and $T\cap P_j^1$, $b$, the \emph{upper right corner}, is the end of $T\cap P_j^1$ and the starting point of $Q_{i+2d+1}$, $c$, the \emph{lower left corner}, is the common end of $T\cap Q_i$ and $T\cap P_{j+2d+1}^2$, finally $d$, the \emph{lower right corner}, is the end of $T\cap Q_{i+2d+1}$ and the starting point of $T\cap P_{j+2d+1}^2$.
	
	The \emph{centre} of $T$ is the boundary of the unique brick $C_T$ of $W$ whose boundary consists of vertices from $Q_{i+d+1}$, $Q_{i+d+2}$, $P^2_{j+d+1}$, and $P^1_{j+d+2}$.
	All vertices of $T$ which are not in the centre and not on the perimeter of $T$ are called \emph{internal}.
	See \cref{fig:atile} for an illustration of a tile.
\end{definition}

Please note that by this definition, only bricks that lie between $P_i^1$ and $P_i^2$ for some $i\in[k]$ can be centre of a tile.
However, if we were to take the mirror image of our currently fixed embedding along a straight vertical line between $Q_{\Floor{\frac{k}{2}}}$ and $Q_{\Floor{\frac{k}{2}}+1}$, we obtain a new embedding for which we then can reapply our parametrisation.
By doing so, every path $P_i^1$ now becomes a path $P_{i'}^2$, and $P_i^2$ becomes $P_{i''}^1$ for $i,i',i''\in[k]$.
This means that we can define for every brick $F$ of $W$ a tile $T_F$ such that $F$ is the centre of $T_F$.

The notion of \hyperref[def:tile]{tiles} allows us to add another layer of parametrisation on top of a cylindrical wall.
In some sense a tile can be seen as a generalisation of a brick that also contains a small acyclic wall inside to allow for additional routing.
The next few definitions are used to add further details to how our walls are divided into different regions and how tiles are used to achieve this additional layer of parametrisation.

\begin{definition}[Triadic partitions]\label{def:triadicpartition}
	Let $k\in\N$ be a positive integer and $W=\Brace{Q_1,\dots,Q_{3k},\hat{P}_1,\dots,\hat{P}_{3k}}$ be a cylindrical $3k$-wall.
	The \emph{triadic partition} of $W$ is the tuple
	\begin{align*}
		\Triade=\Brace{W,k,W_1,W_2,W_3,W^1,W^2,W^3}
	\end{align*}
	such that for each $i\in[3]$, $W_i$ denotes the \hyperref[def:slice]{slice} of $W$ between $Q_{k\Brace{i-1}+1}$ and $Q_{ik}$, and $W^i$ denotes the strip of $W$ between the rows $k\Brace{i-1}+1$ and $ik$.
\end{definition}

\begin{definition}[Tiling]
	Let $k\in\N$ be a positive integer and $W=\Brace{Q_1,\dots,Q_{3k},\hat{P}_1,\dots,\hat{P}_{3k}}$ be a cylindrical $3k$-wall with its \hyperref[def:triadicpartition]{triadic partition} $\Triade=\Brace{W,k,W_1,W_2,W_3,W^1,W^2,W^3}$.
	
	A \emph{tiling} is a family of pairwise disjoint \hyperref[def:tile]{tiles} $\Tiling$.
	Let $W'$ be a \hyperref[def:slice]{slice} of $W$, we say that $\Tiling$ \emph{covers} $W'$ if every vertex $v\in\V{W'}$ with $\DegG{W}{v}=3$ is contained in some tile of $\Tiling$.
	The tiling $\Tiling$ is said to \emph{cover} $\Triade$ if it covers $W_2$ and the \emph{width} of a tiling is the minimum over the width of its tiles.
	
	In most places we will use the following family of tilings:
	Let $f_w\colon\N\rightarrow\N$ be some function, $t\in\N$ a positive integer, and $\xi,\xi'\in[\Fkt{f_w}{t}+1]$.
	We define the \emph{column function} $\ColumnFunction\colon\N\rightarrow\N$ and the \emph{row function} $\RowFunction\colon\N\rightarrow\N$ as follows\footnote{Please note that $\ColumnFunction$ and $\RowFunction$ do also depend on $f_w$, $\xi$, $\xi'$, and $k$. However, it is more convenient to make these dependencies implicit.}:
	\begin{align*}
		\Fkt{\ColumnFunction}{p}&\coloneqq \Brace{k+2-\xi}+\Brace{p-1}\Brace{2\Fkt{f_w}{t}+1}\text{, and}\\
		\Fkt{\RowFunction}{q}&\coloneqq \xi'+\Brace{q-1}\Brace{2\Fkt{f_w}{t}+1}.
	\end{align*}
	We can now define our standard tiling for fixed $f_w$, $\xi$, and $\xi'$.
	\begin{align*}
		\Tiling_{W,k,\Fkt{f_w}{t},\xi,\xi'}\coloneqq\left\{T_{\Fkt{\ColumnFunction}{p},\Fkt{\RowFunction}{q},\Fkt{f_w}{t}} \mid p\in\left[\Ceil{\frac{k+\xi-1}{2\Fkt{f_w}{t}+1}}+1\right]\text{, }q\in\left[\Ceil{\frac{3k-\xi'-1}{2\Fkt{f_w}{t}+1}}+1\right]\right\}
	\end{align*}
\end{definition}

Note that every tiling $\Tiling_{W,k,\Fkt{f_w}{t}.\xi,\xi'}$ covers $\Triade$.
Moreover, every brick of $W_2$ that lies between the two paths of $\hat{P}_i$ for some $i\in[3k]$ is the centre of some tile $T'$ of some tiling $\Tiling'\in\Tiling_{W,k,\Fkt{f_w}{t}.\xi,\xi'}$.
Hence if we perform the mirror-image operation as described after the definition of tiles, we are able to find in total $2\Brace{\Fkt{f_w}{t}+1}^2$ many tilings that cover $W_2$, such that every brick of $W_2$ is the centre of some tile in one of these tilings.

We will commonly want to partition our tilings in a certain way such that tiles from different partition classes are not neighbouring in our wall.

\begin{definition}[Colouring of a tiling]
    Let $f_w\colon\N\rightarrow\N$ be some function, $t,k\in\N$ two positive integers, and $\xi,\xi'\in[\Fkt{f_w}{t}+1]$.
    Moreover, let $W=\Brace{Q_1,\dots,Q_{3k},\hat{P}_1,\dots,\hat{P}_{3k}}$ be a cylindrical $3k$-wall with its \hyperref[def:triadicpartition]{triadic partition} $\Triade=\Brace{W,k,W_1,W_2,W_3,W^1,W^2,W^3}$ and \hyperref[def:tiling]{$\Tiling=\Tiling_{W,k,\Fkt{f_w}{t},\xi,\xi'}$}.
    A \emph{colouring} of $\Tiling$ is a partition of $\Tiling$ into four classes, namely $\Class_1,\Class_2,\Class_3$, and $\Class_4$ as follows:
    For every $i\in\left[\Ceil{\frac{k+\xi-1}{2\Fkt{f_w}{t}+1}}+1\right]$ and every $j\in \left[\Ceil{\frac{3k-\xi'-1}{2\Fkt{f_w}{t}+1}}+1\right]$ we assign to $T_{\Fkt{\ColumnFunction}{i},\Fkt{\RowFunction}{j},\Fkt{f_w}{t}}$ the colour ${ \Brace{i\mod 2}+2\Brace{j\mod 2}+1 }$.
\end{definition}
    
This means that tiles where $\Fkt{\ColumnFunction}{i}$ and $\Fkt{\RowFunction}{j}$ are even get colour $1$, the tiles where $\Fkt{\RowFunction}{j}$ is even but $\Fkt{\ColumnFunction}{i}$ is odd get $3$, and so on.
Hence every column is two-chromatic, every row is so as well, and between each pair of tiles from the same colour that share a row or a column, there is a tile of a different colour that separates those tiles in their respective row or column.
Additionally, if $T$ is some tile, then the eight tiles surrounding $T$ are all of different colour than $T$ itself.

We will use tilings and colourings in several different ways, and sometimes it is necessary to `zoom out' of our current wall, i.e.\ to forget about some of the horizontal paths and vertical cycles in order to obtain a more streamlined version of our wall.

\begin{definition}[Walls from a tiling]\label{def:tiling}
	Let $k,d\in\N$ be positive integers, $W=\Brace{Q_1,\dots,Q_{3k},\hat{P}_1,\dots,\hat{P}_{3k}}$ be a cylindrical $3k$-wall with its triadic partition $\Triade=\Brace{W,k,W_1,W_2,W_3,W^1,W^2,W^3}$, and $\Tiling$ be a tiling of width $d$ that covers $W_2$.
	Moreover, let $\widetilde{W}$ be some slice of $W_2$ and let $I_Q$ be the largest set of integers such that $Q_i$ contains vertices of a tile from $\Tiling$ which intersects $\widetilde{W}$ for every $i\in I_Q$.
	Let $\Fkt{\widetilde{W}}{\Tiling}$ be the union of the cycles $Q_i$, $i\in I_Q$, and the paths $\InducedSubgraph{P_j^i}{\CondSet{Q_h}{h\in I_Q}}$ for every $\Brace{j,i}\in [3k]\times[2]$.
	We call $\Fkt{\widetilde{W}}{\Tiling}$ the \emph{extension of $\widetilde{W}$ that covers $\Tiling$}.
	
	Now, let $\Set{\Class_1,\dots,\Class_4}$ be a four colouring of $\Tiling$ and let $i\in[4]$ be a fixed colour.
	Then let $J_Q\subseteq[3k]$ be the largest set of integers such that for all $j\in J_Q$ the vertical cycle $Q_j$ of $\Fkt{\widetilde{W}}{\Tiling}$ does not contain a vertex of the interior of some tile from $\Class_i$.
	Similarly, let $J_P\subseteq[3k]$ be the largest set of integers such that for every $j\in J_P$, none of the two paths from $\hat{P}_j$ contains a vertex of the interior of a tile from $\Class_i$.
	
	By $\InducedSubgraph{\widetilde{W}}{\Tiling,i}$ we denote the subgraph of $W$ induced by the union of the cycles $Q_\ell$, $\ell \in J_Q$, and the paths $\InducedSubgraph{P_k^q}{\CondSet{Q_h}{h\in J_Q}}$ for every $j \in J_P$ and $q \in [2]$.
	We call $\InducedSubgraph{\widetilde{W}}{\Tiling,i}$ the \emph{$i$th $\Tiling$-slice of $\widetilde{W}$}.
\end{definition}

Note that in $\InducedSubgraph{\widetilde{W}}{\Tiling,i}$ we essentially cut out the tiles of $\Class_i$.
This operation gives us a slice $W'$ of some cylindrical wall for which the perimeter of every tile in $\Class_i$ is a brick in $W'$.
Next, we are going to find a tiling of $W'$ such that every tile of $\Class_i$ that belongs to $\widetilde{W}$ is captured by the centre of some tile in the new tiling.

\begin{definition}[Tier II tiling]\label{def:tilingII}
	Let $t,k,k'\in\N$ be positive integers and $f\colon\N\rightarrow\N$ be some function where $k\geq k'$.
	Let $W$ be a cylindrical $3k$-wall with its \hyperref[def:triadicpartition]{triadic partition} $\Triade=\Brace{W,k,W_1,W_2,W_3,W^1,W^2,W^3}$, and \hyperref[def:tiling]{$\Tiling=\Tiling_{W,k,f,\xi,\xi'}$} for some $\xi,\xi'\in[\Fkt{f}{t}+1]$, as well as $\Set{\Class_1,\dots,\Class_4}$ be a four colouring of $\Tiling$, and $i\in[4]$ be a fixed colour.
	Moreover, let $\widetilde{W}$ be a \hyperref[def:slice]{slice} of $W_2$ of width $k'$ such that no \hyperref[def:tile]{tile} of $\Class_i$ contains a vertex of $\Perimeter{\widetilde{W}}$ and $\widetilde{\Tiling}$ be the collection of all tiles from $\Tiling$ that contain a vertex of $\widetilde{W}$.
	
	The \emph{tier II tiling for $\widetilde{W}$ and $i$ obtained from $\Tiling$} is defined as the unique tiling $\TierIITiling{\Tiling,i,f}{\widetilde{W}}$ of width $f(t)$ of $\InducedSubgraph{\widetilde{W}}{\Tiling,i}$ such that every $T\in\Class_i\cap\widetilde{\Tiling}$ is in the interior of the centre of some tile in $\TierIITiling{\Tiling,i,f}{\widetilde{W}}$.
\end{definition}

Since every tile in $\Tiling$ consists of $2\Fkt{f}{t}+2$ path pairs, $\TierIITiling{\Tiling,i,f}{\widetilde{W}}$ is well defined and does in fact cover all of $\InducedSubgraph{\widetilde{W}}{\Tiling,i,f}$.

Let $W=(Q_1,\dots,Q_k,\hat{P}_1,\dots,\hat{P}_k)$ be a cylindrical $k$-wall.
We say that $W$ \emph{grasps} a butterfly minor model $\mu$ of $\Bidirected{K_t}$ if for every $v\in\V{\Bidirected{K_t}}$ there exist a pair $i_v,j_v\in[k]$ together with some $h_v\in[2]$ such that $\V{Q_{i_v}}\cap\V{P_{j_v}^{h_v}}\subseteq \V{\Fkt{\mu}{v}}$.  

\subsection{Three results from the directed Flat Wall Theorem}

In this section we present the partial results from \cite{Giannopoulou2020DirectedFlatWall} which make up the corner stones of our proof of \cref{thm:oddwall}.


\begin{theorem}[Giannopoulou, Kawarabayashi, Kreutzer, and Kwon \cite{Giannopoulou2020DirectedFlatWall}]\label{thm:xpaths}
	Let $D$ be a digraph and $X\subseteq\V{D}$.
	For all positive $k\in\N$, there are $k$ pairwise vertex disjoint directed $X$-paths in $D$, or there exists a set $S\subseteq\V{D}$ of size at most $2k$ such that every directed $X$-path in $D$ contains a vertex of $S$.
	
	Furthermore, there is a polynomial time algorithm which, given $D$ and $X\subseteq\V{D}$ as input, outputs $k$ pairwise disjoint directed $X$-paths, or a set $S\subseteq\V{D}$ of size at most $2k$ as above.
\end{theorem}


\begin{theorem}[Giannopoulou, Kawarabayashi, Kreutzer, and Kwon \cite{Giannopoulou2020DirectedFlatWall}]\label{thm:halfintegral}
	Let $k\in\N$ be a positive integer, $D$ be a digraph, and $X,Y\subseteq\V{D}$.
	If $\mathcal{P}$ is a half-integral $X$-$Y$-linkage of order $2k$ in $D$, then there exists a family $\mathcal{J}$ of $k$ pairwise disjoint $X$-$Y$-paths such that $\V{\mathcal{J}}\subseteq\V{\mathcal{P}}$.
\end{theorem}

Let $k,w\in\N$ be positive integers, $W$ be a cylindrical $k$-wall and $W'$ be a slice of $W$.
A directed $\V{W'}$-path $P$ is called a \emph{jump over $W'$} if $\E{P}\cap \E{W'}=\emptyset$.
We say that a directed $\V{W'}$-path $P$ is a \emph{$w$-long jump over $W'$} if for all $\xi,\xi'\in[w+1]$ the endpoints of $P$ belong to distinct tiles $T_1$ and $T_2$ of the tiling $\Tiling_{W,k,w,\xi,\xi'}$.
If $w$ is clear from the context, we simply call be a $P$ \emph{long jump}.

The following is a combination of Lemmas 4.3 to 4.8 from \cite{Giannopoulou2020DirectedFlatWall} and a proof can be found in the proof of Lemma 4.9 in \cite{Giannopoulou2020DirectedFlatWall}.
The only difference of Lemma 4.9 from \cite{Giannopoulou2020DirectedFlatWall} and the statement below is, that we extract the last subcase of Case 1 in its proof as a potential outcome.

\begin{lemma}[Giannopoulou, Kawarabayashi, Kreutzer, and Kwon \cite{Giannopoulou2020DirectedFlatWall}]\label{lemma:directedlongjumps}
	There exist functions $f_w\colon\N\rightarrow\N$, $f_P\colon\N\rightarrow\N$, and $f_W\colon\N\rightarrow\N$ such that for every $t\in\N$ the following holds:
	Let 
	\begin{itemize}
		\item $D$ be a digraph,
		\item $W$ be a cylindrical $3k$-wall with $k\geq\Fkt{f_W}{t}$ in $D$,
		\item $\Triade=\Brace{W,k,W_1,W_2,W_3,W^1,W^2,W^3}$ be the \hyperref[def:triadicpartition]{triadic partition} of $W$, and
		\item \hyperref[def:tiling]{$\Tiling=\Tiling_{W,k,\Fkt{f_w}{t},\xi,\xi'}$} for some $\xi,\xi'\in[\Fkt{f_w}{t}+1]$.
	\end{itemize}
	If there exists a subfamily $\Tiling'$ of $\Tiling$ and a family $\Jumps$ of pairwise disjoint directed paths in $D$ with the following properties:
	\begin{enumerate}
		\item Every member of $\Jumps$ is internally disjoint from $W$ but has both endpoints on $W$,
		\item $\Abs{\Tiling'}=\Abs{\Jumps}=\Fkt{f_P}{t}$,
		\item for every $T_{\Fkt{\ColumnFunction}{p},\Fkt{\RowFunction}{q},\Fkt{f_w}{t}}\neq T_{\Fkt{\ColumnFunction}{p'},\Fkt{\RowFunction}{q'},\Fkt{f_w}{t}}\in\Tiling'$ we have $\max\Set{\Abs{p-p'},\Abs{q-q'}}\geq 2$,
		\item there exists a bijection $\Start\colon\Tiling'\rightarrow\Jumps$ ($\End\colon\Tiling'\rightarrow\Jumps$) such that the starting point (endpoint) of the path $\Fkt{\Start}{T}$ ($\Fkt{\End}{T}$) belongs to the \hyperref[def:tile]{centre} of $T$,
		\item $\V{\Fkt{\Start}{T}}\cap\V{\Tiling'}$ ($\V{\Fkt{\End}{T}}\cap\V{\Tiling'}$) contains exactly the endpoint of $\Fkt{\Start}{T}$ ($\Fkt{\End}{T}$) where $\V{\Tiling'}=\bigcup_{T'\in\Tiling'}\V{T'}$, and finally
		\item the endpoints (starting points) of the paths in $\Jumps$ are of mutual \hyperref[def:Wdistance]{$W$-distance} at least $4$. 
	\end{enumerate}
	Then one of the following is true.
	\begin{enumerate}
		\item[a)] $D$ has a $\Bidirected{K_t}$-butterfly minor grasped by $W$,
		\item[b)] there exists a family of tiles $\Tiling''\subseteq\Tiling'$ all contained in a single \hyperref[def:strip]{strip}  $S\subseteq W$ of height equal to the height of the tiles in $\Tiling$ such that
		\begin{itemize}
			\item we can number $\Tiling''=\Set{T_1,\dots,T_h}$ such that for each $i\in[h]$ the graph $S-T_i$ has one component containing exactly the tiles $T_1,\dots,T_{i-1}$,
			\item $\Abs{\Tiling''}\geq \Fkt{f_P}{t}^{\frac{1}{4}}$,
			\item for every $i\in[h-1]$ the tiles $T_i$ and $T_{i+1}$ are separated in $W$ by a \hyperref[def:slice]{slice} of width equal to the width of the tiles in $\Tiling$, and
			\item there is a family $\Jumps'\subseteq\Jumps$ with $\Abs{\Jumps'}=\Abs{\Tiling''}$ such that for each $T\in\Tiling''$ we have $\Fkt{\Start}{T}\in\Jumps'$ ($\Fkt{\End}{T}\in\Jumps'$), and
			\item for each $i\in[h]$ the endpoint of $\Fkt{\Start}{T}$ (starting point of $\Fkt{\End}{T}$) lies in the component of $S-T_i$ that contains no tiles of $\Tiling''$ if $i=1$, or in the splice of $S$ separating $T_{i-1}$ and $T_i$.
		\end{itemize}
		\item[c)] there exists a family of tiles $\Tiling''\subseteq\Tiling'$ all contained in a single \hyperref[def:strip]{strip}  $S\subseteq W$ of height equal to the height of the tiles in $\Tiling$ such that $\Tiling''$ and $S$ meet the properties of outcome b) after taking the mirror image of $W$ along a vertical line.
	\end{enumerate}
\end{lemma}

By using the bounds obtained from the proofs in \cite{Giannopoulou2020DirectedFlatWall}, we get the following rough estimates for the functions $f_w$, $f_P$, and $f_W$:
\begin{enumerate}
	\item $\Fkt{f_w}{t}=2^9t^{10}$,
	\item $\Fkt{f_P}{t}=2^7t^8$, and
	\item $\Fkt{f_W}{t}=2^{32+t^{30}}$.
\end{enumerate}

Note that in case \cref{lemma:directedlongjumps} yields the existence of a $\Bidirected{K_{3t}}$ butterfly minor, this minor contains $t$ pairwise disjoint weak odd bicycles and thus, by \cref{obs:oddbicycleevendicycle}, it must contain $t$ pairwise disjoint even dicycles.
If this does not hold, then \cref{lemma:directedlongjumps} produces many pairwise disjoint slices of the wall $W$, mutually far apart from each other, and each of them having a neighbouring slice together with a fairly long jump into this neighbour.
Such a slice, together with its neighbouring slice and the jump forms a digraph in which one can easily find a cylindrical grid of order $3$ together with a perimeter jump as a butterfly minor.
Hence \cref{lem:gridplusjumpmakesoddbicycle} yields the existence of many pairwise disjoint weak odd bicycles, each providing an even dicycle.
Thus \cref{lemma:directedlongjumps} is already powerful enough to guarantee a large even dicycle packing in case we find many long jumps over a sufficiently large wall.

\subsection{Removing long jumps}

We will now make our observations from above more explicit to obtain two auxiliary results providing ways to find bounded size hitting sets for long jumps over our wall.


\begin{definition}[Auxiliary digraph type I]\label{def:auxiliarydigraphI}
	Let $t,k,k',w\in\N$ be positive integers such that $k\geq k'\geq2\Fkt{f_W}{t}+4\Fkt{f_P}{t}\Brace{2w+1}$, $w\geq 2\Fkt{f_w}{t}$, and $\xi,\xi'\in\left[w+1\right]$.
	Let $D$ be a digraph containing a cylindrical $3k$-wall $W=\Brace{Q_1,\dots,Q_{3k},\hat{P}_1,\dots,\hat{P}_{3k}}$ with its \hyperref[def:triadicpartition]{triadic partition} $\Triade=\Brace{W,k,W_1,W_2,W_3,W^1,W^2,W^3}$, and a tiling \hyperref[def:tiling]{$\Tiling=\Tiling_{W,k,w,\xi,\xi'}$}.
	Let $i\in[4]$, $\Set{\Class_1,\dots,\Class_4}$ be a four colouring of $\Tiling$ and $W'\subseteq W$ be a \hyperref[def:slice]{slice} of width $k'$ of $W_2$.
	At last, let us denote by $\Tiling'$ the family of tiles from $\Tiling$ that share a vertex with $W'$.
	Similarly let $\Class_i'\coloneqq\Tiling'\cap\Class_i$.
	Then $\Fkt{D^1_i}{W'}$ is the digraph obtained from $D$ by performing the following construction steps for every $T\in\Class_i'$:
	\begin{enumerate}
		\item add new vertices $x_T^{\text{in}}$ and $x_T^{\text{out}}$,
		\item for every vertex $u$ in the centre of $T$ introduce the edges $\Brace{u,x_T^{\text{in}}}$ and $\Brace{x_T^{\text{out}},u}$, and then
		\item delete all internal vertices of $T$.
	\end{enumerate}
\end{definition}

\begin{lemma}\label{lemma:longjumps1}
	Let $t,k,k',w\in\N$ be positive integers such that $k\geq k'\geq2\Fkt{f_W}{3t}+2^{16}\Fkt{f_P}{3t}+2$, $w\geq 2\Fkt{f_w}{3t}+2^7\Fkt{f_P}{3t}$, $3k$ is divisible by $4$, and $\xi,\xi'\in\left[w+1\right]$.
	Let $D$ be a digraph containing a cylindrical $3k$-wall $W=\Brace{Q_1,\dots,Q_{3k},\hat{P}_1,\dots,\hat{P}_{3k}}$ with its \hyperref[def:triadicpartition]{triadic partition} $\Triade=\Brace{W,k,W_1,W_2,W_3,W^1,W^2,W^3}$, and a tiling \hyperref[def:tiling]{$\Tiling=\Tiling_{W,k,w,\xi,\xi'}$}.
	Let $i\in[4]$, $\Set{\Class_1,\dots,\Class_4}$ be a four colouring of $\Tiling$ and $W'\subseteq W$ be a \hyperref[def:slice]{slice} of width $k'$ of $W_2$.
	
	Now let $\Tiling'$ be the family of all tiles of $\Tiling$ that are completely contained in $W'$ and let $\widetilde{W}$ be the smallest slice of $W$ that contains all tiles from $\Tiling'$.
	
	Consider the \hyperref[def:auxiliarydigraphI]{auxiliary digraph of type I} $\Fkt{D_i^1}{\widetilde{W}}$ and let $\Class_i'$ be as in the definition of $\Fkt{D_i^1}{\widetilde{W}}$.
	Define the sets 
	\begin{align*}
		X_{\text{I}}^{\text{out}}&\coloneqq\CondSet{x_T^{\text{out}}}{T\in\Class_i'}\text{, and}\\
		X_{\text{I}}^{\text{in}}&\coloneqq\CondSet{x_T^{\text{in}}}{T\in\Class_i'}.
	\end{align*}
	Further, we construct the set $Y_{\text{I}}$ as follows:
	Let $Q$ and $Q'$ be the two cycles of $\Perimeter{W_2}$.
	For every $j\in\left[\frac{3k}{4} \right]$, $Y_{\text{I}}$ contains exactly one vertex of $Q\cap P^1_{4j}$, $Q\cap P^2_{4j+2}$, $Q'\cap P^1_{4j}$, and $Q'\cap P^2_{4j+2}$ each.
	
	If there exists a family $\mathcal{L}$ of pairwise disjoint directed paths in $\Fkt{D_i^1}{\widetilde{W}}$ with $\Abs{\mathcal{L}}= 2^7\Fkt{f_P}{3t}$ such that either
	\begin{itemize}
		\item $\mathcal{L}$ is a family of directed $X_{\text{I}}^{\text{out}}$-$Y_{\text{I}}$-paths, or
		\item $\mathcal{L}$ is a family of directed $Y_{\text{I}}$-$X_{\text{I}}^{\text{in}}$-paths,
	\end{itemize}
	then $D$ has an integral packing of $t$ even dicycles.
\end{lemma}

\begin{proof}
	This proof is divided into two steps.
	First we show that the second and third outcome of \cref{lemma:directedlongjumps} yields a large enough even dicycle packing.
	Afterwards our goal is to construct a cylindrical wall $U\subseteq W$ of sufficient size, together with a family of $\Fkt{f_P}{t}$ directed $U$-paths that meet the requirements of \cref{lemma:directedlongjumps}.
	
	\textbf{Claim 1:} \emph{In the situation of \cref{lemma:directedlongjumps} any outcome yields an integral packing of $\Floor{\frac{1}{3}t}$ even dicycles.}
	
	\emph{Proof of Claim 1}
	If the outcome is case a), then we can immediately find $\Floor{\frac{1}{3}t}$ pairwise disjoint odd bicycle models within the $\Bidirected{K_t}$ butterfly minor model.
	By \cref{obs:oddbicycleevendicycle} this yields the desired outcome.
	Hence we may assume case b) or c) holds.
	Since these two cases are symmetric it suffices to only consider case b).
	Here we find a family $\mathcal{T}''\subseteq\mathcal{T}'$ of $h\geq \Fkt{f_P}{t}^{\frac{1}{4}}\geq \frac{1}{3}t$ contained in a single strip $S\subseteq W$ which is of the same height as the tiles in $\mathcal{T}''$.
	Moreover, we can number the tiles in $\mathcal{T}''=\Set{T_1,\dots,T_h}$ such that for each $i\in[h]$ there exists a slice $S_i\subseteq W$ with width equal to the width of the tiles in $\mathcal{T}''$ whose intersection with $S$ is exactly the tile $T_i$.
	Moreover, for each $i\in[h]$ we can find a slice $H_i$ of $W$ containing $S_i$ and both endpoints of the jump $\Fkt{\Start}{T_i}$ such that $H_i$ and $H_j$ are disjoint if $i\neq j$. 
	Let $J_i$ denote the path $\Fkt{\Start}{T_i}$.
	Then for each $i\in[h]$ there exists a slice $S_i'\subseteq S_i$ of width $3$ whose left perimeter coincides with the left perimeter of $S_i$.
	Moreover, $S_i'$ separates the two endpoints of $J_i$ within $W$.
	We can now find a path from the right perimeter of $S_i'$ to the tail of $J_i$ and a path from the head of $J_i$ to the left perimeter of $S_i'$ within $H_i$.
	The resulting digraph contains, as a butterfly minor, a cylindrical grid of order $3$ with a perimeter jump.
	As we can find $h$ such graphs, by \cref{lem:gridplusjumpmakesoddbicycle} and \cref{obs:oddbicycleevendicycle} we can find $h\geq \frac{1}{3}t$ pairwise disjoint even dicycles in $D$.	\hfill$\blacksquare$
	
	In the statement of our assertion we have chosen to insert $3t$ as the arguments for all functions.
	This means any application of \cref{lemma:directedlongjumps} yields, under Claim 1, an integral packing of $t$ even dicycles.
	With this in mind we can now proceed with the main body of our proof.
	Without loss of generality, let us assume $\mathcal{L}$ is a family of directed $X_{\text{I}}^{\text{out}}$-$Y_{\text{I}}$-paths.
	The other case can be seen using similar arguments.
	
	Let $L$ and $P$ be directed paths.
	We say that $P$ is a \emph{long jump of $L$} if $P$ is a $w$-long jump over $W$ and $P\subseteq L$.
	We also say that $P$ is a \emph{jump of $L$} if $P$ is a directed $W$-path.
	
	Towards our goal, we first show that we can use $\mathcal{L}$ to construct a half-integral $X_{\text{I}}^{\text{out}}$-$Y_{\text{I}}$-linkage $\mathcal{L}_1$ such that
	\begin{enumerate}
		\item $\Abs{\mathcal{L}_1}=2^7\Fkt{f_P}{3t}$,
		\item there exists a family $\mathcal{F}\subseteq\Tiling'$ with $\Abs{\mathcal{F}}\leq 2^7\Fkt{f_P}{3t}$, and
		\item for every $L\in\mathcal{L}_1$, every endpoint $u$ of a jump of $L$ with $u\in\V{\widetilde{W}}$ belongs to a tile from $\Class_i'\cup\mathcal{F}$.
	\end{enumerate}
	Once this is achieved, we use \cref{thm:halfintegral} to obtain a family $\mathcal{L}_2$ of pairwise disjoint directed $X_{\text{I}}^{\text{out}}$-$Y_{\text{I}}$-paths of size $2^6\Fkt{f_P}{3t}$ from $\mathcal{L}_1$.
	Afterwards, we remove the cycles and paths of $\widetilde{W}$ that meet tiles from $\mathcal{F}$ and obtain a new slice $\widetilde{W}'$ of some cylindrical wall.
	For this slice, we construct a \hyperref[def:tiling]{tiling} and a \hyperref[def:tilingII]{tier II tiling} as well as a half-integral linkage $\mathcal{L}_3'$ of size $2^6\Fkt{f_P}{3t}$ from $\mathcal{L}_2$
    that connects the \hyperref[def:tile]{centres} of some tiles in the tier II tiling to vertices of $\widetilde{W}'$, such that their endpoints are mutually far enough apart.
    We then further refine $\mathcal{L}_3'$ into the linkage $\mathcal{L}_4$, such that every path in $\mathcal{L}_4$ is internally disjoint from a large subwall of $\widetilde{W}'$.
	Another application of \cref{thm:halfintegral} then yields the family of long jumps necessary for an application of \cref{lemma:directedlongjumps}.
	
	We start out with the construction of $\mathcal{L}_1$ and $\mathcal{F}$.
	For this, let $\mathcal{L}'\coloneqq \mathcal{L}$, $\mathcal{L}_1\coloneqq \emptyset$, and $\mathcal{F}\coloneqq\emptyset$.
	As long as $\mathcal{L}'$ is non-empty, perform the following actions:
	
	Select some path $L\in\mathcal{L}$.
	In case $L$ is internally disjoint from $\widetilde{W}$, add $L$ to $\mathcal{L}_1$ and remove it from $\mathcal{L}'$.
	Otherwise let $s_L$ be its starting point and let $v_L$ be the first vertex of $L$ that belongs to $\widetilde{W}$, but not to a tile from $\Class_i'$, when traversing along $L$ starting from $s_L$.
	\begin{enumerate}
		\item If $v_L$ does not belong to a tile from $\mathcal{F}$, let $T\in\Tiling'\setminus\Class_i$ be the tile that contains $v_L$ and add $T$ to $\mathcal{F}$.
		Let $R$ be a shortest directed path from $v_L$ to $Y_{\text{I}}$ in $W$ such that $R$ avoids all vertices of $W$ that are contained in two different paths of $\mathcal{L}_1$ and that is internally disjoint from $Lv_L$.
		Now add $Lv_LR$ to $\mathcal{L}_1$ and remove $L$ from $\mathcal{L}'$.
		Note that such a path $R$ must exist since the paths in $\mathcal{L}$ are pairwise disjoint, we never used $T$ for such a re-routing before, and $w$ and $k'$ are chosen sufficiently large in proportion to $2^7\Fkt{f_P}{3t}$.
		Also note that the path $R$ is exactly the part, where we might go from integral to half-integral, but since our paths were pairwise disjoint to begin with, we can be sure that $R$ never meets a vertex contained in two distinct paths.
		
		\item So now suppose $v_L$ belongs to a tile $T$ from $\mathcal{F}$.
		Let us follow along $v_LL$ until the first time we encounter a vertex $u_L$ for which one of the following is true:
		\begin{enumerate}
			\item $u_L$ belongs to a tile $T$ from $\Tiling'\setminus\Brace{\Class_i'\cup\mathcal{F}}$, or
			\item every internal vertex of $u_LL$ belongs to $W-\widetilde{W}$ or to some tile from $\Class_i'\cup\mathcal{F}$. 
		\end{enumerate}
		If a) is the case, repeat the instruction from i) but replace $v_L$ by $u_L$.
		In this case $T$ is added to $\mathcal{F}$.
		Otherwise b) must hold and here we may simply remove $L$ from $\mathcal{L}'$ and add it to $\mathcal{L}_1$.
	\end{enumerate}
	Now for every $L\in\mathcal{L}$ we added at most one tile to $\mathcal{F}$ and thus $\Abs{\mathcal{F}}\leq\Abs{\mathcal{L}}$.
	Moreover, from the construction it is clear that $\mathcal{L}_1$ is indeed a half-integral linkage from $X_{\text{I}}^{\text{out}}$ to $Y_{\text{I}}$.
	Also, please note that we may assume that every $L$ meets each tile in $\mathcal{F}$ in at most $2^7\Fkt{f_P}{3t}+1$ horizontal path pairs and vertical cycles, since otherwise one could find a short cut through $W$ itself.
	
	Next, we may apply \cref{thm:halfintegral} to obtain a family $\mathcal{L}_2$ of pairwise disjoint directed $X_{\text{I}}^{\text{out}}$-$Y_{\text{I}}$-paths with $\V{\mathcal{L}_2}\subseteq\V{\mathcal{L}_1}$ and $\Abs{\mathcal{L}_2}=2^6\Fkt{f_P}{3t}$.
	This completes the second step.
	
	For the third step, let us consider $W''\coloneqq \InducedSubgraph{\widetilde{W}}{\Tiling,i}$ together with the tiling $\Tiling''\coloneqq\TierIITiling{\Tiling,i,f}{\widetilde{W}}$ and a four-colouring $\Set{\widetilde{\Class}_1,\dots,\widetilde{\Class}_4}$ of $\Tiling''$.
	Note that by choice of $k'$ this means that $W''$ is a slice of width $k''\geq \Fkt{f_W}{3t}+2^7\Fkt{f_P}{3t}\Brace{2w+1}+1$ of some cylindrical $3k''$-wall that is completely contained in $W$.
	For each $L\in\mathcal{L}_2$ let $T^1_L\in \Class_i$ such that the successor $w_L$ of the starting point $s_L$ of $L$ belongs to $T^1_L$.
	Choose any vertex $s'_L$ of degree three in $W''$ that is not contained in any path of $\mathcal{L}_2$, and let $R_L$ be a directed path from $s_L'$ to $w_L$ within $T^1_L$.
	Let $\mathcal{L}_3'$ be the resulting, and potentially now again half-integral, family of directed paths $s_L'R_Lt_lL$.
	Now there must exist $j\in[4]$ such that at least $2^4\Fkt{f_P}{3t}$ of the paths from $\mathcal{L}_3'$ start at the centre of a tile from $\widetilde{\Class}_j$.
	Let $\mathcal{L}_3\subseteq\mathcal{L}_3'$ be a family of exactly $2^4\Fkt{f_P}{3t}$ such paths.
    This concludes the third step.
    
	Next let us consider the family $\mathcal{F}$.
	Let $W'''$ be the subgraph of $W''$ induced by all vertical cycles and horizontal path pairs in $W''$ that do not contain a vertex of some tile in $\mathcal{F}$ that belongs to a path in $\mathcal{L}_3$.
	Since $\Abs{\mathcal{F}}\leq 2^7\Fkt{f_P}{3t}$ and each tile in $\mathcal{F}$ meets a path in $\mathcal{L}_3$ in at most $2^7\Fkt{f_P}{3t}+1$ such cycles and pairs of horizontal paths, it follows that $W'''$ is a slice of width $k'''\geq\Fkt{f_W}{3t}+2$ of some cylindrical $3k'''$wall $W^*\subseteq W$.
	Moreover, $W^*$ can be partitioned into three \hyperref[def:slice]{slices} of width $k'''$ as in its \hyperref[def:triadicpartition]{triadic partition}, such that $W'''$ is the slice in the middle.
	Let us rename the paths and cycles of $W^*$ such that $W^*=\Brace{Q^*_1,\dots,Q^*_{3k'''},\hat{P^*}_1,\dots,\hat{P^*}_{3k'''}}$, and we construct the set $Y^*$ as follows:
	Let $Q^*$ and $'Q^*$ be the two cycles of $\Perimeter{W'''}$.
	For every $j\in\left[ \frac{3k'''}{4} \right]$, $Y^*$ contains exactly one vertex of $Q^*\cap P^{*1}_{4j}$, $Q\cap P^{*2}_{4j+2}$, $'Q^*\cap P^{*1}_{4j}$, and $'Q^*\cap P^{*2}_{4j+2}$ each.
	Let $L\in\mathcal{L}_3$ be any path and $t_L$ be the first vertex after its starting point $L$ shares with either $W'''$ or $W^*-W'''$.
	In case $t_L\in\V{W'''}$, simply add $Lt_L$ to $\mathcal{L}_3'''$.
	Otherwise, let $b_L$ be the endpoint of $L$ in $W^*-W'''$.
	Then we can find a path $R_L$ in $W$ from $b_L$ to a vertex $t^*_L$ of $Y^*$ such that $t^*_L$ is of $W^*$-distance at least $4$ to every endpoint of every path already in $\mathcal{L}_3'''$, $R_L$ is internally disjoint from $L$, and $R_L$ does not contain a vertex that is contained in two distinct paths from $\mathcal{L}_3'''$.
	Add $LR_L$ to $\mathcal{L}_3'''$.
	Finally, $\mathcal{L}_3'''$ is a half-integral linkage from the set of starting points $S^*$ of the paths in $\mathcal{L}_3$ to $Y^*$ of size $2^4\Fkt{f_P}{3t}$, and thus by \cref{thm:halfintegral} we can find a family $\mathcal{L}_4$ of pairwise disjoint directed paths from $S^*$ to $Y^*$ with $\V{\mathcal{L}_4}\subseteq\V{\mathcal{L}_3'''}$ that is of size $2^3\Fkt{f_P}{3t}$.
	It follows that all paths in $\mathcal{L}_4$ are internally disjoint from $W'''$.
	This concludes the fourth step.
	
	Let us consider the tiles of $\widetilde{\Class}_i$ whose centres contain a vertex of $S^*$.
	Since $W'''$ might be a proper subgraph of $W''$, $\Tiling''$ is not necessarily a tiling of $W'''$.
	Each such tile $T$, however, contains a tile $T'$ of width $\Fkt{f_w}{3t}$ with the same \hyperref[def:tile]{centre}.
	Since $T$ can be surrounded by at most $8$ tiles from $\mathcal{F}$ in $W'$, we may find, among the $2^3\Fkt{f_P}{3t}$ many such tiles, a family $\Jumps$ of $\Fkt{f_P}{3t}$ tiles that are pairwise disjoint and thus, since they all are construcuted from the family $\widetilde{\Class}_i$, they meet the distance requirements of the tiles in \cref{lemma:directedlongjumps}.
	Hence we may apply \cref{lemma:directedlongjumps} and by Claim 1 we are done.
\end{proof}

The above lemma allows us to deal with long jumps that all start (end) in the centre of tiles from a single colour class but end (start) in tiles from the remaining three classes.
By utilising a second, slightly altered auxiliary digraph we can make use of \cref{thm:xpaths} to also deal with long jumps between tiles of the same colour.

\begin{definition}[Auxiliary digraph type II]\label{def:auxiliarydigraphII}
	Let $t,k,k',w\in\N$ be positive integers such that $k\geq k'\geq2\Fkt{f_W}{t}$, $w\geq 2\Fkt{f_w}{t}$, and $\xi,\xi'\in\left[w+1\right]$.
	Let $D$ be a digraph containing a cylindrical $3k$-wall $W=\Brace{Q_1,\dots,Q_{3k},\hat{P}_1,\dots,\hat{P}_{3k}}$ with its \hyperref[def:triadicpartition]{triadic partition} $\Triade=\Brace{W,k,W_1,W_2,W_3,W^1,W^2,W^3}$, and a \hyperref[def:tiling]{tiling} $\Tiling=\Tiling_{W,k,w,\xi,\xi'}$.
	Let $i\in[4]$, $\Set{\Class_1,\dots,\Class_4}$ be a four colouring of $\Tiling$, and $W'\subseteq W$ be a \hyperref[def:slice]{slice} of width $k'$ of $W_2$ such that no \hyperref[def:tile]{tile} of $\Class_i$ contains a vertex of the perimeter of $W'$.
	Then $\Fkt{D^2_i}{W'}$ is the digraph obtained from $D$ by performing the following construction steps:
	
	for every $T\in\Class_i$, such that $T$ contains a vertex of $W'$, we do the following:
	\begin{enumerate}
		\item add a new vertex $x_T$, and
		\item for every vertex $v$ that belongs to the interior or the centre of $T$, introduce the edges $\Brace{x_T,v}$ and $\Brace{v,x_T}$.
	\end{enumerate}
	Once this is done, delete all vertices of $W'$ that do not belong to tiles of $\Class_i$.
	Let $X_{\text{II}}^i$ be the collection of all newly introduced vertices $x_T$.
\end{definition}

\begin{lemma}\label{lemma:longjumpspahse2}
	Let $t,k,k',w\in\N$ be positive integers, and $\xi,\xi'\in\left[w+1\right]$ where $w\geq2\Fkt{f_w}{3t}$.
	Let $D$ be a digraph containing a cylindrical $3k$-wall $W_0=\Brace{Q_1,\dots,Q_{3k},\hat{P}_1,\dots,\hat{P}_{3k}}$, where $k\geq k'\geq4\Fkt{f_W}{3t}^2$, with its \hyperref[def:triadicpartition]{triadic partition} $\Triade=\Brace{W_0,k,W_1,W_2,W_3,W^1,W^2,W^3}$, a \hyperref[def:slice]{slice} $W\subseteq W_2$ of width $k'$, a \hyperref[def:tiling]{tiling} $\Tiling=\Tiling_{W,k,w,\xi,\xi'}$, a four colouring $\Set{\Class_1,\dots,\Class_4}$, and a fixed colour $i\in[4]$.
	
	Then $D$ has an integral packing of $t$ even dicycles, or there exists a set $\MarkedTiles{2}_{i,\xi,\xi'}\subseteq\Tiling$ with $\Abs{\MarkedTiles{2}_{i,\xi,\xi'}}\leq 8\Fkt{f_P}{3t}$ and a set $Z^2_{i,\xi,\xi'}\subseteq\V{D-W}$ with $\Abs{Z^2_{i,\xi,\xi'}}\leq 8\Fkt{f_P}{3t}$ such that every directed $\V{W_0}$-path $P$ in $D-Z^2_{i,\xi,\xi'}$ whose endpoints belong to different tiles of $\Class_i$ contains a vertex of some tile in $\MarkedTiles{2}_{i,\xi,\xi'}$.
\end{lemma}

\begin{proof}
	Let $W'$ be the largest \hyperref[def:slice]{slice} of $W$ such that no tile of $\Class_i$ contains a vertex of $\Perimeter{W'}$.
	Let us consider the \hyperref[def:auxiliarydigraphI]{auxiliary digraph} $\Fkt{D^2_i}{W'}$ with the set $X_{\text{II}}^i$ of newly added vertices.
	By applying \cref{thm:xpaths} to the set $X_{\text{II}}^i$ in $\Fkt{D^2_i}{W'}$, we either find a set $Z$ of size at most $8\Fkt{f_P}{3t}$ that hits all directed $X_{\text{II}}^i$-paths, or there exists a family $\mathcal{J}'$ of $4\Fkt{f_P}{3t}$ pairwise disjoint directed $X_{\text{II}}^i$-paths in $\Fkt{D^2_i}{W'}$.
	
	Let us assume the latter.
	Then, by construction of $\Fkt{D^2_i}{W'}$, no path in $\mathcal{J}'$ contains a vertex of $W_2$.
	Back in the digraph $D$, let us consider the \hyperref[def:tilingII]{tier II tiling} $\Tiling''\coloneqq \TierIITiling{\Tiling,i,w}{W'}$ of width $w$ of $W''\coloneqq \InducedSubgraph{W'}{\Tiling,i}$.
	Notice that, by choice of $k$, $W''$ still contains a cylindrical $\Fkt{f_W}{3t}$-wall $W'''$ such that the perimeter of every tile $T\in\Class_i$, for which $x_T$ is an endpoint of some path in $\mathcal{J}'$, bounds a cell of $W'''$.
	Let $\Tiling'''$ be a tiling of $W'''$ such that the perimeter of every $T\in\Class_i$, for which $x_T$ is an endpoint of a path in $\mathcal{J}'$, is the centre of some tile in $\Tiling'''$.
	We now consider a four colouring $\Set{\Class_1',\dots,\Class_4'}$ of $\Tiling'''$.
	Then there must exist $j\in[4]$ and a family of paths $\mathcal{J}''$ of size $\Fkt{f_P}{3t}$ such that the starting points of every path in $\mathcal{J}''$ belongs to a tile of $\Class_i$ whose perimeter is the centre of a tile in $\Class_j'$.
	For every $J''\in\mathcal{J}''$, let $T_1,T_2\in\Class_i$ be the two tiles such that $J''$ is a directed $x_{T_1}$-$x_{T_2}$-path.
	We can now find a directed path $J$ that starts on the perimeter of $T_1$, ends on the perimeter of $T_2$, and is internally disjoint from $W'''$.
	Hence we find a family $\mathcal{J}$ of pairwise disjoint directed $W'''$-paths whose endpoints all lie on the centres of distinct tiles of $\Tiling'''$ and that all start at the centres of tiles from $\Class_j'$.
	So we may apply \cref{lemma:directedlongjumps} and by Claim 1 from the proof of \cref{lemma:longjumps1} we have found the desired packing of even dicycles.
	
	Therefore we may assume that we find a set $Z$ of size at most $8\Fkt{f_P}{3t}$ that hits all directed $X_{\text{II}}^i$-paths.
	Let $Z^2_{i,\xi,\xi'}\coloneqq Z\cap \V{D}$, and $\MarkedTiles{2}_{i,\xi,\xi'}\coloneqq\CondSet{T\in\Tiling}{x_T\in Z}$.
	Since $\Abs{Z}\leq 8\Fkt{f_P}{3t}$, the bounds on the two sets follow immediately.
	Moreover, since $Z$ meets every directed $X_{\text{II}}^i$-path in $D'$, every directed path with endpoints in distinct tiles of $\Class_i$ which is otherwise disjoint from $W_0$ must contain a vertex from $Z^2_{i,\xi,\xi'}$ or meet a tile from $\MarkedTiles{2}_{i,\xi,\xi'}$.
\end{proof}

\subsection{Proof of \Cref{thm:oddwall}}

A final piece we require to complete the proof is a local version of Menger's Theorem.

\begin{theorem}[Menger's Theorem \cite{menger1927allgemeinen}]\label{thm:directedlocalmenger}
	Let $D$ be a digraph and $X,Y\subseteq\V{D}$ be two sets of vertices, then the maximum number of pairwise disjoint directed $X$-$Y$-paths in $D$ equals the minimum size of a set $S\subseteq\V{G}$ such that every directed $X$-$Y$-path in $D$ contains a vertex of $S$.
\end{theorem}

Throughout the proof we will collect set of tiles $\mathcal{F}_i$  from different tilings.
In general we will say that a tile is \emph{marked} if it either contains a vertex of some separator obtained from \cref{thm:directedlocalmenger} and \cref{lemma:longjumps1}, or through \cref{lemma:longjumpspahse2}, or one of these separators contains the fresh vertex created for it in the construction of some auxiliary graph (type I or II).
This means that the overall number of tiles that are marked exceeds the size of the sets $\mathcal{F}_i$.
Still, the number of columns that can contain marked tiles is bounded, which is enough for our purposes.
A vertex of our wall $W$ is said to be \emph{marked} if it belongs to a separator obtained from \cref{thm:directedlocalmenger} and \cref{lemma:longjumps1}, or through \cref{lemma:longjumpspahse2}.
In each of the steps we introduce families of marked tiles and vertices and a slice of $W$ is said to be \emph{clear} if it does not contain vertices that are marked.

\begin{proof}[Proof of \Cref{thm:oddwall}]
Let $r,t\in\N$ be positive integers, $D$ be a digraph and $W$ be a cylindrical $\OddWallOrder{t,r}$-wall, where $\OddWallOrder{t,r}$ will be determined throughout the proof.
To do this we will introduce a constant $d_i$ for each step, for which we will make more and more assumptions in the form of lower bounds.
\begin{align*}
	\text{Let us assume $\OddWallOrder{t,r}\geq 3d_1$.}
\end{align*}
Then $W$ is the cylindrical $3d_1$-wall $U=\Brace{Q_1,\dots,Q_{3k},\hat{P}_1,\dots,\hat{P}_{3k}}$ with its triadic partition $\mathcal{W}=\Brace{W,d_1,\widetilde{W}_1,\widetilde{W}_2,\widetilde{W}_3,\widetilde{W}^1,\widetilde{W}^2,\widetilde{W}^3}$.
Throughout the proof let us fix
\begin{align*}
	w\coloneqq2\Fkt{f_w}{3t}+2^7\Fkt{f_P}{3t}.
\end{align*}

\paragraph{Phase I}
Phase I is divided into $2^{11+8\Fkt{f_P}{3t}}\Fkt{f_P}{3t}$ \emph{rounds}, each of which produces a \hyperref[def:slice]{slice} $W_i$ of $W_0\coloneqq \widetilde{W}_2$ which is clean with respect to all tiles that have not been marked up to this point.
We also obtain sets $F_i\subseteq\V{D-\bigcup_{j=1}^{i-1}F_j}$ and $\mathcal{F}_i$ of marked vertices and tiles respectively for each round $i\in[2^{11+8\Fkt{f_P}{3t}}\Fkt{f_P}{3t}]$, and in each round $i$ we will work on the digraph $D_i\coloneqq D_{i-1}-F_i$.
In this context, whenever we ask for a clean slice of the current slice $W_{i-1}$ or $W'_i$ we ask for a slice $W'$ such that there do not exist $\xi,\xi'\in[w+1]$ whose corresponding \hyperref[def:tiling]{tiling} of $W_0$ has a \hyperref[def:tile]{tile} $T$, that is marked or contains a vertex of any separator set found so far, which satisfies $\V{T}\cap\V{W'}\neq\emptyset$.

Each round is divided into two steps, \emph{Step I} and \emph{Step II}.
Let $i\in[2^{11+8\Fkt{f_P}{3t}}\Fkt{f_P}{3t}]$.
From the previous round or the initial state, we may assume to have obtained the following sets and graphs, which serve as the \emph{input} for round $i$:
\begin{itemize}
	\item $F_{I,i-1}\subseteq\V{D}$ such that $\Abs{F_{I,i-1}}\leq \Brace{i-1}\Brace{2^{11}\Brace{w+1}^2\Fkt{f_P}{3t}+2^5\Brace{w+1}^2\Fkt{f_P}{3t}}$,
	\item $D_{i-1}\coloneqq D-F_{I,i-1}$,
	\item $\mathcal{F}_{I,i-1}\subseteq \bigcup_{\xi,\xi'\in[w+1]}\Tiling_{W_0,d_1,w,\xi,\xi'}$ of size at most $\Brace{i-1}\Brace{2^{11}\Brace{w+1}^2\Fkt{f_P}{3t}+2^5\Brace{w+1}^2\Fkt{f_P}{3t}}$, and
	\item a slice $W_{i-1}$ of $W_0$ of width $\Brace{\Brace{2^{12}\Brace{w+1}^3\Fkt{f_P}{3t}+1}\Brace{2^6\Brace{w+1}^3\Fkt{f_P}{3t}+1}}^{2^{11+8\Fkt{f_P}{3t}}\Fkt{f_P}{3t}-i+1}d_2$ that is clean with respect to $F_{I,i-1}$ and $\mathcal{F}_{I,i-1}$ such that, in case $i-1\geq 1$, every long jump over $W_{i-1}$ in $D_{i-1}$ contains a vertex of some tile in $\mathcal{F}_{I,j}\setminus\mathcal{F}_{I,j-1}$ for every $j\in[i-1]$.
\end{itemize}
For $i=1$ the inputs are the graphs $D_0\coloneqq D$, $W_0$, and two empty sets.

After round $i$ is complete we require the following sets and graphs as its \emph{output}:
\begin{itemize}
	\item $F_{I,i}\subseteq\V{D}$ such that $\Abs{F_{I,i}}\leq i\Brace{2^{11}\Brace{w+1}^2\Fkt{f_P}{3t}+2^5\Brace{w+1}^2\Fkt{f_P}{3t}}$,
	\item $D_{i}\coloneqq D-F_{I,i}$,
	\item $\mathcal{F}_{I,i}\subseteq \bigcup_{\xi,\xi'\in[w+1]}\Tiling_{W_0,d_1,w,\xi,\xi'}$ of size at most $i\Brace{2^{11}\Brace{w+1}^2\Fkt{f_P}{3t}+2^5\Brace{w+1}^2\Fkt{f_P}{3t}}$, and
	\item a slice $W_i$ of $W_0$ of width $\Brace{\Brace{2^{12}\Brace{w+1}^3\Fkt{f_P}{3t}+1}\Brace{2^6\Brace{w+1}^3\Fkt{f_P}{3t}+1}}^{2^{11+8\Fkt{f_P}{3t}}\Fkt{f_P}{3t}-i}d_2$ that is clean with respect to $F_{I,i}$ and $\mathcal{F}_{I,i}$ such that every long jump over $W_{i}$ in $D_{i}$ contains a vertex of some tile in $\mathcal{F}_{I,j}\setminus\mathcal{F}_{I,j-1}$ for every $j\in[i]$.
\end{itemize}

To be able to find a slice of width $d_2$ in the end, we therefore must fix
\begin{align*}
	d_1\geq \Brace{\Brace{2^{12}\Brace{w+1}^3\Fkt{f_P}{3t}+1}\Brace{2^6\Brace{w+1}^3\Fkt{f_P}{3t}+1}}^{2^{11+8\Fkt{f_P}{3t}}\Fkt{f_P}{3t}}d_2,
\end{align*}
and we further assume $d_2\geq 2^{16}\Fkt{f_W}{3t}^2$ to make sure we can apply \cref{lemma:longjumps1} and \cref{lemma:longjumpspahse2} in every round.
Note that this is not our final lower bound on $d_2$, just an intermediate assumption.

Next we describe the steps we perform in every round.
Let $i\in[2^{11+8\Fkt{f_P}{3t}}\Fkt{f_P}{3t}]$ and suppose we are given sets $F_{I,i-1}$ and $\mathcal{F}_{I,i-1}$, and graphs $D_{i-1}$ and $W_{i-1}$ as required by the input conditions for round $i$.

\emph{Step I}:
Let $k_i\coloneqq\Brace{\Brace{2^{12}\Brace{w+1}^3\Fkt{f_P}{3t}+1}\Brace{2^6\Brace{w+1}^3\Fkt{f_P}{3t}+1}}^{2^{11+8\Fkt{f_P}{3t}}\Fkt{f_P}{3t}-i+1}d_2$.
For each of the two possible parametrisations of $W$, we consider for every possible choice of $\xi,\xi'\in[w+1]$, the \hyperref[def:tiling]{tiling} $\Tiling\coloneqq\Tiling_{W_{i-1},k_i,w,\xi,\xi'}$ together with its four colouring $\Set{\Class_1,\dots,\Class_4}$.
For each $j\in[4]$, we consider the smallest slice $W'$ of $W$ that contains all vertices which belong to some tile of $\Tiling$.
Consider the auxiliary digraph $\Fkt{D_j^1}{W'}$ of type I obtained from $D_{i-1}$ and define the sets
\begin{align*}
	X_{\text{I}}^{\text{out}}&\coloneqq\CondSet{x_T^{\text{out}}}{T\in\Class_j}\text{, and}\\
	X_{\text{I}}^{\text{in}}&\coloneqq\CondSet{x_T^{\text{in}}}{T\in\Class_j}.
\end{align*}
Additionally we construct the set $Y_{\text{I}}$ as follows:
Let $Q$ and $Q'$ be the two cycles of $\Perimeter{W_0}$.
For every $j\in\left[ \frac{3d_1}{4} \right]$, $Y_{\text{I}}$ contains exactly one vertex of $Q\cap P^1_{4j}$, $Q\cap P^2_{4j+2}$, $Q'\cap P^1_{4j}$, and $Q'\cap P^2_{4j+2}$ each.
Then remove all vertices of $Y_{\text{I}}$ that do not belong to $D_{i-1}$.
Note that, by choice of $d_1$ and the bound on $F_{I,i-1}$, this does not decrease the size of $Y_{\text{I}}$ dramatically.

Then, if there is a family of $2^7\Fkt{f_P}{3t}$ pairwise disjoint directed $X_{\text{I}}^{\text{out}}$-$Y_{\text{I}}$-paths in $\Fkt{D_j^1}{W'}$, \cref{lemma:longjumps1} provides us with the existence of an integral packing of $t$ even dicycles and we are done.
So we may assume that there does not exist such a family and thus we may find a set $Z_1\subseteq\V{\Fkt{D_j^1}{W'}}$ of size at most $2^7\Fkt{f_P}{3t}$ that meets all such paths by \cref{thm:directedlocalmenger}.
With a similar argument, we either find an integral packing of $t$ even dicycles, or a set $Z_2\subseteq\V{\Fkt{D_j^1}{W'}}$ of size at most $2^7\Fkt{f_P}{3t}$ that meets all directed $Y_{\text{I}}$-$X_{\text{I}}^{\text{in}}$-paths in $\Fkt{D_j^1}{W'}$.
Let $\pi\in[2]$ indicate which of the two parametrisations of $W$ we are currently considering.
Then we can define the following two sets:
\begin{align*}
	Z_{\pi,\xi,\xi',j}\coloneqq& \Brace{Z_1\cup Z_2}\cap\V{D}\text{, and}\\
	\mathcal{Z}_{\pi,\xi,\xi',j}\coloneqq& \CondSet{T\in\Tiling}{\Brace{\V{T}\cup\Set{x_T^{\text{out}},x_T^{\text{in}}}}\cap\Brace{Z_1\cup Z_2}\neq\emptyset}.
\end{align*}
Note that $\max\Set{\Abs{Z_{\pi,\xi,\xi',j}},\Abs{\mathcal{Z}_{\pi,\xi,\xi',j}}}\leq 2^8\Fkt{f_P}{3t}$.

If we do not find an integral packing of $t$ even dicycles at any point, the sets $Z_{\pi,\xi,\xi',j}$ and $\mathcal{Z}_{\pi,\xi,\xi',j}$ are well defined for every possible choice of $\pi\in[2]$, $\xi,\xi'\in[w+1]$, and $j\in[4]$.
We collect these sets into the following to sets, which make up the first step towards the output of round $i$:
\begin{align*}
	F'_{I,i}\coloneqq& \bigcup_{\pi\in[2]}\bigcup_{\xi,\xi'\in[w+1]}\bigcup_{j\in[4]}Z_{\pi,\xi,\xi',j}\text{, and}\\
	\mathcal{F}'_{I,i}\coloneqq& \bigcup_{\pi\in[2]}\bigcup_{\xi,\xi'\in[w+1]}\bigcup_{j\in[4]}\mathcal{Z}_{\pi,\xi,\xi',j}.
\end{align*}
Consequently we have that $\max\Set{\Abs{F'_{I,i}},\Abs{\mathcal{F}'_{I,i}}}\leq 2^{11}\Brace{w+1}^2\Fkt{f_P}{3t}$.

We may now find a clear slice $W'_i\subseteq W_{i-1}$ of width $\Brace{2^{12}\Brace{w+1}^3\Fkt{f_P}{3t}+1}^{2^{11+8\Fkt{f_P}{3t}}\Fkt{f_P}{3t}-i}\Brace{2^6\Brace{w+1}^3\Fkt{f_P}{3t}+1}^{2^{11+8\Fkt{f_P}{3t}}\Fkt{f_P}{3t}-i+1}d_2$ which does not contain a marked vertex.
Note that we lose the additional factor of $2\Brace{w+1}$ since we remove whole tiles of width $w$ from $W_{i-1}$.
Let $D_i'\coloneqq D_{i-1}-F'_{I,i}$.
This concludes Step I of round $i$.

\begin{claim}\label{claim:stepIclaim}
	Every long jump $J$ over $W'_i$ in $D'_i$ whose endpoints belong to tiles of different colour contains a vertex of a tile from $\mathcal{F}_{I,j}\setminus \mathcal{F}_{I,j-1}$ for every $j\in[i-1]$, and it contains a vertex of a tile from $\mathcal{F}'_{I,i}$.
\end{claim}

\emph{Proof of \Cref{claim:stepIclaim}:} Suppose $J$ is also a long jump over $W_{i-1}$ in $D_{i-1}$, then, as $J$ still exists in $D'_i$, the tile in whose centre $J$ starts, or the tile in whose centre $J$ ends for some choices of $\pi\in[2]$, $\xi,\xi'\in[w+1]$, and $j\in[4]$, must be marked and therefore cannot belong to $W'_i$.

Hence $J$ must contain some vertex of $W_{i-1}$ as an internal vertex.
Let $T_{\mathsf{s}}$ be the tile of $W'_i$ in whose centre $J$ starts, and let $T$ be the first tile from the same tiling of $W_{i-1}$, that $J$ meets after $T_{\mathsf{s}}$.
Let $J'$ be the shortest subpath of $J$ with endpoints in $T_{\mathsf{s}}$ and $T$.
Then $J'$ is a long jump over $W_{i-1}$ in $D_{i-1}$.
Therefore, by our assumptions on the input of the $i$th round of Phase I, the first part of our claim is satisfied.
Moreover, if $T$ has a different colour than $T_{\mathsf{s}}$, $T$ must be marked.
So suppose $T$ has the same colour as $T_{\mathsf{s}}$.
Nonetheless, since $T_{\mathsf{s}}$ and the tile $T_{\mathsf{t}}$ which contains the endpoint of $J$ in the current tiling have different colours, $J$ contains a directed subpath $J''$ which is a long jump over $W_{i-1}$ and attaches to tiles of different colour.
Hence our claim follows.\hfill$\blacksquare$

With this we are ready for Step II of round $i$.

\emph{Step II}:
For this step let $k_i\coloneqq \Brace{2^{12}\Brace{w+1}^3\Fkt{f_P}{3t}+1}^{2^{11+8\Fkt{f_P}{3t}}\Fkt{f_P}{3t}-i}\Brace{2^6\Brace{w+1}^3\Fkt{f_P}{3t}+1}^{2^{11+8\Fkt{f_P}{3t}}\Fkt{f_P}{3t}-i+1}d_2$.
We are mainly concerned with the digraph $D'_i$.
In Step II it suffices to fix one parametrisation of $W$ since the construction of the \hyperref[def:auxiliarydigraphII]{type II auxiliary digraph} leaves the complete interior of same-colour-tiles intact instead of only their \hyperref[def:tile]{centres}.
For every pair of $\xi,\xi'\in[w+1]$ let us consider the \hyperref[def:tiling]{tiling} $\Tiling\coloneqq\Tiling_{W'_i,k_i,w,\xi,\xi'}$ together with its four colouring $\Set{\Class_1,\dots,\Class_4}$.
Then, for every $j\in[4]$ we call upon \cref{lemma:longjumpspahse2} which either provides us with an integral packing of $t$ even dicycles, and therefore closes the proof, or it produces two sets
\begin{align*}
	Z^{2}_{\xi,\xi',j}\subseteq&\V{D'_i}\text{ of size at most }2^3\Fkt{f_P}{3t}\text{, and}\\
	\mathcal{Z}^{2}_{\xi,\xi',j}\subseteq&\Tiling\text{ of size at most }2^3\Fkt{f_P}{3t},
\end{align*}
such that every directed $\V{W'_i}$-path whose endpoints belong to different tiles of $\Class_j$, contains a vertex of some tile in $\mathcal{Z}^{2}_{\xi,\xi',j}$.
This allows us to form the sets for the second step of round $i$:
\begin{align*}
	F''_{I,i}\coloneqq& \bigcup_{\xi,\xi'\in[w+1]}\bigcup_{j\in[4]}Z^2_{\xi,\xi',j}\text{, and}\\
	\mathcal{F}''_{I,i}\coloneqq& \bigcup_{\xi,\xi'\in[w+1]}\bigcup_{j\in[4]}\mathcal{Z}^2_{\xi,\xi',j}.
\end{align*}
As a result we obtain $\max\Set{\Abs{F''_{I,i}},\Abs{\mathcal{F}''_{I,i}}}\leq 2^5\Brace{w+1}^2\Fkt{f_P}{3t}$, and we are able to produce the two sets which will be passed on to the next round.
\begin{align*}
	F_{I,i}\coloneqq& F'_{I,i}\cup F''_{I,i}\cup F_{I,i-1}\text{, and}\\
	\mathcal{F}_{I,i}\coloneqq& \mathcal{F}'_{I,i}\cup \mathcal{F}''_{I,i}\cup \mathcal{F}_{I,i-1}.
\end{align*}
The bounds on $F_{I,i}$ and $\mathcal{F}_{I,i}$ follow immediately from the bounds on $F'_{I,i}$ and $F''_{I,i}$, $\mathcal{F}'_{I,i}$ and $\mathcal{F}''_{I,i}$, and the assumptions on the input of round $i$ respectively.

The pigeon hole principle allows us to find a clear slice $W_i\subseteq W'_i$ of width $\Brace{\Brace{2^{12}\Brace{w+1}^3\Fkt{f_P}{3t}+1}\Brace{2^6\Brace{w+1}^3\Fkt{f_P}{3t}+1}}^{2^{11+8\Fkt{f_P}{3t}}\Fkt{f_P}{3t}-i}d_2$ which does not contain a marked vertex.
Similar to Step I, we lose the additional factor of $2\Brace{w+1}$ since we remove whole tiles of width $w$ from $W'_i$.
Let $D_i\coloneqq D'_i-F''_{I,i}$.
This concludes Step II of round $i$.

\begin{claim}\label{claim:stepIIclaim}
	Every long jump over $W_i$ in $D_i$ contains a vertex of some tile in $\mathcal{F}_{I,j}\setminus\mathcal{F}_{I,j-1}$ for every $j\in[i]$.
\end{claim}

\emph{Proof of \Cref{claim:stepIIclaim}:}
Let $J$ be a long jump over $W_i$ in $D_i$, and let $\Tiling$ be a tiling of $W_0$ defined by $w$ and some $\xi,\xi'\in[w+1]$ such that $J$ starts at the centre of some tile $T_{\mathsf{s}}\in\Tiling$.
Suppose all tiles of $\Tiling$ which contain vertices of $J$ belong to the same colour.
Then $J$ must have existed during the corresponding part of Step II of round $i$ and thus either $T_{\mathsf{s}}$ or $T_{\mathsf{t}}\in\Tiling$, which is the tile that contains the endpoint of $J$, must have been newly marked in this round.
Therefore $J$ must contain at least one tile of a colour different from the one of $T_{\mathsf{s}}$.
Moreover, we may assume $T_{\mathsf{s}}$ and $T_{\mathsf{t}}$ to be of the same colour as otherwise we would be done by \cref{claim:stepIclaim}.
Next suppose $J$ is also a long jump over $W_{i-1}$, then again $J$ would have been considered during Step II as a long jump connecting two tiles of the same colour and thus $T_{\mathsf{s}}$ or $T_{\mathsf{t}}$ would have been newly marked in this round.
Therefore $J$ must contain a vertex of some tile from $W_{i-1}$.
Let $J'$ be a shortest subpath from $T_{\mathsf{s}}$ to some tile $T$ of $W_{i-1}$, then $J'$ is a long jump over $W_{i-1}$ and thus $J$ contains a vertex of some tile of $\mathcal{F}_{I,j}\setminus\mathcal{F}_{I,j-1}$ for every $j\in[i-1]$ by our assumptions on the input of round $i$.
If $T$ has a different colour than $T_{\mathsf{s}}$, then $T$ would have been marked in Step I of round $i$, and if $T$ shares the colour of $T_{\mathsf{s}}$ it must have been marked in Step II of round $i$.
Either way our claim follows. \hfill$\blacksquare$

From \cref{claim:stepIIclaim} it follows that we satisfy all requirements for the output of round $i$ and thus round $i$ is complete.
We continue until we finish round $i=2^{11+8\Fkt{f_P}{3t}}\Fkt{f_P}{3t}$ and obtain the following four objects as its output:
\begin{itemize}
	\item a slice $W_I\coloneqq W_{2^{11+8\Fkt{f_P}{3t}}\Fkt{f_P}{3t}}$ of width $d_2$,
	\item a set $A\coloneqq F_{I,2^{11+8\Fkt{f_P}{3t}}\Fkt{f_P}{3t}}$ of size at most $(3t)^{28}2^{60+2^{10}(3t)^8}$,
	\item a digraph $D_I\coloneqq D_{2^{11+8\Fkt{f_P}{3t}}\Fkt{f_P}{3t}}=D-A$, and
	\item a sequence $\mathcal{F}_{I,1}\subseteq \mathcal{F}_{I,2}\subseteq\dots\subseteq \mathcal{F}_{I,2^{11+8\Fkt{f_P}{3t}}\Fkt{f_P}{3t}}$ such that every long jump over $W_I$ in $D_I$ contains a vertex of some tile in $\mathcal{F}_{I,i}\setminus\mathcal{F}_{I,i-1}$ for every $i\in[2^{11+8\Fkt{f_P}{3t}}\Fkt{f_P}{3t}]$.
\end{itemize}
This brings us to the final claim of Phase I.

\begin{claim}\label{claim:nolongjumps}
	If there is a long jump over $W_I$ in $D_I$, then there exists an integral packing of $t$ even dicycles in $D$.
\end{claim}

\emph{Proof of \Cref{claim:nolongjumps}:}
Let $J$ be a long jump over $W_I$ in $D_I$.
We fix a parametrisation of $W$, $\xi,\xi'\in[w+1]$, and $c\in[4]$ such that there exists a tile $T_{\mathsf{s}}\in\Tiling\coloneqq\Tiling_{W_0,d_1,w,\xi,\xi'}$ of colour $c$ whose centre contains the starting point of $J$.
Let $T_{\mathsf{t}}\in\Tiling$ be the tile that contains the endpoint of $J$.
As $J$ is a long jump, note that $T_{\mathsf{s}}\neq T_{\mathsf{t}}$.

We now create a family $\mathcal{L}_0$ of $2^{9+8\Fkt{f_P}{3t}}\Fkt{f_P}{3t}$ pairwise disjoint subpaths of $J$ with the following properties:
\begin{enumerate}
	\item for every $L\in\mathcal{L}_0$, let $T_{L,1}$ and $T_{L,2}$ be the tiles of $\Tiling$ that contain the starting point $s_L$ and the endpoint $t_L$ of $L$ respectively, then there exist distinct $i_{L,1},i_{L,2}\in[2^{11+8\Fkt{f_P}{3t}}\Fkt{f_P}{3t}]$ such that $s_L$ is a vertex of a tile from $\mathcal{F}_{I,i_{L,1}}\setminus\mathcal{F}_{I,i_{L,1}-1}$, and $t_L$ is a vertex of some tile in $\mathcal{F}_{I,i_{L,2}}\setminus\mathcal{F}_{I,i_{L,2}-1}$, and
	\item if $L,L'\in\mathcal{L}_0$ are distinct, then $\Set{i_{L,1},i_{L,2}}\cap\Set{i_{L',1},i_{L',2}}=\emptyset$.
\end{enumerate}
Let us initialise $\mathcal{I}_0\coloneqq[2^{11+8\Fkt{f_P}{3t}}\Fkt{f_P}{3t}]$ and for every subset $\mathcal{I}'\subseteq\mathcal{I}_0$, we define the family $\mathcal{F}_{\mathcal{I}'}\coloneqq \bigcup_{i\in\mathcal{I}'}\mathcal{F}_{I,i}\setminus\mathcal{F}_{I,i-1}$.
Please note that every internal vertex of $J$ that belongs to $W$ must belong to some tile from $\mathcal{F}_{I,2^{11+8\Fkt{f_P}{3t}}\Fkt{f_P}{3t}}$, since otherwise we could find a directed path from the centre of $T_{\mathsf{s}}$ to the perimeter of $W_0$ contradicting the construction in Step I of Phase I, or we would have a directed path between two tiles of the same colour, where both of them are unmarked.
This second outcome contradicts the construction in Step II of Phase I.

Consider the shortest subpath of $J$ that shares its starting point with $J$ and is a long jump over $W$.
Let $T_{L_1,1}$ be the tile of $\Tiling$ where this path ends and let $s_{L_1}$ be the first vertex of $J$ for which its successor along $J$ does not belong to $T_{L_1,1}$.
Note that there exists $i_{L_1,1}\in\mathcal{I}_0$ such that $T_{L_1,1}\in\mathcal{F}_{I,i_{L_1,1}}\setminus \mathcal{F}_{I,i_{L_1,1}-1}$ by the discussion above.
Let $L_1$ be the shortest subpath of $J$ that starts in $s_{L_1}$ and ends in a vertex $t_{L_1}$ for which $i_{L_1,2}\in\mathcal{I}_0\setminus\Set{i_{L_1,1}}$ exists such that $t_{L_1}$ belongs to a tile from $\mathcal{F}_{I,i_{L_1,2}}\setminus \mathcal{F}_{I,i_{L_1,2}-1}$.
Let $T_{L_1,2}$ be the tile of $\Tiling$ that contains $t_{L_1}$.
We add $L_1$ to $\mathcal{L}_0$ and set $\mathcal{I}_1\coloneqq \mathcal{I}_0\setminus\Set{i_{L_1,1},i_{L_1,2}}$.
By our choice of $i_{L_1,1}$ and $i_{L_1,2}$, the path $t_{L_1}J$ still contains a vertex of a tile from $\mathcal{F}_{I,j}\setminus\mathcal{F}_{I,j-1}$ for every $j\in\mathcal{I}_1$.

Now let $q\in[2,2^{9+8\Fkt{f_P}{3t}}\Fkt{f_P}{3t}]$ and assume that the paths $L_1,\dots,L_{q-1}$ together with the tiles, indices and the set $\mathcal{I}_{q-1}$ have already been constructed.
Follow along $J$, starting from $t_{L_{q-1}}$, until the next time we encounter the last vertex $s_{L_q}$ of some tile from $\mathcal{F}_{\mathcal{I}_{q-1}}$ before $J$ leaves said tile again.
Let $T_{L_q,1}\in\Tiling$ be the tile that contains $s_{L_q}$, and let $i_{L_q,1}\in\mathcal{I}_{q-1}$ be the integer such that $s_{L_q}$ belongs to a tile of $\mathcal{F}_{I,i_{L_q,q}}\setminus \mathcal{F}_{I,i_{L_q,q}-1}$.
Then let $L_q$ be the shortest subpath of $J$ that starts in $s_{L_q}$ and ends in a vertex $t_{L_q}$ which belongs to a tile from $\mathcal{F}_{I,i_{L_q,2}}\setminus \mathcal{F}_{I,i_{L_q,2}-1}$, where $i_{L_q,2}\in\mathcal{I}_{q-1}\setminus\Set{i_{L-q,1}}$.
We choose $T_{L_q,2}\in\Tiling$ to be the tile that contains $t_{L_q}$ and set $\mathcal{I}_q\coloneqq \mathcal{I}_{q-1}\setminus\Set{i_{L_q,1},i_{L_q,2}}$.
As before notice that $t_{L_q}J$ still contains a vertex from some tile in $\mathcal{F}_{I,j}\setminus\mathcal{F}_{I,j-1}$ for every $j\in\mathcal{I}_q$.
Add $L_q$ to $\mathcal{L}_0$.

With every iteration we remove exactly two members from $\mathcal{I}$ and, as $\Abs{\mathcal{I}}=2^{11+8\Fkt{f_P}{3t}}\Fkt{f_P}{3t}$, this means that by the time we reach some $q$ for which $\mathcal{I}_q=\emptyset$, we have indeed constructed $2^{10+8\Fkt{f_P}{3t}}\Fkt{f_P}{3t}$ paths as required.

Now there must exist $c'\in[4]$ and a family $\mathcal{L}_1\subseteq \mathcal{L}_0$ of size $2^{8+8\Fkt{f_P}{3t}}\Fkt{f_P}{3t}$ such that each path $L\in\mathcal{L}_1$ has at least one endpoint in $\Class_{c'}$.
And, in an immediate second step, we can find a family $\mathcal{L}_2\subseteq\mathcal{L}_1$ of size $2^{7+8\Fkt{f_P}{3t}}\Fkt{f_P}{3t}$ such that every path in $\mathcal{L}_2$ starts in a tile of $\Class_{c'}$, or every path in $\mathcal{L}_2$ ends in a tile of $\Class_{c'}$.
Without loss of generality we may assume that every path in $\mathcal{L}_2$ starts in a tile of $\Class_{c'}$ since the other case follows with similar arguments.

Let $\widetilde{W}'$ be the smallest slice of $W$ that contains all tiles from $\Class_{c'}$, but no tile from $\Class_{c'}$ meets the perimeter of $\widetilde{W}$.
Then let $\widetilde{\Tiling}\coloneqq\TierIITiling{\Tiling,c',w}{\widetilde{W}'}$ be the \hyperref[def:tilingII]{tier II tiling} of $\widetilde{W}\coloneqq\InducedSubgraph{\widetilde{W}'}{\Tiling,c',w}$.
Since the paths in $\mathcal{L}_2$ are pairwise disjoint, we can extend each $L\in\mathcal{L}_2$ such that it starts on the centre of the tile of $\widetilde{\Tiling}$ which encloses its endpoint in $W$, while making sure that the resulting family of paths is still at least half-integral.
Similarly, wherever necessary, we may extend the paths through $W$ such that each of them also ends in a tile of $\widetilde{\Tiling}$.
Indeed, we can even guarantee that the endpoints of the resulting paths are mutually at $\widetilde{W}$-distance at least $4$.
Let $\mathcal{L}_3$ be the resulting half integral linkage.

Next consider the four colouring $\Set{\widetilde{\Class_1},\dots,\widetilde{\Class_4}}$ of $\widetilde{\Tiling}$.
Then there exists $\widetilde{c}\in[4]$ and a family $\mathcal{L}_4\subseteq\mathcal{L}_3$ of size $2^{5+8\Fkt{f_P}{3t}}\Fkt{f_P}{3t}$ such that every path in $\mathcal{L}_4$ starts at the centre of some tile from $\Class_{\widetilde{c}}$.
It follows from the construction of $\mathcal{L}_0$ that no two paths in $\mathcal{L}_4$ start in the same tile.

By a similar argument, there exists a family $\mathcal{L}_5\subseteq\mathcal{L}_4$ of size $2^{4+8\Fkt{f_P}{3t}}\Fkt{f_P}{3t}$ such that either none, or all paths in $\mathcal{L}_5$ end in tiles of $\widetilde{\Class_{\widetilde{c}}}$.

In the first case we can extend every path in $\mathcal{L}_5$ towards the perimeter of $\widetilde{W}$ such that the resulting family $\mathcal{L}_6$ of paths remains at worst half-integral, and the endpoints of the resulting paths are mutually at $\widetilde{W}$-distance at least $4$.
Now \cref{thm:halfintegral} provides us with a family $\mathcal{L}_7$ of size $2^{3+8\Fkt{f_P}{3t}}\Fkt{f_P}{3t}$ such that $\V{\mathcal{L}_7}\subseteq\V{\mathcal{L}_6}$, and the paths in $\mathcal{L}_7$ are pairwise vertex disjoint.
Hence \cref{lemma:longjumps1} yields the existence of an integral packing of $t$ even dicycles and our claim follows.

In the second case we consider two subcases.
Let $\mathcal{X}$ be the family of all tiles of $\widetilde{\Tiling}\setminus\widetilde{\Class_{\widetilde{c}}}$ which contain an internal vertex of some path in $\mathcal{L}_5$ but no endpoint of any path in $\mathcal{L}_5$.

Suppose $\Abs{\mathcal{X}}\geq 2^8\Fkt{f_P}{3t}$ and recall that, if $L$ and $P$ are directed paths, we say that $P$ is a \emph{long jump of $L$} if $P$ is a $w$-long jump over $W$ and $P\subseteq L$.
We also say that $P$ is a \emph{jump of $L$}, if $P$ is a directed $\V{W}$-path.
Then we can use the technique from the first part of the proof of \cref{lemma:longjumps1} to construct a half-integral family $\mathcal{L}_6$ such that
\begin{enumerate}
	\item $\Abs{\mathcal{L}_6}=2^{4+8\Fkt{f_P}{3t}}\Fkt{f_P}{3t}$, and
	\item for every $L\in\mathcal{L}_6$, every endpoint $u$ of a jump of $L$ with $u\in\V{\widetilde{W}}$ belongs to a tile from $\widetilde{\Class_{\widetilde{c}}}\cup\mathcal{X}$.
\end{enumerate}
We can then apply \cref{thm:halfintegral} to obtain a family $\mathcal{L}_7$ of size $2^{3+8\Fkt{f_P}{3t}}\Fkt{f_P}{3t}$ with $\V{\mathcal{L}_7}\subseteq\V{\mathcal{L}_6}$ such that the paths in $\mathcal{L}_7$ are pairwise disjoint and link the same two sets of vertices as the paths in $\mathcal{L}_6$ do.
Finally \cref{lemma:longjumps1} yields the existence of an integral packing of $t$ even dicycles.

So we may assume $\Abs{\mathcal{X}}<2^8\Fkt{f_P}{3t}$.
In this case, we may find a subwall $\widetilde{W}'$ of $\widetilde{W}$ of order $d_1-2^8\Fkt{f_P}{3t}\Brace{2w+1}$ that does not contain a vertex of any tile in $\mathcal{X}$ by removing, for every tile $T\in\mathcal{X}$, all edges and vertices of the horizontal cycles and vertical paths of $T$, which are not used by other cycles or paths.
For each tile we remove during this procedure, we remove a row and a column of tiles and thereby might reduce the number of distinct tiles which contain starting vertices of paths in $\mathcal{L}_5$ by a factor of $\frac{1}{2}$.
However, since $\Abs{\mathcal{X}}<2^8\Fkt{f_P}{3t}$, we can still find, after potentially expanding the start and end sections of some paths to again reach the slightly shifted perimeters of their tiles, a half-integral family $\mathcal{L}_6$ of size $>2^4\Fkt{f_P}{3t}$ of paths that start and end in tiles of $\widetilde{\Class_{\widetilde{c}}}$ and that are otherwise disjoint from $\widetilde{W}'$.
By using \cref{thm:halfintegral} we can transform this family into an integral family $\mathcal{L}_7$ of size $>2^3\Fkt{f_P}{3t}$ and thus an application of \cref{lemma:longjumpspahse2} yields an integral packing of $t$ even dicycles. 

Concluding Phase I, \cref{claim:nolongjumps} either yields an integral packing of $t$ even dicycles and therefore closes the proof, or $W_I$ is in fact clean, meaning that $W_I$ has no long jump in $D_I$.
Consequently we may bound the function $\ApexBound$ from the statement of \cref{thm:oddwall} as follows:
\begin{align*}
	\Fkt{\ApexBound}{t}\leq (3t)^{28}2^{61+2^{10}(3t)^8}.
\end{align*}

\paragraph{Phase II}

The second phase is a fairly easy and short argument.
Essentially, with $W_I$ we have found a wall of still sufficient size but without any long jumps.
Therefore me may now find $t$ slices of $W_I$, mutually still far enough apart from each other within $W_I$, and can ask if among them there is one without an even dicycle.
If the answer is yes we are done and if the answer is no, the absence of long jumps ensures that these $t$ even dicycles are pairwise disjoint and thus we have yet again found an integral packing of $t$ even dicycles.
The only technical part that remains is to provide sufficient definitions for these slices and their mutual distance.

To meet the requirements from Phase I and have enough space left in $W_I$, let us make the following assumption:
\begin{align*}
	d_2\geq t\Brace{r+6+2^{32+(3t)^{30}}}.
\end{align*}

Let us partition $W_I$ into $3t$ many slices as follows:
First we can partition $W_I$ into $t$ slices $S_i$ of width $r+2^{32+(3t)^{30}}$.
Each $S_i$ is then partitioned into a slice $H_i$ of width $r+6+2^{31+(3t)^{30}}$ that contains the left perimeter cycle of the of $S_i$, and a slice $G_i$ of width $2^{31+(3t)^{30}}$ containing the right cycle of the perimeter of $S_i$.
For every $i\in[t]$ we may now further partition $H_i$.
Let $N_{i,L}\subseteq H_i$ be the slice of width $2^{30+t^{30}}$ containing the left cycle of $\Perimeter{H_i}$, let $N_{i,R}$ be the slice of width $2^{30+t^{30}}$ containing the right cycle of $\Perimeter{H_i}$, and let $N_i'\coloneqq H_i-N_{L,i}-N_{R,i}$ be the remaining slice of width $r+6$.
Finally let $N_i$ be the slice obtained from $N_i'$ be removing the three left most and the three right most vertical cycles.
Then $N_i$ is a slice of width $r$.

Let us first assume that for every $i\in[t]$ there exists a directed path $P_i$ with one endpoint in $N_i$, the other endpoint in $W_I-N_i'$, and which is internally disjoint from $W_I$.
In this case there exists a slice $X$ of width three contained in one of the two components of $N_i'-N_i$ which separates the two endpoints of $P_i$ in $W_I$.
Since there does not exist a long jump over $W_I$ in $D_I$ we can further be sure that there exists $Y\in\Set{L,R}$ such that the endpoint of $P_i$ not in $N_i$ must lie in $N_{Y,i}$.
We can now extend the path $P_i$ within $S_i$ until it forms a perimeter jump over $X$.
By \cref{lem:gridplusjumpmakesoddbicycle} and \cref{obs:oddbicycleevendicycle} this means that we can find an even dicycle locally bound to $S_i$.
Let $J\subseteq[t]$ be the set of integers where this is the case.
If $\Abs{J}=t$ we have found an integral packing of $t$ even dicycles and are done.
Hence there exists the set $J'\coloneqq[t]\setminus J$.

Now suppose for each $i\in J'$ the strong component of $D_I-\Perimeter{N_i'}$ that contains $N_i$ has an even dicycle.
Since there is no long jump over $W_I$ these components are pairwise disjoint and also disjoint from the odd bicycles found in the $S_j$, $j\in J$, and thus we have again found an integral packing of $t$ even dicycles.
Thus there must exist some $i\in J'$ for which $N_i'$ has no even dicycle.
In particular this means that the component of $D_I-\Perimeter{N_i}$ containing the remaining vertices of $N_i-\Perimeter{N_i}$ must be free of even dicycles.
Hence $N_i$ is the desired odd wall of order $r$ in $D_I=D-A$ and the proof is complete.
At last let us combine all assumptions on the $d_j$ to obtain the following bound on $\OddWallOrder{t,r}$:
\[ \OddWallOrder{t,r}\leq\Brace{2^{140}t^{72}}^{2^{10}t^8+2^{12}}\Brace{r+6+2^{32+t^{30}}} . \qedhere \] 
\end{proof}

Please note that this proof is very close to the one for the original directed flat wall theorem.
Moreover, asking for an integral packing of even dicycles might even be a bit overkill.
It should be possible to significantly reduce the order of the functions $\OddWallApexNoArg $ and $\OddWallOrderNoArg $ through more careful analysis, especially if we are content with a half- or quarter-integral packing of even dicycles.

Finally, we briefly argue that the process of finding the odd wall can be performed efficiently.
Recently Campos, Lopes, Maia, and Sau \cite{campos2022adapting} provided an adaptation of the proofs in \cite{KawarabayashiKreutzer2015DirectedGrid} and \cite{Johnson2001DirectedTreewidth}, leading to an polynomial time algorithm for a given $k$ that either finds a directed tree decomposition of small width in terms of $k$, or a large cylindrical $k$-wall.

\begin{theorem}[Campos, Lopes, Maia, and Sau \cite{campos2022adapting}]\label{thm:fptdirectedwall}
    There exist computable functions $f, g \colon \N \rightarrow \N$ such that for every $k \in \N$ and every digraph $D$, one can either find a cylindrical $k$-wall in $D$ or a directed tree decomposition of width at most $f(k)$ for $D$ in time $g(k) n^{\mathcal{O}(1)}$.
\end{theorem}

For later use, we can thus combine \Cref{cor:directedwall}, \Cref{thm:fptdirectedwall}, and \Cref{thm:oddwall} into a single statement.

\begin{corollary}\label{cor:findoddwall}
    There exist computable functions $\FindOddWallOrderNoArg \colon \N\times\N\rightarrow\N$, $\FindOddWallApexNoArg \colon \N\rightarrow\N$, and $f \colon \N\times\N\rightarrow\N$, and an algorithm such that for all integers $r,t \geq 1$ and all digraphs $D$ we can find an integral packing of $t$ even dicycles in $D$, a directed tree decomposition of $D$ with width less than $\FindOddWallOrder{r,t}$, or a set $A \subseteq \V{D}$ with $\Abs{A}\leq \FindOddWallApex{t}$ and an $r$-wall $W'\subseteq W-A$, which is an odd wall under $A$, in time at most $f(t,r) n^{\mathcal{O}(1)}$.
\end{corollary}

\section{Shifting even dicycles}\label{sec:shifting}

In this section we prove five key results for many arguments to come, namely \Cref{lem:planarshifting}, \Cref{lem:planartraceshifting}, \Cref{lem:circlenonplanarshifting}, \Cref{lem:diamondnonplanarshifting}, and \Cref{lem:refineddiamondnonplanarshifting}.
These statements will allow us to shift the occurrence of even dicycles into specific areas, if we are provided enough infrastructure.
With \Cref{lem:detour}, we will also prove a convenient utility lemma that lets us reroute a linkage through another linkage if both of them live in a strong planar decomposition.

\subsection{Shifting a planar even dicycle}\label{subsec:planarshifting}

Finding an even dicycle will often not be enough on its own, and we will have to shift this dicycle into a specific part of the graph.
These parts of the graph will be separated by a collection of concentric odd dicycles in the upcoming arguments, which will allow us to use one of the \emph{shifting lemmas} we present in the remainder of this section.

\begin{lemma}[Planar Shifting]\label{lem:planarshifting}
	Let $C_{e}$ be an even and $C_{o}$ be an odd dicycle, together with a drawing of $C_{e} \cup C_{o}$ in the plane.
	There exists an even dicycle $C' \subset C_{e} \cup C_{o}$ which is drawn in the closure of one of the two disks bounded by $C_{o}$.
\end{lemma}
\begin{proof}
	Let $C_e$ and $C_o$ be chosen such that they form a smallest counterexample to our statement.
	It is easy to see that the drawing of $C_{e} \cup C_{o}$ must have at least one face whose boundary is a dicycle $C$, since we can build the graph from a dicycle using directed ears.
	If $C$ itself is even, we are already done, since $C$ bounds a face and thus lies on one side of $C_{o}$.
	
	Thus, we can assume that $C$ is odd and furthermore $C \neq C_{o}$, as we would otherwise again be done already.
	Accordingly, $C$ can be partitioned into paths $P_1, Q_1, P_2, Q_2, \ldots , P_k, Q_k$, with $k \geq 1$, such that each $P_i$ is maximal with $P_i \subseteq C \cap C_o$ and the tail of $P_i$ is the head of $Q_i$, whose tail is in turn the head of $P_{i+1}$.
	Note that since $C_{e} \cup C_{o}$ is a drawing on the plane, for each $i \in [k]$ there exists a subpath $Q'_i$ of $C_o$ with the same endpoints as $Q_i$, which is edge-disjoint and internally vertex-disjoint from $C$, such that $C$ is drawn on one side of $\bigcup_{i \in [k]} (P_i \cup Q'_i)$.
	Since $C_o$ is a dicycle, $Q'_i$ and $Q'_j$ do not intersect for $i \neq j$ and due to the maximality of $P_i$ the path $Q_i$ has length at least one.
	
	We now observe that if there exist $Q_i$ and $Q'_i$ of the same parity, then we can delete the interior of $Q_i$, disrupting $C_e$, and replace this part of the even dicycle with $Q'_i$, yielding a smaller even cycle $C'_e$.
	Due to the minimality of $C_{e} \cup C_{o}$ there now exists an even cycle $C'$ which lies entirely on one side of $C_o$.
	Thus, all $Q_i$ and $Q'_i$ must have differing parities.
	
	If we now suppose that $k \geq 2$, then we can replace both the interior of $Q_1$ and $Q_2$ with $Q'_1$ and $Q'_2$ respectively, again creating another smaller even dicycle, which yields our desired outcome.
	Thus, $k = 1$, but in this case it is easy to see that $C = C_o$, since all that distinguishes $C_o$ and $C_e$ is whether they use $Q_1$ or $Q'_1$.
	This leads to a contradiction and concludes our proof.
\end{proof}

We can in fact be a little more specific as to how we can find this even dicycle.

\begin{corollary}\label{cor:planarshifting}
    Let $C_{e}$ be an even and $C_{o}$ be an odd dicycle, together with a drawing of $C_{e} \cup C_{o}$ in the plane.
	There exists a path $P_e \subseteq C_e$, that is drawn in the closure of one of the two disks bounded by $C_{o}$, such that $C_{o} \cup P_e$ contains an even dicycle.
\end{corollary}
\begin{proof}
    According to \Cref{lem:planarshifting}, there exists an even dicycle $C \subseteq C_e \cup C_o$ that is drawn entirely in the closure of one of the two disks bounded by $C_o$.
    Let $d$ be the disk bounded by $C_o$ that $C$ is drawn on and let $\mathcal{P}$ be the set of components of $C_e \cap C$, namely all $C_o$-paths in $C$ that are drawn on $d$.
    Furthermore, we take the direction of $C_o$ to define a circular order on $V(C_o)$.
    Clearly, if $|\mathcal{P}| = 1$, we are done.
    We may therefore suppose that $C$ is chosen within $C_e \cup C_o$ such that the number of elements in $\mathcal{P}$ is minimal.

    Suppose that $|\mathcal{P}| > 1$ and note that if there exists a $u$-$v$-path $P \in \mathcal{P}$, then if the unique directed $v$-$u$-subpath $P'$ in $C_o$ forms a dicycle together with $P$.
    If $P' \cup P$ is even, we are again done.
    Thus we may assume that all such pairs of paths have opposing parities.
    However, if $P \in \mathcal{P}$ is a $u$-$v$-path such that the $u$-$v$-path $P''$ on $C_o$ does not contain the endpoints of any vertex of a path in $\mathcal{P} \setminus \{ P \}$, then we may replace $P$ with $P''$ in $C$, which results in an even dicycle contradicting the minimality of $\mathcal{P}$.
    Such a path must exist since $C_{e} \cup C_{o}$ have a planar drawing and in particular $C \cup C_o$ are drawn in a planar fashion on $d$.
\end{proof}

\subsection{A short primer on matching theory in bipartite graphs}\label{subsec:matchingprimer}

Our next goal will be to prove that we can shift even dicycles in a similar fashion to \Cref{lem:planarshifting} if we are instead provided a planar decomposition for both dicycles.
This will require an in-depth discussion of matching theory, as our main hurdle will lie in trying to understand how the drawing of our two dicycles behaves in small cycle sums, which correspond to a specific sum operation in the matching setting that is easier to analyse than small cycle sums.

The reader may be relieved to know that our explicit use of matching theory is constrained to this section.
We begin by establishing the link between bipartite graphs with perfect matchings and digraphs.

If $B$ is a bipartite undirected graph, we denote the two colour classes into which $\V{B}$ can be partitioned as $V_1,V_2$.
Should there be several valid choices for $V_1$ and $V_2$, we assume that one of them is chosen at random.

\begin{definition}\label{def:mdir}
	Let $B$ be a bipartite graph with a perfect matching $M = \{ a_1b_1, \ldots , a_{|M|}b_{|M|} \}$, with $a_i \in V_1$ and $b_i \in V_2$, for $i \in [|M|]$.
	The \emph{$M$-direction} $\DirM{B}{M}$ of $B$ is defined such that
	\begin{enumerate}
		\item $\V{\DirM{B}{M}} = \{ v_1, \ldots , v_{|M|} \}$ and
		
		\item $\E{\DirM{B}{M}} = \{ (v_i,v_j) \ | \ a_ib_j \in \E{B} \text{ and } i \neq j \}$.
	\end{enumerate}
	In the other direction, we call $B$ the \emph{split of $\DirM{B}{M}$}.
	Note that every digraph has a split.
\end{definition}

A graph $G$ with $|\V{G}| \geq 2k+2$ is called \emph{$k$-extendable} for a $k \in \N$ if for every subset $F \subseteq \E{G}$ such that $|F| \leq k$ there exists a perfect matching $M$ with $F \subseteq M$.
In particular, a $1$-extendable graph is also called \emph{matching covered}.

The following theorem is a central feature of the matching theory of bipartite graphs and, as McCuaig points out in \cite{mccuaig2000even} (see Theorem 16), it and extensions of this statement can be found in the literature since at least the 1960s.

\begin{theorem}[folklore]\label{thm:kexttostrcon}
    Let $B$ be a bipartite graph with a perfect matching $M$ and let $D$ be a digraph such that $D = \DirM{B}{M}$.
    Then $B$ is $k$-extendable if and only if $D$ is $k$-strongly connected.
\end{theorem}

Since matching theory concerns itself with graphs having perfect matchings, basic notions such as the definition of a minor in the undirected setting, have to be adjusted such that they preserve the existence of perfect matchings, or more specific properties such as being matching covered.
A \emph{bicontraction} at a vertex $v$ of degree two in an undirected graph $G$ is the result of contracting both edges incident to $v$.
We call a subgraph $H \subseteq G$ \emph{conformal} if $G - \V{H}$ has a perfect matching, and $H$ is called a \emph{matching minor} of $G$, if $H$ can be derived from a conformal subgraph of $G$ by repeated bicontractions.
Matching minors and butterfly minors are essentially the same notion expressed in different settings.

\begin{lemma}[McCuaig \cite{mccuaig2000even}]\label{lem:bminortomatminor}
	Let $B$ and $H$ be bipartite, matching covered graphs.
	$H$ is a matching minor of $B$ if and only if there exists a perfect matching $M$ of $B$ and a perfect matching $M'$ of $H$ such that $\DirM{H}{M'}$ is a butterfly minor of $\DirM{B}{M}$.
\end{lemma}

Another way in which bipartite graphs and their $M$-directions are tied together is their drawings.

\begin{observation}\label{obs:strongplanarityisplanarity}
        A bipartite graph $B$ with a perfect matching is planar if and only if $\DirM{B}{M}$ is strongly planar for every perfect matching $M$ of $B$.
\end{observation}

A central reason for the interest in matching theory is the problem of counting perfect matchings.
In the study of this problem, Kasteleyn \cite{kasteleyn1967graph} introduced a notion that allowed him to count perfect matchings on planar graphs in polynomial time.
We call an orientation $\Oriented{G}$ -- that is an assignment of directions to each edge -- of an undirected graph $G$ \emph{Pfaffian} if each even, conformal cycle in $G$ has an odd number of edges agreeing with either direction of traversal in $\Oriented{G}$.
A graph with a Pfaffian orientation is called \emph{Pfaffian}.

As we illustrated in \Cref{sec:notation}, Pfaffian graphs are very closely tied to our problem since they correspond to non-even digraphs.
For a survey on Pfaffian graphs see \cite{Thomas2006PfaffianSurvey} and for a very thorough treatment on the relation of several seemingly disparate problems in combinatorics and other fields of mathematics to the problem of recognizing Pfaffian graphs see \cite{McCuaig2004Polya}.

In the following theorem, the equivalence of the first two items is due to Little \cite{Little1975Convertible} and the equivalence of the second two items is a consequence of combining \Cref{thm:nonevendigraphs} and a result by McCuaig \cite{McCuaig2004Polya}.

\begin{theorem}\label{thm:pfafftooddbicycle}
	Let $B$ be a bipartite graph with a perfect matching $M$.
	The following statements are equivalent:
	\begin{enumerate}
		\item $B$ is not Pfaffian.
		
		\item $B$ contains $K_{3,3}$ as a matching minor.
		
		\item $\DirM{B}{M}$ contains a weak odd bicycle.
	\end{enumerate}
\end{theorem}

We will also need the matching theoretic analogues for directed tight cuts and dibraces.
In a matching covered graph $G$ with a perfect matching, an edge cut $\Cut{X}$ around a vertex set $X$ is called \emph{tight}, if for all perfect matchings $M$ of $G$ we have $|M \cap \Cut{X}| = 1$.
Let $\overline{X} = \V{G} \setminus X$.
Then $\ContractXinGtoV{X}{G}{c}$ and $\ContractXinGtoV{\overline{X}}{G}{c}$ are constructed by contracting $X$ and respectively $\overline{X}$ into a single vertex.
The resulting graphs are the two \emph{tight cut contractions} corresponding to $\Cut{X}$.
If $X$ or $\overline{X}$ contain a single vertex, we call the associated tight cut \emph{trivial}.
Another link between directed graphs and bipartite graphs is given by the following observations, which also illustrates that the term directed tight cut derives from the notions in matching theory.

\begin{observation}\label{obs:dirtightcutandtightcut}
    Let $B$ be a bipartite matching covered graph with a tight cut $\Cut{X}$.
    Then for every perfect matching $M$ of $B$ the vertex in $\DirM{B}{M}$ corresponding to the unique edge $e \in M \cap \Cut{X}$ induces a directed tight cut in $\DirM{B}{M}$.

    Vice versa, for any digraph $D$ with $D = \DirM{B'}{M'}$, any directed tight cut $(Y,Z)$ in $D$ induces a tight cut in $B'$.
\end{observation}

By repeatedly applying tight cut contractions for non-trivial tight cuts to a graph $G$ until we have reached only graphs that lack non-trivial tight cuts, we perform a \emph{tight cut decomposition}.
A bipartite graph without non-trivial tight cuts is called a \emph{brace} and any non-bipartite graph without non-trivial tight cuts is called a \emph{brick}.
As Lov{\'a}sz \cite{lovasz1987matching} showed when he introduced this procedure, the tight cut decomposition of any given graph is independent of the order in which we contract tight cuts.
Again, we note the connection between the terms in the directed setting and the matching setting.
In particular, thanks to \Cref{obs:dirtightcutandtightcut} we now know that the big vertices of a digraph $D$ in a plane decomposition will correspond to tight cuts in the split of $D$.
We can also conclude the following.

\begin{observation}\label{obs:dibracesandbraces}
    The splits of the dibraces of a strongly connected digraph $D$ are exactly the braces of the split of $D$.
\end{observation}

This begs the question of what the small cycle sums turn into, once we look at the split of a digraph.
Let $B_1$ and $B_2$ be bipartite graphs, such that their intersection is a cycle $C$ of length 4 and we have $V(B_i) \setminus V(C) \neq \emptyset$ for all $i \in [2]$.
Further, let $S \subseteq E(C)$ be some subset of the edges of $C$.
A \emph{4-cycle sum of $B_1$ and $B_2$ (at $C$)} is a graph $( B_1 \cup B_2 ) - S$.
Small cycle sums can be translated to 4-cycle sums and vice versa, similar to \Cref{obs:dirtightcutandtightcut}.

\begin{lemma}\label{lem:smallcycsumandtrisum}
    Let $D$ be the small cycle sum of $D_1$ and $D_2$, then the split of $D$ is a 4-cycle sum of the splits of $D_1$ and $D_2$.

    Vice versa, let $B$ be the 4-cycle sum of $B_1$ and $B_2$ at $C$, then for every perfect matching $M$ of $B$ there exist two perfect matchings $M_1$ and $M_2$ of $B_1$ and $B_2$ respectively, such that $\DirM{B}{M}$ is a small cycle sum of $\DirM{B_1}{M_1}$ and $\DirM{B_2}{M_2}$ and $M_i \setminus E(C) \subseteq M$ for each $i \in [2]$.
\end{lemma}
\begin{proof}
    For the first part of the statement, it is easy to observe that all of the graphs in \Cref{fig:smallcyclesum} contain a cycle of length 4 in their split.
    In particular, the split of a digon corresponds to a $C_4$ with a perfect matching.
    The split of the transitively oriented triangle corresponds to a $C_4$ of which one edge is found in the matching and the other two matching edges for the split are found outside of the $C_4$.
    And the split of the anti-directed orientation of the $C_4$ found in \Cref{fig:smallcyclesum} corresponds to a $C_4$ together with four matching edges outside of the $C_4$, each covering one vertex of the $C_4$.
    From these observations it is easy to observe that we can simply perform a 4-cycle-sum of the splits of $D_1$ and $D_2$ at the $C_4$ found in the subgraphs corresponding to the graphs depicted in \Cref{fig:smallcyclesum}.

    When dealing with the second part of the statement, we must consider how the perfect matching $M$ interacts with $C$.
    Let $V_1$ and $V_2$ be the two partitions of $V(B)$ according to which we will build the $M$-direction of $B$.
    If $C$ is conformal, then it is easy to see that we can perform a 2-sum of $\DirM{B_1}{M_1}$ and $\DirM{B_2}{M_2}$, with $M_i = M \cap E(B_i)$ for both $i \in [2]$.
    
    Should $|E(C) \cap M| = 1$, we claim that we can perform a 3-sum at $\DirM{B_1}{M_1}$ and $\DirM{B_2}{M_2}$.
    Consider $B_1$ and $B_2$ and note that by definition $B_1 \cap B_2 = C$.
    Let $e, f \in M \setminus E(C)$ be the two distinct edges covering vertices of $C$ that are not covered by the single edge in $E(C) \cap M|$.
    Since $B$ is bipartite and $M$ is a perfect matching, we must have either $e,f \in E(B_1)$ or $e,f \in E(B_2)$.
    W.l.o.g.\ we suppose that $e,f \in E(B_1)$ holds and we let $M_1 = M \cap E(B_1)$, which is clearly a perfect matching of $B_1$.
    Furthermore, we let $D_1 = \DirM{B_1}{M_1}$.
    Consider $M \cap E(B_2)$ and note that this is almost a perfect matching $B_2$, except that there exist exactly two adjacent vertices $a,b \in V(C)$ which are not covered.
    W.l.o.g.\ we may suppose that $a \in V_1$ and $b \in V_2$.
    
    If we consider the $(M \cap E(B_2)) \cup \{ e,f \}$-direction of the graph we would obtain from simply adding $e$ and $f$, together with their endpoints to $B_2$, we could therefore note that the vertex $a$ was contracted into has no incoming edge and the vertex $b$ was contracted into has no outgoing edge.
    Thus, this graph would be our desired $D_2'$ from \Cref{def:smallcyclesum}.
    To construct $D_2$, we can therefore simply define $M_2 = (M \cap E(B_2)) \cup \{ ab \}$.
    (Note how this corresponds to the contraction of the edge $(w,v)$ in \Cref{def:smallcyclesum}.)
    It is now easy to confirm that $D$ is in fact a 3-sum at $D_1$ and $D_2$.
    The case in which $|E(C) \cap M| = 0$ can then be argued for in an analogous fashion.
\end{proof}

Note in particular the stark contrast in the size and complexity of the definition of small cycle sums (see \Cref{def:smallcyclesum}) and the definition of the 4-cycle sum.
This difference is also felt in a pronounced way when trying to argue about these concepts in the context of a proof, which is the main reason why we are switching entirely to the matching theoretic setting for this subsection.

As a consequence of this lemma and the fact that a 4-cycle sum at $C$ of two planar graphs for which $C$ is a face is planar, we can deduce that forming a small cycle sum of two strongly planar digraphs at a cycle in those digraphs that forms a face yields another strongly planar digraph.
The two lemmas also allow us to translate \cref{thm:nonevenstructure} to the world of bipartite graphs with perfect matchings.

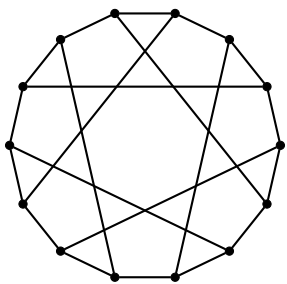
\begin{figure}[h!]
	\begin{center}
					\begin{tikzpicture}[scale=0.6]
			
			\node (V0) at (0:0) [draw=none] {};
			
			\foreach\i in {1,2,3,4,5,6,7,8,9,10,11,12,13,14}
			{
				\node (V\i) at ($(V0)+({(360/14 * \i)}:3)$) [draw, circle, scale=0.3, fill, label={}] {};
			}
			
			\foreach\i in {1,2,3,4,5,6,7,8,9,10,11,12,13}
			{
				\pgfmathtruncatemacro\iplus{\i+1}
				\path (V\i) edge[thick] (V\iplus);
			}
			\path (V14) edge[thick] (V1);
			
			\path
				(V1) edge[thick] (V6)
				(V2) edge[thick] (V11)
				(V3) edge[thick] (V8)
				(V4) edge[thick] (V13)
				(V5) edge[thick] (V10)
				(V7) edge[thick] (V12)
				(V9) edge[thick] (V14)
			;
			
			\end{tikzpicture}
	\end{center}
	\caption{The Heawood graph.
            Note that this is the split of the graph in \Cref{fig:F7}.}
	\label{fig:heawood}
\end{figure}

\begin{theorem}[Robertson, Seymour, and Thomas \cite{Robertson1999PermanentsPfaffianOrientations} and McCuaig \cite{McCuaig2004Polya}]\label{thm:pfaff4cyc}
A bipartite, matching covered graph is Pfaffian if and only if all of its braces are either the Heawood graph (see \Cref{fig:heawood}) or can be obtained from planar braces by means of the 4-cycle-sum.
\end{theorem}

Related to this theorem and also the Two Paths Theorem mentioned earlier, we will make use of a recent result concerning a specific type of cross over a $C_4$.
Though we only make limited use of this definition here, we will use it more extensively later on.

\begin{definition}[Conformal cross]\label{def:conformalcross}
    Let $G$ be an undirected graph with a perfect matching and let $C \subseteq G$ be a conformal cycle.
    Two disjoint $C$-paths $L$ and $R$, where $L$ has the endpoints $l_1,l_2$ and $R$ has two endpoints $r_1,r_2$, constitute a \emph{conformal cross on $C$} if $l_1,r_1,l_2,r_2$ occur on $C$ in the given order and $C \cup L \cup R$ is a conformal and matching covered subgraph of $G$.
\end{definition}

\begin{theorem}[Giannopoulou and Wiederrecht \cite{giannopoulou2021two}]\label{lem:crossoverc4}
    Let $B$ be a brace and $C$ be a 4-cycle in $B$.
    There exists a conformal cross over $C$ in $B$ if and only if $C$ is contained in a conformal bisubdivision of $K_{3,3}$.
\end{theorem}

\subsection{Shifting an even dicycle in a planar decomposition}\label{subsec:strongplanarshifting}

To give an explanation for the necessity of the following lemma, note that if we have a plane decomposition of a digraph containing two dicycles, their traces may in fact induce a planar drawing, but that does not mean that the actual drawing of these two dicycles is planar, since we have to fudge the drawing of an object to find its trace, and thus we cannot directly apply \Cref{lem:planarshifting}.
To prove the lemma we therefore have to either find a planar drawing for these two dicycles, or somehow construct an even dicycle from the fact that certain parts of the drawing of these two dicycles cannot be made planar.

\begin{lemma}[Planar Trace Shifting]\label{lem:planartraceshifting}
Let $D$ be a digraph, let $C_1, C_2 \subseteq D$ respectively be an even and odd dicycle, and let $\delta = (\Gamma, \mathcal{V}, \mathcal{D})$ be a plane decomposition of $D$ with a maelstrom $m \in \mathcal{D}$, such that $C_1 \cap \sigma(m) \neq \emptyset$, the drawing of $C_1 \cap \sigma(m)$ in $\Gamma$ is cross-free, and $C_2$ is grounded in $\delta$.

Then either
\begin{enumerate}
    \item there exists a directed path $P \subseteq C_1$, such that $P$ is drawn by $\Gamma$ on the closure of the $m$-tight disk of $C_2$ and $C_2 \cup P$ contains an even dicycle,
    
    \item there exists a conjunction cell in $C(\delta)$ that contains an even dicycle, or
    
    \item there exists a big vertex $\mathsf{v}$ that contains an even dicycle.
\end{enumerate}
\end{lemma}

The setting of this lemma is not much different from the one in \Cref{cor:planarshifting}.
In fact our main goal will be to prove that we can adjust the drawing such that we are allowed to apply \Cref{cor:planarshifting}, with only a few cases resulting in us actually finding an even dicycle in a big vertex or a conjunction cell.
Our attention thus centers on the behaviour of the two dicycles passing through the big vertices and conjunction cells.
We first argue that the former are in fact not an obstacle.

\begin{lemma}\label{lem:planartraceshiftingnobigvertices}
Let $D$ be a digraph, let $C_1, C_2 \subseteq D$ respectively be an even and odd dicycle, and let $\delta = (\Gamma, \mathcal{V}, \mathcal{D})$ be a plane decomposition of $D$ with a maelstrom $m \in \mathcal{D}$ and without conjunction cells, such that $C_1 \cap \sigma(m) \neq \emptyset$, the drawing of $C_1 \cap \sigma(m)$ in $\Gamma$ is cross-free, and $C_2$ is grounded in $\delta$.

Then there either exists a directed path $P \subseteq C_1$, such that $P$ is drawn by $\Gamma$ on the closure of the $m$-tight disk of $C_2$ and $C_2 \cup P$ contains an even dicycle, or there exists a big vertex $\mathsf{v}$ that contains an even dicycle.
\end{lemma}
\begin{proof}
    As in \Cref{cor:planarshifting}, we let $C_1$ and $C_2$ be a smallest counterexample to our claim.
    Furthermore, we may suppose that no big vertex $\mathsf{v} \in \mathcal{V}$, with $\pi_\delta(\mathsf{v}) = (v,X)$, exists such that $D[X \cup \{ v \}]$ contains an even dicycle, as we are otherwise done.

    Suppose that there exists a $\mathsf{v} \in \mathcal{V}$.
    Let $\pi_\delta(\mathsf{v}) = (v,X)$.
    Recall that by definition $(V(D-X), X \cup \{ v \})$ is a directed tight cut in $D$.
    We let $\OutNeighbours{}{\mathsf{v}}$ be the set of vertices which are the head of an outgoing edge of $\mathsf{v}$ and suppose w.l.o.g.\ that all edges crossing $\Boundary{\mathsf{v}}$ are outgoing.

    We note that since $C_1 \cup C_2$ consist of two dicycles and $v$ is a directed separator, $\OutNeighbours{}{\mathsf{v}}$ contains at most two vertices $u_1, u_2$, with $u_1 \in V(C_1)$ and $u_2 \in V(C_2)$, if the respectively exist.
    Thus the interior of the big vertex can be replaced with a planar graph consisting of at most four distinct $v$-$\OutNeighbours{}{\mathsf{v}}$-paths, two for $u_1$ and $u_2$ each, one for each of the two possible parities, depending on what sort of $v$-$\OutNeighbours{}{\mathsf{v}}$-paths pass through the big vertex.
    This procedure can be performed repeatedly until no more big vertices remain, letting us apply \Cref{cor:planarshifting} once we have gotten rid of all big vertices, which lets us find the desired path.
\end{proof}

As a consequence of this lemma, we can now focus on dealing with conjunction cells, as we can always first contract all directed tight cuts associated with the big vertices of the decomposition.
For this purpose, we define the following setting for the remainder of this subsection based on the assumptions and objects laid out in the statement of \Cref{lem:planartraceshifting}:
Let $C_1$ and $C_2$ be a smallest counterexample to \Cref{lem:planartraceshifting}, which must mean that $C_1$ and $C_2$ intersect.
Let $\mathcal{C} \subseteq \mathcal{D}$ be the set of conjunction cells into which a part of $C_1 \cup C_2$ is drawn and let $H = (C_1 \cup C_2) \cup \bigcup_{c \in \mathcal{C}} \widetilde{c}$, which ensures that we have all vertices used in the small cycle sum, even if they are not found on $C_1$ or $C_2$.

Further, we let $S$ be the subgraph of $H - \sigma(m)$ formed by the union of the strongly connected component of $H - \sigma(m)$ that contains $C_2$ in its entirety and all vertices of all $c \in \mathcal{C}$ into which a vertex of this component is drawn.
As argued in \Cref{lem:planartraceshiftingnobigvertices}, we can assume that $S$ and the restriction of $\delta$ to $S$ are already reduced and thus we can assume that the only non-trivial component of $S$ is in fact a dibrace, or a subdivision of a dibrace.
Note that, since $C_1 \cup C_2$ simply consists of two dicycles, the set of weakly connected components of $(D - S) \cap (C_1 \cup C_2)$ solely contains subpaths of $C_1$.
As a consequence, our goal from now until the end of this subsection will be to provide a planar drawing of $S$ that can be extended to a planar drawing of $H$, as this allows us to apply \Cref{cor:planarshifting} and thus prove the lemma.
The properties of the involved traces of those parts of the even dicycle outside of $m$ then follow immediately from \Cref{def:trace}.

Note that due to $S$ already being reduced, the only parts of $S$ that are not drawn in a planar way by $\Gamma$ lie in the conjunction cells of $\delta$.
Let $s \in C(\delta)$ be a conjunction cell containing a part of $S$ in its interior.
To find drawings for the parts of $S$ drawn in conjunction cells, we will switch to the split $\DirM{B}{M}$ of $S$.
As the only non-trivial component of $S$ is a dibrace, we know that its split is a brace.
According to \Cref{obs:strongplanarityisplanarity}, if we can provide a planar drawing for $B$, then $S$ possesses a strongly planar drawing.
Note that now, instead of having to deal with small cycle sums, we must concern ourselves with 4-cycle-sums, which involve exactly four vertices.
Let $a,b,c,d$ be the four vertices involved in the 4-cycle-sum which corresponds to the small cycle sum in $S$ at $s$.
Further, let $B_0', B_1', \ldots , B_k'$ be the graphs which form $B$ as a 4-cycle-sum, such that $a,b,c,d$ form a 4-cycle $C = (a,b,c,d)$ in each $B_i'$ and let $B_i = B_i' - (E(C) \setminus E(B))$ for $i \in [0,k]$.
Note that due to \Cref{def:maelstrom}, we know that $\sigma(s)$ is non-even and thus $\bigcup_{i=1}^k B_i$, and in fact $\bigcup_{i=1}^k B_i'$, is Pfaffian thanks to \Cref{thm:pfafftooddbicycle} and \Cref{thm:nonevendigraphs}.

We now briefly consider the implications of \Cref{thm:pfaff4cyc} for the brace $B$ in particular.
Since we can assume $B$ to not be planar, as we would otherwise be done, we know that $B$ is built from planar braces via 4-cycle sums.
Due to this, we can imagine $B$ as being decomposable in a tree-like fashion into several planar braces stuck together on 4-cycles, from which a few edges might be removed to actually form $B$ itself.
In particular, due to the tree-like nature of this decomposition, there exists 4-cycle sums at which all but one of the involved graphs are planar braces that are only attached to one other brace in the decomposition.
We can be even more specific, thanks to a lemma by McCuaig.

\begin{lemma}[McCuaig \cite{mccuaig2000even, McCuaig2004Polya}]\label{lem:bracec4takeapart}
    Let $G_1$ and $G_2$ be bipartite graphs, such that $G_1 \cap G_2$ is a 4-cycle, and $G_1$ and $G_2$ are proper subgraphs of $G_1 \cup G_2$.
    If $G_1 \cup G_2$ is a brace, then $G_1$ and $G_2$ are braces.
\end{lemma}

Thus the braces we find are not only planar, but we can actually assume that they are all attached via a facial 4-cycle of theirs.
This allows us to conclude that, if we can prove that the 4-cycle-sum of these braces has a planar drawing, we can then propagate this property by induction through the rest of the decomposition.

We therefore make the following assumption:
$B_1', \ldots , B_k'$ all have planar drawings in which $C$ bounds a face and the split of $S - \sigma(s)$ is a subgraph of $B_0$.
Our aim will then be to find a planar drawing of the 4-cycle sum of $B_1', \ldots , B_k'$ that can be placed into $s$ and is otherwise compatible with the drawing of $B_0$.
Of course, this assumption is not guaranteed to hold automatically.
However should some $B_i'$ with $i \in [k]$ not be planar with $C$ bounding a face, we may consider the tree-like decomposition of $B_k'$ and proceed with the arguments we are about to present to progressively find planar drawings that we can stick together to form a planar drawing of $B_i'$, following the inductive argument we have just laid out above.

We call $B_0$ the \emph{base} of the 4-cycle-sum at $C$.
In other words, $B_0$ represents the part of the graph that is drawn outside of the interior of the conjunction cell, including the vertices on the boundary of $s$.
The graphs $B_1, \ldots , B_k$ will be called \emph{flaps}.

To properly analyze how $B$ behaves inside of $s$, we will however also need to know how $M$ interacts with $U = V(C)$ and especially in which parts of $B$ the edges of $M$ that cover vertices of $C$ are drawn.
In particular, we distinguish between edges containing two vertices of $U$, edges drawn in the interior of $s$ by $\Gamma$, and edges that are found in the base.
The last type of edge is relevant, since we want to preserve the drawing of these edges.
It is easy to observe that there are six types of ways in which $M$ and $C$ can interact (see \Cref{fig:4sumtypes}):

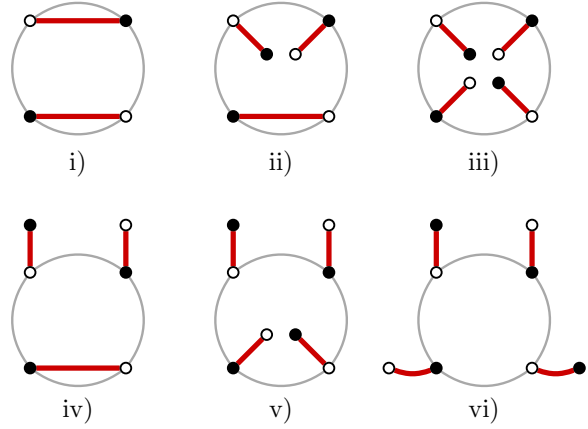
\begin{figure}
    \centering

    \scalebox{0.9}{
    \begin{tikzpicture}[scale=1]

        \pgfdeclarelayer{background}
		\pgfdeclarelayer{foreground}
			
		\pgfsetlayers{background,main,foreground}
			
        \begin{pgfonlayer}{main}
        \node (T) [v:ghost] {};
        \node (B) [v:ghost,position=270:35mm from T] {};

        \node (LT) [v:ghost,position=180:30mm from T] {};
        \node (MT) [v:ghost,position=0:0mm from T] {};
        \node (RT) [v:ghost,position=0:30mm from T] {};

        \node (LTlabel) [v:ghost,position=270:14mm from LT] {i)};
        \node (MTlabel) [v:ghost,position=270:14mm from MT] {ii)};
        \node (RTlabel) [v:ghost,position=270:14mm from RT] {iii)};

        \node (LB) [v:ghost,position=180:30mm from B] {};
        \node (MB) [v:ghost,position=0:0mm from B] {};
        \node (RB) [v:ghost,position=0:30mm from B] {};

        \node (LBlabel) [v:ghost,position=270:15.5mm from LB] {iv)};
        \node (MBlabel) [v:ghost,position=270:15.5mm from MB] {v)};
        \node (RBlabel) [v:ghost,position=270:15.5mm from RB] {vi)};

        \node (LTpicture) [v:ghost,position=0:0mm from LT] {
            \begin{tikzpicture}[scale=1]

                \pgfdeclarelayer{background}
		          \pgfdeclarelayer{foreground}
			
		          \pgfsetlayers{background,main,foreground}
			
                \begin{pgfonlayer}{main}

                    \node (C) [v:ghost] {};

                    \node (a_1) [v:main,position=45:10mm from C] {};
                    \node (b_1) [v:mainempty,position=135:10mm from C] {};
                    \node (a_2) [v:main,position=225:10mm from C] {};
                    \node (b_2) [v:mainempty,position=315:10mm from C] {};
                    
                \end{pgfonlayer}{main}

                \begin{pgfonlayer}{background}

                    \draw [e:main,color=Gray,bend right=35] (a_1) to (b_1);
                    \draw [e:main,color=Gray,bend right=35] (b_1) to (a_2);
                    \draw [e:main,color=Gray,bend right=35] (a_2) to (b_2);
                    \draw [e:main,color=Gray,bend right=35] (b_2) to (a_1);

                    \draw [e:main,line width=2.2pt,color=BostonUniversityRed] (a_1) to (b_1);
                    \draw [e:main,line width=2.2pt,color=BostonUniversityRed] (b_2) to (a_2);
                
                \end{pgfonlayer}{background}

                \begin{pgfonlayer}{foreground}

                \end{pgfonlayer}{foreground}
                
            \end{tikzpicture}
        };

        \node (MTpicture) [v:ghost,position=0:0mm from MT] {
            \begin{tikzpicture}[scale=1]

                \pgfdeclarelayer{background}
		          \pgfdeclarelayer{foreground}
			
		          \pgfsetlayers{background,main,foreground}
			
                \begin{pgfonlayer}{main}

                    \node (C) [v:ghost] {};

                    \node (a_1) [v:main,position=45:10mm from C] {};
                    \node (b_1) [v:mainempty,position=135:10mm from C] {};
                    \node (a_2) [v:main,position=225:10mm from C] {};
                    \node (b_2) [v:mainempty,position=315:10mm from C] {};

                    \node (x_1) [v:mainempty,position=225:7mm from a_1] {};
                    \node (y_1) [v:main,position=315:7mm from b_1] {};
                    
                \end{pgfonlayer}{main}

                \begin{pgfonlayer}{background}

                    \draw [e:main,color=Gray,bend right=35] (a_1) to (b_1);
                    \draw [e:main,color=Gray,bend right=35] (b_1) to (a_2);
                    \draw [e:main,color=Gray,bend right=35] (a_2) to (b_2);
                    \draw [e:main,color=Gray,bend right=35] (b_2) to (a_1);

                    \draw [e:main,line width=2.2pt,color=BostonUniversityRed] (a_1) to (x_1);
                    \draw [e:main,line width=2.2pt,color=BostonUniversityRed] (b_1) to (y_1);
                    \draw [e:main,line width=2.2pt,color=BostonUniversityRed] (b_2) to (a_2);
                
                \end{pgfonlayer}{background}

                \begin{pgfonlayer}{foreground}

                \end{pgfonlayer}{foreground}
                
            \end{tikzpicture}
        };

        \node (RTpicture) [v:ghost,position=0:0mm from RT] {
            \begin{tikzpicture}[scale=1]

                \pgfdeclarelayer{background}
		          \pgfdeclarelayer{foreground}
			
		          \pgfsetlayers{background,main,foreground}
			
                \begin{pgfonlayer}{main}

                    \node (C) [v:ghost] {};

                    \node (a_1) [v:main,position=45:10mm from C] {};
                    \node (b_1) [v:mainempty,position=135:10mm from C] {};
                    \node (a_2) [v:main,position=225:10mm from C] {};
                    \node (b_2) [v:mainempty,position=315:10mm from C] {};

                    \node (x_1) [v:mainempty,position=225:7mm from a_1] {};
                    \node (y_1) [v:main,position=315:7mm from b_1] {};
                    \node (x_2) [v:mainempty,position=45:7mm from a_2] {};
                    \node (y_2) [v:main,position=135:7mm from b_2] {};
                    
                \end{pgfonlayer}{main}

                \begin{pgfonlayer}{background}

                    \draw [e:main,color=Gray,bend right=35] (a_1) to (b_1);
                    \draw [e:main,color=Gray,bend right=35] (b_1) to (a_2);
                    \draw [e:main,color=Gray,bend right=35] (a_2) to (b_2);
                    \draw [e:main,color=Gray,bend right=35] (b_2) to (a_1);

                    \draw [e:main,line width=2.2pt,color=BostonUniversityRed] (a_1) to (x_1);
                    \draw [e:main,line width=2.2pt,color=BostonUniversityRed] (b_1) to (y_1);
                    \draw [e:main,line width=2.2pt,color=BostonUniversityRed] (a_2) to (x_2);
                    \draw [e:main,line width=2.2pt,color=BostonUniversityRed] (b_2) to (y_2);
                
                \end{pgfonlayer}{background}

                \begin{pgfonlayer}{foreground}

                \end{pgfonlayer}{foreground}
                
            \end{tikzpicture}
        };

        \node (LBpicture) [v:ghost,position=0:0mm from LB] {
            \begin{tikzpicture}[scale=1]

                \pgfdeclarelayer{background}
		          \pgfdeclarelayer{foreground}
			
		          \pgfsetlayers{background,main,foreground}
			
                \begin{pgfonlayer}{main}

                    \node (C) [v:ghost] {};

                    \node (a_1) [v:main,position=45:10mm from C] {};
                    \node (b_1) [v:mainempty,position=135:10mm from C] {};
                    \node (a_2) [v:main,position=225:10mm from C] {};
                    \node (b_2) [v:mainempty,position=315:10mm from C] {};

                    \node (x_1) [v:mainempty,position=90:7mm from a_1] {};
                    \node (y_1) [v:main,position=90:7mm from b_1] {};
                    
                \end{pgfonlayer}{main}

                \begin{pgfonlayer}{background}

                    \draw [e:main,color=Gray,bend right=35] (a_1) to (b_1);
                    \draw [e:main,color=Gray,bend right=35] (b_1) to (a_2);
                    \draw [e:main,color=Gray,bend right=35] (a_2) to (b_2);
                    \draw [e:main,color=Gray,bend right=35] (b_2) to (a_1);

                    \draw [e:main,line width=2.2pt,color=BostonUniversityRed] (a_1) to (x_1);
                    \draw [e:main,line width=2.2pt,color=BostonUniversityRed] (b_1) to (y_1);
                    \draw [e:main,line width=2.2pt,color=BostonUniversityRed] (b_2) to (a_2);
                
                \end{pgfonlayer}{background}

                \begin{pgfonlayer}{foreground}

                \end{pgfonlayer}{foreground}
                
            \end{tikzpicture}
        };

        \node (MBpicture) [v:ghost,position=0:0mm from MB] {
            \begin{tikzpicture}[scale=1]

                \pgfdeclarelayer{background}
		          \pgfdeclarelayer{foreground}
			
		          \pgfsetlayers{background,main,foreground}
			
                \begin{pgfonlayer}{main}

                    \node (C) [v:ghost] {};

                    \node (a_1) [v:main,position=45:10mm from C] {};
                    \node (b_1) [v:mainempty,position=135:10mm from C] {};
                    \node (a_2) [v:main,position=225:10mm from C] {};
                    \node (b_2) [v:mainempty,position=315:10mm from C] {};

                    \node (x_1) [v:mainempty,position=90:7mm from a_1] {};
                    \node (y_1) [v:main,position=90:7mm from b_1] {};

                    \node (x_2) [v:mainempty,position=45:7mm from a_2] {};
                    \node (y_2) [v:main,position=135:7mm from b_2] {};
                    
                \end{pgfonlayer}{main}

                \begin{pgfonlayer}{background}

                    \draw [e:main,color=Gray,bend right=35] (a_1) to (b_1);
                    \draw [e:main,color=Gray,bend right=35] (b_1) to (a_2);
                    \draw [e:main,color=Gray,bend right=35] (a_2) to (b_2);
                    \draw [e:main,color=Gray,bend right=35] (b_2) to (a_1);

                    \draw [e:main,line width=2.2pt,color=BostonUniversityRed] (a_1) to (x_1);
                    \draw [e:main,line width=2.2pt,color=BostonUniversityRed] (b_1) to (y_1);

                    \draw [e:main,line width=2.2pt,color=BostonUniversityRed] (a_2) to (x_2);
                    \draw [e:main,line width=2.2pt,color=BostonUniversityRed] (b_2) to (y_2);
                
                \end{pgfonlayer}{background}

                \begin{pgfonlayer}{foreground}

                \end{pgfonlayer}{foreground}
                
            \end{tikzpicture}
        };

        \node (RBpicture) [v:ghost,position=0:0mm from RB] {
            \begin{tikzpicture}[scale=1]

                \pgfdeclarelayer{background}
		          \pgfdeclarelayer{foreground}
			
		          \pgfsetlayers{background,main,foreground}
			
                \begin{pgfonlayer}{main}

                    \node (C) [v:ghost] {};

                    \node (a_1) [v:main,position=45:10mm from C] {};
                    \node (b_1) [v:mainempty,position=135:10mm from C] {};
                    \node (a_2) [v:main,position=225:10mm from C] {};
                    \node (b_2) [v:mainempty,position=315:10mm from C] {};

                    \node (x_1) [v:mainempty,position=90:7mm from a_1] {};
                    \node (y_1) [v:main,position=90:7mm from b_1] {};

                    \node (x_2) [v:mainempty,position=180:7mm from a_2] {};
                    \node (y_2) [v:main,position=0:7mm from b_2] {};
                    
                \end{pgfonlayer}{main}

                \begin{pgfonlayer}{background}

                    \draw [e:main,color=Gray,bend right=35] (a_1) to (b_1);
                    \draw [e:main,color=Gray,bend right=35] (b_1) to (a_2);
                    \draw [e:main,color=Gray,bend right=35] (a_2) to (b_2);
                    \draw [e:main,color=Gray,bend right=35] (b_2) to (a_1);

                    \draw [e:main,line width=2.2pt,color=BostonUniversityRed] (a_1) to (x_1);
                    \draw [e:main,line width=2.2pt,color=BostonUniversityRed] (b_1) to (y_1);

                    \draw [e:main,line width=2.2pt,color=BostonUniversityRed,bend left=20] (a_2) to (x_2);
                    \draw [e:main,line width=2.2pt,color=BostonUniversityRed,bend right=20] (b_2) to (y_2);
                
                \end{pgfonlayer}{background}

                \begin{pgfonlayer}{foreground}

                \end{pgfonlayer}{foreground}
                
            \end{tikzpicture}
        };

        \end{pgfonlayer}{main}
        
        \begin{pgfonlayer}{foreground}
        \end{pgfonlayer}{foreground}

        \begin{pgfonlayer}{background}
        \end{pgfonlayer}{background}
        
    \end{tikzpicture}
    }

    \caption{The types of interactions between $M$ and $U$ at a 4-cycle sum.}
    \label{fig:4sumtypes}
\end{figure}
\begin{enumerate}
    \item two edges of $M$ cover all of $U$,
    
    \item one edge of $M$ covers two vertices of $U$ and the other two are matched outside of $B_0$,
    
    \item all of $U$ is matched outside of $B_0$,
    
    \item one edge of $M$ covers two vertices of $U$ and the other two are matched into $B_0$,
    
    \item two vertices of $M$ are matched into $B_0$ and the other two are matched outside of $B_0$, and
    
    \item all of $U$ is matched into $B_0$.
\end{enumerate}
Knowing the type of $B_0$, we can determine the type of $B_1, \ldots , B_k$, if we were decomposing the graph with them as a base and narrow down the structure of the sum as follows:
\begin{itemize}
    \item If $B_0$ has type i), then all flaps have type i).
    
    \item If $B_0$ has type ii), then one flap has type iv) and the others have type ii) as well.
    
    \item If $B_0$ has type iii), either there are exactly two flaps with type v) or there is one flap with type vi), and all other flaps have type iii) in both cases. 
    
    \item If $B_0$ has type iv), then all flaps have type ii).
    
    \item If $B_0$ has type v), then one flap has type v) and all others have type iii).
    
    \item If $B_0$ has type vi), then all flaps have type iii).
\end{itemize}
Let $\mathcal{P}$ be the set of components of $C_1 \cap \bigcup_{i=1}^k B_i$ and let $\mathcal{Q}$ be the set of components of $C_2 \cap \bigcup_{i=1}^k B_i$.
Since $C_1$ and $C_2$ are cycles and $U$ is separator in $G$, the sets $\mathcal{P}$ and $\mathcal{Q}$ each contain at most two components.
Furthermore, all components in $\mathcal{P} \cup \mathcal{Q}$ are $M$-conformal paths, since they are splits of directed paths in $S$.
We will use $P,P'$ as the names for the up to two elements of $\mathcal{P}$ and similarly, we use $Q,Q'$ for the paths in $\mathcal{Q}$.
Furthermore, we will say that a flap $B_i$ \emph{hosts} a part of a path $R \in \mathcal{P} \cup \mathcal{Q}$, if $(V(R) \cap V(B_i)) \setminus U$ is non-empty.

If either $\mathcal{P}$ or $\mathcal{Q}$ are empty, then there is nothing to prove.
Similarly, if $k = 1$, then there is again nothing to prove, since all $B_i'$ with $i \in [k]$ have a planar drawing with $C$ bounding a face.
Should $k \geq 2$ and for all $i \in [2,k]$ any edge $e \in E(B_i - C)$ can be separated from $U$ in $B_i'$ via a separator consisting of two adjacent vertices in $C$, we can draw $B_1$ in the center of $s$ and, since the only information provided by the other flaps can be represented via drawings whose connected components only touch two vertices of $U$ that are adjacent in $C$, we can easily draw the other flaps around the drawing of $B_1$ in $s$.

The above observations tell us that we can additionally assume that $k \geq 2$, that both $\mathcal{P}$ and $\mathcal{Q}$ are non-empty, and that there exist distinct $i,j \in [k]$, such that both $E(B_i - C)$ and $E(B_j - C)$ contain an edge that cannot be separated from the vertices of $U$ in $B_i'$, respectively $B_j'$, via a separator consisting of wo adjacent vertices in $C$.
Furthermore, since $\bigcup_{i=1}^k B_i'$ is Pfaffian, we know that $\bigcup_{i=1}^k B_i'$ cannot contain a conformal cross on $C$, as otherwise \Cref{lem:crossoverc4} yields a witness that $\bigcup_{i=1}^k B_i'$ is not Pfaffian.
We now move on to proving \Cref{lem:planartraceshifting} for the six cases suggested by \Cref{fig:4sumtypes}.

\begin{lemma}\label{lem:planartraceshiftingtypei}
    If $B_0$ has type i), there exists a planar drawing for the split of $\sigma(s)$ that is compatible with the drawing of the split of $S - \sigma(s)$.    
\end{lemma}
\begin{proof}
Note that $| \mathcal{P} | = | \mathcal{Q} | = 1$, since by definition all paths in these two sets must begin and end on an edge from $M$.
If $P$ is hosted by a flap $H$, then $H$ is the only flap that can host $P$, as $P$ must use two non-matching edges incident with two distinct vertices in $U$ to enter and exit $H$ and the matching edges covering $U$ force $P$ leave the flaps entirely after exiting $H$.
Thus, since $k \geq 2$, $P$ and $Q$ must be hosted by different flaps.
Since these two paths cannot form a conformal cross on $C$, we can therefore easily find a planar drawing for the split of $\sigma(s)$ on $s$ that attaches appropriately to drawing of $B_0$.
\end{proof}

We call a path $P$ in a graph $G$ with a perfect matching $M$ \emph{internally $M$-conformal}, if $P$ minus its endpoints is $M$-conformal.
Note that an (internally) $M$-conformal path in a bipartite graph is of odd length if and only if its endpoints belong to opposite colour classes.
In the upcoming arguments we will often use this implicitly when narrowing down the types of potential $U$-paths. 

\begin{lemma}\label{lem:planartraceshiftingtypeii}
    If $B_0$ has type ii), there exists a planar drawing for the split of $\sigma(s)$ that is compatible with the drawing of the split of $S - \sigma(s)$.    
\end{lemma}
\begin{proof}
Let $ab \in M$ be the edge that covers two vertices in $U$, with $cc', dd' \in M$ being the two edges whose other endpoints are found in the interior of the flap of type iv).
Assume that there exists a path $R \in \mathcal{P} \cup \mathcal{Q}$ with $a$ or $b$ as an endpoint.
W.l.o.g.\ let one endpoint of $R$ be $a$ and note that $ab \in M$ must be the single edge incident with $a$ on $R$, as $R$ is $M$-conformal.
If $R$ solely uses the edge $ab$ and then leaves the flaps again, we can easily represent $R$ by drawing the edge $ab$ near the boundary of $s$ and focusing on the rest of the paths in $\mathcal{P} \cup \mathcal{Q}$.
Thus we can assume that $R$ is hosted by one of the flaps.
This however necessitates, that $R$ includes at least a third vertex of $U$, either $c$ or $d$, to be able to leave the flap again.
In this case, $R$ is the sole path in either $\mathcal{P}$ or $\mathcal{Q}$ and it is hosted by only a single flap.
Accordingly, we may suppose that $| \mathcal{P} | = | \mathcal{Q} | = 1$ and can thus find an appropriate drawing in an analogous fashion to \Cref{lem:planartraceshiftingtypei}.
\end{proof}

Before we delve into the next case, we note that, since we only need some planar drawing of the split of $C_1 \cup C_2$ that agrees with $\Gamma$ outside the conjunction cells, the four vertices on the boundary of $s$ may in fact be moved into the interior of $s$, should no edge incident to them in $C_1 \cup C_2$ be drawn outside of $s$.

\begin{lemma}\label{lem:planartraceshiftingtypeiiipartone}
    If $B_0$ has type iii) and $B_1, B_2$ are type v) flaps, there exists a planar drawing for the split of $\sigma(s)$ that is compatible with the drawing of the split of $S - \sigma(s)$.    
\end{lemma}
\begin{proof}
We let $a,b$ be the two vertices in $U$ such that $aa', bb' \in M \cap E(B_1)$ and let $cc', dd' \in M \cap E(B_2)$.
Note that if any $M$-conformal path $R \in \mathcal{P} \cup \mathcal{Q}$ has $a$ or $b$ as an endpoint, then some part of $R$ is hosted by $B_1$ and an analogous restriction holds for $B_2$ and $c$, as well as $d$.

Suppose some $R \in \mathcal{P} \cup \mathcal{Q}$ has an internally $M$-conformal $b$-$c$-subpath contained in $B_1$.
Since $R$ must contain the edge $bb'$, there must either exist an $M$-conformal $a$-$b$-subpath of $R$ found entirely in $B_1$ or there exists a $b$-$d$-subpath of $R$ in $B_1$.
In fact, as $R$ must similarly contain $cc'$, there must either exist an $M$-conformal $c$-$d$-subpath of $R$ found entirely in $B_2$ or there exists an $a$-$c$-subpath of $R$ in $B_2$.
Note that the existence of a $b$-$d$-subpath similarly forces $R$ to enter $B_2$ again.
Since it must then contain $dd'$ and $R$ is a path, this leads to a contradiction, since we cannot connect $d$ and $c$ here, without turning $R$ into a cycle.
Thus we can narrow down our focus to the case in which $B_1$ contains an $M$-conformal $a$-$b$-subpath of $R$ and $B_2$ an $M$-conformal $c$-$d$-subpath of $R$.
Thus we know that the endpoints of $R$ are $a$ and $d$ and we may further assume w.l.o.g.\ that $R \in \mathcal{P}$, which implies $| \mathcal{P} | = 1$.

Should $| \mathcal{Q} | = 2$, we note that the two paths in $\mathcal{Q}$ must have $a$ and $b$, and respectively $c$ and $d$ as their endpoints, which makes finding the planar drawing easy, as $B_1$ is the only flap that has the potential to not be separable with two two adjacent vertices of $C$.
Thus, we may suppose further that $| \mathcal{Q} | = 1$ and let $Q$ be the unique path in $\mathcal{Q}$.
If $Q$ contains an $a$-$c$-subpath, we note that this subpath must be found in either $B_1$ or $B_2$, as it contains $aa'$ or $cc'$.
In case the $a$-$c$-subpath is found in $B_1$, we know that there also exists a subpath $Q''$ of $Q$ in $B_2$ that has $c$ as an endpoint and starts on the edge of $M$ covering $c$.
The other endpoint of $Q''$ cannot be $b$, as any path between $b$ and $c$ must be of odd length in $B_2$, due to their colours.
Thus $Q''$ is a $c$-$d$-path and $d$ must be an endpoint of $Q$, as there is no way to reach $b$ from this point.
Similarly, $a$ must be an endpoint of $Q$, as there is no way to reach $b$ from $a$ which still allows us to exit the flaps afterwards.
Since $\{ a, b \}$ separates $\sigma(s) - C$ and $B_0 - C$, we can therefore find a planar drawing as desired.

Thus we may assume that $B_2$ contains the $a$-$c$-subpath, implying that $B_1$ hosts an $a$-$b$-subpath of $Q$.
Note that the $a$-$b$-subpath of $Q$ cannot be extended to include $d$, for parity reasons, and $c$ and $d$ must be the endpoints of $Q$, which can be argued for as in the previous paragraph.
We conclude that $k=2$, as well as the fact that $\{ a,c,d \}$ separates $b$ from $B_2 - C$ in $B_2'$ and $\{ a,b,c \}$ separates $d$ and $B_1 - C$ in $B_1'$.
This allows us to find another planar drawing for the split of $\sigma(s)$.
Analogous arguments let us conclude that a $b$-$d$-subpath of $Q$ would also be beneficial to us.

Still under the assumption that $\mathcal{P}$ contains a single path that is comprised of a $c$-$d$-path in $B_2$ followed by an $a$-$c$-path in $B_1$ that has $b$ as an internal vertex, we may now assume that $| \mathcal{Q} | = 1$ and that $Q \in \mathcal{Q}$ visits the vertices of $U$ in the order provided by $C$.
For this reason, we know that the endpoints of $Q$ must also be adjacent on $C$.
Since any $B_i$ with $i \geq 3$ would only be visited by $Q$, we may also assume that $k = 2$ once more.
Further, if $Q$ only visits $B_1$, or only $B_2$, we can also easily find a planar drawing.
Should the endpoints of $Q$ be $a$ and $b$, then $Q$ may only enter $B_2$ via $c$ and $d$ and thus again make it easy to find the desired drawing.
If $c$ and $d$ are the endpoints of $Q$, we note that $\{ a,c,d \}$ separate the split of $\sigma(s)$ and $B_0$, which allows us to find a planar drawing, since we can move $b$ into $s$, as $a$ and $d$ are the endpoints of $P$, and $B_1$ can be separated from $d$ via $\{ a,b,c \}$.
The case in which $a$ and $d$ are the endpoints of $Q$ can be resolved similarly, as $\{ a,d \}$ separate the split of $\sigma(s)$ and $B_0$.
This means we must consider the case in which $b$ and $c$ are the endpoints of $Q$.
However, similar to many previous cases, if $B_2$ does not have a separator of size two consisting of two adjacent vertices on $C$, then
$\{ a,c,d \}$ separates $b$ from $B_2 - C$ in $B_2'$ and $\{ a,b,c \}$ separates $d$ and $B_1 - C$ in $B_1'$ and we are again done.
Thus no path in $\mathcal{P} \cup \mathcal{Q}$ has an internally $M$-conformal $b$-$c$-subpath contained in $B_1$ or, via an analogous argument, in $B_2$.
And using the same methods, we may also exclude the existence of an internally $M$-conformal $a$-$d$-subpath in $B_1$ or $B_2$.

Note that an $M$-conformal $b$-$c$- or $a$-$d$-subpath of any path in $\mathcal{P} \cup \mathcal{Q}$ cannot exist due to the arrangement of the edges of $M$ that cover the vertices in $U$.
Since any other $b$-$c$- or $a$-$d$-subpath of any path in $\mathcal{P} \cup \mathcal{Q}$ would have to be of even length, despite the fact that the endpoints of these paths have opposite colours, no such subpaths may exist at all.

Next, we narrow our assumptions to $k = 2$.
Towards the contrary, suppose there exists a $B_3$, which must have type iii), then we let $T \in \mathcal{P} \cup \mathcal{Q}$ be a path of which some part is hosted by $B_3$.
W.l.o.g.\ we can assume that $a$ is an endpoint of $T$.
Thus some part of $T$ is hosted by $B_1$.
To exit $B_1$ again, and then enter and exit $B_3$ afterwards, $T$ already needs to visit at least three vertices in $U$.
This means that we can assume w.l.o.g.\ that $T$ is the sole member of $\mathcal{P}$.
Note that $d$ cannot be the next vertex of $U$ visited after $T$ and suppose that instead $T$ visits $c$ next.
This forces $T$ to enter $B_2$, which it cannot exit via $b$, as that forces it back into $B_1$, depriving it of the chance to ever enter $B_3$.
Thus $T$ exits $B_2$ via $d$ and must then enter $B_3$, as it is its last chance to do so.
However, the only endpoint left-over from $U$ is $b$, forcing $T$ to enter $B_1$ again, which is impossible.
Thus $T$ must visit $a$ and then $b$.
It cannot stay within $B_1$ and visit $c$ or $d$, as that again deprives it of any chance to ever visit $B_3$, and thus it enters either $B_2$ or $B_3$.
Should $T$ enter $B_2$, it would have to visit $d$ next, then enter $B_3$ and finally exit via $c$, forcing it to enter $B_2$ once more, which is once more impossible.
Thus $T$ enters $B_3$ via $b$.
If $T$ then visits $c$, it must next visit $d$ and we note that $B_3$ contains an internally $M$-conformal $b$-$c$-subpath of $T$ and the situation is overall analogous to what we have previously discussed.

Therefore we have narrowed our situation to the case in which $T$ starts on $a$, moves through $B_1$ to find $b$, moves through $B_3$ to enter $B_2$ via $d$, forcing it to exit all flaps via $c$, which it reaches via an $M$-conformal $c$-$d$-path in $B_2$.
However, the $b$-$d$-subpath of $T$ in $B_3$ would have to be internally $M$-conformal and thus of odd length.
This is a contradiction to the fact that $b$ and $d$ have the same colour.
Thus no such path $T$ can exists, and we can assume that $k=2$, leaving us only with $B_1$ and $B_2$.

Suppose that some path $R \in \mathcal{P} \cup \mathcal{Q}$ contains an $a$-$c$-subpath in $B_1$.
This immediately implies that $R$ must have $a$ as an endpoint, as $R$ cannot next visit $c$ or $d$ in $B_2$ and visiting $b$ next would cause it to enter $B_1$ again in a position from which it cannot exit.
For similar reason, we can deduce that, after visiting $c$, $R$ can only visit $d$ next and then exit all flaps.
Thus the set containing $R$ only contains $R$ and, if the other set contains two paths, it is again easy to find a planar drawing of the contents of $s$.

We may therefore assume that $R \in \mathcal{P}$ and $Q \in \mathcal{Q}$ is the only representative of its set.
If $Q$ starts on $a$ as well, then it cannot visit $d$ next and visiting $c$ causes it to take the same shape as $R$, which makes finding the drawing easy.
Should $Q$ visit $a$ and then $b$, it must move towards $c$ next, as it cannot visit $d$, and thus create a conformal cross on $C$.
In case $Q$ has $b$ as an endpoint, the story is much the same, with $c$ being impossible to visit next, $b$ offering a simple way of finding the drawing and the remaining possibility forces $Q$ to visit $a$, then $b$, move to $c$ within $B_2$ and exit all flaps.
In this last case, we note that $\{ a,b,c \}$ separates $d$ and $B_1 - C$ and $\{ a,c,d \}$ separates $b$ and $B_2 - C$, a type of case we have discussed previously.
The remaining cases, in which $Q$ starts either on $c$ or $d$, are analogous.
Thus we conclude that by symmetry, we can assume that no $a$-$c$- or $b$-$d$-subpath of any path in $\mathcal{P} \cup \mathcal{Q}$ occurs in $B_1$ or $B_2$.

The remaining situation however only allows $a$-$b$- and $c$-$d$-subpaths for any and all paths in $\mathcal{P} \cup \mathcal{Q}$.
Since any path starting on $a$ or $b$ is immediately forced to enter $B_1$, and $c$ and $d$ similarly force the path to enter $B_2$, we are therefore done.
\end{proof}

The next lemma in this series will be the first of two cases in which it is possible to find an even dicycle.
These arguments will be a little more involved than the previous ones, leading us to introduce a few more terms and a lemma that will help us find the even dicycles in both cases.
To avoid speaking of $M$-directions at several points, we call an $M$-alternating path or cycle \emph{$M$-even} if it uses an even number of edges from $M$ and \emph{$M$-odd} otherwise.
More generally, we refer to this as the \emph{$M$-parity} of an $M$-alternating path or cycle.
Clearly, an $M$-alternating cycle is $M$-even if and only if its $M$-direction is an even dicycle.

We are now specifically interested in two $U$-paths hosted by the same flap, that cover the endpoints of $U$, and intersect somewhere in the flap.
First, we need to specify the setting.
Let $B'$ be a flap of type iii) or type vi), let $P_{ab}$ be an (internally) $M$-conformal path with the endpoints $a$ and $b$ hosted by $B'$ and let $Q_{cd}$ be an (internally) $M$-conformal path with the endpoints $c$ and $d$ hosted by $B'$.
Further, we assume that $P_{ab}$ and $Q_{cd}$ intersect within $B'$ and that $P_{ab}$ is the subpath of some path in $\mathcal{P}$, with $Q_{cd}$ being the subpath of some path in $\mathcal{Q}$.

Whether both paths are $M$-conformal or internally $M$-conformal depends on whether $B'$ has type iii) or type vi), but for our arguments we can suppose w.l.o.g.\ that $B'$ has type iii).
Since $P_{ab}$ and $Q_{cd}$ intersect, we can also also define the following paths:
Let $P_a$ be the $a$-$Q_{cd}$-subpath within $P_{ab}$, with $P_b$ being the $b$-$Q_{cd}$-subpath, and let $Q_c$ be the $c$-$P_{ab}$-subpath within $Q_{cd}$, with $Q_d$ being the $d$-$P_{ab}$-subpath within $Q_{cd}$.
Furthermore, let $u$ be the endpoint of $Q_c$ outside of $U$.

As is true for all flaps, $B'$ has a planar drawing inside of $s$ with $a$, $b$, $c$, and $d$ being the only vertices of $B'$ found on the boundary, in lexicographic order.
In the drawing of $B'$ the path $P_{ab}$ describes a curve that splits the disk $s$ into two disks $\Delta$ and $\Delta'$, one of which, say $\Delta$, contains both both $c$ and $d$.
Note that $Q_c$ must be drawn within $\Delta$.
Thanks to the many restrictions of our setting, we can provide a very specific description of the structure these two paths exhibit.

\begin{lemma}\label{lem:twopathstechnicalstructure}
    We can number the distinct components $H_1, \ldots , H_\ell$ of $P_{ab} \cap Q_{cd}$, for some even $\ell \in \N$ with $n \geq 2$, such that for each $i,j \in [\ell]$ with $i < j$ the path $H_i$ appears before $H_j$ when traversing $P_{ab}$ from $b$ to $a$, and each component is $M$-conformal.
    Further, for each $i \in [\ell - 1]$, the graph $P_{ab} \cap Q_{cd}$ contains an $M$-odd $M$-conformal cycle using $H_i$ and $H_{i+1}$.
\end{lemma}
\begin{proof}
    Let $w$ be the first vertex such that $u'Q_{cd}w$ is a $P_{ab}$-path.
    If $w$ lies on $bP_{ab}u'$, we note that $wP_{ab}u'$ and $u'Q_{cd}w$ are internally $M$-conformal and can be used to respectively replace each other in $C_1$ and $C_2$.
    Due to the minimality of $C_1 \cup C_2$, $wP_{ab}u'$ and $u'Q_{cd}w$ must therefore have the same $M$-parity.
    However, this means we can reduce the number of edges in $C_1 \cup C_2$, contradicting the minimality of this pair.
    
    Thus $w$ must lie on $uP_{ab}a$, with the matching edge that covers $w$ lying on $bP_{ab}w$.
    Note that this requires that $u'Q_{cd}w$ is drawn within $\Delta'$ and it also means that $u'P_{ab}w \cup u'Q_{cd}w$ forms an $M$-alternating cycle, that must be odd, since the $M$-direction of $B$ does not contain even dicycles inside of conjunction cells.

    Starting from these arguments, we observe that, if we were to traverse $Q_{cd}$ starting from $c$, we would start on $Q_c$, travel towards $b$ on $P_{ab}$, leave $P_{ab}$ on a non-matching edge and follow a subpath of $Q_{cd}$ drawn in $\Delta'$ that meets $P_{ab}$ on a vertex that is closer to $a$ on $P_{ab}$ than any vertex we have seen so far.
    This then either leads us directly to $Q_d$, or we next take a subpath of $Q_{cd}$ drawn in $\Delta$ that has an endpoint on $P_{ab}$, travel along $P_{ab} \cap Q_{cd}$ towards $b$, leave $P_{ab}$ followed by a subpath drawn in $\Delta'$, and then repeat the pattern stated in this sentence, until we finally reach $d$.
    This yields the desired structure.
\end{proof}

Note that $u$ must be an endpoint of $H_1$ and let $u'$ be the other endpoint.
Let $v$ and $v'$ be the endpoints of $H_\ell$ such that $v'$ lies closer to $b$ on $P_{ab}$ than $v$.
The following lemma will only become helpful in our last case, but relies still lives within the same setting as the previous lemma.
Thus it is convenient to state it here.

\begin{lemma}\label{lem:twopathstechnicalparitycondition}
    The paths $H_1$ and $H_\ell$ have the same $M$-parity if and only if $uP_{ab}v'$ and $u'Q_{cd}v$ have opposite $M$-parities.
\end{lemma}
\begin{proof}
    Let us prove this by induction over $\ell$.
    For $\ell = 2$, this is clearly true, since $H_1$ and $H_2$ are only involved in a single $M$-alternating cycle, which must be $M$-odd.
    In general, note that we can use the induction hypothesis on $H_3$ and $H_\ell$, and thus we only have to concern ourselves with two $M$-alternating cycles: the one involving $H_1$ and $H_2$, and the other involving $H_2$ and $H_3$.
    We let $u_2$ and $u_2'$ be the two endpoints of $H_2$, with $u_2'$ lying closer to $b$ on $P_{ab}$ than $u_2$, and we let $u_3$ and $u_3'$ be the endpoints of $H_3$, with analogous restrictions.
    Suppose that $H_1$ and $H_\ell$ have the same parity.
    Then, analogously to the start of the induction, $uP_{ab}u_2'$ and $u'Q_{cd}u_2'$ have the same $M$-parity if and only if $H_1$ and $H_2$ have different $M$-parity.
    The same is true for $u_2P_{ab}u_3'$ and $u_2'Q_{cd}u_3$, and $H_2$ and $H_3$.
    Since $H_2$ and $H_3$ are contained in both $uP_{ab}v'$ and $u'Q_{cd}v$, their $M$-parity can be disregarded, when determining whether the two paths have opposite parity.
\end{proof}

With these lemmas in place, we can tackle the next lemma, in which finding an even dicycle is an option.

\begin{lemma}\label{lem:planartraceshiftingtypeiiiparttwo}
    If $B_0$ has type iii) and $B_1$ is a type vi) flap, there exists a planar drawing for the split of $\sigma(s)$ that is compatible with the drawing of the split of $S - \sigma(s)$ or $\sigma(s)$ contains an even dicycle.    
\end{lemma}
\begin{proof}
Any path in $\mathcal{P} \cup \mathcal{Q}$ must host some subpath within $B_1$, since all matching edges covering the vertices of $U$ lie in $B_1$.
If there are two paths in $\mathcal{P}$ or $\mathcal{Q}$, then according to the restrictions laid out, both paths reside solely in $B_1$, which easily allows us to find a drawing again.
We can therefore assume that $| \mathcal{P} | = | \mathcal{Q} | = 1$.

This implies that $k \leq 3$ and of course, if either of these paths is only found in $B_1$, we are again done.
Note that no flap can host $a$-$c$- or $b$-$d$-subpaths of either of the two paths, since these would have to be of odd length, contradicting the fact that the endpoints of the subpaths have the same colour.
Thus if $B_2$ and $B_3$ each respectively only host parts of one of the paths, we can easily find the desired drawing.

Let $P \in \mathcal{P}$ and w.l.o.g.\ assume that $P$ starts in $a$ and visits $b$ next.
If $P$ then exits all flaps, we are done.
Thus $P$ must also visit $c$, which must occur via a subpath that w.l.o.g.\ lies in $B_2$.
Since visiting $c$ forces $P$ back into $B_1$, the remaining vertex $d$ must also be visited, which prompts $P$ to exit all flaps afterwards.
Therefore $P$ is an $a$-$d$-path.

Let $Q \in \mathcal{Q}$ and note that no matter where $Q$ starts, it must have exactly two $M$-conformal subpaths in $B_1$ that connect two pairs of vertices from $U$ that are adjacent in $C$.
In particular, since $Q$ operates under the same constraints as $P$, we know that $Q$ can only enter one flap other than $B_1$ and if this flap is not $B_2$, we are done.
Thus we can assume that $k=2$ and that $Q$ enters $B_2$.
If $Q$ has the same endpoints as $P$, we can easily find a planar drawing.
Should $Q$ instead share only one endpoint with $P$, say $a$, with $b$ being the other, then we note that $c$ is not an endpoint of either of the two paths, again allowing us to arrange the two flaps inside of $s$ in such a way that we can draw them in a planar fashion.

Finally, if $Q$ instead has the endpoints $b$ and $c$, we let $P_{ab}, P_{cd}$ be the two subpaths of $P$ comprising $P \cap B_1$, such that $P_{ab}$ has the endpoints $a$ and $b$, and let $P_{bc} = P \cap B_2$.
Similarly, we let $Q_{ab}, Q_{cd}$ be the two subpaths of $Q$ in $Q \cap B_1$, with $Q_{ab}$ being the path with the endpoints $a$ and $b$, and let $Q_{ad} = Q \cap B_2$.
Note that if $P_{bc}$ and $Q_{ad}$ do not intersect, we can easily find a drawing.
Analogously, if no pair of paths from $P_{ab}, P_{cd}, Q_{ab}, Q_{cd}$ intersects, we are also done.
W.l.o.g.\ we assume that $Q_{cd}$ intersects $P_{ab}$.

Note that both $B_1$ and $B_2$ contain the appropriate setting to let us apply \Cref{lem:twopathstechnicalstructure}.
First, let us fix $H_1, \ldots , H_\ell$ as the names for the components of $P_{ab} \cap Q_{cd}$ and let $I_1, \ldots , I_r$ be the components of $P_{bc} \cap Q_{ad}$, respectively numbered in accordance with the statement of \Cref{lem:twopathstechnicalstructure}.
We let $u$ and $u'$ be the endpoints of $H_1$, such that $Q_c = cQ_{cd}u$ is a $c$-$P_{ab}$-path and $P_b = bP_{ab}u'$ is a $b$-$Q_{cd}$-path.
Similarly, we let $v$ and $v'$ be the endpoints of $H_\ell$, with $Q_d = dQ_{cd}v'$ and $P_a = aP_{ab}v$ being a $d$-$P_{ab}$- and an $a$-$Q_{cd}$-path respectively.

Furthermore, there exists an $a$-$d$-path $Y$ such that $Y \cup P = C_1$ and a $b$-$c$-path $X$ with $X \cup Q = C_2$.
Both of these paths are found entirely in $B_0$.
Together with $Q_cuH_1u'P_b$, $X$ forms an $M$-conformal cycle that must be $M$-odd due to the minimality of $C_1$ and $C_2$.
The same is true for $Q_dv'H_\ell vP_a$ and $Y$, letting us conclude that both of these pairs have opposing $M$-parities.

Since $Q_{cd}$ intersects $P_{ab}$ and $B_1$ has a planar drawing on $s$, with the vertices $a,b,c,d$ found in the given order on the boundary, we know that $Q_{ab}$ and $P_{cd}$ do not intersect.
Thus $Q_{ab}$ and $P_{ab}$, respectively $Q_{cd}$ and $P_{cd}$, must have the same $M$-parity.

Suppose $Y$ has even $M$-parity.
Then the $M$-parity of $Q_dv'H_\ell vP_a$ must be odd.
Should the $M$-parity of $P_{ab}$ and $P_{cd}$ be opposite, then in total $bP_{ab}aYdP_{cd}c$ has odd $M$-parity, which forces $P_{bc}$ to have odd $M$-parity as well, since $bP_{ab}aYdP_{cd}c \cup P_{bc} = C_1$.
Thus $X$ must have odd $M$-parity, as otherwise $Q_cuH_1u'P_b$ has odd $M$-parity, which would mean that $cQ_cuH_1u'P_bb \cup P_{bc}$ is a $M$-conformal $M$-even cycle.
As the $M$-parities of $P_{ab}$ and $P_{cd}$ are opposite, the same must be true for $Q_{ab}$ and $Q_{cd}$.
This lets us conclude that $Q_{ad}$ has odd $M$-parity, as $Q_{ab} \cup Q_{ad} \cup Q_{cd} \cup X = C_2$.
However, now $Q_{ad} \cup dQ_dv'H_\ell vP_aa$ is an $M$-conformal $M$-even cycle.

Thus the $M$-parity of $P_{ab}$ and $P_{cd}$ must be equal.
Analogously to our previous argument this implies that $P_{bc}$ is $M$-even, $X$ is $M$-even, and $Q_cuH_1u'P_b$ is $M$-odd.
Again, we can then derive that $Q_dv'H_\ell vP_a$ is $M$-odd and $Q_{ad} \cup dQ_dv'H_\ell vP_aa$ is an $M$-conformal $M$-even cycle.
If we instead suppose that $Y$ has odd $M$-parity the exact same reasoning also always yields an $M$-conformal $M$-even cycle contained in the flaps.
This completes the proof.
\end{proof}

We now turn to a pair of simpler cases.

\begin{lemma}\label{lem:planartraceshiftingtypeiv}
    If $B_0$ has type iv), there exists a planar drawing for the split of $\sigma(s)$ that is compatible with the drawing of the split of $S - \sigma(s)$.    
\end{lemma}
\begin{proof}
    W.l.o.g.\ we let $ab \in M \cap E(C)$.
    Once again any path using only the edge $ab$ is not particularly interesting to us and if a path in $\mathcal{P}$ or $\mathcal{Q}$ starts on $a$ or $b$ and uses more than one edge, then it is the only element in its set.
    Analogously to the arguments in \Cref{lem:planartraceshiftingtypeii}, it is therefore easy to find a drawing for the insides of $s$.
\end{proof}

\begin{lemma}\label{lem:planartraceshiftingtypev}
    If $B_0$ has type v), there exists a planar drawing for the split of $\sigma(s)$ that is compatible with the drawing of the split of $S - \sigma(s)$.    
\end{lemma}
\begin{proof}
W.l.o.g.\ we let $a$ and $b$ be the two vertices in $U$ which are covered by edges in $E(B_0) \cap M$ and let $B_1$ be the type v) flap.
Note that any flap other than $B_1$ cannot contain an $a$-$c$-subpath for any of the paths in $\mathcal{P} \cup \mathcal{Q}$, since this subpath would have to be internally $M$-conformal and thus have odd length, contradicting the colours of $a$ and $c$.
The same holds for $b$-$d$-subpaths by symmetry.
Let $R \in \mathcal{P} \cup \mathcal{Q}$ be a path containing an $a$-$c$- or $b$-$d$-subpath $R'$, then we can assume that $R \neq R'$, as otherwise $R$ is the only path of its set and occupies only a single flap, which makes finding the drawing we search for easy.

W.l.o.g.\ assume that $R'$ is an $a$-$c$-path and note that $a$ must also be an endpoint of $R$, as the edge in $M$ covering $a$ lies outside of $\sigma(s)$.
Therefore, after leaving $B_1$, the path $R$ must visit a type iii) flap $B_2$ and from there either use a $c$-$d$- or a $b$-$c$-subpath $R''$ of $R$ hosted in $B_2$.
Since $R''$ being a $c$-$d$-path would force us to enter $B_1$ again, we conclude that $R''$ must end on $c$ and $b$, which forces us to leave the flaps and enter $B_0$.
Thus in this case $R$ is an $a$-$b$-path and also the only path in its set.
Note that the same argument also holds if $R'$ was a $b$-$d$-path, which would still make $R$ an $a$-$b$-path.

If both sets contain paths which have $a$-$c$- or $b$-$d$-subpaths respectively, we can therefore draw the interior of $s$, since the drawing must only attach to the drawing of $B_0$ at $a$ and $b$, allowing us to drawn $B_1$ in the center of $s$ and $B_2$ surrounding $B_1$, but still within $s$.
(Of course if there also exists a $B_3$, this is made even simpler.)
Suppose therefore that only $\mathcal{P}$ contains such a path and assume w.l.o.g.\ that this path contains an $a$-$c$-subpath hosted in $B_1$.
If $\mathcal{Q}$ contains only one path $T$, then we can assume that $T$ is not an $a$-$b$-path, as we can otherwise find a drawing as we just described.
Suppose that $T$ has $c$ as one endpoint and note that $T$ must therefore immediately enter $B_1$.
If $T$ then visits $a$ or $b$, we are done.
Thus $B_1$ contains a $c$-$d$-subpath of $T$, which is followed by an $a$-$d$-subpath in $B_2$ or $B_3$.
Should be $B_3$ exist, we are done, thus the $a$-$d$-subpath of $T$ lies in $B_2$.
However, now $\{ a,b,c \}$ separates $d$ and $B_1 - C$, allowing us to find a planar drawing for the insides of $s$.
In case $d$ is an endpoint of $T$, the situation is even simpler.

This allows us to assume that $\mathcal{Q}$ contains two paths.
Note that, for parity reasons, $B_1$ cannot contain $a$-$d$- or $b$-$c$-paths and thus $\mathcal{Q}$ contains an $a$-$b$-path and a $c$-$d$-path, providing no novel type of structure and thus again yielding a drawing as desired.
By symmetry, we conclude that no path in $\mathcal{P} \cup \mathcal{Q}$ contains an $a$-$c$- or $b$-$d$-subpath, as we are otherwise done.
However, as we just argued, this severely restricts the available structure to such a degree that finding the drawing is now trivial.
\end{proof}

\begin{lemma}\label{lem:planartraceshiftingtypevi}
    If $B_0$ has type vi), there exists a planar drawing for the split of $\sigma(s)$ that is compatible with the drawing of the split of $S - \sigma(s)$ or $\sigma(s)$ contains an even dicycle.    
\end{lemma}
\begin{proof}
Suppose first that $\mathcal{P}$ contains only a singular path $R$ and suppose w.l.o.g.\ that $R$ has $a$ as an endpoint.
We note immediately that when traversing $R$ starting from $a$, we must leave the flaps and enter back into $B_0$ once we see the second vertex of $U$.
Thus if either $\mathcal{P}$ or $\mathcal{Q}$ contain only a singular vertex, we can find the desired drawing.

Therefore we can assume that $| \mathcal{P} | = | \mathcal{Q} | = 2$ and according to the observation we have just made, each path is entirely hosted by a single flap.
Thus, if $k=1$ or $k \geq 3$, we are again done.
This leaves us in the case in which $k = 2$.
If two paths from the same set are hosted by the same flap or if in $B_1$, respectively $B_2$, the two paths from different sets do not intersect, we are again done.

This narrows the case down to the following setting: $k=2$, both sets contain two paths, and $B_1$, as well as $B_2$, each contain two intersecting $M$-alternating paths, which together cover all vertices in $U$.

Due to the symmetry of our setting, we can make the following assumptions:
Both $\mathcal{P}$ and $\mathcal{Q}$ contain an $a$-$b$-path $P_{ab} \subseteq B_1$, respectively $Q_{ab} \subseteq B_2$ and a $c$-$d$-path $P_{cd} \subseteq B_2$, respectively $Q_{cd} \subseteq B_1$.
For both sets this in particular means that within $C_1 \cup C_2$ there exists a $b$-$c$-path $W \subseteq C_1$, respectively $Y \subseteq C_2$, and an $a$-$d$-path $X \subseteq C_1$, respectively $Z \subseteq C_2$, both entirely found outside of $(B_1 \cup B_2) - U$, such that 
\[ W \cup P_{ab} \cup X \cup P_{cd} = C_1 \text{ and } Y \cup Q_{ab} \cup Z \cup Q_{cd} = C_2 . \]
Clearly, we are once again in the right territory to use \Cref{lem:twopathstechnicalstructure}.
For $B_1$, let $H_1, \ldots , H_\ell$ be the distinct components of $P_{ab} \cap Q_{cd}$, numbered according to the statement of \Cref{lem:twopathstechnicalstructure}.
Furthermore, we let $u$ and $u'$ be the endpoints of $H_1$, such that $Q_c = cQ_{cd}u$ is a $c$-$P_{ab}$-path and $P_b = bP_{ab}u'$ is a $b$-$Q_{cd}$-path.
Similarly, we let $v$ and $v'$ be the endpoints of $H_\ell$, with $Q_d = dQ_{cd}v'$ and $P_a = aP_{ab}v$ being a $d$-$P_{ab}$- and an $a$-$Q_{cd}$-path respectively.

For $B_2$, let $I_1, \ldots , I_r$ be the distinct components of $P_{cd} \cap Q_{ab}$, numbered according to the statement of \Cref{lem:twopathstechnicalstructure}.
We define $Q_a$, $Q_b$, $P_c$, and $P_d$ analogously to the paths in $B_1$ and label the endpoints of the just defined paths such that $x'$ is an endpoint of $I_1$ and $P_d$, $x$ is an endpoint of $I_1$ and $Q_a$, $y'$ is an endpoint of $I_r$ and $Q_b$, and $y$ is an endpoint of $I_r$ and $P_c$.
See \Cref{fig:typevidiagram} as a reference.

\begin{figure}
    \centering

    \scalebox{0.9}{
    \begin{tikzpicture}[scale=1]

        \pgfdeclarelayer{background}
		\pgfdeclarelayer{foreground}
			
		\pgfsetlayers{background,main,foreground}
			
        \begin{pgfonlayer}{main}
        \node (C) [v:ghost] {};

        \node (u) [v:main,position=135:8mm from C] {};
        \node (uprime) [v:mainempty,position=225:8mm from C] {};
        \node (v) [v:main,position=315:8mm from C] {};
        \node (vprime) [v:mainempty,position=45:8mm from C] {};

        \node (c) [v:mainempty,position=135:10mm from u] {};
        \node (b) [v:main,position=225:10mm from uprime] {};
        \node (a) [v:mainempty,position=315:10mm from v] {};
        \node (d) [v:main,position=45:10mm from vprime] {};

        \node (x_1) [v:main,position=225:6mm from c] {};
        \node (x_2) [v:mainempty,position=135:6mm from b] {};
        \node (x_3) [v:main,position=45:6mm from a] {};
        \node (x_4) [v:mainempty,position=315:6mm from d] {};

        \node (y_1) [v:main,position=135:16mm from c] {};
        \node (y_2) [v:mainempty,position=225:16mm from b] {};
        \node (y_3) [v:main,position=315:16mm from a] {};
        \node (y_4) [v:mainempty,position=45:16mm from d] {};

        \node (u_label) [v:ghost,position=75:2.5mm from u] {$u$};
        \node (uprime_label) [v:ghost,position=285:2.5mm from uprime] {$u'$};
        \node (v_label) [v:ghost,position=255:3mm from v] {$v$};
        \node (vprime_label) [v:ghost,position=105:3mm from vprime] {$v'$};

        \node (c_label) [v:ghost,position=45:2.5mm from c] {$c$};
        \node (b_label) [v:ghost,position=315:2.5mm from b] {$b$};
        \node (a_label) [v:ghost,position=225:3mm from a] {$a$};
        \node (d_label) [v:ghost,position=135:3mm from d] {$d$};

        \node (H1_label) [v:ghost,position=180:9mm from C] {$H_1$};
        \node (Hell_label) [v:ghost,position=0:9mm from C] {$H_{\ell}$};

        \node (y_label) [v:ghost,position=315:5.5mm from x_1] {\textcolor{Gray}{$Y$}};
        \node (z_label) [v:ghost,position=225:5.5mm from x_4] {\textcolor{Gray}{$Z$}};

        \node (w_label) [v:ghost,position=110:7mm from x_2] {$W$};
        \node (x_label) [v:ghost,position=70:7mm from x_3] {$X$};

        \node (I1_label) [v:ghost,position=185:10.5mm from x_1] {$I_1$};
        \node (I2_label) [v:ghost,position=355:10.5mm from x_4] {$I_2$};
        
        \end{pgfonlayer}{main}
        
        \begin{pgfonlayer}{foreground}
        \end{pgfonlayer}{foreground}

        \begin{pgfonlayer}{background}

        \draw [e:main,line width=3pt,color=BostonUniversityRed] (u.center) to (uprime.center);
        \draw [e:main,line width=3pt,color=BostonUniversityRed] (v.center) to (vprime.center);
        \draw [e:main,line width=3pt,color=BostonUniversityRed] (c.center) to (x_1.center);
        \draw [e:main,line width=3pt,color=BostonUniversityRed] (b.center) to (x_2.center);
        \draw [e:main,line width=3pt,color=BostonUniversityRed] (a.center) to (x_3.center);
        \draw [e:main,line width=3pt,color=BostonUniversityRed] (d.center) to (x_4.center);
        \draw [e:main,line width=3pt,color=BostonUniversityRed] (y_1.center) to (y_2.center);
        \draw [e:main,line width=3pt,color=BostonUniversityRed] (y_3.center) to (y_4.center);

        \draw [e:main,line width=2.2pt,color=CornflowerBlue] (u.center) to (vprime.center);
        \draw [e:main,line width=2.2pt,color=CornflowerBlue] (uprime.center) to (b.center);
        \draw [e:main,line width=2.2pt,color=CornflowerBlue] (v.center) to (a.center);
        \draw [e:main,line width=2.2pt,color=CornflowerBlue] (c.center) to (y_1.center);
        \draw [e:main,line width=2.2pt,color=CornflowerBlue] (d.center) to (y_4.center);
        \draw [e:main,line width=2.2pt,color=CornflowerBlue] (y_2.center) to (y_3.center);

        \draw [e:main,bend left=30,color=Gray] (x_1.center) to (x_2.center);
        \draw [e:main,bend left=30,color=Gray] (x_3.center) to (x_4.center);

        \draw [e:main] (x_1.center) to (x_2.center);
        \draw [e:main] (x_3.center) to (x_4.center);

        \draw [e:main,line width=2.6pt,color=Amber,dashed] (y_1.center) to (y_4.center);
        \draw [e:main,line width=2.6pt,color=Amber,dashed] (u.center) to (c.center);
        \draw [e:main,line width=2.6pt,color=Amber,dashed] (vprime.center) to (d.center);
        \draw [e:main,line width=2.6pt,color=Amber,dashed] (uprime.center) to (v.center);
        \draw [e:main,line width=2.6pt,color=Amber,dashed] (b.center) to (y_2.center);
        \draw [e:main,line width=2.6pt,color=Amber,dashed] (a.center) to (y_3.center);
        
        \end{pgfonlayer}{background}
        
    \end{tikzpicture}
    }
    \caption{A reference for the arguments in \Cref{lem:planartraceshiftingtypevi}, schematically depicting the interaction of two cycles in a trisum. In this lemma, the interactions of the parities of the individual parts of the cycles in the trisum are particularly involved.}

    \label{fig:typevidiagram}
\end{figure}
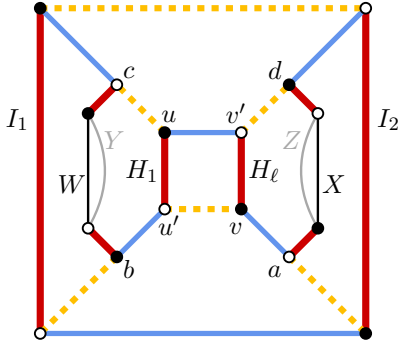

Note that that $Q_c$ is internally $M$-conformal and thus $bP_{ab}u'H_1uQ_cc$ is internally $M$-conformal, resulting in $bP_{ab}u'H_1uQ_cc \cup W$ and $bP_{ab}u'H_1uQ_cc \cup Y$ both being $M$-conformal cycles.
Clearly this implies that both of these cycles are $M$-odd, as otherwise $C_1 \cup C_2$ was not a smallest counterexample.
Thus the parity of the number of edges of $M$ in $W$ and $Y$ is equal.
The same argument holds for $X$ and $Z$.

Since outside of the flaps the part of $C_1$ and $C_2$ attached to $b$ and $c$, and $a$ and $d$ have equal parity, the difference between $C_1$ and $C_2$ must be found within the flaps.
For this reason, the paths $P_{ab}$ and $P_{cd}$ have the same $M$-parity if and only if $Q_{ab}$ and $Q_{cd}$ have opposite $M$-parity.
Therefore there must exist a flap in which the two involved paths also have opposite $M$-parity.
W.l.o.g.\ let $B_1$ be this flap and assume that the $M$-parity of $P_{ab}$ and $P_{cd}$ is the same.
Note that therefore the two paths in $B_2$ have the same $M$-parity.

Suppose that $W$ and $X$ have opposite $M$-parities.
We further suppose that $H_1$ and $H_\ell$ also have opposite $M$-parities, which means that the $M$-parities of $uP_{ab}v'$ and $u'Q_{cd}v$ are equal, according to \Cref{lem:twopathstechnicalparitycondition}.
Should $H_1$ have the same parity as $W$, we note that $bP_{ab}uQ_cc = P_bu'H_1uQ_cc$, and thus $P_b$ and $Q_c$ must have opposing $M$-parities, with the same holding for $P_a$ and $Q_d$ for symmetric reasons, because we assumed $C_1$ and $C_2$ to be a minimal counterexample.
This however implies that $P_{ab}$ and $Q_{cd}$ have the same parity, contradicting our assumptions.
If instead $H_1$ has the opposite parity of $W$, then $P_b$ and $Q_c$ must have the same $M$-parities, with the same also holding for $P_a$ and $Q_d$, again leading to the same contradiction.
We are therefore lead to assume that $H_1$ and $H_\ell$ indeed have the same $M$-parity, but this leads to a contradiction in an analogous fashion.

This means we can suppose instead that $W$ and $X$ have the same $M$-parity.
We now consider $B_2$ and suppose that $I_1$ and $I_r$ have opposite parities, which tells us that $xP_{cd}y'$ and $x'Q_{ab}y$ have the same $M$-parity by \Cref{lem:twopathstechnicalparitycondition}.
W.l.o.g.\ we can assume that $I_1$ and $X$ have the same parity, which implies that $P_d$ and $Q_a$ have opposing $M$-parities and in turn $P_c$ and $Q_b$ share their $M$-parity, which is a contradiction to the fact that $P_{cd}$ and $Q_{ab}$ have the same $M$-parity.
So we must instead consider the case in which $I_1$ and $I_r$ share their $M$-parity, but this can be resolved in a contradiction as all previous cases were.

This means that there are no valid $M$-parities for $W$ and $X$, which in turn implies that $B_1$ and $B_2$ cannot both contain intersecting paths.
Since we were able to find drawings in all other cases, our proof is therefore complete.
\end{proof}

With this we have checked all possible cases for \Cref{lem:planartraceshifting} and the lemma therefore follows.

\subsection{Shifting a non-planar even dicycle}\label{subsec:nonplanarshifting}

For the proof of the non-planar shifting lemmas, we will need a matching theoretic version of the famous two paths theorem.
In the undirected setting, this theorem was proven several times over the course of several decades (see \cite{giannopoulou2021two} for an overview).
A variant of this theorem in the matching theoretic setting was recently proven by Giannopoulou and Wiederrecht \cite{giannopoulou2021two}.

In directed graphs there is no true equivalent of the Two Paths Theorem.
This makes it somewhat surprising that there exists a statement that provides a version of this theorem in the matching theory of bipartite graphs (see \Cref{lem:crossoverc4}), as this result can at least somewhat be translated into the directed setting.
However, to really find crosses or useful embeddings for a specific digraph, we need much more infrastructure around the cycle that we want to view as the outer face.
Recall the definition of a conformal cross from \Cref{def:conformalcross}.

\begin{definition}[Reinforced society]\label{def:reinforcedsociety}
    A \emph{reinforced society} is a tuple $\mathcal{S} = ( D, H, \Psi , \Omega )$ where
    \begin{itemize}
        \item $D$ and $H$ are strongly connected digraphs with $H$ being a subgraph of $D$,
        
        \item $(D, \Psi)$ and $(D - ( H - V(\Omega) ) , \Omega )$ are societies,
        
        \item $H$ contains a butterfly minor model of a cylindrical grid of order 3 with the pairwise disjoint, homogeneous dicycles $C_1$, $C_2$, and $C_3$, such that $V(\Psi) = V(C_1)$ and $V(\Omega) = V(C_3)$,
        
        \item $( V(D - ( H - V(\Omega) ) ) , H )$ is an undirected separation with $V(\Omega)$ being the separator, and
        
        \item $H$ has a maelstrom-free odd decomposition $\rho$ within a disk $\Delta$ that is simultaneously an odd rendition of $(H, \Psi , \Omega)$ in $\Delta$.
    \end{itemize}
    Let $B$ be the split of $D$ and let $C \subseteq B$ be the undirected cycle corresponding to the split of $C_1$.
    We say that $\mathcal{S}$ has a \emph{cross in the split}, if $B$ contains a conformal cross on $C$.
\end{definition}

Having a non-even rendition is equivalent to the notion of the split of the digraph being matching flat, as defined in \cite{Giannopoulou2023ExcludingSingleCrossing}.
Therefore we can translate Theorem 4.6 in the full version \cite{giannopoulou2022excluding} of \cite{Giannopoulou2023ExcludingSingleCrossing} as follows.

\begin{theorem}[Giannopoulou, Thilikos, and Wiederrecht \cite{Giannopoulou2023ExcludingSingleCrossing}]\label{thm:reinforcedmatcross}
    Let $\mathcal{S} = ( D, H, \Psi , \Omega )$ be a reinforced society.
    The society $(D, \Psi)$ has a non-even rendition if and only if $\mathcal{S}$ does not have a cross in the split.
\end{theorem}

To apply \Cref{thm:reinforcedmatcross}, we need a directed analogue of a technical result, which appears in the proof of Theorem 4.6 in \cite{giannopoulou2022excluding}.
This lemma allows us to get a non-even rendition for a reinforced society, if the graph in it is non-even.
Translating this to the world of digraphs, in \cite{giannopoulou2022excluding} the problem is first reduced to the dibraces of the graph (Claim 1 to 3) and then resolved using \hyperref[thm:nonevenstructure]{the decomposition theorem for non-even dibraces} (Claim 4), which on its own already almost provides the desired decomposition.

\begin{lemma}\label{lem:nonevenhasnonevenrendition}
    Let $D$ be a non-even digraph and let $( D, H, \Psi, \Omega )$ be a reinforced society.
    Then $(D, \Psi)$ has a non-even rendition.
\end{lemma}

We will also need a tool that lets us actually find a cylindrical grid of order three, which is a requirement for building a reinforced society.
For this purpose, we will again lean on a type of grid theorem.
However, this result relies on a graph that is already very structured.
In particular, due to the way we will later build our outlines, which starts from a large wall and then attaches large transactions to said wall, we can additionally assume that there exist two large linkages, connecting the innermost and outermost dicycles of a set of homogeneous concentric dicycles.
The following definition captures an idealised version of this structure.

\begin{definition}[Segregated grid]
    Let $k$ be a positive integer.
    A \emph{segregated grid} of order $k$ is a digraph consisting of the cycles $C_1, \ldots , C_k$, with 
    \[ C_i = \Brace{v_0^i, e^i_0, v_1^i, e_1^i, \ldots , e_{2k-3}^i, v_{2k-2}^i, e_{2k-2}^i, v_{2k-1}^i, e_{2k-1}^i, v_0^i} \]
    for each $i \in [k]$, by adding the directed paths $P_i = v_i^1 v_i^2 \ldots v_i^{k-1} v_i^k$ and $Q_i = v_{i+k}^k v_{i+k}^{k-1} \ldots v_{i+k}^2 v_{i+k}^1$ for every $i \in [0, k-1]$.
\end{definition}

In the more abstract setting of a set of homogeneous dicycles, we will often require a particular pair of linkages to be present.
Note that this pair can easily be found in a large cylindrical wall.

\begin{definition}[Segregated $\mathcal{C}$-pair]\label{def:segregatedpair}
Let $D$ be a digraph, let $k , t$ be integers, let $\delta = (\Gamma, \mathcal{V}, \mathcal{D})$ be a strong plane decomposition of $D$, let $\mathcal{C} = \{ C_1, \ldots , C_k \}$ be a homogeneous set of dicycles in $D$, and let $R$ be the connected region that is bounded by the traces of both $C_1$ and $C_k$.
We call a pair $(\mathcal{L}_1, \mathcal{L}_2)$ consisting of a $C_1$-$C_k$-linkage $\mathcal{L}_1$ and a $C_k$-$C_1$-linkage $\mathcal{L}_2$, each of order $t$, a \emph{segregated $\mathcal{C}$-pair} of order $t$, if all paths in $\mathcal{L}_1 \cup \mathcal{L}_2$ are drawn in $R$, the paths in $\mathcal{L}_1 \cup \mathcal{L}_2$ are pairwise disjoint, both $\mathcal{L}_1$ and $\mathcal{L}_2$ have order $t$, and $C_1$ and $C_k$ can each be partitioned into two paths such that the endpoints of $\mathcal{L}_1$ lie on one path and the endpoints of $\mathcal{L}_2$ lie on the other.
\end{definition}

The proof of the planar directed grid theorem tells us that, given a large set of homogeneous dicycles and a large segregated pair for these dicycles, we can easily find a large segregated grid.
In particular, for our specific setting, the proof can be modified with little effort to yield a segregated grid that separates the two disks that are separated by all the homogeneous dicycles in the original set.

\begin{theorem}[Hatzel, Kreutzer, and Kawarabayashi \cite{hatzel2019polynomial}]\label{lem:findseggrid}
    Let $k$ be an integer.
    There exists a function $\PolyPlanarGridNoArg : \N \rightarrow \N$ with $\PolyPlanarGridNoArg \in \mathcal{O}(k^6)$, such that if $D$ is a digraph, with a strong plane decomposition $\delta = (\Gamma, \mathcal{V}, \mathcal{D})$, containing a homogeneous set of dicycles $\mathcal{C} = \{ C_1, \ldots , C_{\PolyPlanarGrid{k}} \}$, together with a segregated $\mathcal{C}$-pair $(\mathcal{L}_1, \mathcal{L}_2)$ of order $\PolyPlanarGrid{k}$, and with $\Delta_1$, respectively $\Delta_{\PolyPlanarGrid{k}}$, being the disk bounded by the trace of $C_1$, respectively $C_{\PolyPlanarGrid{k}}$, whose interior does not contain any cycle in $\mathcal{C}$.

    Then $D$ contains a segregated grid $G$ of order $k$ as a butterfly minor and the trace of each dicycle in $G$ separates $\Delta_1$ and $\Delta_{\PolyPlanarGrid{k}}$.
\end{theorem}

For the reinforced society we wish to build, we actually need a cylindrical grid.
This can be done directly using what we are given by the above lemma, via a construction from \cite{hatzel2019polynomial}.

\begin{lemma}[Hatzel, Kreutzer, and Kawarabayashi \cite{hatzel2019polynomial}]\label{lem:seggridtoseparatingcylgrid}
    Let $k$ be a positive integer, let $H$ be a segregated grid of order $k^2$, and let $\Delta_1$, respectively $\Delta_{k^2}$, be the disk bounded by the trace of $C_1$, respectively $C_{k^2}$, whose interior does not contain $H$.
    
    Then $H$ contains a cylindrical grid $G$ of order $k$ as a butterfly minor and the trace of each dicycle in $G$ separates $\Delta_1$ and $\Delta_{k^2}$.
\end{lemma}

\begin{corollary}\label{cor:getseparatingcylindricalgrid}
    Let $k$ be an integer and let $\PlanarGridNoArg : \N \rightarrow \N$ be defined such that $\PlanarGrid{k} = \PolyPlanarGrid{k^2}$.
    If $D$ is a digraph, with a strong plane decomposition $\delta = (\Gamma, \mathcal{V}, \mathcal{D})$, containing a homogeneous set of dicycles $\mathcal{C} = \{ C_1, \ldots , C_{\PlanarGrid{k}} \}$, together with a segregated $\mathcal{C}$-pair $(\mathcal{L}_1, \mathcal{L}_2)$ of order $\PlanarGrid{k}$, and with $\Delta_1$, respectively $\Delta_{\PlanarGrid{k}}$, being the disk bounded by the trace of $C_1$, respectively $C_{f(k)}$, whose interior does not contain any cycle in $\mathcal{C}$.

    Then $D$ contains a cylindrical grid $G$ of order $k$ as a butterfly minor and the trace of each dicycle in $G$ separates $\Delta_1$ and $\Delta_{\PlanarGrid{k}}$.
\end{corollary}

This means there exists a concrete constant ensuring that we can find a cylindrical grid of order three in a circle rim together with an segregated $\mathcal{C}$-pair of order at least this constant.
For (refined) diamond rims we will need to take twice that value.
Thus we define the following constant to be used throughout the rest of the paper.

\begin{definition}[Shifting constant]\label{def:rootingconstant}
    We define $\SC = 2 \PlanarGrid{3} + 1$.
\end{definition}

Having gathered all the necessary tools, we now prove the first of our properly non-planar shifting lemmas on a much stricter version of outlines called \emph{rims}.
In the case of circle rims, the big distinction from an outline is that we get rid of all transactions in the infrastructure of our outline.
We will provide the definitions for diamond rims and refined diamond rims later on in this section.

\begin{definition}[Rim of a circle maelstrom]\label{def:circlerim}
Let $D$ be a digraph, let $\delta = (\Gamma, \mathcal{V}, \mathcal{D})$ be an odd decomposition of $D$, let $m \in C(\delta)$ be a maelstrom of $\delta$, and let $\theta$ be a positive integer.
A tuple $\mathfrak{R} = (H, \mathcal{C})$, where $\mathcal{C} = \{ C_1 , C_2, \ldots , C_\theta \}$ is a homogeneous family of pairwise disjoint dicycles, is a \emph{circle $\theta$-rim} of $m$, if the points \textbf{\textsf{M1}}, \textbf{\textsf{M2}}, and \textbf{\textsf{M7}} of \Cref{def:outline} hold for $\mathfrak{R}$, provided we set $t = 0$, and $C_1$ has an $m$-disk in $\delta$.

We call $(H_\mathfrak{R}, \Omega_\mathfrak{R})$ the \emph{$\mathfrak{R}$-society of $m$}, where $V(\Omega_\mathfrak{R}) = V(C_\theta)$ and $H_\mathfrak{R} = H \cup \sigma(m)$.
\end{definition}

To be able to apply these lemmas on our outlines we of course first need to argue that we can find these rims inside normal outlines.
This will be done in \cref{sec:localising}.

We can now apply \Cref{lem:nonevenhasnonevenrendition} for the case in which our circle rim contains a cylindrical grid of order three and we have a cross on the outermost cycle to find a weak odd bicycle.
\begin{lemma}\label{lem:keylemmashifting}
    Let $\theta$ be an integer, let $D$ be a digraph, let $\delta = (\Gamma, \mathcal{V}, \mathcal{D})$ be an odd decomposition of $D$, and let $m \in C(\delta)$ be a maelstrom of $\delta$.
    If $(H, \mathcal{C} = \{ C_1 , \ldots , C_\theta \} )$ is a circle $\theta$-rim of $m$ such that $L,R \subseteq H$ form a cross over $C_\theta$ and $H$ contains a butterfly minor model $G$ of a cylindrical grid of order three, whose three pairwise disjoint concentric dicycles $C_1',C_2',C_3'$ are all found in $\mathcal{C}$, then $H$ contains a weak odd bicycle. 
\end{lemma}
\begin{proof}
    Without loss of generality, we may suppose that $H$ is non-even, that the $m$-tight disk of $C_3'$ contains $C_2'$ and $C_1'$, and that the $m$-tight disk of $C_1'$ does not contain $C_2'$. 
    Consider one of the two natural cyclical orderings $\Psi$ induced by $C_1'$ and an analogous ordering $\Omega$ induced by $C_3'$, and let $H' \subseteq H$ be the digraph drawn by $\delta$ into the region bounded by the traces of $C_1'$ and $C_3'$.
    Using \Cref{def:mdir}, we observe that $( D, H', \Psi, \Omega )$ is a reinforced society with a cross in the split, since we can shorten $L$ and $R$ to start on $C_3'$.
    Thus there does not exist a non-even rendition of $( D, \Psi )$, according to \Cref{thm:reinforcedmatcross}.
    Note that $D$ cannot be isomorphic to $F_7$, since it contains a butterfly minor model of the cylindrical grid of order 3 and thus at least 18 vertices.
    Since $D$ is non-even, this leads to a contradiction, as \Cref{lem:nonevenhasnonevenrendition} guarantees us a non-even rendition for $(D, \Psi)$ given our assumptions.
\end{proof}

With the help of all of these tools, we can now prove a version of the non-planar shifting lemma for circle rims with a segregated $\mathcal{C}$-pair.

Let $(D, \Omega)$ be a society and for each $i \in [2]$, let $P_i \subseteq D$ be a directed path with tail $s_i$ and head $t_i$.
The pair $P_1, P_2$ is called a \emph{crossing pair} for $(D, \Omega)$ if $P_1$ and $P_2$ are disjoint and the vertices $s_1, s_2, t_1, t_2$ occur on $\Omega$ in the order listed.
\begin{lemma}[Circle non-planar shifting]\label{lem:circlenonplanarshifting}
    Let $D$ be a digraph, let $\delta = (\Gamma, \mathcal{V}, \mathcal{D})$ be an odd decomposition of $D$, let $m \in C(\delta)$ be a maelstrom of $\delta$, and let $\mathfrak{R} = (H, \mathcal{C} )$ be a circle $\lfloor \nicefrac{\SC}{2} \rfloor$-rim of $m$, with a segregated $\mathcal{C}$-pair of order $\lfloor \nicefrac{\SC}{2} \rfloor$ and a crossing pair $L,R$ for the $\mathfrak{R}$-society $(H_\mathfrak{R}, \Omega_\mathfrak{R})$ of $m$.
    Then the digraph $H_\mathfrak{R}$ contains a weak odd bicycle.
\end{lemma}
\begin{proof}
    Note that $\lfloor \nicefrac{\SC}{2} \rfloor$ is large enough that \Cref{cor:getseparatingcylindricalgrid} nets us a butterfly minor of the cylindrical grid of order three within $H_\mathfrak{R}$ that surrounds the maelstrom $m$.
    Thus \Cref{lem:keylemmashifting} tells us that $H_\mathfrak{R}$ contains a weak odd bicycle.
\end{proof}

We can immediately derive a powerful corollary from this lemma.

\begin{corollary}\label{cor:circleapplyshifting}
    Let $t$ be an integer, let $D$ be a digraph, let $\delta = (\Gamma, \mathcal{V}, \mathcal{D})$ be an odd decomposition of $D$, let $m \in C(\delta)$ be a maelstrom of $\delta$, and let $\mathfrak{R} = (H, \mathcal{C}, \mathcal{E} )$ be a circle $t \lfloor \nicefrac{\SC}{2} \rfloor$-rim of $m$, with an segregated $\mathcal{C}$-pair $(\mathcal{L}_1, \mathcal{L}_2)$ of order $\lfloor \nicefrac{\SC}{2} \rfloor$ and a transaction $\mathcal{L} = \{ L_1, \ldots , L_{2t} \}$ on the $\mathfrak{R}$-society $(H_\mathfrak{R}, \Omega_\mathfrak{R})$, such that for all $i \in [t]$ the paths $L_{2i-1}, L_{2i}$ are a crossing pair for $\mathfrak{R}$-society.
    
    Then the digraph $H_\mathfrak{R}$ contains a half-integral packing of $t$ weak odd bicycles.
\end{corollary}
\begin{proof}
    We can simply divide $\mathcal{C}$ up into $t$ set of dicycles $\mathcal{C}_1, \ldots, \mathcal{C}_t$ such that $\mathcal{C}_i = \{ C_{i\lfloor \nicefrac{\SC}{2} \rfloor} , \ldots , C_{(i+1)\lfloor \nicefrac{\SC}{2} \rfloor - 1 } \}$ for $i \in [t]$.
    For all $i \in [t]$, we can find an segregated $\mathcal{C}_i$-pair within the linkages $\mathcal{L}_1$ and $\mathcal{L}_2$.
    Applying \Cref{lem:circlenonplanarshifting} thus yields $t$ weak odd bicycles, which must constitute a half-integral packing since the elements of $\mathcal{C}$, respectively $\mathcal{L}$, are pairwise disjoint.
\end{proof}

Before we proceed to consider (refined) diamond outlines, we will first concern ourselves with a technical statement claiming that two linkages that cross each other in an almost planar way, can be used to reroute the paths of one linkage along the paths of the other linkage without any loss in the order of the resulting linkage.
Similar statements can be found in Section 4 of \cite{kawarabayashi2020quickly} (see in particular Lemma 4.6).
Our approach however often requires us to preserve the size of the linkages exactly.

\begin{lemma}[Rerouting]\label{lem:detour}
    Let $k \in \N$ be some positive integer and let $\mathcal{P} = \{ P_i \}_{i \in [k]}$ and $\mathcal{Q} = \{ Q_i \}_{i \in [k]}$ be linkages with $V(P_i) \cap V(Q_j) \neq \emptyset$ for all $i,j \in [k]$ and such that the heads of the paths in $\mathcal{P}$ are all found in $V(Q_k)$.
    Furthermore, let $t_i$ be the tail of $P_i$ and $h_i$ be the head of $Q_i$, for all $i \in [k]$.
    Let $\delta = (\Gamma, \mathcal{V}, \mathcal{D})$ be a pure, maelstrom-free odd decomposition of $H := \bigcup_{i \in [k]} (P_i \cup Q_i)$ in a disk $\Delta$, such that all paths in $\mathcal{P} \cup \mathcal{Q}$ are grounded, the endpoints of all paths in $\mathcal{P} \cup \mathcal{Q}$ are drawn on the boundary of $\Delta$, and for all $i \in [2,k-1]$ the trace of $P_i$ separates $P_{i-1}$ and $P_{i+1}$, and the trace of $Q_i$ separates $Q_{i-1}$ and $Q_{i+1}$.

	Then there exists a $\{ t_i \}_{i \in [k]}$-$\{ h_i \}_{i \in [k]}$-linkage $\{ R_1, \ldots , R_k \} $ of order $k$ in $H$, such that for each $i \in [k]$ the path $R_i$ has $t_i$ as its tail and $h_i$ as its head, and for all $i \in [2, k-1]$, the trace of $R_i$ separates $R_{i-1}$ and $R_{i+1}$.
\end{lemma}
\begin{proof}
    In this proof, we want to treat the members of $\mathcal{P} \cup \mathcal{Q}$ as if they were drawn in a planar way.
    For this purpose, we first make a few observations.
    
    For the sake of rerouting, we can ignore big vertices and the directed tight cuts associated with them, since at most two paths may intersect there, one from $\mathcal{P}$ and the other from $\mathcal{Q}$.
    Furthermore, we observe that at most two paths from $\mathcal{P}$, and respectively $\mathcal{Q}$, may be partially drawn by $\Gamma$ on some conjunction cell $c \in \mathcal{D}$, since the boundary of $c$ contains at most four vertices.
    In fact, if $\Boundary{c}$ contains at most three vertices, then the paths interacting with $c$ behave essentially like planar paths for our purposes.

    We may therefore suppose that the boundary of $c$ contains four vertices.
    Clearly, since the traces of the elements of $\mathcal{P}$ and respectively $\mathcal{Q}$ separate each other, we just need to discuss the following situation.
    Suppose that some subpath $P$ of a path in $\mathcal{P}$, respectively $\mathcal{Q}$, is drawn on $c$ such that the two components of $\Boundary{c}- V(P)$ each contain a vertex.
    In this case, we note that $c$ must contain an even dicycle, which can be deduced by applying \Cref{lem:smallcycsumandtrisum} and \Cref{lem:crossoverc4}, noting that the conformal $K_{3,3}$ corresponds to an odd bicycle and then applying \Cref{obs:oddbicycleevendicycle}.
    This contradicts the fact that $\delta$ is pure.
    We can therefore treat the paths in $\mathcal{P}$ and $\mathcal{Q}$ as if they were drawn by $\Gamma$ in an essentially planar way.

    W.l.o.g.\ we can assume that $\mathcal{P}$ is a $\{ t_i \}_{i \in [k]}$-$V(Q_k)$-linkage (otherwise we can shorten the paths appropriately).
	We will prove by induction over all integers $k$ and all $i \in [k]$ that we can choose $w_i \in V(P_i)$ for each $i \in [k]$ such that
	\begin{enumerate}
		\item\label{item:pathsorti}$w_i \in \V{Q_i}$,
		
		\item\label{item:pathsortii}$w_iQ_i$ is disjoint from $P_j$ with $j \in [i+1, k]$, and
		
		\item\label{item:pathsortiii}$w_iQ_i$ is disjoint from $P_{i-1}w_{i-1}$ if $i > 1$.
	\end{enumerate}
	In the case $k=1$, the vertex $w_1$ is simply the head of $P_1$, since $P_1$ is a $t_1$-$Q_k$-path.
    For $k > 1$ and $i \in [k]$ with $i < k$, we proceed backwards from $Q_k$ and first use the induction hypothesis on $\bigcup_{j \in [i+1, k]} (P_j \cup Q_j)$ to find $w_{i+1}, w_{i+2}, \ldots , w_k$.
    We let $w_i \in (\V{Q_i} \cap \V{P_i})$ be the first vertex of $P_i$ on $Q_i$ such that $w_iQ_i$ is disjoint from $P_{i+1}$, which must exist since the paths in $\mathcal{P}$ behave in a planar way.
    This satisfies \ref{item:pathsorti} and \ref{item:pathsortii}.
    
    Suppose now that $w_{i+1}Q_{i+1}$ intersects $P_iw_i$.
    Then, since the paths in $\mathcal{P}$ behave in an essentially planar way, there must exist a $w_i' \in V(P_iw_i) \cap V(Q_iw_i) \setminus \{ w_i \}$ maximising $P_iw_i'$, as some part of $P_i$ must properly cross through $Q_i$ in the drawing given by $\Gamma$ to touch $Q_{i+1}$.
    However, according to the choice of $w_i$, there must exist some $P_j$ with $j \in [i+1, k]$ such that $w_i'Q_iw_i$ is intersected by $P_j$.
    As $P_j$ cannot intersect $P_iw_i'$ and similarly cannot intersect $w_i'P_iw_i$, it must somehow reach $w_i'Q_iw_i$, requiring it to cross through $w_iQ_i$ in the drawing, which is a clear contradiction to our choice of $w_i$.
    Thus $w_{i+1}Q_{i+1}$ does not intersect $P_iw_i$ and we can also satisfy \ref{item:pathsortiii}.

    Let $R_i := P_i w_i Q_i$ for each $i \in [k]$.
    We note that if $R_i$ and $R_{i+1}$ do not intersect, then $R_i$ and $R_j$, for $j \in [k]$ with $j > i+1$, do not intersect either, because $\delta$ is a strong plane decomposition without maelstroms.
    We can suppose that $R_i$ and $R_{i+1}$ do in fact intersect for some $i \in [k-1]$, since otherwise $\{ R_i \}_{i \in [k]}$ is the linkage we are searching for.
    Thanks to assumptions on $\mathcal{P}$ and $\mathcal{Q}$, we know that the paths $P_iw_i$ and $P_{i+1}w_{i+1}$, and similarly $w_iQ_i$ and $w_{i+1}Q_{i+1}$, cannot intersect.
    Our choice of $w_i$ also tells us that $w_iQ_i$ and $P_{i+1}w_{i+1}$ are disjoint thanks to \ref{item:pathsortii}.
    Thus $w_{i+1}Q_{i+1}$ and $P_iw_i$ must intersect, but this directly contradicts \ref{item:pathsortiii} and concludes our proof.
\end{proof}


Unlike the last lemma, the next shifting argument will only require us to first define diamond rims.
Whilst diamond outlines require several transactions to describe the structure around a maelstrom, a diamond rim only requires a single transaction.

\begin{definition}[Diamond rim]\label{def:diamondrim}
Let $D$ be a digraph, let $\delta = (\Gamma, \mathcal{V}, \mathcal{D})$ be an odd decomposition of $D$, let $m \in C(\delta)$ be a maelstrom of $\delta$, and let $\theta$ be a positive integer.
A tuple $\mathfrak{R} = (H, \mathcal{C}, \mathcal{E})$, where $\mathcal{C} = \{ C_1, \ldots , C_\theta \}$ is a homogeneous family of pairwise disjoint dicycles and $\mathcal{E} = \{ E_1,  \ldots , E_\theta \}$ is a $C_\theta$-linkage, is called a \emph{diamond $\theta$-rim} of $m$, if the following requirements are met.

\begin{description}
    \item[DR1] The points \textbf{\textsf{M1}}, \textbf{\textsf{M2}}, \textbf{\textsf{M4}}, \textbf{\textsf{M5}}, \textbf{\textsf{M6}}, and \textbf{\textsf{M7}} of \Cref{def:outline} hold for $\mathfrak{R}$, provided we set $t = 1$, let $\mathcal{E}_1 = \mathcal{E}$, set $E_i^1 = E_i$ for all $i \in [\theta]$, and treat $\mathfrak{R}$ as a diamond outline.
    
    \item[DR2] The graph $C_1 \cup E_1$ contains an $m$-tight cycle in $\delta$.

\end{description}

We call $(H_\mathfrak{R}, \Omega_\mathfrak{R})$ the \emph{$\mathfrak{R}$-society of $m$}, where $V(\Omega_\mathfrak{R}) = V(Q_1)$ and $H_\mathfrak{R} = H \cup \sigma(m)$.
\end{definition}

It turns out that we can find plenty of different somewhat smaller circle rims in a large enough diamond rim.
Our goal will be to move any cross on the society of the diamond rim onto the society of one of these smaller circle rims, which lets us use \Cref{lem:circlenonplanarshifting}.
For this purpose, we will want to take a layer of $\SC$ cycles and paths from the infrastructure of a large diamond rim and then show that we can find a weak odd bicycle within them, if we also have a cross to go along with it.
This way we can get an analogous result to \Cref{cor:circleapplyshifting} for diamond rims.
Since this set of cycles and paths sadly does not form another diamond rim, we need to define this slightly more general structure.

\begin{definition}[Layers of a rim]\label{def:diamondlayer}
    Let $\theta$, $\ell$ and, $k$ be positive integers, let $D$ be a digraph, let $\delta = (\Gamma, \mathcal{V}, \mathcal{D})$ be an odd decomposition of $D$, let $m \in C(\delta)$ be a maelstrom of $\delta$, and let $\mathfrak{R} = (H, \mathcal{C} = \{ C_1, \ldots , C_\theta \} , \mathcal{E} = \{ E_1, \ldots , E_\theta \} )$ be a $\theta$-rim of $m$.
    We call the tuple $(H', \{ C_\ell , \ldots C_{\ell + k - 1} \}, \{ E_\ell , \ldots E_{\ell + k - 1} \} )$ a \emph{$k$-layer at $\ell$ in $\mathfrak{R}$}, where $H'$ is the union of $\bigcup_{j = \ell}^{\ell +k - 1} E_j$ and the subgraph of $H$ drawn on the intersection of the largest $m$-disk bounded by a cycle in $C_{\ell + k - 1} \cup E_{\ell + k - 1}$ and the largest non-$m$-disk within $C_\ell \cup E_\ell$.
    We note that this yields a definition for the layers of a circle rim by ignoring the parts of the definition that have to do with $\mathcal{E}$.
\end{definition}

We first prove that we can find a weak odd bicycle if we have a cross on any $\SC$-layer.
\begin{lemma}\label{lem:diamondlayershifting}
    Let $\theta$ be an integer with $\theta \geq \SC$, let $D$ be a digraph, let $\delta = (\Gamma, \mathcal{V}, \mathcal{D})$ be an odd decomposition of $D$, let $m \in C(\delta)$ be a maelstrom of $\delta$, and let $\mathfrak{R} = (H, \mathcal{C}, \mathcal{E} )$ be a diamond $\theta$-rim of $m$, with a crossing pair $L_1, L_2 \subseteq H_\mathfrak{R}$ for the $\mathfrak{R}$-society $(H_\mathfrak{R}, \Omega_\mathfrak{R})$.
    Then for every $\SC$-layer $(H', \mathcal{C}', \mathcal{E}')$ there exists a weak odd bicycle in $H' \cup L_1 \cup L_2$.
\end{lemma}
\begin{proof}
    Let $\ell = \lfloor \SC \rfloor$.
    We note that the value of $\ell$ stems from the proof of \Cref{lem:circlenonplanarshifting}.
    Our strategy will be to define two types of circle $\ell$-rims within the $\SC$-layer and then show that we can find a cross on one of these two types of circle rims, allowing us to use \Cref{lem:circlenonplanarshifting} to find the weak odd bicycle.
    To ease notation, we will assume that we have chosen the innermost $\SC$ cycles and paths in $\mathfrak{R}$ as our $\SC$-layer.
    
    The first type of circle rim is defined as follows.
    We observe that for each $i \in [\ell , \SC]$ the paths in the linkage $\Set{E_1, E_2, \ldots , E_\ell } \subseteq \mathcal{E}$ contain both a $C_{i - (\ell-1) }$-$C_i$-linkage and a $C_i$-$C_{i - (\ell -1)}$-linkage, each of order $\ell$.
    Let the corresponding circle $\ell$-rim induced by $\{ C_{i-(\ell -1)}, \ldots , C_i \}$, together with the segregated pair of linkages, be called $\mathfrak{U}_i$ and note that $C_i$ is the dicycle that defines the $\mathfrak{U}_i$-society.
    We will call these type of circle rims the \emph{type I rims}.

    For the next type, we first require a general observation on how the paths in $\mathcal{E}$ can be assumed to behave with respect to the contents of $\mathcal{C}$.
    Let $C$ be the $m$-tight cycle in $C_{\ell+1} \cup E_\ell$ and let $P_1, \ldots , P_r$ be the $r$ distinct maximal components contained in $C \cap E_\ell$.
    It is possible that there exists some $i \in [r]$ such that the $m$-tight cycle in $P_i \cup C_{\ell+1}$ is a dicycle, but in this case $P_i$ must be part of a linkage $\{ Z_1 , \ldots , Z_\ell \}$ consisting of $C_{\ell+1}$-paths with $Z_j \subseteq E_j$ for $j \in [\ell]$ and $Z_\ell = P_i$, such that for each $j \in [\ell]$ the $m$-tight cycle of $Z_j \cup C_{\ell+1}$ is a dicycle and for each $j \in [2,\ell-1]$, the trace of $Z_j$ separates the traces of $Z_{j-1}$ and $Z_{j+1}$ in the $m$-disk of $C_{\ell+1}$.
    Let $h \in [\ell + 1]$ be minimal such that $Z_{h-1} \cap C_h \neq \emptyset$.
    This allows us to use \Cref{lem:detour} route $\{ C_h, \ldots ,  C_{\ell+1} \}$ via $\{ Z_{h-1}, \ldots , Z_\ell \}$, which results in the homogeneous dicycles $\{ C_h', \ldots , C_{\ell+1}' \}$, all of which are found in $H$.
    This means we could change the initial diamond rim $\mathfrak{R}$ such that the new homogeneous family of dicycles is $\{ C_1, \ldots , C_{h-1}, C_h', \ldots , C_{\ell+1}', C_{\ell+2}, \ldots , C_\SC \}$.
    Since this is always possible, provided that the $m$-tight cycle in $P_i \cup C_{\ell+1}$ is a dicycle, we will assume that no such $P_i$ exists.
    We make the same assumption for all pairs of $E_i$ and $C_j$ with $i,j \in [\SC]$.

    Thus for each $i \in [r]$ the $m$-tight cycle $P_i \cup C_{\ell+1}$ is a diamond.
    For each $i \in [r]$, we define a circle $\ell$-rim $\mathfrak{W}_i = ( H_i, \mathcal{C}_i = \{ C_1^i, \ldots , C_\ell^i \} ) )$ as follows.
    Let $C_\ell^i$ be the dicycle in $C_{\ell+1} \cup P_i$ that is distinct from $C_{\ell+1}$ and therefore contains $P_i$ in its entirety.
    Furthermore, let $c_i$ be the disk defined by the trace of $C_\ell^i$ that does not contain $m$.
    Due to $\delta$ being an odd decomposition, we can construct the rest of $\mathcal{C}_i$ by demanding that for $j \in [\ell-1]$ the dicycle $C_j^i$ is the subgraph in $E_j \cup C_{\SC-j}$ such that the trace of $C_j^i$ separates $m$ and $c_i$.
    Note that $C_1^i$ is found in $E_1 \cup C_{\SC-1}$ and $C_\SC$ is not used in the construction at all.
    By definition $\mathcal{C}_i$ is therefore a homogeneous set of pairwise disjoint dicycles.
    $H_i$ is defined by restricting $H$ to the disk defined by the $m$-tight cycle of $C_1^i$ that does not contain $m$.
    For each $\mathfrak{W}_i$, the dicycle $C_\ell^i$ is what defines the $\mathfrak{W}_i$ society.
    We can find both a $C_\ell^i$-$C_1^i$- and a $C_1^i$-$C_\ell^i$-linkage, each of order $\ell$, by considering the dicycles $\{ C_1, \ldots , C_\ell \} $ which intersect all dicycles in $\mathcal{C}_i$.

    Analogously, we define $C'$ to be the $m$-tight cycle in $C_{\ell+1} \cup E_{\ell+1}$, with $P_1', \ldots , P_s'$ being the components of $C' \cap E_{\ell+1}$, and for each $i \in [s]$, let $\mathfrak{X}_i = ( H_i' , \mathcal{F}_i = \{ F_1^i, \ldots , F_\ell^i \} )$ be a circle $\ell$-rim, such that for $j \in [\ell]$ the dicycle $F_j^i$ is a subgraph in $E_{j+1} \cup C_{\SC-j}$.
    Thus $F_1^i$ is found in $E_2 \cup C_{\ell+1}$ and $F_\ell^i$ is found in $E_{\ell+1} \cup C_{\SC-1}$, with $C_\SC$ again going unused.
    We call the circle rims $\mathfrak{W}_i$ and $\mathfrak{X}_j$ for $i \in [r]$ and $j \in [s]$, together with their segregated pairs of linkages, the \emph{type II rims}.

    Let $A$ be the $m$-tight cycle in $C_\SC \cup E_{\ell+1}$, note that $P_i \subseteq A$ for all $i \in [r]$, and that the trace of $A$, $C$, and $C'$ all separate $m$ from the vertices on the $\mathfrak{R}$-society, respectively for $A$ those vertices of the $\mathfrak{R}$-society that do not lie on $C_\SC$.
    Thus there exists a pair of $A$-paths $Q_1$ and $Q_2$ forming a cross on $A$, with $Q_i \subseteq L_i$ for $i \in [2]$.
    Since the trace of $C$ also separates $m$ and the trace of $C_\SC$, there also exists a pair of $C$-paths $R_1$ and $R_2$ forming a cross on $C$, with $R_i \subseteq Q_i$ for $i \in [2]$.
    
    We will now reroute the endpoints of $R_1$ and $R_2$, and in some cases $Q_1$ and $Q_2$ or other pairs, onto the society of $\mathfrak{U}_i$, $\mathfrak{W}_j$, or $\mathfrak{X}_h$ for some $i \in [\ell,\SC]$, $j \in [r]$, or $h \in [s]$.
    In general, we will say that we \emph{reroute a path $\alpha$ onto some object $\beta$} or say that we \emph{move the endpoints of a path $\alpha$ onto some object $\beta$} to state that we find a path $\gamma$ with its tail, respectively its head, on $\alpha$ and its head, respectively its tail, on $\beta$, allowing us to replace $\alpha$ with a combination of $\alpha$ and $\gamma$.
    
    Clearly, if all four endpoints of $R_1$ and $R_2$ land either on $C_{\ell+1}$, we have found a cross for a type I rim.
    Alternatively, if all four endpoints land on $P_i$, for some $i \in [r]$, we have found a cross for a type II rim.
    In fact, if all four endpoints land on mutually distinct $P_i$, we note that each $P_i$ is a $C_{\ell+1}$-path by definition and we can therefore reroute each endpoint to lie on $C_{\ell+1}$ via $P_i$.
    Of course, the same holds true if three or less endpoints land on mutually distinct $P_i$ and the rest lie on $C_{\ell+1}$.
    If three endpoints are placed on the same $P_i$, then we can reroute the remaining endpoint first to $C_{\ell+1}$, if it does not lie there already, and then onto $C_\ell^i$, providing us with a cross for a type II rim.
    
    This means we can suppose that exactly two endpoints of $R_1$ and $R_2$ land on the same $P_i$.
    In particular, we can also assume that neither of these endpoints also lies on $C_{\ell+1}$.
    Thus neither of them is also an endpoint of $P_i$.
    Note that these endpoints cannot both be endpoints of $R_1$, or respectively of $R_2$, since the two paths form a cross.
    W.l.o.g.\ we can assume that these endpoints lie on $P_1$.
    For both $i \in [2]$, this implies that $Q_i$ starts or ends on $R_i$, since $P_1 \subseteq A \cap C \subseteq E_{\ell}$.
    Clearly, if exactly one of the other endpoints lands on another $P_i$, we can reroute this end of the path to land on $C_{\ell+1}$, and we can proceed similarly if the other two endpoints find themselves on distinct $P_i$.
    Thus we can further narrow our focus to the two cases in which either both of the other endpoints of $R_1$ and $R_2$ land on $C_{\ell+1}$, or alternatively, they land on some $P_i$ with $i \in [2, r]$.

    First, let us suppose that $P_1$ contains a head and a tail of $R_1$ and $R_2$.
    If the tail occurs before the head on $P_1$, then we can clearly use $P_1$ to route both of them onto $C_{\ell+1}$, giving us a cross on a type I rim.
    In fact, if the other two endpoints land on some $P_i$ with $i \in [2,r]$, then on $P_i$ the head must occur before the tail.
    This leads us to the corresponding case for $P_1$.
    We let $h_i$ be the head and $t_i$ be the tail of $R_i$, for both $i \in [2]$.
    W.l.o.g.\ we can assume that $h_1$ occurs before $t_2$ on $P_1$.
    
    Let $h_i'$ be the head of $Q_i$ and let $t_i'$ be the tail of $Q_i$, for both $i \in [2]$.
    According to our previous observations, we know that $h_1' = h_1$ and $t_2' = t_2$.
    Let $J_1$ be the maximal subpath of $E_\ell$ within $A$ that contains $P_1$ and note that if $h_1'$ or $t_1'$ lie on $A \setminus (C_\SC \cup J_1)$, we can route these endpoints onto $C_\SC$ according to our previous arguments.
    Furthermore, since $Q_1$ and $Q_2$ form a cross, $t_1'$ and $h_2'$ are found on $A - h_1P_1t_2$ in the given order, when traversing this path from $t_2$ to $h_1$.
    
    Suppose that $t_1'$ is an internal vertex of $J_1$.
    Then it must be found on the $t_2$-$C_\SC$-path in $J_1$.
    In particular, $t_2 J_1 t_1'$ contains at least one vertex $v \in V(P_1)$ that is distinct from $t_2$, since we know that $t_2$ is not an endpoint of $P_1$.
    We can therefore use $vP_1t_2''Q_1h_1$ as one path which has both endpoints on $C_\ell^1$.
    And $h_2'$ can then be first rerouted to the tail of $J_1$, if it does not lie on $J_1$ to begin with, and then we can also extend $Q_2$ to have both endpoints on $C_\ell^1$.
    This can even be done if $h_2'$ is found on $t_2 J_1$, since it must then also lie on $t_1' J_1$.
    An analogous argument resolves the case in which $h_2'$ is found on $J_1$.
    Thus we can further suppose that $t_1'$ and $h_2'$ lie on $C_\SC$.

    Note that $J_1 \cup C_\SC$ consists of two dicycles that intersect in a single segment, one of them being $C_\SC$.
    Let $P$ be the directed $J_1$-path in $C_\SC$ that forms a dicycle $P \cup J_1$ that separates $m$ and the trace of $C_\SC - (P \cup J_1)$ and note that $t_1'$ appears after $h_2'$ on $P$.
    Therefore $Q_1 \cup Q_2 \cup h_1P_1t_2 \cup h_2'Pt_1'$ forms a dicycle which intersect both of the dicycles in $J_1 \cup C_\SC$ in a single segment, which is disjoint from their intersection.
    This is a weak odd bicycle and we are therefore done.
    (We note here that this case is in fact the one that forces us to live with being third-integral.)  

    We now move on to our last two cases, which happen to be analogous to one another.
    Suppose that $h_1$ and $h_2$ are the two endpoints of $R_1$ and $R_2$ that are found on $P_1$ in the given order.
    Our intermediate goal will be to first reduce to the case in which $t_1$ and $t_2$ both lie on $C_{\ell+1}$.
    Of course if both $t_1$ and $t_2$ lie on distinct $P_i$, we can simply reroute them onto $C_{\ell+1}$ and we can reduce similarly if only one of them is found on some $P_i$ with $i \in [2,r]$.
    Thus we can assume that $t_1$ and $t_2$ are both found on some component of $C - (C_{\ell+1} \cup P_1)$.

    Consider a cross on $C' \subseteq C_{\ell+1} \cup E_{\ell+1}$ consisting of a pair of paths $S_1$ and $S_2$ with $S_i \subseteq L_i$ for $i \in [2]$.
    Analogous to our previous arguments, we can assume that $S_1$ and $S_2$ have both of their heads, or both of their tails, on some $P_p'$ for $p \in [s]$.
    If the other two endpoints of $S_1$ and $S_2$ can be routed onto $C_{\ell+1}$, we consider this as having reached our intermediate goal and adress this later on.
    Thus, we can assume that both tails of $S_1$ and $S_2$ lie on $P_p'$ for some $p \in [s]$.
    
    Since the trace of $C$ separates $m$ and the trace of $C' - C_{\ell+1}$ there exists at least one $C$-cross $T_1,T_2$ with $T_i \subseteq S_i$ for $i \in [2]$.
    According to our previous arguments, any such $C$-cross must have at least two endpoints on some $P_i$ with $i \in [r]$ and in particular, both of these endpoints must either both be the tails, or respectively heads, of $T_1$ and $T_2$.
    Furthermore, we can again assume that each $C$-cross found in $S_1$ and $S_2$ has both heads on some $P_q$ and both tails on some $P_{q'}$ with $p,q \in [r]$ and $p \neq q$, as we have otherwise reached our intermediate goal and can proceed from there.

    Let $a_i$ be the head and $b_i$ be the tail of $S_i$ for $i \in [2]$.
    We can choose a $C$-cross $T_1$ and $T_2$ within $S_1$ and $S_2$, where $a_i'$ is the head of $T_i$ for $i \in [2]$, such that $b_1S_1a_1'$ and $b_2S_2a_2'$ contain only a single pair of paths forming a cross on $C$, namely $T_1$ and $T_2$.
    Suppose that $T_1$ and $T_2$ have their tails $b_1'$ and $b_2'$ on $P_i$ and their heads on $P_1$ for some $i \in [2,r]$.
    Since $b_1S_1b_1'$ and $b_2S_2b_2'$ do not contain a cross on $C$ according to our assumptions, and as $\delta$ is an odd decomposition, which forces the parts of these paths drawn between $C$ and $C'$ to be part of a strong decomposition, we can use \Cref{lem:detour} to reroute these two paths onto $C_{\ell+1}$ using $P_i$ and $P_p'$.
    
    Thus we have reached our intermediate goal and the pair $T_1$ and $T_2$ can be assumed w.l.o.g.\ to have their heads on $P_1$ and their tails on $C_{\ell+1}$.
    Consider the $m$-tight cycle $F$ in $P_1 \cup C_\ell$.
    Clearly, there exists a cross on $F$ within $T_1$ and $T_2$ consisting of two paths $T_1'$ and $T_2'$ with $T_i' \subseteq T_i$ for $i \in [2]$ and if three or more endpoints of this pair lie on $F \cap C_\ell$ we can find a cross on a type I rim.
    Thus the heads of $T_1'$ and $T_2'$ also lie on $P_1$ and their tails lie on $F \cap C_\ell$.
    Let $P$ be $P_1$-path in $C_{\ell+1}$ that forms a cycle together with $P_1$ whose $m$-disk separates $m$ and the trace of $C_\SC$ and let $P' \subseteq C_\ell$ be defined analogously.
    Note that both tails of $T_1$ and $T_2$ lie on $P$, allowing us to again apply \Cref{lem:detour}, this time to $T_1$ and $T_2$, and $P$ and $P'$ to reroute the tails of $T_1$ and $T_2$ onto $P_1$ guaranteeing us another cross on a type II rim.
    (Reaching the intermediate goal with $S_1$ and $S_2$ can be resolved analogously by considering the $m$-tight cycle in $P_1' \cup C_\ell$, where we assume that both heads, respectively tails, of $S_1$ and $S_2$ lie on $P_1'$.
    The result will then be a cross on $\mathfrak{X}_1$, which is also a type II rim.)
    This completes our proof.
\end{proof}

By splitting up a large diamond rim into several layers with disjoint infrastructure, we can then find a large third-integral packing of weak odd bicycles analogous to \Cref{cor:circleapplyshifting}.
\begin{lemma}[Diamond non-planar shifting]\label{lem:diamondnonplanarshifting}
    Let $\theta, t$ be integers with $\theta$ being a multiple of $t\SC$, let $D$ be a digraph, let $\delta = (\Gamma, \mathcal{V}, \mathcal{D})$ be an odd decomposition of $D$, let $m \in C(\delta)$ be a maelstrom of $\delta$, and let $\mathfrak{R} = (H, \mathcal{C}, \mathcal{E} )$ be a diamond $\theta$-rim of $m$, with a transaction $\mathcal{L} = \{ L_1, \ldots , L_{2t} \}$ on the $\mathfrak{R}$-society $(H_\mathfrak{R}, \Omega_\mathfrak{R})$, such that for all $i \in [t]$ the paths $L_{2i-1}, L_{2i}$ are a crossing pair for $(H_\mathfrak{R}, \Omega_\mathfrak{R})$.
    
    Then the digraph $H_\mathfrak{R}$ contains a third-integral packing of $t$ weak odd bicycles.
\end{lemma}
\begin{proof}
    By splitting up $\mathfrak{R}$ into $t$ layers $\mathfrak{R}_1, \ldots , \mathfrak{R}_t$, with $\mathfrak{R}_i$ being an $\SC$ layer at $(i-1)\SC +1$ for $i \in [t]$, we can find $t$ disjoint layers.
    Applying \Cref{lem:diamondlayershifting} on each layer together with one of our $t$ pairs of paths forming crosses, we can find $t$ weak odd bicycles.
    The third-integrality of the resulting packing is then derived analogously to \Cref{cor:circleapplyshifting}, since we have three sets $\mathcal{C}$, $\mathcal{E}$, and $\mathcal{L}$, each individually containing mutually disjoint objects.
\end{proof}


The proof for refined diamond outlines proceeds in a way that is reminiscent of the proof for diamond outlines, but the technical details ultimately yield quarter-integral.
Thus this represents the first part of our overall approach where we have to resign ourselves to quarter- instead of third- or even half-integraltiy.
To begin in this direction, we define what constitutes a refined diamond rim.
\begin{definition}[Refined diamond rim]\label{def:refineddiamondrim}
Let $\theta, k, g, d$ be positive integers, let $D$ be a digraph, let $\delta = (\Gamma, \mathcal{V}, \mathcal{D})$ be an odd decomposition of $D$, and let $m \in C(\delta)$ be a maelstrom of $\delta$.
A \emph{refined diamond $\theta$-rim} of $n \subseteq m$ with \emph{eddy degree} $d$ and \emph{roughness} $g$ is a tuple $\mathfrak{R} = (H, \mathcal{C}, \mathcal{E}, \mathfrak{B}, \mathfrak{D})$, where $\mathcal{C} = \{ C_1, \ldots , C_\theta \}$ is a homogeneous family of pairwise disjoint dicycles, $\mathcal{E} = \{ E_1, \ldots , E_\theta \}$ is a $C_\theta$-linkage, $\mathfrak{B} = \{ \mathcal{B}_1, \ldots , \mathcal{B}_d \}$ is a set of linkages with $\mathcal{B}_i = \{ B_1^i, \ldots , B_{\theta (g+1)}^i \}$, for all $i \in [d]$, $\mathfrak{D}$ is a collection of $d$ closed disks, each disjoint from $n$, and the following requirements are met.

\begin{description}
    \item[RDR1] $\mathfrak{R}' = (H', \mathcal{C}, \mathcal{E})$ is a diamond $\theta$-rim of $m$, where $H'$ is the restriction of $H$ to $\mathcal{C}$ and $\mathcal{E}$.

    \item[RDR2] The points \textbf{\textsf{E1}}, \textbf{\textsf{E2}}, \textbf{\textsf{E4}}, \textbf{\textsf{E7}}, \textbf{\textsf{E8}}, \textbf{\textsf{E9}}, and \textbf{\textsf{E10}} of \Cref{def:refinedoutline} hold for $\mathfrak{R}$, provided we let $\mathfrak{V}^* = \{ n \}$, $\mathcal{E}_1 = \mathcal{E}$, as well as $E_i^1 = E_i$ for all $i \in [\theta]$, and we let the object $\mathfrak{M}$ be $\mathfrak{R}'$.
    Following \textbf{\textsf{E4}}, we let $Q_{d+1}$ be the $n$-tight cycle of $C_\theta \cup E_\theta \cup \bigcup_{i=1}^d B_{\theta (g+1)}^i$ and define $\Delta_{i+1}$ for all $i \in [d]$.
    
    \item[RDR3] For all $i \in [d]$, $\mathcal{B}_i$ is a choppy transaction on the $\mathfrak{R}'$-society, such that for all $j \in [d]$ with $j > i$, the entirety of $\mathcal{B}_j$ is drawn in $\Delta_{i+1}$.
    
    \item[RDR4] The graph $C_1 \cup E_1 \cup \bigcup_{i=1}^d B_1^i$ has an $n$-tight cycle in $\delta$.

\end{description}

We call $(H_\mathfrak{R}, \Omega_\mathfrak{R})$ the \emph{$\mathfrak{R}$-society of $n$}, where $V(\Omega_\mathfrak{R}) = V(Q_{d+1})$ and $H_\mathfrak{R} = H \cup \sigma(n)$.
\end{definition}

We will again work with layers to prove that we can find weak odd bicycles within refined diamond rims that host a cross.
\begin{definition}[Layers of a refined diamond rim]\label{def:refineddiamondlayer}
    Let $\theta, \ell , k$ be positive integers, let $d,g$ be integers, let $D$ be a digraph, let $\delta = (\Gamma, \mathcal{V}, \mathcal{D})$ be an odd decomposition of $D$, let $m \in C(\delta)$ be a maelstrom of $\delta$, and let $\mathfrak{R} = (H, \mathcal{C} = \{ C_1, \ldots , C_\theta \} , \mathcal{E} = \{ E_1, \ldots , E_\theta \}, \mathfrak{B} = \{ \mathcal{B}_1, \ldots , \mathcal{B}_d \} )$ be a refined diamond $\theta$-rim of $n \subseteq m$ of eddy degree $d$ and roughness $g$, with $\mathcal{B}_i = \{ B_1^i, \ldots , B_{\theta (g+1)}^i \}$, for all $i \in [d]$.
    We call the tuple $(H', \{ C_\ell , \ldots C_{\ell + k - 1} \}, \{ E_\ell , \ldots E_{\ell + k - 1} \}, \{ \{ B_\ell^i , \ldots B_{\ell + k - 1}^i \} \mid i \in [d] \} )$ a \emph{$k$-layer at $\ell$ in $\mathfrak{R}$}, where $H'$ is the union of $\bigcup_{i=1}^d \bigcup_{j = \ell}^{\ell +k -1} ( E_j \cup B_j^i ) $ and the subgraph of $H$ drawn on the intersection of the largest $m$-disk bounded by a cycle in $C_{\ell + k - 1} \cup E_{\ell + k - 1} \cup \bigcup_{i= 1}^d B_{\ell + k - 1}^i$ and the largest non-$m$-disk within $C_\ell \cup E_\ell  \cup \bigcup_{i= 1}^d B_\ell^i$.
\end{definition}

We now prove that we can find weak odd bicycles in refined diamond layers with crosses on them by reducing the problem to an application of the methods discussed in the proof of \Cref{lem:diamondlayershifting}.
\begin{lemma}\label{lem:refineddiamondlayershifting}
    Let $\theta$ be an integer with $\theta \geq \SC$, let $D$ be a digraph, let $\delta = (\Gamma, \mathcal{V}, \mathcal{D})$ be an odd decomposition of $D$, let $m \in C(\delta)$ be a maelstrom of $\delta$, and let $\mathfrak{R} = (H, \mathcal{C}, \mathcal{E}, \mathfrak{B} )$ be a refined diamond $\theta$-rim of $m$, with a crossing pair $L_1, L_2 \subseteq H_\mathfrak{R}$ for the $\mathfrak{R}$-society $(H_\mathfrak{R}, \Omega_\mathfrak{R})$.
    Then for every $\SC$-layer $(H', \mathcal{C}', \mathcal{E}', \mathfrak{B}')$ there exists a weak odd bicycle in $H' \cup L_1 \cup L_2$.
\end{lemma}
\begin{proof}
    Let $\ell = \lfloor \SC \rfloor$, with $\ell$ again corresponding to the value we used in \Cref{lem:circlenonplanarshifting}.
    Let $d$ be the eddy degree of $\mathfrak{R}$ and note that since the roughness $g$ of $\mathfrak{R}$ is a non-negative integer, we know that each member of $\mathfrak{B}$ has order at least $\SC$.
    We will not need many paths and thus we can neglect the actual roughness of $\mathfrak{R}$ and assume that $g = 0$.
    As in \Cref{lem:circlenonplanarshifting}, we will want to route each cross onto an $\SC$-rim for which we have already proven a shifting lemma.
    We note that since $\mathcal{C}$ and $\mathcal{E}$ form a diamond $\SC$-rim, the definitions for $\mathfrak{U}_i$, $\mathfrak{W}_j$, and $\mathfrak{X}_h$ given in the proof of \Cref{lem:diamondnonplanarshifting} can still be used for this proof.
    Once more, we will assume that we have chosen the $\SC$ innermost cycles and paths as our layer to ease notation.

    We make some observations about the structure of $\mathfrak{R}$ before we delve deeper into the proof.
    As in the proof of \Cref{lem:diamondnonplanarshifting}, we assume that within the $m$-tight cycle $C$ of $C_i \cup E_j$, for $i,j \in [\SC]$, the $m$-tight cycle in $C \cup P$ formed by each of the components $P$ of $C \cap E_j$ is not a dicycle.
    Similarly, within the $m$-tight cycle $C'$ of $C_{\SC} \cup E_i \cup B_j^h$, for $i,j \in [\SC]$ and $h \in [d]$, the $m$-tight cycle in $C \cup B$, formed by each of the components $B$ of $C' \cap B_j^h$ is not a diamond.
    Analogous to our arguments in the proof of \Cref{lem:diamondnonplanarshifting} we could otherwise reroute parts of $\mathcal{E}$ to more closely encompass $n$.
    However, we must note that this does not actually leave us with a valid refined diamond rim, since the choppy transaction $\mathcal{B}_h$ might separate a disk $d \in \mathcal{D}$ that partily or entirely contains an even dicycle.
    This is however not a problem for our constructions, since they do not interact with the contents of the disk from $\mathcal{D}$, and we only make these modifications to reduce the number of cases we have to consider.
    We additionally make the assumption that, if $C''$ is the $m$-tight of $C_i \cup E_{\SC} \cup B_j^h$, for $i,j \in [\SC]$ and $h \in [d]$, the $m$-tight cycle in $C \cup B$, formed by each of the components $B$ of $C' \cap B_j^h$ is not a dicycle.
    Now we can proceed with the proof, where we will again be working in layers, each allowing us to find one weak odd bicycle.
    
    Let $S$ be the $m$-tight dicycle in $C_{\ell+1} \cup E_{\ell+1} \cup \bigcup_{i=1}^d B_{\ell+1}^i$.
    Clearly $L_1$ and $L_2$ contain a cross on $S$, consisting of the paths $R_1$ and $R_2$, with $R_i \subseteq L_i$ for $i \in [2]$.
    If all four endpoints of $R_1$ and $R_2$ are on $S' = S \cap (C_{\ell+1} \cup E_{\ell+1})$, we can find a cross and thus a weak odd bicycle by following the constructions in the proof of \Cref{lem:circlenonplanarshifting}.
    Should all four endpoints of $R_1$ and $R_2$ lie on $C_{\ell+1}$ and components of $S - (C_{\ell+1} \cup E_{\ell+1})$ that have both endpoints on $C_{\ell+1}$, the same applies.
    (Note that the paths $P_1, \ldots , P_r$ and the associated transactions from the proof of \Cref{lem:circlenonplanarshifting} behave like choppy transactions with respect to $C_{\ell+1}$.)
    
    Suppose at most two endpoints of $R_1$ and $R_2$ lie on $E_{\ell+1}$, such that they are not a head followed by a tail on $E_{\ell+1}$.
    We can directly route them onto $C_{\ell+1}$, if there is at most one, or a tail appears before a head on $E_{\ell+1}$.
    Otherwise we can use \Cref{lem:detour} on $R_1$ and $R_2$, and $E_\ell$ and $E_{\ell+1}$ to either find three or more endpoints of a cross on $E_\ell$, allowing us to reroute all of them onto $E_\ell$, or we can reroute the two endpoints on $E_{\ell+1}$ either onto $C_{\ell+1}$.
    With either case reducing to one of the simple cases discussed prior to this.

    Suppose that $R_1$ and $R_2$ together have two endpoints on $E_{\ell+1}$ such that there is a head followed by a tail on $E_{\ell+1}$.
    If these two endpoints lie on distinct $C_{\ell+1}$-paths within $E_{\ell+1}$, we can proceed as above.
    Otherwise, we note that due to $R_1$ and $R_2$ forming a cross on $S$, if the other two endpoints of $R_1$ and $R_2$ lie on the same $C_{\ell+1}$-path $B \subseteq B_{\ell+1}^i$, for some $i \in [d]$, the tail on $B$ must appear before the head.
    This allows us to route these two endpoints onto $C_{\ell+1}$.
    Of course all other placements of these two endpoints permit us to do the same.
    Thus we can find a weak odd bicycle within $E_{\ell+1} \cup C_{\ell+1} \cup R_1 \cup R_2$ as in the proof of \Cref{lem:circlenonplanarshifting}, where we assume $R_1$ and $R_2$ to have two of their endpoints on $C_{\ell+1}$.

    Another fortunate situation for us occurs when at least three endpoints of $R_1$ and $R_2$ are found on the same component $B \subseteq B_{\ell+1}^i$, for some $i \in [d]$, in $S - (C_{\ell+1} \cup E_{\ell+1})$.
    We may route the potential four endpoints onto the $C_{\ell+1} \cup E_{\ell+1}$-path $B'$ within $B_{\ell+1}$ that contains $B$.
    If the endpoints of $B'$ lie on $C_{\ell+1}$ it is easy to see that we can construct a circle $\ell$-rim akin to $\mathfrak{W}_i$.
    The same kind of construction also works if $B'$ has both endpoints on $E_{\ell+1}$.
    We note that we can assume that we are in one of these two cases as $\mathcal{B}_i$ is a choppy transaction for $\mathcal{C}$ or $\mathcal{E}$ and thus we could optionally simply follow $B_{\ell+1}^i$ until we have found a $E_{\ell+1}$- or a $C_{\ell+1}$-path containing $B$.
    In either case we are therefore again done.

    Next, let us suppose that three or more endpoints of $R_1$ and $R_2$ are found on $E_{\ell+1}$ and components of $S - (C_{\ell+1} \cup E_{\ell+1})$ that have both endpoints on $E_{\ell+1}$.
    By again invoking the constructions in the proof of \Cref{lem:circlenonplanarshifting}, it is clear that if we can route all endpoints onto $E_{\ell+1}$ we are done.
    Thus at least two endpoints of $R_1$ and $R_2$ lie on some $E_{\ell+1}$-path $P \subseteq B_{\ell+1}^i$, for some $i \in [d]$.
    Using $B_\ell^i$, \Cref{lem:detour}, and other rerouting strategies we have presented prior to this, we can narrow our focus on the case in which $P$ contains a tail appearing before a tail and the other two endpoints lie on $E_{\ell+1}$.
    
    Let $h_i$ be the head and $t_i$ be the tail of $R_i$ for $i \in [2]$ and assume w.l.o.g.\ that $t_1$ appears before $h_2$ on $B$.
    Note that there exists a dicycle $C \subseteq B' \cup E_{\ell+1}$ with $B' \subseteq C$ and $E_{\ell+1} - C$ contains the other two endpoints of $R_1$ and $R_2$.
    If $b_1$ is the head and $b_2$ is the tail of $B'$, then this lets us conclude that $t_2, h_1, b_1, b_2$, or $h_1, b_1, b_2, t_2$, or $b_1, b_2, t_2, h_1$ is the order of the given vertices on $E_{\ell+1}$.
    The second order allows us to immediately reroute all endpoints onto $B'$.
    For the first and third order, we can route one endpoint onto $B'$ on $E_{\ell+1}$ and for the other, we can follow $E_{\ell+1}$ until we meet $C_{\ell+1}$, which we can then follow until we reach the $C_{\ell+1}$-path in $E_{\ell+1}$ that contains $b_1$, which in turn allows us to route this endpoint onto $B'$ as well.

    With this we have narrowed our view to the case in which two endpoints of $R_1$ and $R_2$ lie on $C_{\ell+1}$ and components of $S - (C_{\ell+1} \cup E_{\ell+1})$ that have both endpoints on $C_{\ell+1}$ and the other two endpoints lie on components of $S - (C_{\ell+1} \cup E_{\ell+1})$ that have both endpoints on $E_{\ell+1}$.
    Of course if the two endpoints can be rerouted onto $E_{\ell+1}$ we are again done and by using $B_\ell^i$, $E_\ell^i$, \Cref{lem:detour}, and our other rerouting strategies, we can further narrow this case down.
    Thus we can assume that $t_1$ appears before $h_2$ on a path $Q \subseteq B_{\ell+1}^i$, for some $i \in [d]$, that has both of its endpoints on $E_{\ell+1}$.
    Let $Q'$ be the $E_{\ell+1}$-path $Q$ within $B_{\ell+1}$ that contains $Q$.

    Suppose first that the endpoints of $Q'$ both lie on a common $C_{\ell+1}$-path $F$ within $E_{\ell+1}$.
    We note that $F \cup E_{\ell+1}$ contains a dicycle that intersects the unique dicycle in $E_{\ell+1} \cup C_{\ell+1}$, that is distinct from $C_{\ell+1}$, in a single subpath.
    If the endpoints $t_2$ and $h_1$ lie on $C_{\ell+1}$ or can be rerouted onto $C_{\ell+1}$ via some combination of our rerouting methods, we note that $C_{\ell+1} \cup R_1 \cup R_2 \cup t_1Qh_2$ contains a dicycle distinct from $C_{\ell+1}$ that intersects both of the previously mentioned dicycles in single subpaths, yielding a weak odd bicycle.
    Otherwise $t_2$ and $h_1$ must lie on some $C_{\ell+1}$-subpath $S \subseteq B_{\ell+1}^j$, with $j \in [d]$ and $j \neq i$.
    This means that $R_1 \cup R_2 \cup t_2Sh_1 \cup t_1Qh_2$ contains a dicycle intersecting the dicycle in $F \cup E_{\ell+1}$ and a dicycle in $S \cup C_{\ell+1}$, which together with $C_{\ell+1}$ and the other dicycle in $C_{\ell+1} \cup E_{\ell+1}$ forms a weak odd bicycle.

    If the endpoints of $Q'$ lie on two distinct $C_{\ell+1}$-paths $F_1$ and $F_2$ within $E_{\ell+1}$, we first require some observations about $F_1$ and $F_2$, where we let $F_1$ contain the head of $Q'$.
    Let $u_i$ be the head and $v_i$ be the tail of $F_i$ for $i \in [2]$.
    Our observations about the structure we can assume to find within the rim, we know that on $C_{\ell+1}$ the vertices $u_1,v_1,u_2,v_2$ lie in the given order when traversing $C_{\ell+1}$ according to its orientation.
    Furthermore, if we let $T$ be the $v_1$-$u_2$-path within $C_{\ell+1}$ then $T \cup Q' \cup F_1 \cup F_2$ contains a unique dicycle $J$.
    
    Should we be able to route $t_2$ and $h_1$ onto $C_{\ell+1}$, we do so and note that $t_2$ must occur before $h_1$ on the path $C_{\ell+1} - T$.
    This allows us to route $h_1$ onto $F_1$ and from there onto $Q'$ and we can route $t_2$ onto $Q'$ via $F_2$ in an analogous fashion.
    Thus all we have left to consider is the case in which $h_1$ and $t_2$ occupy the same $C_{\ell+1}$-path $W \subseteq B_{\ell+1}^j$, with $j \in [d]$ and $j \neq i$ and $h_1$ occurs before $t_2$ on $W$.
    In this situation we can find a dicycle within $Q' \cup F_1 \cup F_2 \cup C_{\ell+1}$ that contains $Q'$ and intersects $J$ in a single path.
    A third dicycle, intersecting both the dicycle we have just found and $J$ in a single path, is constituted by $R_1 \cup h_1Wt_2 \cup R_2 \cup h_2Qt_1$.
    This yields another weak odd bicycle and completes the proof.
\end{proof}

This leads us to the conclusion of the section and an analogous result to \Cref{lem:diamondnonplanarshifting}.
\begin{lemma}[Refined diamond non-planar shifting]\label{lem:refineddiamondnonplanarshifting}
    Let $\theta, t$ be integers with $\theta$ being a multiple of $t\SC$, let $D$ be a digraph, let $\delta = (\Gamma, \mathcal{V}, \mathcal{D})$ be an odd decomposition of $D$, let $m \in C(\delta)$ be a maelstrom of $\delta$, and let $\mathfrak{R} = (H, \mathcal{C}, \mathcal{E}, \mathfrak{B}, \mathfrak{D} )$ be a refined diamond $\theta$-rim of $m$, with a transaction $\mathcal{L} = \{ L_1, \ldots , L_{2t} \}$ on the $\mathfrak{R}$-society $(H_\mathfrak{R}, \Omega_\mathfrak{R})$, such that  for all $i \in [t]$ the paths $L_{2i-1}, L_{2i}$ are a crossing pair for $(H_\mathfrak{R}, \Omega_\mathfrak{R})$.
    
    Then the digraph $H_\mathfrak{R}$ contains a quarter-integral packing of $t$ weak odd bicycles.
\end{lemma}
\begin{proof}
    As in the proof of \Cref{lem:refineddiamondnonplanarshifting}, we simply split our rim up into $t$ layers and apply \Cref{lem:refineddiamondlayershifting}.
    The result is a quarter-integral packing of $t$ weak odd bicycles since, for each set $\mathcal{C}$, $\mathcal{E}$, $\mathcal{L}$, and $\bigcup_{i \in [d]} \mathcal{B}_i$, the objects contained within are pairwise disjoint, analogous to our arguments in \Cref{lem:diamondnonplanarshifting} and \Cref{cor:circleapplyshifting}.
\end{proof}

\section{Finding the rim of an outline}\label{sec:localising}

In \cref{sec:shifting} we have introduced a wide array of tools that help us to shift even dicycles given a large rim.
We dedicate this section to proving that these rims can easily be found within the outlines we have defined earlier.
This then allows us to establish higher level versions of the findings from \cref{sec:shifting} that are tailored to be used within the context of maelstrom outlines.


We start with the easiest type of rim to find.

\begin{lemma}\label{lem:circleoutlinererouting}
    Let $\theta$ be an integer, let $D$ be a digraph, let $\delta=(\Gamma,\mathcal{V},\mathcal{D})$ be an odd decomposition of $D$, and let $m \in C(\delta)$ be a maelstrom of $\delta$, with a circle $\theta$-outline $\mathfrak{M}=(H,\mathcal{C},\mathfrak{E},\mathfrak{V})$ of positive degree.

    There exists a circle $\nicefrac{\theta}{2}$-rim $\mathfrak{R} = ( H', \mathcal{C}')$ of $m$, such that the $\mathfrak{R}$-society of $m$ is the $\mathfrak{M}$-society of $m$ and there exists an segregated $\mathcal{C}'$-pair of order $\nicefrac{\theta}{2}$ for $\mathcal{C}'$.
\end{lemma}
\begin{proof}
Let $t$ be the degree of $\mathfrak{M}$.
To start, suppose that $t = 1$.

Note that in this case, clearly any selection of $\nicefrac{\theta}{2}$ paths from $\mathcal{E}_1 \in \mathfrak{E}$ can be used to find a segregated $\mathcal{C}$-pair of order $\nicefrac{\theta}{2}$.
Let us reserve the paths $E_1^1, \ldots , E_{\nicefrac{\theta}{2}}^1$ for this purpose.
Consider the point \textbf{\textsf{M3}} of \Cref{def:outline} and let $O_1$ be the dicycle in $C_{\nicefrac{\theta}{2}} \cup E_{\nicefrac{\theta}{2}}^1$ whose $m$-disk is $m$ itself.
Note that the existence of this cycle, combined with the point \textbf{\textsf{M6}} of \Cref{def:outline} tells us that there exist $\nicefrac{\theta}{2}$ dicycles $C_1', \ldots , C_{\nicefrac{\theta}{2}}'$ with $C_i' \subseteq C_{\nicefrac{\theta}{2} + i} \cup E_{\nicefrac{\theta}{2} + i}^1$ for $i \in [\nicefrac{\theta}{2}]$, such that $\{ C_1', \ldots , C_{\nicefrac{\theta}{2}}', O_1 \}$ is a homogeneous set of concentric dicycles.
In particular, $\bigcup_{i \in [\nicefrac{\theta}{2}]} C_i'$ is drawn entirely in the non-$m$-disk of $O_1$.

Again using \textbf{\textsf{M6}} of \Cref{def:outline}, we can now observe that the paths $E_1^1, \ldots , E_{\nicefrac{\theta}{2}}^1$, which we initially set aside, contain a segregated pair of order $\nicefrac{\theta}{2}$ for $\{ C_1', \ldots , C_{\nicefrac{\theta}{2}}' \}$.
Thus we have found our $\nicefrac{\theta}{2}$-rim $\mathfrak{R}$ together with an appropriate segregated pair.

For $t \geq 2$, we may iteratively continue this process.
To see that this works, note that the linkage $\mathcal{E}_2$ would be a transaction on the $\mathfrak{R}$-society, since this is also the $\mathfrak{M}$-society if $t = 1$.
Or more accurately, the cycle $C_{\nicefrac{\theta}{2}}'$ can be chosen to be $Q_1$ from point \textbf{\textsf{M4}} of \Cref{def:outline}, which is used to define the society on which $\mathcal{E}_2$ is a transaction in point \textbf{\textsf{M5}}.
We can therefore proceed in an analogous fashion to the construction for $t=1$ and find our segregated pair in $\mathcal{E}_t$.
\end{proof}

Of course a similar conclusion may be drawn for circle $\theta$-outlines of degree zero.
We separate this statement out, since in this case we require an segregated $\mathcal{C}$-pair to already exist.

\begin{lemma}\label{lem:circleoutlinereroutingdegreezero}
    Let $\theta$ be an even integer, let $D$ be a digraph, let $\delta=(\Gamma,\mathcal{V},\mathcal{D})$ be an odd decomposition of $D$, and let $m \in C(\delta)$ be a maelstrom of $\delta$, with a circle $\theta$-outline $\mathfrak{M}=(H,\mathcal{C},\mathfrak{E},\mathfrak{V})$ of degree 0 and a segregated $\mathcal{C}$-pair $(\mathcal{L}_1, \mathcal{L}_2)$ of order $\nicefrac{\theta}{2}$.

    There exists a circle $\nicefrac{\theta}{2}$-rim $\mathfrak{R} = ( H', \mathcal{C}')$ of $m$, such that the $\mathfrak{R}$-society of $m$ is the $\mathfrak{M}$-society of $m$ and there exists an segregated $\mathcal{C}'$-pair of order $\nicefrac{\theta}{2}$ for $\mathcal{C}'$.
\end{lemma}

As noted in \Cref{sec:shifting}, we will only really make use of refined diamond rims.
But similar to our process in that section, we first show that we can find a diamond rim in a diamond outline and then show that we can find a refined diamond rim in a refined diamond outline.

\begin{lemma}\label{lem:diamondoutlinererouting}
    Let $\theta$ be an integer, let $D$ be a digraph, let $\delta = (\Gamma, \mathcal{V}, \mathcal{D})$ be an odd decomposition of $D$, and let $m \in C(\delta)$ be a maelstrom of $\delta$, with a diamond $\theta$-outline $\mathfrak{M} = (H, \mathcal{C}, \mathfrak{E}, \mathfrak{V})$.
    
    There exists a diamond $\nicefrac{\theta}{2}$-rim $\mathfrak{R}$ of $m$, such that the $\mathfrak{R}$-society of $m$ is the $\mathfrak{N}$-society of $m$.
\end{lemma}
\begin{proof}
Let $t$ be the degree of $\mathfrak{M}$.
Note that since $\mathfrak{M}$ is a diamond maelstrom, $t$ must be positive.
Let us first consider the case in which $t = 1$.
Here confirming that \textbf{\textsf{DR1}} and \textbf{\textsf{DR2}} of \Cref{def:diamondrim} hold for $(H', \{ C_{\nicefrac{\theta}{2} +1 }, \ldots, C_\theta \}, \{ E_{\nicefrac{\theta}{2} +1 }^1, \ldots, E_\theta^1 \} )$ is simple.

Thus we can assume that $t \geq 2$ and there exists an $\mathcal{E}_2 \in \mathfrak{E}$.
We note that it is possible for all transactions $\mathcal{E}_i \in \mathfrak{E}$ with $i < t$ to form circle outlines, but in this case, we can simply first employ the construction from the proof of \Cref{lem:circleoutlinererouting} to find an appropriate set of homogeneous dicycles to start from.
Therefore, we additionally suppose that $Q_1$ (from point \textbf{\textsf{M4}} of \Cref{def:outline}) forms a diamond.

Unsurprisingly, we consider $\{ E_{\nicefrac{\theta}{2} +1 }^2, \ldots, E_\theta^2 \}$ and want to modify $\{ E_{\nicefrac{\theta}{2} + 1 }^1, \ldots, E_\theta^1 \}$ or $\{ C_{\nicefrac{\theta}{2} + 1 }, \ldots, C_\theta \}$ such that we can otherwise forget the contents of $\mathcal{E}_2$.
If the endpoints of $E_\theta^2$ both lie on $E_\theta^1 - C_\theta$, we note that according to the point \textbf{\textsf{M8}} of \Cref{def:outline} the $m$-tight cycle $Q_2$ in $C_\theta \cup E_\theta^1 \cup E_\theta^2$ must either be a diamond or a dicycle.
Should $Q_2$ be a dicycle, we may proceed as in the proof of \Cref{lem:circleoutlinererouting}.
Otherwise the $m$-tight cycle in $C_\theta \cup E_\theta^1 \cup E_\theta^2$ is a diamond.
Using \Cref{lem:detour} it is now easy to confirm that we can find a linkage $\{ E_1', \ldots , E_{\nicefrac{\theta}{2}}' \}$ such that $E_i' \subseteq \bigcup_{j \in [2]} E_{i+\nicefrac{\theta}{2}}^j$ for all $i \in [\nicefrac{\theta}{2}]$ and together with the dicycles $\{ C_{\nicefrac{\theta}{2} + 1 }, \ldots, C_\theta \}$, this forms a diamond $\nicefrac{\theta}{2}$-rim.

Should both endpoints of $E_\theta^2$ both lie on $C_\theta$, then $Q_2$ may either form a dicycle, in which case we can proceed with the construction from the proof of \Cref{lem:circleoutlinererouting}, or $Q_2$ forms a diamond, in which case we may proceed in an analogous fashion to the previous case.
This leaves us to consider the case in which one endpoint lies on $C_\theta$ and the other on $E_\theta^1 - C_\theta$.
However, here too we can use point \textbf{\textsf{M8}} of \Cref{def:outline} and \Cref{lem:detour} to modify either the transaction or the dicycles such that we find a $\nicefrac{\theta}{2}$-rim which has a society corresponding to $Q_2$.
We iteratively repeat this process as in the proof of \Cref{lem:circleoutlinererouting} to find the desired diamond $\nicefrac{\theta}{2}$-rim.
\end{proof}

\begin{lemma}\label{lem:refineddiamondoutlinererouting}
    Let $g, h, \theta$ be integers, let $D$ be a digraph, let $\delta = (\Gamma, \mathcal{V}, \mathcal{D})$ be an odd decomposition of $D$, let $m \in C(\delta)$ be a maelstrom of $\delta$, with a diamond $\theta$-outline $\mathfrak{M}$, and let $\NN = ( H, \mathcal{C}, \mathfrak{E}, \mathfrak{B}, \mathfrak{V}, \mathfrak{D})$ be a refined diamond $\theta$-outline of $\MM$ around $n \subseteq m$ of eddy degree $h$ and roughness $g$, where $(H_\NN, \Omega_\NN)$ is the $\NN$-society of $n$.
    
    There exists a refined diamond $\nicefrac{\theta}{2}$-rim $\mathfrak{R}$ of $n$ with eddy degree $h$ and roughness $g$, such that the $\mathfrak{R}$-society is the $\mathfrak{N}$-society.
\end{lemma}
\begin{proof}
Let $d$ be the surplus degree of $\mathfrak{N}$ and assume that $d$ is positive, as otherwise $\NN$ is actually a diamond outline and we may simply apply \Cref{lem:diamondoutlinererouting}.
We first note that we can already find a diamond $\nicefrac{\theta}{2}$-rim $\mathfrak{R} = (H', \mathcal{C}' , \mathcal{E}' )$ in $\mathfrak{M}$, again thanks to \Cref{lem:diamondoutlinererouting}, and we thus only have to concentrate on finding appropriate choppy transactions for all of our eddy segments.
This means we are already assured to satisfy \textbf{\textsf{RDR1}}.
Let $\mathcal{C}' = \{ C_1', \ldots , C_\theta' \}$ and $\mathcal{E}' = \{ E_1', \ldots , E_\theta' \}$, such that the $\mathfrak{R}$-society corresponds to a diamond $C'$ in $C_\theta' \cup E_\theta'$.

We will construct one choppy transaction for each of the eddy segments $S$, such that the segment itself will be a path of that transaction and thus be found on the society of the eventual refined diamond rim.
The other transactions can be constructed in an analogous manner, proving that the eddy degree of $\mathfrak{R}$ is $h$.
Note that since $S$ is an eddy segment, it must also be a directed $C'$-path, such that the head and tail both lie on either $C_\theta' \cap C'$ or on $E_\theta' \cap C'$ and the head of $S$ appears before the tail of $S$ on either of these paths, when traversing the relevant path according to its orientation.

We can iteratively construct the set of transactions $\mathfrak{B}_S$ in $\mathfrak{B}$ that belong to the segment $S$ as follows.
$\mathfrak{B}_0 = \{ \mathcal{B}_i \in \mathfrak{B} \mid B_{\theta (g+1)}^i \cap S \neq \emptyset \}$ and for all $j \in [k]$, we set $\mathfrak{B}_j = \{ \mathcal{B}_i \in \mathfrak{B}  \mid B_{\theta (g+1)}^i \cap B_{\theta (g+1)}^h \neq \emptyset, \text{ for some } \mathcal{B}_h \in \mathfrak{B}_{j-1} \}$.
We then set $\mathfrak{B}_S = \mathfrak{B}_k$ and w.l.o.g.\ we can assume that $\mathfrak{B}_S = \{ \mathcal{B}_1, \ldots , \mathcal{B}_\ell \}$ for some $\ell \in [k]$.
Note that according to \Cref{def:refinedoutline}, we can partition $\mathfrak{B}$ in this manner by constructing this type of set for all eddy segments.

Let $C''$ be the $n$-tight dicycle in $S \cup C'$, which contains $S$ in its entirety, and let $\Delta$ be the non-$n$-disk of $C''$.
If there are disks $d_1, \ldots , d_i \in \mathfrak{D}$ that are contained within $\Delta$, we let $\Delta'$ be the smallest closed disk in $\Delta$ that contains $d_1, \ldots , d_i$ and is disjoint from $n$.
We add $\Delta'$ to $\mathfrak{D}^*$.

If $\ell = 1$ we are done and if $\ell \geq 2$ we can proceed as in the previous proofs of this subsection by combining $\mathcal{B}_1$ and $\mathcal{B}_2$ into another transaction $\mathcal{B}_S$ of half their size using \Cref{lem:detour} and point \textbf{\textsf{E7}} from \Cref{def:refinedoutline}.
Due to \textbf{\textsf{E5}} and \textbf{\textsf{E6}} from \Cref{def:refinedoutline}, we know that both $\mathcal{B}_1$ and $\mathcal{B}_2$ are choppy, in an extended sense of the definition given in \textbf{\textsf{E6}}, and this allows us to construct $\mathcal{B}_S$ such that it is again choppy using the methods presented in the proofs of \Cref{lem:circleoutlinererouting} and \Cref{lem:diamondoutlinererouting}.
The points \textbf{\textsf{RDR2}} up to \textbf{\textsf{RDR4}} then all follow either directly from \Cref{def:refinedoutline}, or they are assured by the type of construction we have seen in the proof of the previous two lemmas.
We therefore omit the details of this procedure.
\end{proof}

These lemmas now allow us to establish the following notational convention.

\begin{notation}[Rims of outlines]\label{not:rims}
    From now on, we associate a rim $\mathfrak{R}$ with every outline $\mathfrak{M}$ in the natural way implied by \Cref{lem:circleoutlinererouting}, \Cref{lem:circleoutlinereroutingdegreezero}, \Cref{lem:diamondoutlinererouting}, and \Cref{lem:refineddiamondoutlinererouting}.
    To make future statements easier, we may also give the rim as a tuple containing all possible objects needed for the statement in question, e.g.\ $\mathfrak{R} = ( H , \mathcal{C}, \mathcal{E}, \mathfrak{B}, \mathfrak{D} )$, where all objects are non-empty if $\mathfrak{M}$ is a refined diamond outline with positive surplus degree, the last two objects are empty if $\mathfrak{M}$ is a diamond outline, and the last three objects are empty if $\mathfrak{M}$ is a circle outline.
\end{notation}

Now that we have cleaned up the infrastructure of an outline, we will want to show that transactions on the societies of maelstrom outlines cannot have many non-planar parts.

We can now make use of \Cref{cor:circleapplyshifting}, \Cref{lem:diamondnonplanarshifting}, and \Cref{lem:refineddiamondnonplanarshifting}, and respectively \Cref{lem:circleoutlinererouting}, \Cref{lem:circleoutlinereroutingdegreezero}, \Cref{lem:diamondoutlinererouting}, and \Cref{lem:refineddiamondoutlinererouting}, to obtain half-, third-, or quarter-integral packings of even dicycles from an outline that is thick enough and hosts a large transaction with many crossing pairs.

\begin{corollary}\label{lem:applyshifting}
	Let $\theta, \ell$ be integers with $\theta$ being a multiple of $2 \ell \SC$, let $D$ be a digraph, let $\delta = (\Gamma, \mathcal{V}, \mathcal{D})$ be an odd decomposition of $D$, let $m \in C(\delta)$ be a maelstrom of $\delta$, let $n \subseteq m$, let $\mathfrak{M}$ be a (refined) $\theta$-outline $\mathfrak{M}$ of $m$ (or $n$), and let $\mathcal{L} = \{ L_1, L_2, \ldots , L_{2\ell} \}$ be a transaction on the $\MM$-society $(H_\mathfrak{M}, \Omega_\mathfrak{M})$, such that for all $i \in [\ell]$ the paths $L_{2i-1}, L_{2i}$ form a crossing pair on $(H_\mathfrak{M}, \Omega_\mathfrak{M})$.
    \begin{enumerate}
        \item If $\MM$ is a circle outline and there exists an segregated $\mathcal{C}$-pair of order $\nicefrac{\theta}{2}$, then $H_\mathfrak{M}$ contains a half-integral packing of $2\ell$ even dicycles,

        \item if $\MM$ is a diamond outline, then $H_\mathfrak{M}$ contains a third-integral packing of $\ell$ even dicycles, or
        
        \item if $\MM$ is a refined diamond outline, then $H_\mathfrak{M}$ contains a quarter-integral packing of $\ell$ even dicycles.
    \end{enumerate}
\end{corollary}

\section{Calming a transaction}\label{sec:transaction}

Before taking a full dive into a maelstrom, let us discuss why and how any large transaction on the maelstrom society of a (refined) $\theta$-outline, with sufficiently large $\theta$, either yields a large quarter-integral packing of even dicycles, or contains a large transaction that can be integrated into our odd decomposition.

Let $\mathcal{P}$ be a transaction on a society $(D,\Omega)$.
An element $P \in \mathcal{P}$ is \emph{peripheral} if one of the two intervals within $\Omega$ defined by the endpoints of $P$ does not contain any vertices from $\mathcal{P} \setminus \{ P \}$.
We call $\mathcal{P}$ \emph{crooked} if it contains no peripheral element and \emph{planar} if no $\mathcal{P}' \subseteq \mathcal{P}$ is crooked.
Clearly every crooked transaction contains a crossing pair.

A first step towards our goal is the following simple lemma inspired by Lemma 3.1 in \cite{kawarabayashi2020quickly}.

\begin{lemma}\label{lem:planarorcrossing}
    If $\mathcal{P}$ is a transaction on a society $(D,\Omega)$ with $\mathcal{P} \geq p+2q-1$, such that $p,q \geq 1$ are integers, then $\mathcal{P}$ contains $q$ pairwise disjoint crossing pairs or a planar transaction of order $p$.
\end{lemma}
\begin{proof}
	We iteratively build two sets $\mathcal{P}', \mathcal{C} \subseteq \mathcal{P}$ with the objective that $\mathcal{P}'$ is a large planar transaction and $\mathcal{C}$ contains many pairwise disjoint crossing pairs.
	If $\mathcal{P}$ contains a peripheral element $P$, we add this element to $\mathcal{P}'$ and remove it from $\mathcal{P}$.
	Should $\mathcal{P}$ contain no such element, then it must contain a crossing pair $Q,Q'$, which we can add to $\mathcal{C}$ and remove from $\mathcal{P}$.
	We can repeat this process until $\mathcal{P}' \cup \mathcal{C} = \mathcal{P}$.
	Clearly we can ensure $\mathcal{P}' \cap \mathcal{C} = \emptyset$.
	
	If $|\mathcal{P}'| < p$, the set $\mathcal{C}$ contains at least $2q$ elements and according to our construction, at least $q$ pairwise disjoint crossing pairs.
	On the other hand, if $|\mathcal{C}| < 2q$, then $\mathcal{P}'$ must contain at least $p$ elements, all of which where peripheral at the point at which we added them and thus $\mathcal{P}'$ does not contain a crossing pair and is therefore planar.
\end{proof}

By invoking \cref{lem:applyshifting} we obtain the following corollary which guarantees that we always find a large planar transaction on a maelstrom society whose depth is unbounded.

\begin{corollary}\label{cor:manycyclesorplanartransaction}
Let $t,\theta,$ and $p$ be integers such that $\theta\geq 2t\SC$ is even, let $D$ be a digraph with an odd decomposition $\delta$.
Moreover, let $m$ be a maelstrom of $\delta$ together with a $\theta$-outline $\mathfrak{M}$.
Then the following two statements are true.
\begin{enumerate}
	\item If $\mathfrak{M}$ is a circle outline and $(H_{\mathfrak{M}},\Omega_{\mathfrak{M}})$ is its maelstrom society, then either $H_{\mathfrak{M}}$ contains a half-integral packing of $t$ even dicycles, or every transaction of order at least $p+2t-1$ on $(H_{\mathfrak{M}},\Omega_{\mathfrak{M}})$ contains a planar transaction of order $p$.
 
	\item If $\mathfrak{M}$ is a diamond outline and $\mathfrak{N}$ is a refined diamond $\theta$-outline in $\mathfrak{M}$, then either $H_{\mathfrak{M}}$ contains a third-integral packing of even dicycles, or every transaction of order at least $p+2t-1$ on the $\mathfrak{N}$-society contains a planar transaction of order $p$.
\end{enumerate}
\end{corollary}

\subsection{Many even graphs attaching to an outline}

We will need a few tools for the main theorem of this section and also require a few more definitions.

Let $\theta$ and $p$ be integers, let $D$ be a digraph, let $\delta$ be an odd decomposition of $D$, and let $m$ be a maelstrom of $\delta$, with a $\theta$-rim $\mathfrak{R}$.
Let $H'$ and $G$ be a subgraphs of $H_\mathfrak{R}$, such that $H' \subseteq G$ and $H'$ contains a transaction $\mathcal{P} = \{ P_1, \ldots , P_p \}$ on the $\mathfrak{R}$-society.
Further, let $X$ and $Y$ be the two minimal intervals on $\Omega_\mathfrak{R}$ containing the heads, respectively the tails, of all paths of $\mathcal{P}$.
We call $B = X \cup Y \cup V(P_1) \cup V(P_p)$ the \emph{strip boundary} of $\mathcal{P}$.
We say that \emph{$\mathcal{P}$ frames $H'$ in $G$}, if there exists no $H'$-$(G - V(H'))$- and no $(G - V(H'))$-$H'$-path in $G - B$.

Using this notion, we provide a technical statement that argues that, within an outline that has a planar transaction to which a non-even graph is attached, we can find a non-even decomposition of the outline together with most of the transaction and the graph attached to it.
We note that the following statement is formulated for outlines but can clearly also be applied to rims.
The broader statement will be useful later on.

\begin{lemma}\label{lem:integratenoneventransaction}
    Let $p$ be positive integers with $p \geq 2$, let $\theta \geq 2\SC$ be an even integer, and let $\ell = p + 4$.
    Let $D$ be a digraph, let $\delta=(\Gamma,\mathcal{V},\mathcal{D})$ be an odd decomposition of $D$, let $m$ be a maelstrom of $\delta$, with a $\theta$-outline $\mathfrak{M} = ( H, \mathcal{C}, \mathfrak{E}, \mathfrak{V} )$, with $\mathfrak{N} = ( H^*, \mathcal{C}, \mathfrak{E}, \mathfrak{B}, \mathfrak{V}^*, \mathfrak{D} )$ being a refined $\theta$-outline, if $\mathfrak{M}$ is a diamond outline, and an organised $\mathcal{C}$-pair of order $\lfloor \SC \rfloor$, if $\mathfrak{M}$ is a circle outline.
    Furthermore, let $(G, \Omega)$ be the $\mathfrak{M}$-society, or the $\mathfrak{N}$-society, if $\mathfrak{N}$ exists and let $\mathcal{P} = \{ P_1, \ldots , P_\ell \}$ be a planar transaction on $(G, \Omega)$.
    Finally, let $\mathfrak{R}$ be a $\theta$-rim of $\mathfrak{M}$, or respectively $\mathfrak{N}$, and let $( H_k, \mathcal{C}_k, \mathcal{E}_k , \mathfrak{B}_k )$ be an $\SC$-layer in $\mathfrak{R}$ for $k \in [ \nicefrac{\theta}{2} - (\SC - 1)]$.

    If there exists a subgraph $H_\mathcal{P} \subseteq G$ such that $H_\mathcal{P}$ contains all paths in $\mathcal{P}$ and a subgraph $H_\mathcal{P}' \subseteq H_\mathcal{P}$ that is framed by $\mathcal{Q} = \{ P_{\ell + 3}, \ldots , P_{p + 2} \}$ in $H$ then if $H_\mathcal{P} \cup H_k$ is non-even, $H \cup H_\mathcal{P}'$ has a non-even decomposition.
\end{lemma}
\begin{proof}
Let $k \in [ \nicefrac{\theta}{2}  - (\SC - 1)]$, with $( H_k, \mathcal{C}_k, \mathcal{E}_k , \mathfrak{B}_k )$ being the $\SC$-layer at $k$ in $\mathfrak{R}$.
Since $\mathcal{C}$ has order at least $\lfloor \SC \rfloor$ and we are either provided with an organised $\mathcal{C}$-pair of order $\lfloor \SC \rfloor$, or we can derive such a pair by using the paths in $\mathcal{E}$, we can find a butterfly minor model $B$ of a cylindrical grid of order three in $H_k$ using \Cref{cor:getseparatingcylindricalgrid}.
Based on $B$, we can construct a reinforced society, which cannot have a cross in the split, as we are otherwise not non-even as proven in \Cref{lem:keylemmashifting}.
This allows us to apply \Cref{thm:reinforcedmatcross} to find a non-even rendition $\rho'$ of the subgraph $H'$ of $H_\mathcal{P} \cup H_k$ that is drawn on the inclusion-wise largest $m$-disk $d_B$ of any dicycle found inside of $B$.

Note that due to $\rho'$ being a rendition, the subgraph of $H'$ drawn inside $d_B$, but outside of the interior of the inclusion-wise minimal $m$-disk of of any dicycle in $B$ is already in agreement with $\delta$.
Thus we can concentrate on how $\rho'$ decomposes $H_\mathcal{P}$.

We note that, even though $H_\mathcal{P} \cap H$ has a plane decomposition in $\delta$ as well, this part of the graph may be decomposed in slightly different ways in $\rho'$ and $\delta$, as much of what is drawn in the proximity of $P_3$ and $P_{p + 2}$ within $H$ is part of $H_\mathcal{P}$.
This may cause the appearance of big vertices and conjunction cells involving vertices and edges of $P_1$, $P_2$, $P_{p + 3}$, and $P_{p + 4}$ within $\rho'$ that did not occur in $\delta$.
However, due to directed tight cuts not allowing two disjoint dicycles passing through them and small cycle sums similarly having a limit of at most two disjoint cycles passing through them, the vertices and edges of any path in $\mathcal{Q}$ cannot be involved in such big vertices or conjunction cells.
(There is a small exception for big vertices or conjunction cells whose interior in fact has a cross-free drawing.
But such big vertices and conjunction cells do not have to be drawn in a separate disk, allowing us to just keep the existing drawing.)

A similar phenomenon may occur near $m$ when considering the parts of $H_\mathcal{P}$ and $H$ drawn in the immediate proximity of $C_{\nicefrac{\theta}{2}}$.
Here it could be the case that in $\rho'$ some part of $\sigma_\delta(m)$ is decomposed together with vertices or edges of $C_{\nicefrac{\theta}{2}}$ inside of a big vertex or a small cycle sum.
This would mean that $\delta$ and $\rho'$ do not agree in this part of the decomposition, since $\delta$ has put this part of the drawing into the maelstrom $m$.
However, due to plane decompositions allowing us to have non-planarity inside of conjunction cells and big vertices, we can adjust $\delta$ such that it also uses these same big vertices and conjunction cells of $\rho'$ in its decomposition.

According to the arguments we just laid out, we can adjust $\delta$ such that $\rho'$ and $\delta$ agree on the decomposition of $H_\mathcal{P}' \cup H_k - \sigma(m)$ and in particular, the decomposition of $H_\mathcal{P} \cup (H_k \cap \sigma(m))$ in $\rho'$ is compatible with the decomposition $\delta$ in such a way that we can derive a non-even decomposition of $H \cup H_\mathcal{P}$.
\end{proof}

Let $t$, $\theta$, $d$, and $g$ be positive integers with $\theta > t \geq 1$, let $\mathfrak{R} = (H, \mathcal{C}, \mathcal{E}, \mathfrak{B}, \mathfrak{D})$ be a $\theta$-rim, with eddy degree $d$ and roughness $g$, if $\mathfrak{R}$ is a refined diamond rim, and $\mathcal{C} = \{ C_1, \ldots , C_\theta \}$, $\mathcal{E} = \{ E_1, \ldots , E_\theta \}$, $\mathfrak{B} = \{ \mathcal{B}_1, \ldots , \mathcal{B}_d \}$, and $\mathcal{B}_i = \{ B^i_1, \ldots , B^i_{\theta(g+1)} \}$, for each $i \in [d]$.
The \emph{$t$-retreat of $\mathfrak{R}$} is defined as the $(\theta - t)$-rim $\mathfrak{R}' = ( H' , \mathcal{C}' = \{ C_{t+1}, \ldots , C_\theta \}, \mathcal{E}' = \{ E_{t+1}, \ldots , E_\theta \}, \mathfrak{B}' = \{ \mathcal{B}_1', \ldots , \mathcal{B}_d' \} )$, with $\mathcal{B}'_i = \{ B^i_{t(g+1) +1} , \ldots , B^i_{\theta(g+1)} \}$, for each $i \in [d]$, and $H'$ is the natural restriction of $H$ to the area described by the infrastructure of $\mathfrak{R}'$.
Note that $\mathfrak{R}'$ retains its eddy degree and roughness, if $\mathfrak{R}$ was a refined diamond rim.

The proof of the main theorem of this section will essentially consist of two phases.
First, we try to gather $t$ disjoint subgraphs, each hosting a large planar transaction, and then we want to either find a quarter-integral packing of $t$ even dicycles in these subgraphs together with the outline, or find one subgraph among them that is odd with respect to the $t$-retreat of the rim of the original outline.
This intermediate step is quite technical and will be dealt with in the next two lemmas.
Afterwards, if we find such an odd subgraph in the proof of the main theorem, we will have to ensure that it is decomposable in a way that is compatible with our outline.
This will be dealt with in the proof proper.

\begin{lemma}\label{lem:dicyclesthroughatransaction}
    Let $t, p$ be positive integers with $p \geq 2$, let $\theta \geq 2(\SC + 1)t$ be an even integer, and let $\ell = p + 4$.
    Let $D$ be a digraph, let $\delta$ be an odd decomposition of $D$, and let $m$ be a maelstrom of $\delta$, with a $\theta$-outline $\mathfrak{M}$, whose $\mathfrak{M}$-society is $(H_\mathfrak{M}, \Omega_\mathfrak{M})$.
    Moreover, let $\mathfrak{R} = (H, \mathcal{C}, \mathcal{E})$ be the $(\SC +1)t$-rim of $\mathfrak{M}$, with $(\mathcal{L}_1, \mathcal{L}_2)$ being an organised $\mathcal{C}$-pair of order $\lfloor \SC \rfloor$, if $\mathfrak{R}$ is a circle rim, let $\mathcal{P} = \{ P_1, \ldots , P_{t\ell} \}$ be a planar transaction on the $\mathfrak{M}$-society, and let $\mathfrak{R}^* = ( H^* , \mathcal{C}^*, \mathcal{E}^* )$ be the $t$-retreat of $\mathfrak{R}$.

    If there exist $t$ pairwise disjoint subgraphs $H_1, \ldots, H_t \subseteq H_\mathfrak{M}$ such that $H_i$ contains all paths in $\mathcal{P}_i = \{ P_{(i-1)\ell + 1}, \ldots , P_{i \ell} \}$, the subgraph $H_i' \subseteq H_i$ is framed by $\mathcal{Q}_i = \{ P_{(i-1)\ell + 3}, \ldots , P_{(i-1)\ell + p + 2} \}$ in $H \cup H_i$, and $H \cup H_i'$ contains an even dicycle for each $i \in [t]$, then either $D$ contains a third-integral packing of $t$ even dicycles or there exists an $i \in [t]$ such that $H^* \cup H_i$ is non-even and $H^* \cup H_i'$ is odd.
\end{lemma}
\begin{proof}
Let $\mathcal{H} = \{ H_1 , \ldots , H_t \}$ and $\mathcal{H}' = \{ H_1' , \ldots , H_t' \}$.
We note that since $\mathfrak{M}$ is an outline and $\delta$ is an odd decomposition, any even dicycle in any graph $H_i' \cup H$ with $i \in [t]$ must be hosted by $m$.
To find the even dicycles we desire, we will prove via induction on $t$ that we can always either find a graph $H_i \in \mathcal{H}$, with $i \in [t]$, such that $D_i = H_i \cup \bigcup_{j=(i-1)\SC + 1 + t}^{i\SC + t} (C_j \cup E_j) \cup U_i$ contains an even dicycle, where $U_i$ is the graph in the $\SC$-layer at $(t-1)\SC + 1 + t$ in $\mathfrak{R}$, or we find an even dicycle entirely in $H_i' \cup C_j \cup E_j$, with $i,j \in [t]$ and $j \leq i$.
Since $\theta \geq \SC t$, this then allows us to apply the induction hypothesis by removing the graphs from $\mathcal{H}$ and $\mathcal{H}'$ and the involved infrastructure from $\mathfrak{R}$.
The resulting packing of even dicycles is third-integral since we choose disjoint parts of $\mathcal{C}$ and $\mathcal{E}$ and the graphs in $\mathcal{H}$ are disjoint.
For $t \leq 3$ the statement trivially holds.

Consider the graph $D_t$.
If $D_t \cup H$ is not non-even it contains a weak odd bicycle according to \Cref{thm:nonevendigraphs} and therefore also an even dicycle and we are done.
Thus \Cref{lem:integratenoneventransaction} yields a non-even decomposition $\rho_t$ of $G_t = H_t' \cup \bigcup_{j=(t-1)\SC + 1 + t}^{t\SC + t} (C_j \cup E_j) \cup U_t$.

According to our assumptions, within $H \cup H_t'$ there exists an even dicycle $C$.
If $H_t' \cup C_i \cup E_i$ is even, for any given $i \in [t]$, we are done and should $C$ solely be contained in $G_t$, we are also done.
Accordingly, we may suppose that $G_t$ and $H_t' \cup C_1 \cup E_1$ are odd.
Thus $C$ must have some edges outside of $G_t$ and, since $\mathcal{Q}_t$ frames $H_t'$ in $H \cup H_t$, the cycle $C$ therefore intersects $P_1' = P_{(t-1)\ell + 3}$ or $P_2' = P_{(t-1)\ell + 2 + p}$ in at least two vertices.
We note that all such edges must lie in $H$, as $(H \cup H_t') - G_t \subseteq H$.
Our goal will therefore be to show that we can find an even dicycle in $H$, leading to a contradiction.

Suppose that $\mathfrak{R}$ is a circle rim.
Let $C_1'$ be an odd dicycle in $P_1' \cup C_{t\SC+1+t}$, such that the trace of $C_1'$ bounds a disk that does not contain the traces of any path in $\mathcal{Q}_t \setminus \{ P_1' \}$.
Similarly, we let $C_2'$ be an odd dicycle in $P_2' \cup C_{t\SC+1+t}$, such that the trace of $C_2'$ bounds a disk that does not contain the traces of any path in $\mathcal{Q}_t \setminus \{ P_2' \}$.
Let $P$ be the path corresponding to the intersection of $P_1'$ and $C_1'$.

First, let us note how $P$ interacts with $C_1 \in \mathcal{C}$.
Each $C_1$-path $P'$ within $P$ forms a dicycle $C_1''$ with $C_1$ that is distinct from $C_1$.
Since $C_1''$ is found in $H_t' \cup C_1$, it must be odd, or we are again done, and the unique $C_1''$-path in $C_1 \cup P$ must therefore have the same parity as $P$.

We can now use \Cref{lem:planartraceshifting} to shift $C$ using $C_1'$ and thus find a new even dicycle $C' \subseteq C_1' \cup C$.
We note here that if \Cref{lem:planartraceshifting} yields an even dicycle in some conjunction cell, we are also done, since this conjunction cell must be found in the decomposition of $H_i$ and therefore $H_i$, since we ensured that $\delta$ and $\rho_t$ agree on the decomposition of $H$.
This yields two cases, either $C'$ is found on the disk bounded by $C_1'$ that does not contain the traces of any path in $\mathcal{Q}_t \setminus \{ P_1' \}$, or it is found outside the interior of this disk.
We will start with the former case.

Clearly, $C$ cannot be found entirely in $H$, as this contradicts $H$ being odd.
Therefore some edge of $C$ must lie in $E(\sigma(m)) \setminus E(C_1)$.
However, if we consider $C \cap \sigma(m)$, we note that this graph is entirely contained in $C_1 \cup P$ and thus another application of \Cref{lem:planartraceshifting} will yield an even dicycle outside of the interior of $m$, a contradiction to the fact that $H$ is odd.

Therefore, $C'$ must instead be found on the disk bounded by $C_1'$ that contains the traces of all paths in $\mathcal{Q}_t$.
Here we may again apply \Cref{lem:planartraceshifting} to $C'$ and $C_2'$ to get a new even dicycle $C''$ that lies on either side of the trace of $C_2'$.
If $C''$ is found on the side of the trace of $C_2'$ that contains the traces of the paths in $\mathcal{Q}_t \setminus \{ P_2' \}$, we are done, since $C'' \subseteq H_i'$.
Otherwise, we can again shift $C''$ using $C_1$ to derive a contradiction to the fact that $H$ is odd.

Thus we may suppose instead that $\mathfrak{R}$ is a diamond rim.
We reuse names from the previous section here for convenience.
Similarly to the last case, we let $C_1'$ be an odd dicycle in $P_1' \cup C_{t\SC+1+t} \cup E_{t\SC+1+t}$, such that the trace of $C_1'$ bounds a disk that does not contain the intersection of $m$ and the traces of the paths in $\mathcal{Q}_t \setminus \{ P_1' \}$.
Analogously, we let $C_2'$ be an odd dicycle in $P_2' \cup C_{t\SC+1+t} \cup E_{t\SC+1+t}$, such that the trace of $C_2'$ bounds a disk that does not contain the intersection of $m$ and traces of the paths in $\mathcal{Q}_t \setminus \{ P_2' \}$.
Let $P$ again be the path corresponding to the intersection of $P_1'$ and $C_1'$.

We note that aside from $C_1$, there also exists another $R_1$ in $C_1 \cup E_1$ distinct from $C_1$, such that $R_1 \cup C_1$ together contain the diamond whose trace bounds $m$.
Now we can proceed as in the previous case, first shifting $C$ using $C_1'$ to yield a new even dicycle $C'$, which we first suppose to be located on the disk bounded by $C'$ that does not contain the intersection of $m$ and the traces of the paths in $\mathcal{Q}_t \setminus \{ P_1' \}$.
Again, we note that $C' \cap \sigma(m)$ has to be odd and thus we may shift $C'$ into $H$ via either $C_1$ or $R_1$, yielding a contradiction.
Thus $C'$ lies on the other side of $C_1'$ and we can shift $C'$ into $H_i'$ using $C_2'$.
This concludes our proof.
\end{proof}

Of course an analogous statement can be proven for refined diamond outlines.
We chose to separate this from the previous lemma, since the statement and proof for refined diamond outlines becomes even more technical, obscuring the core of the statement.
This lemma can be proven in an analogous fashion to \Cref{lem:dicyclesthroughatransaction}.

\begin{lemma}\label{lem:dicyclesthroughatransactionrefined}
    Let $t, p$ be positive integers with $p \geq 2$, let $\theta \geq 2(\SC + 1)t$ be an even integer, and let $\ell = p + 4$.
    Let $D$ be a digraph, let $\delta$ be an odd decomposition of $D$, and let $m$ be a maelstrom of $\delta$, with a $\theta$-outline $\mathfrak{M}$, let $\mathfrak{N}$ be a refined diamond $\theta$-outline around $n \subseteq m$, and let $(H_\mathfrak{N}, \Omega_\mathfrak{N})$ be the $\mathfrak{N}$-society.
    Moreover, let $\mathfrak{R} = (H, \mathcal{C}, \mathcal{E}, \mathfrak{B}, \mathfrak{D})$ be the $(\SC +1)t$-rim of $\mathfrak{N}$, let $\mathcal{P} = \{ P_1, \ldots , P_{t\ell} \}$ be a planar transaction on the $\mathfrak{N}$-society, and let $\mathfrak{R}^* = ( H^* , \mathcal{C}^*, \mathcal{E}^*, \mathfrak{B}^*, \mathfrak{D}^* )$ be the $t$-retreat of $\mathfrak{R}$.

    If there exist $t$ pairwise disjoint subgraphs $H_1, \ldots, H_t \subseteq H_\mathfrak{M}$ such that $H_i$ contains all paths in $\mathcal{P}_i = \{ P_{(i-1)\ell + 1}, \ldots , P_{i \ell} \}$, the subgraph $H_i' \subseteq H_i$ is framed by $\mathcal{Q}_i = \{ P_{(i-1)\ell + 3}, \ldots , P_{(i-1)\ell + 2 + p} \}$ in $H \cup H_i$, and $H \cup H_i'$ contains an even dicycle for each $i \in [t]$, then either $D$ contains a quarter-integral packing of $t$ even dicycles or there exists an $i \in [t]$ such that $H^* \cup H_i$ is non-even and $H^* \cup H_i'$ is odd.
\end{lemma}


\subsection{Finding an odd transaction}

As explained earlier, we will have to ensure that the transaction we are given, if we cannot find a quarter-integral packing of $t$ even dicycles, is not only odd, but also decomposable in a way that lets us integrate it into the existing outline.
We formalise the correct notion for this in the next definition.

\begin{definition}[Odd transaction]\label{def:oddtransaction}
	Let $t$, $\theta$, and $p$ be integers such that $\theta$ is even and let $D$ be a digraph with an odd decomposition $\delta$.
	Moreover, let $m$ be a maelstrom of $\delta$ together with a $\theta$-outline $\mathfrak{M}$, in case $\mathfrak{M}=(H,\mathcal{C},\mathfrak{E},\mathfrak{V})$ is not a circle outline, let $\mathfrak{N}=(H^*,\mathcal{C},\mathfrak{E},\mathfrak{B},\mathfrak{V}^*,\mathfrak{D})$ be a refined diamond $\theta$-outline around $n \subseteq m$ in $\mathfrak{M}$.
	In case $\mathfrak{M}$ is a circle outline, let $(K,\Omega)$ be its maelstrom society, otherwise let $(K,\Omega)$ be its $\mathfrak{N}$-society.
	Assume that $K$ is strongly connected.

	Let $\mathcal{T}=\{ T_1,\dots,T_{\ell}\}$ be a planar transaction on $(K,\Omega)$ with the strip boundary $B$.
	Let $A\subseteq V(D)$ be a set of vertices with $A\cap (V(\mathcal{T})\cup V(\Omega))=\emptyset$.
	The \emph{strip} of $\mathcal{T}$ \emph{under $A$} is a subgraph $S$ of the strong component $Q$ of $K-A$ that contains $V(\Omega)$, where $S$ is induced by $B$ together with the collection of all vertices $v \in V(Q)$ such that there exists a directed $v$-$V(\mathcal{T})$-path or a directed $V(\mathcal{T})$-$v$-path in $Q-B$.

	Let $\mathcal{T}$ be a planar transaction on $(K,\Omega)$, a set of vertices $A\subseteq V(D)$ with $A\cap(V(\mathcal{T})\cup V(\Omega))=\emptyset$, and let $S$ be the strip of $\mathcal{T}$ under $A$.
	We say that $\mathcal{T}$ is \emph{odd} with the \emph{apex set} $A$ if either, $\mathfrak{M}$ is a circle outline and the graph $(H\cup S)-A$ is odd, or $\mathfrak{M}$ is a diamond outline and the graph $(H^*\cup S)-A$ is odd.
\end{definition}

Beyond the correct definition for an odd transaction, we will also have to define our goal, which is a decomposition which combines the decomposition of an existing outline with the decomposition of an odd transaction.

\begin{definition}[Expansion of an odd decomposition]\label{def:compatibleoddtransaction}
    Let $\theta , t, k , p$ be integers such that $\theta$ is even and $p$ is divisible by 4, and let $D$ be a digraph with an odd decomposition $\delta$.
	Moreover, let $m$ be a maelstrom of $\delta$ together with a $\theta$-outline $\mathfrak{M} = ( H , \mathcal{C}, \mathfrak{E}, \mathfrak{V} )$ of degree $t$ and a refined diamond $\theta$-outline $\mathfrak{N}= ( H^*, \mathcal{C}, \mathfrak{E}, \mathfrak{B}, \mathfrak{V}^*, \mathfrak{D} )$ with surplus degree $k$ for some $n \subseteq m$, if $\mathfrak{M}$ is a diamond outline, with $(H_\mathfrak{M},\Omega_\mathfrak{M})$ being the $\mathfrak{M}$-society.
 
    We call an odd transaction $\mathcal{P} = \{ P_1, \ldots , P_p \}$ on the $\mathfrak{M}$-, or respectively the $\mathfrak{N}$-society, and an apex set $A \subseteq V(H_\mathfrak{M})$, with $S$ being the strip of $\mathcal{P}$ under $A$, \emph{$\delta$-compatible} for $\mathfrak{M}$, or $\mathfrak{N}$, if
    \begin{itemize}
        \item there exists an odd decomposition $\delta' = (\Gamma', \mathcal{V}', \mathcal{D}')$ of $H_\mathfrak{M} - A$ that agrees with $\delta$ on $H-A$,

        \item the society $((H \cup S) - A, \Omega_\mathfrak{M})$ has a non-even rendition $\rho$, extending the non-even rendition of $(H,\Omega_\mathfrak{M})$ associated with $\mathfrak{M}$, in which all paths of $\mathcal{P}$ are grounded,

        \item for each $i \in [p]$, the edges of $P_i$ are drawn by $\Gamma'$ outside of any maelstrom of $\delta'$,

        \item for all $i \in [2, p-1]$, the trace of $P_i$ in $\delta'$ separates the traces of $P_{i-1}$ and $P_{i+1}$ on the $m$-disk $\Delta_m$, or respectively the $n$-disk $\Delta_n$, of the cycle in $D$ that defines the $\mathfrak{M}$-, or respectively the $\mathfrak{N}$-society, and

        \item the restriction of $\rho$ to the subgraph of $S$ framed by $\{ P_{\nicefrac{p}{4}}, \ldots , P_{\nicefrac{3p}{4} + 1} \}$ together with all nodes drawn on and all conjunction cells of $\delta$ intersecting $\Delta_\mathfrak{M} \setminus (m \setminus \Boundary{m})$, or respectively $\Delta_\mathfrak{N} \setminus (m \setminus \Boundary{m} \cup \bigcup_{d \in \mathfrak{D} (d - \Boundary{d}) )} )$, is an odd rendition that is a restriction of $\delta$.
    \end{itemize}
    The odd decomposition $\delta'$ is called a \emph{$\mathcal{P}$-$A$-expansion of $\delta$}.
    Furthermore, an apex set $A \subseteq V(H_\mathfrak{M})$ is called \emph{pleasant}, if it contains no vertices of elements in $\mathcal{C} \cup \bigcup_{i \in [t]} \mathcal{E}_i \cup \bigcup_{i \in [k]} \mathcal{B}_i$.
\end{definition}

Sadly, we cannot ensure that the apex set $A$ we will find in the next theorem is necessarily pleasant, but this can be remedied by simply removing all parts of the infrastructure that are hit by an apex vertex.
The following statement is technical, but its proof is not surprising.
For this reason, we only provide a rough sketch of it.

\begin{lemma}\label{lem:flowalteredoutline}
    Let $\theta , g$ be integers such that $\theta$ is even, and let $D$ be a digraph with an odd decomposition $\delta$.
	Moreover, let $m$ be a maelstrom of $\delta$ together with a $\theta$-outline $\mathfrak{M} = ( H , \mathcal{C}, \mathfrak{E}, \mathfrak{V} )$, a refined diamond $\theta$-outline $\mathfrak{N}= ( H^*, \mathcal{C}, \mathfrak{E}, \mathfrak{B}, \mathfrak{V}^*, \mathfrak{D} )$ with roughness $g$ for some $n \subseteq m$, if $\mathfrak{M}$ is a diamond outline, with $(G,\Omega)$ being the $\mathfrak{M}$-, or respectively the $\mathfrak{N}$-society, and an organised $\mathcal{C}$-pair of order $\theta$, if $\mathfrak{M}$ is a circle outline.
 
    Let $A \subseteq V(G)$ be a set of vertices with $|A| < \nicefrac{\theta}{6}$, then there exists a $(\theta - 6|A|)$-outline $\mathfrak{M}'$ for $m$, with cycle set $\mathcal{C}'$, or respectively a refined diamond $(\theta - 6|A|)$-outline $\mathfrak{N}'$ with roughness $g$ for $n$, such that $A$ is pleasant for $\mathfrak{M}'$, respectively $\mathfrak{N}'$, and an organised $\mathcal{C}'$-pair of order $\theta$, if $\mathfrak{M}'$ is a circle outline.
\end{lemma}
\begin{proof}[Rough Sketch]
At worst we have to remove $3|A|$ different elements of $\mathcal{C} \cup \bigcup_{\mathcal{E} \in \mathfrak{E}} \mathcal{E} \cup \bigcup_{\mathcal{B} \in \mathfrak{B}} \mathcal{B}$ from the sets in the infrastructure, since each vertex may hit a vertex that lies in the intersection of a cycle from $\mathcal{C}$, a path from a $\mathcal{E} \in \mathfrak{E}$, and a path from $\mathcal{B} \in \mathfrak{B}$.
Note that we do not delete the infrastructure, but instead simply remove these objects from the sets of the outline.
We do not know whether the vertices of $A$ and the parts of the infrastructure we remove lies inside our outside $m$, so to preserve the balance of the outline and ensure that the properties of the outline are preserved, we can at worst delete $3|A|$ more parts of the outline, to ensure that there is as much infrastructure outside of $m$ as there is inside of it.
This also allows us to derive all properties of the new outline from the old one.

One remaining thing to remark is that if some vertex in $A$ hits the society, then we must give up this part of the infrastructure.
This will only move the society further inwards and does not cause any issues with the decomposition.
However, we will have to truncate linkages in $\bigcup_{\mathcal{E} \in \mathfrak{E}} \mathcal{E} \cup \bigcup_{\mathcal{B} \in \mathfrak{B}} \mathcal{B}$ which may cause both the degree and the eddy degree of the outline to change, though the roughness can be maintained.
We will also have to mildly adjust $m$ and the elements of $\mathfrak{V}$, $\mathfrak{V}^*$, and $\mathfrak{D}$ to fit the new (not necessarily directed) cycle that lies at the halfway point in our infrastructure.
As $|A| < \nicefrac{\theta}{6}$ and we balanced out our removal process, the new society has no chance of touching the modified maelstrom and thus this will not produce any new outlines or split $m$ into new maelstroms.
\end{proof}

We now finally proceed to the proof of the main theorem of this section.

\begin{theorem}\label{thm:oddtransaction}
	Let $t$, $\theta$, and $p$ be integers.
    There exist functions $\OddTransactionApexNoArg \colon \N \rightarrow \N$ and $\OddTransactionOrderNoArg \colon \N \times \N \rightarrow \N$ with $\OddTransactionApexNoArg \in \mathcal{O}(t^3)$ and $\OddTransactionOrderNoArg \in \mathcal{O}(t^5, p)$ such that $\theta \geq 2t\SC + 6\OddTransactionApex{t}$ is even and $p \geq t > 1$, and let $D$ be a digraph with an odd decomposition $\delta$.
	Moreover, let $m$ be a maelstrom of $\delta$ together with a $\theta$-outline $\mathfrak{M}$.
	\begin{enumerate}
		\item If $\mathfrak{M}$ is a circle outline, with an organised $\mathcal{C}$-pair of order $\lfloor \SC \rfloor + \OddTransactionApex{t}$, and $(H,\Omega)$ is its maelstrom society, then either $H$ contains a half-integral packing of $t$ even dicycles, or for every transaction $\mathcal{T}$ of order at least $\OddTransactionOrder{t,p}$ on $(H,\Omega)$, there exists an apex set $A$ with $|A| \leq \OddTransactionApex{t}$, a circle $(\theta - 6|A|)$-outline $\mathfrak{M}'$, with cycle set $\mathcal{C}'$, and an organised $\mathcal{C}'$-pair of order $\lfloor \SC \rfloor$ and the society $(H',\Omega')$, and an odd transaction $\mathcal{P}$ of order $p$ such that $A$ and $\mathcal{P}$ on $(H',\Omega')$ are $\delta$-compatible for $\mathfrak{M}'$.
  
		\item If $\mathfrak{M}$ is a diamond outline and $\mathfrak{N}$ is a refined diamond $\theta$-outline in $\mathfrak{M}$ with refined maelstrom society $(H,\Omega)$, then either $H$ contains a quarter-integral packing of $t$ even dicycles, or for every transaction $\mathcal{T}$ of order at least $\OddTransactionOrder{t,p}$ on $(H,\Omega)$, there exists an apex set $A$ with $|A| \leq \OddTransactionApex{t}$, a diamond $(\theta - 6|A|)$-outline, together with a refined diamond $(\theta - 6|A|)$-outline $\mathfrak{N}'$, with the $\mathfrak{N}$-society $(H',\Omega')$, and an odd transaction $\mathcal{P}$ on $(H',\Omega')$ of order $p$ such that $A$ and $\mathcal{P}$ are $\delta$-compatible for $\mathfrak{N}'$.
	\end{enumerate}
\end{theorem}
\begin{proof}
    In the following, let $(G,\Psi)$ denote the society $(H_{\mathfrak{M}},\Omega_{\mathfrak{M}})$ in case $\mathfrak{M}$ is a circle outline, and otherwise let it denote the society $(H_{\mathfrak{N}},\Omega_{\mathfrak{N}})$.
    Moreover, we set
    \[ q \coloneqq 72t^3 - 104t^2 + 38t + 2 = (18t+1) (2t - 1) (2t - 2) + 8t , \ p' \coloneqq p + (36t + 2)q, \text{ and } p_0 \coloneqq (2t - 1) p' . \]

    Using \cref{cor:manycyclesorplanartransaction} we may assume that every transaction of order $p_0 + 2t - 1$ contains a planar transaction $\mathcal{P}$ of order $p_0$, as otherwise we immediately find the required packing of even dicycles.
    Note that $p_0 + 2t - 1 = (2 t - 1) (p + (36 t + 2) (72 t^3 - 104 t^2 + 38 t + 2)) + 2 t - 1$.
	Hence, it suffices to only consider planar transactions from now on and assume that such a transaction $\mathcal{P}$ exists.

    In the first half of the proof we will progressively build a selection of $t$ transactions with associated graphs that are disjoint, allowing us to apply \Cref{lem:dicyclesthroughatransaction} and \Cref{lem:dicyclesthroughatransactionrefined} to either find a quarter-integral packing of $t$ even dicycles, or a single transaction among them which is non-even and an associated restriction of it which is odd.
    In the former case, we are immediately done and in the latter case, we must argue that we can refine the odd part sufficiently to find the odd transaction promised in the statement of the lemma.

    We let $\mathcal{P} = \{ P_1 , \ldots , \mathcal{P}_{p_0} \}$, such that for all $i \in [2, p_0 - 1]$ the endpoints of $P_i$ separate the endpoints of $P_{i-1}$ from the endpoints of $P_{i+1}$ on $\Psi$.
    For each $i \in [2t - 1]$, we let $\mathcal{Q}_i = \{ Q_1^i, \ldots , Q_{p'}^i \}$, with $Q_j^i = P_{(i-1)p' + j}$ for all $j \in [p']$.
    For each $i \in [2t - 1]$ and $j \in [18t+1]$, we set 
    \[ \mathcal{S}^i_{j,R} = \{ Q_{(j-1)q + 1}^i, \ldots Q_{jq}^i \}, \ \mathcal{S}^i_{j,L} = \{ Q_{p' - (18t+2 - j)q + 1}^i, \ldots Q_{p' - ( (18t+1 - j) q )}^i \}, \text{ and }  \mathcal{S}^i_j = \mathcal{S}^i_{j,R} \cup \mathcal{S}^i_{j,L} . \]
    We let $X_j^i$ be the union of the two maximal segments of $\Psi$ that contain all heads, respectively all tails, of the paths in $\mathcal{Q}_i \setminus \bigcup_{h=1}^j \mathcal{S}^i_j$ and no endpoints of the remaining paths in $\mathcal{Q}_i$.
    Furthermore, for $j \neq 1$, we let $Y_j^i$ be the union of the two maximal segments of $\Psi$ that together contain all endpoints of the paths in $\bigcup_{h = 1}^{j-1} \mathcal{S}^i_j$ and for $j = 1$, we let $Y_j^i$ be $V(\Psi)$ minus the vertices of the union of the two minimal segments contain all heads, respectively all tails, of $\mathcal{Q}_i$. 
    
    \begin{claim}\label{claim:findsafetyapex}
        For all $i \in [2t - 1]$ and $j \in [18t+1]$, there either exists a quarter-integral packing of $t$ even dicycles in $G$, or there exists a set of vertices $A_j^i \subseteq V(G) \setminus V( \mathcal{S}_j^i )$ of size at most $2(t-1)$ such that there are no $X_j^i$-$Y_j^i$- and no $Y_j^i$-$X_j^i$-paths in $G - ( V(\mathcal{S}_j^i) \cup A_j^i )$.
	\end{claim}
    \emph{Proof of \Cref{claim:findsafetyapex}:}
    Suppose that there does not exist a quarter-integral packing of $t$ even dicycles in $G$.
	We apply \Cref{thm:directedlocalmenger}, to either get a set $A$ of size at most $t-1$ disrupting all $X_j^i$-$Y_j^i$-paths in $G - V( \mathcal{S}_j^i )$, or we find at least $t$ pairwise disjoint $X_j^i$-$Y_j^i$-paths $P_1, \ldots , P_t$ in $G - V( \mathcal{S}_j^i )$.
    Suppose we find the paths, then by construction, these must be disjoint from the paths in $\mathcal{S}_j^i$ and each path in $\{ P_1, \ldots , P_t \}$ forms a cross with each path of either $\mathcal{S}^i_{j,R}$ or $\mathcal{S}^i_{j,L}$.
    By applying \cref{lem:applyshifting}, we can thus find a quarter-integral packing of $t$ even dicycles in $G$.
    Therefore, we must find the set $A$.
    Via an analogous argument, we find a similar set $A'$ that destroys all $Y_j^i$-$X_j^i$-paths in $G - V( \mathcal{S}_j^i )$.
    By setting $A_j^i = A \cup A'$, we have found our desired set.
	\hfill$\blacksquare$

    We let $A_1 = \bigcup_{i=1}^{2t-1} \bigcup_{j=1}^{18t+1} A_j^i$ and note that
    \[ |A_1| \leq (18t+1) (2t-1) (2t - 2) = 72t^3 - 104t^2 + 30t + 2 = q - 8t . \]
    In particular, we observe that, as a consequence of $t \geq 1$ and the bound on the size of $A_1$ for each $i \in [2t - 1]$, amongst both $Q_{(18t+1)q + 1}^i, \ldots , Q_{(18t+1)q + 8t}^i$ and $Q_{p' - ((18t+1)q + 8t - 1)}^i , \ldots , Q_{p' - (18t+1)q}^i$ there are $t$ paths that are disjoint from $A_1$, making for at least $2t$ paths that are not hit by $A_1$ and are also not found in $\mathcal{S}^i \coloneqq \bigcup_{j = 1}^{(18t+1)} \mathcal{S}_j^i$.
    We let $R_i$ be the path amongst these that has the lowest index and let $L_i$ be the path with the highest index.
    For each $i \in [2t-1]$ we let $X_i'$ be the union of the two minimal segments of $\Psi$ that contain both heads, respectively tails of $R_i$ and $L_i$, let $X_i = X_i' \setminus V( \{ R_i , L_i \} )$, and let $Y_i = Y_1^i$, which we defined earlier.
    We also let $Z_i = V(\Psi) \setminus ( X_i \cup Y_i )$ and start with a basic observation.

    \begin{claim}\label{claim:movethroughsafety}
        For each $i \in [2t-1]$, for all $X_i$-$Y_i$- and all $Y_i$-$X_i$-paths $P$ in $G - (A_1 \cup V( \{ R_i , L_i \} ) )$, we have $J \cap V( \mathcal{S}_j^i ) \neq \emptyset$ for all $j \in [18t+1]$.
	\end{claim}
    \emph{Proof of \Cref{claim:movethroughsafety}:}
    This is an immediate consequence of \Cref{claim:findsafetyapex} and the definition of $A_1$.
	\hfill$\blacksquare$

    The following claim is the technical core of the proof.
    Its idea is similar to \Cref{claim:nolongjumps} from the proof of \Cref{thm:oddwall}.
    We have cleaned up several sections around the part of the transaction that is separated from the outside by $R_i$ and $L_i$ and we will use \Cref{claim:movethroughsafety} to take an $X_i$-$Y_i$- or a $Y_i$-$X_i$-path and split it up into several jumps over paths in our transaction, which give us usable crosses for \Cref{lem:applyshifting}.
    Sadly, unlike \Cref{claim:nolongjumps}, this does not work in full generality, but if it does not, we can guarantee the existence of a path that does not jump out very far from the part of the transaction we isolated and with some work, we can make this benefit it us as well.

    From now on, we will say that a path $Q \in \mathcal{P}$ with $V(Q) \cap A_1 = \emptyset$ has been \emph{spared by $A_1$}.

    \begin{claim}\label{claim:farjumpgivespacking}
        For each $i \in [2t-1]$, if there exists a path $J$ in $G - A_1$ with 
        \begin{itemize}
            \item $V(J) \cap V(\mathcal{S}_j^i) \neq \emptyset$ for all $j \in [18t+1]$, and

            \item $V(J) \cap V(Q) = \emptyset$ for any path $Q \in \mathcal{S}_j^i$ spared by $A_1$, for all but at most one $j \in [18t+1]$, 
        \end{itemize}
        then $G$ contains a quarter-integral packing of $t$ even dicycles.
	\end{claim}
    \emph{Proof of \Cref{claim:farjumpgivespacking}:}
    If there exists a $j \in [18t+1]$ and a $Q \in \mathcal{S}_j^i$ that is spared by $A_1$ and $V(Q) \cap V(J) \neq \emptyset$, we set $h \coloneqq j$ and otherwise let $h \coloneqq 0$.
    We now iteratively identify $6t$ disjoint subpaths $J_1, \ldots , J_{6t}$ of $J$ and $6t$ distinct paths $P_1, \ldots , P_{6t}$ from $\mathcal{S}^i$ that are spared by $A_1$, such that each $J_j$ can be extended to form a cross with $P_j$.
    These then let us find a quarter-integral packing of $t$ even dicycles, after some work on disentangling the mess we made.
    To initialise our search, we set $I_0 = [18t+1] \setminus \{ h \}$, let $s_0$ be the tail of $J$, and start with $j = 1$.

    Let $r_1 \in I_{j-1}$ and $u \in V(\mathcal{S}_{r_1}^i)$  be chosen such that $s_{j-1}Ju$ is minimal.
    Further, We choose $r_2 \in I_{j-1} \setminus \{ r_1 \}$ and $w \in V(\mathcal{S}_{r_2}^i)$ such that $uJw$ is minimal and choose $r_3 \in I_{j-1} \setminus \{ r_1, r_2 \}$ and $v \in V(\mathcal{S}_{r_3}^i)$ such that $wJv$ is minimal.
    We set $I_j = I_{j-1} \setminus \{ r_1, r_2, r_3 \}$.
    To find a third subpath, we choose $s_j \in V(\mathcal{S}_{r_3}^i)$ such that $wJs_j$ does not intersect $\bigcup_{k \in I_j} V(\mathcal{S}_k^i)$ and $wJs_j$ is maximal.
    
    We first consider $uJw$.
    Let $I_{R,1}$ be the set of indices $k$ from $\mathcal{I}_{R,1} \coloneqq [(r_1-1)q + 1, r_1 q]$ such that $uJw$ and $Q_k^i \in \mathcal{S}_{r_1, R}^i$ intersect and let $I_{L,1}$ be defined analogously for $\mathcal{S}_{r_1, L}^i$.
    Suppose there exist three indices $k, k' \in I_{R,1}$ and $\ell \in \mathcal{I}_{R,1} \setminus I_{R,1}$ with $k < \ell < k'$ and $Q_\ell^i$ is spared by $A_1$.
    Then there must exist a $V(Q_k^i)$-$V(Q_{k'}^i)$- or a $V(Q_{k'}^i)$-$V(Q_k^i)$-subpath $J' \subseteq uJw$ that is disjoint from $Q_\ell^i$.
    Within $Q_k^i \cup Q_{k'}^i \cup J'$ we can therefore find a path $J_j$ that forms a cross with $Q_\ell^i$, which we can then choose as $P_j$.
    In this case, we therefore proceed to the next round.

    Otherwise we know that $I_{R,1}$ is found in a containment-wise minimal subinterval $I_{R,1}'$ of $\mathcal{I}_{R,1}$ such that for all $\ell \in I_{R,1}$ the path $Q_\ell^i$ is not spared by $A_1$.
    Of course this is also true for the indices in $I_{L,1}$ and thus we find an analogous interval $I_{L,1}'$.
    The same arguments also apply to $wJv$ and $vJs_j$, allowing us to find the sets $I_{R,2}, I_{L,2}, I_{R,3}, I_{L,3}$ and the associated intervals $I_{R,2}', I_{L,2}', I_{R,3}', I_{L,3}'$.
    In particular, if any of the sets is empty, we choose the associated intervals to be empty.

    Let $u'$ be the last vertex of $uJw$ that is contained in $V(\mathcal{S}_{r-1}^i)$, when traversing the path starting from the tail.
    Note that the index of the path of $\mathcal{S}_{r_1}^i$ that $u'$ belongs to is contained either in $I_{R,1}'$ or in $I_{L,1}'$.
    W.l.o.g.\ we suppose that $u' \in I_{R,1}'$ and note that if $w$ is contained in $I_{L,2}'$, then either $u'Jw$ contains a path $J'$ that jumps over a path of $\mathcal{S}_{r_1,L}^i$ that is spared by $A_1$, or $J'$ jumps over a path of $\mathcal{S}_{r_2,R}^i$ that is spared by $A_1$.
    This allows us to define $J_j$ and $P_j$ analogously to our previous argument and proceed to the next round.
    Thus $w$ must be contained in $I_{R,2}'$.

    Suppose that $r_1 < r_3 < r_2$ or $r_2 < r_3 < r_1$ holds.
    Then we can find a jump $J'$ within $wJv$ over a path of $\mathcal{S}_{r_3}^i$ that is not hit by $A_1$ and again find our desired paths.
    For analogous reasons, we can exclude the possibility that $r_2 < r_1 < r_3$ or $r_3 < r_1 < r_2$ holds.
    Thus we either have $r_1 < r_2 < r_3$ or $r_3 < r_2 < r_1$ and w.l.o.g.\ we can assume that we are in the case in which $r_1 < r_2 < r_3$ is true.

    As a consequence, since we know that $w \in I_{R,2}'$, we also know that $I_{L,2}' = \emptyset$.
    Repeating our arguments we can similarly conclude that $v \in I_{R,3}'$ and $I_{L,3}' = \emptyset$.
    Suppose now that there exists an $\ell$ in $\mathcal{I}_{R,2} \coloneqq [(r_2-1)q + 1, r_2 q]$ such that $k < \ell < k'$ for any $k \in I_{R,2}'$ and $k' \in I_{R,1}'$, and $Q_\ell^i$ is spared by $A_1$.
    Then we can find a jump $J'$ over $Q_\ell^i$ and are again done.
    
    However, if no such $\ell$ exists, there must exists some $\ell \in \mathcal{I}_{R,2}$ such that $k'' < \ell < k'''$ for any $k'' \in I_{R,3}'$ and $k''' \in I_{R,2}'$ such that $Q_\ell^i$ is spared by $A_1$, as $I_{R,2}' \neq \mathcal{I}_{R,2}$ due to its minimality and the fact that $|I_{R,2}'| < q$.
    Hence we can again find a jump over $Q_\ell^i$ and finish the round.
    After $6t$ rounds we thus find our desired paths $P_1, \ldots , P_{6t}$ and $J_1, \ldots , J_{6t}$.

    It is easy to observe that any two paths we have found are pairwise disjoint and in particular, for each $j \in [6t]$, the path $J_j$ ends on two distinct paths $Q_{k_j}^i$ and $Q_{\ell_j}^i$, such that $Q_{k_1}^i, Q_{\ell_j}^i, \ldots , Q_{k_{6t}}^i, Q_{\ell_{6t}}^i$ is another set of pairwise disjoint paths.
    Furthermore each $J_j$ is internally disjoint from $Q_{k_j}^i$ and $Q_{\ell_j}^i$.
    We can therefore find a path $J_j'$ in $Q_{k_j}^i \cup Q_{\ell_j}^i \cup J_j$ that forms a cross with $P_j$.

    Note that the paths $J_1', \ldots , J_{6t}'$ may now be a half-integral linkage.
    We want to apply \Cref{thm:halfintegral} at this point to clean this up, but there is a case in which the resulting linkage may no longer form crosses with the paths $J_1, \ldots , J_{6t}$, or even any other paths of $\mathcal{S}^i$ that is spared by $A_1$.
    Thus we must first get rid of potential troublemakers.

    Consider the paths $Q_{k_j}^i, Q_{\ell_j}^i$ for some $j \in [6t]$ and suppose w.l.o.g.\ that $k_j < \ell_j$.
    Let $Q_x^i$ be a path from $Q_{k_1}^i, Q_{\ell_j}^i, \ldots , Q_{k_{6t}}^i, Q_{\ell_{6t}}^i$ such that $\ell_j < x$ and there is no $Q_y^i$ with $\ell_j < y < x$ amongst $Q_{k_1}^i, Q_{\ell_j}^i, \ldots , Q_{k_{6t}}^i, Q_{\ell_{6t}}^i$.
    If there exists no $Q_z^i \in \mathcal{Q}_i$ with $\ell_j < x$ that is spared by $A_1$, then \Cref{thm:halfintegral} may in fact return a path that ends in the endpoints of $Q_{\ell_j}^i$ and $Q_x^i$, but does not jump over anything.
    Of course an analogous problem may occur next to $Q_{k_j}^i$.
    Let $j'$ and $j''$ be the two indices fitting the two problematic paths we may find this way.

    We note that there cannot be a third such path, since we choose all $J_j$ and thus all $Q_{k_j}^i, Q_{\ell_j}^i$ from pairwise disjoint triples of transactions from $\mathcal{S}_1^i, \ldots , \mathcal{S}_q^1$.
    To proceed we greedily search for triples of indices from $[18t+1]$ such that there are less than $t$ paths spared by $A_1$ between the paths from $Q_{k_1}^i, Q_{\ell_j}^i, \ldots , Q_{k_{6t}}^i, Q_{\ell_{6t}}^i$ we consider.
    This at worst partitions our $6t$ selected paths into $2t$ triples, for each of which we drop the two selected paths on the outside.
    As a result our selection shrinks down to at least $2t$ paths $J_1'', \ldots , J_{2t}''$, which now however have the property that each pair of endpoints is separated in $\Psi$ by the endpoints of at least $t$ paths in $\mathcal{S}^i$ that are spared by $A_1$.
    Clearly, via an application of \Cref{thm:halfintegral}, we find $t$ disjoint paths that can be combined with $t$ paths from $\mathcal{S}^i$ that are spared by $A_1$ to apply \Cref{lem:applyshifting}.
    \hfill$\blacksquare$

    We will not always find such a path, but we can use the fact that these paths do not exist to find a different kind of useful path that does not stray very far from the area we want to isolate.

    \begin{claim}\label{claim:isolatetransactions}
        For each $i \in [2t-1]$, at least one of the following is true:
        \begin{itemize}
            \item There exists a quarter-integral packing of $t$ even dicycles in $G$.

            \item There exist no $X_i$-$(Y_i \cup Z_i)$- and no $(Y_i \cup Z_i)$-$X_i$-paths in $G - (A_1 \cup V( \{ R_i , L_i \} ) )$.
            
            \item There exists an $X_i$-$Z_i$- or $Z_i$-$X_i$-path $J$ in $G - (A_1 \cup V( \{ R_i , L_i \} ) )$ and some $j \in [18t+1]$ such that $V(J) \cap V( \mathcal{S}_j^i ) = \emptyset$.
        \end{itemize}
	\end{claim}
    \emph{Proof of \Cref{claim:isolatetransactions}:}
    Suppose the first two statements of the claim do not hold and let $i \in [2t-1]$ be arbitrary.
    Let $J$ be an $X_i$-$Y_i$-path in $G - (A_1 \cup V( \{ R_i , L_i \} ) )$, or an $X_i$-$Z_i$-path such that $V(J) \cap V( \mathcal{S}_j^i ) \neq \emptyset$ for all $j \in [18t+1]$.
    In the first case, we also know that $J \cap V( \mathcal{S}_j^i ) \neq \emptyset$ holds for each $j \in [18t+1]$ thanks to \Cref{claim:movethroughsafety}.
    
    Note that if there exists $j \in [18t+1]$ and a $Q \in \mathcal{S}_j^i$ that is spared by $A_1$ and $V(Q) \cap V(J) \neq \emptyset$, then we could reroute $J$ via $Q$ to find an $X_i$-$Z_i$-path in $G - (A_1 \cup V( \{ R_i , L_i \} ) )$.
    If there exists a $j \in [18t+1]$ with this property, we can choose this $j$ and a $v \in V( \mathcal{S}_j^i )$ such that $Jv$ is as short as possible.
    We note that if there exists a $j \in [18t+1]$ such that $V(Jv) \cap V( \mathcal{S}_j^i ) = \emptyset$ we are done.
    Thus we can assume that this is not the case and redefine $J \coloneqq Jv$, as well as setting $h \coloneqq j$.
    If no such $j$ exists, we leave $J$ as is and set $h = 0$.

    We conclude that for all $j \in [18t+1] \setminus \{ h \}$ we know that if $Q \in \mathcal{S}_j^i$ is spared by $A_1$ then we have $V(J) \cap V(Q) \neq \emptyset$ allowing us to apply \Cref{claim:farjumpgivespacking} and find a quarter-integral packing of $t$ even dicycles, contradicting our initial assumptions.
    For an $Y_i$-$X_i$-path in $G - (A_1 \cup V( \{ R_i , L_i \} ) )$, or a $Z_i$-$X_i$-path $P$ such that $V(P) \cap V( \mathcal{S}_j^i ) \neq \emptyset$ for all $j \in [18t+1]$, we can argue analogously.
	\hfill$\blacksquare$

    If we find a quarter-integral packing of $t$ even dicycles, we are of course happy and finding $t$ different transactions which we can isolate, corresponding to the second option of \Cref{claim:isolatetransactions}, is actually the outcome we want to proceed with.
    Thus we must show that the third statment in \Cref{claim:isolatetransactions} cannot hold for many transactions.

    \begin{claim}\label{claim:isolationorpacking}
        Either there exists a quarter-integral packing of $t$ even dicycles in $G$, or there exists a set $I \subseteq [2t - 1]$ with $|I| \geq t$ such that for each $i \in I$ there exist no $X_i$-$(Y_i \cup Z_i)$- and no $(Y_i \cup Z_i)$-$X_i$-paths in $G - (A_1 \cup V( \{ R_i , L_i \} ) )$.
	\end{claim}
    \emph{Proof of \Cref{claim:isolationorpacking}:}
    As always, we suppose that there does not exist a quarter-integral packing of $t$ even dicycles in $G$.
    Using \Cref{claim:isolatetransactions}, we thus know that for each $i \in [2t-1]$ either there exist no $X_i$-$(Y_i \cup Z_i)$- and no $(Y_i \cup Z_i)$-$X_i$-paths in $G - (A_1 \cup V( \{ R_i , L_i \} ) )$, or there exists an $X_i$-$Z_i$- or a $Z_i$-$X_i$-path $J_i$ in $G - (A_1 \cup V( \{ R_i , L_i \} ) )$ and some $j \in [18t+1]$ such that $V(J_i) \cap V( \mathcal{S}_j^i ) = \emptyset$.
    
    Of course, if there happen to be $t$ indices $i \in [2t-1]$ such that there exist no $X_i$-$(Y_i \cup Z_i)$- and no $(Y_i \cup Z_i)$-$X_i$-paths in $G - (A_1 \cup V( \{ R_i , L_i \} ) )$, we are done.
    Thus we suppose that this is not the case for at least $t$ indices $i \in [2t-1]$.
    W.l.o.g.\ we will assume that $J_1, \ldots , J_t$ exist to make the notation easier.
    We can choose these $J_i$ such that for each $i \in [t]$ there exists at most one $Q \in \mathcal{S}^i$ that has been spared by $A_1$ such that $J_i$ and $Q$ intersect, and in particular, we can choose $J_i$ such that $J_i \cap Q$ is a single subpath of $J_i$ that contains either the tail or the head of $J_i$, if such a $Q$ exists.
    
    Furthermore, we note that for any distinct $i, j \in [t]$, the path $J_i$ is disjoint from $L_j \cup R_j$, as both these paths have been spared by $A_1$ and thus if these paths intersect, $J_i \cup L_j \cup R_j$ contains an $X_i$-$Y_i$- or a $Y_i$-$X_i$-path $J_i'$ and a $k \in [18t+1]$ with $J_i' \cap V(\mathcal{S}_k^i) = \emptyset$ exists, contradicting \Cref{claim:movethroughsafety}.
    
    If $J_1, \ldots , J_t$ are pairwise disjoint, we note that each $J_i$ is also disjoint from $L_i \cup R_i$ and forms a cross with one of these two paths.
    Thus an application of \Cref{lem:applyshifting} yields a quarter-integral packing of $t$ even dicycles in $G$, contradicting our initial assumption.

    Therefore there exist two distinct $i, j \in [t]$ such that $J_i$ and $J_j$ intersect.
    W.l.o.g.\ we suppose that $J_i$ is an $X_i$-$Z_i$-path.
    Let $u \in V(J_i) \cap V(J_j)$ be the unique vertex minimising the length of $J_iu$.
    Due to our choice of $u$, we know that $J \coloneqq J_iuJ_j$ is an $X_i$-$Y_i$-path, since $X_j, Z_j \subseteq Y_i$.
    As we noted prior, we also know that both $J_i$ and $J_j$ are disjoint from $L_i$, $R_i$, $L_j$, and $R_j$ and in particular, both are spared by $A_1$.
    Thus $J$ is an $X_i$-$Y_i$-path in $G - (A_1 \cup V( R_i \cup L_i ) )$.
    \Cref{claim:movethroughsafety} therefore tells us that we have $V(J) \cap V( \mathcal{S}_k^i )$ for all $k \in [18t+1]$.

    Suppose that $J$ intersects some path $Q \in \mathcal{S}_k^i$ for some $k \in [18t+1]$ that has been spared by $A_1$.
    Then by our assumption, we know that $J_j$ intersects $Q$, as $J_i$ could only intersect $Q$ in a subpath containing its head.
    Thus $J_j \cup Q$ contains an $X_j$-$Y_j$- or a $Y_j$-$X_j$-path $J'$, since $X_i, Z_i \subseteq Y_j$, in $G - (A_1 \cup V( \{ R_j , L_j \} ) )$ and, in contradiction to \Cref{claim:movethroughsafety}, we have $V(J') \cap V(\mathcal{S}_k^j) = \emptyset$ for some $k \in [18t+1]$.

    We conclude that $J$ does not intersect any path $Q \in \mathcal{S}_k^i$ that has been spared by $A_1$ for any $k \in [18t+1]$.
    This means we can use \Cref{claim:farjumpgivespacking} to find a quarter-integral packing of $t$ even dicycles in $G$, contradicting our initial assumption and completing the proof of the claim.
	\hfill$\blacksquare$

    Thus, since we supposed that $G$ does not contain a quarter-integral packing of $t$ even dicycles, we can find a set $I \subseteq [2t-1]$ of size at least $t$ such that there exist no $X_i$-$(Y_i \cup Z_i)$- and no $(Y_i \cup Z_i)$-$X_i$-paths in $G - (A_1 \cup V( \{ R_i , L_i \} ) )$ for all $i \in I$.
    For the sake of simplicity, we assume that $[t] \subseteq I$ and proceed to only work with these $t$ values.
    
    For any $i \in [t]$ note that $\mathcal{Q}_i \setminus \mathcal{S}^i$ contains $p' - 36tq = p + 2q$ paths of which at least $p + q$ are spared by $A_1$, including $R_i$ and $L_i$.
    Thus we can find a transaction $\mathcal{P}_i \subseteq \mathcal{Q}_i \setminus \mathcal{S}^i$ of order $p'' = p + 4(t + 1)$ consisting of paths spared by $A_1$, such that $R_i, L_i \in \mathcal{P}_i$.
    For each $i \in [t]$ we let $\mathcal{P}_i = \{ P_1^i, \ldots , P_{p''}^i \}$ such that for all $j \in [2, p'' -1]$ the endpoints of $P_j^i$ separate the endpoints of $P_{j-1}^i$ and the endpoints of $P_{j+1}^i$ on $\Psi$.
    
    We let $G_i$ be the union of the strip of $\mathcal{P}_i$ under $A_1$ in $G$, the graph $R_i \cup L_i$, and all $X_i$-$V(R_i \cup L_i)$-, as well as all $V(R_i \cup L_i)$-$X_i$-paths in the union of the strong components of $G - A_1$ that contain vertices of $V(\Psi)$, for all $i \in [t]$.
    Notably, for all distinct $i, j \in [t]$, the graphs $G_i$ and $R_j \cup L_j$ are disjoint, as we could otherwise find an $X_i$-$(Y_i \cup Z_i)$- or a $(Y_i \cup Z_i)$-$X_i$-path in $G - (A_1 \cup V( \{ R_i , L_i \} ) )$.

    \begin{claim}\label{claim:disjointstrips}
        For all distinct $i, j \in [t]$ the graphs $G_i$ and $G_j$ are disjoint.
	\end{claim}
    \emph{Proof of \Cref{claim:disjointstrips}:}
    Suppose this is not true and there exist distinct $i, j \in [t]$ such that there exists some $v \in V(G_i) \cap V(G_j) \neq \emptyset$.
    Thus there also exists a $v$-$(V(\mathcal{P}_i) \cup X_i)$- or a $(V(\mathcal{P}_i) \cup X_i)$-$v$-path $J_i$ in $G_i$ and a $v$-$(V(\mathcal{P}_j) \cup X_j)$- or a $(V(\mathcal{P}_j) \cup X_j)$-$v$-path $J_j$ in $G_j$.
    
    By definition $v$ is part of some strong component of $G - A_1$ that also contains a vertex $u \in V(\Psi)$.
    Let $T$ be a $v$-$u$-path in $G - A_1$.
    We note that, since $L_i$, $R_i$, $L_j$, and $R_j$ have been spared, we can choose $T$ such that it only intersects one of these four paths.
    W.l.o.g.\ we assume that $T$ only intersects $L_j \cup R_j$ and furthermore, we will assume that $v$ is the head of $J_i$.
    (If $v$ is the tail, we can choose $T$ to be a $u$-$v$-path instead.)
    Given our choices, we observe that $J_i \cup T$ contains a $(V(\mathcal{P}_i) \cup X_i)$-$u$-path in $G - A_1$ that avoids both $R_i$ and $L_i$.
    If $J_i$ has its tail outside of $V( R_i \cup L_i )$ this immediately leads to a contradiction.
    Thus our case narrows and we can additionally assume that $J_i$ has its tail in $V( R_i \cup L_i )$.

    Should $v$ be the tail of $J_j$ and should the head of $J_j$ lie outside of $(R_j \cup L_j)$, we note that $J_i \cup J_j$ can be extended to an $Y_j$-$X_j$-path in $G - (A_1 \cup V( \{ R_j , L_j \} ) )$, which leads to another contradiction.
    Thus if $v$ is the tail of $J_j$, the head of $J_j$ lies in $(R_j \cup L_j)$.
    As a consequence, by construction of $G_i$ and $G_j$, there must therefore be an $X_j$-$v$-path in $G_j$ and a $v$-$X_i$-path in $G_i$.
    Hence there also exists an $X_j$-$X_i$-path in $G_j$ that does not intersect $R_j \cup L_j$, leading to a contradiction, since $X_i \subseteq Y_j \cup Z_j$.

    Thus $v$ must be the head of $J_j$ and we can immediately narrow down this case as we did in the beginning for $J_i$, leading us to the scenario in which $J_i$ is a $V(R_i \cup L_i)$-$v$-path in $G_i$ and $J_j$ is a $V(R_j \cup L_j)$-$v$-path in $G_j$.
    As a consequence there exists a $v$-$X_j$-path $P$ in $G_j$ that is disjoint from $R_j \cup L_j$ and therefore $J_i \cup P$ is a $(Y_j \cup Z_j)$-$X_j$-path in $G - (A_1 \cup V( \{ R_j , L_j \} ) )$.
    This completes the proof of our claim.
	\hfill$\blacksquare$

    We now let $H_i$ be the strip of $\mathcal{P}_i \setminus \{ P_1^i, P_2^i, P_{p'' -1}^i, P_{p''}^i \}$ under $\emptyset$ in $G_i$ for each $i \in [t]$, and let $\mathfrak{R} = ( H^* , \mathcal{C}^*, \mathcal{E}^*, \mathfrak{B}^*, \mathfrak{D}^* )$ be the $t$-retreat of the $\nicefrac{\theta}{2}$-rim of either $\mathfrak{M}$, if $\mathfrak{M}$ is a circle outline, or $\mathfrak{N}$, if $\mathfrak{M}$ is a diamond outline.
    Using \Cref{lem:dicyclesthroughatransaction} or \Cref{lem:dicyclesthroughatransactionrefined}, we can now either find some $i \in [t]$, such that $H^* \cup H_i$ is odd and $H^* \cup G_i$ is non-even, or find a quarter-integral packing of $t$ even dicycles.

    Thus we may suppose that $H^* \cup H_i$ is odd and $H^* \cup G_i$ is non-even for some $i \in [t]$.
    W.l.o.g.\ we assume that $i = 1$ to simplify notation.
    Let $X$ be the union of the two maximal segments of $\Psi$ that contain both heads, respectively both tails, of $P_{t+2}^1$ and $P_{p'' - (t + 1)}^1$, but no endpoints of $P_{t+1}^1$ or $P_{p'' - t}^1$, and let $Y$ be $V(\Psi)$ minus the union of the two minimal segments containing all heads, respectively all tails, of $\mathcal{P}_1 \setminus \{ R_i, L_i \}$.
    We now find another apex set that disrupts jumps into $H_i$.

    \begin{claim}\label{claim:newapex}
        Either there exists a quarter-integral packing of $t$ even dicycles in $G$, or there exists a set of vertices $A_2$ such that $G - ( A_1 \cup A_2 \cup V(\{ P_2^1, \ldots, P_{t+1}^1, P_{p'' - t}^1, \ldots , P_{p'' - 1}^1 \}) )$ contains no $X$-$Y$- and no $Y$-$X$-paths with $|A_2| \leq 2(t-1)$.
	\end{claim}
    \emph{Proof of \Cref{claim:newapex}:}
    The proof of this claim is entirely analogous to the proof of \Cref{claim:findsafetyapex}.
	\hfill$\blacksquare$

    Under the assumption that there does not exist a quarter-integral packing of $t$ even dicycles in $G$, we let $A = A_1 \cup A_2$ and note that $|A| = 72t^3 - 104t^2 + 32t = 72t^3 - 104t^2 + 30t + 2 + 2(t-1)$.
    Let $\mathcal{P}_1'$ be the transaction that contains exactly those paths from $\mathcal{P}_1 \setminus \{ P_1^1, \ldots, P_{t+1}^1, P_{p'' - t}^1, \ldots , P_{p''}^1 \}$ that were spared by $A_2$.
    Note that since $|\mathcal{P}_1| = p + 4(t + 1)$ and $A_2 \leq 2(t-1)$, we know that $\mathcal{P}_1'$ contains at least $p$ paths.
    Furthermore, all paths in $\mathcal{P}_1'$ were spared by $A$.
    Let $R \in \mathcal{P}_1'$ be the path whose endpoints separate the endpoints of $R_1$ from the endpoints of the remaining paths in $\mathcal{P}_1'$ on $\Psi$ and let $L$ be defined analogously for $L_1$.

    \begin{claim}\label{claim:stripisodd}
        The strip $S$ of $\mathcal{P}_1'$ under $A$ in $G$ is odd.
	\end{claim}
    \emph{Proof of \Cref{claim:stripisodd}:}
    Suppose there exists some even dicycle $C \subseteq S$ and let $v \in V(C)$ be arbitrary.
    Since $v \in V(S)$, there exists a $v$-$V(\mathcal{P}_1')$- or a $V(\mathcal{P}_1')$-$v$-path in $G - A$ that avoids $V(R \cup L)$ and is contained in the same strong component as some vertex of $V(\Psi)$ in $G - A$.
    W.l.o.g.\ we assume that there exists a $v$-$V(\mathcal{P}_1')$-path $J$.
    If $J \subseteq G_1$, then by definition we must also have $v \in V(S)$.
    Otherwise we know that $J$ intersects $R_1 \cup L_1$, since in this case some part of $J$ may lie outside of $G_1$.
    We let $u \in V(G_1)$ be the last vertex from $R_1 \cup L_1$ that is seen on $J$ when traversing the path starting from the tail.
    However, $uJ$ is contained in $G_1$ by definition, since it can be extended to a $V(R_1 \cup L_1)$-$X_1$-path and since $J$ avoids $V(R \cup L)$, this lets us find a cross in $G_1$, contradicting the fact that it is non-even.
	\hfill$\blacksquare$

    We can now apply \Cref{lem:flowalteredoutline} to $\mathcal{P}_1'$ and $A$ to find the desired new outline.
    After possibly truncating the transaction $\mathcal{P}_1'$ to fit the new society, we have thus completed the proof.
\end{proof}

Since the proof makes these numbers explicit, we can note that 
\[ \OddTransactionApex{t} = 72t^3 - 104t^2 + 32t \text{ and } \OddTransactionOrder{t,p} = (2 t - 1) (p + (36 t + 2) (72 t^3 - 104 t^2 + 38 t + 2)) + 2 t - 1 \]
are sufficient to satisfy the statement of \Cref{thm:oddtransaction}.

\section{Building new outlines}\label{sec:buildoutline}

\subsection{Subdividing a circle outline}

Now that we can find a transaction that neatly fits the existing decomposition of a given outline, provided that the maelstrom described by the outline does not have bounded depth, we need to show how this transaction lets us find new outlines.
To simplify this setting somewhat, we show that we can sacrifice a small part of the transaction to find a transaction whose paths only move against the direction of our infrastructure in the area containing $H$.
As we will demonstrate, this is step is not strictly necessary, but will allow us to focus on the actually interesting parts of the proof of the main results for this section.
We only use this tool for circle outlines, and thus will only define it for them.
For (refined) diamond outlines we will rely on the fact that the reader can get the gist of the ideas from the simplified circle outline case and concentrate on explaining the technical challenges of the (refined) diamond outline case.

\begin{definition}
Let $\theta , t , p$ be integers such that $\theta \geq 2$ is even, and let $D$ be a digraph with an odd decomposition $\delta$.
Moreover, let $m$ be a maelstrom of $\delta$ together with a circle $\theta$-outline $\mathfrak{M} = ( H , \mathcal{C}, \mathfrak{E}, \mathfrak{V} )$ of degree $t$, with $(H_\mathfrak{M},\Omega_\mathfrak{M})$ being the $\mathfrak{M}$-society.
If $\mathcal{P} = \{ P_1, \ldots , P_p \}$ is an odd transaction on the $\mathfrak{M}$-society, together with an apex set $A \subseteq V(H_\mathfrak{M})$ such that $\mathcal{P}$ and $A$ are $\delta$-compatible, we say that $\mathcal{P}$ is \emph{caught in the flow of $\mathfrak{M}$} if there exists a $P \in \mathcal{P}$ and
\begin{itemize}
    \item a $C \in \mathcal{C}$ such that $P$ has a $C$-subpath $Q \subseteq H$ for which $Q \cup C$ contains a dicycle $C'$ with $Q \subseteq C'$, whose trace separates $m$ and elements of $V(\Omega_\mathfrak{M})$, and that is not the $m$-tight dicycle of $Q \cup C$, or

    \item an $E \in \bigcup_{i \in [t]} \mathcal{E}_i$ such that $P$ has an $E$-subpath $Q \subseteq H$ for which $Q \cup E \cup C_\theta$ contains a dicycle $C'$ with $Q \subseteq C'$, whose trace separates $m$ and the elements of $V(\Omega_\mathfrak{M})$, and that is not the $m$-tight dicycle of $Q \cup E \cup C_\theta$.
\end{itemize}
If $\mathcal{P}$ is not caught in the flow of $\mathfrak{M}$, we say that $\mathcal{P}$ \emph{flows against} $\mathfrak{M}$.
\end{definition}

We now show that we can always find a transaction that is not caught in the flow of the circle outline, if we sacrifice $\theta - 1$ paths from the initial transaction.
Analogous statements and definitions can also be provided for (refined) diamond outlines, though we will not need them.

\begin{lemma}\label{lem:caughtupintheoutline}
    Let $\theta$ and $p$ be integers such that $\theta$ is even and $p \geq \theta \geq 2$, and let $D$ be a digraph with an odd decomposition $\delta$.
    Moreover, let $m$ be a maelstrom of $\delta$ together with a circle $\theta$-outline $\mathfrak{M} = ( H , \mathcal{C}, \mathfrak{E}, \mathfrak{V} )$, with $(H_\mathfrak{M},\Omega_\mathfrak{M})$ being the $\mathfrak{M}$-society, and let $\mathcal{P} = \{ P_1, \ldots , P_p \}$ be an odd transaction on the $\mathfrak{M}$-society, together with a pleasant apex set $A \subseteq V(H_\mathfrak{M})$, such that $\mathcal{P}$ and $A$ are $\delta$-compatible.
    
    Then there exists a $\theta$-outline $\mathfrak{M}'$ for $m$, such that $(H_\mathfrak{M},\Omega_\mathfrak{M})$ is the $\mathfrak{M}'$-society, and the odd transaction $\mathcal{P}' = \{ P_\theta, \ldots , P_p \}$ on the $\mathfrak{M}'$-society flows against $\mathfrak{M}'$ and is $\delta$-compatible together with $A$.
\end{lemma}
\begin{proof}
    Let $P \in \mathcal{P}$ and let $C \in \mathcal{C}$ such that $P$ has a $C$-subpath $Q \subseteq H$ for which $Q \cup C$ contains a dicycle $C'$ with $Q \subseteq C'$ whose trace separates $m$ and elements of $V(\Omega_\mathfrak{M})$.
    Since $A$ is pleasant and $\mathcal{P}$ together with $A$ is $\delta$-compatible, we know that therefore $P_1$, or $P_p$ depending on the way the elements of $\mathcal{P}$ are labelled, behaves in the same way towards $C_1$.
    We may therefore assume that $P \subseteq P_1$.
    Let $d$ be the closure of the disk defined by the diamond in $P \cup C_1$ that neither contains $m$ nor $V(\Omega_\mathfrak{M})$.
    Further, let $H_d$ be the subgraph of $H$ drawn on $d$ by $\Gamma$, and let $s \in [\theta -1]$ be the largest value for which $P_s \cap H_d \neq \emptyset$.
    Note that several subpaths of $P_1$ may behave just like $P$ and thus several disjoint regions like $H_d$ may exist.
    Our suggested modifications can be performed on all such regions simultaneously.
    Thus we will concentrate on only $H_d$, assuming the other regions are dealt with in analogous fashion using the same paths from $\mathcal{P}$.

    Our plan is now to adjust each cycle in $\mathcal{C}$ that is intersected by some path in $\mathcal{P}$ on $H_d$, such that the cycle ends up a bit more close to $C_\theta$, which itself cannot be intersect since it constitutes part of the society of $\mathfrak{M}$.
    For $P_1$ let $j$ be the largest value in $[\theta -1]$ such that $(P_1 \cap C_j) \cap H_d$ is non-empty.
    There must exist a set of $C_j$-subpaths $\mathcal{R}$ of $P$ such that $C_j \cup \bigcup_{R \in \mathcal{R}} R$ contains a dicycle $C_j'$ whose trace separates not only $m$ and elements of $V(\Omega_\mathfrak{M})$, but also separates all remaining parts of $P$ and the elements of $V(\Omega_\mathfrak{M})$.
    We can now replace $C_j$ with $C_j'$ in $\mathcal{C}$ and since $j$ was largest value in $[\theta -1]$ for which $P_1$ intersects $C_j$, at least on $H_d$, we know that $\mathcal{C}$ remains a homogeneous set of pairwise disjoint dicycles.
    This argument can be repeated for all paths $P_2, \ldots , P_s$, as long as they still intersect elements of $\mathcal{C}$.
    This removes all interactions between the paths $P_{s+1}, \ldots , P_p$ and the elements of $\mathcal{C}$ drawn on $H_d$, since $\delta$ is an odd decomposition.
    To ensure that no member of $\mathcal{C}$ intersects a path in $\mathcal{P}'$, we therefore at most need to sacrifice the paths $P_1, \ldots , P_{\theta -1}$.

    This argument can analogously be extended for any situation in which there exists a path $P \in \mathcal{P}$ and a path $E \in \mathcal{E}_j \in \mathfrak{E}$ such that $P$ has an $E$-subpath $Q \subseteq H$ for which $Q \cup E \cup C_\theta$ contains a dicycle $C'$ with $Q \subseteq C'$ whose trace separates $m$ and elements of $V(\Omega_\mathfrak{M})$.
\end{proof}

\begin{remark}
    In the main results of this section we consider how the odd transaction splits the maelstrom into several smaller maelstroms.
    As part of this we define new outlines for these smaller maelstroms based on the old one.
    In our proofs we will pretend that all parts of the infrastructure of the old outline are relevant for the infrastructure of each new outline, but this is very unlikely to be the case.
    Ignoring this possibility does not make our proofs easier, but avoids certain technicalities when defining the new outlines.
    Therefore we simply ignore this and always construct new outlines using the entirety of the old infrastructure.
    Irrelevant transactions from $\mathfrak{E}$ and $\mathfrak{B}$ can be dropped after the outline is built.
\end{remark}

Before we can move on to proving the first main result, we define how an odd transaction splits a maelstrom into several disks.

\begin{definition}\label{def:maelstromoutlinenewdisks}
Let $\theta, t, s, g, p$ be integers such that $\theta$ is even, and let $D$ be a digraph with an odd decomposition $\delta$.
Moreover, let $m$ be a maelstrom of $\delta$ together with a $\theta$-outline $\mathfrak{M} = ( H , \mathcal{C}, \mathfrak{E}, \mathfrak{V} )$ with degree $t$ and, if $\mathfrak{M}$ is a diamond outline, a refined diamond $\theta$-outline $\mathfrak{N} = ( H^*, \mathcal{C}, \mathfrak{E}, \mathfrak{B}, \mathfrak{V}^*, \mathfrak{D} )$ with surplus degree $s$ and roughness $g$ for some $n \subseteq m$, with $(H_\mathfrak{M},\Omega_\mathfrak{M})$ being the $\mathfrak{M}$-society, and let $\Delta$ be the disk corresponding to the trace of the cycle in $H$, respectively $H^*$, that defines the $\mathfrak{M}$-, respectively the $\mathfrak{N}$-society.
Let $\mathcal{P} = \{ P_1, \ldots , P_p \}$ be an odd transaction on the $\mathfrak{M}$-, or respectively the $\mathfrak{N}$-society, together with an apex set $A \subseteq V(H_\mathfrak{M})$ such that $\mathcal{P}$ and $A$ are $\delta$-compatible, with $\delta'$ being the $\mathcal{P}$-$A$-expansion of $\delta$.
We define the following three sets:
    \begin{itemize}
        \item Let $\Delta_1$, respectively $\Delta_p$, be the closure of the disk in $\Delta - \Trace_{\delta'}(P_1)$, respectively $\Delta - \Trace_{\delta'}(P_p)$, that does not contain the drawing of any vertex of a path in $\mathcal{P}$.
        The set $\mathcal{D}(\mathcal{P}, m)$ consists of the closures of the disks in $(m \cap \Delta_1) \cup (m \cap \Delta_p)$ and, if $\mathfrak{N}$ exists, we define $\mathcal{D}(\mathcal{P}, n)$ analogously.

        \item For each $i \in [t]$ and each $k \in [\theta - 1]$, we define $d^{\mathcal{E}_i}_k$ as the closure of the disk in $\Delta - \Trace_{\delta'}(E^i_k)$ that does not contain $\Trace_{\delta'}(E^i_{k+1})$.
        Analogously, for each $k \in [\theta - 1]$, we define $d^\mathcal{C}_k$, and respectively $d^\mathcal{P}_k$, as the closure of the disk in $\Delta - \Trace_{\delta'}(C_k)$, respectively $\Delta - \Trace_{\delta'}(P_k)$, that does not contain $\Trace_{\delta'}(C_{k+1})$, respectively $\Trace_{\delta'}(P_{k+1})$.
        For any disk $d \subseteq \Delta$ and $k \in [\theta - 1]$, we let $\mathcal{D}(k, d)$ be the set of the closures of the disks in $d \cap \bigcap_{i=1}^t d^{\mathcal{E}^i}_k \cap d^\mathcal{C}_k \cap d^\mathcal{P}_k$, if we are working with a circle outline.

        \item For each $i \in [s]$, and each $k \in [\theta - 1]$, we define $d^{\mathcal{B}_i}_{k(g+1)}$ as the closure of the disk in $\Delta - \Trace_{\delta'}(B^i_{k(g+1)})$ that does not contain $\Trace_{\delta'}(B^i_{k(g+1)+1})$.
        For any disk $d \subseteq \Delta$ and $k \in [\theta - 1]$, we let $\mathcal{D}(k, d)$ be the set of the closures of the disks in $d \cap \bigcap_{i=1}^t d^{\mathcal{E}^i}_k \cap d^\mathcal{C}_k \cap d^\mathcal{P}_{k(g+1)} \cap \bigcap_{j=1}^{s} d^{\mathcal{B}_j}_{k(g+1)}$, if we are working with a (refined) diamond outline.
    \end{itemize}
\end{definition}

Additionally, we will need a small lemma that tells us that we can integrate a maelstrom that is already odd into an existing odd decomposition without issue.

\begin{lemma}\label{lem:oddmaelstromdecomposition}
    Let $D$ be an odd digraph with an odd decomposition $\delta$, such that $\delta$ uses a single maelstrom $m$.
    Then $D$ has an odd, maelstrom-free decomposition $\delta'$, such that $\delta$ and $\delta'$ agree on the subgraph of $D$ drawn outside of the interior of $m$.
\end{lemma}
\begin{proof}
    We begin by contracting all directed tight cuts in $D$, resulting in the dibrace $D'$.
    If we find an odd, maelstrom-free decomposition of $D'$, we may expand all of these directed tight cuts and simply surround them with some appropriate disks to find an odd, maelstrom-free decomposition of $D$.
    Let $C'$ be the (not necessarily directed) cycle in $D'$ whose trace bounds $m$.
    Note that according to \textbf{\textsf{PD2}} in \Cref{def:planedecomposition}, $C'$ separates the part of $D$ decomposed in the interior of $m$ from the remainder of the graph.
    
    Let $D' = \DirM{B}{M}$ and note that $B$ is a brace, since $D'$ is a dibrace.
    In particular, $B$ is also Pfaffian, since $D'$ is odd and therefore non-even.
    Note that within $B$ there exists an undirected cycle $C$ in the split of $C'$ that separates the split of the interior of $m$ from the remainder of $B$.
    Let $B'$ be the subgraph of the split of $\sigma(m)$, together with $C$, that is separated by $C$ from the rest of $B$.
    Note that according to \Cref{thm:pfaff4cyc}, we may construct the entirety of $B$ from planar braces via 4-cycle sums.
    This in particular allows us to undo the 4-cycle sums that are attached to the cycles of length 4 in $B'$.
    We may therefore reduce $B'$ to a planar graph $G$ that has a planar embedding in which $C$ is a face.

    Using this planar embedding, we can draw disks around each cycle of length 4 in $G$ used for a 4-cycle sum and attach the parts of $B'$ we have just removed.
    During this process, we may drawn these parts of the graph onto the new disks in whichever way we please.
    We can now attach this decomposition of $B'$ to the decomposition of the remainder of $B$ at $C$, yielding what is essentially a planar embedding of $B$, with the exception of some disks whose trace contains only the vertices of a cycle of length four in $B$.
    According to \Cref{lem:smallcycsumandtrisum}, we can translate this embedding into an odd decomposition of $D'$, since the 4-cycle sums on these disks turn into small cycle sums.
    We may then expand the directed tight cuts again and add the appropriate disks back into the decomposition to find an odd decomposition $\delta'$ of $D$, as argued in the beginning.
    Since the only part of $B$ whose embedding we had to change was $B'$, we conclude that $\delta$ and $\delta'$ agree on the subgraph of $D$ drawn outside of the interior of $m$.
\end{proof}

Since we have proven \Cref{lem:caughtupintheoutline}, we can assume that we are working with an odd transaction flowing against our outline in the next statement.

\begin{lemma}\label{lem:buildnewmaelstromscircle}
    Let $\theta , p, \ell$ be positive integers such that $\theta$ is a multiple of $2^\ell \omega $ and $p \geq 2\theta$, and let $D$ be a digraph with an odd decomposition $\delta$.
    Moreover, let $m$ be a maelstrom of $\delta$, together with a circle $\theta$-outline $\mathfrak{M}  = ( H , \mathcal{C}, \mathfrak{E}, \mathfrak{V} )$, with $(H_\mathfrak{M},\Omega_\mathfrak{M})$ being the $\mathfrak{M}$-society, and let $\mathcal{P} = \{ P_1, \ldots , P_p \}$ be an odd transaction on the $\mathfrak{M}$-society, flowing against $\mathfrak{M}$, together with a pleasant apex set $A \subseteq V(H_\mathfrak{M})$ such that $\mathcal{P}$ and $A$ are $\delta$-compatible.
    Let $\delta'$ be the $\mathcal{P}$-$A$-expansion of $\delta$, let $H_\mathcal{P} = (H \cup S) - A$, where $S$ is the strip of $\mathcal{P}$ under $A$.

    Either $H_\mathfrak{M}$ is odd, there exists an integral packing of $\ell$ even dicycles in $H_\mathfrak{M}$, or there exists an $h < \ell$ and pairwise disjoint disks $d_1, \ldots , d_h \subseteq m$ such that for each $i \in [h]$ there exists a $\theta_i$-outline $\mathfrak{M}_i = ( H_i , \mathcal{C}_i, \mathfrak{E}_i, \mathfrak{V}_i )$ for $d_i' \subseteq d_i$, with $\theta_i \geq \nicefrac{\theta}{2^{h-1}}$, $H_i$ is a subgraph of $H_\mathcal{P}$, the graphs $H_j$ and $H_{j'}$ are disjoint for any distinct $j,j' \in [h]$, the graph $\sigma_{\delta'}(d_i') \cup H_\mathcal{P}$ contains an even dicycle, $d_i$ is bounded by the trace of the cycle corresponding to the $\mathfrak{M}_i$-society $(H_{\mathfrak{M}_i}, \Omega_{\mathfrak{M}_i})$ in $\delta'$, and $H_{\mathfrak{M}_i}$ is a proper subgraph of $H_\mathfrak{M}$.
\end{lemma}
\begin{proof}
    Let $t$ be the degree of $\mathfrak{M}$.
    Suppose that $H_\mathfrak{M}$ contains an even dicycle and that there does not exist a integral packing of $\ell$ even dicycles in $H_\mathfrak{M}$.
    Let $Q$ be the dicycle corresponding to $\Omega_\mathfrak{M}$.

    We begin by establishing in which way we assume $\mathcal{P}$ to be labelled.
    W.l.o.g.\ we assume that $P_1$ is the unique path in $\mathcal{P}$ such that the trace of the unique diamond in $Q \cup P_1$ bounds a disk that contains the traces of all paths in $\mathcal{P} \setminus \{ P_1 \}$, but does not contain the entire trace of $Q$.

    We set $\mathcal{P}' = \{ P_{\nicefrac{\theta}{2}}, \ldots , P_{p - \nicefrac{\theta}{2} + 1} \}$.
    Consider the elements $e_1, \ldots , e_h$ of $\mathcal{D}(\mathcal{P}', m)$, labelled such that, if $e_i$ has a vertex of $V(P_{\nicefrac{\theta}{2}})$ on its boundary that occurs before any vertex of $V(P_{\nicefrac{\theta}{2}})$ on the boundary of $e_j$, then $i < j$ for all $i,j \in [h]$.
    This of course only considers the side of the transaction that is closer to $P_1$ and it is possible that some $e_i$ with $i \in [h]$ does not contain such a vertex.
    However, we claim that, should such an $i \in [h]$ exist, then it is unique and we can therefore let $e_1, \ldots , e_h$ be labelled such that $i = h$, if it exists.

    To prove this claim, consider the dicycles $Q_t$ and $O_t$ from $\mathfrak{M}$ and their associated $m$-disks $\Delta$ and $d$.
    Furthermore, consider the disk $d'$ bounded by the trace of the unique diamond in $P_{\nicefrac{\theta}{2}} \cup Q_t$ that does not contain the entirety of the trace of $Q_t$.
    Clearly we have $e_i \subseteq d' \cap d$.
    Suppose there exists some $j \in [h]$ with $i \neq j$ such that $e_j$ does not contain a vertex of $P_1$ either.
    Then we must also have $e_j \subseteq d' \cap d$ and in particular, $e_i$ and $e_j$ must be found in distinct components of $d' \cap d$.
    However, for these two disks to be disjoint, there must exist some $O_t$-subpath of $P_{\nicefrac{\theta}{2}}$ whose trace is drawn outside of the interior of $d$.
    This subpath is therefore a witness for the fact that $\mathcal{P}$ is caught in the flow of $\mathfrak{M}$, a contradiction to our assumptions.

    Suppose that $\sigma(e_h)$ does not contain a vertex of $V(P_{\nicefrac{\theta}{2}})$, then we note that $\{ P_{\nicefrac{p}{2} - 1} , \ldots , P_p \}$ is a transaction of order at least $\theta$ on the $\mathfrak{M}$-society and, since this transaction goes against the flow of $\mathfrak{M}$, we may in fact easily build a diamond $\theta$-outline for a disk encompassing $e_h$.
    Thus if $\sigma(e_h) \cup H_\mathcal{P}$ is even, we can find any desired diamond outline using $\theta$ or less dicycles as their base.
    We may therefore concentrate our attention on those $e_i$ for which $\sigma(e_i)$ contains a vertex of $V(P_{\nicefrac{\theta}{2}})$.

    Suppose that there exists some $i \in [h]$, such that $H_\mathcal{P} \cup \sigma(e_1)$ is odd.
    Then we can refine $\delta'$, $\Gamma$, and $H_\mathcal{P}$ as follows.
    Let $H_\mathcal{P} := H_\mathcal{P} \cup \sigma(e_1)$ and note that the restriction of $\delta'$ to $H_\mathcal{P}$ by definition only contains maelstroms on $e_1$.
    \Cref{lem:oddmaelstromdecomposition} tells us that we can find a maelstrom-free, pure odd decomposition of $H_\mathcal{P}$ that agrees with $\delta$ on the decomposition of $H_\mathcal{P} - ( \sigma(e_1 - \Boundary{e_1} ) $, though this may require modifying $\Gamma$ on $e_1$.
    This allows us to refine both $\delta'$ and $\Gamma$, such that $\delta'$ restricted to $H_\mathcal{P}$ is a maelstrom-free, pure odd decomposition.

    At this point, we note that for at least one $i \in [h]$ the graph $H_\mathcal{P} \cup \sigma(e_i)$ must contain an even dicycle, since $H_\mathfrak{M}$ is not odd.
    We iterate over $e_1, \ldots , e_h$ and refine $\delta'$, $\Gamma$, and $H_\mathcal{P}$ as just laid out, though we may skip $e_h$, if $\sigma(e_h)$ does not contain a vertex of $V(P_{\nicefrac{\theta}{2}})$.
    Let $I \subseteq [h]$ be the maximal set such that $\sigma(e_i) \cup H_\mathcal{P}$ is even, for all $i \in I$.

    \textbf{Case 1:} Suppose $|I| = 1$.
    We are already done if $h \in I$ and $\sigma(e_h)$ does not contain a vertex of $V(P_{\nicefrac{\theta}{2}})$.
    Otherwise, let $i \in I$.
    We note that since $e_i \in \mathcal{D}(\mathcal{P}', m)$, there exists a dicycle $C_{\nicefrac{\theta}{2}}'$ in $C_{\nicefrac{\theta}{2}} \cup P_{\nicefrac{\theta}{2}} \cup \bigcup_{j=1}^t E^j_{\nicefrac{\theta}{2}}$ whose trace bounds $e_i$.
    Thus we can easily satisfy the requirements of \Cref{def:outline} and therefore our statement, by providing a circle $\theta$-outline for $e_i$ via the outline $( H_\mathcal{P}, \mathcal{C}, \mathfrak{E} \cup \{ \mathcal{E}_{t+1} \}, \mathfrak{V}' )$, where $\mathcal{E}_{t+1} = \{ P_1, \ldots , P_\theta \}$ and $\mathfrak{V}'$ is defined through $H_\mathcal{P}$, $\delta'$, and $e_i$.

    \textbf{Case 2:} Suppose instead that $|I| \geq \ell \geq 2$.
    Without loss of generality, let $e_1, \ldots , e_\ell \in I$ to ease notation.
    We note that by definition $e_i$ and $e_j$ are disjoint for all distinct $i,j \in [\ell]$.
    Since $e_1, \ldots , e_\ell \in \mathcal{D}(\mathcal{P}', m)$, for each $i \in [\ell]$ there exists a dicycle $R_i$ in $C_{\nicefrac{\theta}{2}} \cup P_{\nicefrac{\theta}{2}} \cup \bigcup_{j=1}^t E^j_{\nicefrac{\theta}{2}}$ whose trace bounds $e_i$ and furthermore there exists an even dicycle $Z_i \subseteq H_\mathcal{P} \cup \sigma(e_i)$, such that $V(\sigma(e_i)) \cap V(Z_i) \neq \emptyset$.
    Note that therefore, if $Z_i$ and $R_i$ do not intersect, $Z_i$ is found entirely in $\sigma(e_i)$.
    Accordingly, we also have $V(R_i) \cap V(R_j) = \emptyset$ for all distinct $i,j \in [\ell]$.

    For any $i \in [\ell]$, if there exists a choice for $Z_i$, such that $R_i \cup Z_i$ has a plane decomposition, then we may use \Cref{lem:planartraceshifting} to find an even dicycle in $\sigma(e_i)$, since $H_\mathcal{P}$ is odd.
    We adjust each $Z_i$ to lie in $\sigma(e_i)$ if this is the case.
    Clearly, if $Z_i \subseteq \sigma(e_i)$ for all $i \in [\ell]$, we have found an integral packing of even dicycles in $H_\mathfrak{M}$.

    Suppose there exists an $i \in [\ell]$ for which there does not exists a choice for $Z_i$, such that $R_i \cup Z_i$ has a plane decomposition, then we claim that we can still find an even dicycle $R_i \subseteq \sigma(e_i)$.
    For this purpose, we need to further refine the disk $e_i$.


    Let $e_1^i, \ldots e_r^i$ be the pairwise disjoint disks found in $\mathcal{D}(\nicefrac{\theta}{4}, e_i)$.
    We iterate over $[r]$ to refine $\delta'$, $\Gamma$, and $H_\mathcal{P}$, until we find an $e_j^i$ such that $H_\mathcal{P} \cup \sigma(e_j^i)$ contains an even dicycle.
    This $j \in [r]$ must exist, since $H_\mathcal{P} \cup \sigma(e_i)$ is even.
    We may choose $Z_i \subseteq H_\mathcal{P} \cup \sigma(e_j^i)$ and suppose that $R_i \cup Z_i$ does not have a plane decomposition, as we are otherwise done according to \Cref{lem:planartraceshifting}.
    Thus there exists a cross on $R_i$ in $Z_i$ consisting of the paths $Y_1,Y_2 \subseteq Z_i$.

    We note that $\nicefrac{\theta}{4} \geq \lceil \nicefrac{\omega}{2} \rceil$.
    Since $e_j^i \in \mathcal{D}(\nicefrac{\theta}{4}, e_i)$, there exists an $e_j^i$-tight dicycle $C_{\nicefrac{\theta}{4}}'$ in $P_{\nicefrac{\theta}{4}} \cup C_{\nicefrac{\theta}{4}} \cup \bigcup_{j=1}^t E^j_{\nicefrac{\theta}{4}}$.
    Due to $C_{\nicefrac{\theta}{4}}'$ and $R_i$ being concentric on $e_i$, the conditions on $\mathfrak{M}$ from \Cref{def:outline}, the fact that $\mathcal{P}$ is an odd transaction, and $\delta'$ being the $\mathcal{P}$-$A$-expansion of $\delta$, we conclude that there exists $\nicefrac{\theta}{4} + 1$ odd homogeneous dicycles $C_{\nicefrac{\theta}{4}}', C_{\nicefrac{\theta}{4} + 1}', \ldots , C_{\nicefrac{\theta}{2}}'$ in $\sigma(e_i)$, such that $C_{\nicefrac{\theta}{2}}' = R_i$ and $C_k' \subseteq P_k \cup C_k \cup \bigcup_{i=1}^t E^i_k$ for all $k \in [\nicefrac{\theta}{4}, \nicefrac{\theta}{2}]$.
    We note that $\nicefrac{\theta}{2} \geq \omega$, since we have $|I| \geq 2$ in this case.
    Using these dicycles it is easy to observe that there must exist a circle $\nicefrac{\theta}{2}$-outline for $e_j^i$ entirely contained in $\sigma(e_i)$.
    Observe that this new outline has degree at least 1 and thus we can find a segregated pair for its cycle set of order $\lfloor \nicefrac{\SC}{2} \rfloor$ for this outline and apply \Cref{lem:circlenonplanarshifting} to it and the crossing pair $Y_1, Y_2$ on $C_{\nicefrac{\theta}{2}}'$ to find an even dicycle on $\sigma(e_i)$.

    We can therefore always find an integral packing of $\ell$ even dicycles in this case.
    
    \textbf{Case 3:} Finally, we suppose that $2 \leq |I| < \ell$.
    Let $s = |I|$ and suppose w.l.o.g.\ that $e_1, \ldots , e_s \in I$.
    For any $i \in [s]$, consider the contents of $\mathcal{D}(\nicefrac{\theta}{4}, e_i)$.
    If $|\mathcal{D}(\nicefrac{\theta}{4}, e_i)| = 1$, we may proceed as in \textbf{Case 2} to find a circle $\nicefrac{\theta}{2}$-outline for the sole element of this set.
    On the other hand, if $|\mathcal{D}(\nicefrac{\theta}{4}, e_i)| \geq \ell - 1$, we may argue as in \textbf{Case 2} to find a packing of $\ell - 1$ even dicycles and find another even dicycle in some $e_j$ with $j \in [s] \setminus \{ i \}$, since $s \geq 2$.
    Should $2 \leq |\mathcal{D}(\nicefrac{\theta}{4}, e_i)| < \ell - 1$, we can iterate the procedure we have just outlined.

    Clearly, this procedure must terminate in at most $\ell - 1$ steps, yielding either an integral packing of $\ell$ even dicycles, or we find $h$ circle outlines embedded on $h$ pairwise disjoint disks on $\Delta$, each hosting an even dicycle.
    According to our arguments, we have $h \geq 1$ and we also know that $h < \ell$, as we can otherwise find an integral packing of $\ell$ even dicycles.
    Thus, each of the outlines must be at least a circle $\nicefrac{\theta}{2^{h-1}}$-outline, except possibly for the diamond outline provided to $e_h$, but even here we may choose all dicycles of $\mathcal{C}$ used to construct this outline disjoint from the infrastructure of the other outlines and thus guarantee, that the outlines and their infrastructure are mutually disjoint, satisfying the requirements of the lemma.
\end{proof}

Whilst we found new maelstroms surrounded by new outlines, we have not actually separated the contents of these new maelstroms, as the odd decomposition allows for edges with endpoints in two different maelstroms.
This will be the goal of the next lemma, which either yields a large half-integral packing of even dicycles, or a small set of vertices which separates all of the new maelstroms from one another.
Note that $\mathsf{g}$ is the function required for finding a cylindrical grid in \Cref{cor:getseparatingcylindricalgrid}.

\begin{lemma}\label{lem:buildnewmaelstromscircleapexset}
    Let $\theta , p, \ell$ be positive integers such that $\theta$ is a multiple of $2^{\ell+1} ( \mathsf{g}(3\ell) + \omega )$ and $p \geq 2\theta$, and let $D$ be a digraph with an odd decomposition $\delta$.
    Moreover, let $m$ be a maelstrom of $\delta$, together with a circle $\theta$-outline $\mathfrak{M}  = ( H , \mathcal{C}, \mathfrak{E}, \mathfrak{V} )$, with $(H_\mathfrak{M},\Omega_\mathfrak{M})$ being the $\mathfrak{M}$-society, and let $\mathcal{P} = \{ P_1, \ldots , P_p \}$ be an odd transaction on the $\mathfrak{M}$-society, flowing against $\mathfrak{M}$, together with a pleasant apex set $A \subseteq V(H_\mathfrak{M})$ such that $\mathcal{P}$ and $A$ are $\delta$-compatible.
    Let $\delta'$ be the $\mathcal{P}$-$A$-expansion of $\delta$ and let $H_\mathcal{P} = (H \cup S) - A$, where $S$ is the strip of $\mathcal{P}$ under $A$.

    Either $H_\mathfrak{M}$ is odd, there exists a half-integral packing of $\ell$ even dicycles in $H_\mathfrak{M}$, or there exist
    \begin{enumerate}
        \item an $h < \ell$ and pairwise disjoint disks $d_1, \ldots , d_h \subseteq m$ such that for each $i \in [h]$ there exists a $\theta_i'$-outline $\mathfrak{M}_i = ( H_i , \mathcal{C}_i, \mathfrak{E}_i, \mathfrak{V}_i )$ for $d_i' \subseteq d_i$, with $\theta_i' \geq \nicefrac{\theta}{2^{h-1}}$, $H_i$ is a subgraph of $H_\mathcal{P}$, the set $\{ H_1, \ldots , H_h \}$ is half-integral, the graph $\sigma_{\delta'}(d_i') \cup H_\mathcal{P}$ contains an even dicycle, and $d_i$ is bounded by the cycle corresponding to the $\mathfrak{M}_i$-society $(H_{\mathfrak{M}_i}, \Omega_{\mathfrak{M}_i})$, and

        \item there exists a vertex set $T \subseteq V(\bigcup_{i \in [h]} H_{\mathfrak{M}_i})$ with $|T| \leq 2\ell^2$, such that for all distinct $i,j \in [h]$ there does not exist a $V(\Omega_{\mathfrak{M}_i})$-$V(\Omega_{\mathfrak{M}_j})$- or a $V(\Omega_{\mathfrak{M}_j})$-$V(\Omega_{\mathfrak{M}_i})$-path in $D[V(H_{\mathfrak{M}_i} \cup H_{\mathfrak{M}_j})] - T$.
    \end{enumerate}
\end{lemma}
\begin{proof}
    We may assume that $H_\mathfrak{M}$ is in fact not odd and that we cannot find a half-integral packing of $\ell$ even dicycles in $H_\mathfrak{M}$.
    Furthermore, via the methods employed in the proof of \Cref{lem:buildnewmaelstromscircle}, we can find the desired disks $d_1, \ldots , d_h$ and the associated circle outlines $\mathfrak{M}_1, \ldots , \mathfrak{M}_h$.
    Note that by construction all of the outlines $\mathfrak{M}_i$ with $i \in [h]$ have degree at least 1.
    According to \Cref{lem:circleoutlinererouting}, we can therefore also find a $\nicefrac{\theta_i}{2}$-rim $\mathfrak{R}_i = ( H_i', \mathcal{C}_i' )$ for each $\mathfrak{M}_i$, such that the $\mathfrak{R}_i$-society is the $\mathfrak{M}_i$-society and there exists an organised $\mathcal{C}_i'$-pair $(\mathcal{L}_1^i, \mathcal{L}_2^i)$ of order $\nicefrac{\theta_i}{2}$.

    Clearly, if $h = 1$, there is nothing to show.
    So we may suppose that $2 \leq h < \ell$.
    For the moment, we also suppose that $\mathfrak{M}_1, \ldots , \mathfrak{M}_h$ are all circle outlines and thus each $\sigma(d_i)$ for $i \in [h]$ contains some vertex of $P_1 \in \mathcal{P}$ by construction, allowing us to assume that $d_1, \ldots , d_h$ are labelled in order of their occurrence along $P_1$ when traversed according to the direction of its edges.
    For each $i \in [h]$, we let $R_i$ be the dicycle corresponding to $V(\Omega_{\mathfrak{M}_i})$.
    
    Let $i, j \in [h]$ with $i < j$, let $t_i$ be the degree of $\mathfrak{M}_i$, and let $\mathfrak{E}_i = \{ \mathcal{E}^i_1, \ldots \mathcal{E}^i_{t_i} \}$.
    Suppose there exists a $V(\Omega_{\mathfrak{M}_j})$-$V(\Omega_{\mathfrak{M}_i})$-linkage $\mathcal{L} = \{ L_1, \ldots , L_\ell \}$ in $D[V(H_{\mathfrak{M}_i} \cup H_{\mathfrak{M}_j})]$, then we note that since $H_{\mathfrak{M}_i}$ and $H_{\mathfrak{M}_j}$ are disjoint and $\delta'$ provides an odd decomposition for $H_i$ and $H_j$, there must exists at least $\ell$ edges $e_1, \ldots , e_\ell$ with their head in $V(\sigma(d_i))$ and their tail in $V(\sigma(d_j))$, such that $e_s \in E(L_s)$ for $s \in [\ell]$.
    In particular, due to $\delta'$, we know that each $L_s$ must intersect all members of both $\mathcal{C}_i'$ and $\mathcal{C}_j'$ for all $s \in [h]$, as we could otherwise find an edge jumping from the middle of the decomposition of $H_\mathcal{P}$ to somewhere inside the part of $m$ that is split into new maelstroms.
    Specifically, when traversed in the direction of their edges, each $L_s$ must intersect the members $\mathcal{C}_i'$ in descending order of their indices and the members of $\mathcal{C}_j'$ in ascending order of their indices, though the path is of course allowed to return to any cycle in $\mathcal{C}_i' \cup \mathcal{C}_j'$ once it has intersected it for the first time.

    Several applications of \Cref{lem:detour}, along with the fact that $i < j$, let us reroute the heads of the paths in $\mathcal{L}$ first via the members of $\mathcal{C}_i \cup \bigcup_{s \in [t_i]} \mathcal{E}^i$ onto $\mathcal{P}$, if they did not intersect these paths already, and from there, we can then reroute the heads of $\mathcal{L}$ to all reside on $R_j$.
    Since $\theta_i \geq \nicefrac{\theta}{2^{h-1}}$, as well as $2 \leq h < \ell$, and $\theta$ being a multiple of $2^\ell ( \mathsf{g}(3\ell) + 4\ell^2 + \omega )$, we conclude that 
    \[ \nicefrac{\theta_i}{2} \geq \nicefrac{\theta}{2^h} \geq \mathsf{g}(3\ell) . \]
    Thus, we can find a planar grid $G$ of order $3\ell$ within $H_j$ whose trace separates $d_j$ and the trace of $R_j$ using \Cref{lem:findseggrid} and \Cref{lem:seggridtoseparatingcylgrid}.
    We may now cut off the paths in $\mathcal{L}$ such that they all start inside of $G$ and may also further move the heads of $\mathcal{L}$ along the paths in $\mathcal{P}$, allowing us to find $\ell$ disjoint grids $G_1, \ldots , G_\ell$ of order three within $G$ such that each has a perimeter jump and no three grids together with their perimeter jump intersect in a common vertex.
    This allows us to use \Cref{lem:gridplusjumpmakesoddbicycle} to find a half-integral packing of even dicycles within $H_\mathfrak{M}$.
    Therefore no such linkage $\mathcal{L}$ may exist, which implies the existence of a directed $V(\Omega_{\mathfrak{M}_j})$-$V(\Omega_{\mathfrak{M}_i})$-separator of size at most $\ell$ according to Menger's theorem.

    If there exists a $V(\Omega_{\mathfrak{M}_i})$-$V(\Omega_{\mathfrak{M}_j})$-linkage of order $\ell$ or larger, we may proceed analogously, though we must extend the paths first off of $\mathcal{P}$ onto the infrastructure of $\mathfrak{M}$ and reach the society of $\mathfrak{M}_j$ by moving in the opposite direction compared to the previous case.

    Should some $\mathfrak{M}_i$ with $i \in [h]$ be a diamond outline, then it is the sole diamond outline amongst $\mathfrak{M}_1, \ldots \mathfrak{M}_h$, as discussed in the proof of \Cref{lem:buildnewmaelstromscircle}.
    If there exists some $V(\Omega_{\mathfrak{M}_j})$-$V(\Omega_{\mathfrak{M}_i})$, or respectively a $V(\Omega_{\mathfrak{M}_i})$-$V(\Omega_{\mathfrak{M}_j})$-linkage with $\ell$ or more paths, for some $j \in [h] \setminus \{ i \}$, we can again find a circle rim $\mathfrak{R}_j$ corresponding to $\mathfrak{M}_j$ and are able to move the tails, respectively the heads of $\mathcal{L}$ along $\mathcal{P}$ to first meet the infrastructure of $\mathfrak{M}$ and then meet the $\mathfrak{R}_j$-society using \Cref{lem:detour}.
    Analogous to our previous arguments in this proof, this allows us to find a small directed $V(\Omega_{\mathfrak{M}_j})$-$V(\Omega_{\mathfrak{M}_i})$-, and respectively a small directed $V(\Omega_{\mathfrak{M}_i})$-$V(\Omega_{\mathfrak{M}_j})$-separator.

    Clearly, since we cannot find a half-integral packing of $\ell$ even dicycles, we must therefore be able to combine all of these separators into a set of vertices $T \subseteq V(H_\mathfrak{M})$ with $|T| \leq 2\ell^2$, such that for all distinct $i,j \in [h]$ there does not exist a $V(\Omega_{\mathfrak{M}_i})$-$V(\Omega_{\mathfrak{M}_j})$- or a $V(\Omega_{\mathfrak{M}_j})$-$V(\Omega_{\mathfrak{M}_i})$-path in $D[V(H_{\mathfrak{M}_i} \cup H_{\mathfrak{M}_j})] - T$.
\end{proof}

\subsection{Subdividing a (refined) diamond outline}

Now we first prove a version of \Cref{lem:buildnewmaelstromscircle} for (refined) diamond outlines.
However, in this version of the statement we forgo asking the odd transaction to go against the flow, though we do restrict ourselves to calm transactions instead.
Note that since we allowed the surplus degree of a refined diamond outline to be 0, this also takes care of regular diamond outlines.
As promised earlier, the fundamental ideas of the proofs here are not much different from the proofs for the circle outlines, but many parts of the proofs get considerably more technical.
We note that the case in which the condition $\mathfrak{B}_i = \mathfrak{V}_i = \emptyset$ of our statement is true, corresponds to $H_{\mathfrak{M}_i}$ being a diamond or circle outline.

\begin{lemma}\label{lem:buildnewmaelstromsrefineddiamond}
    Let $\theta , p, \ell$ be positive integers and let $g$ be non-negative, such that $\theta$ is a common multiple of $\ell^{\ell+1} \omega$ and $p \geq 2\theta(g+1)$, and let $D$ be a digraph with an odd decomposition $\delta$.
    Moreover, let $m$ be a maelstrom of $\delta$, together with a diamond $\theta$-outline $\mathfrak{M}$, with $(H_\mathfrak{M},\Omega_\mathfrak{M})$ being the $\mathfrak{M}$-society, and a refined $\theta$-outline $\mathfrak{N} = ( H , \mathcal{C}, \mathfrak{E}, \mathfrak{B}, \mathfrak{V}, \mathfrak{D} )$ of $n \subseteq m$ with roughness $g$.
    Further, let $\mathcal{P} = \{ P_1, \ldots , P_p \}$ be a calm, odd transaction on the $\mathfrak{N}$-society $(H_\mathfrak{N},\Omega_\mathfrak{N})$, together with a pleasant apex set $A \subseteq V(H_\mathfrak{N})$ such that $\mathcal{P}$ and $A$ are $\delta$-compatible.
    Let $\delta'$ be the $\mathcal{P}$-$A$-expansion of $\delta$, let $H_\mathcal{P} = (H \cup S) - A$, where $S$ is the strip of $\mathcal{P}$ under $A$.

    Either $H_\mathfrak{N}$ is odd, there exists a quarter-integral packing of $\ell$ even dicycles in $H_\mathfrak{N}$, or there exist an $h < \ell$ and pairwise disjoint disks $d_1, \ldots , d_h \subseteq n$ such that for each $i \in [h]$ there exists a $\theta_i$-outline $\mathfrak{M}_i = ( H_i , \mathcal{C}_i, \mathfrak{E}_i, \mathfrak{B}_i, \mathfrak{V}_i, \mathfrak{D}_i )$ (of roughness $g$) for $d_i' \subseteq d_i$, with $\theta_i \geq \nicefrac{\theta}{\ell^{h-1}}$, the graph $\sigma_{\delta'}(d_i') \cup H_\mathcal{P}$ contains an even dicycle, $\{ H_1, \ldots , H_h \}$ is a third-integral set of graphs, $d_i$ is bounded by the cycle corresponding to the $\mathfrak{M}_i$-society $(H_{\mathfrak{M}_i}, \Omega_{\mathfrak{M}_i})$, and either $H_{\mathfrak{M}_i} \subset H_\mathfrak{M}$, if $\mathfrak{B}_i = \mathfrak{V}_i = \emptyset$, or $H_{\mathfrak{M}_i} \subset H_\mathfrak{N}$.
\end{lemma}
\begin{proof}
    Let $t$ be the degree and $t'$ be the surplus degree of $\mathfrak{N}$.
    Suppose that $H_\mathfrak{N}$ contains an even dicycle and that there does not exist a quarter-integral packing of $\ell$ even dicycles in $H_\mathfrak{N}$.
    As a consequence, we may also assume that $\ell \geq 2$.
    Further, we let $Q$ be the cycle corresponding to $\Omega_\mathfrak{N}$.

    We set $\mathcal{P}_{\nicefrac{\theta(g+1)}{\ell}-1} = \{ P_{\nicefrac{\theta(g+1)}{\ell}-1}, \ldots , P_{p - (\nicefrac{\theta(g+1)}{\ell} + 1)} \}$ and consider the elements $e_1 , \ldots , e_h$ of $\mathcal{D}(\mathcal{P}_{\nicefrac{\theta(g+1)}{\ell}-1}, m)$.
    If there exists some $i \in [h]$ such that $H_\mathcal{P} \cup \sigma(e_i)$ is odd, then with the help of \Cref{lem:oddmaelstromdecomposition}, we refine $\delta'$, $\Gamma$, and $H_\mathcal{P}$ as in the proof of \Cref{lem:buildnewmaelstromscircle}.
    Clearly, this cannot be true for all $i \in [h]$.
    We iterate over all $e_1, \ldots , e_h$ and refine $\delta'$, $\Gamma$, and $H_\mathcal{P}$ whenever possible.

    Let $I \subseteq [h]$ be the maximal set such that for all $i \in I$ the graph $H_\mathcal{P} \cup \sigma(e_i)$ contains an even dicycle $Z_i$.
    As in the proof of \Cref{lem:buildnewmaelstromscircle}, we distinguish three cases.

    \textbf{Case 1:} Suppose $|I| = 1$.
    Then let $i \in I$ and note that since $e_i \in \mathcal{D}(\mathcal{P}_{\nicefrac{\theta(g+1)}{\ell}-1}, m)$ there exists a cycle in $C_{\nicefrac{\theta}{\ell}} \cup P_{\nicefrac{\theta(g+1)}{\ell}} \cup \bigcup_{j=1}^t E^j_{\nicefrac{\theta}{\ell}} \cup \bigcup_{j=1}^{t'} B^j_{\nicefrac{\theta(g+1)}{\ell}}$ whose trace bounds $e_i$.
    This allows us to find a $\theta$-outline for $e_i$:
    \begin{enumerate}
        \item Either we use the refined diamond $\theta$-outline $( H_\mathcal{P}, \mathcal{C}, \mathfrak{E}, \mathfrak{B} \cup \{ \mathcal{B}_{t'+1} \}, \mathfrak{V}', \mathfrak{D}' )$ of roughness $g$, if $\mathcal{B}_{t+1}$ satisfies \textbf{\textsf{E6}} of \Cref{def:refinedoutline}, where $\mathcal{B}_{t+1} = \{ P_1, \ldots , P_{\theta(g+1)} \}$, and the sets $\mathfrak{V}$ and $\mathfrak{D}'$ are derived from $\delta'$, $\Gamma$, and $H_\mathcal{P}$, or

        \item since $\{ P_1, \ldots , P_{\theta(g+1)} \}$ does not satisfy \textbf{\textsf{E6}} of \Cref{def:refinedoutline} and $\mathcal{P}$ is calm, we can conclude that the heads and tails of $\mathcal{P}$ either appear on the diamond $Q'$ corresponding to $\Omega_\mathfrak{M}$, or they may be rerouted onto $Q'$ through the paths in the sets of $\mathfrak{B}$ using \Cref{lem:detour}, resulting in a transaction $\mathcal{E}_{t+1}$ of order $\theta$ on the $\mathfrak{M}$-society, with $E_j^{t+1}$ containing large parts of $P_{\nicefrac{\theta(g+1)}{\ell} - (\nicefrac{\theta}{2} - j)}$ for $j \in [\theta]$, which allows us to construct a diamond or circle $\theta$-outline $( H_\mathcal{P}, \mathcal{C}, \mathfrak{E} \cup \{ \mathcal{E}_{t+1} \}, \mathfrak{V}' )$, where the set $\mathfrak{V}'$ is derived from $\delta'$, $\Gamma$, and $H_\mathcal{P}$.
        (We can also construct an organised $\mathcal{C}'$-pair of order $2\omega$ from the paths in the sets within $\mathfrak{E}$ and $\mathcal{P}$ itself provides candidates for these paths.) 
    \end{enumerate}
    We can therefore satisfy either \Cref{def:refinedoutline} with the first item, or \Cref{def:outline} with the second item, and thus find an outline fitting our statement.

    \textbf{Case 2:} Suppose instead that $|I| \geq \ell$.
    Without loss of generality, we let $e_1, \ldots , e_\ell \in I$ to ease notation.
    We note that $e_i$ and $e_j$ are disjoint for distinct $i,j \in [\ell]$ and since $e_1, \ldots , e_\ell \in \mathcal{D}(\mathcal{P}_{\nicefrac{\theta(g+1)}{\ell}-1}, m)$, for each $i \in [\ell]$ there exists a cycle $R_i$ in $C_{\nicefrac{\theta}{\ell}} \cup P_{\nicefrac{\theta(g+1)}{\ell}} \cup \bigcup_{j=1}^t E^j_{\nicefrac{\theta}{\ell}} \cup \bigcup_{j=1}^{t'} B^j_{\nicefrac{\theta(g+1)}{\ell}}$ whose trace bounds $e_i$.
    Accordingly, we also have $V(R_i) \cap V(R_j) = \emptyset$ for all distinct $i,j \in [\ell]$.


    Let $c_i = \nicefrac{\theta}{\ell^2} - \omega (i-1)$ and let $e_1^i, \ldots e_r^i$ be the pairwise disjoint disks found in $\mathcal{D}(c_i, e_i)$ for some fixed but arbitrary $i \in [\ell]$.
    We iterate over the elements of $[r]$ to refine $\delta'$, $\Gamma$, and $H_\mathcal{P}$, until we find an $e_j^i$ such that $H_\mathcal{P} \cup \sigma(e_j^i)$ contains an even dicycle.
    This $j$ must exist, since $H_\mathcal{P} \cup \sigma(e_i)$ is even.

    We note that $\nicefrac{\theta}{\ell^2} \geq \omega \ell$, since $\theta$ is a multiple of $\ell^{\ell+1} \omega$, and we also note that there exists an $e_j^i$-tight cycle in $P_{c_i(g+1)} \cup C_{c_i} \bigcup_{s=1}^t E^s_{c_i} \cup \bigcup_{s=1}^{t'} B^s_{c_i(g+1)}$.
    Based on this, we now have to make quite a technical construction to get the new outlines.
    
    To help the reader understand the intuition behind our choices, we provide an explanation ahead of the construction, whilst referencing it.
    Note that within $\mathcal{C}'$ below there are essentially two batches of dicycles, one with indices higher than $c_i$ and one with indices lower than $c_i$.
    Recall that elements of the infrastructure in an outline with higher indices lie further away from the maelstrom.
    Thus the elements with high indices naturally have to form what becomes the rim of the new outline.
    We can guarantee that this part is sufficiently structured to get an odd decomposition, since we have identified $e_j^i$ as the new maelstrom and have refined $\delta'$, $\Gamma$, and $H_\mathcal{P}$ appropriately throughout.
    The choice of the value $c_i$ then guarantees that for any $i' \in I \setminus \{ i \}$ we will not choose any of the same infrastructure.
    This leaves the different outlines to only overlap when different parts of the infrastructure of two outlines pass through each other. 
    
    Due to the existence of the $e_j^i$-tight cycle and $R_i$, one of the following two options is a $2\omega$-outline:
    \begin{enumerate}
        \item Either we choose the refined diamond $2\omega$-outline $\mathfrak{M}_i' = (H_i' , \mathcal{C}', \mathfrak{E}', \mathfrak{B}' = \mathfrak{B}'' \cup \{ \mathcal{B}_{t'+1}' \}, \mathfrak{V}', \mathfrak{D}' )$, if $\mathcal{B}_{t+1}'$ satisfies \textbf{\textsf{E6}} of \Cref{def:refinedoutline}, with
        \[ \mathcal{C}' = \{ C_{c_i+\nicefrac{\theta}{4}}, C_{c_i+\nicefrac{\theta}{4} - 1}, \ldots , C_{c_i+\nicefrac{\theta}{4} - (\omega - 1)}, C_{c_i}, C_{c_i - 1}, \ldots , C_{c_i - (\omega - 1)} \} , \]
        the $t$ sets in $\mathfrak{E}'$ being analogous restrictions of the transactions $\mathcal{E}_j$ with $j \in [t]$, the set
        \[ \mathcal{B}_{t'+1}' = \{ P_{c_i(g+1) + \nicefrac{\theta}{4}}, \ldots , P_{c_i(g+1)+\nicefrac{\theta}{4} - (\omega - 1)}, P_{c_i(g+1)}, \ldots , P_{c_i (g+1) - (\omega - 1)} \} , \]
        the $t'$ sets in $\mathfrak{B}''$ being analogous restrictions of the transactions $\mathcal{B}_j$ with $j \in [t']$, the graph
        \[ H_i' = (H_\mathcal{P} \cap \sigma(e_i)) \cup \bigcup_{C \in \mathcal{C}'} C \cup \bigcup_{E \in \mathcal{E}' \in \mathfrak{E}'} E \cup \bigcup_{B \in \mathcal{B}' \in \mathfrak{B}'} B , \]
        and the sets $\mathfrak{V}'$ and $\mathfrak{D}'$ are derived from $\delta'$, $\Gamma$, and $H_\mathcal{P}$, while our outline is explicitly meant to reside in $(H_\mathcal{P} \cup \sigma(e_j^i)) \cup H_i'$\footnote{Since our intent in this case is to find a quarter-integral packing of even dicycles, we do not need to decompose the entire graph here.}, or

        \item analogous to the situation in \textbf{Case 1}, we can find a transaction $\mathcal{E}_{t+1}$ of order $\omega$ on the $\mathfrak{M}$-society, with $E_s^{t+1}$ containing large parts of $P_{c_i(g+1) + \nicefrac{\theta(g+1)}{4} - (s - 1)}$ for $s \in [\omega]$, and for $s \in [\omega+1, 2\omega]$, $E_s^{t+1}$ containing much of $P_{c_i(g+1) - (s - 1)}$.
        This allows us to construct a diamond or circle $2\omega$-outline $\mathfrak{M}_i' = ( H_i', \mathcal{C}', \mathfrak{E}'' = \mathfrak{E}' \cup \{ \mathcal{E}_{t+1} \}, \mathfrak{V}' )$, where $\mathcal{C}'$ and $\mathcal{E}''$ are defined analogously to the corresponding sets in the previous item, we let
        \[ H_i' = (H_\mathcal{P} \cap \sigma(e_i)) \cup \bigcup_{C \in \mathcal{C}'} C \cup \bigcup_{E \in \mathcal{E}'' \in \mathfrak{E}''} E , \]
        and the set $\mathfrak{V}'$ is again derived from $\delta'$, $\Gamma$, and $H_\mathcal{P}$, while our outline is explicitly meant to reside in $(H_\mathcal{P} \cup \sigma(e_j^i)) \cup H_i'$
        (We can also construct an organised $\mathcal{C}'$-pair of order $2\omega$ from the paths in the sets within $\mathfrak{E}$.) 
    \end{enumerate}
    For each $i \in [\ell]$, we can find such a $2\omega$-outline $\mathfrak{M}_i'$ around a disk $e_i' \subseteq e_i$, which we may treat as a maelstrom in $\delta'$.
    Thanks to our choices, we know that $\{ H_1', \ldots , H_\ell' \}$ are third-integral, since only the elements of $\mathcal{C}'$, the sets in $\mathfrak{E}'$, and the sets in $\mathfrak{B}'$ may pairwise intersect each other.

    We claim that this allows us to find a quarter-integral packing of $\ell$ even dicycles.
    Recall that each outline $\mathfrak{M}_i'$ must have some even dicycle $Z_i$ with $Z_i \subseteq H_i' \cup \sigma(e_i')$, such that $V(\sigma(e_i')) \cap V(Z_i) \neq \emptyset$.
    Note that therefore, if $Z_i$ and $R_i$ are disjoint, we have $Z_i \subseteq \sigma(e_i)$ and clearly each $Z_i$ for which this is true can be used for our packing.

    For any $Z_i$ which contains a cross $P_i, Q_i$ on the $\mathfrak{M}_i$-society, we can apply \Cref{lem:applyshifting} to find an even dicycle within $H_i' \cup P_i \cup Q_i$.
    Should $Z_i \cup R_i$ instead have a planar decomposition, we consider the graph in which we found $R_i$ and let $X_1, \ldots , X_a$ be the set of mutually distinct dicycles within $C_{\nicefrac{\theta}{\ell}} \cup P_{\nicefrac{\theta(g+1)}{\ell}} \cup \bigcup_{j=1}^t E^j_{\nicefrac{\theta}{\ell}} \cup \bigcup_{j=1}^{t'} B^j_{\nicefrac{\theta(g+1)}{\ell}}$.
    We note that $a$ must be positive and that $Z_i \cup \bigcup_{j = 1}^x X_j$ has a planar decomposition.
    To shift $Z_i$ into $e_i$, we can now repeatedly apply \Cref{lem:planartraceshifting} to $Z_i$ and members of $X_1, \ldots , X_a$ until the resulting even dicycle lies on one side of the trace of each $X_j$ with $j \in [a]$ and thus must be contained within $H_{\mathfrak{M}_i'} \cup H_i'$.
    Thus, since $\{ H_1', \ldots , H_\ell' \}$ are third-integral, we can find a quarter-integral packing of even dicycles within $\bigcup_{i \in [\ell]} (H_i' \cup Z_i)$.

    \textbf{Case 3:} Suppose instead that $\ell > |I| > 1$.
    Let $s = |I|$ and suppose w.l.o.g.\ that $e_1, \ldots , e_s \in I$.
    Consider the contents of $\mathcal{D}(\nicefrac{\theta}{\ell^2}, e_1)$ and let $I_2^1$ be the set of disks $f \in \mathcal{D}(\nicefrac{\theta}{\ell^2}, e_1)$ for which $H_\mathcal{P} \cup \sigma(f)$ contains an even dicycle.
    $I_2^1$ must be non-empty, since $e_1 \in I$.
    
    If $|I_2^1| \geq \ell - 1$, we may proceed as in \textbf{Case 2} to find a quarter-integral packing of $\ell$ even dicycles, since $s \geq 2$.
    Should $|I_2^1| = 1$ instead, we let $f \in I_2^1$ be the sole representative and save this as a tuple $(f,2)$, where 2 is the exponent of the divisor in $\nicefrac{\theta}{\ell^2}$.

    In case $\ell - 1 > |I_2^1| > 1$, we iterate on the elements $f \in I_2^1$ by considering $\mathcal{D}(\nicefrac{\theta}{\ell^3}, f)$.
    Clearly this procedure must terminate, as in each step of the iteration, we are required to find one even dicycle less.
    Thus, if we also perform this procedure over all other members of $I$, we either find a quarter-integral packing of $\ell$ even dicycles in $H_\mathfrak{N}$ or we find at most $\ell - 1$ tuples $(f_1, x_1), \ldots , (f_k, x_k)$, such that for each $i \in [k]$ we have $x_k \leq \ell - 1$ and $f_1, \ldots , f_k$ are pairwise disjoint disks.
    We sort these tuples such that for $i,j \in [k]$ with $i < j$ we have $x_i \geq x_j$.

    Our next goal will be to find outlines for all of the disk we have found in this fashion and choose these outlines such that they are ultimately third-integral.
    For this purpose, for each $(f_i,x_i)$ we use both $i$ and the value of $x_i$ to locate the particular part of the old infrastructure which we draw from.
    This way we can guarantee that our goals are met.
    As described in \textbf{Case 2}, the breakpoint in each of the definitions of the following sets is set such that the objects with higher indices lie in the part of the graph that we can already guarantee an odd decomposition for.

    Let $r \in [k]$ be the largest index such that $x_1, \ldots , x_r$ are all equal, let $x = x_1$, and let $\theta' = \nicefrac{\theta}{\ell^{k - 1}}$.
    For $f_1, \ldots , f_r$ we can then build $\theta'$-outlines for each $i \in [r]$ as follows:
    \begin{enumerate}
        \item Either we can find a refined diamond $\theta'$-outlines $\mathfrak{M}_i = (H_i, \mathcal{C}_i, \mathfrak{E}_i, \mathfrak{B}_i \cup \{ \mathcal{B}^i_{t'+1} \}, \mathfrak{V}_i, \mathfrak{D}_i )$ of roughness $g$, if $\mathcal{B}^i_{t'+1}$ satisfies \textbf{\textsf{E6}} of \Cref{def:refinedoutline}, where
        \begin{align*}
            \mathcal{B}^i_{t'+1} = & \{ P_{\nicefrac{\theta(g+1)}{\ell^{x-1}} - h_i}, P_{\nicefrac{\theta(g+1)}{\ell^{x-1}} - (h_i + 1)}, \ldots , P_{\nicefrac{\theta(g+1)}{\ell^{x-1}} - (h_i + (\nicefrac{\theta'(g+1)}{2} - 1))}, \\
            & \quad P_{\nicefrac{\theta(g+1)}{\ell^x} - h_i}, P_{\nicefrac{\theta(g+1)}{\ell^x} - (h_i + 1)}, \ldots , P_{\nicefrac{\theta(g+1)}{\ell^x} - (h_i + (\nicefrac{\theta'(g+1)}{2} - 1))} \} 
        \end{align*}
        with $h_i = \nicefrac{\theta'(g+1)(i-1)}{2}$.
        We select the transactions within $\mathfrak{B}_i$ as being analogous restrictions of transactions in $\mathfrak{B}$.
        Further, we choose $\mathcal{C}_i \subset \mathcal{C}$ such that
        \[ \mathcal{C}_i = \{ C_{\nicefrac{\theta}{\ell^{x-1}} - h_i'}, \ldots , C_{\nicefrac{\theta}{\ell^{x-1}} - (h_i' + (\nicefrac{\theta'(g+1)}{2} - 1))}, C_{\nicefrac{\theta}{\ell^x} - h_i'}, \ldots , C_{\nicefrac{\theta}{\ell^{x-1}} - (h_i' + (\nicefrac{\theta'(g+1)}{2} - 1))} \} \]
        with $h_i' = \nicefrac{\theta'(i-1)}{2}$, and the transactions in $\mathfrak{E}_i$ are chosen analogously to the set $\mathcal{C}_i$ from within the transactions in $\mathfrak{E}$.
        We then define
        \[ H_i = (H_\mathcal{P} \cap \sigma(f_i)) \cup \bigcup_{C \in \mathcal{C}_i} C \cup \bigcup_{E \in \mathcal{E} \in \mathfrak{E}_i} E \cup \bigcup_{B \in \mathcal{B} \in \mathfrak{B}_i} B  \]
        and the sets $\mathfrak{V}_i$ and $\mathfrak{D}_i$ are derived from through $\delta'$, $\Gamma$, and $H_\mathcal{P}$.

        \item If $\mathcal{B}^i_{t'+1}$, defined as above, does not satisfy \textbf{\textsf{E6}} of \Cref{def:refinedoutline} then, because $\mathcal{P}$ is calm, we can reroute the heads and tails of select paths of $\mathcal{P}$ using \Cref{lem:detour} to form a transaction $\mathcal{E}^i_{t+1}$ of order $\theta'$ on the $\mathfrak{M}$-society, using the paths in the let
        \begin{align*}
            & \{ P_{\nicefrac{\theta(g+1)}{\ell^{x-1}} - h_i'}, P_{\nicefrac{\theta(g+1)}{\ell^{x-1}} - (h_i' + 1)}, \ldots , P_{\nicefrac{\theta(g+1)}{\ell^{x-1}} - (h_i' + (\nicefrac{\theta'(g+1)}{2} - 1))}, \\
            & \quad P_{\nicefrac{\theta(g+1)}{\ell^x} - h_i'}, P_{\nicefrac{\theta(g+1)}{\ell^x} - (h_i' + 1)}, \ldots , P_{\nicefrac{\theta(g+1)}{\ell^x} - (h_i' + (\nicefrac{\theta'(g+1)}{2} - 1))} \}
        \end{align*}
        and the paths with the same indices from the transactions within $\mathfrak{B}$.
        This allows us to construct a diamond or circle $\theta'$-outline $\mathfrak{M}_i = (H_i, \mathcal{C}_i, \mathfrak{E}_i \cup \{ \mathcal{E}^i_{t+1} \}, \mathfrak{V}_i )$, where $\mathcal{C}_i$ and $\mathcal{E}_i$ are defined analogously to the corresponding sets in the previous item, we have
        \[ H_i' = (H_\mathcal{P} \cap \sigma(f_i)) \cup \bigcup_{C \in \mathcal{C}_i} C \cup \bigcup_{E \in \mathcal{E} \in \mathfrak{E}_i} E , \]
        and the set $\mathfrak{V}_i$ is again derived from $\delta'$, $\Gamma$, and $H_\mathcal{P}$.
        (We can also once more construct an organised $\mathcal{C}_i$-pair of order $2\omega$ from the paths in the sets within $\mathfrak{E}$, if $\mathfrak{M}_i$ is a circle outline.)
    \end{enumerate}
    Our choices allow us to conclude that $\{ H_1, \ldots , H_r \}$ are third-integral (see \textbf{Case 2}) and each $\mathfrak{M}_i$ is a $\theta'$-outline of some type.
    We note that thanks to our choices, the paths 
    \[ P_{\nicefrac{\theta(g+1)}{\ell^x} - h_{r+1}}, P_{\nicefrac{\theta(g+1)}{\ell^x} - (h_{r+1} + 1) }, \ldots , P_{\nicefrac{\theta(g+1)}{\ell^{x-1}} + 1} \]
    remain untouched and the dicycles $C_{\nicefrac{\theta}{\ell^x} - h_{r+1}'}, \ldots , C_{\nicefrac{\theta}{\ell^{x-1}} + 1}$, as well as the paths within transactions of $\mathfrak{E}$ with the same indices, are also not used in $\mathfrak{M}_1, \ldots , \mathfrak{M}_r$.
    This allows us to construct $\theta'$-outlines $\mathfrak{M}_{r+1}, \ldots , \mathfrak{M}_k$ for the remaining tuples $(f_{r+1}, x_{r+1}), \ldots , (f_k, x_k)$, if we follow this method for choosing our outlines and adjust indices appropriately.
    This then yields $k$ new $\theta'$-outlines $\mathfrak{M}_1, \ldots , \mathfrak{M}_k$ of differing types.
    Note that $f_i$ is the maelstrom that is outlined by $\mathfrak{M}_i$.
    The disk $d_i$ mentioned in the statement is the disk we found $f_i$ in, or more accurately, the restriction of that disk to fit the society of $\mathfrak{M}_i$.
    Since $\theta' = \nicefrac{\theta}{\ell^{k - 1}}$ and $k < \ell$, we therefore found the desired outlines.
\end{proof}

We can close this section by proving an analogous statement to \Cref{lem:buildnewmaelstromscircleapexset} for refined diamond outlines, again slightly more general than the previous version.

\begin{lemma}\label{lem:buildnewmaelstromsrefineddiamondapexset}
    Let $\theta , p, \ell , g$ be positive integers such that $\theta$ is a multiple of $2 \ell^\ell ( \mathsf{g}(3\ell) + 4\ell^2 + \omega )$ and $p \geq 2\theta(g+1)$, and let $D$ be a digraph with an odd decomposition $\delta$.
    Moreover, let $m$ be a maelstrom of $\delta$, together with a diamond $\theta$-outline $\mathfrak{M}$, with $(H_\mathfrak{M},\Omega_\mathfrak{M})$ being the $\mathfrak{M}$-society, and a refined $\theta$-outline $\mathfrak{N} = ( H , \mathcal{C}, \mathfrak{E}, \mathfrak{B}, \mathfrak{V}, \mathfrak{D} )$ of $n \subseteq m$ with roughness $g$.
    Further, let $\mathcal{P} = \{ P_1, \ldots , P_p \}$ be a calm, odd transaction on the $\mathfrak{N}$-society $(H_\mathfrak{N},\Omega_\mathfrak{N})$, together with a pleasant apex set $A \subseteq V(H_\mathfrak{N})$ such that $\mathcal{P}$ and $A$ are $\delta$-compatible.
    Let $\delta'$ be the $\mathcal{P}$-$A$-expansion of $\delta$, let $H_\mathcal{P} = (H \cup S) - A$, where $S$ is the strip of $\mathcal{P}$ under $A$.

    Either $H_\mathfrak{N}$ is odd, there exists a quarter-integral packing of $\ell$ even dicycles in $H_\mathfrak{N}$, or there exist
    \begin{enumerate}
        \item an $h < \ell$ and pairwise disjoint disks $d_1, \ldots , d_h \subseteq m$ such that for each $i \in [h]$ there exists a $\theta_i$-outline $\mathfrak{M}_i = ( H_i , \mathcal{C}_i, \mathfrak{E}_i, \mathfrak{B}_i, \mathfrak{V}_i, \mathfrak{D}_i )$ (of roughness $g$) for $d_i' \subseteq d_i$, with $\theta_i \geq \nicefrac{\theta}{\ell^{h-1}}$, the graph $\sigma_{\delta'}(d_i') \cup H_\mathcal{P}$ contains an even dicycle, $\{ H_1, \ldots , H_h \}$ is third-integral, $d_i$ is bounded by the cycle corresponding to the $\mathfrak{M}_i$-society $(H_{\mathfrak{M}_i}, \Omega_{\mathfrak{M}_i})$, and either $H_{\mathfrak{M}_i} \subset H_\mathfrak{M}$, if $\mathfrak{B}_i = \mathfrak{V}_i = \emptyset$, or $H_{\mathfrak{M}_i} \subset H_\mathfrak{N}$, and

        \item there exists a vertex set $T \subseteq V(\bigcup_{i \in [h]} H_{\mathfrak{M}_i})$ with $|T| \leq 4\ell^2$, such that for all distinct $i,j \in [h]$ there does not exist a $V(\Omega_{\mathfrak{M}_i})$-$V(\Omega_{\mathfrak{M}_j})$- or a $V(\Omega_{\mathfrak{M}_j})$-$V(\Omega_{\mathfrak{M}_i})$-path in $D[V(H_{\mathfrak{M}_i} \cup H_{\mathfrak{M}_j})] - T$.
    \end{enumerate}
\end{lemma}
\begin{proof}
    Suppose that $H_\mathfrak{M}$ contains an even dicycle and there does not exist a quarter-integral packing of $\ell$ even dicycles.
    We let $\mathfrak{M}_1, \ldots , \mathfrak{M}_h$ be the outlines we find in this setting via \Cref{lem:buildnewmaelstromsrefineddiamond} and let $d_1, \ldots , d_h \subseteq \Delta$ be the associated disks.
    This satisfies the first item.
    For the second item, we may assume that $h \geq 2$, as otherwise there is nothing to show.
    
    Let $i,j \in [h]$ be distinct, let $(H_{\mathfrak{M}_i}, \Omega_{\mathfrak{M}_i})$ be the $\mathfrak{M}_i$-society, let $(H_{\mathfrak{M}_j}, \Omega_{\mathfrak{M}_j})$ be the $\mathfrak{M}_j$-society, and let $R_i, R_j$ be the two cycles associated with $\Omega_{\mathfrak{M}_i}$ and $\Omega_{\mathfrak{M}_j}$ respectively.
    Furthermore, let $t_i, t_j$ be the respective degree of $\mathfrak{M}_i$ and $\mathfrak{M}_j$ and let $t_i'$ and $t_j'$ be their respective surplus degree, if they are refined diamond outlines.

    Suppose there exists a $V(\Omega_{\mathfrak{M}_i})$-$V(\Omega_{\mathfrak{M}_i})$-linkage $\mathcal{L}$ of order $2\ell$ in $D[V(H_{\mathfrak{M}_i}) \cup V(H_{\mathfrak{M}_j})]$.
    According to \Cref{lem:circleoutlinererouting}, \Cref{lem:diamondoutlinererouting}, and \Cref{lem:refineddiamondoutlinererouting}, we can find $\nicefrac{\theta}{2}$-rims for both outlines.
    If either of the two outlines is a refined diamond outline, we can reroute the endpoints of the paths in $\mathcal{L}$ onto the society of their underlying diamond outlines using \Cref{lem:detour}.
    Additionally, if either of the two rims are circle rims, the construction from the proof of \Cref{lem:buildnewmaelstromsrefineddiamond} allows us to find organised $\mathcal{C}_i$-pairs, if we need them.
    
    If $\mathfrak{M}_i$ is a (refined) diamond outline, we make use of the construction used in \Cref{lem:diamondlayershifting} to cover the (underlying) diamond $\nicefrac{\theta_i}{2}$-rim corresponding to $\mathfrak{M}_i$ in two circle $\nicefrac{\theta_i}{2}$-rims and thus ensure that at least $\ell$ of the paths in $\mathcal{L}$ have their endpoints on one of the two rims.
    Let $\mathcal{L}'$ be the selection of paths we just mentioned and let $\mathfrak{R}$ be the rim we selected.
    Should $\mathfrak{M}_i$ be a circle outline, we simply let $\mathfrak{R}$ be its $\nicefrac{\theta_i}{2}$-rim and let $\mathcal{L}' = \mathcal{L}$.
    We can then use \Cref{lem:detour} to route the endpoints of the paths in $\mathcal{L}'$ found on the society of the $\nicefrac{\theta_j}{2}$-rim corresponding to $\mathfrak{M}_j$ towards the society $\mathfrak{R}$ from the outside.
    Analogously to the proof of \Cref{lem:buildnewmaelstromscircleapexset}, we can find a planar grid of order $3\ell$ in $\mathfrak{R}$, together with enough perimeter jumps on said grid to let us construct a quarter-integral packing of even dicycles via \Cref{lem:gridplusjumpmakesoddbicycle}.

    Thus $\mathcal{L}$ cannot exist and Menger's Theorem gives us a set of vertices of size at most $2\ell$ that meets every $V(\Omega_{\mathfrak{M}_i})$-$V(\Omega_{\mathfrak{M}_i})$-path in $D[V(H_{\mathfrak{M}_i}) \cup V(H_{\mathfrak{M}_j})]$.
    We may repeat this argument for any pair of outlines amongst $\mathfrak{M}_1, \ldots , \mathfrak{M}_h$, yielding a set of vertices $T \subseteq \bigcup_{i \in [h]} V(H_{\mathfrak{M}_i})$ disrupting any paths between the societies of the outlines $\mathfrak{M}_1, \ldots , \mathfrak{M}_h$ within $D[\bigcup_{i \in [h]} H_{\mathfrak{M}_i}]$ such that $|T| \leq 4\ell^2$, as desired by the statement.
\end{proof}

Before we close this section, we note that \Cref{lem:buildnewmaelstromscircleapexset} and \Cref{lem:buildnewmaelstromsrefineddiamondapexset} yield apex sets that may hit a part of the newly constructed outline.
This needs to be cleaned up, which can be simply done with an application of \Cref{lem:flowalteredoutline} to each outline we find, giving us an opportunity to unify the two statements into one.
For this purpose and for later reuse, we define the function $\BuildOutlineNoArg \colon \N \rightarrow \N$, which gives a lower bound on the needed order of the outline for the application of the results in this section, and the function $\BuildOutlineApexNoArg \colon \N \rightarrow \N$, which gives an upper bound for the number of vertices we must delete:
\[ \BuildOutline{\ell} = 24\ell^{\ell+2} + 2 \ell^\ell ( \SC + \PlanarGrid{3\ell} + 6 \OddTransactionApex{t} ) \text{ and } \BuildOutlineApex{\ell} = \OddTransactionApex{t} + 4\ell^2 . \]
Further than this, we can now gather up the results of \Cref{sec:transaction} and \Cref{sec:buildoutline} into a single statement that says that, if we have an outline with depth above a certain threshold, then we can find a large odd transaction that splits this maelstrom into several new maelstroms with appropriate outlines.

\begin{corollary}\label{cor:buildnewmaelstroms}
    Let $\theta , d, p, \ell , g$ be positive integers such that $\theta$ is a multiple of $\BuildOutline{\ell}$, we have $p = 2\theta(g+1)$, and $d = \OddTransactionOrder{p}$, and let $D$ be a digraph with an odd decomposition $\delta$.
    Moreover, let $m$ be a maelstrom of $\delta$, together with a $\theta$-outline $\mathfrak{M}$, with $(H_\mathfrak{M},\Omega_\mathfrak{M})$ being the $\mathfrak{M}$-society, and, if $\mathfrak{M}$ is a diamond outline, a refined $\theta$-outline $\mathfrak{N} = ( H , \mathcal{C}, \mathfrak{E}, \mathfrak{B}, \mathfrak{V}, \mathfrak{D} )$ of $n \subseteq m$ with roughness $g$.

    Then, if $\mathfrak{M}$, respectively $\mathfrak{N}$, has depth at least $d$, there exists an odd transaction $\mathcal{P} = \{ P_1, \ldots , P_p \}$ on the $\mathfrak{M}$-society $(H_\mathfrak{M},\Omega_\mathfrak{M})$, if $\mathfrak{M}$ is a circle outline, or a calm, odd transaction $\mathcal{P}$ on the $\mathfrak{N}$-society $(H_\mathfrak{N},\Omega_\mathfrak{N})$, together with an apex set $A \subseteq V(H_\mathfrak{M})$ such that $|A| \leq \BuildOutlineApex{\ell}$.
    Let $H_\mathcal{P} = (H \cup S) - A$, where $S$ is the strip of $\mathcal{P}$ under $A$.

    Additionally, either $H_\mathfrak{M} - A$ is odd, there exists a quarter-integral packing of $\ell$ even dicycles in $H_\mathfrak{M}$, or there exists a positive integer $h$ with $h < \ell$ and
    \begin{itemize}
        \item If $\mathfrak{M}$ is a circle outline, there exist pairwise disjoint disks $d_1, \ldots , d_h \subseteq m$ such that for each $i \in [h]$ there exists a $\theta_i'$-outline $\mathfrak{M}_i = ( H_i , \mathcal{C}_i, \mathfrak{E}_i, \mathfrak{V}_i )$ for $d_i' \subseteq d_i$, with $\theta_i' \geq \nicefrac{\theta}{2^{h-1}}$, $H_i$ is a subgraph of $H_\mathcal{P}$, $\{ H_1, \ldots , H_h \}$ is half-integral, $d_i$ is bounded by the cycle corresponding to the $\mathfrak{M}_i$-society $(H_{\mathfrak{M}_i}, \Omega_{\mathfrak{M}_i})$, and $H_{\mathfrak{M}_i}$ contains an even dicycle, or

        \item If $\mathfrak{M}$ is a diamond outline, there exist pairwise disjoint disks $d_1, \ldots , d_h \subseteq n$ such that for each $i \in [h]$ there exists a $\theta_i$-outline $\mathfrak{M}_i = ( H_i , \mathcal{C}_i, \mathfrak{E}_i, \mathfrak{B}_i, \mathfrak{V}_i, \mathfrak{D}_i )$ (of roughness $g$) for $d_i' \subseteq d_i$, with $\theta_i \geq \nicefrac{\theta}{\ell^{h-1}}$, $H_i$ is a subgraph of $H_\mathcal{P}$, $\{ H_1, \ldots , H_h \}$ is third-integral 
        , $d_i$ is bounded by the cycle corresponding to the $\mathfrak{M}_i$-society $(H_{\mathfrak{M}_i}, \Omega_{\mathfrak{M}_i})$, the graph $H_{\mathfrak{M}_i}$ contains an even dicycle, and either $H_{\mathfrak{M}_i} \subset H_\mathfrak{M}$, if $\mathfrak{B}_i = \mathfrak{V}_i = \emptyset$, or $H_{\mathfrak{M}_i} \subset H_\mathfrak{N}$.
    \end{itemize}
\end{corollary}

It should be noted that there exists the possibility that the odd transaction we find only produces a single new maelstrom with an outline.
The problem with this is that we may need to delete vertices for this outline to be decomposable.
This detail will haunt us in \Cref{sec:localstructure} and presents one of the main problems we have to solve to finish the proof of the local structure theorem.


\section{Killing maelstroms}\label{sec:killmaelstroms}

Now that we have all the appropriate tools to decompose a digraph, we need to show that we can take advantage of those regions where we can no longer continue our approach of finding large transactions and splitting them up into new outlines.
In other words, we need to deal with maelstroms of bounded depth.

\subsection{Killing a circle maelstrom}

Following the example of the previous sections, we will develop the method for solving our problem in the relatively simple case of circle outlines and then get more technical to adapt this method for (refined) diamond outlines.
The principle idea here is that, since our maelstrom has bounded depth, we can split up the society of the maelstrom into several segments that each host an even dicycle, or something that lets us find an even dicycle, if we use the infrastructure of the outline.
If there are many such segments, we can find a large quarter-integral packing of even dicycles, or a half-integral packing in the case of circle outlines.
Otherwise, using Menger's theorem, we can delete the separators for the segments we have found and the resulting graph contains no even dicycles that touch the society of our maelstrom, allowing us to split off the remainder of the maelstrom from the rest of the graph.

For any digraph $D$, any society $(U, \Omega)$ with $U \subseteq D$, any segment $I \subseteq V(\Omega)$, and any subgraph $H \subseteq D$ we define $S_H(I)$ as the union of all strong components in $D - (V(\Omega) \setminus I)$ that contain some vertex of $I$.

\begin{lemma}\label{lem:killcirclemaelstrom}
    Let $\theta , t , d$ be integers with $\theta \geq 2 \lfloor \nicefrac{\SC}{2} \rfloor (t + 1)$, such that $\theta$ is even, and let $D$ be a digraph with an odd decomposition $\delta$.
    Moreover, let $m$ be a maelstrom of $\delta$, together with a circle $\theta$-outline $\mathfrak{M}$ with depth less than $d$, where $(H_\mathfrak{M},\Omega_\mathfrak{M})$ is the $\mathfrak{M}$-society, such that $D - (H_\mathfrak{M} - V(\Omega_\MM))$ is odd.

    Then either $H_\mathfrak{M}$ contains a half-integral packing of even dicycles, or there exists a set $S \subseteq V(H_\mathfrak{M})$, with $|S| \leq 4t d + t$, such that the union of the strong components of $D$ containing the vertices of $V(D - H_\mathfrak{M}) \cup V(\Omega_\mathfrak{M})$ is odd.
\end{lemma}
\begin{proof}
    Let $\mathfrak{R} = (H, \mathcal{C} = \{ C_1 , \ldots , C_{\nicefrac{\theta}{2}} \})$ be the circle $\nicefrac{\theta}{2}$-rim of $\mathfrak{M}$, which exists according to \Cref{lem:circleoutlinererouting}, and let $C = C_{\nicefrac{\theta}{2}}$ be the dicycle corresponding to $\Omega_\mathfrak{M}$.
    By definition $\Omega_\mathfrak{M}$ gives us a circular ordering of the vertices of $V(C)$.
    We fix one direction of traversal in this circular ordering for the remainder of this proof.

    For any cycle $C' \in \mathcal{C}$, we call a $C'$-path $P$ a \emph{parity breaking path} if $C' \cup P$ contains an even dicycle and $P$ is found on the $m$-disk of $C'$.
    Note that any even dicycle $Z$ in $D$ intersecting $C'$ such that $Z \cup C'$ has a planar decomposition contains at least one parity breaking path according to \Cref{lem:planartraceshifting}.

    Let $u \in V(C)$ be some vertex.
    We define the segments $I_1, \ldots , I_r \subseteq V(C)$ as follows:
    \begin{itemize}
        \item $I_1, \ldots , I_r$ partition $V(C)$ and occur on $V(C)$ in the order suggested by their indices.

        \item $u \in I_1$ and $I_1$ contains the successor of $u$.
        Furthermore, $I_1$ may only contain the predecessor of $u$ if $r = 1$.

        \item For each $i \in [r-1]$, the segment $I_i$ is minimal such that, if $T_i$ is the union of an $I_i$-$V(C - I_i)$- and a $V(C - I_i)$-$I_i$-separator in $U_i = H_\mathfrak{M} - \bigcup_{j=1}^{i-1}T_j$, then $S_{U_i - T_i}(I_i)$ contains an even dicycle, a parity breaking path for $C$ with its endpoints in $I_i$, or a cross on $C$ with its endpoints in $I_i$.
        (Note that these even dicycles and paths must intersect $\sigma(m)$.)

        \item The integer $r$ is chosen such that, if $T_r$ is the union of an $I_r$-$V(C - I_r)$- and a $V(C - I_r)$-$I_r$-separator in $U_r = H_\mathfrak{M} - \bigcup_{j=1}^{r-1}T_j$, it is the smallest integer such that $S_{U_r - T_r}(I_r)$ does not contain an even dicycle, a parity breaking path with its endpoints in $I_i$, or a cross on $C$ with its endpoints in $I_i$.
    \end{itemize}
    For each $i \in [r]$, let $H_i = S_{D - T_i'}(I_i)$, where $T_i' = \bigcup_{j=1}^i T_j$.
    Suppose that $r > t$ and note that $T_r$ does not contribute any usable dicycles or paths to allow us to find even dicycles.
    Let $J \subseteq [r]$ be the set of indices $i$ such that $H_i$ contains a pair of paths $P_i, Q_i$ forming a cross on $C_{\nicefrac{\theta}{2} - t}$.
    We choose $J$ in relation to $C_{\nicefrac{\theta}{2} - t}$, since parity breaking paths may also be non-planar, but might only form crosses on cycles in $\mathcal{C}$ that have a lower index and are thus closer to $m$.
    
    Next, let $I \subseteq [r] \setminus J$ be the set of indices such that $H_i$ contains a parity breaking path $P_i'$ for $C$.
    Note that for each $i$ in $[r] \setminus (I \cup J)$ there therefore exists an even dicycle $C_i' \subseteq H_i \cap H_\mathfrak{M}$.
    Clearly, if $|[r] \setminus (I \cup J)| > t$, we have found an integral packing of $t$ even dicycles.

    Suppose that $|I| \geq t$, then we can choose $t$ of these indices $i_1, \ldots, i_t \in I$ and find new even dicycles $C_{i_j}' \subseteq P_{i_j}' \cup C_{\nicefrac{\theta}{2} - (j-1)}$ for each $j \in [t]$ using \Cref{lem:planartraceshifting}.
    This yields a half-integral packing of $t$ even dicycles in $\mathfrak{M}$.
    Thus we instead suppose that $|I| < t$.

    Thus, for the remaining indices $i \in J$, we know that each $H_i$ contains a cross $P_i, Q_i$ on $C_{\nicefrac{\theta}{2} - t}$ drawn on the $m$-disk of this dicycle.
    In fact, each $P_i \cup Q_i$ contains a cross on $C_j$ for each $j \in [\nicefrac{\theta}{2} - t]$.
    We can now first find a half integral packing of $|I|$ even dicycles as discussed previously and then use \Cref{lem:circlenonplanarshifting} for the indices in $J$ on $\min(|J|, t)$ of the pairwise disjoint $\lfloor \nicefrac{\SC}{2} \rfloor$-rims induced by $\{ C_{\theta - \lfloor \nicefrac{\SC}{2} \rfloor \cdot j} , \ldots , C_{\theta - ( \lfloor \nicefrac{\SC}{2} \rfloor \cdot (j+1) - 1 )} \}$ for $j \in [|J|]$.
    If $|I| + |J| \geq t$, then this yields a half-integral packing of $t$ even dicycles.

    Finally, we may suppose that $|I| + |J| < t$.
    This then allows us to simply use our previous methods to find a half-integral packing of $|I| + |J|$ even dicycles and then add the remaining dicycles $C_i'$, which are found entirely in $H_i \cap H_\mathfrak{M}$ for each $i \in [r-1] \setminus (I \cup J)$.
    This packing is still a half-integral packing, as all these $C_i$ can only intersect the other even dicycles we found in the cycles from $\mathcal{C}$ that we used for shifting.

    Thus we may instead assume that there are fewer than $t$ usable segments and thus $r \leq t$.
    For each $i \in [r]$, let $u_i$ be the unique vertex in $I_i$ that does not have a successor in $I_i$.
    We define $S_i$ to be the union of $T_i \cup \{ u_i \}$ and the union of an $I_i \setminus \{ u_i \}$-$V(C - (I_i \setminus \{ u_i \}))$- and a $V(C - (I_i \setminus \{ u_i \}))$-$I_i \setminus \{ u_i \}$-separator in $H_\mathfrak{M} - ( T_i' \cup \bigcup_{j=1}^{i-1} S_j )$.
    If $V(I_r) = \{ u_r \}$, we let $S_r = T_r$.

    Altogether, we set $S = \bigcup_{i=1}^r S_i$ and note that each separator we have chosen has size at most $d$, as $\mathfrak{M}$ has depth less than $d$.
    Hence we have $|S| \leq 4t d + t$.

    Suppose that the union of strong components of $D - S$ that contain vertices of $V(D - H_\mathfrak{M}) \cup V(\Omega_\mathfrak{M})$ contains an even dicycle $C'$ containing some vertices of $V(\Omega_\mathfrak{M})$.
    Since $H$ itself is also odd, we know that $V(C') \cap V(\sigma(m)) \neq \emptyset$.
    Should $C \cup C'$ have a plane decomposition, then we can preserve the drawing of $C$ within this decomposition and use \Cref{lem:planartraceshifting} to prove that $H_i - S_i$, for some $i \in [r]$, contains a parity breaking path for $C$, contradicting the definition of our segments.
    Thus $C'$ must instead contain a cross on $C$ consisting of two paths that must be hosted in $H_i$ for some $i \in [r]$, but this again contradicts the definition of our segments.
    Therefore no such even dicycle exists.
\end{proof}

\subsection{Killing a (refined) diamond maelstrom}

We will now endeavour to replicate a version of \Cref{lem:killcirclemaelstrom} for (refined) diamond outlines with bounded depth, though we will search for a third-integral packing instead.
The main issue will be the fact that bounded depth does not actually prevent the existence of a transaction of large order on the maelstrom-society, as whirls on eddies cannot be excluded (see \Cref{def:calmtransaction}).
Once this issue is dealt with however, we can proceed in a largely analogous fashion to the proof of the previous lemma.
To facilitate this, we first briefly investigate what bounded depth for a refined diamond outline tells us about transactions on the diamond outline that we refined.

Let $\theta , g$ be integers and let $D$ be a digraph with an odd decomposition $\delta$.
Moreover, let $m$ be a maelstrom of $\delta$, together with a diamond $\theta$-outline $\mathfrak{M}$, with the $\mathfrak{M}$-society $(H_\mathfrak{M},\Omega_\mathfrak{M})$, and a refined diamond $\theta$-outline $\mathfrak{N}$ for $n \subseteq m$, with roughness $g$, depth less than $g$, and the $\mathfrak{N}$-society $(H_\mathfrak{N},\Omega_\mathfrak{N})$.
We let $C$ be the diamond corresponding to $\Omega_\mathfrak{M}$ and let $Q$ be the cycle corresponding to $\Omega_\mathfrak{N}$.
If $S_1, \ldots , S_h \subseteq Q$ are the $h$ distinct eddy segments of $Q$, then we note that, for each $i \in [h]$, there exists a directed path $S_i' \subseteq C$, with the same endpoints as $S_i$, whose internal vertices are drawn outside of the $n$-disk of $Q$.
We call $S_1', \ldots , S_h'$ the \emph{shorelines} of the eddy segments.
Given these definitions, we can now easily show that the shorelines behave much like the eddy segments with respect to the depth of their respective societies.

\begin{lemma}\label{lem:shorelinedepth}
    Let $\theta , g, d, h$ be integers with $d \geq (g+1) \theta $ and let $D$ be a digraph with an odd decomposition $\delta$.
    Moreover, let $m$ be a maelstrom of $\delta$, together with a diamond $\theta$-outline $\mathfrak{M}$, with the $\mathfrak{M}$-society $(H_\mathfrak{M},\Omega_\mathfrak{M})$, and a refined diamond $\theta$-outline $\mathfrak{N}$ for $n \subseteq m$, with roughness $g$, depth less than $d$, and the $\mathfrak{N}$-society $(H_\mathfrak{N},\Omega_\mathfrak{N})$.
    Furthermore, let $S_1, \ldots , S_h$ be the eddy segments of the cycle on $V(\Omega_\mathfrak{N})$ and let $S_1', \ldots , S_h'$ be their shorelines.

    Then, for any $i \in [h]$, there does not exists a transaction $\mathcal{T}$ on $(H_\mathfrak{M},\Omega_\mathfrak{M})$ of order $d$ that is a $V(S_i')$-$(V(\Omega_\mathfrak{M}) \setminus V(S_i'))$- or $(V(\Omega_\mathfrak{M}) \setminus V(S_i'))$-$V(S_i')$-linkage.
\end{lemma}
\begin{proof}
    We let $C$ be the diamond corresponding to $\Omega_\mathfrak{M}$ and let $Q$ be the cycle corresponding to $\Omega_\mathfrak{N}$.
    Let $i \in [h]$ and suppose there exists a transaction $\mathcal{T}$ on $(H_\mathfrak{M},\Omega_\mathfrak{M})$ of order $g$ and $\mathcal{T}$ is a $(V(\Omega_\mathfrak{M}) \setminus V(S_i'))$-$V(S_i')$-linkage.
    Note that $S_i' \cup S_i$ forms a cycle whose non-$n$-disk cannot contain the vertices of $C - S_i'$.
    Therefore all paths in $\mathcal{T}$ must intersect $S_i$ and in particular, this allows us to find a transaction $\mathcal{T}'$ of order $g$ that is a $(V(\Omega_\mathfrak{N}) \setminus V(S_i))$-$V(S_i)$-linkage, contradicting the bound on the depth of $\mathfrak{N}$, as $\mathcal{T}'$ is not whirly.
    Of course an analogous argument takes care of $V(S_i')$-$(V(\Omega_\mathfrak{M}) \setminus V(S_i'))$-linkages that are too large.
\end{proof}

This allows us to state the version of \Cref{lem:killcirclemaelstrom} for (refined) diamond maelstroms.
Note that if there are no eddies, then the statement of this lemma can be read exactly like \Cref{lem:killcirclemaelstrom}.
However, any eddy that hosts an even dicycle must be cut off and dealt with slightly differently, since we cannot properly kill the eddy or separate it from the rest of the graph with few vertices.
Instead, we must further decompose the graph at these points, but this will be done once we combine all of our arguments for the proof of the main theorem.

\begin{lemma}\label{lem:killdiamondmaelstrom}
    Let $\theta , t, g, d, h, k$ be integers with $\theta \geq 2 \SC (t + 1)$ and $d \geq (g+1)\theta$, such that $\theta$ even, and let $D$ be a digraph with an odd decomposition $\delta$.
    Moreover, let $m$ be a maelstrom of $\delta$, together with a diamond $\theta$-outline $\mathfrak{M}$, with the $\mathfrak{M}$-society $(H_\mathfrak{M},\Omega_\mathfrak{M})$, and a refined diamond $\theta$-outline $\mathfrak{N}$ for $n \subseteq m$, with roughness $g$, depth less than $d$, and the $\mathfrak{N}$-society $(H_\mathfrak{N},\Omega_\mathfrak{N})$, such that $D - (H_\mathfrak{N} - V(\Omega_\NN))$ is odd.
    Furthermore, let $S_1', \ldots , S_h'$ be the shorelines on the cycle $C$ corresponding to $V(\Omega_\mathfrak{M})$. 

    Then either $H_\mathfrak{M}$ contains a third-integral packing of even dicycles, or there exists a set $S \subseteq V(H_\mathfrak{M})$, with $|S| \leq 4t d + t$, and at most $t-1$ indices $i_1, \ldots , i_k \in [h]$ such that the union of strong components of $D - ( S \cup \bigcup_{j=1}^k V(S_{i_j}') )$ containing vertices $V(D - H_\mathfrak{M}) \cup (V(C) \setminus \bigcup_{j=1}^k V(S_{i_j}'))$ is odd.
\end{lemma}
\begin{proof}
    Let $\mathfrak{R}$ be the refined $\SC (t + 1)$-diamond rim with roughness $g$ corresponding to $\mathfrak{N}$, with $Q$ being the cycle corresponding to $\Omega_\mathfrak{N}$, which exists according to \Cref{lem:refineddiamondoutlinererouting}.
    Furthermore, let $\mathfrak{R}' = ( H, \mathcal{C}, \mathcal{E}, \mathfrak{B}, \mathfrak{D} )$ be the $\SC (t + 1)$-diamond rim corresponding to $\mathfrak{M}$, which exists according to \Cref{lem:diamondoutlinererouting}.
    We let $S_1, \ldots , S_h \subseteq Q$ be the $h$ distinct eddy segments of $Q$.
    
    From \Cref{lem:shorelinedepth} and \Cref{thm:directedlocalmenger}, we conclude that there exists a set $W_i$ that is the union of an $V(S_i')$-$V(C - S_i')$- and a $V(C - S_i')$-$V(S_i')$-separator in $H_\mathfrak{M}$, such that $|W_i| \leq 2d$, for each $i \in [h]$.
    For the remainder of the proof we fix some direction of transversal of $\Omega_{\mathfrak{M}}$.
    We let the notion of a parity breaking path for $C$ be defined as it is in the proof of \Cref{lem:killcirclemaelstrom}.

    If there exist $t$ or more graphs amongst $S_{H_\mathfrak{M} - W_1}(V(S_1')), \ldots , S_{H_\mathfrak{M} - W_h}(V(S_h'))$ that contain an even dicycle, a parity breaking path on $C$, or a cross on $C$, then we can find a third-integral packing of $t$ even dicycles in an analogous fashion to the methods used in the proof of \Cref{lem:killcirclemaelstrom} via the use of \Cref{lem:planartraceshifting} and \Cref{lem:diamondnonplanarshifting}.

    Thus we may suppose that there exists some set $J \subseteq [h]$ comprised of all indices such that $|J| \leq t - 1$ and for each $i \in J$ the graph $S_{H_\mathfrak{M} - W_i}(V(S_i'))$ contains an even dicycle, a parity breaking path on $C$, or a cross on $C$.
    Note that for any $i \in [h] \setminus J$ the graph $S_{H_\mathfrak{M}}(V(S_i'))$ is thus odd.

    We can now build segments $I_1, \ldots , I_r \subseteq V(C) \setminus ( \bigcup_{i \in J} V(S_i') \setminus V(Q) )$, similar to the proof of \Cref{lem:killcirclemaelstrom}, by respecting the following requirements.
    Let $u \in V(C) \cap V(Q)$ be some vertex.
    \begin{itemize}
        \item $I_1, \ldots , I_r$ partition $V(C)$ and occur on $V(C)$ in the order suggested by their indices.

        \item $u \in I_1$ and $I_1$ contains the successor of $u$.
        Furthermore, $I_1$ may only contain the predecessor of $u$ on $C$ if $r = 1$.

        \item For any $i \in [r]$ and $j \in [h]$, we have $V(S_j') \subseteq I_i$ if and only if $V(S_j') \cap I_i \neq \emptyset$.

        \item For each $i \in [r-1]$, the segment $I_i$ is minimal such that, if $T_i$ is the union of an $I_i$-$V(C - I_i)$- and a $V(C - I_i)$-$I_i$-separator in $U_i = H_\mathfrak{M} - \bigcup_{j=1}^{i-1}T_j$, then $S_{U_i - T_i}(I_i)$ contains an even dicycle, a parity breaking path for $C$ with its endpoints in $I_i$, or a cross on $C$ with its endpoints in $I_i$.
        (Note that these even dicycles and paths must intersect $\sigma(m)$.)

        \item The integer $r$ is chosen such that, if $T_r$ is the union of an $I_r$-$V(C - I_r)$- and a $V(C - I_r)$-$I_r$-separator in $U_r = H_\mathfrak{M} - \bigcup_{j=1}^{r-1}T_j$, it is the smallest integer such that $S_{U_r - T_r}(I_r)$ does not contain an even dicycle, a parity breaking path with its endpoints in $I_i$, or a cross on $C$ with its endpoints in $I_i$.
    \end{itemize}
    Now, if $|J| + r \geq t$, we can proceed analogously to the proof of \Cref{lem:killcirclemaelstrom} to find a third-integral packing, using \Cref{lem:planartraceshifting} and \Cref{lem:diamondnonplanarshifting}, and otherwise we can define the desired set $S$ of size at most $4t d + t$, by combining those sets $W_i$ for which $i \in J$ and sets defined analogously to the sets $S_1, \ldots , S_r$ in the proof of \Cref{lem:killcirclemaelstrom}.
\end{proof}

\section{Local structure theorem}\label{sec:localstructure}

In this section we present a structure theorem for a digraph with respect to a large odd wall.
It should be mentioned that this structure theorem does not entirely suffice to prove \Cref{thm:mainthm1}, since \Cref{thm:localstructure} incurs a set of vertices that need to be deleted.
This may split the graph into several distinct strong components and only the strong component containing whatever is left of the wall will be decomposed here.
A resolution to the problem, in the form of a global structure theorem, is given in \Cref{sec:globalstructure}.

We first prove a simple theorem that lets us separate the two sides of an odd wall from each other with a few vertices, which establishes one more pair of functions $\SeparatingWallNoArg \colon \N \rightarrow \N$ and $\SeparatingWallApexNoArg \colon \N \rightarrow \N$ that are defined as follows:
\[ \SeparatingWall{t} = 4(t+1)\SC + 3 \text{ and } \SeparatingWallApex{t} = 4t-2 .\]
Further, we say that an odd wall $W = ( Q_1, \ldots , Q_k, \hat{P}_1, \ldots , \hat{P}_k )$ is \emph{separating} in $D$ under a set of vertices $A$ if there exist undirected separations $(X_1, Y_1), (X_k, Y_k)$ in $K - A$, where $K$ is the $W$-component of $D$ under $A$, with $X_i \cap Y_i = A \cup V(Q_i)$ and every vertex in $Y_i$ reaches or is reachable by a vertex of $V(W)$ for both $i \in \{ 1, k \}$.

\begin{lemma}\label{lem:jumpsoverawall}
    Let $t, k$ be integers such that $k \geq 4(t+1)\SC + 3$, let $D$ be a digraph, $A \subseteq V(D)$ be a set of vertices and let $W = ( Q_1, \ldots , Q_k, \hat{P}_1, \ldots , \hat{P}_k )$ be an odd $k$-wall under $A$ in $D$, where $Q_1$ and $Q_k$ are the two dicycles forming the perimeter of $W$.
    
    Then either $D - A$ contains a half-integral packing of $t$ even dicycles in $D$, or there exists a set $S$ of at most $4t-2$ vertices in $V(D - A)$ such that $W$ is separating in $D$ under $(A \cup S)$.
\end{lemma}
\begin{proof}
    We suppose that there does not exist a half-integral packing of $t$ even dicycles in $D$.
    Since $W$ is an odd wall under $A$ in $D$, by definition there exists a separation $(X,Y)$ of $D$ such that $X \cap Y = A \cup \Perimeter{W}$, $W - \Perimeter{W} \subseteq Y$, and every vertex in $Y$ reaches a vertex of $W - \Perimeter{W}$ or is reachable from it.
    Furthermore, the cylindrical society $(D[Y], \Omega_1, \Omega_2)$ with $V(\Omega_i) = V(Q_i)$ for both $i \in \{ 1, k \}$ has an odd rendition $\rho$ in the disk.

    Using these facts, we note that in $D[Y] - V(Q_{2(t+1)\SC+1} \cup Q_{2(t+1)\SC+2} \cup Q_{2(t+1)\SC+3})$ the dicycles $Q_1, \ldots , Q_{2(t+1)\SC}$ and the dicycles $Q_{2(t+1)\SC+4}, \ldots, Q_k$ lie in different strong components, due to $\rho$ being an odd rendition.
    In particular, there exist no $V(Q_{2(t+1)\SC})$-$V(Q_{2(t+1)\SC+4})$- and no $V(Q_{2(t+1)\SC+4})$-$V(Q_{2(t+1)\SC})$-paths in $D[Y] - V(Q_{2(t+1)\SC+1} \cup Q_{2(t+1)\SC+2} \cup Q_{2(t+1)\SC+3})$.

    Let $D' = D - V(Q_{2(t+1)\SC+1} \cup Q_{2(t+1)\SC+2} \cup Q_{2(t+1)\SC+3})$ and suppose that there exist $t$ disjoint paths $\mathcal{R} = \{ R_1, \ldots , R_{2t} \}$ in $D'$ with their tails in $V(Q_{2(t+1)\SC})$ and their heads in $V(Q_{2(t+1)\SC+4})$.
    Due to the existence of the separation $(X,Y)$, the rendition $\rho$, and our prior observations, we know that, starting from their tails, the paths in $\mathcal{R}$ traverse through the dicycles $Q_1, \ldots , Q_{2(t+1)\SC}$ in ascending order and traverse the dicycles $Q_{2(t+1)\SC+4}, \ldots, Q_k$ in descending order.
    (Though the paths may turn back whilst the are travelling through the dicycles, they still do so in an orderly way, only ever visiting a given dicycle $Q_i$ if their last intersection was with $Q_{i+1}$ or $Q_{i-1}$.)
    
    It is possible that within the paths of $\mathcal{R}$  we can find subpaths that form crosses on some $Q_i$ within $W$.
    Clearly, if this is true for $Q_i$ with $\SC \leq i \leq 2(t+1)\SC$, then there exists a cross on $Q_\SC$ and we may use \Cref{lem:circlenonplanarshifting} to find a dicycle within the first $\SC$ dicycles.
    This at most involves two paths of $\mathcal{R}$ and since we chose at least $(t+1)\SC$ dicycles between $Q_{2(t+1)\SC}$, respectively $Q_{2(t+1)\SC+4}$, and the perimeter of $W$, we can transform $t$ crosses on $Q_{t\SC}$ into a half-integral packing of $t$ even dicycles.
    Thus we may suppose that there is no cross on $Q_{t\SC}$, which means that with respect to the dicycles $Q_{t\SC} , \ldots , Q_{2(t+1)\SC}$ the paths $R_1, \ldots , R_{2t}$ behave like planar paths.

    Recall that each $\hat{P}_i$ consists of a pair of paths $P_i^1$ and $P_i^2$, with $P_i^1$ traversing the dicycles $Q_1, \ldots , Q_k$ in ascending order of their indices and $P_i^2$ traversing them in descending order.
    Since we just observed that we can suppose that the dicycles $Q_{t\SC} , \ldots , Q_{2(t+1)\SC}$ and the paths in $\mathcal{R}$ behave in a planar way, we can apply \Cref{lem:detour} to the paths in $\mathcal{R}$, the dicycles $Q_{t\SC} , \ldots , Q_{2(t+1)\SC}$, and the paths $v_1P_1^2, \ldots , v_tP_{2t}^2$, where $v_i$ is the last vertex in the intersection of $V(P_i^2)$ and $V(Q_{(t+1)\SC})$, when traversing $P_i^2$ starting from its tail, for each $i \in [t]$.
    This allows us to reroute the paths in $\mathcal{R}$ such that $R_i$ starts on $v_i$ for each $i \in [t]$.
    
    Via analogous arguments, we can also reroute the paths in $\mathcal{R}$ to have their heads join up with the paths $P_1^2, \ldots , P_{2}t^2$ at their respective first vertex on $Q_{3(t+1)\SC + 4}$.
    Thus we have found a collection of at least $t$ perimeter jumps on a wall $W'$ of order at least $t\SC$, which can be extended to any particular subwall of $W'$, since they intersect with $v_1P_1^2, \ldots , v_tP_{2t}^2$.
    As $\SC \geq 3$, we can therefore use \Cref{lem:jumpsoverawall} to find a half-integral packing of $t$ even dicycles in $D$, contradicting our assumptions.
    Thus \Cref{thm:directedlocalmenger} gives us a separator $S_1$ for all $V(Q_{2(t+1)\SC})$-$V(Q_{2(t+1)\SC+4})$-paths in $D'$ of size at most $2t -1$ and via an analogous argument we can find such a separator $S_2$ for all $V(Q_{2(t+1)\SC+4})$-$V(Q_{2(t+1)\SC})$-paths in $D'$, also of size at most $2t -1$.
    If we now consider $D - (S_1 \cup S_2)$, it becomes easy to observe that the desired separations must exist, using the existence of the separation $(X,Y)$ for $D$ and the rendition $\rho$.
\end{proof}

\subsection{Splitting a maelstrom}\label{subsec:splittingmaelstrom}

We now show that we can resolve the situation mentioned at the end of \Cref{sec:buildoutline}.
The solution will be to prove that we either split the maelstrom into several new maelstroms hosting even dicycles, or we can repeatedly refine the maelstrom to contain less and less of the graph, whilst giving a constant bound on the number of vertices that we need to delete to give the remaining maelstrom an odd outline throughout the entire process.

Before we begin, we establish the following three functions $\InductionStepApexNoArg \colon \N \rightarrow \N$, $\InductionStepOutlineNoArg \colon \N \rightarrow \N$, and $\InductionStepOrderNoArg \colon \N \times \N \rightarrow \N$, which are intended to capture both the size of transaction and outline we search for throughout, and the maximum size of the apex set we need.
For any integers $t , g \in \N$ we define these functions as follows:
\begin{align*}
    \InductionStepApex{t}       \coloneqq \ & 2 ( \BuildOutlineApex{t} + 4(\FindOddWallApex{t} + \SeparatingWallApex{t}) ), \\
    \InductionStepOutline{t}    \coloneqq \ & 4t^t ( \InductionStepOutline{t-1} + 18 \InductionStepApex{t} + 3 \PlanarGrid{ \FindOddWallOrder{ \SeparatingWall{t} , t } } ), \text{ and } \\
    \InductionStepOrder{t,g}    \coloneqq \ & 3 ( g + 1 ) ( \OddTransactionOrder{\InductionStepOutline{t} + 3 \PlanarGrid{ \FindOddWallOrder{ \SeparatingWall{t} , t } } + 18\InductionStepApex{t} + t(2\InductionStepApex{t} + 2)} + 2t - 1.
\end{align*}

Next we give the definition of the formal objects needed for the induction step of proof of the local structure theorem.
Our goal is to split up a maelstrom into several smaller maelstroms, each with their own outlines.
However over the course of this process, we may encounter a maelstrom that refuses to split up.
In this case we must prove that we can continuously keep shrinking the maelstrom whilst keeping the apex set used to perform this procedure below a fixed bound, no matter how often we iterate the process and fail to split the maelstrom.

\begin{definition}[Overflow in an outline]\label{def:overflow}
    Let $\theta , t, d, g, k, s$ be integers such that $\theta = \InductionStepOutline{t}$, $p = \InductionStepOrder{t,g} - (2t - 1)$, and $q = \nicefrac{p}{3}$.
    Moreover, let $D$ be a digraph with an odd decomposition $\delta_0$ and let $m_0$ be a maelstrom of $\delta_0$, together with a $\theta$-outline $\mathfrak{M}_0  = ( H , \mathcal{C}, \mathfrak{E}, \mathfrak{V} )$ of degree $d$ and potentially a refined diamond $\theta$-outline $\mathfrak{N}$ for $n_0 \subseteq m_0$ of roughness $g$ and surplus degree $k$, if $\mathfrak{M}_0$ is a diamond outline, with $(H_0,\Omega_0)$ being the $\mathfrak{M}_0$-society, respectively the $\mathfrak{N}_0$-society.
    We set $m_0 \coloneqq n_0$, if $n_0$ exists, and let $d_0$ be the $m_0$-disk of the cycle corresponding to $\Omega_0$ in $\delta_0$. 

    We say that $\mathfrak{M}$, respectively $\mathfrak{N}$, has \emph{overflow} $s$ if there exist linkages $\mathcal{Q}_1, \ldots , \mathcal{Q}_s$ and for each $i \in [s]$ we have:
    \begin{enumerate}
            \item $\mathcal{Q}_i$ is a planar transaction of order $p$ on $(H_{i-1},\Omega_{i-1})$, which is required to be a calm transaction if $\mathfrak{M}_{i-1}$ is a refined diamond outline,

            \item there exists a set $A_i \subseteq V(H_0)$ with $|A_i| \leq \nicefrac{\InductionStepApex{t}}{2}$, such that there exists an odd decomposition $\delta_i$ of $H_0 - A_i$ that contains a unique maelstrom $m_i$ together with a $(\theta - 6\InductionStepApex{t})$-outline $\mathfrak{M}_i$ (of roughness $g$) such that $H_i$ contains an even dicycle, with $(H_i,\Omega_i)$ being the $\mathfrak{M}_i$-society,

            \item for each $i \in [s]$ the disk $d_i$ is the $m_i$-disk of the cycle corresponding to $\Omega_{\mathfrak{M}_i}$ in $\delta_i$, and the transaction $\mathcal{Q}_i = \{ Q^i_1, \ldots , Q^i_p \}$ is indexed such that for all $j \in [2, p-1]$ the trace of $Q^i_j$ separates the trace of $Q^i_{j-1}$ and $Q^i_{j+1}$ in $d_{i-1}$, and each $\mathcal{Q}_i$ is partitioned into $\mathcal{Q}^i_1 = \{ Q^i_1, \ldots , Q^i_q \}$, $\mathcal{Q}_2^i = \{ Q^i_{q + 1}, \ldots , Q^i_{2q} \}$, and $\mathcal{Q}_3^i = \{ Q^i_{2q + 1}, \ldots , Q^i_{3q} \}$,

            \item for each $i \in [s]$ there exists an integer $s_i \in [0, t-1]$ and $s_i$ crossing pairs $(L_1^i, R_1^i), \ldots , (L_{s_i}^i, R_{s_i}^i)$ on the $\mathfrak{M}_{i-1}$-society\footnote{We explicitly do not care if these crossing pairs are hit by any apex set and in fact expect them to be hit.}, such that all paths in the crossing pairs are pairwise disjoint, for each $j \in [s_i]$ the path $L_j^i$ is disjoint from the paths in $\{ Q^i_1, \ldots , Q^i_{p - (\ell - \nicefrac{1}{2}) (2\InductionStepApex{t} + 2)} \} $ and from $H_{\mathfrak{M}_i}$, and $R_j^i = Q^i_{p - (j - 1)(2\InductionStepApex{t} + 2)}$,

            \item for all $i \in [s]$, we have $V(\Omega_i) \subseteq V(\Omega_0) \cup V( \bigcup_{j = 1}^s Q^j_{ p - s_j(2\InductionStepApex{t} + 2)} )$, and
            
            \item there exists an odd rendition $\rho_i$ for $( H_0 - (V(H_i - V(\Omega_i)) \cup A_i), \Omega_0, \Omega_i )$.
    \end{enumerate}
    We call $\mathcal{Q}_1, \ldots , \mathcal{Q}_s$ the \emph{overflow-linkages} associated with $\mathfrak{M}$, respectively $\mathfrak{N}$, the value $s_i$ is called the \emph{spill-over} of $\mathcal{Q}_i$ for each $i \in [s]$, and we call $\mathfrak{M}_1, \ldots , \mathfrak{M}_s$ the \emph{overflow-outlines} associated with $\mathfrak{M}$, respectively $\mathfrak{N}$.
\end{definition}

The following lemma takes care of all cases that might crop up while splitting a maelstrom.
One case worth remarking on is the one in which we manage to somehow find an even dicycle that lies outside of the new maelstrom and its outline, but it explicitly touches the apex set we use to construct that outline.
This is still useful to us, because we are able to either use this even dicycle or kill all even dicycles associated with the outside of this new outline by removing this apex set of bounded size.
Thus we can reduce the number of even dicycles we are hunting for whilst at the same time only incurring a small number of new vertices that need to be deleted.

\begin{lemma}\label{lem:localstructureinductionstepgeneral}
    Let $\theta , t, g, s$ be positive integers such that $\theta = \InductionStepOutline{t}$.
    Moreover, let $D$ be a digraph with an odd decomposition $\delta$ and let $m$ be a maelstrom of $\delta$, together with a $\theta$-outline $\mathfrak{M}  = ( H , \mathcal{C}, \mathfrak{E}, \mathfrak{V} )$, with $(H_\mathfrak{M},\Omega_\mathfrak{M})$ being the $\mathfrak{M}$-society and a segregated $\mathcal{C}$-pair of order $\nicefrac{\theta}{2}$ if $\mathfrak{M}$ is a circle outline with degree zero, and if $\mathfrak{M}$ is a diamond outline, a refined diamond $\theta$-outline $\mathfrak{N}$ for $n \subseteq m$ of roughness $g$.
    Further, let $\mathcal{Q}_1, \ldots , \mathcal{Q}_s$ be the overflow-linkages associated with $\mathfrak{M}$, respectively $\mathfrak{N}$, and let $\mathfrak{M}_1, \ldots , \mathfrak{M}_s$ be the overflow-outlines associated with $\mathfrak{M}$, respectively $\mathfrak{N}$, where $( H_{\mathfrak{M}_s}, \Omega_{\mathfrak{M}_s} )$ is the $\mathfrak{M}_s$-society.

    Then one of the following holds
    \begin{enumerate}
        \item there exists a quarter-integral packing of $t$ even dicycles in $D$,\label{item:ISpacking}

        \item $\mathfrak{M}_i$ is a $(\theta - 6\InductionStepApex{t})$-outline (of roughness $g$) with depth less than $\InductionStepOrder{t,g}$,\label{item:ISdepth}

        \item there exists a set $A \subseteq V(H_\mathfrak{M})$ with $|A| \leq \InductionStepApex{t}$ such that $H_\mathfrak{M} - A$ is odd,\label{item:ISodd}
        
        \item there exists an $h \in \N$ with $1 < h < t$, pairwise disjoint disks $d_1^*, \ldots , d_h^* \subseteq m$, and a set $S \subseteq V(H_\mathfrak{M})$ with $|S| \leq \InductionStepApex{t}$, such that for each $i \in [h]$ there exists a $\theta_i$-outline $\mathfrak{M}_i^* = ( H_i^* , \mathcal{C}_i^*, \mathfrak{E}_i^*, \mathfrak{B}_i^*, \mathfrak{V}_i^*, \mathfrak{D}_i^* )$ (with roughness $g$) for $d_i^{**} \subseteq d_i^*$ within $H_\mathfrak{M} - S$ with $\theta_i \geq \InductionStepOutline{t - (h-1)}$, the set $\{ H_1^*, \ldots , H_h^* \}$ is third-integral, the trace of $d_i^*$ corresponds to the $\mathfrak{M}_i^*$-society $(H_{\mathfrak{M}_i^*}, \Omega_{\mathfrak{M}_i^*})$, and the graph $H_{\mathfrak{M}_i^*}$ contains an even dicycle,\label{item:ISsplitmaelstrom}
        
        \item there exists a set $A \subseteq V(H_\mathfrak{M})$ with $|A| \leq \InductionStepApex{t}$ such that $H_\mathfrak{M} - A$ has an odd decomposition $\delta'$ with a unique maelstrom $m'$ and an associated $\InductionStepOutline{t-1}$-outline $\mathfrak{M}'$ (of roughness $g$), where $(H_{\mathfrak{M}'},\Omega_{\mathfrak{M}'})$ is the $\mathfrak{M}'$-society, $H_{\mathfrak{M}'}$ contains an even dicycle, there exists an odd rendition $\rho'$ of $( H_\mathfrak{M} - (V(H_{\mathfrak{M}'} - V(\Omega_{\mathfrak{M}'})) \cup A), \Omega_\mathfrak{M}, \Omega_{\mathfrak{M}'} )$, and $H_\mathfrak{M} - V(H_{\mathfrak{M}'} - V(\Omega_{\mathfrak{M}'}))$ contains an even dicycle (which must contain vertices of $A$), or\label{item:ISsplitoffdicycle}

        \item there exists a planar transaction $\mathcal{Q}_{s+1}$ on $( H_{\mathfrak{M}_s}, \Omega_{\mathfrak{M}_s} )$, a set of vertices $A_{s+1} \subseteq V(H_\mathfrak{M})$, and a $(\theta - 6\InductionStepApex{t})$-outline $\mathfrak{M}_{s+1}$ (of roughness $g$) in $H_\mathfrak{M} - A_{s+1}$, which together verify that $\mathfrak{M}$ has overflow $s + 1$.\label{item:ISoverflow}
    \end{enumerate}
\end{lemma}
\begin{proof}
    We let $p$, $q$, $\delta_i$, $m_i$, $d_i$, $A_i$, $\mathfrak{M}_i$, and the transactions $\mathcal{Q}^i_1, \mathcal{Q}^i_2, \mathcal{Q}^i_3$ for each $i \in [0, s]$ be defined as in \Cref{def:overflow}, with some objects undefined for $i = 0$, and suppose that the first two items of our statement do not hold.
    According to \Cref{cor:manycyclesorplanartransaction} there must therefore exist a planar transaction $\mathcal{Q}_{s+1} = \{ Q^{s+1}_1, \ldots , Q^{s+1}_p \}$ of order $p$ on $( H_{\mathfrak{M}_s}, \Omega_{\mathfrak{M}_s} )$ contained in a transaction of order $\InductionStepOrder{t,g}$ on $( H_{\mathfrak{M}_s}, \Omega_{\mathfrak{M}_s} )$.

    We choose indices such that for all $j \in [2, p-1]$ the trace of $Q^{s+1}_j$ separates the trace of $Q^{s+1}_{j-1}$ and $Q^{s+1}_{j+1}$ in $d_s$.
    We also partition $\mathcal{Q}_{s+1}$ into $\mathcal{Q}^{s+1}_1 = \{ Q^{s+1}_1, \ldots , Q^{s+1}_q \}$, $\mathcal{Q}_2^{s+1} = \{ Q^{s+1}_{q + 1}, \ldots , Q^{s+1}_{2q} \}$, and $\mathcal{Q}_3^{s+1} = \{ Q^{s+1}_{2q + 1}, \ldots , Q^{s+1}_{3q} \}$.

    Suppose that the third, fourth, and fifth item of our statement also do not hold, leaving us to prove \cref{item:ISoverflow}.
    To slightly simplify the proof, we neglect the case in which $\mathfrak{M}$ is a circle outline with degree zero, as it does not differ from the other cases, aside from forcing us to bring up the segregated $\mathcal{C}$-pair in a few places.
    
    We first show that we can find an outline for a special type of maelstrom $w_{s+1}$ that is itself not a good candidate for $m_{s+1}$, but provides a useful reference point.
    Let $A_0 = \emptyset$.
    
    \begin{claim}\label{claim:findsmallmaelstrom}
        For each $i \in [s+1]$, there exists an odd transaction $\mathcal{Q}^{i*}_2 \subseteq \mathcal{Q}^i_2$ of order $(g+1)\theta$, a set of vertices $A_i'$ of size at most $\BuildOutlineApex{t}$, an odd decomposition $\delta_i'$ of $H - (A_{i-1} \cup A_i')$ with a unique maelstrom $w_i$ that comes equipped with a $(\theta - 6\InductionStepApex{t})$-outline $\mathfrak{W}_i$ (of roughness $g$), and an odd rendition of $\rho'_i$ of $( H - (V(H_{\mathfrak{W}_i} - V(\Omega_{\mathfrak{W}_i})) \cup A_{i-1} \cup A_i'), \Omega_\mathfrak{M}, \Omega_{\mathfrak{W}_i} )$, where $( H_{\mathfrak{W}_i}, \Omega_{\mathfrak{W}_i} )$ is the $\mathfrak{W}_i$-society.
    \end{claim}
    \emph{Proof of \Cref{claim:findsmallmaelstrom}:}
    Since $q \geq (g+1)\OddTransactionOrder{\InductionStepOutline{t}}$, we can then apply \Cref{cor:buildnewmaelstroms}, whilst using $\mathcal{Q}^i_2$ as the witness surpassing the bound on the depth of $\mathfrak{M}_{i-1}$, to $\mathfrak{M}_{i-1}$ and $\mathcal{Q}^{i*}_2$ to find an apex set $A_i'$ of order at most $\BuildOutlineApex{t}$.
    We note the following about our chosen functions:
    \[ q \geq \OddTransactionOrder{\InductionStepOutline{t}}, \ \InductionStepApex{t} \geq \nicefrac{\InductionStepApex{t}}{2} + \BuildOutlineApex{t}, \text{ and } \nicefrac{\theta}{t^{h-1}} - 6\InductionStepApex{t} \geq \InductionStepOutline{t - (h-1)} . \]
    Thus, as \cref{item:ISsplitmaelstrom} does not hold, \Cref{cor:buildnewmaelstroms} holds with $h = 1$, proving our claim.
	\hfill$\blacksquare$

    Note that $w_i$ will lie on one of the two sides of $\mathcal{Q}^i_2$ in $d_{i-1}$.
    Either a majority of the paths in $\mathcal{Q}^i_1$ will lie in $H_{\mathfrak{W}_i}$ or a majority of the paths in $\mathcal{Q}^i_3$ will lie in $H_{\mathfrak{W}_i}$.
    We speak of the majority here, since the apex set $A_i'$ may in fact separate subpaths of any path $P$ in $\mathcal{Q}^i$ from its endpoints.
    Suppose that $P \in \mathcal{Q}^i_3$ and some subpath $P' \subseteq P$ is separated from the endpoints of $P$ by $A_i'$, then, even if the majority of the paths in $\mathcal{Q}^i_1$ lie in $H_{\mathfrak{W}_i}$ and thus the even maelstrom intuitively lies on the side of $\mathcal{Q}^i_2$ that belongs to $\mathcal{Q}^i_1$, $P'$ may still be so well connected into $H_{\mathfrak{W}_i}$ that it is actually found in $H_{\mathfrak{W}_i}$.
    We mention this explicitly since it may at the moment seem like a minor corner case, but this problem is responsible for most of the complexity of the remainder of this proof.

    In the following we will assume w.l.o.g.\ that a majority of the paths in $\mathcal{Q}^{s+1}_1$ lie in $H_{\mathfrak{W}_{s+1}}$ and we let the other $\mathcal{Q}_i$ be indexed consistently such that the majority of the paths in $\mathcal{Q}^i_1$ lie in $H_{\mathfrak{W}_i}$.
    To be precise, we say that the majority of the paths in $\mathcal{Q}^i_1$ lie in $H_{\mathfrak{W}_i}$ if at least $q - \InductionStepApex{t}$ paths of $\mathcal{Q}^i_1$ are completely contained in $H_{\mathfrak{W}_i}$, which is uniquely determined, since
    \[ q \geq 2 \InductionStepApex{t} + 1 . \] 

    Let $\ell = \nicefrac{\theta - 6\InductionStepApex{t}}{2}$, and let $\mathfrak{R}_i = ( U_i , \mathcal{R}_i, \mathcal{F}_i, \mathfrak{A}_i, \mathfrak{D}_i )$ be the $\ell$-rim of roughness $g$ for $\mathfrak{W}_i$, for $i \in [s+1]$.
    Further, let $\mathcal{R}_i = \{ R_1^i, \ldots , R_\ell^i \}$ and, if $\mathfrak{R}_i$ is a diamond rim, let $\mathcal{F}_i = \{ F_1^i, \ldots , F_\ell^i \}$, otherwise we let $\mathcal{F}_i = \emptyset$.
    If $\mathfrak{R}_i$ is a refined diamond rim of eddy degree $k$, we let $\mathfrak{A}_i = \{ \mathcal{A}^i_1, \ldots , \mathcal{A}^i_k \}$ with $\mathcal{A}^i_j = \{ A_1^{i,j}, \ldots , A_{(g+1)\ell}^{i,j} \}$ for each $j \in [k]$.
    Again giving necessary bounds for our functions, we let $\ell' = \nicefrac{\ell}{2}$ and note that
    \[ \ell' \geq 3 \PlanarGrid{ \FindOddWallOrder{ \SeparatingWall{t} , t } } + 18\InductionStepApex{t} . \]

    Our next goal will be to find one or several large odd walls in the rim $\mathfrak{R}_{s+1}$.
    We split this into three cases, one for circle rims, one for diamond rims, and one for refined diamond rims, each yielding an increasingly large number of walls.

    \textbf{Finding an odd wall in a circle rim:}
    First, we note that, if $\mathfrak{R}_i$ is a circle rim, we can easily find a segregated $\mathcal{R}_i$-pair of order $\ell$, since the initial outline either has positive degree or a segregated $\mathcal{C}$-pair, and thus we can assume that $\mathcal{R}_i$ comes equipped with such a pair $(\mathcal{L}_1^i, \mathcal{L}_2^i)$.    
    Thus, we can use the dicycles $R_{\ell' + \nicefrac{\ell'}{3} + 1}^i, \ldots , R_{\ell' + \nicefrac{2\ell'}{3}}^i$, the segregated $\mathcal{R}_i$-pair, \Cref{cor:findoddwall}, and \Cref{lem:jumpsoverawall} to find a separating odd wall $W_1^i$ of order at least 3 within $H_\mathfrak{M}$, having returned all previously deleted vertices to the graph, under an associated apex set $A_1^i$ of order $\SeparatingWallApex{t}$.
    Let $(X_1^i, Y_1^i), (X_2^i, Y_2^i)$ be the two separations that witness that $W_1^i$ is separating, such that $X_j^i \cap \sigma_{\delta_i'}(w_i) \neq \emptyset$ and thus $Y_j^i \cap \sigma_{\delta_i'}(w_i) = \emptyset$ for both $j \in [2]$ and $W_1^i \subseteq Y_2^i$.

    \textbf{Finding an odd wall in a diamond rim:}
    Within the dicycles $R_{\ell' + \nicefrac{\ell'}{3} + 1}^i, \ldots , R_{\ell' + \nicefrac{2\ell'}{3}}^i$, we can still find $W_1^i$ and the associated set $A_1^i$.
    Additionally, we can find a set of dicycles in $R_{\ell' + \nicefrac{2\ell}{3} + 1}^i, \ldots , R_\ell^i$ and $F_{\ell' + \nicefrac{2\ell}{3} + 1}^i, \ldots , F_\ell^i$ using \Cref{lem:detour} and then apply \Cref{cor:findoddwall} and \Cref{lem:jumpsoverawall} once more to find another separating odd wall $W_2^i$ of order at least 3 in $H_\mathfrak{M}$ under an apex set $A_2^i$ of order $\SeparatingWallApex{t}$.
    Note that we can choose the wall $W_2^i$ and the involved dicycles and paths such that they are disjoint from $A_1^i$ since
    \[ \ell' \geq 3 \PlanarGrid{ \FindOddWallOrder{ \SeparatingWall{t} , t } } + 18\InductionStepApex{t} . \]
    For this wall, we let $(X_3^i, Y_3^i), (X_4^i, Y_4^i)$ be the two separations that witness that $W_2^i$ is separating, such that $X_j^i \cap \sigma_{\delta_i'}(w_i) \neq \emptyset$ and thus $Y_j^i \cap \sigma_{\delta_i'}(w_i) = \emptyset$ for both $j \in \{ 3, 4 \}$ and $W_2^i \subseteq Y_4^i$.

    \textbf{Finding odd walls in a refined diamond rim:}
    From the previous two cases we can take the walls $W_1^i$ and $W_2^i$, and their associated apex sets $A_1^i$ and $A_2^i$.
    Our goal will be to find one odd wall for each transaction in $\mathfrak{A}_i$.
    However, we are not able to incur an additional apex set for each of these transactions, since $k$ is not bounded by any function of $t$.
    
    To resolve this, we can find $k$ mutually disjoint cylindrical $\FindOddWallOrder{ \SeparatingWall{t} , t }$-walls $U^{i*}_1, \ldots , U^{i*}_k$, using \Cref{cor:getseparatingcylindricalgrid}, such that for each $j \in [k]$ the wall $U^{i*}_j$ is found in the dicycles $R_{\ell' + 1}^i, \ldots , R_{\ell' + \nicefrac{\ell'}{3}}^i$, the paths $F_{\ell' + 1}^i, \ldots , F_{\ell' + \nicefrac{\ell'}{3}}^i$, and the paths $A_{g\ell + \ell' + 1}^{i,j}, \ldots , A_{g\ell + \ell' + \nicefrac{\ell'}{3}}^{i,j}$.
    Note that these walls are disjoint since each $\mathcal{A}^i_j$ is an eddy transaction.
    To make these walls odd, we prove a claim which essentially asserts the odd wall theorem again, but for several odd walls which all live in the same outline.

    \begin{claim}\label{claim:quantumoddwall}
        There either exists a half-integral packing of $t$ even dicycles or there exists a set of vertices $A_3^i \subseteq H_\mathfrak{M}$ with $|A_3^i| \leq 2 ( \FindOddWallApex{t} + \SeparatingWallApex{t} )$, such that there exist $k$ odd walls $U^i_1, \ldots , U^i_k$ under $A_3^i$ in $H_\mathfrak{M}$ and each $U^i_j$ is separating in $H_\mathfrak{M}$ under $A_3^i$.
    \end{claim}
    \emph{Proof of \Cref{claim:quantumoddwall} (Rough Sketch):}
    This essentially boils down to retreading the proof of \Cref{thm:oddwall}.
    However, unlike in \Cref{sec:oddwall}, we simultaneously find tilings and try to clean-up all walls $U^{i*}_1, \ldots , U^{i*}_k$ at once.
    This is largely straightforward, since the walls are already disjoint, except whenever we build an auxiliary digraph and apply \Cref{thm:xpaths}, since the paths this theorem yields may have endpoints in the auxiliary digraphs belonging to different walls.
    At this point, we can use the fact that all of these walls reside in the same outline and this allows us to reroute the paths as we did in the proofs of \Cref{lem:buildnewmaelstromscircleapexset} and \Cref{lem:buildnewmaelstromsrefineddiamond}, resulting in many peripheral jumps, possibly on different walls.
    This allows us to apply \Cref{lem:gridplusjumpmakesoddbicycle} to find a half-integral packing of $t$ even dicycles if enough such paths exist.
    If many of the paths have endpoints in the same wall, we can of course then proceed as in the proof of \Cref{thm:oddwall}.
    This necessitates another factor of 2 in the size of the apex set we find in \Cref{thm:oddwall}.

    Afterwards, we reprove \Cref{lem:jumpsoverawall} in the same way, increasing the potential size of the apex set we find by a factor of 2 again for the size of the apex set found in \Cref{lem:jumpsoverawall}, and rerouting paths to allow us to apply \Cref{lem:gridplusjumpmakesoddbicycle}, if we cannot proceed exactly as in the proof of \Cref{lem:jumpsoverawall}.
    Thus we can find the apex set $A_3^i$ with the desired properties.
	\hfill$\blacksquare$

    In particular, due to \Cref{claim:quantumoddwall}, there must exist an undirected separation $(X_5^i, Y_5^i)$ in $H_\mathfrak{M}$ such that $X_5^i \cap Y_5^i$ consists of $A_3^i$ and one dicycle of the perimeter of each wall $U^i_1, \ldots , U^i_k$, such that $U^i_1, \ldots , U^i_k \subseteq Y_5^i$ and every vertex of $Y_5^i$ reaches or is reached by a vertex in $\bigcup_{j=1}^k V(U^i_j)$.
    We note that since all walls have order at least 3, \Cref{thm:reinforcedmatcross} tells us that $Y_2^i$, $Y_4^i$, and $Y_5^i$ each have a non-even rendition.
    Further, we point out that for each $i \in [s+1]$ we have
    \[ \InductionStepApex{t} \geq \nicefrac{\InductionStepApex{t}}{2} + \BuildOutlineApex{t} + 4 \FindOddWallApex{t} + 4 \SeparatingWallApex{t} \geq |A_{i-1} \cup A_i' \cup A_1^i \cup A_2^i \cup A_3^i| . \]
    Let $A_{s+1}'' = A_s \cup A_{s+1}' \cup A_1^{s+1} \cup A_2^{s+1} \cup A_3^{s+1}$
    Consider $\mathfrak{W}_{s+1}$ and $\mathfrak{R}_{s+1}$ again and note that within $H_\mathfrak{M} - A_{s+1}''$ we can easily shrink $\mathfrak{W}_{s+1}$ into a $\InductionStepOutline{t-1}$-outline $\mathfrak{W}_{s+1}'$ (of roughness $g$) around $w_{s+1}$ and an associated $\nicefrac{\InductionStepOutline{t-1}}{2}$-rim $\mathfrak{R}_{s+1}'$ chosen such that $\mathfrak{R}_{s+1}'$ uses the dicycles $R_1^i, \ldots , R_{\ell'}^i$, the paths $F_1^i, \ldots , F_{\ell'}^i$, if they exists, and for each $j \in [k]$, the paths $A_1^{i,j}, \ldots , A_{(g+1)\ell'}^{i,j}$, since
    \[ \nicefrac{\InductionStepOutline{t}}{2} - 18\InductionStepApex{t} \geq \InductionStepOutline{t-1} . \]
    Let $( H_{\mathfrak{W}_{s+1}'}, \Omega_{\mathfrak{W}_{s+1}'} )$ be the $\mathfrak{W}_{s+1}'$-society.
    Note that our choices for $W_1^{s+1}$, $W_2^{s+1}$, $U^i_1, \ldots , U^i_k$, and $\mathfrak{W}_{s+1}'$ therefore imply that $W_1^{s+1} \cup W_2^{s+1} \cup \bigcup_{j=1}^k U^i_j$ and $H_{\mathfrak{W}_{s+1}'}$ are disjoint.
    In particular, $Y_2^{s+1} \cup Y_4^{s+1} \cup Y_5^{s+1}$ and $H_{\mathfrak{W}_{s+1}'}$ are disjoint.

    We let $U = Y_2^{s+1} \cup Y_4^{s+1} \cup Y_5^{s+1}$.
    Suppose $U$ contains an even dicycle $Z$, then \Cref{claim:findsmallmaelstrom} implies that $H_{\mathfrak{W}_{s+1}} - (A_1^{s+1} \cup A_2^{s+1} \cup A_3^{s+1})$ is disjoint from $Z$.
    However, according to \Cref{claim:findsmallmaelstrom}, there exists an odd rendition of $\rho'_{s+1}$ of $( H - (V(H_{\mathfrak{W}_{s+1}} - V(\Omega_{\mathfrak{W}_{s+1}})) \cup A_s \cup A_{s+1}'), \Omega_\mathfrak{M}, \Omega_{\mathfrak{W}_{s+1}} )$ and there exists another odd rendition for $H_{\mathfrak{R}_{s+1}}$, where $( H_{\mathfrak{R}_{s+1}}, \Omega_{\mathfrak{R}_{s+1}} )$ is the $\mathfrak{R}_{s+1}$-society.
    The existence of these two compatible renditions tells us that $U - (A_s \cup A_{s+1}')$ is odd and thus $Z$ must intersect $A_s \cup A_{s+1}'$.
    Therefore, the set $A_{s+1}''$, the $\InductionStepOutline{t-1}$-outline $\mathfrak{W}_{s+1}'$ around $w_{s+1}$, and the even dicycle $Z$ satisfy \cref{item:ISsplitoffdicycle} of our statement.
    Since we assumed that this item does not hold, we therefore now know that $U$ is odd.
    We set $A_{s+1} = A_1^{s+1} \cup A_2^{s+1} \cup A_3^{s+1}$ and note that $U \subseteq H_\mathfrak{M} - A_{s+1}$.
    
    Our goal is to now find the new odd transaction and the associated outline to verify the last point of our statement.    
    Recall that for all $i \in [s+1]$ we have chosen indices such that the majority of the paths in $\mathcal{Q}^i_1$ lie in $H_{\mathfrak{W}_i}$.
    We also note that for the disks $d_1, \ldots , d_s$ from \Cref{def:overflow} we have $d_i \subseteq d_j$ for each pair $i,j \in [s]$ with $i > j$ and of course $d_i \subseteq m$ for all $i \in [s]$.
    An immediate consequence of both of these facts is that for each $i \in [s+1]$ the majority of the paths in $\mathcal{Q}^i_3$ are found in $U$.
    We are now at the part of the proof in which the fact that $A_{s+1}$ may separate parts of the paths in $\mathcal{Q}^i_3$ and move them into what will become the maelstrom makes life difficult.

    Let us first outline the core problem we seek to remedy in the second part of the proof and describe our approach to resolving it.
    It is straightforward to show that given $A_{s+1}$ we can now construct an outline that basically satisfies the last item of our statement, since each $\mathcal{Q}^i_3$ lies in $U$ and the transactions are large enough to allow for the loss of $\InductionStepApex{t}$ vertices and still leave enough paths to build an odd transaction of order $\BuildOutline{t}$.
    This will give us an outline that is very similar to $\mathfrak{M}_s$, though the difference in which vertices are contained in $A_{s+1}$ and those contained in $A_s$ may have shifted things somewhat, and from here we simply need to choose enough paths in $\mathcal{Q}_{s+1}$ to find the outline that confirms that $\mathfrak{M}$ has overflow $s+1$.
    
    However, the complication is that the overflow of an outline is not necessarily bounded in $t$.
    Suppose that in the first iteration of this process $A_1$ hits the path $Q^1_p$ and thus the society of the outline $\mathfrak{M}_1$ we find has moved inwards slightly.
    In the second iteration, we choose a new transaction which has its endpoints on the $\mathfrak{M}_2$-society.
    From this point forward, we can therefore not continue building our outlines with $Q^1_p$ since the paths of the second transaction we find do not interact with this path.
    The next apex set we find may now contain a vertex of $Q^1_{p-1}$ and this process may repeat itself in each iteration, which ultimately would chip away so much of $\mathcal{Q}_1$ that it becomes too small to apply the lemmas we use in the other parts of the proof.

    Thus we must somehow argue that the number of times this can happen is bounded in some function of $t$.
    We will achieve this through either finding a cross consisting of the path that is hit by the apex set and another path we find, which correspond to the crosses associated with the transactions in \Cref{def:overflow}, or alternatively, if we cannot find such a path, we will show that the remainder of this path in $U$ still separates certain parts of $U$.
    The latter option allows us to still use the remainder of the path to build a new society, since for our purposes a society does not necessarily have to form a cycle, but instead needs to be a set of vertices with a cyclic permutation that describes a closed curve separating two parts of an odd decomposition.
    Note that for this reason we do not need to concern ourselves with the case in which $A_{s+1}$ intersects $V(\Omega_\mathfrak{M})$, since the remaining vertices of the society will still have the properties stated above.

    We show this by iterating over the transactions we found in ascending order of their index, starting with $i = 1$ and thus $\mathcal{Q}_1$.
    Once we have shown this for $i=1$ the procedure can be repeated analogously for all other transactions.
    Note that according to \Cref{def:overflow} the transaction $\mathcal{Q}_1$ already has an associated spill-over $s_1 \in [0, t-1]$ and if $s_1 \neq 0$ there exist $s_1$ crossing pairs $(L_1^1, R_1^1), \ldots , (L_{s_1}^1, R_{s_1}^1)$ on the $\mathfrak{M}$-society, such that all paths in the crossing pairs are pairwise disjoint, for each $r \in [s_1]$ the path $L_r^1$ is disjoint from the paths in $\{ Q^1_1, \ldots , Q^1_{p - (r - \nicefrac{1}{2}) (2\InductionStepApex{t} + 2)} \}$ and $H_{\mathfrak{M}_1}$, and $R_r^1 = Q^1_{p - (r - 1)(2\InductionStepApex{t} + 2)}$.
    
    We now might need to update the spill-over value and the outline $\mathfrak{M}_1$ as follows.
    Consider the path $Q = Q^1_{p - s_1(2\InductionStepApex{t} + 2)}$.
    Recall that, as previously laid out, $\mathcal{Q}^1_3$ resides in $U$, which has an odd rendition, and note that $Q \in \mathcal{Q}^1_3$.
    At this point, we again need to distinguish cases according to the type of the outline $\mathfrak{M}_{i-1}$, though here the two types of diamond outlines only require small adjustments to the constructions.
    
    \textbf{Determining the spill-over of $\mathcal{Q}_1$ starting from a circle outline:}
    If $A_{s+1}$ does not contain any vertices of $Q$, we can choose $Q$ and the $\InductionStepOutline{t} - (6\InductionStepApex{t} + 1)$ other paths with the highest indices in $\mathcal{Q}^1_3 \setminus \{ Q^1_{p - s_1(2\InductionStepApex{t} + 2)} , \ldots , Q^1_p \}$ that are also not hit by $A_{s+1}$ and continue with the next transaction, provided that we choose our functions such that
    \[ \nicefrac{\InductionStepOrder{t,g}}{3} \geq \OddTransactionOrder{\InductionStepOutline{t} + 6\InductionStepApex{t} + t(2\InductionStepApex{t} + 2)} . \]
    We may therefore suppose that $A_{s+1}$ contains at least one vertex of $Q$.
    According to \Cref{def:overflow} the $\InductionStepApex{t} + 1$ paths $Q^1_{(p + 1) - s_1(2\InductionStepApex{t} + 2)}, \ldots , Q^1_{p - (r - \nicefrac{1}{2}) (2\InductionStepApex{t} + 2)}$ do not intersect any of the crossing pairs already associated with $\mathcal{Q}_1$.
    Thus there exists some path $R$ among these which is not intersected by $A_{s+1}$ and similarly, since $Q$ is hit by $A_{s+1}$, there exists another path $L$ amongst the paths $ Q^1_{p - (r + \nicefrac{1}{2})(2\InductionStepApex{t} + 2)}, \ldots , Q^1_{p - s_1(2\InductionStepApex{t} + 2)}$ that is distinct from $Q$ and is also spared by $A_{s+1}$.

    We note that the endpoints of the paths $Q$, $L$, and $R$ all reside in $V(\Omega_\mathfrak{M})$.
    Since both $V(L)$ and $V(R)$ are disjoint from $A_{s+1}$ and $U$ has an odd rendition, the removal of the endpoints of $L$, or respectively $R$, separates $\Omega_\mathfrak{M}$ into two segments such that there are no paths with endpoints in both segments in $U - V(L)$, or $U - V(R)$ respectively.
    As a consequence we can choose a subset of vertices of $V(\Omega_\mathfrak{M}) \cup V(L) \cup V(R)$ with a cylindrical order $\Omega$ such that $V(L) \cup V(R) \subseteq V(\Omega)$ and both endpoints of $Q$ reside in $V(\Omega)$.
    (We note that the endpoints of $Q$ may reside in $A_{s+1}$ but for the definition of $\Omega$ we can return the vertices of $V(\Omega) \cap A_{s+1}$ to the graph.)
    Thus $Q$ similarly divides $\Omega$ into two segments $S, S'$.
    If there exist no $S$-$S'$- and no $S'$-$S$-paths in $U - V(Q)$, we may proceed by leaving $s_1$ as is and using $Q$ as part of the $\mathfrak{M}_1$-society as explained above, whilst adjusting the remainder of $\mathfrak{M}_1$ by replacing some paths in $\mathcal{E}^1_{d+1}$ with paths in $\mathcal{Q}^1_3 \setminus \{ Q = Q^1_{p - s_1(2\InductionStepApex{t} + 2)} , \ldots , Q^1_p \}$ that are not hit by $A_{s+1}$, which is possible since we choose our functions such that
    \[ \nicefrac{\InductionStepOrder{t,g}}{3} \geq \OddTransactionOrder{\InductionStepOutline{t} + 18\InductionStepApex{t} + t(2\InductionStepApex{t} + 2)} . \]
    This results in a new outline $\mathfrak{M}_1'$ that we can continue our construction with, whilst $\mathcal{Q}_1$ retains its spill-over value and thus in particular has not accumulated $t$ or more crosses.

    We can therefore suppose that there exists a $S_1$-$S_1'$- or a $S_1'$-$S_1$-path $P$ in $U - V(Q)$.
    Note that if both endpoints of $P$ are already in $V(\Omega_\mathfrak{M})$, we have found a cross, consisting of $R_1^{s_1 + 1} = Q$ and $L^1_{s_1 + 1} = P$, that confirms that the spill-over of $\mathcal{Q}$ has increased.
    Otherwise, we note that by our choice of $L$ and $R$, since $L$ has a strictly lower index than $p - s_1(2\InductionStepApex{t} + 2)$ and $R$ has a strictly higher index than $p - s_1(2\InductionStepApex{t} + 2)$, we have $V(L) \subseteq V(S)$ and $V(R) \subseteq V(S')$ or vice versa.
    This means, since $L$ and $R$ are directed paths that are not hit by $A_{s+1}$, we find a directed path $P'$ in $L \cup R \cup P$ with both endpoints in $V(\Omega_\mathfrak{M})$ that forms a crossing pair with $Q$.
    
    Thus we can increase the spill-over of $\mathcal{Q}_1$ and this process can be repeated until we either find some new value $s_1' \in [0, t-1]$ such that $Q^1_{p - s_1'(2\InductionStepApex{t} + 2)}$ is disjoint from $A_{s+1}$, which allows us to replace $\mathfrak{M}_1$ as laid out above, or we find $t$ pairwise disjoint crossing pairs on the $\mathfrak{M}$-society, which clearly allows us to find a quarter-integral packing of $t$ even dicycles via \Cref{lem:applyshifting}.

    \textbf{Determining the spill-over of $\mathcal{Q}_1$ starting from a diamond outline:}
    Note that by our choice of $W_1^1$ and $W_2^1$ the dicycles $R^1_{\nicefrac{\ell}{2} + 1}, \ldots , R^1_\ell$ all lie in $U$.
    In particular, all paths in $\mathcal{Q}^1_3$ that are not hit by $A_{s+1}$ intersect those dicycles in $\{ R^1_{\nicefrac{\ell}{2}}, \ldots , R^1_\ell \}$ spared by $A_{s+1}$ in a planar way, since $U$ has an odd rendition.
    Thus, since we can choose our functions such that
    \[ \ell' \geq t\SC + 18\InductionStepApex{t} \text{ and } \nicefrac{\InductionStepOrder{t,g}}{3} \geq t\SC + 18\InductionStepApex{t} + t(2\InductionStepApex{t} + 2) , \]
    we can use \Cref{lem:detour} to find a set of $t\SC$ homogeneous dicycles $R_1' , \ldots , R_{t\SC}'$ within $\bigcup_{k = \nicefrac{\ell}{2} + 1 }^\ell R^1_k \cup \bigcup_{k = 2q + 1}^{p - t(2\InductionStepApex{t} + 2)} Q^1_k$.
    These dicycles define a $t\SC$-rim $\mathfrak{R}$ positioned such that each crossing pair associated with $\mathcal{Q}_1$ contains a crossing pair on the $\mathfrak{R}$-society.
    Thus we can use $\mathfrak{R}$ in the same way we used $\mathfrak{M}$ and its associated rim in the previous case.
    The remainder of this case can then be resolved analogously to the previous one, either adjusting $\mathfrak{M}_1$ and $s_1$ or finding a quarter-integral packing of $t$ even dicycles.

    \textbf{Determining the spill-over of $\mathcal{Q}_1$ starting from a refined diamond outline:}
    For this case we first need to recall that we choose $\mathcal{Q}_i$ to be a calm transaction if $\mathfrak{M}_{i-1}$ is a refined diamond outline.
    As a consequence, we can reroute any crossing pair we find onto the society of the diamond outline underlying $\mathfrak{N}$ via \Cref{lem:detour}.
    This allows us to then proceed analogously to the previous two cases and use the rim described in the case for diamond outlines, should we find $t$ crossing pairs.
    With this we have finished the description of a step in our iterative procedure.

    As mentioned above, for each $i \in [2,s+1]$, we can repeat this procedure for $\mathcal{Q}_i$, where we replace $\mathfrak{M}$ with $\mathfrak{M}_{i-1}$ in the arguments.
    We also note that there does not exist a value $s_{s+1}$, but we can simply initialise this value with 0.
    This proves that we can use $\mathcal{Q}_{s+1}$ to show that $\mathfrak{M}$ in fact has at least overflow $s+1$, as demanded in \cref{item:ISoverflow}, or find a quarter-integral packing of $t$ even dicycles in the process.
\end{proof}

Via a simple induction we can derive the following powerful corollary from \Cref{lem:localstructureinductionstepgeneral}.

\begin{corollary}\label{cor:localstructureinductionstepgeneral}
    Let $\theta , t, g$ be positive integers such that $\theta = \InductionStepOutline{t}$.
    Moreover, let $D$ be a digraph with an odd decomposition $\delta$ and let $m$ be a maelstrom of $\delta$, together with a $\theta$-outline $\mathfrak{M}  = ( H , \mathcal{C}, \mathfrak{E}, \mathfrak{V} )$, with $(H_\mathfrak{M},\Omega_\mathfrak{M})$ being the $\mathfrak{M}$-society and a segregated $\mathcal{C}$-pair of order $\nicefrac{\theta}{2}$ if $\mathfrak{M}$ is a circle outline with degree zero, and if $\mathfrak{M}$ is a diamond outline, a refined diamond $\theta$-outline $\mathfrak{N}$ for $n \subseteq m$ of roughness $g$.

    Then one of the following holds
    \begin{enumerate}
        \item there exists a quarter-integral packing of $t$ even dicycles in $D$,\label{item:packing}

        \item there exists a set $A \subseteq V(H_\mathfrak{M})$ with $|A| \leq \InductionStepApex{t}$ such that $H_\mathfrak{M} - A$ has an odd decomposition $\delta'$ with a unique maelstrom $m'$ and an associated $(\InductionStepOutline{t} - 6\InductionStepApex{t})$-outline $\mathfrak{M}'$ (of roughness $g$) with depth less than $\InductionStepOrder{t,g}$,\label{item:depth}

        \item there exists a set $A \subseteq V(H_\mathfrak{M})$ with $|A| \leq \InductionStepApex{t}$ such that $H_\mathfrak{M} - A$ is odd,\label{item:odd}
        
        \item there exists an $h \in \N$ with $1 < h < t$, pairwise disjoint disks $d_1^*, \ldots , d_h^* \subseteq m$, and a set $S \subseteq V(H_\mathfrak{M})$ with $|S| \leq \InductionStepApex{t}$, such that for each $i \in [h]$ there exists a $\theta_i$-outline $\mathfrak{M}_i^* = ( H_i^* , \mathcal{C}_i^*, \mathfrak{E}_i^*, \mathfrak{B}_i^*, \mathfrak{V}_i^*, \mathfrak{D}_i^* )$ (with roughness $g$) for $d_i^{**} \subseteq d_i^*$ within $H_\mathfrak{M} - S$ with $\theta_i \geq \InductionStepOutline{t - (h - 1)}$, the set $\{ H_1^*, \ldots , H_h^* \}$ is third-integral, the trace of $d_i^*$ corresponds to the $\mathfrak{M}_i^*$-society $(H_{\mathfrak{M}_i^*}, \Omega_{\mathfrak{M}_i^*})$, and the graph $H_{\mathfrak{M}_i^*}$ contains an even dicycle,\label{item:splitmaelstrom} or
        
        \item there exists a set $A \subseteq V(H_\mathfrak{M})$ with $|A| \leq \InductionStepApex{t}$ such that $H_\mathfrak{M} - A$ has an odd decomposition $\delta'$ with a unique maelstrom $m'$ and an associated $\InductionStepOutline{t-1}$-outline $\mathfrak{M}'$ (of roughness $g$), where $(H_{\mathfrak{M}'},\Omega_{\mathfrak{M}'})$ is the $\mathfrak{M}'$-society, $H_{\mathfrak{M}'}$ contains an even dicycle, there exists an odd rendition $\rho'$ of $( H_\mathfrak{M} - (V(H_{\mathfrak{M}'} - V(\Omega_{\mathfrak{M}'})) \cup A), \Omega_\mathfrak{M}, \Omega_{\mathfrak{M}'} )$, and $H_\mathfrak{M} - V(H_{\mathfrak{M}'} - V(\Omega_{\mathfrak{M}'}))$ contains an even dicycle (which must contain vertices of $A$).\label{item:splitoffdicycle}
    \end{enumerate}
\end{corollary}
\begin{proof}
    We assign $\mathfrak{M}$, respectively $\mathfrak{N}$, overflow 0 and then iteratively apply \Cref{lem:localstructureinductionstepgeneral}, which must terminate in an item of \Cref{lem:localstructureinductionstepgeneral} other than \cref{item:ISoverflow}, since $D$ is finite.
\end{proof}

\subsection{Finding the bottom of the abyss}

There is an unfortunate edge case in the proof of \Cref{thm:localstructure} that we need to discuss at some length to show that we can resolve it.
In particular, once we reach a refined diamond outline of bounded depth, \Cref{lem:killdiamondmaelstrom} provides a set of vertices which we can delete that almost separates all even dicycles from the infrastructure outside of the maelstrom, with the exception of a select few shorelines corresponding to eddies, which may still be both substantially well connected to the outside of the maelstrom and contain even dicycles.
If we can at least find two separable even dicycles in this refined diamond outline, we are already happy and if all even dicycles we find happen to be attached to the outside of the maelstrom through a part of the outline that is not a shoreline, we can also proceed easily.
Thus we identify the following situation as problematic.

\begin{definition}[Undertow]\label{def:undertow}
    Let $\theta , d, g$ be positive integers.
    Moreover, let $D$ be a digraph with an odd decomposition $\delta$ and let $m$ be a maelstrom of $\delta$, together with a $\theta$-outline $\mathfrak{M}$, with $(H_\mathfrak{M},\Omega_\mathfrak{M})$ being the $\mathfrak{M}$-society, and a refined diamond $\theta$-outline $\mathfrak{N}$ for $n \subseteq m$ of roughness $g$ and bounded depth $d$.
    Let $S_1', \ldots , S_h'$ be the shorelines on the cycle $C$ corresponding to $V(\Omega_\mathfrak{M})$, and let $W_1, \ldots, W_h$ be the undirected $V(C - S_i')$-$V(S_i')$-separator of size at most $2d$ implied to exist by \Cref{lem:shorelinedepth} for each $i \in [h]$.
    If there exists a unique $i \in [h]$ such that $U_i = S_{H_\mathfrak{M} - W_i}(V(S_i'))$ contains an even dicycle and each even dicycle in $H_\mathfrak{M}$ contains a vertex of $V(U_i) \cup W_i$, we say that $\mathfrak{N}$ has \emph{undertow at $S_i'$}.
\end{definition}

At this point we must pull back the curtain and reveal to the reader that we have been working with a simplified version of what an outline can be.
The reason we define the depth of refined diamond outline so strangely is that it is in fact perfectly normal for the outline to keep growing further on the eddies, resulting in an overall fractal-like cactus of interlocking transactions.
In almost all situations it suffices to consider this cactus only up to the first or second level of additional transactions.
However, there is an exception and this subsection is dedicated to it.

Therefore, should we encounter a refined diamond outline with undertow during the proof of \Cref{thm:localstructure} beckoning us deeper into the maelstrom, we will heed its call and dive in.
For this purpose, we first refine the definition of an outline by building the structure further in the direction of a fixed eddy.
Afterwards we rework the definition of overflow to be compatible with this new outline definition.
Using these definitions, our goal is to prove that when we find an outline with undertow, we can ensure that we ultimately find a small set of vertices which makes the graph odd, or we reach one of the other favourable options presented in \Cref{cor:localstructureinductionstepgeneral}.

\begin{definition}[Abyssal outline]\label{def:abyssaloutline}
    Let $\theta , t, g, k$ be positive integers and let $p, s$ be integers, with $a_i$ being a positive integer for each $i \in [p]$.
    Moreover, let $D$ be a digraph with an odd decomposition $\delta$ and let $m$ be a maelstrom of $\delta$, together with a $\theta$-outline $\mathfrak{M}$ of degree $t$ and a refined diamond $\theta$-outline $\mathfrak{N}$ for $n \subseteq m$ with surplus degree $k$ and roughness $g$ with undertow at the shoreline $S'$.
    Let $S \subseteq \Omega_\mathfrak{N}$ be the eddy segment corresponding to the shoreline $S'$, with $(H_\mathfrak{N}, \Omega_\mathfrak{N})$ being the $\mathfrak{N}$-society.

    We define the \emph{abyssal outline} $\mathfrak{G} = (H, \mathcal{C}, \mathfrak{E}, \mathfrak{B}, \mathfrak{P}, \mathfrak{S}, \mathfrak{V}, \mathfrak{D} )$ with \emph{profundity} $p$, roughness $g$, and \emph{scrutiny} $s$ for $w \subseteq n$, such that
    \begin{itemize}
        \item $\mathfrak{M} = ( H', \mathcal{C}, \mathfrak{E}, \mathfrak{V}' )$ with $H' \subseteq H$,
        
        \item $\mathfrak{N} = (H'', \mathcal{C}, \mathfrak{E}, \mathfrak{B}, \mathfrak{A}, \mathfrak{V}'', \mathfrak{D}'' )$ with $H'' \subseteq H$,
        
        \item the set $\mathfrak{P} = \{ \mathfrak{A}_1, \ldots , \mathfrak{A}_p \}$, with each $\mathfrak{A}_i = \{ \mathcal{A}^i_1, \ldots , \mathcal{A}^i_{a_i} \}$ consisting of linkages $\mathcal{A}^i_j = \{ A_1^{i,j}, \ldots , A_{(g+1)\theta}^{i,j} \}$ for each $i \in [p]$ and each $j \in [a_i]$,

        \item the set $\mathfrak{S} = \{ \mathcal{S}_1, \ldots , \mathcal{S}_s \}$ consists of linkages $\mathcal{S}_i = \{ S_1^i, \ldots , S_{(g+1)\theta}^i \}$ for each $i \in [s]$,

        \item $\mathfrak{D}$ is a collection of disk that are disjoint from $n$,
    \end{itemize}
    and the following hold:
    \begin{description}
        \item[A1] The graph $H'' \cup \bigcup_{i=1}^p \bigcup_{j = 1}^{a_i} \bigcup_{h=1}^{\theta(g+1)} A^{i,j}_h \cup \bigcup_{i=1}^s \bigcup_{j=1}^{\theta(g+1)} S^i_j $ is a subgraph of $H$.

        \item[A2] Let $\Delta$ be the $m$-disk of the cycle corresponding to the $\mathfrak{M}$-society and let $\rho'$ be the non-even rendition from the definition of $\mathfrak{N}$.
        Let $H_0 \coloneq H$.
        The society $(H_0, \Omega_0)$ has a non-even rendition $\rho = (\Gamma_{H_0}, \mathcal{V}_{H_0}, \mathcal{D}_{H_0})$ in $\Delta$ which contains the rendition $\rho'$, such that all paths of the linkages in $\mathfrak{A}$ and $\mathfrak{B}$ are grounded in $\rho$.

        \item[A3] For all $i \in [p]$ and $j \in [a_i]$, let $p_{i,j} = \sum_{h \in [i-1]} a_h + j$, and let $O_{p_{i,j}}$ be the $w$-tight cycle of $C_{\nicefrac{\theta}{2}} \cup \bigcup_{j\in[t]}E^j_{\nicefrac{\theta}{2}} \cup \bigcup_{j \in [k]}B^j_{\nicefrac{\theta(g+1)}{2}} \cup \bigcup_{j \in [i]} \bigcup_{h \in [a_j]} A^{j,h}_{\nicefrac{\theta(g+1)}{2}}$ in $\rho'$.
        Let $p' = \sum_{i \in [p]} a_i$.
        Furthermore, for all $i \in [p'+1, p'+s]$, let $O_i$ be the $w$-tight cycle of $C_{\nicefrac{\theta}{2}} \cup \bigcup_{j\in[t]}E^j_{\nicefrac{\theta}{2}} \cup \bigcup_{j \in [k]}B^j_{\nicefrac{\theta(g+1)}{2}} \cup \bigcup_{j \in [p]} \bigcup_{h \in [a_j]} A^{j,h}_{\nicefrac{\theta(g+1)}{2}} \cup \bigcup_{j \in [i - a]}S^j_{\nicefrac{\theta(g+1)}{2}}$ in $\rho'$.
        For each $i \in [p' + a]$ we denote the $w$-disk of $O_i$ by $c_i$.
        Notice that $O_i$ is grounded in $\delta$ and $\rho'$, and we require that $c_{p'+s} = w$.

        \item[A4] For all $i \in [p]$ and $j \in [a_i]$, let $Q_{p_{i,j}}$ be the $w$-tight cycle of $C_\theta \cup \bigcup_{j\in[t]}E^j_{\theta} \cup \bigcup_{j \in [k]}B^j_{\theta(g+1)} \cup \bigcup_{j \in [i]} \bigcup_{h \in [a_j]} A^{j,h}_{\theta(g+1)}$ in $\rho'$.
        For all $i \in [p'+1, p'+s]$, let $Q_i$ be the $w$-tight cycle of $C_\theta \cup \bigcup_{j\in[t]}E^j_{\theta} \cup \bigcup_{j \in [k]}B^j_{\theta(g+1)} \cup \bigcup_{j \in [p]} \bigcup_{h \in [a_j]} A^{j,h}_{\theta(g+1)} \cup \bigcup_{j \in [i - a]}S^j_{\theta(g+1)}$ in $\rho'$.
        For each $i \in [p' + a]$ we denote the $w$-disk of $Q_i$ by $\Delta_i$.
        For $i \in [p'+s]$, we define the society $(H_i,\Omega_i)$, such that $V(\Omega_i) = V(Q_i)$ and $H_i \subseteq H$ is the graph induced by all vertices and big vertices of $H$ drawn on $\Delta_i$ by $\Gamma_{H_0}$.

        \item[A5] Let $R_0 = S$.
        For all $i \in [p]$ and all $j \in [a_i - 1]$, the linkage $\mathcal{A}^{i,j}$ is a transaction on $(H_{p_{i,j}}, \Omega_{p_{i,j}})$, or $(H_\mathfrak{N}, \Omega_\mathfrak{N})$ if $i = j = 1$, with all endpoints on $R_{p_{i,j} -1}$, with the tails of $\mathcal{A}^{i,j}$ appearing before the heads when traversing $R_{p_{i,j} -1}$ from its head to its tail and we let $R_{p_{i,j}}$ be the unique directed path in $R_{p_{i,j} -1} \cup A^{i,j}_{(g+1)\theta}$ that contains $A^{i,j}_{(g+1)\theta}$ in its entirety.
        
        Further, for all $i \in [p]$, the linkage $\mathcal{A}^{i,a_i}$ is a transaction on $(H_{p_{i,a_i}}, \Omega_{p_{i,a_i}})$, or $(H_\mathfrak{N}, \Omega_\mathfrak{N})$ if $i = a_i = 1$, with all endpoints of $R_{p_{i,a_i} -1}$, with the heads of $\mathcal{A}^{i,a_i}$ appearing before the tails when traversing $R_{p_{i,a_i} -1}$ from its head to its tail and we let $R_{p_{i,a_i}} \coloneq A^{i,a_i}_{(g+1)\theta}$.

        \item[A6] For all $i \in [p]$ let $Z_i$ be the unique subpath of $A^{i,a_i}_{(g+1)\theta}$ that forms a dicycle with $R_{p_{i,a_i} - 1}$ in $R_{p_{i,a_i} - 1} \cup A^{i,a_i}_{(g+1)\theta}$.
        Let $(J^i, \Omega^i) = (H_{p_{i,a_i}},\Omega_{p_{i,a_i}})$ for ease of notation.
        There do not exist transactions on $(J^i, \Omega^i)$ of order at least $\InductionStepOrder{t,g}$ with their heads in $V(\Omega^i) \setminus V(Z_i)$ and their tails in $V(Z_i)$, or vice versa.
        Let $W_i$ be the undirected $(V(\Omega^i) \setminus V(Z_i))$-$V(Z_i)$-separator of order $2\InductionStepOrder{t,g}$, implied to exist by \Cref{thm:directedlocalmenger}, then $U_i = S_{J_i - W_i}(V(Z_i))$ contains an even dicycle and all even dicycles in $J_i$ contain a vertex of $V(U_i) \cup W_i$\footnote{Thus in the extended sense of \Cref{def:undertow}, this partial outline again has undertow and pulls the infrastructure one level deeper into the abyss.}.

        \item[A7] Let $P_0 = R_{p'}$ with $u$ being its tail and $v$ being its head.
        For all $i \in [s]$ the linkage $\mathcal{S}_i$ is a transaction on $(H_{p' + i-1}, \Omega_{p' + i-1})$ with all endpoints on $V(P_{i-1})$, the heads of $\mathcal{S}_i$ occurring before the tails of $\mathcal{S}_i$ on $P_{i-1}$ when traversing $P_{i-1}$ from $u$ to $v$.
        The path $P_i$ is then defined as the unique (undirected) path in $P_{i-1} \cup S^i_{(g+1)\theta}$ that has the endpoints $u$ and $v$ and contains $S^i_{(g+1)\theta}$ in its entirety\footnote{Note that the endpoints $u$ and $v$ can remain fixed, since no linkage of $\mathfrak{S}$ is allowed to cut them out of the society.}.

        \item[A8] For every $i \in [p]$, $j \in [a_i]$, and $h \in [2,\theta(g+1)-1]$ the trace of $A_h^{i,j}$ separates $A_{h-1}^{i,j}$ from $A_{h+1}^{i,j}$ in $\Delta_{p_{i,j}}$.

        \item[A9] For every $i \in [s]$, $j\in[2,\theta(g+1)-1]$ the trace of $S_j^i$ separates $S_{j-1}^i$ from $S_{j+1}^i$ in $\Delta_{p'+i}$.

        \item[A10] The restriction of $\rho'$ to all nodes drawn on common and conjunction cells intersecting $\Delta - ( (w - \Boundary{w}) \cup \bigcup_{d \in \mathfrak{D}} (d - \Boundary{d} ) )$ is an odd rendition that is a restriction of $\delta$.

        \item[A11] The vertices of $\sigma(w) - V(H)$ are drawn by $\Gamma$ into the disks from $\mathfrak{V}$ and the vertices of $\sigma(m - (V(H) \cup \sigma(w)))$ are drawn by $\Gamma$ into the disks from $\mathfrak{D}$.
		Finally, no vertex of $H$ is drawn in the interior of some disk from $\mathfrak{V} \cup \mathfrak{D}$.

        \item[A12]For every pair $d, d' \in \mathfrak{D} \cup \{ w \}$, every (not necessarily directed) path $P$ in $D$ with one endpoint in $\sigma(d)$ and the other in $\sigma(d')$ contains a vertex of $H - (\sigma(d) \cup \sigma(d'))$.
    \end{description}
    We call $(H_\mathfrak{G},\Omega_\mathfrak{G})$ the \emph{$\mathfrak{G}$-society} of $w$, where $H_\mathfrak{G} = H_{p'+k} \cup \sigma(w)$ and $\Omega_\mathfrak{G} = \Omega_{p'+k}$.
    Furthermore, we call $P_s$ the \emph{brink} of $\mathfrak{G}$ with the \emph{quasi-tail} $u$ and the \emph{quasi-head} $v$, and let $R_{p'}$ be the \emph{proto-brink}.
\end{definition}

Note that according to \Cref{def:abyssaloutline} any refined diamond outline with undertow can be cast as an abyssal outline.
We also want to mention that the definition for abyssal outlines still does not capture the full complexity of what an outline might look like, since the incredible specificity of the case we deal with here, allows us to focus on growing a very specific region of the outline.

We will also need the notion of a rim for abyssal outlines, though this definition is straightforward along the lines of the other rim definitions.
\begin{definition}[Abyssal rim]\label{def:abyssalrim}
Let $\theta, k, g, d, p, s$ be positive integers, let $D$ be a digraph, let $\delta = (\Gamma, \mathcal{V}, \mathcal{D})$ be an odd decomposition of $D$, and let $m \in C(\delta)$ be a maelstrom of $\delta$ with $n \subseteq m$.
An \emph{abyssal $\theta$-rim} of $w \subseteq n$ with \emph{profundity} $p$, \emph{scrutiny} $s$, and \emph{roughness} $g$ is a tuple $\mathfrak{R} = (H, \mathcal{C}, \mathcal{E}, \mathfrak{B}, \mathfrak{P} , \mathfrak{S}, \mathfrak{D})$, where
\begin{itemize}
    \item $\mathcal{C} = \{ C_1, \ldots , C_\theta \}$ is a homogeneous family of pairwise disjoint dicycles,

    \item $\mathcal{E} = \{ E_1, \ldots , E_\theta \}$ is a $C_\theta$-linkage,

    \item $\mathfrak{B} = \{ \mathcal{B}_1, \ldots , \mathcal{B}_d \}$ is a set of linkages with $\mathcal{B}_i = \{ B_1^i, \ldots , B_{\theta (g+1)}^i \}$, for all $i \in [d]$,

    \item $\mathfrak{P} = \{ \mathcal{A}_1, \ldots , \mathcal{A}_p \}$ is a set of linkages with $\mathcal{A}_i = \{ A_1^i, \ldots , A_{\theta (g+1)}^i \}$, for all $i \in [p]$,

    \item $\mathfrak{S} = \{ \mathcal{S}_1, \ldots , \mathcal{S}_s \}$ is a set of linkages with $\mathcal{S}_i = \{ S_1^i, \ldots , S_{\theta (g+1)}^i \}$, for all $i \in [s]$,

    \item $\mathfrak{D}$ is a collection of $d$ closed disks, each disjoint from $w$, and
\end{itemize}
the following requirements are met.

\begin{description}
    \item[AR1] $\mathfrak{R}' = (H', \mathcal{C}, \mathcal{E}, \mathfrak{B}, \mathfrak{D})$ is a refined diamond $\theta$-rim of $n$, where $H'$ is the restriction of $H$ to $\mathcal{C}$, $\mathcal{E}$, and $\mathfrak{B}$.

    \item[AR2] The points \textbf{\textsf{A1}}, \textbf{\textsf{A2}}, \textbf{\textsf{A4}}, \textbf{\textsf{A5}}, \textbf{\textsf{A6}}, \textbf{\textsf{A7}}, \textbf{\textsf{A8}}, \textbf{\textsf{A9}}, \textbf{\textsf{A10}}, \textbf{\textsf{A11}}, and \textbf{\textsf{A12}} of \Cref{def:abyssaloutline} hold for $\mathfrak{R}$, provided that
    \begin{itemize}
        \item we let $\mathfrak{V} = \{ w \}$, $\mathcal{E}_1 = \mathcal{E}$, as well as $E_i^1 = E_i$ for all $i \in [\theta]$,

        \item we let $\mathfrak{A}_i = \{ \mathcal{A}_i \}$ with $A^{i,1}_j = A^i_j$ for all $i \in [s]$ and $j \in [\theta (g+1)]$, and

        \item we let the object $\mathfrak{M}$ be $\mathfrak{R}'$.
    \end{itemize}
    Following \textbf{\textsf{A4}}, we let $Q_{p+s}$ be the $w$-tight cycle of $C_\theta \cup E_\theta \cup \bigcup_{i=1}^d B_{\theta (g+1)}^i \cup \bigcup_{j \in [p]} A^j_{\theta(g+1)} \cup \bigcup_{j \in [s]}S^j_{\theta(g+1)}$ and define $\Delta_i$ for all $i \in [p+s]$.
    
    \item[AR3] The graph $C_1 \cup E_1 \cup \bigcup_{i=1}^d B_1^i \cup \bigcup_{j \in [p]} A^j_1 \cup \bigcup_{j \in [s]}S^j_1$ has a $w$-tight cycle in $\delta$.
\end{description}

We call $(H_\mathfrak{R}, \Omega_\mathfrak{R})$ the \emph{$\mathfrak{R}$-society of $w$}, where $V(\Omega_\mathfrak{R}) = V(Q_{p+s})$ and $H_\mathfrak{R} = H \cup \sigma(w)$.
\end{definition}

It is again straightforward to prove that we can find abyssal rims inside of abyssal outline analogous to the results in \Cref{sec:localising}.
\begin{lemma}\label{lem:abyssaloutlinererouting}
    Let $g, h, \theta$ be integers, let $D$ be a digraph, let $\delta = (\Gamma, \mathcal{V}, \mathcal{D})$ be an odd decomposition of $D$, let $m \in C(\delta)$ be a maelstrom of $\delta$, with a diamond $\theta$-outline $\mathfrak{M}$, a refined diamond outline $\mathfrak{M}$ around $n \subseteq m$, and an abyssal outline $\mathfrak{G}$ around $w \subseteq n$ with profundity $p$, scrutiny $s$, and roughness $g$, where $(H_\mathfrak{G}, \Omega_\mathfrak{G})$ is the $\mathfrak{G}$-society of $w$.
    
    There exists an abyssal $\nicefrac{\theta}{2}$-rim $\mathfrak{R}$ of $w$ with profundity $p$, scrutiny at most $s$, and roughness $g$, such that the $\mathfrak{R}$-society is the $\mathfrak{N}$-society.
\end{lemma}

To prove the main result of this subsection, we also need a refinement of the definition for overflow that focuses on refining a given outline with undertow towards the abyss.

\begin{definition}[Riptide within an outline]\label{def:riptide}
    Let $\theta , t, d, g, k, r$ be integers such that $\theta = \InductionStepOutline{t}$, $g$ is positive, $p = \InductionStepOrder{t,g} - (2t - 1)$, and $q = \nicefrac{p}{3}$.
    Moreover, let $D$ be a digraph with an odd decomposition $\delta_0$ and let $m_0$ be a maelstrom of $\delta_0$, together with a $\theta$-outline $\mathfrak{M}_0  = ( H , \mathcal{C}, \mathfrak{E}, \mathfrak{V} )$ of degree $d$ and potentially a refined diamond $\theta$-outline $\mathfrak{N}$ for $n_0 \subseteq m_0$ of roughness $g$ and surplus degree $k$, if $\mathfrak{M}_0$ is a diamond outline, with $(H_0,\Omega_0)$ being the $\mathfrak{M}_0$-society, respectively the $\mathfrak{N}_0$-society.
    We set $m_0 \coloneqq n_0$, if $n_0$ exists, and let $d_0$ be the $m_0$-disk of the cycle corresponding to $\Omega_0$ in $\delta_0$. 

    We say that $\mathfrak{M}$, respectively $\mathfrak{N}$, has \emph{riptide} $r$ if there exist linkages $\mathcal{Q}_1, \ldots , \mathcal{Q}_r$ and for each $i \in [r]$
    \begin{enumerate}
                    \item $\mathcal{Q}_i$ is a planar transaction of order $p$ on $(H_{i-1},\Omega_{i-1})$, such that
                    \begin{enumerate}
                        \item if $\mathfrak{M}_{i-1}$ is a refined diamond outline $\mathcal{Q}_i$ is calm, and

                        \item if $\mathfrak{M}_{i-1}$ is an abyssal outline, with $P$ being its brink, then all of the endpoints of $\mathcal{Q}_i$ lie in $V(P)$,
                    \end{enumerate}
            
                    \item there exists a set $A_i \subseteq V(H_0)$ with $|A_i| \leq \nicefrac{\InductionStepApex{t}}{3}$, such that there exists an odd decomposition $\delta_i$ of $H_0 - A_i$ that contains a unique maelstrom $m_i$ together with a $(\theta - 6\InductionStepApex{t})$-outline $\mathfrak{M}_i$ (of roughness $g$) such that $H_i$ contains an even dicycle, with $(H_i,\Omega_i)$ being the $\mathfrak{M}_i$-society,

                    \item for each $i \in [r]$ the disk $d_i$ is the $m_i$-disk of the cycle corresponding to $\Omega_{\mathfrak{M}_i}$ in $\delta_i$, and the transaction $\mathcal{Q}_i = \{ Q^i_1, \ldots , Q^i_p \}$ is indexed such that for all $j \in [2, p-1]$ the trace of $Q^i_j$ separates the trace of $Q^i_{j-1}$ and $Q^i_{j+1}$ in $d_{i-1}$, and each $\mathcal{Q}_i$ is partitioned into $\mathcal{Q}^i_1 = \{ Q^i_1, \ldots , Q^i_q \}$, $\mathcal{Q}_2^i = \{ Q^i_{q + 1}, \ldots , Q^i_{2q} \}$, and $\mathcal{Q}_3^i = \{ Q^i_{2q + 1}, \ldots , Q^i_{3q} \}$,

                    \item for each $i \in [r]$ there exists an integer $r_i \in [0, t-1]$ and $r_i$ crossing pairs $(L_1^i, R_1^i), \ldots , (L_{r_i}^i, R_{r_i}^i)$ on the $\mathfrak{M}_{i-1}$-society, such that all paths in the crossing pairs are pairwise disjoint, for each $j \in [r_i]$ the path $L_j^i$ is disjoint from the paths in $\{ Q^i_1, \ldots , Q^i_{p - (\ell - \nicefrac{1}{2}) (2\InductionStepApex{t} + 2)} \} $ and from $H_{\mathfrak{M}_i}$, and $R_j^i = Q^i_{p - (j - 1)(2\InductionStepApex{t} + 2)}$,

                    \item for all $i \in [r]$, we have $V(\Omega_i) \subseteq V(\Omega_0) \cup V( \bigcup_{j = 1}^r Q^j_{ p - r_j(2\InductionStepApex{t} + 2)} )$, and
            
                    \item there exists an odd rendition $\rho_i$ for $( H_0 - (V(H_i - V(\Omega_i)) \cup A_i), \Omega_0, \Omega_i )$, or
    \end{enumerate}
    We call $\mathcal{Q}_1, \ldots , \mathcal{Q}_r$ the \emph{riptide-linkages} associated with $\mathfrak{M}$, respectively $\mathfrak{N}$, the value $s_i$ is called the \emph{spill-over} of $\mathcal{Q}_i$ for each $i \in [r]$, and we call $\mathfrak{M}_1, \ldots , \mathfrak{M}_r$ the \emph{riptide-outlines} associated with $\mathfrak{M}$, respectively $\mathfrak{N}$.
\end{definition}

We now prove a lemma that is in many ways analogous to \Cref{lem:localstructureinductionstepgeneral}, but of course deals with abyssal outlines.
This involves revisiting most sections of the paper and giving arguments as to why their contents translate to this specific setting.
Since this is only relevant for this single lemma, we restrict the discussion of these matters to its proof.

\begin{lemma}\label{lem:exploringtheabyss}
    Let $\theta , t, g, r$ be positive integers such that $\theta = \InductionStepOutline{t}$.
    Moreover, let $D$ be a digraph with an odd decomposition $\delta$ and let $m$ be a maelstrom of $\delta$, together with a $\theta$-outline $\mathfrak{M}  = ( H , \mathcal{C}, \mathfrak{E}, \mathfrak{V} )$, with $(H_\mathfrak{M},\Omega_\mathfrak{M})$ being the $\mathfrak{M}$-society and a segregated $\mathcal{C}$-pair of order $\nicefrac{\theta}{2}$ if $\mathfrak{M}$ is a circle outline with degree zero, and if $\mathfrak{M}$ is a diamond outline, a refined diamond $\theta$-outline $\mathfrak{N}$ for $n \subseteq m$ of roughness $g$.
    Further, let $\mathcal{Q}_1, \ldots , \mathcal{Q}_r$ be the riptide-linkages associated with $\mathfrak{M}$, respectively $\mathfrak{N}$, and let $\mathfrak{M}_1, \ldots , \mathfrak{M}_r$ be the riptide-outlines associated with $\mathfrak{M}$, respectively $\mathfrak{N}$, where $( H_{\mathfrak{M}_r}, \Omega_{\mathfrak{M}_r} )$ is the $\mathfrak{M}_r$-society.

    Then one of the following holds
    \begin{enumerate}
        \item there exists a quarter-integral packing of $t$ even dicycles in $D$,\label{item:ABYSSpacking}

        \item there exists an $A \subseteq V(H_\mathfrak{M})$ with $|A| \leq \max(\InductionStepApex{t}, 4t \InductionStepOrder{t,g} + t)$ such that $H_\mathfrak{M} - A$ is odd,\label{item:ABYSSodd}

        \item there exists a set $S \subseteq V(H_\mathfrak{M})$, with $|S| \leq \max(\InductionStepApex{t}, 4t \InductionStepOrder{t,g} + t)$, and there exist $h_1, h_2 \in \N$ with $h_1$ being positive and $1 < h_1+h_2 < t$, such that there exist disjoint even dicycles $C_1, \ldots , C_{h_2}$ and for each $i \in [h_1]$ there exists a circle $\theta_i$-outline $\mathfrak{M}_i^* = ( H_i^* , \mathcal{C}_i^*, \mathfrak{E}_i^* )$ within $H_\mathfrak{M} - S$ with $\theta_i \geq \InductionStepOutline{t - (h-1)}$, the set $\{ H_{\mathfrak{M}_1^*}, \ldots , H_{\mathfrak{M}_h^*} \}$ is half-integral, where $(H_{\mathfrak{M}_i^*}, \Omega_{\mathfrak{M}_i^*})$ is the $\mathfrak{M}_i^*$-society, the graph $H_{\mathfrak{M}_i^*} - V(H_i^*)$ is disjoint from all $C_1, \ldots , C_{h_2}$, and $H_{\mathfrak{M}_i^*}$ contains an even dicycle,\label{item:ABYSSeddies}
        
        \item there exists an $h \in \N$ with $1 < h < t$, pairwise disjoint disks $d_1^*, \ldots , d_h^* \subseteq m$, and a set $S \subseteq V(H_\mathfrak{M})$ with $|S| \leq \InductionStepApex{t}$, such that for each $i \in [h]$ there exists a circle, diamond, or refined diamond $\theta_i$-outline $\mathfrak{M}_i^* = ( H_i^* , \mathcal{C}_i^*, \mathfrak{E}_i^*, \mathfrak{B}_i^*, \mathfrak{V}_i^*, \mathfrak{D}_i^* )$ (with roughness $g$) for $d_i^{**} \subseteq d_i^*$ within $H_\mathfrak{M} - S$ with $\theta_i \geq \InductionStepOutline{t - (h-1)}$, the set $\{ H_1^*, \ldots , H_h^* \}$ is third-integral, the trace of $d_i^*$ corresponds to the $\mathfrak{M}_i^*$-society $(H_{\mathfrak{M}_i^*}, \Omega_{\mathfrak{M}_i^*})$, and the graph $H_{\mathfrak{M}_i^*}$ contains an even dicycle,\label{item:ABYSSsplitmaelstrom}
        
        \item there exists a set $A \subseteq V(H_\mathfrak{M})$ with $|A| \leq \InductionStepApex{t}$ such that $H_\mathfrak{M} - A$ has an odd decomposition $\delta'$ with a unique maelstrom $m'$ and an associated circle, diamond, or refined diamond $\InductionStepOutline{t-1}$-outline $\mathfrak{M}'$ (of roughness $g$), where $(H_{\mathfrak{M}'},\Omega_{\mathfrak{M}'})$ is the $\mathfrak{M}'$-society, $H_{\mathfrak{M}'}$ contains an even dicycle, there exists an odd rendition $\rho'$ of $( H_\mathfrak{M} - (V(H_{\mathfrak{M}'} - V(\Omega_{\mathfrak{M}'})) \cup A), \Omega_\mathfrak{M}, \Omega_{\mathfrak{M}'} )$, and $H_\mathfrak{M} - V(H_{\mathfrak{M}'} - V(\Omega_{\mathfrak{M}'}))$ contains an even dicycle (which must contain vertices of $A$), or\label{item:ABYSSsplitoffdicycle}

        \item there exists a planar transaction $\mathcal{Q}_{r+1}$ of order $\InductionStepOrder{t,g} - (2t - 1)$ on $( H_{\mathfrak{M}_r}, \Omega_{\mathfrak{M}_r} )$, a set of vertices $A_{r+1} \subseteq V(H_\mathfrak{M})$, and a $(\theta - 6\InductionStepApex{t})$-outline $\mathfrak{M}_{r+1}$ (of roughness $g$) in $H_\mathfrak{M} - A_{r+1}$, which together verify that $\mathfrak{M}$ has riptide $r + 1$.\label{item:ABYSSoverflow}
    \end{enumerate}
\end{lemma}
\begin{proof}
    We let $q$, $\delta_i$, $m_i$, $d_i$, $A_i$, $\mathfrak{M}_i$, and the transactions $\mathcal{Q}^i_1, \mathcal{Q}^i_2, \mathcal{Q}^i_3$ for each $i \in [0, s]$ be defined as in \Cref{def:riptide}, with some objects undefined for $i = 0$, and suppose that the first two items of our statement do not hold.
    Note that unlike in the proof of \Cref{lem:localstructureinductionstepgeneral}, this does not immediately mean that the $\mathfrak{M}_r$-society has depth large enough to find a new transaction and proceed.

    Should $\mathfrak{M}_r$ be a circle, diamond, or refined diamond outline with depth at least $\InductionStepOrder{t,g}$, we can follow the proof strategy of \Cref{lem:localstructureinductionstepgeneral} in its entirety.
    Suppose instead that $\mathfrak{M}_r$ is a circle outline with bounded depth, then \Cref{lem:killcirclemaelstrom} implies that there exists a set of vertices $X$ such that $H_\mathfrak{M} - X$ is odd with $|X| \leq 4t \InductionStepOrder{t,g} + t$.
    Suppose instead that $\mathfrak{M}_r$ is an abyssal outline.
    Note that we explicitly include refined diamond outlines here.
    Let $P$ be the brink of $\mathfrak{M}_r$ with the quasi-tail $u$ and the quasi-head $v$.
    Furthermore, let $P'$ be the proto-brink of $\mathfrak{M}_r$ and let $\mathfrak{M}_r'$ be the abyssal outline corresponding to the society that $P'$ is a part of.
    Consider the maximal family $\{ P_1', \ldots , P_k' \}$ of subpaths of $P'$ which are internally disjoint from $P$, thus having both endpoints in $V(P)$.
    These correspond in a natural way to eddy segments in a refined diamond outline, due to \textbf{\textsf{A6}} guaranteeing that $P'$ cannot be refined into another directed path that embeds more of the maelstrom and the condition on the transaction on the $\mathfrak{M}_r$-society corresponding to calm transactions.
    We distinguish two cases, with the first corresponding naturally to a refined diamond outline with bounded depth, and the second ultimately leading to us showing that we can increase the riptide.
    
    \textbf{Case 1: The brink cannot be extended.}
    Suppose that there does not exist a transaction of order at least $\InductionStepOrder{t,g}$ with all endpoints in $V(P)$ with its heads occurring before its tails when traversing $P$ from $u$ to $v$.
    Due to \textbf{\textsf{A6}} and the fact that the initial refined diamond outline that $\mathfrak{M}_r$ is built on has undertow (and bounded depth), we know that if there exist transaction of order at least $\InductionStepOrder{t,g}$ with all endpoints in $V(P)$, then all of its endpoints are found in some $P'_i$, with $i \in [k]$.
    This now corresponds to a refined diamond outline having bounded depth, allowing us to apply the proof strategy from \Cref{lem:killdiamondmaelstrom} (essentially just ignoring the part of the society that is not $P'$) to find either find a quarter-integral packing of $t$ even dicycles, which would mean we are done, or a set of vertices $X$ with $|X| \leq 4t \InductionStepOrder{t,g} + t$, such that there exists a selection $P_{i_1}', \ldots , P_{i_\ell}'$, with $i_1, \ldots , i_\ell \in [k]$ and the union of strong components of $H_\mathfrak{M} - ( X \cup \bigcup_{j=1}^\ell P_{i_j}' )$ containing vertices of $V(H_\mathfrak{M}) \cup (V(\Omega_{\mathfrak{M}_r'}) \setminus \bigcup_{j=1}^\ell V(P_{i_j}') ) $ is odd.
    
    If $\ell = 0$, we know $H_\mathfrak{M} - X$ is odd and we are done.
    Should $\ell \geq 1$ and there exists an even dicycle in $H_{\mathfrak{M}_r'}$ that is disjoint from $S_{H_{\mathfrak{M}_r'} - W_1}(P_{i_1}') \cup W_1$\footnote{This implies that this outline does not have undertow in the general sense.}, where $W_1$ is the undirected $V(\Omega_{\mathfrak{M}_r'}) \setminus V(P_{i_1}')$-$V(P_{i_1}')$-separator of size at most $2\InductionStepOrder{t,g}$ which must exist in this case, we can proceed as follows.
    Let $\mathfrak{R}_r = (H, \mathcal{C}, \mathcal{E}, \mathfrak{B}, \mathfrak{P} , \mathfrak{S}, \mathfrak{D})$ be the abyssal $\nicefrac{\theta}{2}$-rim of $\mathfrak{M}_r$ and note that $\mathfrak{R}_r$ has scrutiny $k$ since for each $j \in [k]$ we have $P_{i_j} \in \mathcal{S}_\ell$, for some $\ell \in [k]$.
    W.l.o.g.\ we assume that $P_{i_j}' \in \mathcal{S}_j$ for each $j \in [k]$.

    Let $\theta' = \nicefrac{\theta}{2}$ and let $p$ be the profundity of $\mathfrak{M}_r$, and thus also $\mathfrak{R}_r$.
    We let $\mathcal{S}_i = \{ S^i_1, \ldots , S^i_{\theta'(g+1)} \}$ for each $i \in [k]$ and $\mathcal{S}_i \in \mathfrak{S}$, and let $\mathcal{A}_p = \{ A^p_1, \ldots , A^p_{\theta'(g+1)} \}$.
    Note that for each $i \in [k]$ the head of $S^i_{\theta'(g+1)}$ appears before its tail on $A^p_{\theta'(g+1)}$.
    We let $\theta'' = \nicefrac{\theta'}{2^{t - (k - 1)}}$.
    For each $i \in [k]$, we can find a set of $\theta''$ homogeneous dicycles $\mathcal{R}_i = \{ R^i_1, \ldots , R^i_{\theta''(g+1)} \}$ in the union of the paths $A^p_{\theta' - ((k - (i+1)) \theta'')  + 1}, \ldots , A^p_{\theta' - ((k - i) \theta'')} \in \mathcal{A}_p$ and the paths $S^i_{\theta' - ((k - (i+1)) \theta'')  + 1}, \ldots , S^i_{\theta' - ((k - i) \theta'')} \in \mathcal{S}_i$.
    We note that this ultimately requires that for all $\ell' \in [t-1]$ we have $\nicefrac{\InductionStepOutline{t}}{2^{\ell' + 1}} \geq \InductionStepOutline{t - \ell'}$ which is guaranteed by the functions chosen in the previous subsection.
    
    Let $J_i = S_{H_{\mathfrak{M}_r'} - X}(P_{i_j}') \cup \bigcup_{h = 1}^{\theta''} R^i_h$ for each $i \in [k]$.
    Using $\mathcal{R}_i$, we can construct the circle $\theta''$-outline $\mathfrak{R}_i = ( J_i', \mathcal{R}_i, \emptyset )$, where we let $J_i'$ be a restriction of $J_i'$ to the infrastructure of $\mathcal{R}_i$ for each $i \in [k]$.
    We note that $J_1 \cup J_1', \ldots, J_k \cup J_k'$ are half-integral, with $J_i$ and $J_j'$ intersecting for distinct $i,j \in [k]$.
    Therefore, the circle outlines we found satisfy \cref{item:ABYSSeddies}.

    We are left with having to argue for the case in which $\ell = 1$ and all even dicycles in $H_{\mathfrak{M}_r'}$ use a vertex from $V(S_{H_{\mathfrak{M}_r'} - W_1}(P_{i_1}')) \cup W_1$.
    If there exists no transaction $\mathcal{Q}_{r+1}'$ of order at least $\InductionStepOrder{t,g}$ with all endpoints in $V(P)$ with its tails occurring before its heads when traversing $P$ from $u$ to $v$, we can use \textbf{\textsf{A6}} and the fact that the initial refined diamond outline that $\mathfrak{M}_r$ is built on has undertow (and bounded depth) to apply the proof strategy from \Cref{lem:killdiamondmaelstrom} to find a set of vertices $X$ with $|X| \leq 4t \InductionStepOrder{t,g} + t$ such that $H_{\mathfrak{M}_r'} - X$ is odd.
    Thus we instead suppose that $\mathcal{Q}_{r+1}'$ exists.

    This now leads us into the part of the proof in which we have to retread some sections of the proof.
    In particular, we need to argue that we can shift analogously to the results in \Cref{sec:shifting}, can find an odd transaction like in \Cref{sec:transaction}, and build new outlines like in \Cref{sec:buildoutline}.
    We will write each of these steps out as a claim.
    For the first claim on shifting, we note that we make a much more strict assumption on where the endpoints of the crosses land.

    \begin{claim}\label{claim:abyssshifting}
        Let $\theta, t$ be integers with $\theta$ being a multiple of $2t\SC$, let $D$ be a digraph, let $\delta = (\Gamma, \mathcal{V}, \mathcal{D})$ be an odd decomposition of $D$, let $m \in C(\delta)$ be a maelstrom of $\delta$, and let $\mathfrak{G}$ be an abyssal $\theta$-outline of $m$, with a transaction $\mathcal{L} = \{ L_1, \ldots , L_{2t} \}$ on the $\mathfrak{G}$-society $(H_\mathfrak{G}, \Omega_\mathfrak{G})$, such that for all $i \in [t]$ the paths $L_{2i-1}, L_{2i}$ are a crossing pair for the $\mathfrak{G}$-society and all heads or all tails of the paths in $\mathcal{L}$ lie on the brink of $\mathfrak{G}$.
    
        Then the digraph $H_\mathfrak{G}$ contains a quarter-integral packing of $t$ weak odd bicycles.
    \end{claim}
    \emph{Proof of \Cref{claim:abyssshifting} (Rough Sketch):}
    To see that this holds, we first note that we just showed that we can find a large circle rim involving the pseudo-brink of $\mathfrak{G}$ and using the linkages in $\mathcal{S}$ of the abyssal rim of $\mathfrak{G}$, we can modify this circle rim into a refined diamond $t\SC$-rim which involves both the pseudo-brink and the brink.
    Suppose now that w.l.o.g.\ the tails of the paths in $\mathcal{L}$ lie on the brink of $\mathfrak{G}$, then the heads of $\mathcal{L}$ through the infrastructure of the $\nicefrac{\theta}{2}$-rim of $\mathfrak{G}$ towards the refined diamond rim using \Cref{lem:detour} as we have done countless times before.
    \hfill$\blacksquare$

    The more strict assumption is warranted since, we actually only ever find this kind of linkage producing crosses in the proof of \Cref{thm:oddtransaction}, which allows us to easily deduce the next claim that we present without proof since the proof is analogous once we have made that observation.
    Note that we only try to find transactions on the brink of $\mathfrak{G}$ in the following claim.

    \begin{claim}\label{claim:abyssoddtransaction}
        Let $t$, $\theta$, and $p$ be integers with $\theta \geq 2t\SC + 6\OddTransactionApex{t}$ being even and $p \geq t > 1$.
        Let $D$ be a digraph with an odd decomposition $\delta$ and let $m$ be an abyssal maelstrom of $\delta$ together with a $\theta$-outline $\mathfrak{G}$ with the maelstrom society $(H,\Omega)$.
	
        Either $H$ contains a half-integral packing of $t$ even dicycles, or for every transaction $\mathcal{T}$ of order at least $\OddTransactionOrder{t,p}$ with all endpoints on the brink of $\mathfrak{G}$, there exists an apex set $A$ with $|A| \leq \OddTransactionApex{t}$, an abyssal $(\theta - 6|A|)$-outline $\mathfrak{G}'$ and an odd transaction $\mathcal{P}$ of order $p$ with all endpoints on the brink such that $A$ and $\mathcal{P}$ are $\delta$-compatible for $\mathfrak{G}'$.
    \end{claim}

    With both of these claims dealt with, we must now argue that an odd transaction on the brink of $\mathfrak{G}$ actually helps us split the abyssal maelstrom up in a useful way.

    \begin{claim}\label{claim:abysssplitting}
        Let $\theta , p, t$ be positive integers and let $g$ be non-negative, such that $\theta$ is a common multiple of $\BuildOutline{t}$, and let $D$ be a digraph with an odd decomposition $\delta$.
        Moreover, let $m$ be a maelstrom of $\delta$, together with a diamond $\theta$-outline $\mathfrak{M}$, a refined $\theta$-outline $\mathfrak{N}$ of $n \subseteq m$ with undertow and roughness $g$, and an abyssal $\theta$-outline $\mathfrak{G} = (H, \mathcal{C}, \mathfrak{E}, \mathfrak{B}, \mathfrak{P}, \mathfrak{S}, \mathfrak{V}, \mathfrak{D} )$ of $w \subseteq n$ with roughness $g$.
        Further, let $\mathcal{P} = \{ P_1, \ldots , P_p \}$ be an odd transaction with all endpoints on the $\mathfrak{G}$-society, together with a pleasant apex set $A \subseteq V(H_\mathfrak{G})$ such that $\mathcal{P}$ and $A$ are $\delta$-compatible.
        Let $\delta'$ be the $\mathcal{P}$-$A$-expansion of $\delta$, let $H_\mathcal{P} = (H \cup S) - A$, where $S$ is the strip of $\mathcal{P}$ under $A$.

        Either $H_\mathfrak{G}$ is odd, there exists a quarter-integral packing of $t$ even dicycles in $H_\mathfrak{G}$, or there exist 
        \begin{enumerate}
            \item a positive $h \in \N$ with $h < t$ and pairwise disjoint disks $d_1, \ldots , d_h \subseteq n$ such that for each $i \in [h]$ there exists a $\theta_i$-outline $\mathfrak{M}_i = ( H_i , \mathcal{C}_i, \mathfrak{E}_i, \mathfrak{B}_i, \mathfrak{P}_i, \mathfrak{A}_i, \mathfrak{V}_i, \mathfrak{D}_i )$ (of roughness $g$) for $d_i' \subseteq d_i$, with $\theta_i \geq \nicefrac{\theta}{t^{h-1}}$, the graph $\sigma_{\delta'}(d_i') \cup H_\mathcal{P}$ contains an even dicycle, $\{ H_1, \ldots , H_h \}$ is a third-integral set of graphs, $d_i$ is bounded by the cycle corresponding to the $\mathfrak{M}_i$-society $(H_{\mathfrak{M}_i}, \Omega_{\mathfrak{M}_i})$, and at most one of these is an abyssal outline, and

            \item there exists a vertex set $T \subseteq V(\bigcup_{i \in [h]} H_{\mathfrak{M}_i})$ with $|T| \leq 4t^2$, such that for all distinct $i,j \in [h]$ there does not exist a $V(\Omega_{\mathfrak{M}_i})$-$V(\Omega_{\mathfrak{M}_j})$- or a $V(\Omega_{\mathfrak{M}_j})$-$V(\Omega_{\mathfrak{M}_i})$-path in $D[V(H_{\mathfrak{M}_i} \cup H_{\mathfrak{M}_j})] - T$.
        \end{enumerate}
    \end{claim}
    \emph{Proof of \Cref{claim:abysssplitting} (Rough Sketch):}
    Since $\mathfrak{N}$ has undertow, which includes $\mathfrak{N}$ having bounded depth, and thanks to \textbf{\textsf{A6}}, we know that $\mathcal{P}$ cannot intersect with more than half of any transaction in $\mathfrak{B}$ or $\mathfrak{P}$, since $g$ is non-negative and this would otherwise contradict that $\mathfrak{N}$ has bounded depth, as we can use \Cref{lem:detour} to reroute a large transaction onto the $\mathfrak{N}$-society and this transaction would not be a whirl on any of the eddy segments of $\mathfrak{N}$.
    Thanks to this structural observation, we know that whilst $\mathcal{P}$ starts on the brink, it will not split our maelstrom in a way that grows several new outlines whose infrastructure leads deeper than refined diamond outlines do.
    We can thus adapt the arguments used to prove \Cref{lem:buildnewmaelstromscircle} and \Cref{lem:buildnewmaelstromsrefineddiamond} to prove the first part of this claim and separating the maelstroms can then be accomplished using arguments from \Cref{lem:buildnewmaelstromscircleapexset} and \Cref{lem:buildnewmaelstromsrefineddiamondapexset}.
    \hfill$\blacksquare$

    Now, with all of these claims under our belt, we can find an odd transaction $\mathcal{Q}_{r+1}$ within $\mathcal{Q}_{r+1}'$ of order $\InductionStepOrder{t,g} - (2t - 1)$.
    Using $\mathcal{Q}_{r+1}$, we now proceed through the proof in an analogous fashion to the proof for \Cref{lem:localstructureinductionstepgeneral}.
    We let $p$ and $q$ be defined as in \Cref{def:riptide} and choose indices such that for all $j \in [2, p-1]$ the trace of $Q^{r+1}_j$ separates the trace of $Q^{r+1}_{j-1}$ and $Q^{r+1}_{j+1}$ in $d_r$.
    We also partition $\mathcal{Q}_{r+1}$ into $\mathcal{Q}^{r+1}_1 = \{ Q^{r+1}_1, \ldots , Q^{r+1}_q \}$, $\mathcal{Q}_2^{r+1} = \{ Q^{r+1}_{q + 1}, \ldots , Q^{r+1}_{2q} \}$, and $\mathcal{Q}_3^{r+1} = \{ Q^{r+1}_{2q + 1}, \ldots , Q^{r+1}_{3q} \}$.
    
    \begin{claim}\label{claim:findsmallmaelstromABYSS}
        For each $i \in [r+1]$, there exists an odd transaction $\mathcal{Q}^{i*}_2 \subseteq \mathcal{Q}^i_2$ of order $(g+1)\theta$, a set of vertices $A_i'$ of size at most $\BuildOutlineApex{t}$, an odd decomposition $\delta_i'$ of $H - (A_{i-1} \cup A_i')$ with a unique maelstrom $w_i$ that comes equipped with a $\theta - 6\InductionStepApex{t}$-outline $\mathfrak{W}_i$ (of roughness $g$), and an odd rendition of $\rho'_i$ of $( H - (V(H_{\mathfrak{W}_i} - V(\Omega_{\mathfrak{W}_i})) \cup A_{i-1} \cup A_i'), \Omega_\mathfrak{M}, \Omega_{\mathfrak{W}_i} )$, where $( H_{\mathfrak{W}_i}, \Omega_{\mathfrak{W}_i} )$ is the $\mathfrak{W}_i$-society.
    \end{claim}
    \emph{Proof of \Cref{claim:findsmallmaelstromABYSS}:}
    This proof proceeds analogously to the proof of \Cref{claim:findsmallmaelstrom} in the proof of \Cref{lem:localstructureinductionstepgeneral}, using \Cref{claim:abysssplitting}.
	\hfill$\blacksquare$

    We will again assume w.l.o.g.\ that a majority of the paths in $\mathcal{Q}^{r+1}_1$ lie in $H_{\mathfrak{W}_{r+1}}$ and we let the other $\mathcal{Q}_i$ for $i \in [r]$ be indexed consistently such that the majority of the paths in $\mathcal{Q}^i_1$ lie in $H_{\mathfrak{W}_i}$.
    Our next step will be finding a collection of odd walls within our abyssal outline $\mathfrak{M}_r$.
    Recall that we chose $\mathfrak{R}_r$ as the abyssal $\nicefrac{\theta}{2}$-rim of $\mathfrak{M}_r$, let $p'$ be the profundity of $\mathfrak{R}_r$, and let $\ell' = \nicefrac{\theta - 6\InductionStepApex{t}}{2}$.
    We observe that for each $i \in [2, p']$ we can find a cylindrical $\FindOddWallOrder{ \SeparatingWall{t} , t }$-wall $U^*_i$ in $A^i_{\ell' + 1}, \ldots , A^i_{\ell' + \nicefrac{\ell'}{3}}$ and $A^{i-1}_{\ell' + 1}, \ldots , A^{i-1}_{\ell' + \nicefrac{\ell'}{3}}$ in an analogous fashion in which we found the walls in the eddies in the proof of \Cref{lem:localstructureinductionstepgeneral}.
    For $i=1$, we can find a wall in $A^1_{\ell' + 1}, \ldots , A^1_{\ell' + \nicefrac{\ell'}{3}}$ and the eddy of the refined diamond rim that had undertow.
    Given distinct $i,j \in [p']$ of the same parity, we can choose $U^*_i$ and $U^*_j$ such that they are disjoint and in particular, all of these walls are disjoint from the infrastructure of the underlying circle and diamond rims.
    Thus using the same arguments that we gave for \Cref{claim:quantumoddwall}, we can prove the following claim.

    \begin{claim}\label{claim:quantumoddwallABYSS}
        There either exists a half-integral packing of $t$ even dicycles or there exists a set of vertices $A_3^i \subseteq H_\mathfrak{M}$ with $|A_3^i| \leq 2 ( \FindOddWallApex{t} + \SeparatingWallApex{t} )$, such that there exist $k$ odd walls $U^i_1, \ldots , U^i_k$ under $A_3^i$ in $H_\mathfrak{M}$ and each $U^i_j$ is separating in $H_\mathfrak{M}$ under $A_3^i$.
    \end{claim}

    After this point, the arguments to first exclude \cref{item:ABYSSsplitoffdicycle} and then control the spill-over of the riptide-outlines are straightforward adaptations of the methods used in the proof of \Cref{lem:localstructureinductionstepgeneral} and we therefore omit this part.

    \textbf{Case 2: The brink can be extended.} 
    In case there does exist a transaction of order at least $\InductionStepOrder{t,g}$ with all endpoints in $V(P)$ with its heads occurring before its tails when traversing $P$ from $u$ to $v$.
    We can in fact simply use this transaction to repeat our arguments from above, with the slight difference that we do not need to discuss how to kill the maelstrom, since we do not have bounded depth in any sense.
    Note that the claims are stated in a way that is agnostic towards the direction of the transaction we start with.
\end{proof}

As in the previous section, we can iteratively use \Cref{lem:exploringtheabyss} to derive a useful corollary.

\begin{corollary}\label{cor:exploringtheabyss}
    Let $\theta , t, g, r$ be positive integers such that $\theta = \InductionStepOutline{t}$.
    Moreover, let $D$ be a digraph with an odd decomposition $\delta$ and let $m$ be a maelstrom of $\delta$, together with a $\theta$-outline $\mathfrak{M}  = ( H , \mathcal{C}, \mathfrak{E}, \mathfrak{V} )$, with $(H_\mathfrak{M},\Omega_\mathfrak{M})$ being the $\mathfrak{M}$-society and a segregated $\mathcal{C}$-pair of order $\nicefrac{\theta}{2}$ if $\mathfrak{M}$ is a circle outline with degree zero, and if $\mathfrak{M}$ is a diamond outline, a refined diamond $\theta$-outline $\mathfrak{N}$ for $n \subseteq m$ of roughness $g$.

    Then one of the following holds
    \begin{enumerate}
        \item there exists a quarter-integral packing of $t$ even dicycles in $D$,\label{item:ABYSSFINISHpacking}

        \item there exists an $A \subseteq V(H_\mathfrak{M})$ with $|A| \leq \max(\InductionStepApex{t}, 4t \InductionStepOrder{t,g} + t)$ such that $H_\mathfrak{M} - A$ is odd,\label{item:ABYSSFINISHodd}

        \item there exists a set $S \subseteq V(H_\mathfrak{M})$, with $|S| \leq \max(\InductionStepApex{t}, 4t \InductionStepOrder{t,g} + t)$, and there exist $h_1, h_2 \in \N$ with $h_1$ being positive and $1 < h_1+h_2 < t$, such that there exist disjoint even dicycles $C_1, \ldots , C_{h_2}$ and for each $i \in [h_1]$ there exists a circle $\theta_i$-outline $\mathfrak{M}_i^* = ( H_i^* , \mathcal{C}_i^*, \mathfrak{E}_i^* )$ within $H_\mathfrak{M} - S$ with $\theta_i \geq \InductionStepOutline{t - (h-1)}$, the set $\{ H_{\mathfrak{M}_1^*}, \ldots , H_{\mathfrak{M}_h^*} \}$ is half-integral, where $(H_{\mathfrak{M}_i^*}, \Omega_{\mathfrak{M}_i^*})$ is the $\mathfrak{M}_i^*$-society, the graph $H_{\mathfrak{M}_i^*} - V(H_i^*)$ is disjoint from all $C_1, \ldots , C_{h_2}$, and $H_{\mathfrak{M}_i^*}$ contains an even dicycle,\label{item:ABYSSFINISHeddies}
        
        \item there exists an $h \in \N$ with $1 < h < t$, pairwise disjoint disks $d_1^*, \ldots , d_h^* \subseteq m$, and a set $S \subseteq V(H_\mathfrak{M})$ with $|S| \leq \InductionStepApex{t}$, such that for each $i \in [h]$ there exists a circle, diamond, or refined diamond $\theta_i$-outline $\mathfrak{M}_i^* = ( H_i^* , \mathcal{C}_i^*, \mathfrak{E}_i^*, \mathfrak{B}_i^*, \mathfrak{V}_i^*, \mathfrak{D}_i^* )$ (with roughness $g$) within $H_\mathfrak{M} - S$ with $\theta_i \geq \InductionStepOutline{t - (h-1)}$, the set $\{ H_1^*, \ldots , H_h^* \}$ is third-integral, the trace of $d_i^*$ corresponds to the $\mathfrak{M}_i^*$-society $(H_{\mathfrak{M}_i^*}, \Omega_{\mathfrak{M}_i^*})$, and $H_{\mathfrak{M}_i^*}$ contains an even dicycle, and\label{item:ABYSSFINISHsplitmaelstrom}
        
        \item there exists a set $A \subseteq V(H_\mathfrak{M})$ with $|A| \leq \InductionStepApex{t}$ such that $H_\mathfrak{M} - A$ has an odd decomposition $\delta'$ with a unique maelstrom $m'$ and an associated circle, diamond, or refined diamond $\InductionStepOutline{t-1}$-outline $\mathfrak{M}'$ (of roughness $g$), where $(H_{\mathfrak{M}'},\Omega_{\mathfrak{M}'})$ is the $\mathfrak{M}'$-society, $H_{\mathfrak{M}'}$ contains an even dicycle, there exists an odd rendition $\rho'$ of $( H_\mathfrak{M} - (V(H_{\mathfrak{M}'} - V(\Omega_{\mathfrak{M}'})) \cup A), \Omega_\mathfrak{M}, \Omega_{\mathfrak{M}'} )$, and $H_\mathfrak{M} - V(H_{\mathfrak{M}'} - V(\Omega_{\mathfrak{M}'}))$ contains an even dicycle (which must contain vertices of $A$).\label{item:ABYSSFINISHsplitoffdicycle}
    \end{enumerate}
\end{corollary}

\subsection{The proof of the local structure theorem}

We can now finally gather up all of our tools and prove \Cref{thm:localstructure}, which we restate here for the convenience of the reader.
After the proof we give concrete definitions for $\LocalStructureApexNoArg$ and $\LocalStructureOutlineNoArg$.

\begin{theorem}[Local structure theorem]\label{thm:localstructure}
    There exist two functions $\LocalStructureOutlineNoArg \colon \N \rightarrow \N$ and $\LocalStructureApexNoArg \colon \N \rightarrow \N$ such that for any integer $t$ and any digraph $D$ containing a cylindrical $\LocalStructureOutline{t}$-wall $W'$ one of the following holds:
    \begin{itemize}
        \item $D$ contains a quarter-integral packing of $t$ even dicycles, or

        \item there exists a set of vertices $A \subseteq V(D)$, with $|A| \leq \LocalStructureApex{t}$, such that the $W'$-component of $D$ under $A$ has a maelstrom-free, pure odd decomposition and is thus odd.
    \end{itemize}
\end{theorem}

\begin{proof}[Proof of \Cref{thm:localstructure}]
    Let $W'$ be the $\LocalStructureOutline{t}$-wall promised in the statement.
    For the purpose of this proof we set $g = 1$ as the roughness of any potential refined diamond outlines.
    To give a first bound on $\LocalStructureOutlineNoArg$, we demand that
    \[ \LocalStructureOutline{t} \geq \FindOddWallOrder{2 \InductionStepOutline{t} + \SeparatingWallApex{t}} . \]
    We begin by using \Cref{cor:findoddwall} and subsequently \Cref{lem:jumpsoverawall} on $W'$, thus yielding either a quarter-integral packing of $t$ even dicycles, or a separating, odd $2 \InductionStepOutline{t}$-wall $W$ under an apex set $A_W$ of size at most $\FindOddWallApex{t} + \SeparatingWallApex{t}$, which tells us that we must demand
    \[ \LocalStructureApex{t} \geq \FindOddWallApex{t} + \SeparatingWallApex{t} . \]
    Since $W$ is a cylindrical wall, we have $W = ( Q_1, \ldots , Q_{2\InductionStepOutline{t}}, \hat{P}_1, \ldots , \hat{P}_{2\InductionStepOutline{t}} )$.
    Furthermore, according to \Cref{def:oddwall}, there is an undirected separation $(X,Y)$ of $D$ such that $X\cap Y=A\cup\Perimeter{W}$, $W \subseteq Y$, and every vertex in $Y$ reaches a vertex of $W-\Perimeter{W}$ or is reachable from it, and the cylindrical society $(\InducedSubgraph{D}{Y - A_W},\Omega_1,\Omega_2)$ with $\V{\Omega_1}=\V{Q_1}$ and $\V{\Omega_2}=\V{Q_{2\InductionStepOutline{t}}}$ has an odd rendition in the disk.
    In particular, we can use $Q_1$ and $Q_k$ to respectively define disks $m$ and $m'$ into which we draw the remainder of the graph, namely $X-A_W$, which yields an odd decomposition $\delta$ of $D - A_W$ with two maelstroms $m$ and $m'$.

    At this point, we note that both $Q_1, \ldots , Q_{\InductionStepOutline{t}}$ and $Q_{\InductionStepOutline{t} + 1}, \ldots , Q_{2\InductionStepOutline{t}}$ form disjoint sets of homogeneous dicycles living in the same odd rendition.
    Thus we can define two disjoint circle $\InductionStepOutline{t}$-outlines $\mathfrak{M} = (H, \mathcal{C}, \mathcal{E})$ and $\mathfrak{M}' = (H', \mathcal{C}', \mathcal{E}')$, where $\mathfrak{M}$ belongs to $m$ and $\mathfrak{M}'$ belongs to $m'$.
    By definition we therefore have $V(\Omega_\mathfrak{M}) = V(Q_{\InductionStepOutline{t}})$ and $V(\Omega_{\mathfrak{M}'}) = V(Q_{\InductionStepOutline{t} + 1})$, where $(H_\mathfrak{M}, \Omega_\mathfrak{M})$ is the $\mathfrak{M}$-society and $(H_{\mathfrak{M}'}, \Omega_{\mathfrak{M}'})$ is the $\mathfrak{M}'$-society.
    The path pairs $\hat{P}_1, \ldots , \hat{P}_{2\InductionStepOutline{t}}$ additionally give us two segregated $\nicefrac{\InductionStepOutline{t}}{2}$-pairs for the two outlines.

    If either $H_\mathfrak{M}$ or $H_{\mathfrak{M}'}$ is odd, we can integrate the associated maelstrom into the existing odd decomposition according to \Cref{lem:oddmaelstromdecomposition}.
    Should both be odd, we are already done and have found a pure odd decomposition.
    On the other hand, if both contain an even dicycle, we already know that we can at least find two disjoint even dicycles.
    Thus we can assume that at least one of the graphs $H_\mathfrak{M}$ and $H_{\mathfrak{M}'}$ contains an even dicycle.
    At this point our intermediate goal will be to either find the apex set of bounded size giving us a pure, odd decomposition for the graph described by the outline we are currently working on or to end up with $t-1$ or less outlines with bounded depth.
    Note that we are still working under the assumption that there does not exist a quarter-integral packing of $t$ even dicycles in $D$.

    Suppose w.l.o.g.\ that $H_\mathfrak{M}$ contains an even dicycle.
    We start by applying \Cref{cor:localstructureinductionstepgeneral} to $\mathfrak{M}$.
    Of course, we are very happy if \cref{item:packing} of \Cref{cor:localstructureinductionstepgeneral} is the result.
    If either \cref{item:depth} or \cref{item:odd} are the outcome, we have reached our intermediate goal, if $H_{\mathfrak{M}'}$ is odd, or we proceed with the same approach for $\mathfrak{M}'$.
    Should we instead get the even dicycle and the $\InductionStepOutline{t-1}$-outline promised in \cref{item:splitoffdicycle}, we can continue our approach in the new outline, whilst working with $t-1$ from now on.
    Note that this case cannot occur more than $t-2$ times in total and thus the union of the apex sets we incur here combine to contain at most $(t-2) \InductionStepApex{t}$ vertices, providing another lower bound for $\LocalStructureApex{t}$.

    Thus the last case left to check is encapsulated in \cref{item:splitmaelstrom}, which means that there exists an $h \in \N$ with $1 < h < t$, pairwise disjoint disks $d_1, \ldots , d_h \subseteq m$, and a set $S \subseteq V(H_\mathfrak{M})$ with $|S| \leq \InductionStepApex{t}$, such that for each $i \in [h]$ there exists a $\theta_i$-outline $\mathfrak{M}_i = ( H_i , \mathcal{C}_i, \mathfrak{E}_i, \mathfrak{B}_i, \mathfrak{V}_i, \mathfrak{D}_i )$ (with roughness $g$) for $d_i' \subseteq d_i$ within $H_\mathfrak{M} - S$ with $\theta_i \geq \InductionStepOutline{t - (h - 1)}$, the set $\{ H_1, \ldots , H_h \}$ is third-integral, the trace of $d_i$ corresponds to the $\mathfrak{M}_i$-society $(H_{\mathfrak{M}_i}, \Omega_{\mathfrak{M}_i})$, and the graph $H_{\mathfrak{M}_i}$ contains an even dicycle.
    We note that since $d_1, \ldots , d_h$ are pairwise disjoint, we also know that $H_{\mathfrak{M}_1}, \ldots , H_{\mathfrak{M}_h}$ are pairwise disjoint graphs.
    Therefore we can continue our procedure independently on $\mathfrak{M}_1, \ldots, \mathfrak{M}_h$, whilst adjusting our parameter such that we start with $t - (h-1)$ for each outline.
    Since $t > t - (h-1)$, as $h > 1$, this process will therefore terminate after a finite number of steps.

    Note that if at any point of the procedure we have found $t$ or more outlines hosting an even dicycle, we can find a third-integral packing of $t$ even dicycles, since $\InductionStepOutline{1} \geq \SC$ and we can simply apply \Cref{lem:applyshifting} to each.
    (This may result in a third-integral packing instead of an integral packing, since the support graphs $H_1, \ldots , H_h$ produced in \cref{item:splitmaelstrom} are third-integral, even if $H_{\mathfrak{M}_1}, \ldots , H_{\mathfrak{M}_h}$ are pairwise disjoint.)
    Thus we can assume that we do not reach \cref{item:odd} more than $t - 2$ times.
    We can also clearly continue with the same procedure if $H_{\mathfrak{M}'}$ contains an even dicycle.
    
    In total, we will accumulate an apex set $A'$ whose total size gives another lower bound for $\LocalStructureApex{t}$:
    \[ \LocalStructureApex{t} \geq \FindOddWallApex{t} + \SeparatingWallApex{t} + (t-2)\InductionStepApex{t} + (t-2)\InductionStepApex{t} + ({\log t})\InductionStepApex{t} + (t-1)\InductionStepApex{t} , \]
    which already includes the bound incurred by $A_W$,\footnote{Note that this is a very rough estimate, which simply sums up the obvious upper bounds for the number of times each item of \Cref{cor:localstructureinductionstepgeneral} can occur.} and unless $D - A'$ already has a pure odd decomposition, we know that there exists an $\ell \in \N$ with $1 \leq \ell < t$, pairwise disjoint disks $c_1, \ldots , c_\ell \subseteq m$, such that for each $i \in [\ell]$ there exists a $\theta_i'$-outline $\mathfrak{M}_i' = ( H_i' , \mathcal{C}_i', \mathfrak{E}_i', \mathfrak{B}_i', \mathfrak{V}_i', \mathfrak{D}_i' )$ (with roughness $g$) and depth less than $\InductionStepOrder{t - (\ell - 1),g}$ for $c_i' \subseteq c_i$ within $H_\mathfrak{M} - A'$ with $\theta_i' \geq \InductionStepOutline{t - (\ell - 1)}$, the set $\{ H_1', \ldots , H_\ell' \}$ is third-integral, the trace of $c_i$ corresponds to the $\mathfrak{M}_i'$-society $(H_{\mathfrak{M}_i'}, \Omega_{\mathfrak{M}_i'})$, and the graph $H_{\mathfrak{M}_i'}$ contains an even dicycle.

    Note that the condition of \Cref{lem:killcirclemaelstrom} that demands everything outside of the maelstroms to be odd can be easily realised by first removing all graphs associated with our outlines and then iteratively reintroducing and turning them odd via \Cref{lem:killcirclemaelstrom}.
    For any circle outline amongst $\mathfrak{M}_1', \ldots, \mathfrak{M}_\ell'$, we can thus use \Cref{lem:killcirclemaelstrom}.
    If this results in us finding a quarter-integral packing of $t$ even dicycles, we are done.
    Otherwise, we incur at most $(t-1) (4t\InductionStepOrder{t,0} + t)$ new apex vertices, which provides another update to the lower bound for $\LocalStructureApex{t}$.

    All that remains is for us to deal with refined diamond outlines of bounded depth.
    This will, as so often, be the most arduous part of our proof, requiring us to delve deeply into the structure of an outline one last time.
    Suppose that $\mathfrak{M}_1'$ is a refined diamond $\theta_i'$-outline of roughness $g$ and let $\mathfrak{W}$ be the diamond $\theta_i'$-outline that forms the basis of $\mathfrak{M}_1'$, with $( H_\mathfrak{W}, \Omega_\mathfrak{W} )$ being the $\mathfrak{W}$-society.
    Let $t' = t - (\ell - 1)$.
    Applying \Cref{lem:killdiamondmaelstrom} to $\mathfrak{M}_1'$ could give us a quarter-integral packing of $t$ even dicycles, which means we are done.
    Otherwise, we find both a set of apices $S$ of order at most $4 t \InductionStepApex{t,g} + t $ and a set of $k$ shorelines $S_1', \ldots , S_k'$ on the diamond $C$ corresponding to $\Omega_\mathfrak{W}$, with $k \leq t'$, such that the union of strong components of $D - (A' \cup S \cup \bigcup_{j = 1}^k V(S_j') \cup \bigcup_{j = 2}^\ell V(H_{\mathfrak{M}_j'} - V(\Omega_{\mathfrak{M}_j'})) )$ containing vertices of $V(D - H_\mathfrak{W}) \cup (V(\Omega_\mathfrak{W}) \setminus \bigcup_{j = 1}^k V(S_j'))$ is odd.
    We distinguish three cases.

    \textbf{Case 1:} We do not find any noteworthy shorelines, meaning $k = 0$.
    
    In this case, the outcome of \Cref{lem:killdiamondmaelstrom} is analogous to \Cref{lem:killcirclemaelstrom} and we can proceed further.

    \textbf{Case 2:} There is exactly one noteworthy shoreline, meaning $k = 1$.

    In this case $\mathfrak{M}_1'$ has undertow and we can apply \Cref{cor:exploringtheabyss}.
    If this results in \cref{item:ABYSSFINISHpacking} or \cref{item:ABYSSFINISHodd}, we are done (at least with this maelstrom) and in all other cases we iterate our procedure.

    \textbf{Case 3:} There are at least two noteworthy shorelines, meaning $k \geq 2$.

    We recall the definition of $S_H(I)$ from the proofs of \Cref{lem:killcirclemaelstrom} and \Cref{lem:killdiamondmaelstrom}, which is defined such that for a subgraph $H \subseteq D$ and a segment $I \subseteq V(\Omega_\mathfrak{W})$ the graph $S_H(I)$ is the union of strong components in $H$ that contain vertices of $I$.
    For each $i \in [k]$, we let $U_i = S_{H_\mathfrak{W} - (A' \cup S)}(S_i')$ and note that the graphs $U_1, \ldots , U_k$ are pairwise disjoint.

    Consider the $\nicefrac{\theta_1'}{2}$-rim $\mathfrak{R}_1' = ( U_1' , \mathcal{R}_1', \mathcal{F}_1', \mathfrak{A}_1', \mathfrak{D}_1' )$ of $\mathfrak{M}_1'$, with degree $h'$, $\theta' = \nicefrac{\theta_1'}{2}$, $\mathcal{R}_1' = \{ R_1, \ldots , R_{\theta'} \}$, $\mathcal{F}_1' = \{ F_1, \ldots , F_{\theta'} \}$, $\mathfrak{A}_1' = \{ \mathcal{A}_1^1, \ldots , \mathcal{A}_{\theta'}^1 \}$, and $\mathcal{A}_j^1 = \{ A_1^1, \ldots , A_{(g + 1)\theta'}^1 \}$ for all $j \in [h']$.
    We also note that $C$, which is the diamond corresponding to $\Omega_\mathfrak{W}$, consists of two directed paths $P,Q$ who share their endpoints.
    These paths can be chosen such that $P \subseteq R_{\theta'}$ and $Q \subseteq F_{\theta'}$.
    Note that for each $i \in [k]$ we either have $S_i' \subseteq P$ or $S_i' \subseteq Q$.

    For each $i \in [k]$ we construct a new circle outline as follows.
    Suppose first that $S_i' \subseteq P$, then we let $\mathcal{R}_i^1$ be the set of dicycles consisting of $R_{\theta' - ((k-(i+1))  \theta'') + 1} , \ldots , R_{\theta' - ((k-i) \theta'')}$ and note that, if we let $U_i = \bigcup_{j = \theta' - ((k-(i+1))  \theta'') + 1}^{\theta' - ((k-i) \theta'')} R_j \cup (H_\mathfrak{W} \cap U_1)$, we can define the circle $\theta''$-outline $\mathfrak{W}_i^1 = ( U_i, \mathcal{R}_i^1, \emptyset )$ for $m$ (or expansions of this disk) in $U_i \cup U_i' \subseteq H_\mathfrak{W} \cup H_i'$.
    If $S_i' \subseteq Q$, we first note that we can find a homogeneous set of dicycles $R_{\theta' - ((k-(i+1))  \theta'') + 1}' , \ldots , R_{\theta' - ((k-i) \theta'')}'$ in $\bigcup_{j = \theta' - ((k-(i+1))  \theta'') + 1}^{\theta' - ((k-i) \theta'')} (R_j \cup F_j)$ via the use of \Cref{lem:detour}, which form the set $\mathcal{R}_i^1$ and allow us to construct $U_i$ as $\bigcup_{j = \theta' - ((k-(i+1))  \theta'') + 1}^{\theta' - ((k-i) \theta'')} R_j' \cup (H_\mathfrak{W} \cap U_1)$ and we can analogously define the circle $\theta''$-outline $\mathfrak{W}_i^1 = ( U_i, \mathcal{R}_i^1, \emptyset )$ for $m$ (or expansions of this disk) in $U_i \cup U_i' \subseteq H_\mathfrak{W} \cup H_i'$.
    Note that we can use the odd decomposition of $\mathfrak{M}_1'$ to verify the properties we need and find odd decompositions associated with $\mathfrak{W}_1^1, \ldots , \mathfrak{W}_k^1$.

    By definition we conclude that $U_1 \cup U_1', \ldots , U_k \cup U_k'$ is half-integral, with intersections occurring for pairs $U_i$ and $U_j'$ with $i \neq j$.
    On the outlines $\mathfrak{W}_1^1, \ldots , \mathfrak{W}_k^1$ we can now again apply \Cref{lem:localstructureinductionstepgeneral} with the parameter $t'$ and continue our procedure, due to our previous bounds on $\LocalStructureOutlineNoArg$ and the fact that  for all $\ell' \in [t-1]$ we have
    \[ \nicefrac{\InductionStepOutline{t}}{2^{\ell' + 1}} \geq \InductionStepOutline{t - \ell'} . \]
    Note that, if we find a total of $t$ outlines whose infrastructure is at worst quarter-integral, our previous bounds on $\LocalStructureApexNoArg$ still largely suffice, with only the apex set from \Cref{lem:killdiamondmaelstrom} getting added as follows:
    \[ \LocalStructureApex{t} \geq \FindOddWallApex{t} + \SeparatingWallApex{t} + (3t-5 + {\log t})\InductionStepApex{t} + (t-1) (4t\InductionStepOrder{t,g} + t) . \]
    Thus in all possible cases, we can keep our procedure going and only incur additional vertices that must be deleted when we are able to increase the potential size of a quarter-integral packing of even dicycles.
    This completes the proof of our theorem.
\end{proof}
As a consequence of the proof and in particular the fact that $g$ was chosen as a constant, we can give the following definition for the two functions $\LocalStructureOutlineNoArg \colon \N \rightarrow \N$ and $\LocalStructureApexNoArg \colon \N \rightarrow \N$:
\begin{align*}
    \LocalStructureApex{t}      \coloneqq \ & \FindOddWallApex{t} + \SeparatingWallApex{t} + (3t-5 + {\log t})\InductionStepApex{t} + (t-1) (4t\InductionStepOrder{t,g} + t) \text{ and } \\
    \LocalStructureOutline{t}   \coloneqq \ & \FindOddWallOrder{2 \InductionStepOutline{t} + \SeparatingWallApex{t}}.
\end{align*}


\section{Global structure theorem}\label{sec:globalstructure}
\newcommand{\drop}[1]{}
\drop{
\newcommand{\TTT}{\mathcal{T}}
\newcommand{\SSS}{\mathcal{S}}
\newcommand*{\GlobalStructureNoArg}{\mathsf{f}_{\ref{thm:globalstructure}}}
\newcommand*{\GlobalStructure}[1]
}
\subsection{More tools from directed structure theory}

We adapt the concept of directed tangles from \cite{giannopoulou2020canonical}.
Giannopoulou, Kawarabayashi, Kreutzer, and Kwon define directed separations in a slightly more liberal way, that is they do not let the tuple notation encode in which direction the edges avoiding the separator are allowed to go.
This makes for a more compact representation for tangles.
In this paper we will instead employ the notation of Erde \cite{erde2020directed} which was used to describe the cousin of directed tangles for directed pathwidth: \textit{diblockages}.

Let $D$ be a digraph and $k$ be a positive integer.
By $\mathcal{S}_k(D)$ we denote the collection of all directed separations of order $<k$ for $D$.
In case the digraph is understood from the context, we drop the argument and simply write $\mathcal{S}_k$.
An \emph{orientation} of $\mathcal{S}_k$ is a bipartition of $\mathcal{S}_k$ into the sets $\mathcal{O}^{\mathsf{out}}$ and $\mathcal{O}^{\mathsf{in}}$.
Given some orientation $\mathcal{O}=(\mathcal{O}^{\mathsf{out}},\mathcal{O}^{\mathsf{in}})$ we associate to it two functions which assign to every $(A,B)\in\mathcal{S}_k$ a \textit{small side} and a \textit{big side}.
\begin{align*}
    \mathsf{big}_{\mathcal{O}}((A,B))&\coloneqq 
    \begin{cases}
    A\text{, if }(A,B)\in \mathcal{O}^{\mathsf{out}}\text{, and}\\
    B\text{ otherwise.}
    \end{cases}\\
    \mathsf{small}_{\mathcal{O}}((A,B))&\coloneqq 
    \begin{cases}
    B\text{, if }(A,B)\in \mathcal{O}^{\mathsf{out}}\text{, and}\\
    A\text{ otherwise.}
    \end{cases}
\end{align*}
If the orientation is understood from the context we drop the subscript for both functions.
Notice that $\mathsf{big}_{\mathcal{O}}$ fully determines $\mathsf{small}_{\mathcal{O}}$ and vice versa.
Moreover, the orientation $\mathcal{O}=(\mathcal{O}^{\mathsf{out}},\mathcal{O}^{\mathsf{in}})$ is also completely represented by $\mathsf{big}_{\mathcal{O}}$.
This gives us the right to fully treat an orientation as a mapping $\mathsf{big}$ which maps each $(A,B)\in\mathcal{S}_k$ to exactly one of its sides and to assume that the corresponding function $\mathsf{small}$ is given implicitly.
If $\mathsf{big}((A,B))=X$ we refer to it as the \emph{big side}, and to $Y\in\{ A,B\}\setminus \{X\}$ we refer to as the \emph{small side}.

\begin{definition}[Tangle]\label{def:directedtangle}
    Let $D$ be a digraph and $k$ be a positive integer.
    An orientation $\mathsf{big}$ of $\mathcal{S}_k$ is called a \emph{tangle} of \emph{order $k$}, if for any three directed separations $(A_1, B_1), (A_2, B_2), (A_3, B_3) \in \mathcal{S}_k$ it holds that
    \[ \mathsf{small}((A_1,B_1))\cup \mathsf{small}((A_2,B_2)) \cup \mathsf{small}((A_3,B_3)) \neq V(D). \]   

    Let $\mathsf{big}_1$ be a tangle of order $k_1$ and $\mathsf{big}_2$ be a tangle of order $k_2<k_1$ in $D$.
    We say that $\mathsf{big}_2$ is a \emph{truncation} of $\mathsf{big}_1$ if for every $(A,B)\in\mathcal{S}_{k_2}$
    \[ \mathsf{big}_1((A,B))=\mathsf{big}_2((A,B)). \]
\end{definition}

To show that tangles act as obstructions to directed treewidth and how to use them to identify specific other objects, we adopt the strategies from \cite{giannopoulou2020canonical}.
However, for our purposes we need slightly stronger statements and thus provide the necessary proofs for the sake of completeness.

Let $k$ be a positive integer, $D$ be a digraph and $X\subseteq V(D)$ be a set of vertices.
A set $S\subseteq V(D)$ is a \emph{balanced separator} for $X$ if for every strong component $K$ of $D-S$ we have that $|X\cap V(K)|<\frac{2}{3}|X|$.
We say that $X$ is \emph{$k$-linked} if there is no balanced separator of order less than $k$ in $D$.

\begin{lemma}\label{lemma:linkedtotangle}
Let $k$ be a positive integer, $D$ be a digraph and $X\subseteq V(D)$ be a $k$-linked set.
Then there exists a tangle $\mathsf{big}_X$ of order $k$ in $D$ such that for every $(A,B)\in \mathcal{S}_k$, the unique component $K$ of $D-(A\cap B)$ with $|X\cap V(K)|\geq \frac{2}{3}|X|$ belongs to $\mathsf{big}_X((A,B))$.
\end{lemma}

\begin{proof}
Let us simply define $\mathsf{big}_X$ to satisfy the condition.
That is, let $\mathsf{big}_X$ be the orientation of $\mathcal{S}_k$ such that for every $(A,B)\in \mathcal{S}_k$, the unique component $K$ of $D-(A\cap B)$ with $|X\cap V(K)|\geq \frac{2}{3}|X|$ belongs to $\mathsf{big}_X((A,B))$.
We now just have to check whether three small sides can cover the entire vertex set of $D$.
For this let $(A_1,B_1),(A_2,B_2),(A_3,B_3)\in\mathcal{S}_k$ be any three directed separations and for each $i\in[3]$ let $Z_i$ denote $A_i\cap B_i$ and $L_i$ denote the small side of $(A_i,B_i)$.
Observe that $|X\cap V(J)|<\frac{1}{3}|X|$ for any component $J$ of $D[L_i]-Z_i$.
Indeed, it even holds that $|L_i\cap X|<\frac{1}{3}|X|$.
This, however, means that $|(L_1\cup L_2\cup L_3)\cap X|<|X|$ and thus there must exist some vertex of $X$, and therefore of $V(D)$, which does not belong to the union of the three small sides.
\end{proof}

We use a definition for well-linked sets in digraphs by Reed \cite{Reed1999Introducing}.

\begin{definition}[Well-linked set]\label{def:welllinkedset}
    A \emph{well-linked set} of order $m$ in a digraph $D$ is a set of vertices $W \subseteq V(D)$ with $|W| = m$ such that for all subsets $A, B \subseteq W$ with $|A| = |B|$ there are $|A|$ disjoint paths from $A$ to $B$ in $D - (W \setminus (A \cup B))$.
\end{definition}

Similar to the way a tangle can be obtained from a $k$-linked set, we can associate a tangle with every large enough cylindrical wall.

\begin{lemma}\label{lem:walltangle}
    Let $W = (Q_1,\dots,Q_l,\hat{P}_1,\dots,\hat{P}_l)$ be a cylindrical wall of order $l\geq 3k$. Then for any separation $(A, B) 
    \in \SSS_k$ there is exactly one side $X \in \{ A, B \}$ such that $X \setminus (A\cap B)$ contains an entire cycle $Q_i$ of $W$. 
    
    Let $\mathsf{big}_W$ be the function that maps any separation $(A,B) \in \SSS_k$ to the side containing a cycle of $W$.
    Then $\mathsf{big}_W$ is a tangle.
\end{lemma}
\begin{proof}
    The first statement of the lemma is obvious as any pair $Q_i, Q_j$ of cycles in $W$ are linked by $l$ paths in each direction. Thus, if there was such a pair  $Q_i \subseteq A$ and $Q_j \subseteq B$ then $A \cap B$ would have to contain a vertex from each of the $l \geq 3k > |A\cap B|$ paths from $A$ to $B$, a contradiction to $|A \cap B| < k$.

    Thus the function $\mathsf{big}_W$ is well-defined. What remains is to prove that $\mathsf{big}_W$ is a tangle. Towards a contradiction let $(A_i, B_i) \in \SSS_k$ and let $L_i$ be the small side of $(A_i, B_i)$, for $i=1,2,3$.  If $L = \bigcup_{i=1}^3 L_i = V(D)$ then, in particular, $L$ must contain all cycles $Q_1, \dots, Q_l$. Let $X = \bigcup_{i=1}^3 A_i\cap B_i$. Then $|X| < 3k$ and therefore there is a cycle $Q_i$ such that $V(Q_i) \cap X = \emptyset$. But then $Q_i$ must be in the big side of all three separations and therefore $L \not= V(D)$.  
\end{proof}

The last result we need in preparation for the proof of the global structure theorem is the following lemma stating that with every $k$-linked set $Z$, for some large enough value of $k$, we can associate a wall $W$ generated from $Z$ such that the tangle $T_W$ induced by $W$ is a truncation of the tangle $T_Z$ induced by $Z$.
\begin{lemma}\label{lem:welllinkedsetgiveswallandtangles}
    There is a function $f_{\ref{lem:welllinkedsetgiveswallandtangles}}$ such that for all $k \geq 1$, if $Z \subseteq V(D)$ is a $f_{\ref{lem:welllinkedsetgiveswallandtangles}}(k)$-linked set then $D$ contains a cylindrical wall $W$ of order $k$ such that $T_W$ is a truncation of $T_Z$, where $T_Z, T_W$ are the tangles induced by $Z$ and $W$, respectively. 
\end{lemma}
\begin{proof}
    Let $f'$ be the function of the directed grid theorem as defined in \cref{thm:fptdirectedwall}. In \cite{KawarabayashiKreutzer2015DirectedGrid}, Kawarabayashi and Kreutzer construct a cylindrical wall of order $k$ in any digraph of directed tree-width at least $f'(k)$. The proof starts with an obstruction to directed tree-width called \emph{directed brambles}. A directed bramble $\mathcal{B}$ is a set of strongly connected subgraphs of $D$ which pairwise intersect. Its order is the minimal size of any set $X$, called a \emph{cover} or \emph{hitting set} of $\mathcal{B}$, which intersects all bramble elements. 

    Any $\ell$-linked set $X$ induces a bramble $\mathcal{B}_X$ of order $>\ell$ by adding for each $S \subseteq V(D)$ of order $|S|<\ell$ the strong component of $D-S$ to $\mathcal{B}_X$ that contains the majority of the elements of $X$. It can be shown that a minimum order cover of any bramble is well-linked.
    
    Kawarabayashi and Kreutzer  \cite{KawarabayashiKreutzer2015DirectedGrid} construct the cylindrical wall of order $3k$ from the bramble $\mathcal{B}_Z$ obtained from an $f'(3k)$-linked set $Z$ by dividing the well-linked cover of $\mathcal{B}_Z$ into subsets and taking linkages between pairs of these subsets. The cylindrical wall can then be found within these linkages in one of several possible ways.
    However, it is a direct consequence of their proof that the wall $W$ obtained in this way is well-connected to the set $Z$, that is, there cannot be a separation $(A, B) \in \SSS_k$ which contains an entire cycle of $W$ in one side but the majority of the elements of $Z$ in the other side. 
    But this immediately implies that the tangle $T_W$ induced by $W$ is a truncation of the tangle $T_Z$ induced by $Z$.
\end{proof}

\subsection{The proof of the global structure theorem}

Having proven the relevant lemmas, we can now construct our global decomposition.
We restate \Cref{thm:globalstructure} here for the convenience of the reader.

\newcommand{\fX}[1]{f_X(#1)}
\newcommand{\fI}[1]{\LocalStructureOutline{#1}}
\newcommand{\fII}[1]{\LocalStructureApex{#1}}
\newcommand{\fV}[1]{ f_{\ref{lem:welllinkedsetgiveswallandtangles}}(#1)}

\begin{theorem}\label{thm:globalstructure}
    There exists a function $\GlobalStructureNoArg \colon \N \rightarrow \N$ such that for every digraph $D$ and every $k \in \N$ either
    \begin{enumerate}
        \item $D$ has a quarter-integral packing of $k$ even dicycles or

        \item $D$ has a directed tree-decomposition $(T, \beta , \gamma )$ and a map $\alpha \colon V(T) \rightarrow 2^{V(D)}$ such that
        \begin{enumerate}
            \item $|\gamma(e)| \leq \GlobalStructure{k}$ for all $e \in E(T)$,

            \item $\alpha(t) \subseteq \Gamma(t)$ for all $t \in V(T)$,

            \item $|\alpha(t)| \leq \GlobalStructure{k}$ for all $t \in V(T)$, and

            \item no even dicycle of $D - \alpha(t)$ contains a vertex of $\Gamma(t) \setminus \alpha(t)$.
        \end{enumerate}
    \end{enumerate}
\end{theorem}
\begin{proof}
By induction on $|V(D)|$ we prove the following slightly stronger claim which clearly implies the theorem by setting $\GlobalStructure{k} := 4\cdot \fV{3\cdot \fI{k}}$. 

\begin{claim}\label{claim:strongerglobalstructure}
    For every digraph $D$ and every $k \in \N$ either
    \begin{enumerate}
        \item $D$ has a quarter-integral packing of $k$ even dicycles or

        \item for any set $Z \subseteq V(D)$ of size at most $3 \cdot \fV{3\cdot \fI{k}}$ there exits a directed tree-decomposition $(T, \beta, \gamma)$ of $D$ with root $r$ and a map $\alpha \colon V(T) \rightarrow 2^{V(D)}$ such that
        \begin{enumerate}
        \item $|\gamma(e)| \leq 4\cdot\fV{3\cdot \fI{k}}$ for all $e \in E(T)$,
        \item $\alpha(t) \subseteq \Gamma(t)$ for all $t \in V(T)$,
        \item $|\alpha(t)| \leq 4\cdot\fV{3\cdot \fI{k}}$ for all $t \in V(T)$, 
        \item no even dicycle of $D - \alpha(t)$ contains a vertex of $\Gamma(t) \setminus \alpha(t)$, for all $t\in V(T)$
        \item $Z \subseteq \alpha(r) \subseteq \beta(r)$,
        \item $\gamma((d,t))  \subseteq \alpha(t)$ for all $(d, t) \in E(T)$, and
        \item $\gamma((d, t)) \cap \beta(T_t) = \emptyset$ for all $(d, t) \in E(T)$.
        \end{enumerate}
    \end{enumerate}
\end{claim}

\medskip

    If $|V(D)| < 3\cdot \fV{3\cdot \fI{k}}$ there is nothing to show.
    So we may assume that $|V(D)| \geq 3\cdot \fV{\fI{k}}$.
    If $|Z| < 3\cdot \fV{3\cdot \fI{k}}$ we add arbitrary vertices to $Z$ until $|Z| = 3\cdot \fV{3\cdot \fI{k}}$.

    \medskip

    \noindent\textbf{Case 1. }
    Suppose first that there exists $X \subseteq V(D)$ such that $|X| \leq \fV{3\cdot \fI{k}}$ and $(V(H) \cap Z) \leq \frac23|Z|$ for all strong components of $D-X$. 
    Let $H_1, \dots, H_l$ be the strong components of $D-X$.
    For $1 \leq i \leq l$ let $Z_i := V(H_i) \cap Z$.
    By induction on $H_i$ and $Z_i$ there is a directed tree-decomposition $(T_i, \beta_i, \gamma_i)$ with root $r_i$ and a map $\alpha_i \colon V(T) \rightarrow 2^{V(H_i)}$ satisfying the statement of the theorem.
    
    We construct a new directed tree-decomposition $(T, \beta, \gamma)$ as follows. 
    Let $T$ be obtained from the disjoint union of $T_1, \dots, T_l$ by adding a new vertex $r$ and edges $(r, r_i)$ for all $1 \leq i \leq l$. 
    We define 
    \begin{itemize}
    \item 
    $\beta(r) = \alpha(r) := Z \cup X$,  
    \item  $\gamma((r, r_i)) := Z_i \cup X$ and  
    $\beta(r_i) := \beta_i(r_i) \setminus Z_i$, for all $1 \leq i \leq l$, and 
    \item for all $t \in V(T_i)\setminus \{r_i\}$ and all $e \in E(T_i)$ we define $\beta(t) := \beta_i(t)$ and $\gamma(e) := \gamma_i(e)$. 
    \end{itemize}
    This completes the construction. We claim that $\TTT \coloneqq (T, \beta, \gamma)$ and $\alpha$ satisfy the claim of the theorem.
    It is easily verified that $\TTT$ is indeed a directed tree-decomposition.
    Furthermore, Properties a), b), c), and e) are clearly satisfied.
    We prove Property d) next.
    For $t=r$ there is nothing to show as $\beta(r) \subseteq \alpha(r)$.
    If $t\not=r$, Property d) follows from the induction hypothesis using Property f) to ensure that $\alpha(t)$ contains the guard of the incoming edge of $t$ and therefore $D-\alpha(t)$ has no cycle that contains vertices in $D[\beta(T_t)]$ and vertices not in $D[\beta(T_t)]$. 

    \medskip

    \noindent\textbf{Case 2. }
    So we may assume that there is no such set $X$ as above.
    This implies that $Z$ is $\fV{3\cdot \fI{k}}$-linked.
    By \cref{lem:welllinkedsetgiveswallandtangles} there exists a cylindrical wall $W$ of order $3\cdot \fI{k}$ such that the tangle $\TTT_W$ induced by $W$ is the truncation of the tangle $\TTT_Z$ induced by $Z$. 

    We now apply \cref{thm:localstructure} to $D$ and $W$ and either obtain a quarter-integral packing of $t$ even dicycles ---which is the first outcome of the theorem--- or a set $A \subseteq V(D)$ of order $|A| \leq \fII{k}$ such that the strong component $H_W$ of $W$ under $A$ is odd. 

    Let $H_1, \dots, H_l$ be the strong components of $G-A$ other than $H_W$.
    For $1 \leq i \leq l$ let $Z_i := (V(H_i) \cap Z) \cup A$.
    As $|A| \leq \fII{k}$, which is easily seen to be $\leq\fI{k}$, $|A|$ is smaller than the order of the tangle $\TTT_Z$, this implies that $V(H_i) \cup A$ is always contained in the small side with respect to $\TTT_Z$ of any separation with $V(H_i) \cup A$ on one side and $A \cup V(H_W)$ on the other side.
    Therefore $|V(H_i) \cap Z|  \leq \frac13|Z|$ and thus $|Z_i| \leq \frac13 \fV{3\cdot \fI{k}} + \fII{k} \leq \fV{3\cdot \fI{k}}$. 
    
    For each $1 \leq i \leq l$, we apply induction on $D[H_i \cup A]$ and $Z_i \cup A$ to obtain a directed tree-decomposition $\TTT_i := (T_i, \beta_i, \gamma_i)$ with root $r_i$ and a function $\alpha_i$ satisfying Properties a) - e). 

    We now construct a new directed tree-decomposition $(T, \beta, \gamma)$ as follows.
    Let $T$ be the tree obtained from the disjoint union of $T_1, \dots, T_l$ by adding a new root $r$ and the edges $(r, r_i)$ for all $1  \leq i \leq l$. 

    We define 
    \begin{itemize}
    \item 
    $\alpha(r) := A \cup Z$ and 
    \item $\beta(r) := V(H_W) \cup Z \cup A$.
    \item 
    For $1 \leq i \leq l$ we define 
    $\gamma((r, r_i)) := Z_i \cup A$, 
    $\alpha(r_i) := \alpha_i(r_i)$, and 
    $\beta(r_i) := \beta_i(r_i) \setminus (Z_i \cup A)$. 
    \item 
    For all $t \in V(T_i) \setminus \{ r_i\}$ and all $e \in E(T_i)$ we define 
    $\beta(t) := \beta_i(t)$ and 
    $\gamma(e) := \gamma_i(e)$. 
   
    \end{itemize}
    This completes the construction.

    We now verify that $(T, \beta, \gamma)$ is a directed tree-decomposition which together with $\alpha$ satisfies Properties a) - e) of our assumption.     

    By induction hypothesis, $\TTT_i$ satisfies Property e) and therefore $A \subseteq Z_i \subseteq \alpha_i(r_i)  \subseteq \beta_i(r_i)$ and $\bigcup \{ \beta_i(t) \colon t \in V(T_i) \} = V(H_i) \cup A$.
    Thus, as we remove $Z_i \cup A$ from $\beta_i(r_i)$ to get $\beta(r_i)$, every vertex of $D$ is contained in precisely one bag of the decomposition. 
    
    Now let $e \in E(T)$ be an edge. We need to show that $\gamma(e)$ guards $D[\beta(T_i)]$. 
    If $e = (r, r_i)$  then $D[\beta_i(T_i)]$ is  $H_i - Z_i$.
    But $H_i$ is a strong component of $D-A$ and $A \cup Z_i = \gamma(r, r_i)$. Thus $\gamma(r, r_i)$ guards $D[\beta(T_i)]$. 
    Otherwise, $e = (s, t) \in E(T_i)$ for some $1 \leq i \leq l$.
    By induction hypothesis, $\gamma_i(e)$ guards $\beta_i((T_i)_t)$ in $D[H_i \cup A]$.
    By Property e), $Z_i \cap \beta_i((T_i)_t) = \emptyset$, which implies that $\gamma_i(e) = \gamma(e)$ guards $\beta_i((T_i)_t) = \beta(T_t)$.
    So $\TTT$ is a directed tree-decomposition of $D$. 

    We need to verify that $\TTT$ and $\alpha$ satisfy Properties a) to e).
    Properties a), b), c), e), and g) follow immediately from the construction. 
    
    We verify Property f) next.
    For edges $(s, t) \in E(T_i)$ this follows immediately from the induction hypothesis.
    Now let $e = (r, r_i)$, for some $1\leq i \leq l$.
    Then $\gamma(e) = Z_i \cup A$ and $\alpha(r_i) = \alpha_i(r_i)$.
    But to obtain $\TTT_i$ we applied induction to $D[H_i \cup A]$ and $Z_i \cup A$.
    So Property e) applied to $\TTT_i$, i.e.~where $Z_i \cup A$ takes on the r\^ole of $Z$, guarantees that $Z_i \cup A \subseteq \alpha_i(r_i) = \alpha(r_i)$ which proves Property f).  
    
    It remains to verify Property d).
    That is, we need to show that for all $t \in V(T)$ no even dicycle of $D - \alpha(t)$ contains a vertex of $\Gamma(t) \setminus \alpha(t)$.
    We first consider the case where $t=r$.
    Then, by construction, $D[\beta(r)\setminus \alpha(r)] \subseteq D[\beta(r)\setminus A]$ is odd and therefore does not have any even dicycle.
    Furthermore, $\gamma(r, r_i) \subseteq A \cup Z = \alpha(r)$, for all $1 \leq i \leq l$, and therefore no cycle, even or not, in $D-\alpha(r)$ that contains a vertex outside of $D[\beta(r) \setminus \alpha(r)]$ can contain a vertex of $\Gamma(r)\setminus \alpha(r)$. 
  
    Now suppose $t \not= r$.
    Then Property d) follows from the induction hypothesis using Property d) and the fact that by Property f) the set $\gamma(s, t)$ guarding the incoming edge of $t$ is a subset of $\alpha(t)$.
    This concludes the proof of the global structure theorem.
\end{proof}



\bibliographystyle{acm}
\bibliography{literature}

\begin{thebibliography}{10}

\bibitem{BermondT1981Cycles}
{\sc Bermond, J.~C., and Thomassen, C.}
\newblock Cycles in digraphs - a survey.
\newblock {\em Journal of Graph Theory 5\/} (1981), 1--43.

\bibitem{Bjorklund2022ShortestEvenDicycle}
{\sc Bj\"{o}rklund, A., Husfeldt, T., and Kaski, P.}
\newblock The shortest even cycle problem is tractable.
\newblock In {\em S{TOC} '22---{P}roceedings of the 54th {A}nnual {ACM}
  {SIGACT} {S}ymposium on {T}heory of {C}omputing\/} ([2022] \copyright 2022),
  ACM, New York, pp.~117--130.

\bibitem{campos2022adapting}
{\sc Campos, V., Lopes, R., Maia, A.~K., and Sau, I.}
\newblock Adapting the directed grid theorem into an fpt algorithm.
\newblock {\em SIAM Journal on Discrete Mathematics 36}, 3 (2022), 1887--1917.

\bibitem{Chudnovsky2008APathAlgorithm}
{\sc Chudnovsky, M., Cunningham, W.~H., and Geelen, J.}
\newblock An algorithm for packing non-zero {$A$}-paths in group-labelled
  graphs.
\newblock {\em Combinatorica 28}, 2 (2008), 145--161.

\bibitem{Chudnovsky2006PackingAPaths}
{\sc Chudnovsky, M., Geelen, J., Gerards, B., Goddyn, L., Lohman, M., and
  Seymour, P.}
\newblock Packing non-zero {$A$}-paths in group-labelled graphs.
\newblock {\em Combinatorica 26}, 5 (2006), 521--532.

\bibitem{Conforti2020StableSetGenusOddCycles}
{\sc Conforti, M., Fiorini, S., Huynh, T., Joret, G., and Weltge, S.}
\newblock The stable set problem in graphs with bounded genus and bounded odd
  cycle packing number.
\newblock In {\em Proceedings of the 2020 {ACM}-{SIAM} {S}ymposium on
  {D}iscrete {A}lgorithms\/} (2020), SIAM, Philadelphia, PA, pp.~2896--2915.

\bibitem{Conforti2020ExtendedFormulationsStableSet}
{\sc Conforti, M., Fiorini, S., Huynh, T., and Weltge, S.}
\newblock Extended formulations for stable set polytopes of graphs without two
  disjoint odd cycles.
\newblock In {\em Integer programming and combinatorial optimization},
  vol.~12125 of {\em Lecture Notes in Comput. Sci.} Springer, Cham, [2020]
  \copyright 2020, pp.~104--116.

\bibitem{CurticapeanXia2015ParameterizingThePermanent}
{\sc Curticapean, R., and Xia, M.}
\newblock Parameterizing the permanent: genus, apices, minors, evaluation
  {$\textup{mod}\, 2^k$}.
\newblock In {\em 2015 {IEEE} 56th {A}nnual {S}ymposium on {F}oundations of
  {C}omputer {S}cience---{FOCS} 2015}. IEEE Computer Soc., Los Alamitos, CA,
  2015, pp.~994--1009.

\bibitem{cygan2013planar}
{\sc Cygan, M., Marx, D., Pilipczuk, M., and Pilipczuk, M.}
\newblock {T}he {P}lanar {D}irected $k$-{V}ertex-{D}isjoint {P}aths {P}roblem
  is {F}ixed-{P}arameter {T}ractable.
\newblock In {\em 2013 IEEE 54th Annual Symposium on Foundations of Computer
  Science\/} (2013), IEEE, pp.~197--206.

\bibitem{erde2020directed}
{\sc Erde, J.}
\newblock {D}irected {P}ath-{D}ecompositions.
\newblock {\em SIAM Journal on Discrete Mathematics 34}, 1 (2020), 415--430.

\bibitem{ErdosPosa1965OnIndependetCircuits}
{\sc Erd{\H{o}}s, P., and P\'{o}sa, L.}
\newblock On independent circuits contained in a graph.
\newblock {\em Canadian J. Math. 17\/} (1965), 347--352.

\bibitem{Even1975Timetable}
{\sc Even, S., Itai, A., and Shamir, A.}
\newblock On the complexity of timetable and multi-commodity flow problems.
\newblock In {\em 16th {A}nnual {S}ymposium on {F}oundations of {C}omputer
  {S}cience ({U}niv. {C}alifornia, {B}erkeley, {C}alif., 1975)}. IEEE, Long
  Beach, CA, 1975, pp.~184--193.

\bibitem{Fortune1980DirectedSubgraphHomeomorphism}
{\sc Fortune, S., Hopcroft, J., and Wyllie, J.}
\newblock The directed subgraph homeomorphism problem.
\newblock {\em Theoret. Comput. Sci. 10}, 2 (1980), 111--121.

\bibitem{Geelen2009OddHadwiger}
{\sc Geelen, J., Gerards, B., Reed, B., Seymour, P., and Vetta, A.}
\newblock On the odd-minor variant of {H}adwiger's conjecture.
\newblock {\em J. Combin. Theory Ser. B 99}, 1 (2009), 20--29.

\bibitem{Giannopoulou2020DirectedFlatWall}
{\sc Giannopoulou, A.~C., Kawarabayashi, K.-i., Kreutzer, S., and Kwon, O.-j.}
\newblock The directed flat wall theorem.
\newblock In {\em Proceedings of the 2020 {ACM}-{SIAM} {S}ymposium on
  {D}iscrete {A}lgorithms\/} (2020), SIAM, Philadelphia, PA, pp.~239--258.

\bibitem{giannopoulou2020canonical}
{\sc Giannopoulou, A.~C., Kawarabayashi, K.-i., Kreutzer, S., and Kwon, O.-j.}
\newblock Directed tangle tree-decompositions and applications.
\newblock In {\em Proceedings of the 2022 {A}nnual {ACM}-{SIAM} {S}ymposium on
  {D}iscrete {A}lgorithms ({SODA})\/} (2022), [Society for Industrial and
  Applied Mathematics (SIAM)], Philadelphia, PA, pp.~377--405.

\bibitem{Giannopoulou2024ExcludingPlanar}
{\sc Giannopoulou, A.~C., Kreutzer, S., and Wiederrecht, S.}
\newblock Excluding a planar matching minor in bipartite graphs.
\newblock {\em J. Combin. Theory Ser. B 164\/} (2024), 161--221.

\bibitem{giannopoulou2022excluding}
{\sc Giannopoulou, A.~C., Thilikos, D.~M., and Wiederrecht, S.}
\newblock Excluding single-crossing matching minors in bipartite graphs.
\newblock {\em arXiv preprint arXiv:2212.09348\/} (2022).

\bibitem{Giannopoulou2023ExcludingSingleCrossing}
{\sc Giannopoulou, A.~C., Thilikos, D.~M., and Wiederrecht, S.}
\newblock Excluding single-crossing matching minors in bipartite graphs.
\newblock In {\em Proceedings of the 2023 {A}nnual {ACM}-{SIAM} {S}ymposium on
  {D}iscrete {A}lgorithms ({SODA})\/} (2023), SIAM, Philadelphia, PA,
  pp.~2111--2121.

\bibitem{giannopoulou2021flat}
{\sc Giannopoulou, A.~C., and Wiederrecht, S.}
\newblock A {F}lat {W}all {T}heorem for {M}atching {M}inors in {B}ipartite
  {G}raphs.
\newblock {\em arXiv preprint arXiv:2110.07553\/} (2021).

\bibitem{giannopoulou2021two}
{\sc Giannopoulou, A.~C., and Wiederrecht, S.}
\newblock Two disjoint alternating paths in bipartite graphs.
\newblock {\em arXiv preprint arXiv:2110.02013\/} (2021).

\bibitem{gollin2021unified}
{\sc Gollin, J.~P., Hendrey, K., Kawarabayashi, K.-i., Kwon, O.-j., and Oum,
  S.-i.}
\newblock A unified half-integral {E}rd{\H{o}}s-{P}{\'o}sa theorem for cycles
  in graphs labelled by multiple abelian groups.
\newblock {\em arXiv preprint arXiv:2102.01986\/} (2021).

\bibitem{gollin2022unified}
{\sc Gollin, J.~P., Hendrey, K., Kwon, O.-j., Oum, S.-i., and Yoo, Y.}
\newblock A unified {E}rd{\H{o}}s-{P}{\'o}sa theorem for cycles in graphs
  labelled by multiple abelian groups.
\newblock {\em arXiv preprint arXiv:2209.09488\/} (2022).

\bibitem{GueninThomas2011PackingDirectedCircuits}
{\sc Guenin, B., and Thomas, R.}
\newblock Packing directed circuits exactly.
\newblock {\em Combinatorica 31}, 4 (2011), 397--421.

\bibitem{hatzel2019polynomial}
{\sc Hatzel, M., Kawarabayashi, K.-i., and Kreutzer, S.}
\newblock {P}olynomial {P}lanar {D}irected {G}rid {T}heorem.
\newblock In {\em Proceedings of the Thirtieth Annual ACM-SIAM Symposium on
  Discrete Algorithms\/} (2019), SIAM, pp.~1465--1484.

\bibitem{Huynh2019UnifiedErdosPosaConstrained}
{\sc Huynh, T., Joos, F., and Wollan, P.}
\newblock A unified {E}rd{\H{o}}s-{P}\'{o}sa theorem for constrained cycles.
\newblock {\em Combinatorica 39}, 1 (2019), 91--133.

\bibitem{jaffke2023dynamic}
{\sc Jaffke, L., Morelle, L., Sau, I., and Thilikos, D.~M.}
\newblock Dynamic programming on bipartite tree decompositions.
\newblock {\em arXiv preprint arXiv:2309.07754\/} (2023).

\bibitem{Jansen2021VertexDeletingAndEvenLess}
{\sc Jansen, B. M.~P., de~Kroon, J. J.~H., and W\l~odarczyk, M.}
\newblock Vertex deletion parameterized by elimination distance and even less.
\newblock In {\em S{TOC} '21---{P}roceedings of the 53rd {A}nnual {ACM}
  {SIGACT} {S}ymposium on {T}heory of {C}omputing\/} ([2021] \copyright 2021),
  ACM, New York, pp.~1757--1769.

\bibitem{Johnson2001DirectedTreewidth}
{\sc Johnson, T., Robertson, N., Seymour, P.~D., and Thomas, R.}
\newblock Directed tree-width.
\newblock {\em J. Combin. Theory Ser. B 82}, 1 (2001), 138--154.

\bibitem{Joos2016ParityLinkage}
{\sc Joos, F.}
\newblock Parity linkage and the {E}rd{\H{o}}s-{P}\'{o}sa property of odd
  cycles through prescribed vertices in highly connected graphs.
\newblock In {\em Graph-theoretic concepts in computer science}, vol.~9224 of
  {\em Lecture Notes in Comput. Sci.} Springer, Berlin, 2016, pp.~339--350.

\bibitem{jung1970verallgemeinerung}
{\sc Jung, H.~A.}
\newblock {E}ine {V}erallgemeinerung des $n$-fachen {Z}usammenhangs f{\"u}r
  {G}raphen.
\newblock {\em Mathematische Annalen 187}, 2 (1970), 95--103.

\bibitem{KakimuraKawarabayashi2012PackingDirectedCircuits}
{\sc Kakimura, N., and Kawarabayashi, K.-i.}
\newblock Packing directed circuits through prescribed vertices bounded
  fractionally.
\newblock {\em SIAM J. Discrete Math. 26}, 3 (2012), 1121--1133.

\bibitem{Kakimura2012ErdosPosaParityConstratins}
{\sc Kakimura, N., Kawarabayashi, K.-i., and Kobayashi, Y.}
\newblock Erd{\H{o}}s-{P}\'{o}sa property and its algorithmic
  applications---parity constraints, subset feedback set, and subset packing.
\newblock In {\em Proceedings of the {T}wenty-{T}hird {A}nnual {ACM}-{SIAM}
  {S}ymposium on {D}iscrete {A}lgorithms\/} (2012), ACM, New York,
  pp.~1726--1736.

\bibitem{kasteleyn1967graph}
{\sc Kasteleyn, P.}
\newblock {G}raph {T}heory and {C}rystal {P}hysics.
\newblock {\em Graph Theory and Theoretical Physics\/} (1967), 43--110.

\bibitem{Kawarabayashi2012DirectedCyclesFixedVertexSet}
{\sc Kawarabayashi, K., Kr\'{a}l', D., Kr\v{c}\'{a}l, M., and Kreutzer, S.}
\newblock Packing directed cycles through a specified vertex set.
\newblock In {\em Proceedings of the {T}wenty-{F}ourth {A}nnual {ACM}-{SIAM}
  {S}ymposium on {D}iscrete {A}lgorithms\/} (2012), SIAM, Philadelphia, PA,
  pp.~365--377.

\bibitem{KawarabayashiKreutzer2015DirectedGrid}
{\sc Kawarabayashi, K.-i., and Kreutzer, S.}
\newblock The directed grid theorem.
\newblock In {\em S{TOC}'15---{P}roceedings of the 2015 {ACM} {S}ymposium on
  {T}heory of {C}omputing\/} (2015), ACM, New York, pp.~655--664.

\bibitem{Kawarabyashi2023HalfInegralOddDirectedCycles}
{\sc Kawarabayashi, K.-i., Kreutzer, S., Kwon, O.-j., and Xie, Q.}
\newblock A half-integral {E}rd{\H{o}}s-{P}\'{o}sa theorem for directed odd
  cycles.
\newblock In {\em Proceedings of the 2023 {A}nnual {ACM}-{SIAM} {S}ymposium on
  {D}iscrete {A}lgorithms ({SODA})\/} (2023), SIAM, Philadelphia, PA,
  pp.~3043--3062.

\bibitem{KawarabayashiReed2010AlmostLinearTimeOddCycleTransversal}
{\sc Kawarabayashi, K.-i., and Reed, B.}
\newblock An (almost) linear time algorithm for odd cycles transversal.
\newblock In {\em Proceedings of the {T}wenty-{F}irst {A}nnual {ACM}-{SIAM}
  {S}ymposium on {D}iscrete {A}lgorithms\/} (2010), SIAM, Philadelphia, PA,
  pp.~365--378.

\bibitem{KawarabayashiReed2010OddCyclePacking}
{\sc Kawarabayashi, K.-i., and Reed, B.}
\newblock Odd cycle packing [extended abstract].
\newblock In {\em S{TOC}'10---{P}roceedings of the 2010 {ACM} {I}nternational
  {S}ymposium on {T}heory of {C}omputing\/} (2010), ACM, New York,
  pp.~695--704.

\bibitem{kawarabayashi2020quickly}
{\sc Kawarabayashi, K.-i., Thomas, R., and Wollan, P.}
\newblock {Q}uickly {E}xcluding a {N}on-{P}lanar {G}raph.
\newblock {\em arXiv preprint arXiv:2010.12397\/} (2020).

\bibitem{LaPaughPapadimitriou1984EvenPathProblem}
{\sc LaPaugh, A.~S., and Papadimitriou, C.~H.}
\newblock The even-path problem for graphs and digraphs.
\newblock {\em Networks 14}, 4 (1984), 507--513.

\bibitem{Little1975Convertible}
{\sc Little, C. H.~C.}
\newblock A characterization of convertible (0,1)-matrices.
\newblock {\em J. Combinatorial Theory Ser. B 18\/} (1975), 187--208.

\bibitem{Lokshtanov2020ComplexityOddCycle}
{\sc Lokshtanov, D., Ramanujan, M.~S., Saurabh, S., and Zehavi, M.}
\newblock Parameterized complexity and approximability of directed odd cycle
  transversal.
\newblock In {\em Proceedings of the 2020 {ACM}-{SIAM} {S}ymposium on
  {D}iscrete {A}lgorithms\/} (2020), SIAM, Philadelphia, PA, pp.~2181--2200.

\bibitem{lovasz1987matching}
{\sc Lov{\'a}sz, L.}
\newblock {M}atching {S}tructure and the {M}atching {L}attice.
\newblock {\em Journal of Combinatorial Theory, Series B 43}, 2 (1987),
  187--222.

\bibitem{LucchesiYounger1978Minmax}
{\sc Lucchesi, C.~L., and Younger, D.~H.}
\newblock A minimax theorem for directed graphs.
\newblock {\em J. London Math. Soc. (2) 17}, 3 (1978), 369--374.

\bibitem{Mader1978Kreuzungsfrei}
{\sc Mader, W.}
\newblock \"{U}ber die {M}aximalzahl kreuzungsfreier {$H$}-{W}ege.
\newblock {\em Arch. Math. (Basel) 31}, 4 (1978/79), 387--402.

\bibitem{Masarik2019PackingDirectedCircuits}
{\sc Masa{\v r}{\'i}k, T., Muzi, I., Pilipczuk, M., Rza{\.z}ewski, P., and
  Sorge, M.}
\newblock Packing directed circuits quarter-integrally.
\newblock In {\em 27th {A}nnual {E}uropean {S}ymposium on {A}lgorithms},
  vol.~144 of {\em LIPIcs. Leibniz Int. Proc. Inform.} Schloss Dagstuhl.
  Leibniz-Zent. Inform., Wadern, 2019, pp.~Art. No. 72, 13.

\bibitem{masarik2022constant}
{\sc Masarik, T., Pilipczuk, M., Rzazewski, P., and Sorge, M.}
\newblock Constant congestion brambles in directed graphs.
\newblock {\em SIAM Journal on Discrete Mathematics 36}, 2 (2022), 922--938.

\bibitem{mccuaig2000even}
{\sc McCuaig, W.}
\newblock {E}ven {D}icycles.
\newblock {\em Journal of Graph Theory 35}, 1 (2000), 46--68.

\bibitem{McCuaig2004Polya}
{\sc McCuaig, W.}
\newblock P\'{o}lya's permanent problem.
\newblock {\em Electron. J. Combin. 11}, 1 (2004), Research Paper 79, 83.

\bibitem{mccuaig1997permanents}
{\sc McCuaig, W., Robertson, N., Seymour, P.~D., and Thomas, R.}
\newblock Permanents, pfaffian orientations, and even directed circuits.
\newblock In {\em Proceedings of the twenty-ninth annual ACM Symposium on
  Theory of Computing\/} (1997), pp.~402--405.

\bibitem{menger1927allgemeinen}
{\sc Menger, K.}
\newblock {Z}ur {A}llgemeinen {K}urventheorie.
\newblock {\em Fundamenta Mathematicae 10}, 1 (1927), 96--115.

\bibitem{polya1913aufgabe}
{\sc P{\'o}lya, G.}
\newblock Aufgabe 424.
\newblock {\em Archiv der Mathematik und Physik 20\/} (1913), 271.

\bibitem{Reed1999Introducing}
{\sc Reed, B.}
\newblock Introducing {{Directed Tree Width}}.
\newblock {\em Electronic Notes in Discrete Mathematics 3\/} (May 1999),
  222--229.

\bibitem{Reed1999Fruitsalad}
{\sc Reed, B.}
\newblock Mangoes and blueberries.
\newblock {\em Combinatorica 19}, 2 (1999), 267--296.

\bibitem{Reed1996PackingDirectedCirctuis}
{\sc Reed, B., Robertson, N., Seymour, P., and Thomas, R.}
\newblock Packing directed circuits.
\newblock {\em Combinatorica 16}, 4 (1996), 535--554.

\bibitem{Reed2004FindingOddCycleTransversals}
{\sc Reed, B., Smith, K., and Vetta, A.}
\newblock Finding odd cycle transversals.
\newblock {\em Oper. Res. Lett. 32}, 4 (2004), 299--301.

\bibitem{reed1992finding}
{\sc Reed, B.~A.}
\newblock {F}inding {A}pproximate {S}eparators and {C}omputing {T}ree {W}idth
  {Q}uickly.
\newblock In {\em Proceedings of the twenty-fourth annual ACM symposium on
  Theory of computing\/} (1992), pp.~221--228.

\bibitem{Reed1997Brambles}
{\sc Reed, B.~A.}
\newblock Tree width and tangles: a new connectivity measure and some
  applications.
\newblock In {\em Surveys in combinatorics, 1997 ({L}ondon)}, vol.~241 of {\em
  London Math. Soc. Lecture Note Ser.} Cambridge Univ. Press, Cambridge, 1997,
  pp.~87--162.

\bibitem{RobertsoSeymour1986GraphMinorsV}
{\sc Robertson, N., and Seymour, P.~D.}
\newblock Graph minors. {V}. {E}xcluding a planar graph.
\newblock {\em J. Combin. Theory Ser. B 41}, 1 (1986), 92--114.

\bibitem{robertson1990graph}
{\sc Robertson, N., and Seymour, P.~D.}
\newblock {G}raph {M}inors: {VIII}. {A} {K}uratowski {T}heorem for {G}eneral
  {S}urfaces.
\newblock {\em Journal of Combinatorial Theory, Series B 48}, 2 (1990),
  255--288.

\bibitem{robertson1991graph}
{\sc Robertson, N., and Seymour, P.~D.}
\newblock {G}raph {M}inors: {X}. {O}bstructions to {T}ree-{D}ecomposition.
\newblock {\em Journal of Combinatorial Theory, Series B 52}, 2 (1991),
  153--190.

\bibitem{GraphMinorsXIII}
{\sc Robertson, N., and Seymour, P.~D.}
\newblock Graph minors. {XIII}. {T}he disjoint paths problem.
\newblock {\em J. Combin. Theory Ser. B 63}, 1 (1995), 65--110.

\bibitem{robertson2003graph}
{\sc Robertson, N., and Seymour, P.~D.}
\newblock {G}raph minors. {XVI}. {E}xcluding a non-planar graph.
\newblock {\em Journal of Combinatorial Theory, Series B 89}, 1 (2003), 43--76.

\bibitem{Robertson1999PermanentsPfaffianOrientations}
{\sc Robertson, N., Seymour, P.~D., and Thomas, R.}
\newblock Permanents, {P}faffian orientations, and even directed circuits.
\newblock {\em Ann. of Math. (2) 150}, 3 (1999), 929--975.

\bibitem{schrijver1994finding}
{\sc Schrijver, A.}
\newblock {F}inding $k$ {D}isjoint {P}aths in a {D}irected {P}lanar {G}raph.
\newblock {\em SIAM Journal on Computing 23}, 4 (1994), 780--788.

\bibitem{seymour1980disjoint}
{\sc Seymour, P.~D.}
\newblock {D}isjoint {P}aths in {G}raphs.
\newblock {\em Discrete Mathematics 29}, 3 (1980), 293--309.

\bibitem{Seymour1995PackingDirectedCircuitsFractionally}
{\sc Seymour, P.~D.}
\newblock Packing directed circuits fractionally.
\newblock {\em Combinatorica 15}, 2 (1995), 281--288.

\bibitem{Seymour1996PackingDirectectedCircuitsEulerian}
{\sc Seymour, P.~D.}
\newblock Packing circuits in {E}ulerian digraphs.
\newblock {\em Combinatorica 16}, 2 (1996), 223--231.

\bibitem{SeymourThomassen1987EvenDirectedGraphs}
{\sc Seymour, P.~D., and Thomassen, C.}
\newblock Characterization of even directed graphs.
\newblock {\em J. Combin. Theory Ser. B 42}, 1 (1987), 36--45.

\bibitem{shiloach1980polynomial}
{\sc Shiloach, Y.}
\newblock {A} {P}olynomial {S}olution to the {U}ndirected {T}wo {P}aths
  {P}roblem.
\newblock {\em Journal of the ACM (JACM) 27}, 3 (1980), 445--456.

\bibitem{slivkins2010parameterized}
{\sc Slivkins, A.}
\newblock Parameterized tractability of edge-disjoint paths on directed acyclic
  graphs.
\newblock {\em SIAM Journal on Discrete Mathematics 24}, 1 (2010), 146--157.

\bibitem{steiner2020parameterized}
{\sc Steiner, R., and Wiederrecht, S.}
\newblock {P}arameterized {A}lgorithms for {D}irected {M}odular {W}idth.
\newblock In {\em Conference on Algorithms and Discrete Applied Mathematics\/}
  (2020), Springer, pp.~415--426.

\bibitem{thilikos2022killing}
{\sc Thilikos, D.~M., and Wiederrecht, S.}
\newblock Killing a vortex.
\newblock In {\em 2022 IEEE 63rd Annual Symposium on Foundations of Computer
  Science (FOCS)\/} (2022), IEEE, pp.~1069--1080.

\bibitem{Thomas2006PfaffianSurvey}
{\sc Thomas, R.}
\newblock A survey of {P}faffian orientations of graphs.
\newblock In {\em International {C}ongress of {M}athematicians. {V}ol. {III}}.
  Eur. Math. Soc., Z\"{u}rich, 2006, pp.~963--984.

\bibitem{ThomasYoo2023PackingCyclesUndirected}
{\sc Thomas, R., and Yoo, Y.}
\newblock Packing cycles in undirected group-labelled graphs.
\newblock {\em J. Combin. Theory Ser. B 161\/} (2023), 228--267.

\bibitem{thomassen19802}
{\sc Thomassen, C.}
\newblock 2-{L}inked {G}raphs.
\newblock {\em European Journal of Combinatorics 1}, 4 (1980), 371--378.

\bibitem{Thomassen1992EvenDicycles}
{\sc Thomassen, C.}
\newblock The even cycle problem for directed graphs.
\newblock {\em J. Amer. Math. Soc. 5}, 2 (1992), 217--229.

\bibitem{valiant1979complexity}
{\sc Valiant, L.~G.}
\newblock The complexity of computing the permanent.
\newblock {\em Theoretical Computer Science 8}, 2 (1979), 189--201.

\bibitem{VaziraniYannakakis1989PFaffian}
{\sc Vazirani, V.~V., and Yannakakis, M.}
\newblock Pfaffian orientations, {$0$}-{$1$} permanents, and even cycles in
  directed graphs.
\newblock vol.~25. 1989, pp.~179--190.
\newblock Combinatorics and complexity (Chicago, IL, 1987).

\bibitem{whitney1992congruent}
{\sc Whitney, H.}
\newblock {C}ongruent {G}raphs and the {C}onnectivity of {G}raphs.
\newblock In {\em Hassler Whitney Collected Papers}. Springer, 1992,
  pp.~61--79.

\bibitem{Wollan2010PackingAPaths}
{\sc Wollan, P.}
\newblock Packing non-zero {$A$}-paths in an undirected model of group labeled
  graphs.
\newblock {\em J. Combin. Theory Ser. B 100}, 2 (2010), 141--150.

\bibitem{Yamaguchi2019PackingAPaths}
{\sc Yamaguchi, Y.}
\newblock Packing {$A$}-paths in group-labelled graphs via linear matroid
  parity.
\newblock {\em SIAM J. Discrete Math. 30}, 1 (2016), 474--492.

\bibitem{younger1973graphs}
{\sc Younger, D.~H.}
\newblock Graphs with interlinked directed circuits.
\newblock In {\em Proceedings of the Midwest symposium on circuit theory\/}
  (1973), vol.~2, pp.~XVI--2.

\bibitem{Zaslavsky1989BiasedGraphsI}
{\sc Zaslavsky, T.}
\newblock Biased graphs. {I}. {B}ias, balance, and gains.
\newblock {\em J. Combin. Theory Ser. B 47}, 1 (1989), 32--52.

\end{thebibliography}

\end{document}